UNIVERSITÀ DEGLI STUDI DI BRESCIA
Department of Civil, Environmental, Architectural Engineering and Mathematics

PhD Program in Civil and Environmental engineering, curriculum:
MATHEMATICAL METHODS AND MODELS FOR ENGINEERING

Academic Discipline:
MAT/08 – Numerical Analysis

COURSE XXIX

# TIME INTEGRATION SCHEMES FOR FLUID-STRUCTURE INTERACTION PROBLEMS: NON-FITTED FEM S FOR IMMERSED THIN STRUCTURES

PhD Candidate
Michele Annese

Tutor
Ch.ma Prof.ssa Lucia Gastaldi

Coordinator of the PhD Program
Ch.mo Prof. Baldassarre Bacchi

# Abstract

La presente tesi affronta lo studio del problema di interazione fluido-struttura nel caso di strutture sottili immerse in un fluido incomprimibile. Il problema accoppiato si compone di un modello per il fluido, uno per il solido e condizioni di trasmissione assegnate all'interfaccia tra fluido e solido.

Partendo da una formulazione debole basata sulla tecnica dei moltiplicatori di Lagrange, viene riportata la semidiscretizzazione in tempo e viene dimostrata la buona positura degli step temporali; i risultati presentati in questo contesto sono già presenti in letteratura.

Successivamente si passa alla discretizzazione spaziale usando griglie non allineate all'interfaccia solido-fluido ed elementi finiti stabilizzati.

In questa tesi vengono studiati schemi numerici monolitici, nei quali, ad ogni passo temporale, le incognite vengono calcolate tutte insieme, e schemi numerici partizionati che consentono di affrontare separatamente, ad ogni passo temporale, la soluzione del problema del fluido e del problema nel solido.

Vengono proposti ed analizzati tre schemi numerici; il primo schema analizzato è monolitico ed è stato proposto in letteratura per coppie di spazi stabili per il problema di Stokes, qui viene estesa l'analisi di stabilità e buona positura al caso di spazi di elementi finiti stabilizzati. Inoltre viene dimostrata la convergenza nel caso linearizzato.

Allo scopo di ottenere schemi numerici partizionati sono stati proposti due ulteriori algoritmi. Entrambi gli schemi partizionati si basano sul fatto che le forze che il solido scarica sul fluido sono di due tipi: le forze di tipo inerziale, dovute alla massa del corpo solido che si muove, e le forze associate alla deformazione del solido. Il partizionamento del problema viene ottenuto trattando in forma esplicita i termini di forza legati alla deformazione e in forma implicita quelli legati all'inerzia. Uno degli algoritmi partizionati ottenuti in questo modo prevede un passo correttivo per il calcolo dello spostamento del solido ad ogni step temporale.

Per entrambi gli algoritmi proposti è stata dimostrata la stabilità ed inoltre la buona positura dei problemi associati ai passi temporali. Per l'algoritmo che prevede il passo correttivo per il solido, è stata dimostrata anche la convergenza nel caso linearizzato.

I risultati teorici ottenuti sono stati poi verificati risolvendo numericamente un caso test.

# Contents

# Introduction

The object of this thesis is the study of mechanical systems involving the interaction of an incompressible fluid with an immersed deformable thin-walled structure. This type of problem is very common in a variety of scientific and engineering contexts, as for examples: the dynamics of parachutes, paragliders or sails, the dynamics of wind and sea turbines, the biomechanics of animal cells and physiological flows (see e.g. [109, 111, 110, 14, 74, 104, 78, 3, 83, 44, 65]).

In the field of fluid-structure interaction there are two main families of methods suited for the modeling and the simulation of the coupled problem.

- The first family of methods considers the so-called *fitted methods* where the *solid mesh* and the *fluid mesh* have common nodes on the interface. In this case it is possible to transfer directly the coupling conditions between fluid and solid imposing the force equilibrium of the common nodes and zero relative velocities. These types of methods are usually based on an ALE (Arbitrary Lagrangian-Eulerian) description of the fluid dynamics [108, 77, 46]. Since the solid displacements imply the movement of the nodes on the interface, the *fitted methods* need an algorithm for the movement of the fluid grid in order to avoid large distortions of the mesh. For this reason "fitted methods" are mainly used when the solid displacements are small.

- The second family of methods is referred to as *unfitted methods* since the solid and fluid meshes are not constrained to have common nodes. This is the most convenient context in which to address problems involving large structural displacements, or problems in which different solids can come in contact or can be split up. *Unfitted methods* offer the possibility to use:

  - a. the Eulerian framework in the fluid on a fixed grid,
  - b. the Lagrangian framework for the solid using a mesh independent of the fluid grid.

  In general, these methods are known to be inaccurate in space, mainly because the space discretization does not allow discontinuities across the immersed interface. In this case, a local enrichment of the fluid finite element space can be used to avoid this drawback (see [1, 36]).

In both the fitted and unfitted methods, the solution of the fluid-structure problem involves the simulation of the fluid and of the structural dynamics and the interaction between them. The simulations can be performed using either monolithic methods, partitioned methods or splitting methods that will be described in some details in Chapter 1. Here we recall briefly the main features of the these approaches.

- In monolithic methods a single discrete problem is written taking into account both the sub-domains (the fluid and the structure) at once. The coupling conditions are generally imposed in implicit form.

- Partitioned methods allow for a separate solution of the fluid and the structure equations. The coupling conditions can be fulfilled with the desired accuracy iterating the solution of the two subsystems in a fixed point procedure.

- Splitting methods are intermediate between monolithic and partitioned methods and are generally obtained considering an explicit/implicit treatment of the kinematic and dynamic coupling conditions.

# Main contributions

The main contribution of this thesis is the analysis of three time integration schemes for unfitted methods in fluid structure interaction.

The starting point of the work is given by the results obtained in [20] and [23] about a mathematical formulation using Lagrange multiplier technique. In particular, the continuum model analyzed was proposed in [20] where was also presented an unconditionally stable monolithic algorithm for its solution. In [23] the well-posedness of the stationary problem and consequently of the time-steps of the monolithic algorithm were proved. In the cited papers the spatial discretization of the fluid problem is obtained by using a couple of "stable" finite elements in the sense of the LBB condition.

Here, we extend the stability analysis of the fully discrete monolithic algorithm and the well-posedness of the time-step problem to $\mathbb{P}_1 - \mathbb{P}_1$ stabilized finite elements for the fluid problem (see Propositions 2.4.2 and 2.4.3 respectively).

Besides the stability and well-posedness, we analyze the convergence of the space-time discretization of the linearized monolithic formulation that in the following will be denoted Algorithm 1 (see Theorem 3.3.7). The analysis gives a first order convergence rate in time, as expected since we use the Euler scheme, and convergence rates in space that are in agreement with the supposed regularity of the solution of the continuous problem.

In order to obtain partitioned schemes, we introduce and analyze two split algorithms for the coupled problem with Lagrange multipliers inspired by [60, 56, 1] where similar splitting strategies are applied to the coupled problem written using Nitsche method. In particular we introduce two algorithms that allow for a partitioning of the coupled problem by exploiting an explicit/implicit treatment of the transmission conditions.

Algorithm 2 represents, essentially, a simplification of Algorithm 1 since it simply treat the elastic forces that the solid applies to the fluid in explicit form using two expressions for the extrapolations of the solid displacement and velocities evaluated in the previous time steps.

Instead, Algorithm 3, is really a splitting algorithm that involves the solution of two staggered problems. It splits the forces that solid transfers to fluid in two contributions: the inertial contribution that is treated in implicit form and the elastic contribution that is treated in explicit form. This algorithm involves the computation of a fractional step solid velocity in the first subproblem which is corrected in the second subproblem. This scheme uses three expressions for the extrapolations of the solid displacements and velocities evaluated in the previous time steps.

We perform the stability analysis for both the schemes in Theorems 4.3.1 and 4.3.3. Algorithm 2 results conditionally stable for all the extrapolations considered, instead Algorithm 3 is unconditionally stable, for extrapolations of order zero and one, and conditionally stable for the extrapolation of order two.

Since, owing to the results of the stability analysis, Algorithm 3 is the most promising, in the same spirit of what done for the monolithic algorithm, we perform the convergence analysis in the linearized case (see Theorem 4.4.2) obtaining results in line with those of the monolithic case. In particular the splitting introduced preserves the first order convergence rate in time, as expected for the Euler scheme, and gives convergence rates with respect to the spatial mesh sizes in agreement with the supposed regularity of the continuous solution.

Moreover we report numerical simulations using the finite element library FreeFem++ in order to compare the results of the analyzed schemes.

# Thesis outline

**Chapter 1**

In this Chapter we introduce the problem under study. We start presenting the geometrical setting of the coupled problem, then we introduce the strong form of the solid and fluid problem and the transmission conditions. After that, we present the strong form af the coupled problem and we derive the variational formulation using Lagrange multipliers. In the last part of the chapter we report an analysis regarding the state of art in the context of coupling techniques and time integration schemes.

**Chapter 2**.

In this chapter we recall the results obtained in [20] and [23] concerning the continuous coupled problem and its discretization. This chapter represents the starting point of our work. In particular, we recall the well-posedness of a stationary continuous problem that gives automatically the well-posedness of the problem solved at each time step in the time semi-discrete scheme. Then, we move forward to the analysis of the fully discrete problem; in this context we recall the stability analysis of the monolithic algorithm and the well-posedness of the discrete problem solved at each time-step in the case of inf-sup stable couple of fluid finite elements. In the last part of the chapter we introduce the fully discrete problem using $\mathbb{P}_1 - \mathbb{P}_1$ stabilized fluid finite elements and give the proof of the stability and well-posedness of the time-step problem. In the end, after some heuristic considerations about the expected regularity of the solution of the continuous coupled problem, we present numerical experiments confirming the unconditional stability of the monolithic algorithm.

**Chapter 3**.

In this chapter we report the analysis of convergence of the monolithic scheme in the linearized case. This result is obtained considering the case of small displacements from the reference position of the solid (see Theorem 3.3.7). In order to analyze the convergence of the numerical scheme, we introduce and study the projection operators that map the continuous solution of the coupled system in the discrete spaces; for this analysis we exploit the results presented in Appendix C about the approximations of Stokes and solid problems. The last section of the chapter contains numerical simulations which confirm the results of the theoretical convergence analysis.

**Chapter 4**.

In this chapter we introduce and analyze two splitting schemes in order to solve the coupled problem. The design of the two algorithms is inspired by [60, 56, 1]. We show the stability of both the algorithms for all the extrapolations introduced. The stability analysis of Algorithm 2, gives conditional stability for all the extrapolations considered (see Theorem 4.3.1). Algorithm 3 is proved to be unconditionally stable, for extrapolations of order zero and one, and conditionally stable for extrapolation of order two (see Theorem 4.3.3). Then we discuss the convergence analysis of Algorithm 3, in the linearized case, performed in Theorem 4.4.2. The chapter ends with numerical simulations designed in order to verify the stability and accuracy property of Algorithm 3.

**Appendix A**.

In this Appendix, we present some standard tools of functional analysis and finite element discretization. In particular we collect here basic results about variational problems and their finite dimensional approximation.

**Appendix B**.

In this Appendix we recall basic concepts of continuum mechanics; particularly the continuum kinematic and we discuss the conservation laws. Then, we introduce the constitutive relations for Newtonian fluids, recalling the Navier-Stokes equations and hyper-elastic solids. Since the main topic of this thesis is the study of the interaction between fluids and thin-structures, we introduce appropriate structural models; in particular we recall the weak formulation of general shell theory discussing the mathematical properties of the models.

**Appendix C**.

In this appendix we analyze the static solid problem and the Stokes problem. In particular we give results regarding the finite element approximation of the solid and fluid problems. The outcomes of this appendix are useful for the analysis of the projection operators used in the study of the space-time convergence of the coupled schemes.

**Appendix D**.

In this final Appendix we report the matrix computations needed for the implementation of the numerical codes.

# Chapter 1

# The fluid-structure interaction problem

## 1.1 Introduction

In this chapter, we present the fluid-structure interaction problem in the case of thin elastic solids immersed in incompressible fluids. The aim here is to introduce the solid and fluid models and the coupling conditions. Then, in the spirit of interface problems (see [27, 73, 50, 16]), we rewrite the coupled system using the Lagrange multiplier approach.

## 1.2 Problem setting

The model problem we will be working on is sketched in Figure 1.1; from the geometrical point of view the fluid and the solid fill a fixed domain $\Omega \subset \mathbb{R}^d$, $(d = 2,3)$, with boundary $\Gamma$.

At each time $t \in [0,T]$, the *solid domain* is represented by an oriented (d-1)-manifold $\boldsymbol{\phi}(\cdot,t) : \Sigma \to \Sigma(t) \subset \Omega$, with positive normal direction field $\boldsymbol{n}_{\Sigma(t)}$. We introduce a fixed (d-1)-dimensional set $\Sigma$ that represents the reference configuration of the solid and it is represented in Figure 1.1 by the dashed lines with normal vector field $\boldsymbol{n}_{\boldsymbol{\Sigma}}$. Since we consider solid bodies that are represented by curves in 2D or surfaces in 3D, they could be closed, that is without boundary. When the solid body has a boundary we impose boundary conditions for the displacement field and the traction field. In this thesis we consider only the two types of boundary conditions that are sketched in Figure 1.1, namely the supported edges "a" and "b", where we impose the homogeneous Dirichlet boundary condition on the displacement field, and the unloaded edges "c" and "d" where a zero boundary condition on the traction is imposed.

The *fluid domain* is denoted by $\Omega(t) = \Omega \setminus \Sigma(t)$ and its boundary is $\partial\Omega(t) = \Gamma \cup \Sigma(t)$.

## 1.3 The solid problem

Concerning the solid models for thin linear elastic structures, we assume that the solid behaves as a membrane, i.e., a thin elastic shell for which are valid the following claims:

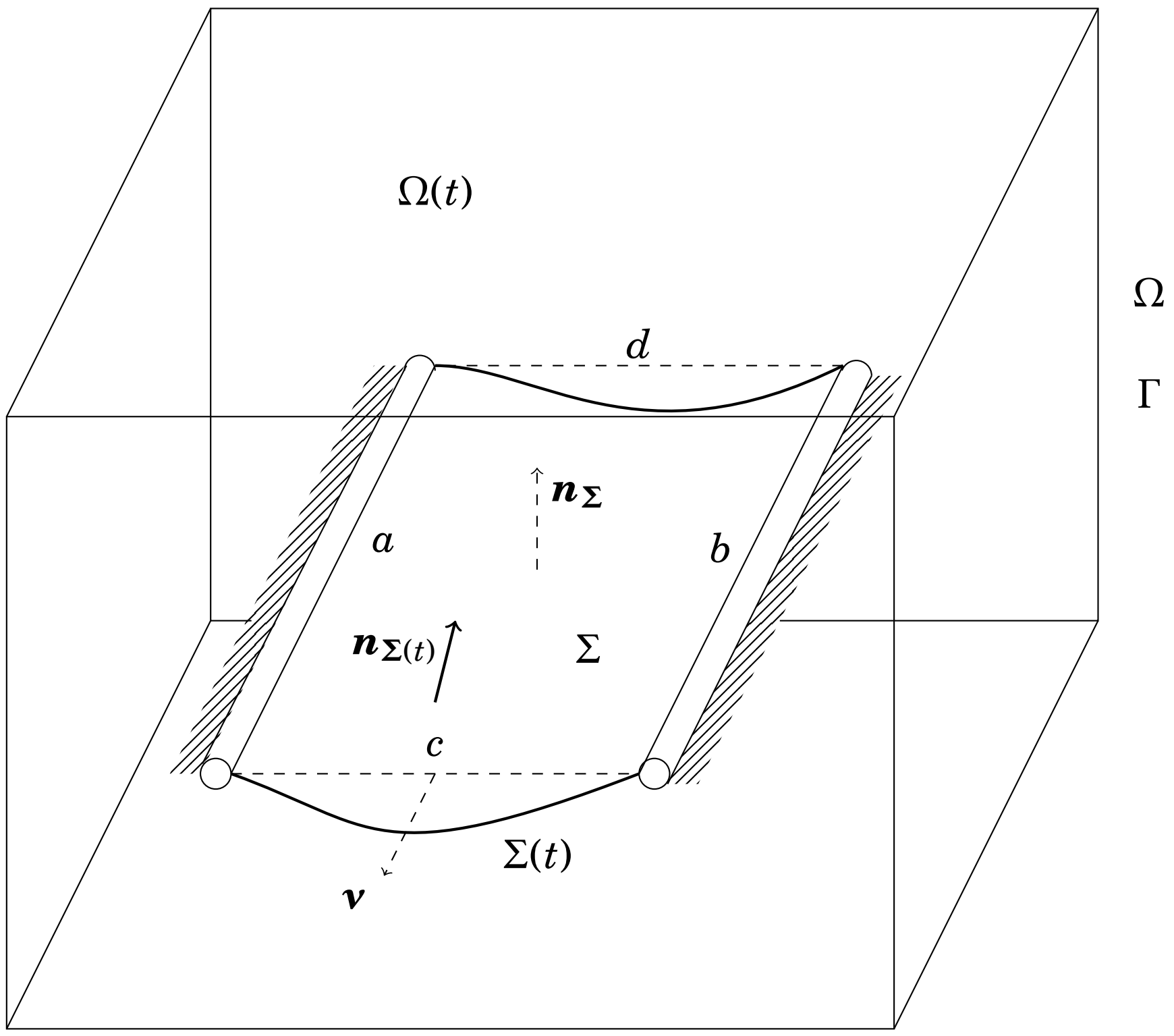

Figure 1.1: Coupled problem

1. linear constitutive stress-strain relation and isotropic-homogeneous material,
2. small deformations,
3. negligible bending terms (membrane deformation),
4. sections normal to the reference configuration $\Sigma$ remain normal to the deformed configuration (Kirchhoff-Love assumption),
5. displacements normal to the reference configuration.

Using these assumptions, it is possible to derive solid models which are widely used in the applications; the derivation of some of the named models is presented in the Appendix B. For the purpose of the theoretical analysis we consider an abstract solid differential operator $\boldsymbol{L}^s : \boldsymbol{D}^s \subset L^2(\Sigma)^d \to L^2(\Sigma)^d$ that takes as argument the solid displacements from the reference configuration, denoted by $\boldsymbol{d}$. If the solid domain has a boundary $\partial\Sigma$, we denote by $\partial\Sigma_t$ the part of the boundary on which it is imposed the traction, while by $\partial\Sigma_d$, the part of boundary on which are imposed Dirichlet boundary conditions for $\boldsymbol{d}$. Moreover we assume that $\partial\Sigma$ is the union of the two disjoint parts $\partial\Sigma_t$ and $\partial\Sigma_d$.

From this abstract point of view, the strong formulation of the elastodynamics for a solid of thickness $\epsilon$ and surface mass density $\rho_{\rm s}$, can be written as (see Appendix B for more

details)

$$\begin{aligned}
&\rho_s \epsilon \partial_{tt} \boldsymbol{d} + \boldsymbol{L}^s \boldsymbol{d} = \boldsymbol{f} && \text{in} \quad \Sigma \times [0,T], \\
&\boldsymbol{d} = \boldsymbol{0} && \text{on} \quad \partial\Sigma_d \times [0,T], \\
&\boldsymbol{B}_s \boldsymbol{d} = \boldsymbol{0} && \text{on} \quad \partial\Sigma_t \times [0,T], \\
&\boldsymbol{d}(0) = \boldsymbol{d}_0 && \text{in} \quad \Sigma, \\
&\partial_t \boldsymbol{d}(0) = \boldsymbol{d}_1 && \text{in} \quad \Sigma,
\end{aligned} \tag{1.1}$$

where the actual expression of the elastic operator $\boldsymbol{L}^s$ depends on the material constants and the geometry of the solid domain $\Sigma$ . The vector field $\boldsymbol{f}$ represents the externally imposed load on $\Sigma$ and $\boldsymbol{B}_s \boldsymbol{d}$ represents the traction on $\partial\Sigma_t$ expressed in terms of an operator $\boldsymbol{B}_s$ whose actual form depends on the operator $\boldsymbol{L}^s$. We recall that in the following we will consider always zero traction on $\partial\Sigma_t$ and zero displacements on $\partial\Sigma_d$.

## 1.4 Fluid problem and coupled problem

Let us consider a Newtonian incompressible fluid flowing, at time $t$, in the deformable domain $\Omega(t)$ as specified in the problem setting (see Figure 1.1). Assuming a low-speed regime, the evolution of the fluid is well described using the linearized Navier-Stokes equations (B.45). Denoting by $\boldsymbol{u}(\cdot,t) : \Omega(t) \to \mathbb{R}^d$ the velocity field and by $p(\cdot,t) : \Omega(t) \to \mathbb{R}$ the pressure field, we have

$$\begin{aligned}
&\rho_f \partial_t \boldsymbol{u} - \operatorname{div} \boldsymbol{\sigma}(\boldsymbol{u}, p) = \boldsymbol{0} && \text{in} \quad \Omega(t), \\
&\operatorname{div} \boldsymbol{u} = 0 && \text{in} \quad \Omega(t), \\
&\boldsymbol{\sigma}(\boldsymbol{u}, p) = -p\boldsymbol{I} + 2\mu \boldsymbol{\epsilon}(\boldsymbol{u}) && \text{in} \quad \Omega(t), \\
&\boldsymbol{\epsilon}(\boldsymbol{u}) = \frac{1}{2}\left(\nabla \boldsymbol{u} + \nabla \boldsymbol{u}^T\right) && \text{in} \quad \Omega(t), \\
&\boldsymbol{u} = \boldsymbol{0} && \text{on} \quad \Gamma \times [0,T], \\
&\boldsymbol{u}(0) = \boldsymbol{u}_0 && \text{in} \quad \Omega \backslash \Sigma(0),
\end{aligned} \tag{1.2}$$

where $\rho_f$ and $\mu$ are positive constants that stand for the density and the dynamic viscosity of the fluid, respectively. We point out that the evolution of $\Omega(t)$ is not known a priori since it is determined by the interaction between the fluid and the structure. The mathematical problem is therefore highly coupled and consists of finding the velocity and the pressure of the fluid, the displacement of the structure and, therefore, the current domains $\Omega(t)$ and $\Sigma(t)$.

### 1.4.1 Coupling conditions

In order to couple the fluid and the solid problems we assign the transmission conditions, that have to be satisfied at the interface. In general the fluid and the solid have different settings, in fact usually, the framework for the fluid is the Eulerian one, while for the solid we use the Lagrangian framework; since we write the interface conditions in the Lagrangian framework, we introduce, also for the fluid, a reference domain $\Omega_0$ which contains the reference configuration of the solid $\Sigma$. Associated to $\Omega_0$ we have also a regular mapping

$$A_t : \Omega_0 \to \Omega \tag{1.3}$$

which transforms $\Omega_0\backslash\Sigma$ into the actual configuration $\Omega(t)$ and such that $A_t(\Sigma) = \boldsymbol{\phi}(\Sigma, t) = \Sigma(t)$. Using this map we can write the fluid velocity and the fluid Cauchy stress tensor in the reference domain $\Omega_0$ as follows

$$\hat{\boldsymbol{u}} = \boldsymbol{u} \circ A_t, \qquad \hat{\boldsymbol{\sigma}}_{\mathrm{f}} = \boldsymbol{\sigma}_{\mathrm{f}} \circ A_t, \qquad \boldsymbol{P}_{\mathrm{f}} = J\hat{\boldsymbol{\sigma}}_{\mathrm{f}}\boldsymbol{F}^{-T}, \tag{1.4}$$

where we used the notations $J = \det(\nabla A_{\mathrm{t}})$ and $\boldsymbol{F} = \nabla \boldsymbol{A_t}$. The tensor $\boldsymbol{P}_{\mathrm{f}}$ is the *first Piola-Kirchhoff Stress Tensor* of the fluid. Using these quantities we can write properly the coupling conditions that are normally imposed in this case. We refer to the case in which the domain $\Omega$ is split into two disjoint parts $\Omega^+$ and $\Omega^-$ by the solid. If this is not the case, as illustrated in Figure 1.2, we can extend the solid domain up to $\partial\Omega$ so that it ideally splits the domain $\Omega$. The solid $\Sigma$ is then denoted as $\Sigma^+$ if it is considered as part of $\partial\Omega^+$ or $\Sigma^-$ if it is seen as part of $\partial\Omega^-$.

- Kinematic condition: the fluid and the solid velocities have to coincide at the interface

$$\hat{\boldsymbol{u}}\big|_{\Sigma^-} = \hat{\boldsymbol{u}}\big|_{\Sigma^+} \stackrel{def}{=} \hat{\boldsymbol{u}}\big|_{\Sigma} = \partial_t \boldsymbol{d} \qquad \text{on} \qquad \Sigma \times [0, T]. \tag{1.5}$$

- Equilibrium condition: the tractions on the interface have to be equilibrated

$$\boldsymbol{P}_{\mathrm{f}}^{+}\boldsymbol{n}_{\boldsymbol{\Sigma}} - \boldsymbol{P}_{\mathrm{f}}^{-}\boldsymbol{n}_{\boldsymbol{\Sigma}} + \boldsymbol{f} = [\![\boldsymbol{P}_{\mathrm{f}}\boldsymbol{n}_{\boldsymbol{\Sigma}}]\!] + \boldsymbol{f} = \boldsymbol{0} \qquad \text{on} \qquad \Sigma \times [0, T], \tag{1.6}$$

where $\boldsymbol{P}_{\mathrm{f}}^{+}$ and $\boldsymbol{P}_{\mathrm{f}}^{-}$ represent the restriction of the Piola-Kirchhoff Stress Tensor of the fluid to $\Sigma$ seen as part of $\partial\Omega^+$ or $\partial\Omega^-$, respectively. The symbol $[\![\cdot]\!]$ represents the jump operator and it is defined by equation (1.6) itself.

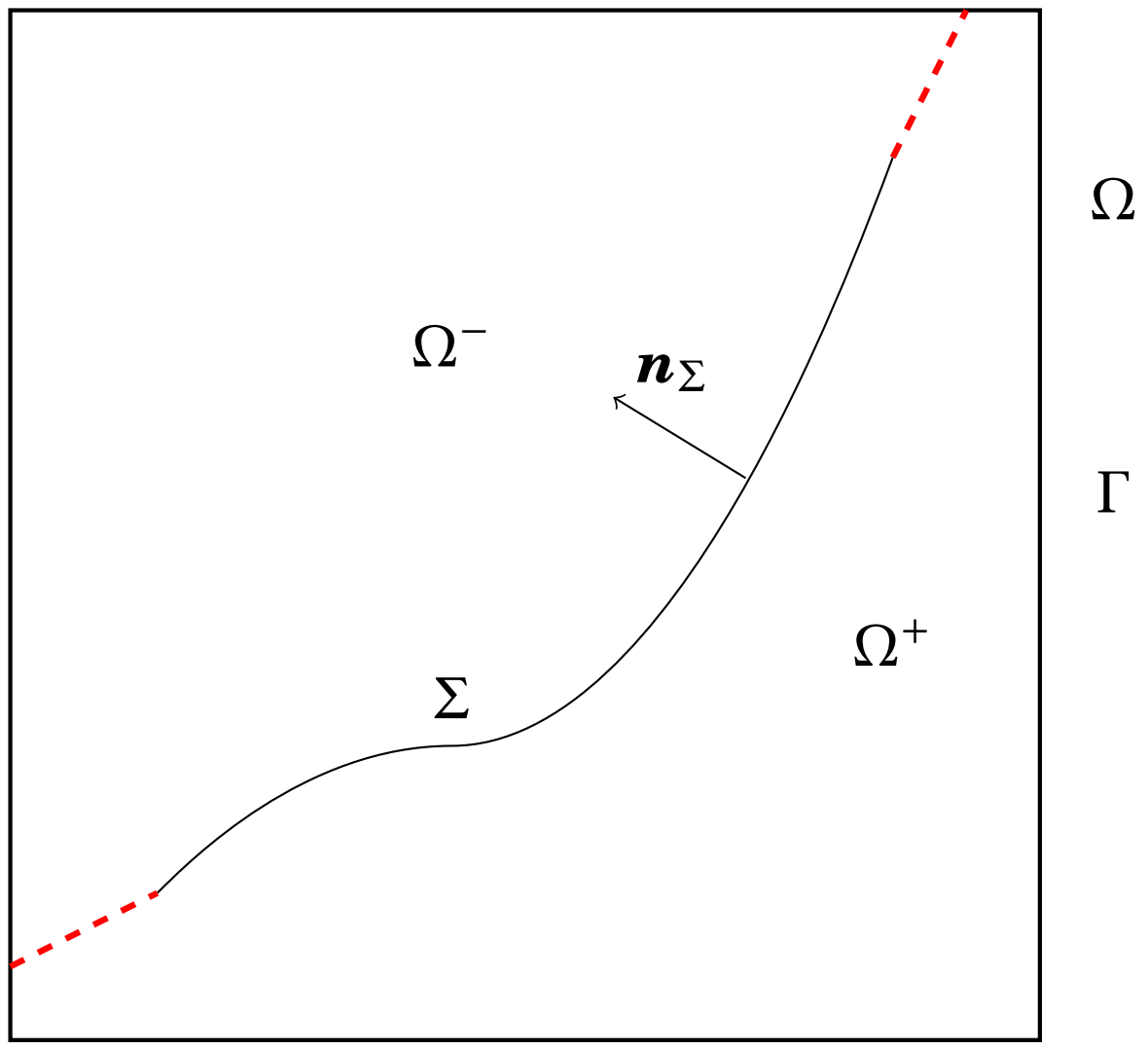

Figure 1.2: Coupled problem

### 1.4.2 Coupled problem

We are lead to consider the coupled problem in the whole domain $\Omega$ seen as a partitioned domain in which we have the solid and the fluid.

**Problem 1.** *Given* $\boldsymbol{u}_0, \boldsymbol{d}_0, \boldsymbol{d}_1$ *such that* $\boldsymbol{u}_0\big|_\Sigma = \boldsymbol{d}_1$ *, find* $\boldsymbol{u}, p, \boldsymbol{d}$ *such that*

- ***Fluid problem***

$$
\begin{aligned}
&\rho_f \partial_t \boldsymbol{u} - \mathrm{div}\boldsymbol{\sigma}_f(\boldsymbol{u}, p) = \boldsymbol{0} && in \quad \Omega(t)\\
&\mathrm{div}\,\boldsymbol{u} = 0 && in \quad \Omega(t)\\
&\boldsymbol{\sigma}_f(\boldsymbol{u}, p) = -p\boldsymbol{I} + 2\mu\boldsymbol{\epsilon}(\boldsymbol{u}) && in \quad \Omega(t)\\
&\boldsymbol{\epsilon}(\boldsymbol{u}) = \frac{1}{2}\left(\nabla\boldsymbol{u} + \nabla\boldsymbol{u}^T\right) && in \quad \Omega(t)\\
&\boldsymbol{u} = \boldsymbol{0} && on \quad \Gamma \times [0, T]\\
&\boldsymbol{u}(0) = \boldsymbol{u}_0 && in \quad \Omega\backslash\Sigma
\end{aligned}
\tag{1.7}
$$

- ***Solid problem***

$$
\begin{aligned}
&\rho_s \epsilon \partial_{tt}\boldsymbol{d} + \boldsymbol{L}^s\boldsymbol{d} = \boldsymbol{f} && in \quad \Sigma \times [0, T]\\
&\boldsymbol{d} = \boldsymbol{0} && on \quad \partial\Sigma_d \times [0, T]\\
&\boldsymbol{B}_s\boldsymbol{d} = \boldsymbol{0} && on \quad \partial\Sigma_t \times [0, T]\\
&\boldsymbol{d}(0) = \boldsymbol{d}_0 && in \quad \Sigma\\
&\partial_t\boldsymbol{d}(0) = \boldsymbol{d}_1 && in \quad \Sigma
\end{aligned}
\tag{1.8}
$$

- ***Coupling conditions***

$$
\begin{aligned}
&\boldsymbol{u} \circ \boldsymbol{\phi} = \partial_t\boldsymbol{d} && on \quad \Sigma \times [0, T],\\
&[\![\boldsymbol{P}_f\boldsymbol{n}_\Sigma]\!] + \boldsymbol{f} = \boldsymbol{0} && on \quad \Sigma \times [0, T]
\end{aligned}
\tag{1.9}
$$

where the actual configuration of the solid $\Sigma(t)$, is given by $\boldsymbol{\phi} = \boldsymbol{\phi}_0 + \boldsymbol{d}$, with $\boldsymbol{\phi}_0$ being the position vector of the solid reference configuration.

## 1.5 Weak form with Lagrange multipliers

In this thesis, we introduce a weak formulation of the Problem 1 involving a suitable Lagrange multiplier. This technique is widely used to couple different problems and its introduction in the finite elements goes back to the early 70' as a method to impose the Dirichlet boundary conditions in weak form [5]. In the following we proceed formally multiplying the equations by regular test functions and then integrating by parts.

- **Fluid**: let $\boldsymbol{v} : \Omega \to \mathbb{R}^d$ and $q : \Omega \to \mathbb{R}$ be regular test functions defined on $\Omega = \Omega(t) \cup \Sigma(t)$ with $\boldsymbol{v}\big|_\Gamma = \boldsymbol{0}$, we multiply the equation $(1.7)_1$ by $\boldsymbol{v}$ and integrate by parts on the two subdomains $\Omega^+(t)$ and $\Omega^-(t)$, then we multiply the equation $(1.7)_2$ by $q$ and integrate. Summing the resulting equations we obtain

$$
\begin{aligned}
&\int_{\Omega(t)} \big(\rho_f \partial_t \boldsymbol{u} - \operatorname{div}\boldsymbol{\sigma}_f(\boldsymbol{u},p)\big)\cdot \boldsymbol{v} + \int_{\Omega(t)} q \operatorname{div}\boldsymbol{u} \\
&\qquad = \rho_f \int_\Omega \partial_t \boldsymbol{u}\cdot\boldsymbol{v} + \int_\Omega \boldsymbol{\sigma}_f(\boldsymbol{u},p):\nabla\boldsymbol{v} + \int_\Omega q \operatorname{div}\boldsymbol{u} - \int_{\Sigma(t)} [\![\boldsymbol{\sigma}_f(\boldsymbol{u},p)\boldsymbol{n}_{\Sigma(t)}]\!]\cdot\boldsymbol{v} \\
&\qquad = \rho_f \int_\Omega \partial_t \boldsymbol{u}\cdot\boldsymbol{v} + \int_\Omega \boldsymbol{\sigma}_f(\boldsymbol{u},p):\nabla\boldsymbol{v} + \int_\Omega q \operatorname{div}\boldsymbol{u} - \int_{\Sigma} [\![\boldsymbol{P}_f(\boldsymbol{u},p)\boldsymbol{n}_{\Sigma}]\!]\cdot\boldsymbol{v}\circ\boldsymbol{\phi} = 0, \quad (1.10)
\end{aligned}
$$

where we used (1.4) to transform for the fluid stress tensor $\boldsymbol{\sigma}_f$ into $\boldsymbol{P}_f$ (see Appendix B section B.2.7).

- **Solid**: we denote by $\boldsymbol{W}$ the space of admissible solid displacements, namely the space of displacements that fulfill the Dirichlet boundary conditions; since we consider only zero Dirichlet boundary conditions $\boldsymbol{w} = \boldsymbol{0}$ on $\partial\Sigma_d$, the space of admissible displacements coincide with the space of virtual displacements. The weak formulation associated to the solid problem is assumed to be given in terms of a bilinear form $a^s : \boldsymbol{W}\times\boldsymbol{W} \to \mathbb{R}$ defined as

$$
a^s(\boldsymbol{d},\boldsymbol{w}) \overset{def}{=} \langle \boldsymbol{L}^s\boldsymbol{d},\boldsymbol{w}\rangle = \langle \boldsymbol{d},\boldsymbol{L}^s\boldsymbol{w}\rangle \qquad \forall \boldsymbol{d},\boldsymbol{w}\in\boldsymbol{W}. \quad (1.11)
$$

Then the weak formulation for the solid is written as

$$
\rho_s\epsilon \int_\Sigma \partial_{tt}\boldsymbol{d}\cdot\boldsymbol{w} + a^s(\boldsymbol{d},\boldsymbol{w}) = \int_\Sigma \boldsymbol{f}\cdot\boldsymbol{w}. \quad (1.12)
$$

Summing the fluid and the solid weak equations and adding and subtracting the quantity $\int_\Sigma \boldsymbol{f}\cdot\boldsymbol{v}\circ\boldsymbol{\phi}$, we obtain the following equation

$$
\begin{aligned}
\rho_f \int_\Omega \partial_t \boldsymbol{u}\cdot\boldsymbol{v} + \int_\Omega \boldsymbol{\sigma}_f(\boldsymbol{u},p):\nabla\boldsymbol{v} + \int_\Omega \operatorname{div}\boldsymbol{u}\, q + \rho_s\epsilon\int_\Sigma \partial_{tt}\boldsymbol{d}\cdot\boldsymbol{w} + a^s(\boldsymbol{d},\boldsymbol{w}) \\
+ \int_\Sigma \boldsymbol{f}\cdot\big(\boldsymbol{v}\circ\boldsymbol{\phi} - \boldsymbol{w}\big) - \int_\Sigma \big([\![\boldsymbol{P}_f(\boldsymbol{u},p)\boldsymbol{n}_\Sigma]\!] + \boldsymbol{f}\big)\cdot\boldsymbol{v}\circ\boldsymbol{\phi} = 0.
\end{aligned}
$$

Using the coupling condition $(1.9)_2$ we can write the following **constrained problem**:

find $\boldsymbol{u},\boldsymbol{d},p$ with $\boldsymbol{u}\circ\boldsymbol{\phi} = \partial_t\boldsymbol{d}$ such that for all test functions $\boldsymbol{v},\boldsymbol{w},q$ with $\boldsymbol{v}\circ\boldsymbol{\phi} = \boldsymbol{w}$ it holds true

$$
\rho_f \int_\Omega \partial_t \boldsymbol{u}\cdot\boldsymbol{v} + \int_\Omega \boldsymbol{\sigma}_f(\boldsymbol{u},p):\nabla\boldsymbol{v} + \int_\Omega \operatorname{div}\boldsymbol{u}\, q + \rho_s\epsilon\int_\Sigma \partial_{tt}\boldsymbol{d}\cdot\boldsymbol{w} + a^s(\boldsymbol{d},\boldsymbol{w}) = 0.
$$

In order to eliminate the constraint and impose the coupling condition $(1.9)_1$, enforcing that the velocities of the solid and the fluid are equal on $\Sigma(t)$, we introduce a Lagrange multiplier $\boldsymbol{\lambda}$ that represents the interaction force between the fluid and the solid. Finally, the **weak formulation with Lagrange multiplier** is written as: find $\boldsymbol{u},\boldsymbol{d},\boldsymbol{\lambda},p$ such that for all test functions $\boldsymbol{v},\boldsymbol{w},\boldsymbol{\mu},q$

$$
\begin{aligned}
&\rho_f \int_\Omega \partial_t \boldsymbol{u}\cdot\boldsymbol{v} + \int_\Omega \boldsymbol{\sigma}_f(\boldsymbol{u},p):\nabla\boldsymbol{v} + \int_\Omega \operatorname{div}\boldsymbol{u}\, q + \rho_s\epsilon\int_\Sigma \partial_{tt}\boldsymbol{d}\cdot\boldsymbol{w} + a^s(\boldsymbol{d},\boldsymbol{w}) + \int_\Sigma \boldsymbol{\lambda}\cdot(\boldsymbol{v}\circ\boldsymbol{\phi} - \boldsymbol{w}) = 0 \\
&\int_\Sigma \boldsymbol{\mu}\cdot(\boldsymbol{u}\circ\boldsymbol{\phi} - \partial_t\boldsymbol{d}) = 0.
\end{aligned} \quad (1.13)
$$

In order to make rigorous the weak formulation, we introduce the following function spaces.

- The **space of fluid velocities**: $\boldsymbol{V} = H_0^1(\Omega)^d$.
- The **space of fluid pressures**: $Q = L_0^2(\Omega) = \left\{ q \in L^2(\Omega) : \frac{1}{|\Omega|} \int_\Omega q\,d\boldsymbol{x} = 0 \right\}$.
- The **space of solid displacements**: $\boldsymbol{W} \subset H^1(\Sigma)^d$, where the proper choice of the functional space $\boldsymbol{W}$ depends on the boundary conditions imposed on the solid displacement.

We observe that the functions $\boldsymbol{v} \in \boldsymbol{V}$ are continuous across the structure and that the traces of that functions have the following properties (see Appendix A)

- the trace of $\boldsymbol{v}$ on $\Sigma(t)$, devoted by $\boldsymbol{v}|_{\Sigma(t)}$, belongs to $H^{\frac{1}{2}}(\Sigma(t))^d$,
- assuming that the map $\boldsymbol{A}_t : \Omega_0 \to \Omega(t)$, introduced in (1.3), is at least Lipschitz for all $t \in [0,T]$, we have

$$\boldsymbol{v} \circ \boldsymbol{\phi} = \hat{\boldsymbol{v}} \in H^{\frac{1}{2}}(\Sigma)^d. \tag{1.14}$$

In order to give meaning to equation $(1.13)_2$, we introduce the **space** $\boldsymbol{\Lambda} = \left(H^{\frac{1}{2}}(\Sigma)^d\right)'$, that is, the dual space of $H^{\frac{1}{2}}(\Sigma)^d$. With this choice, the integrals over $\Sigma$ are regarded as the duality pairing between $H^{\frac{1}{2}}(\Sigma)^d$ and $\boldsymbol{\Lambda}$ that will be denoted as

$$c(\boldsymbol{\mu},\boldsymbol{\eta}) \stackrel{def}{=} {}_{\boldsymbol{\Lambda}}\langle \boldsymbol{\mu},\boldsymbol{\eta}\rangle_{H^{\frac{1}{2}}(\Sigma)^2} \qquad \forall \boldsymbol{\mu} \in \boldsymbol{\Lambda}, \forall \boldsymbol{\eta} \in H^{\frac{1}{2}}(\Sigma)^d. \tag{1.15}$$

Owing to its definition, the biliner form $c(\cdot,\cdot)$ has the following properties

$$c(\boldsymbol{\mu},\boldsymbol{\eta}) \le \|\boldsymbol{\mu}\|_{\boldsymbol{\Lambda}} \|\boldsymbol{\eta}\|_{\frac{1}{2},\Sigma} \qquad \forall \boldsymbol{\mu} \in \boldsymbol{\Lambda}, \forall \boldsymbol{\eta} \in H^{\frac{1}{2}}(\Sigma)^d, \tag{1.16}$$

$$c(\boldsymbol{\mu},\boldsymbol{\eta}) = 0 \qquad \forall \boldsymbol{\mu} \in \boldsymbol{\Lambda} \qquad \Rightarrow \qquad \boldsymbol{\eta} = \boldsymbol{0} \quad \text{in} \quad H^{\frac{1}{2}}(\Sigma)^d. \tag{1.17}$$

We introduce also the following bilinear form

$$\begin{aligned} a^{\mathrm{f}}\big((\boldsymbol{u},p),(\boldsymbol{v},q)\big) &\stackrel{def}{=} \int_\Omega \boldsymbol{\sigma}_{\mathrm{f}}(\boldsymbol{u},p) : \nabla \boldsymbol{v} + \int_\Omega q \operatorname{div} \boldsymbol{u} \\ &= 2\mu\big(\boldsymbol{\epsilon}(\boldsymbol{u}),\boldsymbol{\epsilon}(\boldsymbol{v})\big)_{0,\Omega} - (p,\operatorname{div}\boldsymbol{v})_{0,\Omega} + (q,\operatorname{div}\boldsymbol{u})_{0,\Omega} \qquad \forall (\boldsymbol{u},p),(\boldsymbol{v},q) \in \boldsymbol{V} \times Q, \end{aligned} \tag{1.18}$$

where the second equality is obtained considering the constitutive relation for Newtonian fluids $(1.7)_3$. Then the resulting weak problem with Lagrange multiplier is:

**Problem 2.** *Let $\boldsymbol{\phi}_0$ be the position vector of the solid reference configuration and $\boldsymbol{u}_0, \boldsymbol{d}_0, \boldsymbol{d}_1$ with $\boldsymbol{u}_0 \circ \boldsymbol{\phi}_0 = \boldsymbol{d}_1$ the initial data, then, for $t > 0$, find $\big(\boldsymbol{u}(t), p(t), \boldsymbol{d}(t), \boldsymbol{\lambda}(t)\big) \in \boldsymbol{V} \times Q \times \boldsymbol{W} \times \boldsymbol{\Lambda}$ such that*

- $\boldsymbol{\phi} = \boldsymbol{\phi}_0 + \boldsymbol{d}, \qquad$ *on* $\Sigma$

- *for all* $(\boldsymbol{v}, q, \boldsymbol{w}, \boldsymbol{\mu}) \in \boldsymbol{V} \times Q \times \boldsymbol{W} \times \boldsymbol{\Lambda}$

$$\rho_{\mathrm{f}}\big(\partial_t \boldsymbol{u}, \boldsymbol{v})_{0,\Omega} + a^{\mathrm{f}}((\boldsymbol{u},p),(\boldsymbol{v},q)\big) + c\big(\boldsymbol{\lambda}, \boldsymbol{v}\circ\boldsymbol{\phi} - \boldsymbol{w}\big) \\ - c\big(\boldsymbol{\mu}, \boldsymbol{u}\circ\boldsymbol{\phi} - \partial_t \boldsymbol{d}\big) + \rho_{\mathrm{s}}\epsilon\big(\partial_{tt}\boldsymbol{d}, \boldsymbol{w}\big)_{0,\Sigma} + a^{\mathrm{s}}\big(\boldsymbol{d}, \boldsymbol{w}\big) = 0, \quad (1.19)$$

- $\boldsymbol{u}(0) = \boldsymbol{u}_0$ *in* $\Omega$, $\partial_t \boldsymbol{d}(0) = \boldsymbol{d}_1$, $\boldsymbol{d}(0) = \boldsymbol{d}_0$ *in* $\Sigma$,

where $(\cdot,\cdot)_{0,\Omega}$ and $(\cdot,\cdot)_{0,\Sigma}$ are the $L^2$-scalar product in $L^2(\Omega)^d$ and $L^2(\Sigma)^d$.

### 1.5.1 Mathematical properties of the fluid and the solid problems

The analysis of stability, well-posedness and convergence of the numerical algorithms that we will perform in this thesis, is heavily founded on the mathematical properties of the bilinear forms $a^{\mathrm{s}}(\cdot,\cdot)$ for the solid and $a^{\mathrm{f}}(\cdot,\cdot)$ for the fluid.

**Properties of $a^{\mathrm{s}}: \boldsymbol{W}\times\boldsymbol{W} \to \mathbb{R}$.** The bilinear form $a^{\mathrm{s}}(\cdot,\cdot)$ is assumed **symmetric, continuous and coercive** on the space $\boldsymbol{W}\times\boldsymbol{W}$, then there exist two constants $\beta_{\mathrm{c}} > 0$ and $\beta_{\mathrm{s}} > 0$ for which it holds

$$\begin{aligned} \forall \boldsymbol{d}, \boldsymbol{w} \in \boldsymbol{W}, \qquad & a^{\mathrm{s}}(\boldsymbol{d},\boldsymbol{w}) \le \beta_{\mathrm{s}} \|\boldsymbol{d}\|_{1,\Sigma} \|\boldsymbol{w}\|_{1,\Sigma}. \\ \forall \boldsymbol{w} \in \boldsymbol{W}, \qquad & a^{\mathrm{s}}(\boldsymbol{w},\boldsymbol{w}) \ge \beta_{\mathrm{c}} \|\boldsymbol{w}\|^2_{1,\Sigma} \end{aligned} \quad (1.20)$$

From this assumption we can infer that $a^{\mathrm{s}}(\cdot,\cdot)$ defines an inner product on the space of admissible displacement $\boldsymbol{W}$, that becomes an Hilbert space with the norm $\|\boldsymbol{w}\|_{\mathrm{s}} \stackrel{def}{=} a^{\mathrm{s}}(\boldsymbol{w},\boldsymbol{w})^{\frac{1}{2}}$ which is equivalent to the $H^1$ norm. Moreover we assume that $a^{\mathrm{s}}(\cdot,\cdot)$ commutes with the time derivative operator, that is, for all $\boldsymbol{w}(t) \in \boldsymbol{W}$

$$\frac{d}{dt} a^{\mathrm{s}}(\boldsymbol{w},\boldsymbol{w}) = a^{\mathrm{s}}(\partial_t \boldsymbol{w}, \boldsymbol{w}) + a^{\mathrm{s}}(\boldsymbol{w}, \partial_t \boldsymbol{w}) = 2a^{\mathrm{s}}(\partial_t \boldsymbol{w}, \boldsymbol{w}). \quad (1.21)$$

**Well-posedness of the solid problem** We recall here some properties of the solution of the solid problem associated with the bilinear form $a^{\mathrm{s}}(\cdot,\cdot)$

**Proposition 1.5.1.** *Let* $a^{\mathrm{s}}(\cdot,\cdot)$ *satisfy* (1.20) *then, for all* $\boldsymbol{f} \in \boldsymbol{W}'$ *there exists a unique* $\boldsymbol{d} \in \boldsymbol{W}$ *such that*

$$a^{\mathrm{s}}(\boldsymbol{d},\boldsymbol{w}) = \langle \boldsymbol{f}, \boldsymbol{w}\rangle \qquad \forall \boldsymbol{w} \in \boldsymbol{W}, \quad (1.22)$$

*and there exists* $C > 0$ *independent of* $\boldsymbol{f}$*, such that*

$$\|\boldsymbol{d}\|_{1,\Sigma} \le C \|\boldsymbol{f}\|_{\boldsymbol{W}'}. \quad (1.23)$$

*Proof.* The proof of the proposition is obtained applying the Lax-Milgram Lemma to the variational equation (1.22). □

Let $\boldsymbol{W}_h \subset \boldsymbol{W}$ be finite dimensional, we consider the following problem: find $\boldsymbol{d}_h \in \boldsymbol{W}_h$ such that

$$a^{\mathrm{s}}(\boldsymbol{d}_h, \boldsymbol{w}_h) = \langle \boldsymbol{f}, \boldsymbol{w}_h\rangle \qquad \forall \boldsymbol{w}_h \in \boldsymbol{W}_h. \quad (1.24)$$

The properties of coercivity and continuity transfers directly to any subspace $\boldsymbol{W}_h \subset \boldsymbol{W}$, then applying the previous theorem when $a^s(\cdot,\cdot)$ is restricted to $\boldsymbol{W}_h$, it is possible to state the existence of a solution, its uniqueness and its continuous dependence from the data. Moreover it holds the following optimal error estimate

**Proposition 1.5.2.** *Let $\boldsymbol{d}$ be the unique solution of the continuous problem* (1.22) *and $\boldsymbol{d}_h$ the unique solution of* (1.24)*, then*

$$\|\boldsymbol{d}-\boldsymbol{d}_h\|_{1,\Sigma} \le C \inf_{\boldsymbol{w}_h\in\boldsymbol{W}_h} \|\boldsymbol{d}-\boldsymbol{w}_h\|_{1,\Sigma}. \tag{1.25}$$

*Proof.* The proof of this statement is a direct consequence of the Céa Lemma (see Proposition A.2.3). □

**The fluid bilinear form $a^f(\cdot,\cdot)$.** The fluid bilinear form $a^f(\cdot,\cdot)$ is defined on the product space $\boldsymbol{V}\times Q$ and is associated to the Stokes problem. Find $(\boldsymbol{u},p)\in\boldsymbol{V}\times Q$ such that

$$a^f((\boldsymbol{u},p),(\boldsymbol{v},q)) = 2\mu\big(\boldsymbol{\epsilon}(\boldsymbol{u}),\boldsymbol{\epsilon}(\boldsymbol{v})\big)_{0,\Omega} - (p,\operatorname{div}\boldsymbol{v})_{0,\Omega} + (q,\operatorname{div}\boldsymbol{u})_{0,\Omega} \qquad \forall(\boldsymbol{v},q)\in\boldsymbol{V}\times Q. \tag{1.26}$$

In the following proposition we recall the existence, uniqueness and stability of the solution of the Stokes problem. We refer to Appendix C and its references for the details.

**Proposition 1.5.3.** *For all $\boldsymbol{f}\in\boldsymbol{V}'$, there exists a unique $(\boldsymbol{u},p)\in\boldsymbol{V}\times Q$ such that*

$$a^f((\boldsymbol{u},p),(\boldsymbol{v},q)) = \langle\boldsymbol{f},\boldsymbol{v}\rangle \qquad \forall(\boldsymbol{v},q)\in\boldsymbol{V}\times Q, \tag{1.27}$$

*moreover there exists $C>0$ independent of $\boldsymbol{f}$, such that*

$$\|\boldsymbol{u}\|_{1,\Omega} + \|p\|_{0,\Omega} \le C\|\boldsymbol{f}\|_{\boldsymbol{V}'}$$

Since Stokes problem (1.27) is a saddle point problem, its finite element discretization is more delicate than the finite element discretization of the solid problem; indeed we need to assign a couple of finite element spaces $\boldsymbol{V}_h, Q_h$ that fulfill the inf-sup condition (see Proposition C.1.4 of Appendix C) or we need to add a suitable stabilization term. For a more exaustive discussion see Appendix C and its references.

## 1.6 Alternative coupling techniques

Given the strong form of the coupled problem, there are many ways to obtain a weak formulation; the possible formulations differ principally for the choice of the technique used to impose the coupling conditions. To tackle this question there are essentially two approaches.

- We can leave the coupling conditions in strong form: this choice, from the discrete point of view, forces to assign prescribed values to interface nodal unknowns. This approach involves different techniques for matching and non matching grids.
  - Matching grids for fluid and solid have common nodes on the interface, then the coupling conditions can be imposed directly on these nodes. This type of approach forces to deform the fluid mesh as the solid evolves and typically ALE formulation is used in this case [108, 77, 46].

- In the case of non matching grids, solid and fluid meshes are independent and do not need to be modified if the solid moves. In this case an interpolation/projection step has to be carried out to transfer the information between the two domains. The simplest way is to obtain the information for one variable from the closest point in the other mesh. This approach is known as the nearest neighbour interpolation [101, 42]. A more accurate way of handling the data transfer is by projection. To obtain information from the other mesh, a point can be orthogonally projected on that mesh and the information in that projection point can be used at the original point [49, 38, 85].

  A particular mention in this field has to be reserved to the **Immersed Boundary Method** (IBM) proposed by C.S. Peskin at the beginning of the '70 in order to simulate the the blood flow in the heart during a beat [97, 88]. The IBM, in the original version, applies to the case of solid fibers that move in the fluid, then in order to perform 3D simulations, the solid has to be seen as a bundle of fibers. After its introduction, the IBM method has been applied to three dimensional heart valve simulations [97, 88] and to a wide variety of other biological and non biological problems [98]. The IBM combines a fixed Cartesian grid to compute the fluid motion, expressed in Eulerian coordinates, along with a solid grid on which is computed the solid motion expressed in Lagrangian coordinates. The key idea of this method is to replace the solid material with fluid and to apply to it an appropriate force density term localized along the structure by means of the Dirac delta distribution, then the fluid equations are solved over the entire domain with the influence of the solid mimed by the added forcing term. The original IBM equations are discretized in space using the finite difference method, for both the Eulerian and the Lagrangian variables. In order to transfer informations between the solid and fluid meshes, the IBM, in its original formulation, uses an approximation of the Dirac delta distribution, centered at the solid nodes, that assign to the fluid nodes in its support a force term [98, 89]. The approximation of the Dirac delta distribution represents the main difficulty in the discretization of classical IBM [24]; as shown in [22, 24], in the finite element method, the Dirac delta distribution can be treated weakly, so that its approximation is no more needed. The variational formulation presented in these articles can be applied naturally also to thick solids, eliminating the fiber-like geometry assumption for the solid in the original formulation [26].

- Weak form of the coupling conditions: this choice consists in a modification of the variational formulation in such a way that the transmission conditions are automatically satisfied by the weak solution of the problem. In this category we can collect the methods that we are going to describe briefly.

  - **Penalty method** Penalty method, can be seen as an approximation of Lagrange multiplier method, in the sense that, from a formal point of view we can obtain penalized formulation from the Lagrange multiplier formulation neglecting the constraint equation and assuming that the Lagrange multiplier $\boldsymbol{\lambda}$ associated to the constraint $\boldsymbol{u} \circ \boldsymbol{\phi} - \partial_t \boldsymbol{d} = 0$ has the form

$$\boldsymbol{\lambda} \cong \frac{1}{\epsilon}(\boldsymbol{u} \circ \boldsymbol{\phi} - \partial_t \boldsymbol{d}),$$

where $\epsilon$ is a small parameter that has to be set in such a way that the constraint, enforcing the equality of fluid and solid velocities on the structure, is fulfilled with the desired degree of accuracy. Indeed the penalty method can be seen as a technique to stabilize the saddle point problem obtained using the Lagrange multiplier method. In order to introduce the penalty formulation for fluid-structure problem we consider the following abstract form of the **stationary problem**. Let us introduce the Hilbert space

$$\mathcal{V} = \boldsymbol{V} \times Q \times \boldsymbol{W}, \tag{1.28}$$

product of the spaces of the fluid velocity, the fluid pressure and the solid displacement. $\boldsymbol{\Lambda}$ denotes the space of the Lagrange multiplier. Then, the stationary problem with Lagrange multiplier can be rewritten in operator form as: find $\boldsymbol{U} = (\boldsymbol{u}, p, \boldsymbol{d}) \in \mathcal{V}$ and $\boldsymbol{\lambda} \in \boldsymbol{\Lambda}$ such that

$$\begin{aligned} \boldsymbol{AU} + \boldsymbol{B}^T \boldsymbol{\lambda} &= \boldsymbol{F}, \\ \boldsymbol{BU} &= \boldsymbol{0}, \end{aligned} \tag{1.29}$$

where the operator $\boldsymbol{A} : \mathcal{V} \to \mathcal{V}'$ represents the fluid and solid problem, the operator $\boldsymbol{B} : \mathcal{V} \to \boldsymbol{\Lambda}'$ stands for the coupling term associate to the Lagrange multiplier and $\boldsymbol{F} \in \mathcal{V}'$ represents the external loads.

The penalty method can be obtained perturbing the saddle point problem (1.29) in the following form

$$\begin{aligned} \boldsymbol{AU} + \boldsymbol{B}^T \boldsymbol{\lambda} &= \boldsymbol{F}, \\ \boldsymbol{BU} - \epsilon \boldsymbol{C} \boldsymbol{\lambda} &= \boldsymbol{0}, \end{aligned} \tag{1.30}$$

where $\boldsymbol{C} : \boldsymbol{\Lambda} \to \boldsymbol{\Lambda}'$ is linear continuous and invertible. Then it is possible to eliminate the Lagrange multiplier from the system obtaining

$$\boldsymbol{AU} + \frac{1}{\epsilon} \boldsymbol{B}^T \boldsymbol{C}^{-1} \boldsymbol{BU} = \boldsymbol{0}. \tag{1.31}$$

When the parameter $\epsilon \to 0$, the perturbed problem (1.30) converges to the original problem (1.29) with Lagrange multiplier and, on the other hand, the imposition of the constraint $\boldsymbol{BU} = \boldsymbol{0}$ in the penalized problem becomes "exact". The advantage of this technique with respect to the use of Lagrange multiplier, is that the penalized problem does not introduce a new unknown and does not require an inf-sup property. On the other hand, there are drawbacks because, in order to obtain accurate solutions that respect the constraint, we need small values of $\epsilon$ but the more $\epsilon$ is small, the more the system becomes ill-conditioned [5, 94, 93].

- **Augmented Lagrange multiplier method** The contradictory demands placed on $\epsilon$ by accuracy and efficiency requirements are a serious weakness of the penalty method. The augmented Lagrange multiplier method is a combination of Lagrange multiplier method and the penalty method. As the penalty method, it can be obtained from a perturbation of the saddle point formulation (1.29) as follows

$$\begin{aligned} \left(\boldsymbol{A} + \tfrac{1}{\epsilon}\boldsymbol{B}^T\boldsymbol{C}^{-1}\boldsymbol{B}\right)\boldsymbol{U} + \boldsymbol{B}^T\boldsymbol{\lambda} = \boldsymbol{F}, \\ \boldsymbol{B}\boldsymbol{U} = \boldsymbol{0}. \end{aligned} \tag{1.32}$$

Choosing a pair of conforming subspaces $\mathbb{V}_h \subset \mathbb{V}$ and $\boldsymbol{\Lambda}_h \subset \boldsymbol{\Lambda}$ in such a way that the inf-sup condition is satisfied, the discretization of (1.32) results in a linear system that depends on the parameter $\epsilon$. Thanks to the presence of the Lagrange multiplier, the constraint is enforced "exactly" ad the parameter $\epsilon$ can be chosen to optimize the solution of the linear system. A great advantage of the augmented Lagrange multiplier is that if the correct values of the Lagrange multiplier $\boldsymbol{\lambda}^*$ is known, it can be substituted in the system and we obtain

$$\begin{aligned} \boldsymbol{A}\boldsymbol{U} + \tfrac{1}{\epsilon}\boldsymbol{B}^T\boldsymbol{C}^{-1}\boldsymbol{B}\boldsymbol{U} + \boldsymbol{B}^T\boldsymbol{\lambda}^* = \boldsymbol{F} \\ \boldsymbol{B}\boldsymbol{U} = \boldsymbol{0}, \end{aligned}$$

In this case, for the discretization, we need only the subspaces $\mathbb{V}_h \subset \mathbb{V}$ and the invertibility of the operator $\boldsymbol{A} + \frac{1}{\epsilon}\boldsymbol{B}^T\boldsymbol{C}^{-1}\boldsymbol{B}$. For a complete analysis of this kind of approach to the fluid-structure problem see [81, 82, 79, 12].

Another possibility offered by the augmented Lagrange multiplier method is the implementation of an iterative method in order to compute the solution. In fact, given a guess $\boldsymbol{\lambda}_0$, for $k \in \mathbb{N}$ we can solve the system

$$\begin{aligned} \left(\boldsymbol{A} + \tfrac{1}{\epsilon}\boldsymbol{B}^T\boldsymbol{C}^{-1}\boldsymbol{B}\right)\boldsymbol{U}_k = \boldsymbol{F} - \boldsymbol{B}^T\boldsymbol{\lambda}_k \\ \boldsymbol{\lambda}_{k+1} = \boldsymbol{\lambda}_k + \tfrac{1}{\epsilon}\boldsymbol{B}\boldsymbol{U}_k, \end{aligned}$$

until the difference between $\boldsymbol{\lambda}_k$ and $\boldsymbol{\lambda}_{k+1}$ is smaller that a fixed value.

- **Nitsche method** Nitsche's method [90] was originally proposed to weakly enforce Dirichlet boundary conditions as an alternative to pointwise constraints without using Lagrange multiplier. The Nitsche approach can be obtained from the Lagrange multiplier method replacing the multiplier with its physical representation, namely the normal flux at the interface. In recent years, the method has been reconsidered for interface problems [70, 72, 45], for connecting overlapping meshes [71, 15, 70, 106, 107], for imposing Dirichlet boundary conditions in meshfree methods [61], in immersed boundary methods [105, 9, 47], in fluid mechanics [34, 13] and for contact mechanics [113]. Like the penalty method, the Nitsche method, is very simple to implement with the advantage of having no conditioning problem and optimal convergence properties in $L^2$ and $H^1$. Considering equations (1.10), (1.12) and (1.9) we can derive the following expression of the Nitsche method at the continuous level find $(\boldsymbol{u}, p, \boldsymbol{d}) \in \boldsymbol{V} \times Q \times \boldsymbol{W}$ such that

$$\rho_{\mathrm{f}}(\partial_t \boldsymbol{u}, \boldsymbol{v})_{0,\Omega} + \rho_{\mathrm{s}}(\partial_{tt}\boldsymbol{d}, \boldsymbol{w})_{0,\Sigma} + a^{\mathrm{f}}((\boldsymbol{u}, p);(\boldsymbol{v}, q)) + a_{\mathrm{s}}(\boldsymbol{d}, \boldsymbol{w})$$

$$-\int_{\Sigma} [\![ \boldsymbol{P}_{\mathrm{f}}(\boldsymbol{u},p)\boldsymbol{n}_{\Sigma} ]\!] \cdot (\boldsymbol{v}\circ\boldsymbol{\phi}-\boldsymbol{w}) - \int_{\Sigma} [\![ \boldsymbol{P}_{\mathrm{f}}(\boldsymbol{v},p)\boldsymbol{n}_{\Sigma} ]\!] \cdot (\boldsymbol{u}\circ\boldsymbol{\phi}-\partial_t \boldsymbol{d}) = 0 \quad \forall(\boldsymbol{v},q,\boldsymbol{w}) \in \boldsymbol{V}\times \boldsymbol{Q} \times \boldsymbol{W}.$$

At the discrete level the Nitsche method contains a stabilization term of the form

$$\frac{\gamma}{h_{\mathrm{s}}}\int_{\Sigma}(\boldsymbol{u}_h^n\circ\boldsymbol{\phi}-\partial_\tau \boldsymbol{d}_h^n)\cdot(\boldsymbol{v}_h\circ\boldsymbol{\phi}-\boldsymbol{w}_h),$$

where $\gamma$ is an appropriate constant to be chosen in order to guarantee the coercivity of the formulation and $h_{\mathrm{s}}$ is the solid mesh size [60, 84].

### 1.6.1 Time integration schemes

In the computational treatment of a dynamical coupled system, the design of a time marching scheme is a very important issue. A very rough classification of the possible approaches is based on the method used to solve the fully discrete coupled problem; in this regard we have essentially the following approaches [54, 51]

- *Monolithic or Simultaneous Treatment* The whole problem is treated as a monolithic entity, and all components are advanced simultaneously in time. For the problem of fluid structure interaction, the conceptual scheme of a time-step of the a monolithic algorithm is shown in the following Figure 1.3 where $\boldsymbol{U}_h^n = (\boldsymbol{u}_h^n, p_h^n, \boldsymbol{d}_h^n, \boldsymbol{\lambda}_h^n)$ represents the vector of unknowns at time $t^n$.

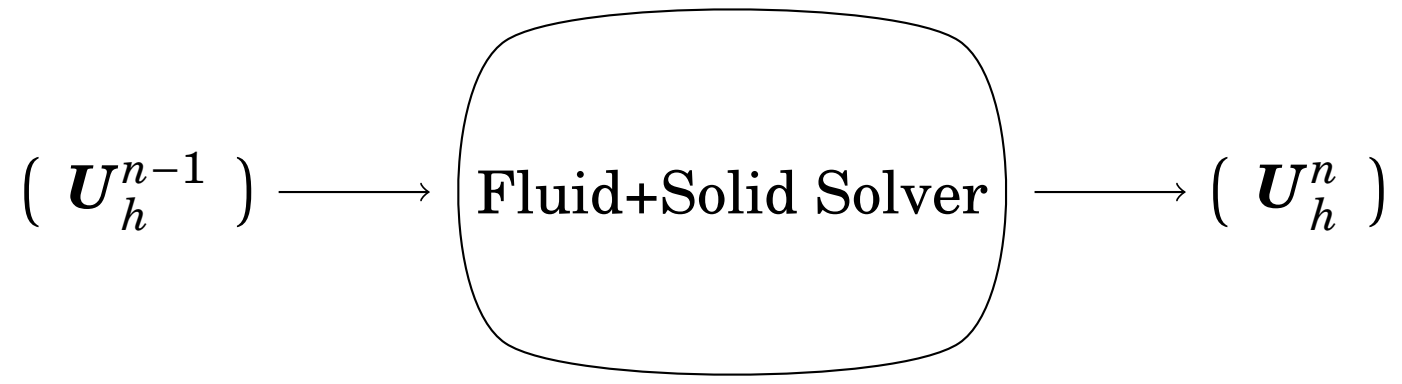

Figure 1.3: Time step of a monolithic algorithm.

  The monolithic approach is usually unconditionally stable, but requires the development of a unified fluid-structure simulator and can be computationally expensive [31].

- *Partitioned Treatment* The field models are computationally treated as isolated entities that are separately stepped in time. Interaction effects are viewed as forcing terms that are communicated between the individual components using prediction, substitution and synchronization techniques. Generally, the solution of a time step begins with a prediction of the unknown of one of the problems, then the other problem is solved using the prediction, after that, the computed unknown are transferred to the first problem whose solution corrects the predicted value. This process can be performed only one time, giving a *staggered scheme*, or can be iterated until the interface conditions are fulfilled with the desired accuracy giving a *sequentially scheme*.

The conceptual scheme of a time-step of a partitioned algorithm for the fluid structure problem is shown in the following Figure 1.4

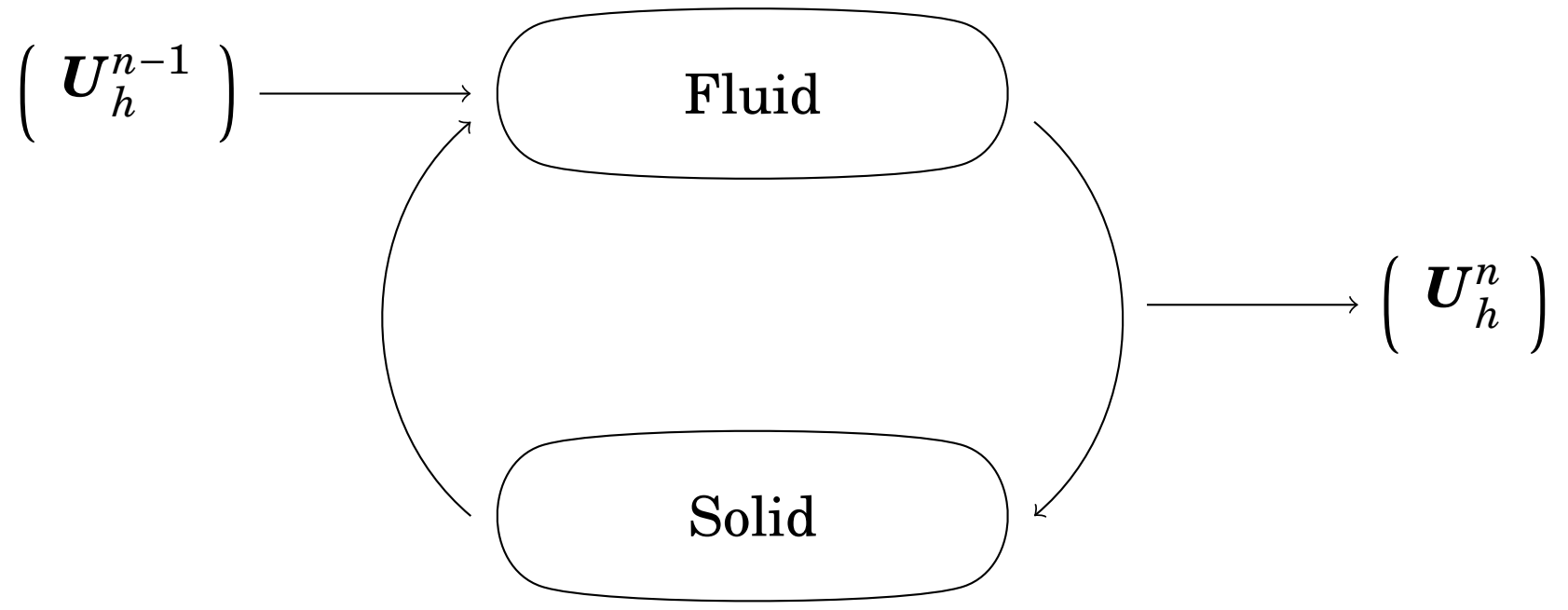

Figure 1.4: Time-step of a partitioned algorithm.

Partition of schemes, sequential or staggered, offer wide flexibility and are highly desirable from a software engineering perspective. Moreover, this type of scheme allows for using specialized numerical methods for the solid and the fluid problem [51]. For a sequential simulation framework to be competitive in engineering practice, it must have stability and convergence properties that are similar to those enjoyed by the monolithic method. For this reason significant efforts to find stable and efficient partitioned schemes for coupled fluid-structure problem have been pursued in the engineering and mathematical communities [100, 50, 99, 92].

Another possible classification of the schemes is based on the way that the coupling conditions are treated at the fully discrete level; in the technical literature we find essentially three type of approaches

- *implicit or strongly coupled schemes* In a strongly coupled scheme there is no time delay between the fluid and solid time-marchings. This features, generally, give unconditional stability and optimal accuracy, but at the price of an high computational cost, since it is necessary to solve the entire coupled system at each time-step. The corresponding solution procedures can be monolithic but also sequentially partitioned.

- *explicit or loosely coupled schemes* Explicit or loosely coupled schemes, uncouple the fluid and the solid time-marchings via appropriate explicit discretization of the interface conditions. The resulting algorithm belongs naturally to the class of partitioned schemes. Loosely coupled schemes are very attractive for the applications since they present a reduced computational cost but, generally, they do not satisfy exactly the coupling transmission conditions. As a consequence, the work exchanged between the two sub-problems is not perfectly balanced and this may induce instabilities in the numerical scheme. For example, it was shown in [37] (see also [63]) that an explicit coupling is unstable in those applications where the added mass effect is important, as in haemodynamics. For many years, stability in explicit coupling has been an open problem that attracted the interest of many researchers. Stabilized explicit coupled schemes have been conceived but at the expense of a degradation of accuracy, which requires suitable correction iterations [34, 35, 67, 56]. Many of these issues are overcome by the explicit Robin-Neumann schemes proposed in [56, 58], which simultaneously deliver unconditional stability and optimal (first-order) time accuracy.

- *semi-implicit or splitting schemes*

  An intermediate possibility are the semi-implicit schemes that enforce an explicit/implicit treatment of the kinematic and dynamic coupling conditions. These schemes are generally stables and less computationally demanding than implicit schemes. They, often, involve a fractional-step time-marching in the fluid (see [57, 8, 102, 4, 2]) or in the solid (see [67, 55, 30, 86]). The implicit part of the coupling can be solved in a monolithic or a partitioned fashion and guarantees stability, while the explicit one reduces computational complexity. The conceptual scheme of the semi-implicit scheme studied in this thesis is reported in Figure 1.5. In practice each time-step is partitioned in two sub-steps; in the first sub-step we obtain the solution of a coupled problem in which the coupling is implicit for the force term related to solid inertia and explicit for the force term related to the solid elasticity. The solution of the first sub-step gives the vector of unknowns $\boldsymbol{U}_h^{n-\frac{1}{2}}$ that is used in the second sub-step where a residual solid problem is solved for the sake of correct the computed values of the solid unknown giving the final value of $\boldsymbol{U}_h^n$. The advantage of this approach is that it allow for a partitioned treatment of the coupled problem preserving in the same time stability and convergence. Indeed the first subproblem is a pure fluid problem with an additional term that is related to the solid inertia, while, the second subproblem is a pure elastodynamic system.

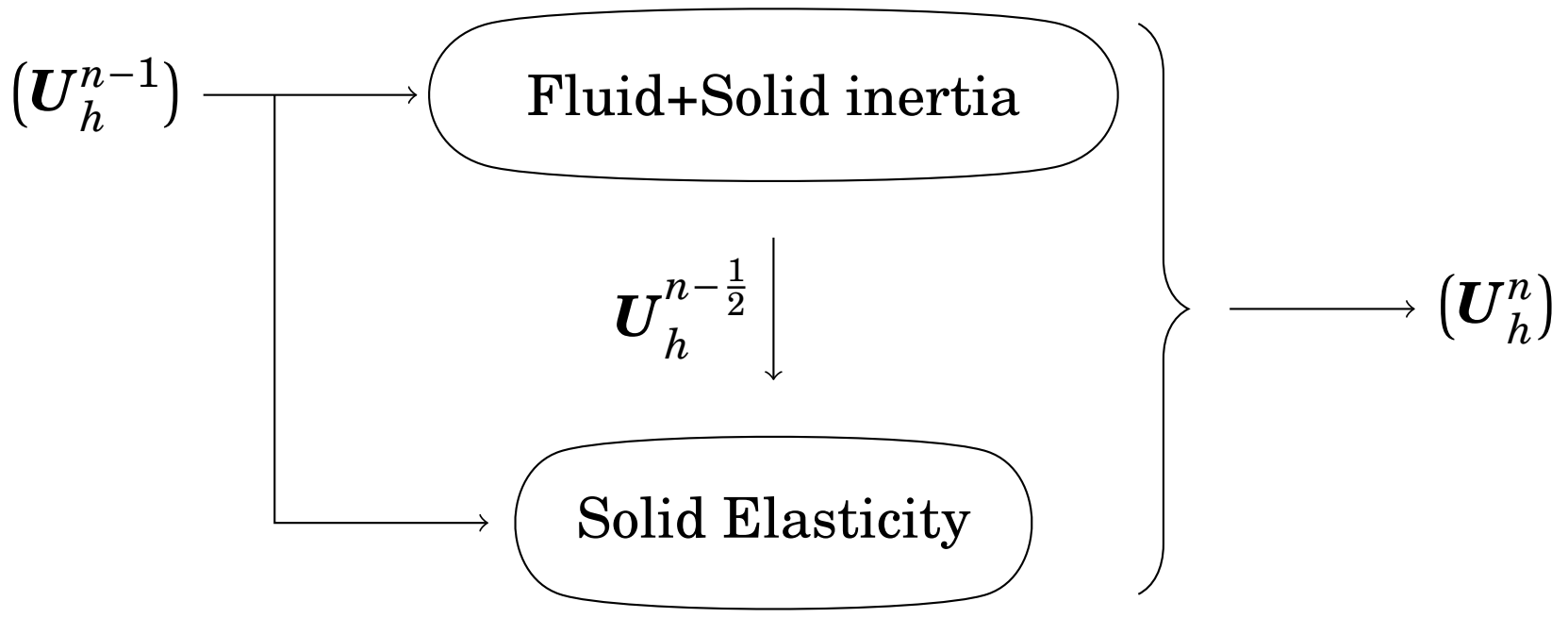

Figure 1.5: Semi-implicit Algorithm

# Chapter 2

# Analysis of the coupled problem

## 2.1 Introduction

In this chapter we start to analyze the numerical solution of the coupled Problem 2 presented in Chapter 1. In the first part of the chapter we present results regarding the stability of the time semi-discrete scheme and the well-posedness of the associated problem. In the second part we introduce the spatial discretization and analyze the well-posedness of the time-steps of the fully discrete problem. In this regard we consider two cases, in the first one the spatial discretization uses Stokes-stable finite element spaces for the fluid velocity and pressure, while in the second case we consider stabilized $\mathbb{P}_1 - \mathbb{P}_1$ elements. The analysis of stability and well-posedness of the time semi-discrete problem and of the fully discrete problem in the case of Stokes-stable finite elements, has been performed in [20, 23] and we report here a summary of the results. In this thesis we show also the well-posedness of the time-steps of the fully discrete problem using the stabilized finite element method, which in our knowledge is a new result. We recall briefly the notation used and the problem under study

$$\boldsymbol{V} \stackrel{def}{=} H^1_0(\Omega)^d, \qquad Q \stackrel{def}{=} L^2_0(\Omega), \qquad \boldsymbol{W} \stackrel{def}{=} H^1(\Sigma)^d \cap \mathscr{B}\mathscr{C}, \qquad \boldsymbol{\Lambda} \stackrel{def}{=} \left(H^{\frac{1}{2}}(\Sigma)^d\right)', \tag{2.1}$$

where $\mathscr{B}\mathscr{C}$ represents the set of functions $\boldsymbol{w} : \Sigma \to \mathbb{R}^d$ that satisfy the essential boundary conditions for the solid displacements.

**Problem 3.** *Let $\boldsymbol{\phi}_0$ be the position vector of the solid reference configuration and let $(\boldsymbol{u}_0, \boldsymbol{d}_0, \boldsymbol{d}_1) \in \boldsymbol{V} \times \boldsymbol{W} \times \boldsymbol{W}$ be the initial values for the fluid velocity, the solid displacement and the fluid velocity. For $t \geq 0$, find $(\boldsymbol{u}(t), p(t), \boldsymbol{d}(t), \boldsymbol{\lambda}(t)) \in \boldsymbol{V} \times Q \times \boldsymbol{W} \times \boldsymbol{\Lambda}$ such that*

$$\begin{aligned}
&\boldsymbol{\phi}(t) = \boldsymbol{\phi}_0 + \boldsymbol{d}(t), \\
&\rho_{\mathrm{f}}\big(\partial_t \boldsymbol{u}, \boldsymbol{v}\big)_{0,\Omega} + 2\mu\big(\boldsymbol{\epsilon}(\boldsymbol{u}), \boldsymbol{\epsilon}(\boldsymbol{v})\big)_{0,\Omega} - (p, \operatorname{div} \boldsymbol{v})_{0,\Omega} + c\left(\boldsymbol{\lambda}, \boldsymbol{v} \circ \boldsymbol{\phi}\right) = 0 \qquad && \forall \boldsymbol{v} \in \boldsymbol{V}, \\
&\rho_{\mathrm{s}}\epsilon\big(\partial_{tt} \boldsymbol{d}, \boldsymbol{w}\big)_{0,\Sigma} + a^{\mathrm{s}}\big(\boldsymbol{d}, \boldsymbol{w}\big) - c\left(\boldsymbol{\lambda}, \boldsymbol{w}\right) = 0 && \forall \boldsymbol{w} \in \boldsymbol{W}, \\
&c\left(\boldsymbol{\mu}, \boldsymbol{u} \circ \boldsymbol{\phi} - \partial_t \boldsymbol{d}\right) = 0 && \forall \boldsymbol{\mu} \in \boldsymbol{\Lambda}, \\
&(\operatorname{div} \boldsymbol{u}, q)_{0,\Omega} = 0 && \forall q \in Q, \\
&(\boldsymbol{u}(0), \boldsymbol{d}(0), \partial_t \boldsymbol{d}(0)) = (\boldsymbol{u}_0, \boldsymbol{d}_0, \boldsymbol{d}_1),
\end{aligned}$$

where we used the notation $(\cdot,\cdot)_{0,\Omega}$ and $(\cdot,\cdot)_{0,\Sigma}$ to indicate the scalar products in $L^2(\Omega)^d$ and $L^2(\Sigma)^d$ respectively. The following energy estimate holds true

**Proposition 2.1.1.** *Energy estimate*
*For $t \geq 0$, suppose that $(\boldsymbol{u}(t), p(t), \boldsymbol{d}(t), \boldsymbol{\lambda}(t)) \in \boldsymbol{V} \times Q \times \boldsymbol{W} \times \boldsymbol{\Lambda}$ satisfy Problem 3, then*

$$\frac{\rho_{\mathrm{f}}}{2}\frac{d}{dt}\|\boldsymbol{u}(t)\|_{0,\Omega}^2 + 2\mu\|\boldsymbol{\epsilon}(\boldsymbol{u}(t))\|_{0,\Omega}^2 + \frac{\rho_{\mathrm{s}}\epsilon}{2}\frac{d}{dt}\|\partial_t\boldsymbol{d}(t)\|_{0,\Sigma}^2 + \frac{1}{2}\frac{d}{dt}\|\boldsymbol{d}(t)\|_{\mathrm{s}}^2 = 0 \tag{2.2}$$

*Proof.* We take $\boldsymbol{v} = \boldsymbol{u}$, $q = p$, $\boldsymbol{w} = \partial_t\boldsymbol{d}$, $\boldsymbol{\mu} = \boldsymbol{\lambda}$ in Problem 3, then summing up the equations we have

$$\begin{aligned}&\rho_{\mathrm{f}}(\partial_t\boldsymbol{u},\boldsymbol{u})_{0,\Omega} + 2\mu(\boldsymbol{\epsilon}(\boldsymbol{u}),\boldsymbol{\epsilon}(\boldsymbol{u}))_{0,\Omega} - (\mathrm{div}\,\boldsymbol{u},p)_{0,\Omega} + (\mathrm{div}\,\boldsymbol{u},p)_{0,\Omega}\\&+ c(\boldsymbol{\lambda},\boldsymbol{u}\circ\boldsymbol{\phi} - \partial_t\boldsymbol{d}) - c(\boldsymbol{\lambda},\boldsymbol{u}\circ\boldsymbol{\phi} - \partial_t\boldsymbol{d}) + \rho_{\mathrm{s}}\epsilon(\partial_{tt}\boldsymbol{d},\partial_t\boldsymbol{d})_{0,\Sigma} + a^{\mathrm{s}}(\boldsymbol{d},\partial_t\boldsymbol{d}) = 0.\end{aligned}$$

Using (1.21) and the following relations

$$(\partial_t\boldsymbol{u},\boldsymbol{u})_{0,\Omega} = \frac{1}{2}\frac{d}{dt}(\boldsymbol{u},\boldsymbol{u})_{0,\Omega} = \frac{1}{2}\frac{d}{dt}\|\boldsymbol{u}\|_{0,\Omega}^2, \qquad (\partial_{tt}\boldsymbol{d},\partial_t\boldsymbol{d})_{0,\Sigma} = \frac{1}{2}\frac{d}{dt}(\partial_t\boldsymbol{d},\partial_t\boldsymbol{d})_{0,\Sigma} = \frac{1}{2}\frac{d}{dt}\|\partial_t\boldsymbol{d}\|_{0,\Sigma}^2,$$

we infer the result. □

## 2.2 Time semi-discrete problem

In this section, we consider the time semi-discretization of Problem 3. The evolutionary problem, in this way, is reduced to a sequence of stationary problems to be solved at each time step. Following [23], we present results of well-posedness for the stationary problem. Given an integer $N > 0$, set $\tau = \frac{T}{N}$ the time step and $t_n = n\tau$. For a given function $g$ depending on $t$ we adopt the following notation

$$g^n = g(t_n), \qquad \partial_\tau g^n = \frac{g^n - g^{n-1}}{\tau}, \qquad \partial_{\tau\tau} g^n = \frac{g^n - 2g^{n-1} + g^{n-2}}{\tau^2}$$

then, approximating the time derivatives with the backward finite differences, we have the following semi-discrete problem in time

**Problem 4.** *Given $\boldsymbol{\phi}_0$, the position vector of the solid reference configuration, and let $(\boldsymbol{u}_0, \boldsymbol{d}_0, \boldsymbol{d}_1) \in \boldsymbol{V} \times \boldsymbol{W} \times \boldsymbol{W}$ be the initial values for the fluid velocity, the solid displacement and the fluid velocity. For $n \geq 1$, find $(\boldsymbol{u}^n, \boldsymbol{d}^n, \boldsymbol{\lambda}^n, p^n) \in \boldsymbol{V} \times \boldsymbol{W} \times \boldsymbol{\Lambda} \times Q$ such that*

$$\begin{aligned}&\rho_{\mathrm{f}}\big(\partial_\tau\boldsymbol{u}^n,\boldsymbol{v}\big)_{0,\Omega} + 2\mu\big(\boldsymbol{\epsilon}(\boldsymbol{u}^n),\boldsymbol{\epsilon}(\boldsymbol{v})\big)_{0,\Omega} - \big(p^n,\mathrm{div}\,\boldsymbol{v}\big)_{0,\Omega} + c\left(\boldsymbol{\lambda}^n,\boldsymbol{v}\circ\boldsymbol{\phi}^{n-1}\right) = 0 && \forall\boldsymbol{v}\in\boldsymbol{V},\\&\rho_{\mathrm{s}}\epsilon\big(\partial_{\tau\tau}\boldsymbol{d}^n,\boldsymbol{w}\big)_{0,\Sigma} + a^{\mathrm{s}}\big(\boldsymbol{d}^n,\boldsymbol{w}\big) - c\left(\boldsymbol{\lambda}^n,\boldsymbol{w}\right) = 0 && \forall\boldsymbol{w}\in\boldsymbol{W},\\&c\left(\boldsymbol{\mu},\boldsymbol{u}^n\circ\boldsymbol{\phi}^{n-1} - \partial_\tau\boldsymbol{d}^n\right) = 0 && \forall\boldsymbol{\mu}\in\boldsymbol{\Lambda},\\&\big(\mathrm{div}\,\boldsymbol{u}^n,q\big)_{0,\Omega} = 0 && \forall q\in Q,\\&\boldsymbol{\phi}^n = \boldsymbol{\phi}_0 + \boldsymbol{d}^n,\end{aligned} \tag{2.3}$$

*with*

$$\boldsymbol{u}^0 = \boldsymbol{u}_0, \qquad \boldsymbol{d}^0 = \boldsymbol{d}_0, \qquad \boldsymbol{d}^{-1} = \boldsymbol{d}_0 - \tau\boldsymbol{d}_1.$$

At each time-step, the semi-discrete problem in time is a stationary problem that we rewrite in the following generalized form

**Problem 5.** *Let* $\overline{\boldsymbol{\phi}} = \boldsymbol{\phi}_0 + \tau\overline{\boldsymbol{d}} \in W^{1,\infty}(\Sigma)^d$ *be invertible with Lipschitz inverse. Given* $\boldsymbol{f} \in L^2(\Omega)^d$ *and* $\boldsymbol{g} \in L^2(\Sigma)^d$, *find* $(\boldsymbol{u}, p, \boldsymbol{d}, \boldsymbol{\lambda}) \in \boldsymbol{V} \times Q \times \boldsymbol{W} \times \boldsymbol{\Lambda}$ *such that*

$$\begin{aligned}
&\alpha_1(\boldsymbol{u},\boldsymbol{v})_{0,\Omega} + 2\mu\big(\boldsymbol{\epsilon}(\boldsymbol{u}),\boldsymbol{\epsilon}(\boldsymbol{v})\big)_{0,\Omega} - (p,\operatorname{div}\boldsymbol{v})_{0,\Omega} + c\left(\boldsymbol{\lambda},\boldsymbol{v}\circ\overline{\boldsymbol{\phi}}\right) = (\boldsymbol{f},\boldsymbol{v})_{0,\Omega} && \forall \boldsymbol{v}\in\boldsymbol{V},\\
&\alpha_2(\boldsymbol{d},\boldsymbol{w}) + a^{\mathrm{s}}(\boldsymbol{d},\boldsymbol{w}) - c(\boldsymbol{\lambda},\boldsymbol{w}) = (\boldsymbol{g},\boldsymbol{w})_{0,\Sigma} && \forall \boldsymbol{w}\in\boldsymbol{W},\\
&c\left(\boldsymbol{\mu},\boldsymbol{u}\circ\overline{\boldsymbol{\phi}} - \boldsymbol{d}\right) = -c(\boldsymbol{\mu},\overline{\boldsymbol{d}}) && \forall \boldsymbol{\mu}\in\boldsymbol{\Lambda},\\
&(\operatorname{div}\boldsymbol{u},q)_{0,\Omega} = 0 && \forall q\in Q.
\end{aligned}$$

Problem 5 gives a time step of Problem 4 if

$$\alpha_1 = \frac{\rho_f}{\tau}, \qquad \alpha_2 = \frac{\rho_s\epsilon}{\tau^2}, \qquad \boldsymbol{f} = \frac{\rho_f}{\tau}\boldsymbol{u}^{n-1}, \qquad \boldsymbol{g} = \frac{\rho_s\epsilon}{\tau^2}(2\boldsymbol{d}^{n-1} - \boldsymbol{d}^{n-2}),$$
$$\boldsymbol{d} = \frac{\boldsymbol{d}^n}{\tau}, \qquad \overline{\boldsymbol{d}} = \frac{\boldsymbol{d}^{n-1}}{\tau}, \qquad \boldsymbol{u} = \boldsymbol{u}^n, \qquad p = p^n, \qquad \boldsymbol{\lambda} = \boldsymbol{\lambda}^n.$$

### 2.2.1 Analysis of the time semi-discrete problem

Problem 5 has the structure of a saddle point problem and can be written naturally in the following matrix operator form

$$\left(\begin{array}{ccc|c}
A_f & 0 & L_f^\top & B^\top\\
0 & A_s & -L_s^\top & 0\\
L_f & -L_s & 0 & 0\\
\hline
B & 0 & 0 & 0
\end{array}\right)
\left(\begin{array}{c}
\boldsymbol{u}\\ \boldsymbol{d}\\ \boldsymbol{\lambda}\\ \hline p
\end{array}\right)
=
\left(\begin{array}{c}
\boldsymbol{f}\\ \boldsymbol{g}\\ \mathbf{m}\\ \hline 0
\end{array}\right),$$

with the following definition for the operators.

$$\begin{aligned}
&A_f : \boldsymbol{V}\to\boldsymbol{V}' && \text{such that} && \langle A_f\boldsymbol{u},\boldsymbol{v}\rangle \overset{def}{=} \alpha_1(\boldsymbol{u},\boldsymbol{v})_{0,\Omega} + 2\mu\big(\boldsymbol{\epsilon}(\boldsymbol{u}),\boldsymbol{\epsilon}(\boldsymbol{v})\big)_{0,\Omega} && \forall \boldsymbol{u},\boldsymbol{v}\in\boldsymbol{V},\\
&A_s : \boldsymbol{W}\to\boldsymbol{W}' && \text{such that} && \langle A_s\boldsymbol{d},\boldsymbol{w}\rangle \overset{def}{=} \alpha_2(\boldsymbol{d},\boldsymbol{w}) + a^{\mathrm{s}}(\boldsymbol{d},\boldsymbol{w}) && \forall \boldsymbol{d},\boldsymbol{w}\in\boldsymbol{W},\\
&L_f : H^{\frac12}(\Sigma)^d\to\boldsymbol{\Lambda}' && \text{such that} && \langle L_f\boldsymbol{u},\boldsymbol{\mu}\rangle \overset{def}{=} c\left(\boldsymbol{\mu},\boldsymbol{u}\circ\overline{\boldsymbol{\phi}}\right) && \forall \boldsymbol{u}\in\boldsymbol{V},\forall\boldsymbol{\mu}\in\boldsymbol{\Lambda},\\
&L_s : H^{\frac12}(\Sigma)^d\to\boldsymbol{\Lambda}' && \text{such that} && \langle L_s\boldsymbol{d},\boldsymbol{\mu}\rangle \overset{def}{=} c\left(\boldsymbol{\mu},\boldsymbol{d}\right) && \forall \boldsymbol{d}\in\boldsymbol{W},\forall\boldsymbol{\mu}\in\boldsymbol{\Lambda},\\
&B : \boldsymbol{V}\to Q' && \text{such that} && \langle B\boldsymbol{u},\boldsymbol{q}\rangle \overset{def}{=} (\operatorname{div}\boldsymbol{u},q)_{0,\Omega} && \forall \boldsymbol{u}\in\boldsymbol{V},\forall q\in Q.
\end{aligned} \tag{2.4}$$

In order to analyze the well-posedness of Problem 5, we consider the following functional setting. Let us introduce the Hilbert space

$$\mathbb{V} = \boldsymbol{V}\times\boldsymbol{W}\times\boldsymbol{\Lambda} \tag{2.5}$$

with norm

$$\|\boldsymbol{U}\|_{\mathbb{V}} := \left(\|\boldsymbol{u}\|_{\boldsymbol{V}}^2 + \|\boldsymbol{d}\|_{\boldsymbol{W}}^2 + \|\boldsymbol{\lambda}\|_{\boldsymbol{\Lambda}}^2\right)^{\frac12} \qquad \forall\, \boldsymbol{U} = (\boldsymbol{u},\boldsymbol{d},\boldsymbol{\lambda})\in\mathbb{V}. \tag{2.6}$$

Then we define the following bilinear forms $\mathbb{A} : \mathbb{V}\times\mathbb{V}\to\mathbb{R}$ and $\mathbb{B} : \mathbb{V}\times\boldsymbol{Q}\to\mathbb{R}$

$$\begin{aligned}\mathbb{A}(\boldsymbol{U},\boldsymbol{V}) \stackrel{def}{=} &\alpha_1(\boldsymbol{u},\boldsymbol{v})_{0,\Omega} + 2\mu(\epsilon(\boldsymbol{u}),\epsilon(\boldsymbol{v})) \\ &+ \alpha_2(\boldsymbol{d},\boldsymbol{w}) + a^{\mathrm{s}}(\boldsymbol{d},\boldsymbol{w}) + c\left(\boldsymbol{\lambda},\boldsymbol{v}\circ\overline{\boldsymbol{\phi}} - \boldsymbol{w}\right) - c\left(\boldsymbol{\mu},\boldsymbol{u}\circ\overline{\boldsymbol{\phi}} - \boldsymbol{d}\right)\end{aligned} \tag{2.7}$$

$$\mathbb{B}(\boldsymbol{U},q) \stackrel{def}{=} (\operatorname{div}\boldsymbol{u},q) \tag{2.8}$$

and the following linear form $\boldsymbol{F}:\mathbb{V}\to\mathbb{R}$

$$\langle \boldsymbol{F},\boldsymbol{V}\rangle \stackrel{def}{=} (\boldsymbol{f},\boldsymbol{v})_{0,\Omega} + (\boldsymbol{g},\boldsymbol{w})_{0,\Sigma} - c(\boldsymbol{\mu},\overline{\boldsymbol{d}}) \qquad \forall \boldsymbol{V} = \left(\boldsymbol{v},\boldsymbol{w},\boldsymbol{\mu}\right)\in\mathbb{V}. \tag{2.9}$$

Then Problem 5 can be written in abstract form as: find $(\boldsymbol{U},p)\in\mathbb{V}\times Q$ such that

$$\begin{aligned}&\mathbb{A}(\boldsymbol{U},\boldsymbol{V}) + \mathbb{B}(\boldsymbol{V},p) = \langle \boldsymbol{F},\boldsymbol{V}\rangle && \forall \boldsymbol{V}\in\mathbb{V}\\ &\mathbb{B}(\boldsymbol{U},q) = 0 && \forall q\in \boldsymbol{Q}.\end{aligned} \tag{2.10}$$

The analysis of well-posedness of problem (2.10) is performed in [23] by checking two *inf-sup* conditions (see [21, Chapter 4] and [23]).

**Proposition 2.2.1.** *Let* $\mathbb{K} \stackrel{def}{=} \{\boldsymbol{V}\in\mathbb{V} : \mathbb{B}(\boldsymbol{V},q) = 0 \quad \forall q\in Q\}$ *be the kernel of* $\mathbb{B}$*. There exist* $\alpha>0$ *and* $\beta>0$ *such that*

$$\inf_{\boldsymbol{U}\in\mathbb{K}}\sup_{\boldsymbol{V}\in\mathbb{K}} \frac{\mathbb{A}(\boldsymbol{U},\boldsymbol{V})}{\|\boldsymbol{U}\|_{\mathbb{V}}\|\boldsymbol{V}\|_{\mathbb{V}}} \geq \alpha, \qquad \inf_{q\in Q}\sup_{\boldsymbol{V}\in\mathbb{V}} \frac{\mathbb{B}(\boldsymbol{V},q)}{\|\boldsymbol{V}\|_{\mathbb{V}}\|q\|_{Q}} \geq \beta \tag{2.11}$$

Owing to proposition 2.2.1, the time-step problem is well posed and moreover we can state the following result that gives a stability bound for the solution of the stationary problem (for the proof see [21, Chapter 4])

**Proposition 2.2.2.** *There exists a unique* $(\boldsymbol{U},p)\in\mathbb{V}\times Q$ *solution of Problem 5, moreover*

$$\|\boldsymbol{U}\|_{\mathbb{V}} + \|p\|_{Q} \leq C\|\boldsymbol{F}\|_{\mathbb{V}'}. \tag{2.12}$$

## 2.3 Fully discrete problem using Stokes-stable spaces

In this section we introduce the finite element spatial discretization of Problem 4. In the following we assume that the **fluid triangulation family** $\{\mathscr{T}_h\}_{h>0}$ is **regular** and the **solid triangulation family** $\{\mathscr{S}_h\}_{h>0}$ is **quasi-uniform** (see Definition 6 in Appendix A). The finite element spaces are chosen as reported in the following.

- By $\boldsymbol{V}_h\subset\boldsymbol{V}$ and $Q_h\subset Q$ we denote the finite element spaces for the fluid velocity and pressure, respectively; these spaces are chosen such that to fulfill the inf-sup stability condition for the Stokes problem, that is, there exists $\beta'>0$ such that

$$\inf_{q_h\in Q_h}\sup_{\boldsymbol{v}_h\in\boldsymbol{V}_h} \frac{(\operatorname{div}\boldsymbol{v}_h,q_h)_{0,\Omega}}{\|\boldsymbol{v}_h\|_{1,\Omega}\|q_h\|_{0,\Omega}} \geq \beta'. \tag{2.13}$$

- For the solid displacements and the Lagrange multiplier, the finite element spaces are constructed considering

$$Y_h = \left\{ w_h \in C^0(\overline{\Sigma}) : \boldsymbol{w}_h\big|_{\boldsymbol{K}} \in \mathbb{P}_1 \quad \forall \boldsymbol{K} \in \mathscr{S}_h \right\},$$

  where $\mathbb{P}_1$ denotes the space of polynomial functions of degree at most one. Then we define the solid displacements space $\boldsymbol{W}_h \subset \boldsymbol{W}$ and the Lagrange multiplier space $\boldsymbol{\Lambda}_h \subset \boldsymbol{\Lambda}$ as follows

$$\boldsymbol{W}_h = \boldsymbol{W} \cap Y_h^d, \qquad \boldsymbol{\Lambda}_h = Y_h^d. \tag{2.14}$$

  We point out that, in general $\boldsymbol{W}_h \subset \boldsymbol{\Lambda}_h$ since the spaces $\boldsymbol{W}_h$ and $\boldsymbol{\Lambda}_h$ may differ because of the essential boundary conditions on the displacement field that are contained in the space $\boldsymbol{W}_h$.

Then the fully discrete problem reads as:

**Problem 6.** *Given $\boldsymbol{\phi}_0$, the position vector of the solid reference configuration, and $(\boldsymbol{u}_{0h}, \boldsymbol{d}_{0h}, \boldsymbol{d}_{1h}) \in \boldsymbol{V}_h \times \boldsymbol{W}_h \times \boldsymbol{W}_h$ the initial fluid velocity, solid displacement and solid velocity of the fully discrete problem, for $n \geq 1$, find $\left(\boldsymbol{u}_h^n, p_h^n, \boldsymbol{d}_h^n, \boldsymbol{\lambda}_h^n\right) \in \boldsymbol{V}_h \times \boldsymbol{Q}_h \times \boldsymbol{W}_h \times \boldsymbol{\Lambda}_h$ such that*

$$\begin{aligned}
&\rho_{\mathrm{f}}\big(\partial_\tau \boldsymbol{u}_h^n, \boldsymbol{v}_h\big)_{0,\Omega} + 2\mu\big(\boldsymbol{\epsilon}(\boldsymbol{u}_h^n), \boldsymbol{\epsilon}(\boldsymbol{v}_h)\big)_{0,\Omega} - \big(\operatorname{div} \boldsymbol{v}_h, p_h^n\big)_{0,\Omega} + c\big(\boldsymbol{\lambda}_h^n, \boldsymbol{v}_h \circ \boldsymbol{\phi}_h^{n-1}\big) = 0 && \forall \boldsymbol{v}_h \in \boldsymbol{V}_h,\\
&\rho_{\mathrm{s}}\epsilon\big(\partial_{\tau\tau} \boldsymbol{d}_h^n, \boldsymbol{w}_h\big)_{0,\Sigma} + a^{\mathrm{s}}\big(\boldsymbol{d}_h^n, \boldsymbol{w}_h\big) - c\big(\boldsymbol{\lambda}_h^n, \boldsymbol{w}_h\big) = 0 && \forall \boldsymbol{w}_h \in \boldsymbol{W}_h,\\
&c\big(\boldsymbol{\mu}_h, \boldsymbol{u}_h^n \circ \boldsymbol{\phi}_h^{n-1} - \partial_\tau \boldsymbol{d}_h^n\big) = 0 && \forall \boldsymbol{\mu}_h \in \boldsymbol{\Lambda}_h,\\
&\big(\operatorname{div} \boldsymbol{u}_h^n, q_h\big)_{0,\Omega} = 0 && \forall q_h \in \boldsymbol{Q}_h,\\
&\boldsymbol{\phi}_h^n = \boldsymbol{\phi}^0 + \boldsymbol{d}_h^n.
\end{aligned}$$

*with*

$$\boldsymbol{u}_h^0 = \boldsymbol{u}_{0h}, \qquad \boldsymbol{d}_h^0 = \boldsymbol{d}_{0h}, \qquad \boldsymbol{d}_h^{-1} = \boldsymbol{d}_{0h} - \tau \boldsymbol{d}_{1h}.$$

We observe that the bilinear form $c(\cdot,\cdot)$, originally defined on the product space $\boldsymbol{\Lambda} \times H^{\frac{1}{2}}(\Sigma)^d$, when restricted to the finite dimensional spaces can be identified with the scalar product in $L^2(\Sigma)^d$.

## 2.3.1 Stability analysis of Problem 6

The stability of the numerical scheme given by Problem 6 has been analyzed in [20] where it is possible to find the proof of the theoretical result and extended numerical simulations showing the stability properties of the algorithm. The numerical scheme has also been compared with similar algorithms that does not make use of Lagrange multiplier and that result only conditionally stable. In the following we report the theoretical result about the stability referring, for the proof, to [20].

**Proposition 2.3.1.** *For $n \geq 1$ let $(\boldsymbol{u}_h^n, \boldsymbol{d}_h^n, \boldsymbol{\lambda}_h^n, p_h^n)$ be the solution of Problem 6, then*

$$\frac{\rho_{\mathrm{f}}}{2}\|\boldsymbol{u}_h^n\|_{0,\Omega}^2 + \frac{\rho_{\mathrm{s}}\epsilon}{2}\|\partial_\tau \boldsymbol{d}_h^n\|_{0,\Sigma}^2 + \frac{1}{2}\|\boldsymbol{d}_h^n\|_{\mathrm{s}}^2 \leq \frac{\rho_{\mathrm{f}}}{2}\|\boldsymbol{u}_h^0\|_{0,\Omega}^2 + \frac{\rho_{\mathrm{s}}\epsilon}{2}\|\partial_\tau \boldsymbol{d}_h^0\|_{0,\Sigma}^2 + \frac{1}{2}\|\boldsymbol{d}_h^0\|_{\mathrm{s}}^2, \tag{2.15}$$

*where $\partial_\tau \boldsymbol{d}_h^0 = \left(\boldsymbol{d}_h^0 - \boldsymbol{d}_h^{-1}\right) \Big/ \tau$.*

### 2.3.2 Well-posedness of the time-steps of Problem 6

The well-posedness of the time-step of the fully discrete problem is studied, as in the time semi-discrete case, considering the following stationary problem

**Problem.** *Let $\overline{\boldsymbol{\phi}} = \boldsymbol{\phi}_0 + \tau \overline{\boldsymbol{d}} \in W^{1,\infty}(\Sigma)^d$ be invertible with Lipschitz inverse. Given $\boldsymbol{f} \in L^2(\Omega)^d$ and $\boldsymbol{g} \in L^2(\Sigma)^d$, find $(\boldsymbol{u}_h, p_h, \boldsymbol{d}_h, \boldsymbol{\lambda}_h) \in \boldsymbol{V}_h \times Q_h \times \boldsymbol{W}_h \times \boldsymbol{\Lambda}_h$ such that*

$$\begin{aligned}
&\alpha_1(\boldsymbol{u}_h, \boldsymbol{v}_h)_{0,\Omega} + 2\mu\big(\boldsymbol{\epsilon}(\boldsymbol{u}_h), \boldsymbol{\epsilon}(\boldsymbol{v}_h)\big)_{0,\Omega} - (p_h, \operatorname{div} \boldsymbol{v}_h)_{0,\Omega} + c\left(\boldsymbol{\lambda}_h, \boldsymbol{v}_h \circ \overline{\boldsymbol{\phi}}\right) = (\boldsymbol{f}, \boldsymbol{v}_h)_{0,\Omega} && \forall \boldsymbol{v}_h \in \boldsymbol{V}_h,\\
&\alpha_2(\boldsymbol{d}_h, \boldsymbol{w}_h) + a^{\mathrm{s}}(\boldsymbol{d}_h, \boldsymbol{w}_h) - c(\boldsymbol{\lambda}_h, \boldsymbol{w}_h) = (\boldsymbol{g}, \boldsymbol{w}_h)_{0,\Sigma} && \forall \boldsymbol{w}_h \in \boldsymbol{W}_h,\\
&c\left(\boldsymbol{\mu}_h, \boldsymbol{u}_h \circ \overline{\boldsymbol{\phi}} - \boldsymbol{d}_h\right) = -c(\boldsymbol{\mu}_h, \overline{\boldsymbol{d}}) && \forall \boldsymbol{\mu}_h \in \boldsymbol{\Lambda}_h,\\
&(\operatorname{div} \boldsymbol{u}_h, q_h)_{0,\Omega} = 0 && \forall q_h \in Q_h.
\end{aligned}$$

The stationary system just introduced has the structure of a saddle point problem and can be written in the following matrix form

$$\left(\begin{array}{ccc|c} A_{fh} & 0 & L_{fh}^\top & B_h^\top \\ 0 & A_{sh} & -L_{sh}^\top & 0 \\ L_{fh} & -L_{sh} & 0 & 0 \\ \hline B_h & 0 & 0 & 0 \end{array}\right) \left(\begin{array}{c} \boldsymbol{u}_h \\ \boldsymbol{d}_h \\ \boldsymbol{\lambda}_h \\ \hline p_h \end{array}\right) = \left(\begin{array}{c} \boldsymbol{f} \\ \boldsymbol{g} \\ \mathbf{m} \\ \hline 0 \end{array}\right), \tag{2.16}$$

with proper definition of the right hand side terms and of the sub-matrices in agreement with the operators defined in (2.4). The well-posedness of the time step problem expressed by the system (2.16) relies on the inf-sup condition (2.13), assumed for the finite element spaces $(\boldsymbol{V}_h, Q_h)$, and on the following proposition. Denoting $\boldsymbol{V}_{h,0}$ the subspace of $\boldsymbol{V}_h$ of divergence free functions, that is

$$\boldsymbol{V}_{h,0} \stackrel{def}{=} \left\{\boldsymbol{v}_h \in \boldsymbol{V}_h : \qquad (\operatorname{div} \boldsymbol{v}_h, q_h)_{0,\Omega} = 0, \qquad \forall q_h \in Q_h\right\}.$$

we have the following

**Proposition 2.3.2.** *Assume that the domain $\Omega$ is convex and that the ratio $h_f/h_{\mathrm{s}}$ is "small", then there exists $\kappa > 0$ such that*

$$\inf_{\boldsymbol{\mu}_h \in \boldsymbol{\Lambda}_h} \sup_{\boldsymbol{v}_h \in \boldsymbol{V}_{h,0}} \frac{c(\boldsymbol{\mu}_h, \boldsymbol{v}_h|_\Sigma)}{\|\boldsymbol{v}_h\|_{1,\Omega}} \geq \kappa. \tag{2.17}$$

The proof of this result is given in [19, Proposition 8]. It is immediate to obtain the following corollary

**Corollary 2.3.3.** *In the same hypotheses of Proposition 2.3.2, there exists $\kappa > 0$ such that*

$$\inf_{\boldsymbol{\mu}_h \in \boldsymbol{\Lambda}_h} \sup_{(\boldsymbol{v}_h, \boldsymbol{w}_h) \in \boldsymbol{V}_{h,0} \times \boldsymbol{W}_h} \frac{c(\boldsymbol{\mu}_h, \boldsymbol{v}_h|_\Sigma - \boldsymbol{w}_h)}{\left(\|\boldsymbol{v}_h\|_{1,\Omega}^2 + \|\boldsymbol{w}_h\|_{1,\Sigma}^2\right)^{\frac{1}{2}}} \geq \kappa. \tag{2.18}$$

Considering the structure of problem 2.16 we recognize that the upper left block constitutes a saddle point problem on the finite dimensional space $\boldsymbol{V}_{h,0}$. Owing to this observation, to the inf-sup (2.18) and to the positive definiteness of the solid and fluid matrices $A_f$ and $A_s$ (this property is related to the coercivity of the corresponding bilinear forms), we can state the following proposition whose proof can be found in [19].

**Proposition 2.3.4.** *Assume that the domain $\Omega$ is convex and the ratio $h_f/h_s$ is "small", then, for any fixed $n \geq 1$, there exist only one solution $(\boldsymbol{u}_h^n, \boldsymbol{d}_h^n, \boldsymbol{\lambda}_h^n, p_h^n)$ to problem (2.16) and, if $(\boldsymbol{u}^n, \boldsymbol{d}^n, \boldsymbol{\lambda}^n, p^n)$ is the solution at time-step n of the stationary problem 5, it holds the following optimal error estimate: for any fixed $n \geq 1$*

$$\begin{aligned}
&\|\boldsymbol{u}^n - \boldsymbol{u}_h^n\|_{1,\Omega} + \|\boldsymbol{d}^n - \boldsymbol{d}_h^n\|_{1,\Sigma} + \|\boldsymbol{\lambda}^n - \boldsymbol{\lambda}_h^n\|_{\Lambda} + \|p^n - p_h^n\|_{0,\Omega} \\
&\quad \leq C \inf_{\boldsymbol{v}_h \in \boldsymbol{V}_h, \boldsymbol{w}_h \in \boldsymbol{W}_h, \boldsymbol{\mu}_h \in \boldsymbol{\Lambda}_h, q_h \in Q_h} \big( \|\boldsymbol{u}^n - \boldsymbol{v}_h\|_{1,\Omega} + \|\boldsymbol{d}^n - \boldsymbol{w}_h\|_{1,\Sigma} + \|\boldsymbol{\lambda}^n - \boldsymbol{\mu}_h\|_{\Lambda} + \|p^n - q_h\|_{0,\Omega} \big).
\end{aligned}$$

## 2.4 Fully discrete problem using $\mathbb{P}_1 - \mathbb{P}_1$ elements

In this section we consider a finite element spatial discretization of Problem 4 that consider a couple of finite dimensional spaces for the fluid velocity and pressure that does not satisfy the basic inf-sup 2.13 but a generalized inf-sup condition (2.23). We consider, as in the case of spatial discretization with Stokes-stable elements, a **fluid triangulation family** $\{\mathscr{T}_h\}_{h>0}$ that is assumed **regular**, and a **solid triangulation family** $\{\mathscr{S}_h\}_{h>0}$ that is assumed **quasi-uniform** (see Definition 6 in Appendix A). We introduce the following piecewise polynomial spaces

$$\begin{aligned}
X_h &= \left\{ v_h \in C^0(\overline{\Omega}) : v_h\big|_{\boldsymbol{K}} \in \mathbb{P}_1 \quad \forall \boldsymbol{K} \in \mathscr{T}_h \right\}, \\
Y_h &= \left\{ w_h \in C^0(\overline{\Sigma}) : \boldsymbol{w}_h\big|_{\boldsymbol{K}} \in \mathbb{P}_1 \quad \forall \boldsymbol{K} \in \mathscr{S}_h \right\},
\end{aligned} \tag{2.19}$$

where $\mathbb{P}_1$ denotes the space of polynomial functions of degree at most one. The finite element spaces for the discrete fluid problem are defined as

$$\boldsymbol{V}_h = \boldsymbol{V} \cap X_h^d, \qquad Q_h = Q \cap X_h. \tag{2.20}$$

In order to discretize the velocity, the displacement and the Lagrange multiplier in the solid domain we use the finite element spaces already introduced in (2.14).

**Fully discrete problem** It is well known that the $\mathbb{P}_1 - \mathbb{P}_1$ couple of fluid spaces $(\boldsymbol{V}_h, Q_h)$ is not inf-sup stable (see Appendix C for details); in order to compensate this drawback we perturb the discrete problem introducing a **pressure stabilization term** of Brezzi-Pitkaranta type [29], namely

$$s_h(p_h, q_h) \stackrel{def}{=} \gamma \sum_{K \in \mathscr{T}_h} h_K^2 \int_K \nabla p_h \cdot \nabla q_h d\boldsymbol{x}, \tag{2.21}$$

where $\gamma > 0$ is a constant to be chosen in a suitable way. We introduce also the notation

$$|q_h|_{s_h}^2 \stackrel{def}{=} s_h(q_h, q_h) \qquad \forall q_h \in Q_h. \tag{2.22}$$

The analysis of the stabilized Stokes problem is performed in Appendix C; the main result concerning the stabilization term that we will use in the following analysis is given in the auxiliary Lemma 2.4.1 whose proof can be found, for example, in [21, Lemma 8.5.1].

**Lemma 2.4.1.** *Let $Q_h \subset C^0(\overline{\Omega})$, then there exist $c_1 > 0$ and $c_2 > 0$ such that for all $q_h \in Q_h$ it holds*

$$\sup_{0 \neq \boldsymbol{v}_h \in \boldsymbol{V}_h} \frac{(\operatorname{div} \boldsymbol{v}_h, q_h)_{0,\Omega}}{\|\boldsymbol{v}_h\|_{1,\Omega}} \geq c_1 \|q_h\|_{0,\Omega} - c_2 |q_h|_{s_h}. \tag{2.23}$$

Using the finite element spaces $\boldsymbol{V}_h$, $Q_h$, $\boldsymbol{W}_h$ and $\boldsymbol{\Lambda}_h$, the fully discrete stabilized problem reads:

**Problem 7.** *Given $\boldsymbol{\phi}_0$, the position vector of the solid reference configuration, and $(\boldsymbol{u}_{0h}, \boldsymbol{d}_{0h}, \boldsymbol{d}_{1h}) \in \boldsymbol{V}_h \times \boldsymbol{W}_h \times \boldsymbol{W}_h$ the initial fluid velocity, solid displacement and solid velocity of the fully discrete problem, for $n \geq 1$, find $\left(\boldsymbol{u}_h^n, p_h^n, \boldsymbol{d}_h^n, \boldsymbol{\lambda}_h^n\right) \in \boldsymbol{V}_h \times Q_h \times \boldsymbol{W}_h \times \boldsymbol{\Lambda}_h$ such that*

$$\begin{aligned}
&\rho_{\mathrm{f}}\big(\partial_\tau \boldsymbol{u}_h^n, \boldsymbol{v}_h\big)_{0,\Omega} + 2\mu\big(\boldsymbol{\epsilon}(\boldsymbol{u}_h^n), \boldsymbol{\epsilon}(\boldsymbol{v}_h)\big)_{0,\Omega} - \big(\operatorname{div} \boldsymbol{v}_h, p_h^n\big)_{0,\Omega} + c\big(\boldsymbol{\lambda}_h^n, \boldsymbol{v}_h \circ \boldsymbol{\phi}_h^{n-1}\big) = 0 && \forall \boldsymbol{v}_h \in \boldsymbol{V}_h,\\
&\rho_{\mathrm{s}}\epsilon\big(\partial_{\tau\tau} \boldsymbol{d}_h^n, \boldsymbol{w}_h\big)_{0,\Sigma} + a^{\mathrm{s}}\big(\boldsymbol{d}_h^n, \boldsymbol{w}_h\big) - c\big(\boldsymbol{\lambda}_h^n, \boldsymbol{w}_h\big) = 0 && \forall \boldsymbol{w}_h \in \boldsymbol{W}_h,\\
&c\big(\boldsymbol{\mu}_h, \boldsymbol{u}_h^n \circ \boldsymbol{\phi}_h^{n-1} - \partial_\tau \boldsymbol{d}_h^n\big) = 0 && \forall \boldsymbol{\mu}_h \in \boldsymbol{\Lambda}_h,\\
&\big(\operatorname{div} \boldsymbol{u}_h^n, q_h\big)_{0,\Omega} + s_h(p_h^n, q_h) = 0 && \forall q_h \in Q_h,\\
&\boldsymbol{\phi}_h^n = \boldsymbol{\phi}^0 + \boldsymbol{d}_h^n.
\end{aligned}$$

*with*

$$\boldsymbol{u}_h^0 = \boldsymbol{u}_{0h}, \qquad \boldsymbol{d}_h^0 = \boldsymbol{d}_{0h}, \qquad \boldsymbol{d}_h^{-1} = \boldsymbol{d}_{0h} - \tau \boldsymbol{d}_{1h}.$$

Also in this case the bilinear form $c(\cdot,\cdot)$, originally defined on the product space $\boldsymbol{\Lambda} \times H^{\frac{1}{2}}(\Sigma)^d$, when restricted to the finite dimensional spaces is identified with the scalar product in $L^2(\Sigma)^d$.

### 2.4.1 Stability analysis of Problem 7

In [20] it has been proved the unconditional stability of Problem 6 as reported in Proposition 2.3.1. In this section we analyze the stability of Problem 7 that differs from Problem 6 because of the pressure stabilization term. Using the notation

$$\dot{\boldsymbol{d}}_h^n \stackrel{def}{=} \partial_\tau \boldsymbol{d}_h^n,$$

we define the energy and the dissipation, respectively, as

$$E_h^n := \frac{\rho_{\mathrm{f}}}{2}\|\boldsymbol{u}_h^n\|_{0,\Omega}^2 + \frac{\rho_{\mathrm{s}}\epsilon}{2}\|\dot{\boldsymbol{d}}_h^n\|_{0,\Sigma}^2 + \frac{1}{2}\|\boldsymbol{d}_h^n\|_{\mathrm{s}}^2, \quad \text{for} \quad n \geq 0, \tag{2.24}$$

and

$$\begin{aligned}
D_h^n :=& \tau \sum_{k=1}^{n} \Big(2\mu \|\epsilon(\boldsymbol{u}_h^k)\|_{0,\Omega}^2 + |p_h^k|_{s_h}^2\Big)\\
&+ \tau^2 \sum_{k=1}^{n} \left(\frac{\rho_{\mathrm{f}}}{2}\|\partial_\tau \boldsymbol{u}_h^k\|_{0,\Omega}^2 + \frac{\rho_{\mathrm{s}}\epsilon}{2}\|\partial_\tau \dot{\boldsymbol{d}}_h^k\|_{0,\Sigma}^2 + \frac{1}{2}\|\partial_\tau \boldsymbol{d}_h^k\|_{\mathrm{s}}^2\right) \quad \text{for} \quad n \geq 1.
\end{aligned} \tag{2.25}$$

**Theorem 2.4.2.** *If* $(\boldsymbol{u}_h^m, p_h^m, \boldsymbol{d}_h^m, \boldsymbol{\lambda}_h^m) \in \boldsymbol{V}_h \times Q_h \times \boldsymbol{W}_h \times \boldsymbol{\Lambda}_h$ *satisfy Problem 7 for all* $m = 1, \dots, n$ *then*

$$E_h^n + D_h^n = E_h^0. \tag{2.26}$$

*Proof.* We take $(\boldsymbol{v}_h, q_h, \boldsymbol{w}_h, \boldsymbol{\mu}_h) = \tau(\boldsymbol{u}_h^m, p_h^m, \partial_\tau \boldsymbol{d}_h^m, \boldsymbol{\lambda}_h^m)$ in Problem 7, sum the equations, and obtain

$$\begin{aligned}\tau\rho_{\mathrm{f}}\big(\partial_\tau \boldsymbol{u}_h^m, \boldsymbol{u}_h^m)_{0,\Omega} + \tau a^{\mathrm{f}}\big((\boldsymbol{u}_h^m, p_h^m),(\boldsymbol{u}_h^m, p_h^m)\big) + s_h(p_h^m, p_h^m) + \tau c\big(\boldsymbol{\lambda}_h^m, \boldsymbol{u}_h^m \circ \boldsymbol{\phi}_h^{m-1} - \partial_\tau \boldsymbol{d}_h^m\big)\\ - \tau c\big(\boldsymbol{\lambda}_h^m, \boldsymbol{u}_h^m \circ \boldsymbol{\phi}_h^{m-1} - \partial_\tau \boldsymbol{d}_h^m\big) + \tau\rho_{\mathrm{s}}\epsilon\big(\partial_{\tau\tau}\boldsymbol{d}_h^m, \partial_\tau \boldsymbol{d}_h^m\big)_{0,\Sigma} + \tau a^{\mathrm{s}}\big(\boldsymbol{d}_h^m, \partial_\tau \boldsymbol{d}_h^m\big) = 0,\end{aligned} \tag{2.27}$$

where $a^{\mathrm{f}}(\cdot,\cdot)$ is given by (1.18) and $s_h(\cdot,\cdot)$ is defined by (2.21). With the notation $\dot{\boldsymbol{d}}_h^m = \partial_\tau \boldsymbol{d}_h^m$ we can write

$$\tau\rho_{\mathrm{s}}\epsilon\big(\partial_{\tau\tau}\boldsymbol{d}_h^m, \partial_\tau \boldsymbol{d}_h^m\big)_{0,\Sigma} = \tau\rho_{\mathrm{s}}\epsilon\big(\partial_\tau \dot{\boldsymbol{d}}_h^m, \dot{\boldsymbol{d}}_h^m\big)_{0,\Sigma}.$$

Using this expression in (2.27) and the following relation valid for all $\boldsymbol{a}, \boldsymbol{b}$ in an Hilbert space

$$2(\boldsymbol{a}, \boldsymbol{b}) = \|\boldsymbol{a}\|^2 + \|\boldsymbol{b}\|^2 - \|\boldsymbol{a} - \boldsymbol{b}\|^2,$$

we have

$$\begin{aligned}\frac{\rho_{\mathrm{f}}}{2}\Big(\|\boldsymbol{u}_h^m\|_{0,\Omega}^2 - \|\boldsymbol{u}_h^{m-1}\|_{0,\Omega}^2 + \|\boldsymbol{u}_h^m - \boldsymbol{u}_h^{m-1}\|_{0,\Omega}^2\Big) + 2\mu\tau\|\boldsymbol{\epsilon}(\boldsymbol{u}_h^m)\|_{0,\Omega}^2 + \tau|p_h^m|_{s_h}^2\\ + \frac{\rho_{\mathrm{s}}\epsilon}{2}\Big(\|\dot{\boldsymbol{d}}_h^m\|_{0,\Sigma}^2 - \|\dot{\boldsymbol{d}}_h^{m-1}\|_{0,\Sigma}^2 + \|\dot{\boldsymbol{d}}_h^m - \dot{\boldsymbol{d}}_h^{m-1}\|_{0,\Sigma}^2\Big)\\ + \frac{1}{2}\big(\|\boldsymbol{d}_h^m\|_{\mathrm{s}}^2 - \|\boldsymbol{d}_h^{m-1}\|_{\mathrm{s}}^2 + \|\boldsymbol{d}_h^m - \boldsymbol{d}_h^{m-1}\|_{\mathrm{s}}^2\big) = 0.\end{aligned} \tag{2.28}$$

Summing over $m = 1, \dots, n$, we have the stability relation reported in (2.26)

$$\begin{aligned}&\underbrace{\frac{\rho_{\mathrm{f}}}{2}\|\boldsymbol{u}_h^n\|_{0,\Omega}^2 + \frac{1}{2}\|\boldsymbol{d}_h^n\|_{\mathrm{s}}^2 + \frac{\rho_{\mathrm{s}}\epsilon}{2}\|\dot{\boldsymbol{d}}_h^n\|_{0,\Sigma}^2}_{E_h^n}\\ &\quad \underbrace{+\tau\sum_{m=1}^n\Big(2\mu\|\boldsymbol{\epsilon}(\boldsymbol{u}_h^m)\|_{0,\Omega}^2 + |p_h^m|_{s_h}^2\Big) + \tau^2\sum_{m=1}^n\Big(\frac{\rho_{\mathrm{f}}}{2}\|\partial_\tau \boldsymbol{u}_h^m\|_{0,\Omega}^2 + \frac{1}{2}\|\boldsymbol{d}_h^n\|_{\mathrm{s}}^2 + \frac{\rho_{\mathrm{s}}\epsilon}{2}\|\partial_\tau \dot{\boldsymbol{d}}_h^m\|_{0,\Sigma}^2\Big)}_{D_h^n}\\ &\quad = \underbrace{\Big(\frac{\rho_{\mathrm{f}}}{2}\|\boldsymbol{u}_h^0\|_{0,\Omega}^2 + \frac{1}{2}\|\boldsymbol{d}_h^0\|_{\mathrm{s}}^2 + \frac{\rho_{\mathrm{s}}\epsilon}{2}\|\dot{\boldsymbol{d}}_h^0\|_{0,\Sigma}^2\Big)}_{E_h^0}.\end{aligned} \tag{2.29}$$

□

### 2.4.2 Well-posedness of time steps of Problem 7

In [23] it is proved that each time-step of Problem 6 is well posed. We will give also the proof of the well-posedness of Problem 7 in the same spirit of the analysis carried out for Problem 6, we analyze the following time independent problem

**Problem 8.** *Let* $\overline{\boldsymbol{\phi}}_h = \boldsymbol{\phi}_0 + \tau \overline{\boldsymbol{d}}_h \in W^{1,\infty}(\Sigma)^d$ *be invertible with Lipschitz inverse. Given* $\boldsymbol{f} \in L^2(\Omega)^d$ *and* $\boldsymbol{g} \in L^2(\Sigma)^d$*, find* $(\boldsymbol{u}_h, p_h, \boldsymbol{d}_h, \boldsymbol{\lambda}_h) \in \boldsymbol{V}_h \times Q_h \times \boldsymbol{W}_h \times \boldsymbol{\Lambda}_h$ *such that*

$$
\begin{aligned}
&\alpha_1(\boldsymbol{u}_h,\boldsymbol{v}_h)_{0,\Omega} + 2\mu\big(\boldsymbol{\epsilon}(\boldsymbol{u}_h),\boldsymbol{\epsilon}(\boldsymbol{v}_h)\big)_{0,\Omega} - \big(\operatorname{div}\boldsymbol{v}_h, p_h\big)_{0,\Omega} + c\big(\boldsymbol{\lambda}_h, \boldsymbol{v}_h\circ\overline{\boldsymbol{\phi}}_h\big) = (\boldsymbol{f},\boldsymbol{v}_h)_{0,\Omega} && \forall \boldsymbol{v}_h \in \boldsymbol{V}_h,\\
&\alpha_2(\boldsymbol{d}_h,\boldsymbol{w}_h) + a^{\mathrm{s}}(\boldsymbol{d}_h,\boldsymbol{w}_h) - c(\boldsymbol{\lambda}_h,\boldsymbol{w}_h) = (\boldsymbol{g},\boldsymbol{w}_h)_{0,\Sigma} && \forall \boldsymbol{w}_h \in \boldsymbol{W}_h,\\
&c\left(\boldsymbol{\mu}_h, \boldsymbol{u}_h\circ\overline{\boldsymbol{\phi}}_h - \boldsymbol{d}_h\right) = -c\left(\boldsymbol{\mu}_h, \overline{\boldsymbol{d}}_h\right) && \forall \boldsymbol{\mu}_h \in \boldsymbol{\Lambda}_h,\\
&\big(\operatorname{div}\boldsymbol{u}_h^n, q_h\big)_{0,\Omega} + s_h(p_h^n, q_h) = 0 && \forall q_h \in Q_h.
\end{aligned}
$$

Problem 8 gives a time step of Problem 7 if

$$
\alpha_1 = \frac{\rho_f}{\tau}, \qquad \alpha_2 = \frac{\rho_s\epsilon}{\tau^2}, \qquad \boldsymbol{f} = \frac{\rho_f}{\tau}\boldsymbol{u}_h^{n-1}, \qquad \boldsymbol{g} = \frac{\rho_s\epsilon}{\tau^2}(2\boldsymbol{d}_h^{n-1} - \boldsymbol{d}_h^{n-2}),
$$
$$
\boldsymbol{d}_h = \frac{\boldsymbol{d}_h^n}{\tau}, \qquad \overline{\boldsymbol{d}}_h = \frac{\boldsymbol{d}_h^{n-1}}{\tau}, \qquad \boldsymbol{u}_h = \boldsymbol{u}_h^n, \qquad p_h = p_h^n, \qquad \boldsymbol{\lambda}_h = \boldsymbol{\lambda}_h^n.
$$

**Well-posedness of Problem 8** We consider the following functional setting. Let us introduce the space

$$
\mathbb{V}_h = \boldsymbol{V}_h \times \boldsymbol{W}_h \times \boldsymbol{\Lambda}_h \times Q_h \tag{2.30}
$$

with norm

$$
\|\boldsymbol{U}_h\|_{\mathbb{V}} \overset{def}{=} \left(\|\boldsymbol{u}_h\|_{1,\Omega}^2 + \|\boldsymbol{d}_h\|_{1,\Sigma}^2 + \|\boldsymbol{\lambda}_h\|_{\boldsymbol{\Lambda}}^2 + \|p_h\|_{0,\Omega}^2\right)^{\frac{1}{2}} \qquad \forall\, \boldsymbol{U}_h = (\boldsymbol{u}_h, \boldsymbol{d}_h, \boldsymbol{\lambda}_h, p_h) \in \mathbb{V}_h, \tag{2.31}
$$

and the following bilinear form $\mathbb{A}_h : \mathbb{V}_h \times \mathbb{V}_h \to \mathbb{R}$,

$$
\begin{aligned}
\mathbb{A}_h(\boldsymbol{U}_h,\boldsymbol{V}_h) \overset{def}{=} {} & \alpha_1(\boldsymbol{u}_h,\boldsymbol{v}_h)_{0,\Omega} + 2\mu\big(\boldsymbol{\epsilon}(\boldsymbol{u}_h),\boldsymbol{\epsilon}(\boldsymbol{v}_h)\big)_{0,\Omega} - \big(\operatorname{div}\boldsymbol{v}_h, p_h\big)_{0,\Omega} + (\operatorname{div}\boldsymbol{u}_h, q_h)_{0,\Omega} + s_h(p_h,q_h)\\
& + \alpha_2(\boldsymbol{d}_h,\boldsymbol{w}_h) + a^{\mathrm{s}}(\boldsymbol{d}_h,\boldsymbol{w}_h) + c\big(\boldsymbol{\lambda}_h, \boldsymbol{v}_h\circ\overline{\boldsymbol{\phi}}_h - \boldsymbol{w}_h\big) - c\left(\boldsymbol{\mu}_h, \boldsymbol{u}_h\circ\overline{\boldsymbol{\phi}}_h - \boldsymbol{d}_h\right)
\end{aligned} \tag{2.32}
$$

Then Problem 8 can be written, in abstract form as

$$
\text{Find} \quad \boldsymbol{U}_h \in \mathbb{V}_h \qquad \text{such that} \qquad \mathbb{A}(\boldsymbol{U}_h,\boldsymbol{V}_h) = \langle \boldsymbol{F}, \boldsymbol{V}_h\rangle \qquad \forall \boldsymbol{V}_h \in \mathbb{V}_h, \tag{2.33}
$$

where $\boldsymbol{F} : \mathbb{V}_h \to \mathbb{R}$ is defined as

$$
\langle \boldsymbol{F}, \boldsymbol{V}_h\rangle \overset{def}{=} (\boldsymbol{f},\boldsymbol{v}_h)_{0,\Omega} + (\boldsymbol{g},\boldsymbol{w}_h)_{0,\Sigma} - c\left(\boldsymbol{\mu}_h, \overline{\boldsymbol{d}}_h\right) \qquad \forall \boldsymbol{V}_h \in \mathbb{V}_h. \tag{2.34}
$$

Then, well-posedness of Problem 8 depend on the following inf-sup condition for the bilinear form $\mathbb{A}(\cdot,\cdot)$.

**Proposition 2.4.3.** *Suppose that the ratio* $\frac{h_f}{h_s}$ *is sufficiently small, then there exists a constant* $\alpha > 0$*, independent of* $h_f$ *and* $h_s$ *such that it holds true*

$$
\inf_{\boldsymbol{U}_h \in \mathbb{V}_h} \sup_{\boldsymbol{V}_h \in \mathbb{V}_h} \frac{\mathbb{A}(\boldsymbol{U}_h,\boldsymbol{V}_h)}{\|\boldsymbol{U}_h\|_{\mathbb{V}} \|\boldsymbol{V}_h\|_{\mathbb{V}}} \geq \alpha. \tag{2.35}
$$

*Proof.* We have

$$\mathbb{A}_h((\boldsymbol{u}_h, \boldsymbol{d}_h, \boldsymbol{\lambda}_h, p_h), (\boldsymbol{u}_h, \boldsymbol{d}_h, \boldsymbol{\lambda}_h, p_h) \geq \alpha_f \|\boldsymbol{u}_h\|_1^2 + \alpha_s \|\boldsymbol{d}_h\|_{1,\Sigma}^2 + |p_h|_{s_h}^2 \tag{2.36}$$

The following inf-sup condition has been proven valid for the couple of spaces considered in this section [21, 64] (see Lemma 2.4.1)

$$\sup_{\boldsymbol{v}_h \in \boldsymbol{V}_h} \frac{(\operatorname{div} \boldsymbol{v}_h, q_h)}{\|\boldsymbol{v}_h\|_1} \geq c_1 \|q_h\|_0 - c_2 |q_h|_{s_h}.$$

Hence there exists $\boldsymbol{v}_{1,h} \in \boldsymbol{V}_h$ such that

$$(\operatorname{div} \boldsymbol{v}_{1,h}, p_h) \geq c_1 \|p_h\|_0^2 - c_2 |p_h|_{s_h}^2 \tag{2.37}$$

$$\|\boldsymbol{v}_{1,h}\|_1 \leq C \|p_h\|_0 \tag{2.38}$$

and

$$\begin{aligned}
&\mathbb{A}_h((\boldsymbol{u}_h, \boldsymbol{d}_h, \boldsymbol{\lambda}_h, p_h), (-\boldsymbol{v}_{1,h}, 0, 0, 0)) \\
&\qquad = -\alpha_1 \left(\boldsymbol{u}_h, \boldsymbol{v}_{1,h}\right)_{0,\Omega} - 2\mu \left(\boldsymbol{\epsilon}(\boldsymbol{u}_h), \boldsymbol{\epsilon}(\boldsymbol{v}_{1,h})\right)_{0,\Omega} - c(\boldsymbol{\lambda}_h, \boldsymbol{v}_{1,h}|_\Sigma) + (\operatorname{div} \boldsymbol{v}_{1,h}, p_h) \\
&\qquad \geq -M_f \|\boldsymbol{u}_h\|_1 \|\boldsymbol{v}_{1,h}\|_1 - M_c \|\boldsymbol{\lambda}_h\|_{\boldsymbol{\Lambda}} \|\boldsymbol{v}_{1,h}\|_1 + c_1 \|p_h\|_0^2 - c_2 |p_h|_{s_h}^2 \\
&\qquad \geq -M_f \|\boldsymbol{u}_h\|_1 \|p_h\|_0 - M_c \|\boldsymbol{\lambda}_h\|_{\boldsymbol{\Lambda}} \|p_h\|_0 + c_1 \|p_h\|_0^2 - c_2 |p_h|_{s_h}^2.
\end{aligned} \tag{2.39}$$

In [19] it has been proved the following inf-sup condition

$$\sup_{\boldsymbol{v} \in \boldsymbol{V}_0} \frac{c(\boldsymbol{\mu}, \boldsymbol{v} \circ \overline{\boldsymbol{\phi}}_h)}{\|\boldsymbol{v}\|_{1,\Omega}} \geq \kappa \|\boldsymbol{\mu}\|_{\boldsymbol{\Lambda}} \qquad \forall \boldsymbol{\mu} \in \boldsymbol{\Lambda}, \tag{2.40}$$

where $\boldsymbol{V}_0$ denotes the subspace of $\boldsymbol{V}$ of divergence free functions

$$\boldsymbol{V}_0 \stackrel{def}{=} \{\boldsymbol{v} \in \boldsymbol{V} : \qquad (\operatorname{div} \boldsymbol{v}, q)_{0,\Omega} = 0, \qquad \forall q \in \boldsymbol{Q}\}.$$

Let $\boldsymbol{v}_2 \in \boldsymbol{V}_0$ be such that the supremum in (2.40) is attained when $\boldsymbol{\mu} = \boldsymbol{\lambda}_h$, that is

$$c(\boldsymbol{\lambda}_h, \boldsymbol{v}_2) \geq \alpha_c \|\boldsymbol{\lambda}_h\|_{\boldsymbol{\Lambda}} \|\boldsymbol{v}_2\|_1, \quad \|\boldsymbol{v}_2\|_1 = \|\boldsymbol{\lambda}_h\|_\Lambda.$$

We construct an approximation of $\boldsymbol{v}_2$ as follows. The pair $(\boldsymbol{v}_2, p_2) \in \boldsymbol{V} \times \boldsymbol{Q}$ is the solution of the following Stokes problem

$$\begin{aligned}
(\nabla \boldsymbol{v}_2, \nabla \boldsymbol{v}) - (\operatorname{div} \boldsymbol{v}, p_2) &= (\nabla \boldsymbol{v}_2, \nabla \boldsymbol{v}) \quad && \forall \boldsymbol{v} \in \boldsymbol{V} \\
(\operatorname{div} \boldsymbol{v}_2, q) &= 0 && \forall q \in \boldsymbol{Q},
\end{aligned} \tag{2.41}$$

with $p_2 = 0$. Let $(\boldsymbol{v}_{2,h}, p_{2,h}) \in V_h \times \boldsymbol{Q}_h$ be the solution of the associated discrete problem

$$\begin{aligned}
(\nabla \boldsymbol{v}_{2,h}, \nabla \boldsymbol{v}_h) - (\operatorname{div} \boldsymbol{v}_h, p_{2,h}) &= (\nabla \boldsymbol{v}_2, \nabla \boldsymbol{v}_h) \quad && \forall \boldsymbol{v}_h \in V_h \\
(\operatorname{div} \boldsymbol{v}_{2,h}, q_h) + s_h(p_{2,h}, q_h) &= 0 && \forall q_h \in \boldsymbol{Q}_h
\end{aligned} \tag{2.42}$$

We have

$$
\begin{aligned}
&\|\boldsymbol{v}_{2,h}\|_1 + \|p_{2,h}\|_0 + |p_{2,h}|_h \le C\|\boldsymbol{v}_2\|_1 \\
&\|\boldsymbol{v}_2 - \boldsymbol{v}_{2,h}\|_0 \le Ch_f\|\boldsymbol{v}_2\|_1 \\
&\|\boldsymbol{v}_{2,h}\|_1 \le C\|\boldsymbol{v}_2\|_1
\end{aligned}
\tag{2.43}
$$

Working as in the proof of [19, Prop. 16] and using an inverse inequality, we obtain that

$$
\begin{aligned}
c(\boldsymbol{\lambda}_h, \boldsymbol{v}_{2,h}) &= c(\boldsymbol{\lambda}_h, \boldsymbol{v}_2) + c(\boldsymbol{\lambda}_h, \boldsymbol{v}_{2,h} - \boldsymbol{v}_2) \\
&\ge \alpha_c\|\boldsymbol{\lambda}_h\|_{\boldsymbol{\Lambda}}\|\boldsymbol{v}_2\|_1 - C\|\boldsymbol{\lambda}_h\|_{0,\Sigma}h_f^{1/2}\|\boldsymbol{v}_2\|_1 \\
&\ge \alpha_c\|\boldsymbol{\lambda}_h\|^2_{\boldsymbol{\Lambda}}\left(\alpha_c - C\left(\frac{h_f}{h_s}\right)^{1/2}\right),
\end{aligned}
\tag{2.44}
$$

so that for $h_f/h_s$ sufficiently small we have

$$
c(\boldsymbol{\lambda}_h, \boldsymbol{v}_{2,h}) \ge \bar{\alpha}_c\|\boldsymbol{\lambda}_h\|_{\boldsymbol{\Lambda}}\|\boldsymbol{v}_2\|_{1,\Omega} = \bar{\alpha}_c\|\boldsymbol{\lambda}_h\|^2_{\boldsymbol{\Lambda}}.
$$

Let us consider $(\boldsymbol{v}_h, \boldsymbol{w}_h, \boldsymbol{\mu}_h, q_h) = (\boldsymbol{v}_{2,h}, \boldsymbol{0}, \boldsymbol{0}, 0)$ in (2.32)

$$
\begin{aligned}
&\mathbb{A}_h((\boldsymbol{u}_h, \boldsymbol{d}_h, \boldsymbol{\lambda}_h, p_h), (\boldsymbol{v}_{2,h}, 0, 0, 0)) \\
&\qquad = \alpha_1\left(\boldsymbol{u}_h, \boldsymbol{v}_{2,h}\right)_{0,\Omega} + 2\mu\left(\boldsymbol{\epsilon}(\boldsymbol{u}_h), \boldsymbol{\epsilon}(\boldsymbol{v}_{2,h})\right)_{0,\Omega} + c(\boldsymbol{\lambda}_h, \boldsymbol{v}_{2,h}) - (\operatorname{div}\boldsymbol{v}_{2,h}, p_h) \\
&\qquad \ge -M_f\|\boldsymbol{u}_h\|_1\|\boldsymbol{\lambda}_h\|_{\boldsymbol{\Lambda}} + \bar{\alpha}_c\|\boldsymbol{\lambda}_h\|^2_{\boldsymbol{\Lambda}} + s_h(p_{2,h}, p_h) \\
&\qquad \ge -M_f\|\boldsymbol{u}_h\|_1\|\boldsymbol{\lambda}_h\|_{\boldsymbol{\Lambda}} + \bar{\alpha}_c\|\boldsymbol{\lambda}_h\|^2_{\boldsymbol{\Lambda}} - |p_{2,h}|_h|p_h|_h \\
&\qquad \ge -M_f\|\boldsymbol{u}_h\|_1\|\boldsymbol{\lambda}_h\|_{\boldsymbol{\Lambda}} + \bar{\alpha}_c\|\boldsymbol{\lambda}_h\|^2_{\boldsymbol{\Lambda}} - C\|\boldsymbol{\lambda}_h\|_{\boldsymbol{\Lambda}}|p_h|_h.
\end{aligned}
\tag{2.45}
$$

Summing up (2.36), (2.39) and (2.45) we obtain

$$
\begin{aligned}
&\mathbb{A}_h((\boldsymbol{u}_h,\boldsymbol{d}_h,\boldsymbol{\lambda}_h,p_h),(\boldsymbol{u}_h+a\boldsymbol{v}_{1,h}+b\boldsymbol{v}_{2,h},\boldsymbol{d}_h,-\boldsymbol{\lambda}_h,-p_h))\\
&\quad\geq \alpha_f\|\boldsymbol{u}_h\|_1^2+\alpha_s\|\boldsymbol{d}_h\|_{1,\Sigma}^2+|p_h|_{s_h}^2\\
&\quad -aM_f\|\boldsymbol{u}_h\|_1\|p_h\|_1-aM_c\|\boldsymbol{\lambda}_h\|_{\boldsymbol{\Lambda}}\|p_h\|_0+ac_1\|p_h\|_0^2-ac_2|p_h|_{s_h}^2\\
&\quad -bM_f\|\boldsymbol{u}_h\|_1\|\boldsymbol{\lambda}_h\|_{\boldsymbol{\Lambda}}+b\bar{\alpha}_c\|\boldsymbol{\lambda}_h\|_{\boldsymbol{\Lambda}}^2-bC\|\boldsymbol{\lambda}_h\|_{\boldsymbol{\Lambda}}|p_h|_h\\
&\quad\geq \alpha_f\|\boldsymbol{u}_h\|_1^2+\alpha_s\|\boldsymbol{d}_h\|_{1,\Sigma}^2+|p_h|_{s_h}^2\\
&\quad -\frac{\delta_1M_f^2}{2}\|\boldsymbol{u}_h\|_1^2-\frac{a^2}{2\delta_1}\|p_h\|_0^2-\frac{\delta_2M_c^2}{2}\|\boldsymbol{\lambda}_h\|_{\boldsymbol{\Lambda}}^2-\frac{a^2}{2\delta_2}\|p_h\|_0^2\\
&\quad +ac_1\|p_h\|_0^2-ac_2|p_h|_{s_h}^2-\frac{\delta_3M_f^2}{2}\|\boldsymbol{u}_h\|_1^2-\frac{b^2}{2\delta_3}\|\boldsymbol{\lambda}_h\|_{\boldsymbol{\Lambda}}^2\\
&\quad +b\bar{\alpha}_c\|\boldsymbol{\lambda}_h\|_{\boldsymbol{\Lambda}}^2-\frac{\delta_4}{2}|p_h|_h^2-\frac{b^2C^2}{2\delta_4}\|\boldsymbol{\lambda}_h\|_{\boldsymbol{\Lambda}}^2\\
&\quad =\left(\alpha_f-\frac{\delta_1M_f^2}{2}-\frac{\delta_3M_f^2}{2}\right)\|\boldsymbol{u}_h\|_1^2\\
&\quad +\alpha_s\|\boldsymbol{d}_h\|_{1,\Sigma}^2+\left(1-ac_2-\frac{\delta_4}{2}\right)|p_h|_h^2\\
&\quad +a\left(c_1-\frac{a}{2\delta_1}-\frac{a}{2\delta_2}\right)\|p_h\|_0^2\\
&\quad +b\left(\bar{\alpha}_c-\frac{\delta_2M_c^2}{2}-\frac{b}{2\delta_3}-\frac{bC^2}{2\delta_4}\right)\|\boldsymbol{\lambda}_h\|_{\boldsymbol{\Lambda}}^2
\end{aligned}
$$

We can choose $\delta_i$ for $i=1,\dots,4$ such that

$$
\begin{aligned}
&\alpha_f-\frac{\delta_1M_f^2}{2}-\frac{\delta_3M_f^2}{2}=\bar{\beta}\\
&1-\frac{\delta_4}{2}=1/2\\
&\bar{\alpha}_c-\frac{\delta_2M_c^2}{2}=\frac{\bar{\alpha}_c}{2}.
\end{aligned}
$$

Then $a$ and $b$ can be taken sufficiently small so that

$$
\begin{aligned}
&c_1-\frac{a}{2\delta_1}-\frac{a}{2\delta_2}>0\\
&\frac{\bar{\alpha}_c}{2}-\frac{b}{2\delta_3}-\frac{bC^2}{2\delta_4}>0
\end{aligned}
$$

□

## 2.5 Regularity assumptions for the solution of Problem 3

The existence of a solution to Problem 3 is still an open question. Nevertheless we assume its existence and its regularity, in order to perform convergence analysis of the finite element solution. In this paragraph we will present some heuristic considerations that are based on the physics intuition. First of all we note that the pressure field has to be discontinuous across the structure, in fact, if not the force orthogonal to the solid is zero. Then, the named discontinuity, lies on a line in $\mathbb{R}^2$ or a surface in $\mathbb{R}^3$; this means that **the pressure field** $p$ can not be better that $p(t) \in L^2_0(\Omega) \cap H^l(\Omega)$ with $0 \le l < \frac{1}{2}$. From this observation, using the equation of momentum conservation of the fluid

$$\rho_{\mathrm{f}} \partial_t \boldsymbol{u} - 2\mu\Delta\boldsymbol{u} + \nabla p = \mathbf{0}$$

it results that **the velocity field** $\boldsymbol{u}$ is expected to be $\boldsymbol{u}(t) \in H^1_0(\Omega)^d \cap H^{1+l}(\Omega)^d$. Moreover we assume that the solid displacements and velocities $\boldsymbol{d}$ and $\partial_t \boldsymbol{d}$ are functions in $\boldsymbol{D}^{\mathrm{s}} \subset H^1(\Sigma)^d \cap H^{1+m}(\Sigma)^d$ with $m > 0$ such that $\boldsymbol{L}^{\mathrm{s}}\boldsymbol{d} \in L^2(\Sigma)^d$. $\boldsymbol{\lambda}$ represents the force acting on the solid and it is proportional to the jump of the symmetric gradient of $\boldsymbol{u}$ across $\Sigma$. Supposing that $\boldsymbol{u}$ restricted to $\Omega^+$ and $\Omega^-$ belongs to $H^{\frac{3}{2}+\delta}(\Omega^+)^d$ and $H^{\frac{3}{2}+\delta}(\Omega^-)^d$, respectively, with $0 \le \delta < \frac{1}{2}$, we conclude that $[\![\boldsymbol{\sigma}\boldsymbol{n}_{\Sigma(t)}]\!] = [\![(-p\boldsymbol{I} + 2\mu\boldsymbol{\epsilon}(\boldsymbol{u}))\boldsymbol{n}_{\Sigma(t)}]\!] \in H^z(\Sigma(t))^d$ with $z \overset{def}{=} \min\{l, \delta\}$, then we infer that $\boldsymbol{\lambda} \in H^z(\Sigma)^d$. In conclusion we assume the following regularity in the space variable

$$\boldsymbol{u}(t) \in H^{l+1}(\Omega)^d, \qquad p(t) \in H^l(\Omega), \qquad \boldsymbol{\lambda}(t) \in H^z(\Sigma)^d, \qquad \boldsymbol{d}(t) \in \boldsymbol{D}^{\mathrm{s}}, \tag{2.46}$$

where $\boldsymbol{D}^{\mathrm{s}} \subset H^{1+m}(\Sigma)^d$ and $0 \le l < \frac{1}{2}$ , $0 \le m$, $0 \le z < \frac{1}{2}$.

## 2.6 Numerical experiments

In this section we report numerical tests in order to explore the stability of the monolithic scheme. The results of these experiments have been presented in [20], where the monolithic scheme with Lagrange multipliers has been compared with a monolithic scheme that does not uses the multipliers technique.

In all simulations reported in this section the velocity and pressure spaces are discretized using the $\mathbb{P}_1 - \mathbb{P}_1$ stabilized couple, whereas in [20] it was used the enhanced Bercovier-Pironneau $\mathbb{P}_1 - iso - \mathbb{P}_2/(\mathbb{P}_1 + \mathbb{P}_0)$ element introduced in [18]. The numerical experiments are performed using the algorithm associated to Problem 7 whose matrix form is presented in the following

Given the initial values $\boldsymbol{\phi}^0_h = \boldsymbol{\phi}_0 + \boldsymbol{d}_{0h}, \boldsymbol{d}^0_h = \boldsymbol{d}_{0h}, \boldsymbol{d}^{-1}_h = d_{0h} - \tau\boldsymbol{d}_{1h}, \boldsymbol{u}^0 = \boldsymbol{u}_{0h}$, for $n \ge 1$ find $(\boldsymbol{u}^n_h, p^n_h, \boldsymbol{d}^n_h, \boldsymbol{\lambda}^n_h) \in \boldsymbol{V}_h \times \boldsymbol{Q}_h \times \boldsymbol{W}_h \times \boldsymbol{\Lambda}_h$ such that

$$\left(\begin{array}{cc|c|c} A_f & B^\top & 0 & L_f^\top(\boldsymbol{\phi}^{n-1}_h) \\ B & S & 0 & 0 \\ \hline 0 & 0 & A_s & -L_s^\top \\ \hline L_f(\boldsymbol{\phi}^{n-1}_h) & 0 & -\frac{1}{\tau}L_s & 0 \end{array}\right) \left(\begin{array}{c} \boldsymbol{u}^n_h \\ p^n_h \\ \hline \boldsymbol{d}^n_h \\ \hline \boldsymbol{\lambda}^n_h \end{array}\right) = \left(\begin{array}{c} \boldsymbol{f} \\ 0 \\ \hline \boldsymbol{g} \\ \hline \mathbf{m} \end{array}\right). \tag{2.47}$$

Denoting by $\varphi, \psi, \chi$ and $\zeta$ the basis functions respectively in $\boldsymbol{V}_h, Q_h, \boldsymbol{W}_h, \boldsymbol{\Lambda}_h$, the sub-matrices in 2.47 have the following expressions

$$
\begin{aligned}
&A_f = \frac{\rho_\mathrm{f}}{\tau} M_\mathrm{f} + K_\mathrm{f}, \quad \text{with} \quad (M_\mathrm{f})_{ij} = \left(\varphi_j, \varphi_i\right)_{0,\Omega}, \quad (K_\mathrm{f})_{ij} = \left(\boldsymbol{\epsilon}\varphi_j, \boldsymbol{\epsilon}\varphi_i\right)_{0,\Omega}, \\
&(B)_{ij} = -\left(\operatorname{div}\varphi_j, \psi_i\right)_{0,\Omega}, \\
&(S)_{ij} = \gamma \sum_{K \in \mathscr{T}_h} h_K^2 \left(\nabla\psi_j, \nabla\psi_i\right)_{0,K}, \\
&A_\mathrm{s} = \frac{\rho_\mathrm{s}\epsilon}{\tau^2} M_\mathrm{s} + K_\mathrm{s}, \quad \text{with} \quad (M_\mathrm{s})_{ij} = \left(\chi_j, \chi_i\right)_{0,\Sigma}, \quad (K_\mathrm{s})_{ij} = a^s\left(\chi_j, \chi_i\right), \\
&(L_f(\boldsymbol{\phi}_h^{n-1}))_{ij} = c(\boldsymbol{\zeta}_j, \boldsymbol{\varphi}_i \circ \boldsymbol{\phi}_h^{n-1}) = (\boldsymbol{\zeta}_j, \boldsymbol{\varphi}_i \circ \boldsymbol{\phi}_h^{n-1})_{0,\Omega}, \qquad (2.48)\\
&(L_\mathrm{s})_{ij} = c(\boldsymbol{\zeta}_j, \boldsymbol{\chi}_i) = (\boldsymbol{\zeta}_j, \boldsymbol{\chi}_i)_{0,\Sigma}, \\
&\boldsymbol{f} = \frac{\rho_\mathrm{f}}{\tau} M_\mathrm{f} \boldsymbol{u}_h^{n-1}, \\
&\boldsymbol{g} = \frac{\rho_\mathrm{s}\epsilon}{\tau^2} M_\mathrm{s}(2\boldsymbol{d}_h^{n-1} - \boldsymbol{d}_h^{n-2}), \\
&\boldsymbol{m} = -\frac{1}{\tau} L_\mathrm{s} \boldsymbol{d}_h^{n-1}.
\end{aligned}
$$

The effective computation of the sub-matrices in (2.47) is performed in Appendix D for the case of an elliptical structures immersed in the fluid.

The goal of the numerical experiments presented here is to confirm the unconditional stability of the monolithic scheme (2.47). To this aim, we perform a series of tests using the classical benchmark problem of an ellipsoidal structure that evolves to a circular equilibrium position. The ellipsoid is centered at the midpoint of the square $[0,1]\times[0,1]$ that is the physical domain and the fluid and the solid have zero initial velocities. The time evolution of the system is represented in Figure 2.2 with the purpose to give an idea of the scheme behavior. In order to verify the stability we perform long time simulations with fixed solid and fluid mesh sizes, varying the time step size. The numerical simulations have been run for a total time $T = 180s$, considering the following discretization parameters:

$$
h = h_\mathrm{s} = 1/64, \qquad \tau \in \{0.04, 0.08, 0.16, 0.32, 0.64, 1.28\}.
$$

Then for each time step $\tau$, we plot the total energy versus the time $t = n\tau < T$

$$
E_{tot} = \|\boldsymbol{u}_h^n\|_{0,\Omega}^2 + \|\dot{\boldsymbol{d}}_h^n\|_{0,\Sigma}^2 + \|\boldsymbol{d}_h^n\|_{s,\Sigma}^2.
$$

The following diagrams show that the total energy remains bounded during the simulation in all the cases, and in fact, it decreases, reflecting the unconditional stability of the monolithic algorithm predicted by the theoretical results given in [20] and reported in Theorem 2.4.2.

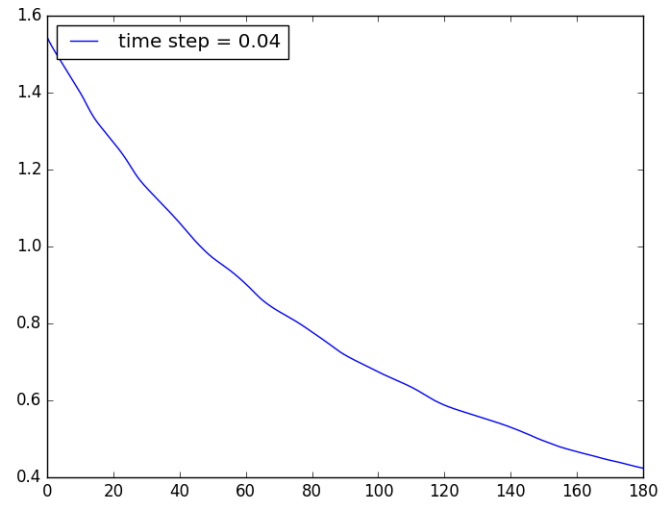

(a) $\tau$ = 0.04 s, $h_f = h_s = 1/64$

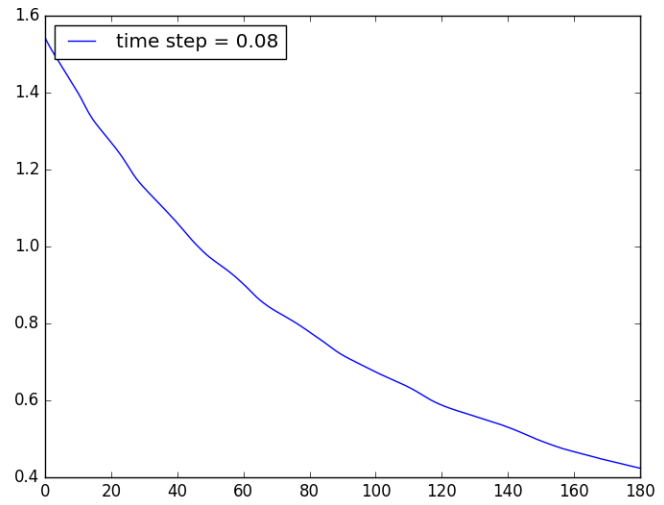

(b) $\tau$ = 0.08 s, $h_f = h_s = 1/64$

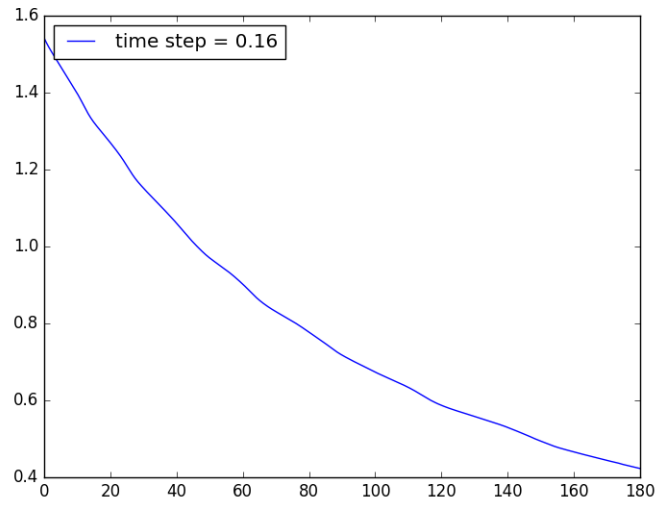

(c) $\tau$ = 0.16 s, $h_f = h_s = 1/64$

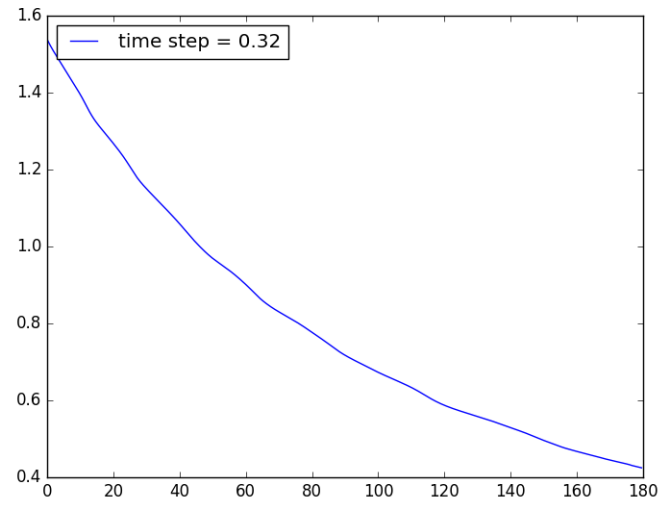

(d) $\tau$ = 0.32 s, $h_f = h_s = 1/64$

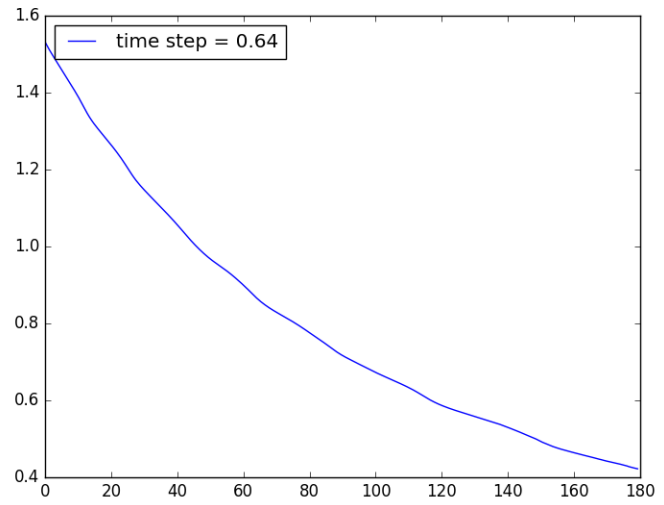

(e) $\tau$ = 0.64 s, $h_f = h_s = 1/64$

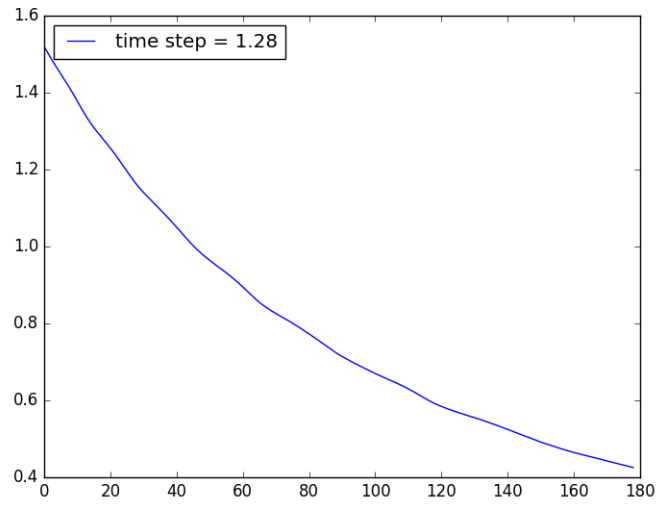

(f) $\tau$ = 1.28 s, $h_f = h_s = 1/64$

Figure 2.1: Stability Analysis of the Monolithic algorithm

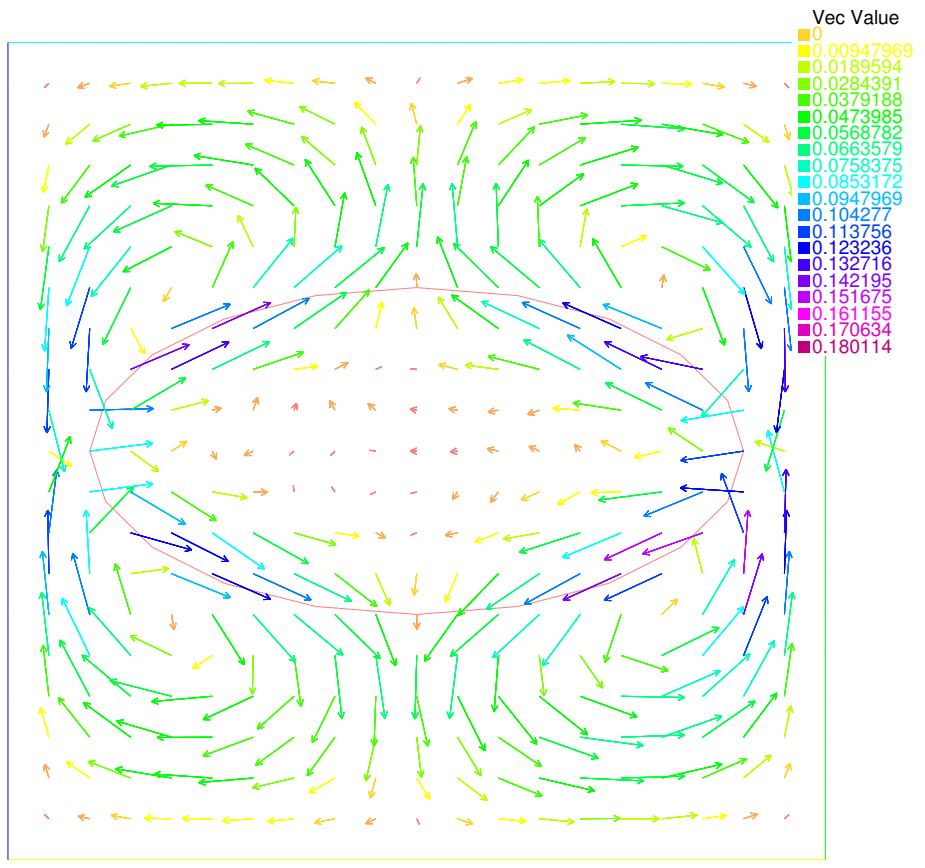

(a) t=0.1 s

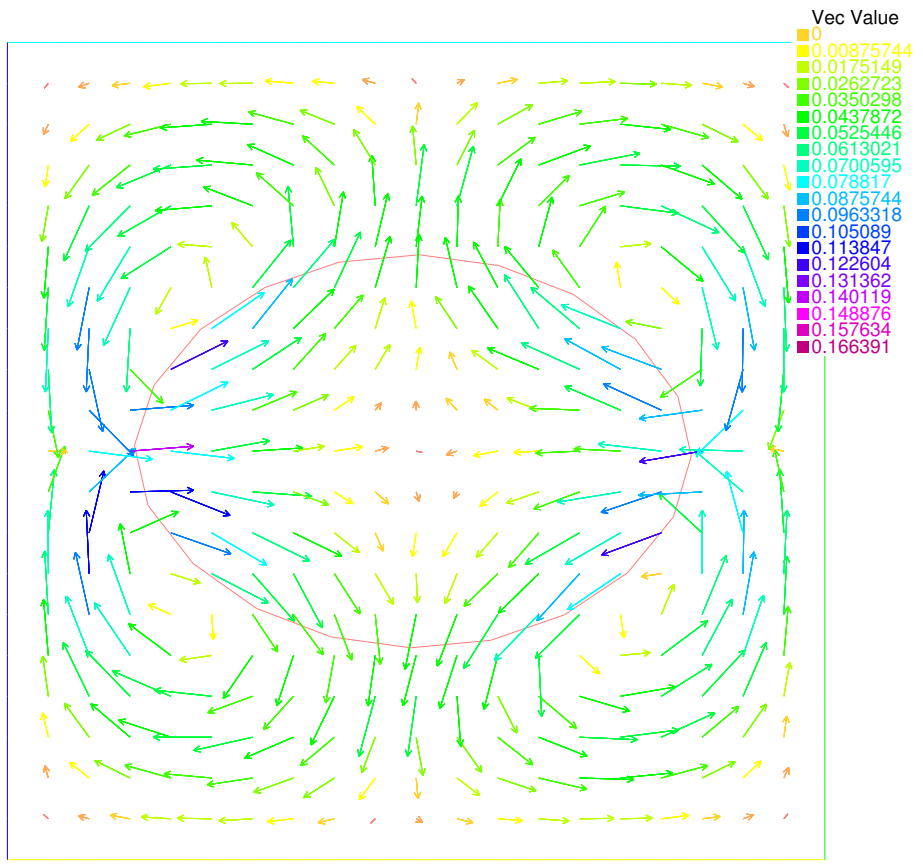

(b) t=0.6 s

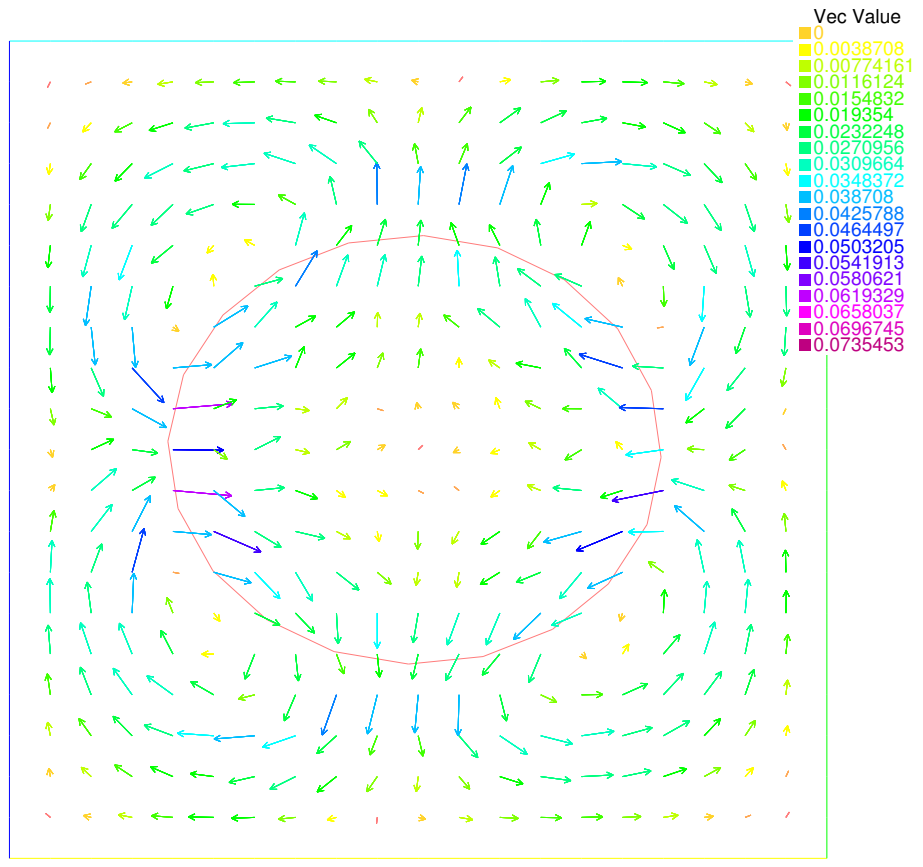

(c) t=1.2 s

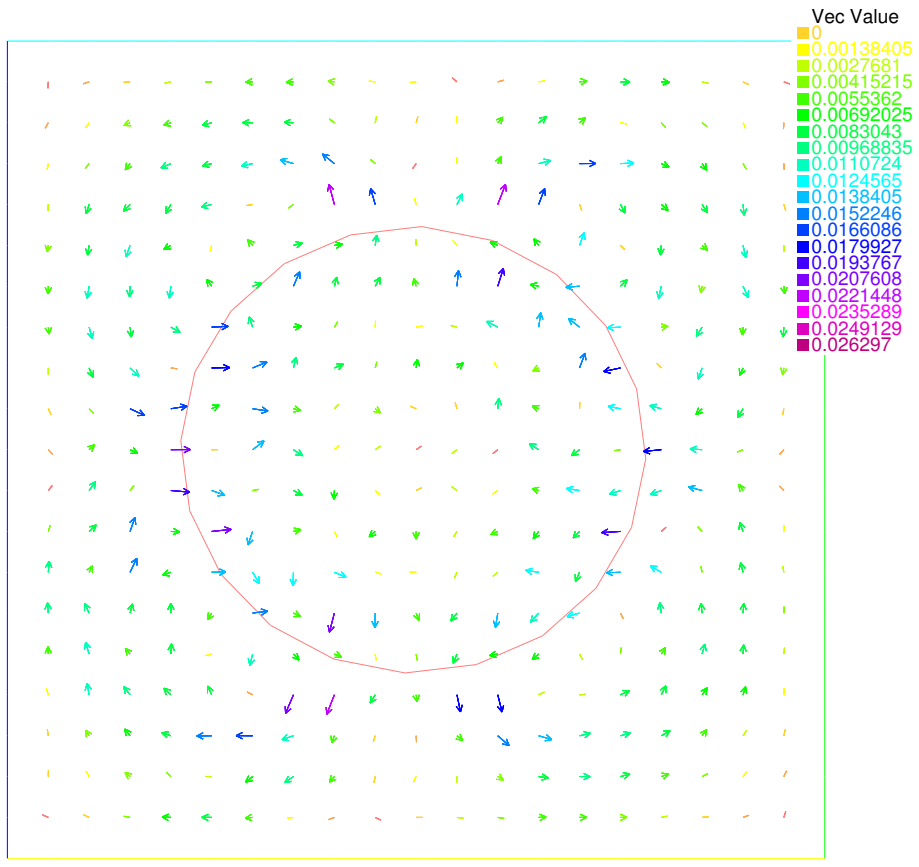

(d) t=1.8 s

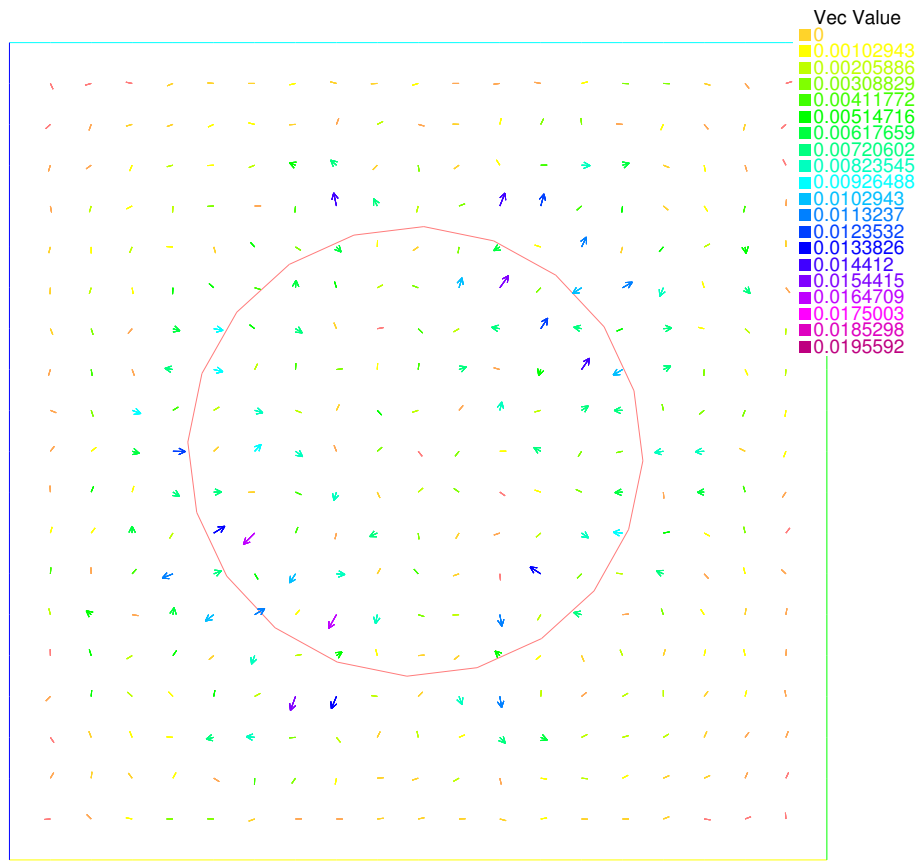

(e) t=2.4 s

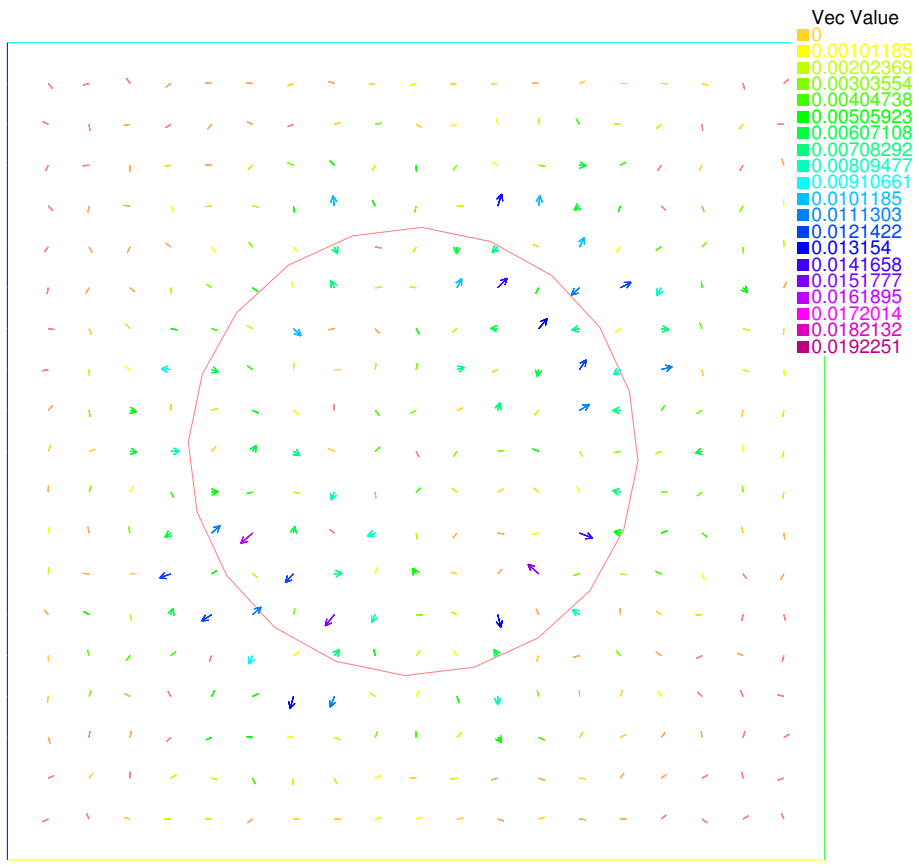

(f) t=3.0 s

Figure 2.2: Test case evolution with the monolithic algorithm

# Chapter 3

# Convergence analysis of Problem 7 - the linearized case

## 3.1 Introduction

In this chapter we analyze the convergence of the space-time scheme given in Problem 7 of Chapter 2. As said before, this scheme has been proposed in [20] where it is proved its unconditional stability and it is showed the well-posedness of the stationary problem solved in each time step in the case of Stokes stable spatial discretization. In Chapter 2 we performed also the analysis of the stability and well-posedness of the problem discretized using $\mathbb{P}_1-\mathbb{P}_1$ stabilized finite elements for the fluid. In this chapter we are interested in the analysis of the convergence of the fully discrete problem in the linearized case. More precisely we assume that the structural displacements are very small so that the solid deformed configuration can be thought as coincident with the undeformed configuration, that is $\Sigma(t)=\Sigma$ for all $t\in[0,T]$.

## 3.2 Linearized Monolithic scheme

The linearized version of the continuous Problem 2, obtained considering small displacements from the reference configuration, is reported in the following Problem 9; the current configuration is identified with the reference configuration, whose position is given by the map $\boldsymbol{\phi}_0$. We adopt the notation $\boldsymbol{v}|_\Sigma \stackrel{def}{=} \boldsymbol{v}\circ\boldsymbol{\phi}_0$ for all $\boldsymbol{v}\in\boldsymbol{V}$

**Problem 9.** *For $t\geq 0$, find $\big(\boldsymbol{u}(t),p(t),\boldsymbol{d}(t),\dot{\boldsymbol{d}}(t),\boldsymbol{\lambda}(t)\big)\in\boldsymbol{V}\times Q\times\boldsymbol{W}\times\boldsymbol{W}\times\boldsymbol{\Lambda}$ such that*

- $\partial_t\boldsymbol{d}=\dot{\boldsymbol{d}}$
- *for all $(\boldsymbol{v},q,\boldsymbol{w},\boldsymbol{\mu})\in\boldsymbol{V}\times Q\times\boldsymbol{W}\times\boldsymbol{\Lambda}$*

$$\rho_{\mathrm{f}}\big(\partial_t\boldsymbol{u},\boldsymbol{v}\big)_{0,\Omega}+a^{\mathrm{f}}\big((\boldsymbol{u},p),(\boldsymbol{v},q)\big)+c\big(\boldsymbol{\lambda},\boldsymbol{v}|_\Sigma-\boldsymbol{w}\big)\\ -c\big(\boldsymbol{\mu},\boldsymbol{u}|_\Sigma-\dot{\boldsymbol{d}}\big)+\rho_{\mathrm{s}}\epsilon\big(\partial_t\dot{\boldsymbol{d}},\boldsymbol{w}\big)_{0,\Sigma}+a^{\mathrm{s}}\big(\boldsymbol{d},\boldsymbol{w}\big)=0 \quad (3.1)$$

- $\boldsymbol{u}(0)=\boldsymbol{u}_0$ *in* $\Omega$, $\dot{\boldsymbol{d}}(0)=\dot{\boldsymbol{d}}_0$, $\boldsymbol{d}(0)=\boldsymbol{d}_0$ *in* $\Sigma$,

where $\boldsymbol{u}_0$ and $\dot{\boldsymbol{d}}_0$ are the fluid and solid initial velocities assigned such that $\boldsymbol{u}_0 \circ \boldsymbol{\phi}_0 = \dot{\boldsymbol{d}}_0$ and $\boldsymbol{d}_0$ is the initial solid displacement. Moreover we recall that

$$a^{\mathrm{f}}((\boldsymbol{u},p);(\boldsymbol{v},q)) \stackrel{def}{=} 2\mu(\boldsymbol{u},\boldsymbol{v})_{0,\Omega} - (\operatorname{div}\boldsymbol{v},p)_{0,\Omega} + (\operatorname{div}\boldsymbol{u},q)_{0,\Omega} \qquad \forall \boldsymbol{u},\boldsymbol{v} \in \boldsymbol{V}, \forall p,q \in Q.$$

The corresponding monolithic scheme is reported in the following Algorithm 1.

**Algorithm 1** Linearized monolithic scheme from [20] (strongly coupled).

Given the position vector of the solid reference configuration $\boldsymbol{\phi}_0$ and the initial data $\boldsymbol{u}_h^0$, $\boldsymbol{d}_h^0$, $\dot{\boldsymbol{d}}_h^0$, assigned such that $\boldsymbol{u}_h^0 \circ \boldsymbol{\phi}_0 = \dot{\boldsymbol{d}}_h^0$;
for $n \geq 1$, find $(\boldsymbol{u}_h^n, p_h^n, \boldsymbol{d}_h^n, \dot{\boldsymbol{d}}_h^n, \boldsymbol{\lambda}_h^n) \in \boldsymbol{V}_h \times Q_h \times \boldsymbol{W}_h \times \boldsymbol{W}_h \times \boldsymbol{\Lambda}_h$ with $\partial_\tau \boldsymbol{d}_h^n = \dot{\boldsymbol{d}}_h^n$ and such that

$$\begin{aligned}\rho_{\mathrm{f}}\big(\partial_\tau \boldsymbol{u}_h^n, \boldsymbol{v}_h)_{0,\Omega} + a^{\mathrm{f}}\big((\boldsymbol{u}_h^n, p_h^n),(\boldsymbol{v}_h,q_h)\big) + s_h(p_h^n,q_h) + c\big(\boldsymbol{\lambda}_h^n, \boldsymbol{v}_h \circ \boldsymbol{\phi}_0 - \boldsymbol{w}_h\big)\\ - c\big(\boldsymbol{\mu}_h, \boldsymbol{u}_h^n \circ \boldsymbol{\phi}_0 - \dot{\boldsymbol{d}}_h^n\big) + \rho_{\mathrm{s}}\epsilon\big(\partial_\tau \dot{\boldsymbol{d}}_h^n, \boldsymbol{w}_h\big)_{0,\Sigma} + a^{\mathrm{s}}\big(\boldsymbol{d}_h^n, \boldsymbol{w}_h\big) = 0 \quad (3.2)\end{aligned}$$

for all $(\boldsymbol{v}_h, q_h, \boldsymbol{w}_h, \boldsymbol{\mu}_h) \in \boldsymbol{V}_h \times Q_h \times \boldsymbol{W}_h \times \boldsymbol{\Lambda}_h$.

We recall that the stabilization term $s_h(\cdot,\cdot)$ and the bilinear form $c(\cdot,\cdot)$ have been introduced in (2.21) and (1.15) respectively. The bilinear form $a^s(\cdot,\cdot)$ represents the weak formulation of the solid problem introduced in (1.11).

## 3.3 Convergence analysis

For the convergence analysis of Algorithm 1 we need some preliminary results from literature. We report here results that will be useful in the following

### 3.3.1 Tools for convergence analysis

We will use extensively a well known discrete version of Gronwall Lemma (see [76, Lemma 5.1])

**Lemma 3.3.1.** *Discrete Gronwall Lemma*
*Let $\tau, B, a_m, b_m, c_m, \gamma_m$ (for integers $m \geq 1$) be non negative numbers such that*

$$a_n + \tau \sum_{m=1}^{n} b_m \leq \tau \sum_{m=1}^{n} \gamma_m a_m + \tau \sum_{m=1}^{n} c_m + B \qquad for \quad n \geq 1$$

*Suppose that $\tau\gamma_m < 1$ for all $m \geq 1$. Then, there holds*

$$a_n + \tau \sum_{m=1}^{n} b_m \leq exp\left(\tau \sum_{m=1}^{n} \frac{\gamma_m}{1-\tau\gamma_m}\right)\left(\tau \sum_{m=1}^{n} c_m + B\right) \qquad for \quad n \geq 1.$$

In the next proposition we report some useful estimates on the finite difference approximation of the time derivatives.

**Proposition 3.3.2.** *Let $(\boldsymbol{X},(\cdot,\cdot)_{\boldsymbol{X}})$ be a real Hilbert space, then*

- *For all $\boldsymbol{v} \in H^2(0,T;\boldsymbol{X})$*

$$\tau\|\partial_\tau \boldsymbol{v}^n - \partial_t \boldsymbol{v}^n\|_{\boldsymbol{X}} \le \tau^{\frac{3}{2}} \|\partial_{tt}\boldsymbol{v}\|_{L^2(t_{n-1},t_n;\boldsymbol{X})}. \tag{3.3}$$

- *For all $\boldsymbol{v} \in H^1(0,T;\boldsymbol{X})$*

$$\tau\|\partial_\tau \boldsymbol{v}^n\|_{\boldsymbol{X}} \le \tau^{\frac{1}{2}} \|\partial_t\boldsymbol{v}\|_{L^2(t_{n-1},t_n;\boldsymbol{X})}. \tag{3.4}$$

*Proof.* Let us begin with the first claim. Let be given a function $\boldsymbol{v} \in H^2(0,T;\boldsymbol{X})$, then, using the Taylor expansions with reminder in integral form, we can write

$$\begin{aligned}
\tau\|\partial_\tau \boldsymbol{v}^n - \partial_t \boldsymbol{v}^n\|_{\boldsymbol{X}} &= \left\| \int_{t_{n-1}}^{t_n} (t-t_{n-1})\partial_{tt}\boldsymbol{v} dt \right\|_{\boldsymbol{X}} \le \int_{t_{n-1}}^{t_n} \left\| (t-t_{n-1})\partial_{tt}\boldsymbol{v} \right\|_{\boldsymbol{X}} dt \\
&\le \left( \int_{t_{n-1}}^{t_n} (t-t_{n-1})^2 dt \right)^{\frac{1}{2}} \left( \int_{t_{n-1}}^{t_n} \|\partial_{tt}\boldsymbol{v}\|_{\boldsymbol{X}}^2 dt \right)^{\frac{1}{2}} \le \tau^{\frac{3}{2}} \|\partial_{tt}\boldsymbol{v}\|_{L^2(t_{n-1},t_n;\boldsymbol{X})}.
\end{aligned} \tag{3.5}$$

The proof of the second claim is a direct consequence of the fundamental Theorem of calculus

$$\tau\|\partial_\tau \boldsymbol{v}^n\|_{\boldsymbol{X}} = \|\boldsymbol{v}^n - \boldsymbol{v}^{n-1}\|_{\boldsymbol{X}} \le \int_{t_{n-1}}^{t_n} \|\partial_t \boldsymbol{v}\|_{\boldsymbol{X}} dt \le \tau^{\frac{1}{2}} \|\partial_t\boldsymbol{v}\|_{L^2(t_{n-1},t_n;\boldsymbol{X})}. \tag{3.6}$$

□

### 3.3.2 Ritz Projection Operators

In this section we define and analyze the Ritz projection operators that we will use in the following for the convergence analysis of the numerical schemes. Essentially the definitions of the projection operators for the fluid velocity and pressure are based on the Stokes Problem; the projection operators for the solid displacement and velocity are obtained using the solid bilinear form $a^{\mathrm{s}}(\cdot,\cdot)$ and we use the $L^2$-projection for the Lagrange multiplier.

**Definition 1.** *Given $\left(\boldsymbol{u},p,\boldsymbol{d},\lambda\right) \in \boldsymbol{V} \times Q \times \boldsymbol{W} \times \boldsymbol{\Lambda}$, we define the following projection operators*

- *Solid displacement and velocity projection*

$$\boldsymbol{\Pi}_{\boldsymbol{W}} : \boldsymbol{W} \longrightarrow \boldsymbol{W}_h$$

*such that for all $\boldsymbol{d} \in \boldsymbol{W}$*

$$a^{\mathrm{s}}(\boldsymbol{\Pi}_{\boldsymbol{W}}\boldsymbol{d}, \boldsymbol{w}_h) = a^{\mathrm{s}}(\boldsymbol{d}, \boldsymbol{w}_h) \qquad \forall \boldsymbol{w}_h \in \boldsymbol{W}_h, \tag{3.7}$$

- *Fluid velocity and pressure projections*

$$\boldsymbol{\Pi}_{\boldsymbol{V}} : \boldsymbol{V} \times Q \longrightarrow \boldsymbol{V}_h, \qquad \boldsymbol{\Pi}_Q : \boldsymbol{V} \times Q \longrightarrow Q_h$$

*such that for all $(\boldsymbol{u},p) \in \boldsymbol{V} \times Q$*

$$\begin{aligned}
a^{\mathrm{f}}\big((\boldsymbol{\Pi}_{\boldsymbol{V}}(\boldsymbol{u},p), \boldsymbol{\Pi}_Q(\boldsymbol{u},p)); (\boldsymbol{v}_h, q_h)\big) &+ s_h(\boldsymbol{\Pi}_Q(\boldsymbol{u},p), q_h) \\
&= a^{\mathrm{f}}((\boldsymbol{u},\boldsymbol{p}); (\boldsymbol{v}_h, q_h)) \quad \forall (\boldsymbol{v}_h, q_h) \in \boldsymbol{V}_h \times \boldsymbol{Q},
\end{aligned} \tag{3.8}$$

*with $a^{\mathrm{f}}(\cdot,\cdot)$ and $s_h(\cdot,\cdot)$ defined in* (1.18)*, and* (2.21)*, respectively.*

- *Lagrange multiplier projection: since, by the regularity assumption,* $\boldsymbol{\lambda} \in H^z(\Sigma)^d$*, we define the projection* $\boldsymbol{\Pi}_{\boldsymbol{\Lambda}}\boldsymbol{\lambda}$ *as the* $L^2$*-projection.*

$$\left(\boldsymbol{\Pi}_{\boldsymbol{\Lambda}}\boldsymbol{\lambda},\boldsymbol{\mu}_h\right)_{0,\Sigma} = \left(\boldsymbol{\lambda},\boldsymbol{\mu}_h\right)_{0,\Sigma} \qquad \forall \boldsymbol{\mu}_h \in \boldsymbol{\Lambda}_h. \tag{3.9}$$

**Remark 1.** *We point out that the solid velocity* $\dot{\boldsymbol{d}}$ *is projected in the finite dimensional space* $\boldsymbol{W}_h$ *using* $\boldsymbol{\Pi}_{\boldsymbol{W}}$*. The existence and uniqueness of solution to problem* (3.7) *is given in Proposition 1.5.1. Problem*(3.8) *has been analyzed in Appendix C where we report also the error estimates for the finite element approximations.*

In order to perform convergence analysis of the coupled problem, we need a preliminary analysis of the projection errors related to the operators just defined. In the following analysis we will assume additional regularity for the solution of the coupled problem (see heuristic considerations reported in paragraph 2.5);

**Proposition 3.3.3.** *Given* $(\boldsymbol{u},p,\boldsymbol{d},\boldsymbol{\lambda}) \in \boldsymbol{V}\times Q\times \boldsymbol{W}\times\boldsymbol{\Lambda}$ *such that*

$$\boldsymbol{u} \in H^{l+1}(\Omega)^d, \qquad p \in H^l(\Omega), \qquad \boldsymbol{\lambda} \in H^z(\Sigma)^d \qquad \boldsymbol{d} \in \boldsymbol{D}^{\mathrm{s}}, \tag{3.10}$$

*where* $\boldsymbol{D}^{\mathrm{s}} \subset H^{1+m}(\Sigma)^d$ *is the domain of the solid operator* $\boldsymbol{L}^s$ *and* $0 \le l < \frac{1}{2}$ , $0 \le m$, $0 \le z < \frac{1}{2}$.

1. *Let* $\boldsymbol{\Pi}_{\boldsymbol{W}}\boldsymbol{d} \in \boldsymbol{W}_h$ *be the solution of Problem* (3.7)*, then*

$$\|\boldsymbol{d} - \boldsymbol{\Pi}_{\boldsymbol{W}}\boldsymbol{d}\|_{1,\Sigma} \le Ch_s^m |\boldsymbol{d}|_{1+m,\Sigma}, \tag{3.11}$$

$$\|\boldsymbol{d} - \boldsymbol{\Pi}_{\boldsymbol{W}}\boldsymbol{d}\|_{0,\Sigma} \le Ch_s^{1+m} |\boldsymbol{d}|_{1+m,\Sigma}, \tag{3.12}$$

   *where* $C > 0$ *are a constants that do not depend on the mesh size* $h_s$.

2. *Let* $\left(\boldsymbol{\Pi}_{\boldsymbol{V}},\boldsymbol{\Pi}_Q\right) \in \boldsymbol{V}_h \times Q_h$ *be the unique solution of Problem* (3.8)*, there exists* $C > 0$*, independent of* $h_f$*, such that*

$$\|\boldsymbol{u} - \boldsymbol{\Pi}_{\boldsymbol{V}}(\boldsymbol{u},p)\|_{1,\Omega} + \|p - \boldsymbol{\Pi}_Q(\boldsymbol{u},p)\|_{0,\Omega} \le Ch_f^l \left(\|\boldsymbol{u}\|_{1+l,\Omega} + \|p\|_{l,\Omega}\right), \tag{3.13}$$

   *moreover, if* $\Omega$ *is a convex polygon (or a domain of class* $C^{1,1}$ *in the case* $d = 3$ *[48, Lemma 4.17])*

$$\|\boldsymbol{u} - \boldsymbol{\Pi}_{\boldsymbol{V}}(\boldsymbol{u},p)\|_{0,\Omega} \le Ch_f^{1+l} \left(\|\boldsymbol{u}\|_{1+l,\Omega} + \|p\|_{l,\Omega}\right). \tag{3.14}$$

3. *Let* $\mathbf{\Pi_\Lambda}\boldsymbol{\lambda} \in \mathbf{\Lambda}_h$ *be the solution of Problem* (3.9)*, then*

$$\|\boldsymbol{\lambda} - \mathbf{\Pi_\Lambda}\boldsymbol{\lambda}\|_\Lambda \leq C h_s^{\frac{1}{2}+z} \|\boldsymbol{\lambda}\|_{z,\Sigma}, \tag{3.15}$$

*where* $C > 0$ *is a constant that does not depend on the mesh size* $h_s$*.*

*Proof.* 1. Since $\boldsymbol{d}$ is in $H^{1+m}(\Sigma)^d$ and $\Sigma$ is of dimension $d-1$, we infer, applying Sobolev immersion theorem ( point 6 of Proposition A.3.2), that $\boldsymbol{d}$ is continuous, in fact, since $d = 2$ or $d = 3$, we always have $1 + m > \frac{d-1}{2}$. Since $\boldsymbol{d}$ is continuous we can consider its Lagrange interpolation $\mathscr{I}(\boldsymbol{d}) \in \boldsymbol{W}_h$. The Lagrange interpolate $\mathscr{I}(\boldsymbol{d})$ fulfils the following approximation estimate (see Theorem A.4.1 and Proposition A.44)

$$\|\boldsymbol{d} - \mathscr{I}(\boldsymbol{d})\|_{1,\Sigma} \leq C h_s^m |\boldsymbol{d}|_{1+m,\Sigma}, \tag{3.16}$$

whit $C > 0$ independent of $h_s$. Then, using the result in (1.25), we obtain

$$\|\boldsymbol{d} - \mathbf{\Pi_W}\boldsymbol{d}\|_{1,\Sigma} \leq C \inf_{\boldsymbol{d}_h \in \boldsymbol{W}_h} \|\boldsymbol{d} - \boldsymbol{d}_h\|_{1,\Sigma} \leq C\|\boldsymbol{d} - \mathscr{I}(\boldsymbol{d})\|_{1,\Sigma} \leq C h_s^m |\boldsymbol{d}|_{1+m,\Sigma}. \tag{3.17}$$

The claim 3.12 follows applying the Aubin-Nitsche Lemma to the solid variational problem.

2. (3.13) is a corollary to Theorem C.1.7; let $\mathscr{S}\mathscr{Z}_h : \boldsymbol{V} \to \boldsymbol{V}_h$ be the Scott-Zhang interpolation for the fluid velocity, owing to inequality (A.41), there exist a constant $C > 0$ uniform in $h_f$ such that

$$\|\boldsymbol{v} - \mathscr{S}\mathscr{Z}_h \boldsymbol{v}\|_{1,\Omega} \leq C h_f^l \|\boldsymbol{v}\|_{1+l,\Omega}, \qquad \forall \boldsymbol{v} \in H^{1+l}(\Omega)^d. \tag{3.18}$$

Here we used Scott-Zhang interpolation because it preserves zero Dirichlet boundary value of $\boldsymbol{v} \in \boldsymbol{V}$. Let $\mathscr{C}_h : \boldsymbol{Q} \to \boldsymbol{Q}_h$ be the Clément interpolation operators for the fluid pressure, owing to (A.38), there exist a constant $C > 0$ uniform in $h_f$ such that, for all $l \in [0, \frac{1}{2})$

$$\|q - \mathscr{C}_h q\|_{0,\Omega} \leq C h_f^l \|q\|_{l,\Omega}, \qquad \forall q \in H_0^l(\Omega). \tag{3.19}$$

Using (3.18) and (3.19) in (C.35), we have

$$\|\boldsymbol{u} - \mathbf{\Pi_V}(\boldsymbol{u}, p)\|_{1,\Omega} + \|p - \mathbf{\Pi_Q}(\boldsymbol{u}, p)\|_{0,\Omega}$$

$$\le C\left[\inf_{\boldsymbol{u}_h\in\boldsymbol{V}_h}\|\boldsymbol{u}-\boldsymbol{u}_h\|_{1,\Omega}+\inf_{p_h\in Q_h}\|p-p_h\|_{0,\Omega}\right]$$
$$\le C\left(\|\boldsymbol{u}-\mathscr{S}\mathscr{I}_h\boldsymbol{u}\|_{1,\Omega}+\|p-\mathscr{C}_hp\|_{0,\Omega}\right)$$
$$\le C\left(h_f^l\|\boldsymbol{u}\|_{1+l,\Omega}+h_f^l\|p\|_{l,\Omega}\right). \quad (3.20)$$

(3.14) comes applying Aubin-Nitsche Lemma.

$$\|\boldsymbol{u}-\boldsymbol{\Pi}_{\boldsymbol{V}}(\boldsymbol{u},p)\|_{0,\Omega}\le Ch_f^{1+l}\left(\|\boldsymbol{u}\|_{l+1,\Omega}+\|p\|_{l,\Omega}\right). \quad (3.21)$$

3. In order to prove the estimate for the Lagrange multiplier error we observe that, since $\boldsymbol{\lambda}\in H^z(\Sigma)^d\subset L^2(\Sigma)^d$, and $H^{\frac{1}{2}}(\Sigma)^d\hookrightarrow L^2(\Sigma)^d\hookrightarrow\boldsymbol{\Lambda}$ with dense inclusions, we can make the following identification

$$c(\boldsymbol{\lambda},\boldsymbol{w})=(\boldsymbol{\lambda},\boldsymbol{w})_{0,\Sigma}\qquad\forall\boldsymbol{w}\in H^{\frac{1}{2}}(\Sigma)^d.$$

From the definition of the projection for $\boldsymbol{\lambda}$ and a duality argument, denoting by $\mathscr{C}_h\boldsymbol{\mu}$ the Clément interpolate of $\boldsymbol{\mu}\in H^{\frac{1}{2}}(\Sigma)^d$, we obtain

$$\|\boldsymbol{\lambda}-\boldsymbol{\Pi}_{\boldsymbol{\Lambda}}\boldsymbol{\lambda}\|_{\boldsymbol{\Lambda}}=\sup_{\boldsymbol{w}\in H^{\frac{1}{2}}(\Sigma)^d}\frac{\langle\boldsymbol{\lambda}-\boldsymbol{\Pi}_{\boldsymbol{\Lambda}}\boldsymbol{\lambda},\boldsymbol{w}\rangle}{\|\boldsymbol{w}\|_{\frac{1}{2},\Sigma}}=\sup_{\boldsymbol{w}\in H^{\frac{1}{2}}(\Sigma)^d}\frac{(\boldsymbol{\lambda}-\boldsymbol{\Pi}_{\boldsymbol{\Lambda}}\boldsymbol{\lambda},\boldsymbol{w}-\mathscr{C}_h\boldsymbol{w})_{0,\Sigma}}{\|\boldsymbol{w}\|_{\frac{1}{2},\Sigma}} \quad (3.22)$$
$$\le\|\boldsymbol{\lambda}-\boldsymbol{\Pi}_{\boldsymbol{\Lambda}}\boldsymbol{\lambda}\|_{0,\Sigma}\frac{\|\boldsymbol{w}-\mathscr{C}_h\boldsymbol{w}\|_{0,\Sigma}}{\|\boldsymbol{w}\|_{\frac{1}{2},\Sigma}}\le C\|\boldsymbol{\lambda}-\boldsymbol{\Pi}_{\boldsymbol{\Lambda}}\boldsymbol{\lambda}\|_{0,\Sigma}h_s^{\frac{1}{2}} \quad (3.23)$$
$$\le Ch_s^{\frac{1}{2}+z}\|\boldsymbol{\lambda}\|_{z,\Sigma}, \quad (3.24)$$

where $C>0$ is a constant that does not depend on the mesh size $h_s$.

□

### 3.3.3 Discrete solid operator

In the analysis of the numerical integration schemes we will use the discrete counterpart of the elastic operator $\boldsymbol{L}^{\mathrm{s}}$ defined as follows: for all $\boldsymbol{w}\in\boldsymbol{W}$, $\boldsymbol{L}_h^{\mathrm{s}}\boldsymbol{w}\in\boldsymbol{W}_h$ is such that

$$\left(\boldsymbol{L}_h^{\mathrm{s}}\boldsymbol{w},\boldsymbol{w}_h\right)_{0,\Sigma}=a^{\mathrm{s}}(\boldsymbol{w},\boldsymbol{w}_h)\qquad\forall\boldsymbol{w}_h\in\boldsymbol{W}_h. \quad (3.25)$$

We summarize some important properties of this operator in the following Lemma (see [56, Lemma 1])

**Lemma 3.3.4.** • *Let $\boldsymbol{w}\in\boldsymbol{D}^{\mathrm{s}}\subset\boldsymbol{W}$ then,*

$$\|\boldsymbol{L}_h^{\mathrm{s}}\boldsymbol{w}\|_{0,\Sigma}\le\|\boldsymbol{L}^{\mathrm{s}}\boldsymbol{w}\|_{0,\Sigma} \quad (3.26)$$

• *Let $\boldsymbol{w}\in\boldsymbol{W}$, then*

$$\boldsymbol{L}_h^{\mathrm{s}}\left(\boldsymbol{\Pi}_{\boldsymbol{W}}(\boldsymbol{w})\right)=\boldsymbol{L}_h^{\mathrm{s}}(\boldsymbol{w}) \quad (3.27)$$

- *Let* $\boldsymbol{w}_h \in \boldsymbol{W}_h$

$$\|\boldsymbol{L}_h^{\mathrm{s}} \boldsymbol{w}_h\|_{\mathrm{s}} \leq \frac{\beta_{\mathrm{s}} C_I^2}{h_s^2} \|\boldsymbol{w}_h\|_{\mathrm{s}} \tag{3.28}$$

$$\|\boldsymbol{L}_h^{\mathrm{s}} \boldsymbol{w}_h\|_{0,\Sigma} \leq \frac{\beta_{\mathrm{s}}^{1/2} C_I}{h_s} \|\boldsymbol{w}_h\|_{\mathrm{s}} \tag{3.29}$$

$$\|\boldsymbol{w}_h\|_{\mathrm{s}}^2 \leq \frac{\beta_{\mathrm{s}} C_I^2}{h_s^2} \|\boldsymbol{w}_h\|_{0,\Sigma}^2 \tag{3.30}$$

*where* $C_I > 0$ *is the constant of the inverse estimate* (A.43) *(we recall that we used the hypothesis of uniformity for the solid mesh family as specified in Section 2.3 and Section 2.4).*

### 3.3.4 Error definition

In order to perform the convergence analysis we will study the behavior of the error $\boldsymbol{e} = \boldsymbol{a} - \boldsymbol{a}_h$, between the generic quantity $\boldsymbol{a}$ satisfying a the weak form of a differential problem and the associated finite element solution $\boldsymbol{a}_h$, when the discretization parameters go to zero. We split the error as follows

$$\boldsymbol{e} = \boldsymbol{a} - \boldsymbol{a}_h = \underbrace{\boldsymbol{a} - \boldsymbol{a}_\Pi}_{\boldsymbol{e}_\Pi} + \underbrace{\boldsymbol{a}_\Pi - \boldsymbol{a}_h}_{\boldsymbol{e}_h},$$

where $\boldsymbol{a}_\Pi$ is a Ritz projection of $\boldsymbol{a}$ in the discrete space. the term $\boldsymbol{e}_\Pi$ is called projection error and $\boldsymbol{e}_h$ is the discrete error. In particular, the projection error is controlled using the classical interpolations properties for Sobolev spaces given in Chapter A, while the behavior of the discrete error is analyzed in the present section. In this spirit we introduce the following notation, given the operators $\boldsymbol{\Pi}_{\boldsymbol{V}}$, $\boldsymbol{\Pi}_Q$, $\boldsymbol{\Pi}_{\boldsymbol{W}}$ and $\boldsymbol{\Pi}_{\boldsymbol{\Lambda}}$ of Definition 1

**Definition 2.** *For* $n \geq 0$ *let* $(\boldsymbol{u}^n, p^n, \boldsymbol{d}^n, \dot{\boldsymbol{d}}^n, \boldsymbol{\lambda}^n) \in \boldsymbol{V} \times Q \times \boldsymbol{W} \times \boldsymbol{W} \times \boldsymbol{\Lambda}$ *and* $(\boldsymbol{u}_h^n, p_h^n, \boldsymbol{d}_h^n, \dot{\boldsymbol{d}}_h^n, \boldsymbol{\lambda}_h^n) \in \boldsymbol{V}_h \times Q_h \times \boldsymbol{W}_h \times \boldsymbol{W}_h \times \boldsymbol{\Lambda}_h$*, then we give the following notations*

$$\boldsymbol{u}_\Pi^n \stackrel{def}{=} \boldsymbol{\Pi}_{\boldsymbol{V}}(\boldsymbol{u}^n, p^n), \qquad \boldsymbol{u}^n - \boldsymbol{u}_h^n = \underbrace{\boldsymbol{u}^n - \boldsymbol{u}_\Pi^n}_{\boldsymbol{\theta}_\pi^n} + \underbrace{\boldsymbol{u}_\Pi^n - \boldsymbol{u}_h^n}_{\boldsymbol{\theta}_h^n},$$

$$p_\Pi^n \stackrel{def}{=} \boldsymbol{\Pi}_Q(\boldsymbol{u}^n, p^n), \qquad p^n - p_h^n = \underbrace{p^n - p_\Pi^n}_{\varphi_\pi^n} + \underbrace{p_\Pi^n - p_h^n}_{\varphi_h^n},$$

$$\boldsymbol{d}_\Pi^n \stackrel{def}{=} \boldsymbol{\Pi}_{\boldsymbol{W}} \boldsymbol{d}^n, \qquad \boldsymbol{d}^n - \boldsymbol{d}_h^n = \underbrace{\boldsymbol{d}^n - \boldsymbol{d}_\Pi^n}_{\boldsymbol{\xi}_\pi^n} + \underbrace{\boldsymbol{d}_\Pi^n - \boldsymbol{d}_h^n}_{\boldsymbol{\xi}_h^n},$$

$$\dot{\boldsymbol{d}}_\Pi^n \stackrel{def}{=} \boldsymbol{\Pi}_{\boldsymbol{W}} \dot{\boldsymbol{d}}^n, \qquad \dot{\boldsymbol{d}}^n - \dot{\boldsymbol{d}}_h^n = \underbrace{\dot{\boldsymbol{d}}^n - \dot{\boldsymbol{d}}_\Pi^n}_{\dot{\boldsymbol{\xi}}_\pi^n} + \underbrace{\dot{\boldsymbol{d}}_\Pi^n - \dot{\boldsymbol{d}}_h^n}_{\dot{\boldsymbol{\xi}}_h^n},$$

$$\boldsymbol{\lambda}_\Pi^n \stackrel{def}{=} \boldsymbol{\Pi}_\Lambda \boldsymbol{\lambda}^n, \qquad \boldsymbol{\lambda}^n - \boldsymbol{\lambda}_h^n = \underbrace{\boldsymbol{\lambda}^n - \boldsymbol{\lambda}_\Pi^n}_{\boldsymbol{\omega}_\pi^n} + \underbrace{\boldsymbol{\lambda}_\Pi^n - \boldsymbol{\lambda}_h^n}_{\boldsymbol{\omega}_h^n}$$

Before introducing the convergence theorem, it is useful to establish the following expression for the error $\dot{\boldsymbol{\xi}}_h^n$ introduced in Definition 2

**Proposition 3.3.5.** *It holds*

$$\dot{\boldsymbol{\xi}}_h^n = \dot{\boldsymbol{d}}_\Pi^n - \partial_\tau \boldsymbol{d}_\Pi^n + \partial_\tau \boldsymbol{\xi}_h^n. \tag{3.31}$$

*Proof.* The claim establish a relation between the error in the solid velocities $\dot{\boldsymbol{\xi}}_h^n = \dot{\boldsymbol{d}}_\Pi^n - \dot{\boldsymbol{d}}_h^n$ and the error in the solid displacements $\boldsymbol{\xi}_h^n = \boldsymbol{d}_\Pi^n - \boldsymbol{d}_h^n$

$$\dot{\boldsymbol{\xi}}_h^n = \dot{\boldsymbol{d}}_\Pi^n - \dot{\boldsymbol{d}}_h^n = \dot{\boldsymbol{d}}_\Pi^n - \partial_\tau \boldsymbol{d}_h^n = \dot{\boldsymbol{d}}_\Pi^n - \partial_\tau \boldsymbol{d}_\Pi^n + \partial_\tau \boldsymbol{d}_\Pi^n - \partial_\tau \boldsymbol{d}_h^n = \dot{\boldsymbol{d}}_\Pi^n - \partial_\tau \boldsymbol{d}_\Pi^n + \partial_\tau \boldsymbol{\xi}_h^n$$

□

In the proof of convergence of Algorithm 1 we will use also the following estimates based on regularity assumptions already discussed in section 2.5.

**Proposition 3.3.6.** *Let* $\big(\boldsymbol{u}(t), p(t), \boldsymbol{d}(t), \dot{\boldsymbol{d}}(t), \boldsymbol{\lambda}(t)\big) \in \boldsymbol{V} \times \boldsymbol{Q} \times \boldsymbol{W} \times \boldsymbol{W} \times \boldsymbol{\Lambda}$ *be such that for a given final time* $T > \tau$*:*

$$\begin{aligned} &\boldsymbol{u} \in H^1(0,T;H^{l+1}(\Omega)^d), && p \in H^1(0,T;H^l(\Omega)),\\ &\boldsymbol{\lambda} \in L^\infty(0,T;H^z(\Sigma)^d), && \boldsymbol{d},\dot{\boldsymbol{d}} \in H^1(0,T;\boldsymbol{D}^{\rm s}), \end{aligned} \tag{3.32}$$

*where* $\boldsymbol{D}^{\rm s} \subset H^{1+m}(\Sigma)^d$ *and* $0 \le l < \frac{1}{2}$ *,* $0 \le m,$ $0 \le z < \frac{1}{2}$*. Then, with the notation of Definition 2, we have*

$$\begin{aligned} &\|\partial_t \boldsymbol{\theta}_\pi\|_{L^2(0,T;L^2(\Omega)^d)} \approx \mathcal{O}(h_f^{1+l}), && \|\boldsymbol{\theta}_\pi\|_{L^\infty(0,T;H^1(\Omega)^d)} \approx \mathcal{O}(h_f^l),\\ &\|\dot{\boldsymbol{\xi}}_\pi\|_{L^\infty(0,T;H^{\frac{1}{2}}(\Sigma)^d)} \approx \mathcal{O}(h_s^{\frac{1}{2}+m}), && \|\partial_t \dot{\boldsymbol{\xi}}_\pi\|_{L^2(0,T;L^2(\Sigma)^d)} \approx \mathcal{O}(h_s^{1+m}),\\ &\|\boldsymbol{\omega}_\pi\|_{L^\infty(0,T;\boldsymbol{\Lambda})} \approx \mathcal{O}(h_s^{z+\frac{1}{2}}). \end{aligned} \tag{3.33}$$

*Proof.* Given the assumptions on the regularity of the functions 3.32, the claims (3.33) are consequence of Proposition 3.3.3.

$$\|\partial_t \boldsymbol{\theta}_\pi\|_{L^2(0,T;L^2(\Omega)^d)} \le C h_f^{1+l} \Big( \|\partial_t \boldsymbol{u}\|_{L^2(0,T;H^{1+l}(\Omega)^d)} + \|\partial_t p\|_{L^2(0,T;H^l(\Omega))} \Big) \approx \mathcal{O}(h_f^{1+l}). \tag{3.34}$$

$$\|\boldsymbol{\theta}_\pi\|_{L^\infty(0,T;H^1(\Omega)^d)} \le C h_f^l \Big( \|\partial_t \boldsymbol{u}\|_{L^\infty(0,T;H^{1+l}(\Omega)^d)} + \|\partial_t p\|_{L^\infty(0,T;H^l(\Omega))} \Big) \approx \mathcal{O}(h_f^l). \tag{3.35}$$

$$\|\dot{\boldsymbol{\xi}}_\pi\|_{L^\infty(0,T;H^{\frac{1}{2}}(\Sigma)^d)} \le C h_s^{\frac{1}{2}+m} \|\dot{\boldsymbol{d}}\|_{L^\infty(0,T;H^{1+m}(\Omega))} \approx \mathcal{O}(h_s^{\frac{1}{2}+m}) \tag{3.36}$$

$$\|\partial_t \dot{\boldsymbol{\xi}}_\pi\|_{L^2(0,T;L^2(\Sigma)^d)} \leq Ch_s^{1+m}\|\partial_t \dot{\boldsymbol{d}}\|_{L^2(0,T;H^{1+m}(\Omega))} \approx \mathcal{O}(h_s^{1+m}) \quad (3.37)$$

$$\|\boldsymbol{\omega}_\pi\|_{L^\infty(0,T;\boldsymbol{\Lambda})} \leq Ch_s^{\frac{1}{2}+z}\|\boldsymbol{\lambda}\|_{L^\infty(0,T;H^z(\Sigma)^d)} \approx \mathcal{O}(h_s^{\frac{1}{2}+z}) \quad (3.38)$$

□

### 3.3.5 Convergence estimates

In this paragraph we give the main result of this section concerning the convergence of Algorithm 1. We assume the following regularity for the solution of the continuous problem

$$\begin{aligned}
&\boldsymbol{u} \in H^1(0,T;H^{l+1}(\Omega)^d), && \boldsymbol{u}|_\Sigma \in H^1(0,T;H^{l+\frac{1}{2}}(\Sigma)^d),\\
&\partial_{tt}\boldsymbol{u} \in L^2(0,T;L^2(\Omega)^d), && \partial_{tt}\boldsymbol{u}|_\Sigma \in L^2(0,T;L^2(\Sigma)^d),\\
&p \in H^1(0,T;H^l(\Omega)), && \boldsymbol{\lambda} \in L^\infty(0,T;H^z(\Sigma)^d)\\
&\boldsymbol{d} \in H^2(0,T;\boldsymbol{D}^{\mathrm{s}}), && \dot{\boldsymbol{d}} \in H^2(0,T;\boldsymbol{D}^{\mathrm{s}}),
\end{aligned} \quad (3.39)$$

then the error estimate follows.

**Theorem 3.3.7.** *Let $(\boldsymbol{u},p,\dot{\boldsymbol{d}},\boldsymbol{d},\boldsymbol{\lambda})$ satisfies Problem 9 with regularity assumptions* (3.39) *and let $\{(\boldsymbol{u}_h^n, p_h^n, \dot{\boldsymbol{d}}_h^n, \boldsymbol{d}_h^n, \boldsymbol{\lambda}_h^n)\}_{n>0}$ be given by the Algorithm 1, then, if $\tau < 1$, we have the following error estimate. For $n > 0$ and $n\tau < T$*

$$\begin{aligned}
\mathcal{E}_h^n &\overset{def}{=} \frac{\rho_{\mathrm{f}}}{2}\|\boldsymbol{\theta}_h^n\|_{0,\Omega}^2 + \frac{\rho_s\epsilon}{2}\|\dot{\boldsymbol{\xi}}_h^n\|_{0,\Sigma}^2 + \frac{1}{2}\|\boldsymbol{\xi}_h^n\|_{\mathrm{s}}^2\\
&\lesssim \exp\left(\sum_{m=1}^n \frac{\tau}{1-\tau}\right)\left[\underbrace{\mathcal{O}(h_f^{2l}) + \mathcal{O}(h_s^{2m+1})}_{A_0} + \tau\left(\underbrace{\mathcal{O}(h_s^{2z+1}) + \mathcal{O}(h_s^{2m+1}) + \mathcal{O}(h_f^{2l})}_{A_1}\right) + \tau^2 A_2\right].
\end{aligned} \quad (3.40)$$

*Proof.* We begin deriving a discrete error equation considering the linearized Problem 9 given at the beginning of this chapter; for $t = t^n$, taking test functions in the discrete spaces and using the notation $\boldsymbol{v}|_\Sigma \overset{def}{=} \boldsymbol{v} \circ \boldsymbol{\phi}_0$ for all $\boldsymbol{v} \in \boldsymbol{V}$, we have

$$\begin{aligned}
\rho_{\mathrm{f}}\big(\partial_t \boldsymbol{u}^n, \boldsymbol{v}_h\big)_{0,\Omega} &+ a^{\mathrm{f}}\big((\boldsymbol{u}^n,p^n),(\boldsymbol{v}_h,q_h)\big) + c\left(\boldsymbol{\lambda}^n, \boldsymbol{v}_h|_\Sigma - \boldsymbol{w}_h\right)\\
&- c\left(\boldsymbol{\mu}_h, \boldsymbol{u}^n|_\Sigma - \dot{\boldsymbol{d}}^n\right) + \rho_s\epsilon\left(\partial_t \dot{\boldsymbol{d}}^n, \boldsymbol{w}_h\right)_{0,\Sigma} + a^{\mathrm{s}}\left(\boldsymbol{d}^n, \boldsymbol{w}_h\right) = 0
\end{aligned} \quad (3.41)$$

for all $(\boldsymbol{v}_h, q_h, \boldsymbol{w}_h, \boldsymbol{\mu}_h) \in \boldsymbol{V}_h \times Q_h \times \boldsymbol{W}_h \times \boldsymbol{\Lambda}_h$.

We can derive an error equation by subtracting from (3.41) the equation (3.2) of Algorithm 1 evaluated on the same test functions; the result is reported in the following, where we add and subtract the terms $\rho_{\mathrm{f}}(\partial_\tau \boldsymbol{u}^n, \boldsymbol{v}_h)_{0,\Omega}$ and $\rho_s\epsilon\left(\partial_\tau \dot{\boldsymbol{d}}^n, \boldsymbol{w}_h\right)_{0,\Sigma}$

$$
\begin{aligned}
&\rho_{\mathrm{f}}\big(\partial_\tau(\boldsymbol{u}^n - \boldsymbol{u}_h^n), \boldsymbol{v}_h\big)_{0,\Omega} + a^{\mathrm{f}}\big((\boldsymbol{u}^n - \boldsymbol{u}_h^n, p^n - p_h^n), (\boldsymbol{v}_h, q_h)\big) + c\left(\boldsymbol{\lambda}^n - \boldsymbol{\lambda}_h^n, \boldsymbol{v}_h|_\Sigma - \boldsymbol{w}_h\right) \\
&\quad - c\left(\boldsymbol{\mu}_h, \boldsymbol{u}^n|_\Sigma - \boldsymbol{u}_h^n|_\Sigma\right) + c\left(\boldsymbol{\mu}_h, \dot{\boldsymbol{d}}^n - \dot{\boldsymbol{d}}_h^n\right) + \rho_s\epsilon\left(\partial_\tau(\dot{\boldsymbol{d}}^n - \dot{\boldsymbol{d}}_h^n), \boldsymbol{w}_h\right)_{0,\Sigma} + a^{\mathrm{s}}\left(\boldsymbol{d}^n - \boldsymbol{d}_h^n, \boldsymbol{w}_h\right) \\
&\qquad = s_h(p_h^n, q_h) + \rho_{\mathrm{f}}\big((\partial_\tau - \partial_t)\boldsymbol{u}^n, \boldsymbol{v}_h\big)_{0,\Omega} + \rho_s\epsilon\left((\partial_\tau - \partial_t)\dot{\boldsymbol{d}}^n, \boldsymbol{w}_h\right)_{0,\Sigma} \qquad (3.42)
\end{aligned}
$$

for all $(\boldsymbol{v}_h, q_h, \boldsymbol{w}_h, \boldsymbol{\mu}_h) \in \boldsymbol{V}_h \times Q_h \times \boldsymbol{W}_h \times \boldsymbol{\Lambda}_h$.

Introducing the error notations given in Definition 2 we have

$$
\begin{aligned}
&\rho_{\mathrm{f}}\big(\partial_\tau\boldsymbol{\theta}_h^n, \boldsymbol{v}_h\big)_{0,\Omega} + a^{\mathrm{f}}\big((\boldsymbol{\theta}_h^n, \varphi_h^n), (\boldsymbol{v}_h, q_h)\big) + \rho_s\epsilon\left(\partial_\tau\dot{\boldsymbol{\xi}}_h^n, \boldsymbol{w}_h\right)_{0,\Sigma} + a^{\mathrm{s}}\left(\boldsymbol{\xi}_h^n, \boldsymbol{w}_h\right) = \\
&\quad \rho_{\mathrm{f}}\big((\partial_\tau - \partial_t)\boldsymbol{u}^n, \boldsymbol{v}_h\big)_{0,\Omega} - \rho_{\mathrm{f}}\big(\partial_\tau\boldsymbol{\theta}_\pi^n, \boldsymbol{v}_h\big)_{0,\Omega} + \rho_s\epsilon\left((\partial_\tau - \partial_t)\dot{\boldsymbol{d}}^n, \boldsymbol{w}_h\right)_{0,\Sigma} - \rho_s\epsilon\left(\partial_\tau\dot{\boldsymbol{\xi}}_\pi^n, \boldsymbol{w}_h\right)_{0,\Omega} \\
&\qquad - c\left(\boldsymbol{\omega}_h^n, \boldsymbol{v}_h|_\Sigma - \boldsymbol{w}_h\right) + c\big(\boldsymbol{\mu}_h, \boldsymbol{\theta}_h^n|_\Sigma - \dot{\boldsymbol{\xi}}_h^n\big) + s_h(p_h^n, q_h) + c\left(\boldsymbol{\mu}_h, \boldsymbol{\theta}_\pi^n|_\Sigma - \dot{\boldsymbol{\xi}}_\pi^n\right) \\
&\qquad\quad - \left[a^{\mathrm{f}}\big((\boldsymbol{\theta}_\pi^n, \varphi_\pi^n), (\boldsymbol{v}_h, q_h)\big) + a^{\mathrm{s}}\left(\boldsymbol{\xi}_\pi^n, \boldsymbol{w}_h\right) + c\left(\boldsymbol{\omega}_\pi^n, \boldsymbol{v}_h|_\Sigma - \boldsymbol{w}_h\right)\right], \qquad (3.43)
\end{aligned}
$$

for all $(\boldsymbol{v}_h, q_h, \boldsymbol{w}_h, \boldsymbol{\mu}_h) \in \boldsymbol{V}_h \times Q_h \times \boldsymbol{W}_h \times \boldsymbol{\Lambda}_h$.

From the definition of the projection operators given in Definition 1 we have

$$
\begin{aligned}
&\text{from (3.7)} \quad a^{\mathrm{s}}\left(\boldsymbol{\xi}_\pi^n, \boldsymbol{w}_h\right) = 0 && \forall \boldsymbol{w}_h \in \boldsymbol{W}_h, \\
&\text{from (3.8)} \quad a^{\mathrm{f}}\big((\boldsymbol{\theta}_\pi^n, \varphi_\pi^n), (\boldsymbol{v}_h, q_h)\big) = s_h(p_\Pi^n, q_h) && \forall (\boldsymbol{v}_h, q_h) \in \boldsymbol{V}_h \times Q_h, \qquad (3.44) \\
&\text{from (3.9)} \quad c\left(\boldsymbol{\omega}_\pi^n, \boldsymbol{v}_h|_\Sigma - \boldsymbol{w}_h\right) = c\left(\boldsymbol{\omega}_\pi^n, \boldsymbol{v}_h|_\Sigma\right) && \forall (\boldsymbol{v}_h, \boldsymbol{w}_h) \in \boldsymbol{V}_h \times \boldsymbol{W}_h.
\end{aligned}
$$

Hence, using (3.44) and testing with $(\boldsymbol{v}_h, q_h, \boldsymbol{\mu}_h, \boldsymbol{w}_h) = \tau(\boldsymbol{\theta}_h^n, \varphi_h^n, \boldsymbol{\omega}_h^n, \dot{\boldsymbol{\xi}}_h^n)$, we obtain the following error relation

$$
\begin{aligned}
&\rho_{\mathrm{f}}\tau\left(\partial_\tau\boldsymbol{\theta}_h^n, \boldsymbol{\theta}_h^n\right)_{0,\Omega} + \tau a^{\mathrm{f}}\big((\boldsymbol{\theta}_h^n, \varphi_h^n), (\boldsymbol{\theta}_h^n, \varphi_h^n)\big) + \rho_s\epsilon\tau\left(\partial_\tau\dot{\boldsymbol{\xi}}_h^n, \dot{\boldsymbol{\xi}}_h^n\right)_{0,\Sigma} \\
&\quad + \tau a^{\mathrm{s}}\left(\boldsymbol{\xi}_h^n, \dot{\boldsymbol{\xi}}_h^n\right) + \tau s_h(\varphi_h^n, \varphi_h^n) = \rho_{\mathrm{f}}\tau\left((\partial_\tau - \partial_t)\boldsymbol{u}^n, \boldsymbol{\theta}_h^n\right)_{0,\Omega} - \rho_{\mathrm{f}}\tau\left(\partial_\tau\boldsymbol{\theta}_\pi^n, \boldsymbol{\theta}_h^n\right)_{0,\Omega} \\
&\quad + \rho_s\epsilon\tau\left((\partial_\tau - \partial_t)\dot{\boldsymbol{d}}^n, \dot{\boldsymbol{\xi}}_h^n\right)_{0,\Sigma} - \rho_s\epsilon\tau\left(\partial_\tau\dot{\boldsymbol{\xi}}_\pi^n, \dot{\boldsymbol{\xi}}_h^n\right)_{0,\Sigma} - \tau c\left(\boldsymbol{\omega}_\pi^n, \boldsymbol{\theta}_h^n|_\Sigma\right) + \tau c\left(\boldsymbol{\omega}_h^n, \boldsymbol{\theta}_\pi^n|_\Sigma - \dot{\boldsymbol{\xi}}_\pi^n\right) \qquad (3.45)
\end{aligned}
$$

Using the result of Proposition 3.3.5 about the representation of the error $\dot{\boldsymbol{\xi}}_h^n$, we can write

$$
\begin{aligned}
&\rho_{\mathrm{f}}\tau\left(\partial_\tau\boldsymbol{\theta}_h^n, \boldsymbol{\theta}_h^n\right)_{0,\Omega} + \tau a^{\mathrm{f}}\big((\boldsymbol{\theta}_h^n, \varphi_h^n), (\boldsymbol{\theta}_h^n, \varphi_h^n)\big) + \rho_s\epsilon\tau\left(\partial_\tau\dot{\boldsymbol{\xi}}_h^n, \dot{\boldsymbol{\xi}}_h^n\right)_{0,\Sigma} + \tau a^{\mathrm{s}}\left(\boldsymbol{\xi}_h^n, \partial_\tau\boldsymbol{\xi}_h^n\right) + \tau s_h(\varphi_h^n, \varphi_h^n) \\
&\quad = \underbrace{-\tau a^{\mathrm{s}}\left(\boldsymbol{\xi}_h^n, \dot{\boldsymbol{d}}_\Pi^n - \boldsymbol{\Pi}_{\boldsymbol{W}}\partial_\tau\boldsymbol{d}^n\right)}_{T_0} + \underbrace{\rho_{\mathrm{f}}\tau\left((\partial_\tau - \partial_t)\boldsymbol{u}^n, \boldsymbol{\theta}_h^n\right)_{0,\Omega} - \rho_{\mathrm{f}}\tau\left(\partial_\tau\boldsymbol{\theta}_\pi^n, \boldsymbol{\theta}_h^n\right)_{0,\Omega}}_{T_1}
\end{aligned}
$$

$$+\underbrace{\rho_s\epsilon\tau\left((\partial_\tau-\partial_t)\dot{\boldsymbol{d}}^n,\dot{\boldsymbol{\xi}}_h^n\right)_{0,\Sigma}-\rho_s\epsilon\tau\left(\partial_\tau\dot{\boldsymbol{\xi}}_\pi^n,\dot{\boldsymbol{\xi}}_h^n\right)_{0,\Sigma}}_{T_2}\underbrace{-\tau c\left(\boldsymbol{\omega}_\pi^n,\boldsymbol{\theta}_h^n|_\Sigma\right)}_{T_3}+\underbrace{\tau c\left(\boldsymbol{\omega}_h^n,\boldsymbol{\theta}_\pi^n|_\Sigma-\dot{\boldsymbol{\xi}}_\pi^n\right)}_{T_4} \tag{3.46}$$

The term $T_0$ can be bounded as follows, using the definition of the projection operator (3.7)

$$\begin{aligned}
T_0 =& -\tau a^{\mathrm{s}}\left(\boldsymbol{\xi}_h^n,\dot{\boldsymbol{d}}_\Pi^n-\partial_\tau\boldsymbol{d}_\Pi^n\right)\\
=&\tau a^{\mathrm{s}}\big(\boldsymbol{\xi}_h^n,\dot{\boldsymbol{d}}^n-\dot{\boldsymbol{d}}_\Pi^n\big)+\tau a^{\mathrm{s}}\big(\boldsymbol{\xi}_h^n,(\partial_\tau-\partial_t)\boldsymbol{d}^n\big)+\tau a^{\mathrm{s}}\big(\boldsymbol{\xi}_h^n,\partial_\tau(\boldsymbol{d}_\Pi^n-\boldsymbol{d}^n)\big)\\
=&\tau a^{\mathrm{s}}\big(\boldsymbol{\xi}_h^n,(\partial_\tau-\partial_t)\boldsymbol{d}^n\big)\le\|\boldsymbol{\xi}_h^n\|_{\mathrm{s}}\|\tau\partial_\tau\boldsymbol{d}^n-\tau\dot{\boldsymbol{d}}^n\|_{\mathrm{s}}\le\tau^{\frac{3}{2}}\|\boldsymbol{\xi}_h^n\|_{\mathrm{s}}\|\partial_{tt}\boldsymbol{d}\|_{L^2(t_{n-1},t_n;\boldsymbol{W})}\\
\le&\frac{\delta_0\tau}{2}\|\boldsymbol{\xi}_h^n\|_{\mathrm{s}}^2+\frac{\tau^2}{2\delta_0}\|\partial_{tt}\boldsymbol{d}\|^2_{L^2(t_{n-1},t_n;\boldsymbol{W})},
\end{aligned} \tag{3.47}$$

where we have used the following property of the projection $\boldsymbol{\Pi}_{\boldsymbol{W}}$, for all $\boldsymbol{w}_h\in\boldsymbol{W}_h$,

$$a^{\mathrm{s}}\left(\dot{\boldsymbol{d}}^n-\dot{\boldsymbol{d}}_\Pi^n,\boldsymbol{w}_h\right)=0,\qquad a^{\mathrm{s}}\left(\partial_\tau(\boldsymbol{d}^n-\boldsymbol{d}_\Pi^n),\boldsymbol{w}_h\right)=0$$

and the notation $\|\partial_{tt}\boldsymbol{d}\|^2_{L^2(t_{n-1},t_n;\boldsymbol{W})}=\int_{t_{n-1}}^{t_n}\|\partial_{tt}\boldsymbol{d}(\tau)\|_s^2d\tau$.

For $T_1$ we obtain the following bounding relation; let $\delta_1>0$, using inequalities (3.3), (3.4)

$$\begin{aligned}
T_1\le&\tau\rho_{\mathrm{f}}\big(\|\partial_\tau\boldsymbol{u}^n-\partial_t\boldsymbol{u}^n\|_{0,\Omega}+\|\partial_\tau\boldsymbol{\theta}_\pi^n\|_{0,\Omega}\big)\|\boldsymbol{\theta}_h^n\|_{0,\Omega}\\
\le&\rho_{\mathrm{f}}\tau\left(\tau^{\frac{1}{2}}\|\partial_{tt}\boldsymbol{u}\|_{L^2(t_{n-1},t_n;L^2(\Omega)^d)}+\tau^{-\frac{1}{2}}\|\partial_t\boldsymbol{\theta}_\pi\|_{L^2(t_{n-1},t_n;L^2(\Omega)^d)}\right)\|\boldsymbol{\theta}_h^n\|_{0,\Omega}\\
\le&\frac{\rho_{\mathrm{f}}\tau^2}{2\delta_1}\|\partial_{tt}\boldsymbol{u}\|^2_{L^2(t_{n-1},t_n;L^2(\Omega)^d)}+\frac{\rho_{\mathrm{f}}}{2\delta_1}\|\partial_t\boldsymbol{\theta}_\pi\|^2_{L^2(t_{n-1},t_n;L^2(\Omega)^d)}+\rho_{\mathrm{f}}\tau\delta_1\|\boldsymbol{\theta}_h^n\|^2_{0,\Omega}
\end{aligned} \tag{3.48}$$

We obtain a similar relation for $T_2$; let $\delta_2>0$, using inequalities (3.3), (3.4)

$$\begin{aligned}
T_2\le&\tau\rho_s\epsilon\left(\|\partial_\tau\dot{\boldsymbol{d}}_h^n-\partial_t\dot{\boldsymbol{d}}_h^n\|_{0,\Sigma}+\|\partial_\tau\dot{\boldsymbol{\xi}}_\pi^n\|_{0,\Sigma}\right)\|\dot{\boldsymbol{\xi}}_h^n\|_{0,\Sigma}\\
\le&\rho_s\epsilon\tau\left(\tau^{\frac{1}{2}}\|\partial_{tt}\dot{\boldsymbol{d}}\|_{L^2(t_{n-1},t_n;L^2(\Sigma)^d)}+\tau^{-\frac{1}{2}}\|\partial_t\dot{\boldsymbol{\xi}}_\pi\|_{L^2(t_{n-1},t_n;L^2(\Sigma)^d)}\right)\|\dot{\boldsymbol{\xi}}_h^n\|_{0,\Sigma}\\
\le&\frac{\rho_s\epsilon\tau^2}{2\delta_2}\|\partial_{tt}\dot{\boldsymbol{d}}\|^2_{L^2(t_{n-1},t_n;L^2(\Sigma)^d)}+\frac{\rho_s\epsilon}{2\delta_2}\|\partial_t\dot{\boldsymbol{\xi}}_\pi\|^2_{L^2(t_{n-1},t_n;L^2(\Sigma)^d)}+\rho_s\epsilon\tau\delta_2\|\dot{\boldsymbol{\xi}}_h^n\|^2_{0,\Sigma}.
\end{aligned} \tag{3.49}$$

The term $T_3$ is estimated using the trace inequality of theorem A.3.3

$$T_3\;=\;\tau c\left(\boldsymbol{\omega}_\pi^n,\boldsymbol{\theta}_h^n|_\Sigma\right)\;\le\;\tau\|\boldsymbol{\omega}_\pi^n\|_{\boldsymbol{\Lambda}}\|\boldsymbol{\theta}_h^n|_\Sigma\|_{\frac{1}{2},\Sigma}\;\le\;C\tau\|\boldsymbol{\omega}_\pi^n\|_{\boldsymbol{\Lambda}}\|\boldsymbol{\theta}_h^n\|_{1,\Omega}\;\le\;\frac{\tau}{\delta_3}\|\boldsymbol{\omega}_\pi^n\|^2_{\boldsymbol{\Lambda}}+\tau\delta_3\|\boldsymbol{\theta}_h^n\|^2_{1,\Omega}.$$

Using the first Korn's inequality reported in Theorem A.3.4 we have

$$\|\boldsymbol{\theta}_h^m\|_{1,\Omega}^2 \le C_K \|\boldsymbol{\epsilon}(\boldsymbol{\theta}_h^m)\|_{0,\Omega}^2.$$

Then, since $\boldsymbol{\omega}_\pi \in L^\infty(0,T;\boldsymbol{\Lambda})$ we have

$$T_3 \le \frac{\tau}{\delta_3} \|\boldsymbol{\omega}_\pi\|_{L^\infty(0,T;\boldsymbol{\Lambda})}^2 + \tau\delta_3 C_K \|\boldsymbol{\epsilon}(\boldsymbol{\theta}_h^n)\|_{0,\Omega}^2. \tag{3.50}$$

The last term to be bounded is $T_4$.

$$T_4 = \tau c\left(\boldsymbol{\omega}_h^n, \boldsymbol{\theta}_\pi^n|_\Sigma - \dot{\boldsymbol{\xi}}_\pi^n\right) \le \tau \|\boldsymbol{\omega}_h^n\|_{\boldsymbol{\Lambda}} \|\boldsymbol{\theta}_\pi^n|_\Sigma - \dot{\boldsymbol{\xi}}_\pi^n\|_{\frac{1}{2},\Sigma}.$$

The estimate is obtained working as in the proof of the well-posedness of the time-step problem (see Proposition 2.4.3). Exploiting the continuous inf-sup property proved in [19], for all $\boldsymbol{\mu}_h \in \boldsymbol{\Lambda}_h$ there exists $\boldsymbol{v}_2 \in \boldsymbol{V}_0 \overset{def}{=} \{\boldsymbol{v} \in \boldsymbol{V} : \quad (\operatorname{div}\boldsymbol{v}, q)_{0,\Omega} = 0, \quad \forall q \in Q\}$, such that

$$c(\boldsymbol{\mu}_h, \boldsymbol{v}_2|_\Sigma) \ge \beta_1 \|\boldsymbol{v}_2\|_{1,\Omega} \|\boldsymbol{\mu}_h\|_{\boldsymbol{\Lambda}},$$

with the couple $(\boldsymbol{v}_2, p_2) \in \boldsymbol{V} \times Q$ solving the Stokes problem (2.41). Let $(\boldsymbol{v}_{2h}, p_{2h}) \in \boldsymbol{V}_h \times Q_h$ be the approximation of $(\boldsymbol{v}_2, p_2)$ given in Proposition 2.4.3, then we conclude that

$$\kappa \|\boldsymbol{\omega}_h^n\|_{\boldsymbol{\Lambda}} \|\boldsymbol{v}_2\|_{1,\Omega} \le c(\boldsymbol{\omega}_h^n, \boldsymbol{v}_{2h}|_\Sigma).$$

Using the error equation (3.42), the relations (3.44) and the splitting $\boldsymbol{\omega}_h^n = (\boldsymbol{\lambda}^n - \boldsymbol{\lambda}_h^n) - (\boldsymbol{\lambda}^n - \boldsymbol{\lambda}_\pi^n) = \boldsymbol{\lambda}^n - \boldsymbol{\lambda}_h^n - \boldsymbol{\omega}_\pi^n$, we can write for $\boldsymbol{v}_h = \boldsymbol{v}_{2h}, q_h = p_{2h}, \boldsymbol{w}_h = \boldsymbol{0}, \boldsymbol{\lambda}_h = \boldsymbol{0}$

$$\begin{aligned} c(\boldsymbol{\omega}_h^n, \boldsymbol{v}_{2h}|_\Sigma) =& c(\boldsymbol{\lambda}^n - \boldsymbol{\lambda}_h^n, \boldsymbol{v}_{2h}|_\Sigma) - c(\boldsymbol{\omega}_\pi^n, \boldsymbol{v}_{2h}|_\Sigma) \\ =& - c(\boldsymbol{\omega}_\pi^n, \boldsymbol{v}_{2h}|_\Sigma) - \rho_{\mathrm{f}}\big((\partial_t - \partial_\tau)\boldsymbol{u}^n + \partial_\tau \boldsymbol{\theta}_\pi^n, \boldsymbol{v}_{2h}\big)_{0,\Omega} - \rho_{\mathrm{f}}\big(\partial_\tau \boldsymbol{\theta}_h^n, \boldsymbol{v}_{2h}\big)_{0,\Omega} \\ & - a^f((\boldsymbol{\theta}_h^n, \varphi_h^n); (\boldsymbol{v}_{2h}, p_{2h})) - s_h(\phi_h^n, p_{2h}). \end{aligned} \tag{3.51}$$

We estimate separately each term.

- By definition of the bilinear form $c(\cdot,\cdot)$ and the trace inequality

$$c(\boldsymbol{\omega}_\pi^n, \boldsymbol{v}_{2h}|_\Sigma) \le M_c \|\boldsymbol{\omega}_\pi^n\|_{\boldsymbol{\Lambda}} \|\boldsymbol{v}_{2h}\|_{\frac{1}{2},\Sigma} \le C \|\boldsymbol{\omega}_\pi\|_{L^\infty(0,T;\boldsymbol{\Lambda})} \|\boldsymbol{v}_{2h}\|_{1,\Omega} \le C \|\boldsymbol{\omega}_\pi\|_{L^\infty(0,T;\boldsymbol{\Lambda})} \|\boldsymbol{v}_2\|_{1,\Omega}.$$

- Using the same type of estimate already used for the term $T_1$ (3.48)

$$\rho_{\mathrm{f}}\big((\partial_t - \partial_\tau)\boldsymbol{u}^n + \partial_\tau \boldsymbol{\theta}_\pi^n, \boldsymbol{v}_{2h}\big)_{0,\Omega}$$

$$\leq \rho_{\mathrm{f}}\Big(\tau^{\frac{1}{2}}\|\partial_{tt}\boldsymbol{u}\|_{L^2(t_{n-1},t_n;L^2(\Omega)^d)} + \tau^{-\frac{1}{2}}\|\partial_t\boldsymbol{\theta}_\pi\|_{L^2(t_{n-1},t_n;L^2(\Omega)^d)}\Big)\|\boldsymbol{v}_{2h}\|_{1,\Omega}$$
$$\leq \rho_{\mathrm{f}}\Big(\tau^{\frac{1}{2}}\|\partial_{tt}\boldsymbol{u}\|_{L^2(t_{n-1},t_n;L^2(\Omega)^d)} + \tau^{-\frac{1}{2}}\|\partial_t\boldsymbol{\theta}_\pi\|_{L^2(t_{n-1},t_n;L^2(\Omega)^d)}\Big)\|\boldsymbol{v}_2\|_{1,\Omega}.$$

- Using the Cauchy-Schwartz inequality we have

$$\rho_{\mathrm{f}}\big(\partial_\tau\boldsymbol{\theta}_h^n, \boldsymbol{v}_{2h}\big)_{0,\Omega} \leq \|\partial_\tau\boldsymbol{\theta}_h^n\|_{0,\Omega}\|\boldsymbol{v}_{2h}\|_{1,\Omega} \leq \|\partial_\tau\boldsymbol{\theta}_h^n\|_{0,\Omega}\|\boldsymbol{v}_2\|_{1,\Omega}.$$

- We recall that for all $(\boldsymbol{u}_h, p_h), (\boldsymbol{v}_h, q_h) \in \boldsymbol{V}_h \times Q_h$ we have

$$a^f((\boldsymbol{u}_h, p_h);(\boldsymbol{v}_h, q_h)) \overset{def}{=} 2\mu\big(\boldsymbol{\epsilon}(\boldsymbol{u}_h), \boldsymbol{\epsilon}(\boldsymbol{v}_h)\big)_{0,\Omega} - \big(\operatorname{div}\boldsymbol{v}_h, p_h\big)_{0,\Omega} + \big(\operatorname{div}\boldsymbol{u}_h, q_h\big)_{0,\Omega},$$

hence, using the Cauchy-Schwartz inequality and the continuity of the solid and fluid bilinear forms, we have, for $\boldsymbol{v}_h = \boldsymbol{v}_{2h}$ and $p_h = p_{2h}$, exploiting properties $(2.42)_2$ and (2.43)

$$\begin{aligned}
&-a^f((\boldsymbol{\theta}_h^n, \varphi_h^n);(\boldsymbol{v}_{2h}, p_{2h})) - s_h(\phi_h^n, p_{2h}) \\
&\quad= -2\mu\big(\boldsymbol{\epsilon}(\boldsymbol{\theta}_h^n), \boldsymbol{\epsilon}(\boldsymbol{v}_{2h})\big)_{0,\Omega} + \big(\operatorname{div}\boldsymbol{v}_{2h}, \phi_h^n\big)_{0,\Omega} - \big(\operatorname{div}\boldsymbol{\theta}_h^n, p_{2h}\big)_{0,\Omega} - s_h(\phi_h^n, p_{2h}) \\
&\quad= -2\mu\big(\boldsymbol{\epsilon}(\boldsymbol{\theta}_h^n), \boldsymbol{\epsilon}(\boldsymbol{v}_{2h})\big)_{0,\Omega} - \big(\operatorname{div}\boldsymbol{\theta}_h^n, p_{2h}\big)_{0,\Omega} - 2s_h(\phi_h^n, p_{2h}) \\
&\quad\leq 2C\mu\|\boldsymbol{\epsilon}(\boldsymbol{\theta}_h^n)\|_{0,\Omega}\|\boldsymbol{v}_{2h}\|_{1,\Omega} + C_K\|\boldsymbol{\epsilon}(\boldsymbol{\theta}_h^n)\|_{0,\Omega}\|p_{2h}\|_{0,\Omega} + 2|\phi_h^n|_{s_h}|p_{2h}|_{s_h} \\
&\quad\leq C\left(\|\boldsymbol{\epsilon}(\boldsymbol{\theta}_h^n)\|_{0,\Omega} + |\phi_h^n|_{s_h}\right)\|\boldsymbol{v}_2\|_{1,\Omega}.
\end{aligned}$$

In sum we have

$$\begin{aligned}
\|\boldsymbol{\omega}_h^n\|_{\boldsymbol{\Lambda}} \leq C\Big(&\tau^{\frac{1}{2}}\|\partial_{tt}\boldsymbol{u}\|_{L^2(t_{n-1},t_n;L^2(\Omega)^d)} + \tau^{-\frac{1}{2}}\|\partial_t\boldsymbol{\theta}_\pi\|_{L^2(t_{n-1},t_n;L^2(\Omega)^d)} + \|\boldsymbol{\omega}_\pi\|_{L^\infty(0,T;\boldsymbol{\Lambda})} \\
&+ \|\boldsymbol{\epsilon}(\boldsymbol{\theta}_h^n)\|_{0,\Omega} + |\phi_h^n|_{s_h}\Big) + C\|\partial_\tau\boldsymbol{\theta}_h^n\|_{0,\Omega}.
\end{aligned}$$

and, exploiting the regularity assumptions and the following Sobolev Inclusions

- since $\boldsymbol{\theta}_\pi \in H^2(0,T;H^{1+l}(\Omega)) \hookrightarrow L^\infty(0,T;H^{1+l}(\Omega))$ we have

$$\|\boldsymbol{\theta}_\pi^n\|_{1,\Omega} \leq \|\boldsymbol{\theta}_\pi\|_{L^\infty(0,T;H^1(\Omega)^d)}, \tag{3.52}$$

- since $\dot{\boldsymbol{\xi}}_\pi \in H^2(0,T;H^{1+m}(\Sigma)) \hookrightarrow L^\infty(0,T;H^{1+m}(\Sigma))$ we have

$$\|\dot{\boldsymbol{\xi}}_\pi^n\|_{\frac{1}{2},\Sigma} \leq \|\dot{\boldsymbol{\xi}}_\pi\|_{L^\infty(0,T;H^{\frac{1}{2}}(\Sigma)^d)}, \tag{3.53}$$

then the estimate of the term $T_4$ follows

$$\begin{aligned}T_4 \leq & C\delta_4 \|\partial_t \boldsymbol{\theta}_\pi\|^2_{L^2(t_{n-1},t_n;L^2(\Omega)^d)} + \frac{C}{\delta_4}\left[\|\dot{\boldsymbol{\xi}}_\pi\|^2_{L^\infty(0,T;H^{\frac{1}{2}}(\Sigma)^d)} + \|\boldsymbol{\theta}_\pi\|^2_{L^\infty(0,T;H^1(\Omega)^d)}\right] \\ & + C\tau\delta_4\|\boldsymbol{\omega}_\pi\|^2_{L^\infty(0,T;\boldsymbol{\Lambda})} + \frac{C\tau}{\delta_4}\left[\|\dot{\boldsymbol{\xi}}_\pi\|^2_{L^\infty(0,T;H^{\frac{1}{2}}(\Sigma)^d)} + \|\boldsymbol{\theta}_\pi\|^2_{L^\infty(0,T;H^1(\Omega)^d)}\right] \\ & + C\tau^2\delta_4\|\partial_{tt}\boldsymbol{u}\|^2_{L^2(t_{n-1},t_n;L^2(\Omega)^d)} + C\delta_4\|\tau\partial_\tau\boldsymbol{\theta}_h^n\|^2_{0,\Omega} + C\tau\delta_4\|\boldsymbol{\epsilon}(\boldsymbol{\theta}_h^n)\|^2_{0,\Omega} + C\tau\delta_4|\phi_h^n|^2_{s_h}.\end{aligned} \tag{3.54}$$

Using the previous estimates (3.47), (3.48), (3.49), (3.50) and (3.54) in (3.46) we obtain

$$\begin{aligned}&\left(\frac{\rho_{\mathrm{f}}}{2}\tau\partial_\tau\|\boldsymbol{\theta}_h^n\|^2_{0,\Omega} + \frac{\rho_s\epsilon}{2}\tau\partial_\tau\|\dot{\boldsymbol{\xi}}_h^n\|^2_{0,\Sigma} + \frac{1}{2}\tau\partial_\tau\|\boldsymbol{\xi}_h^n\|^2_{\mathrm{s}}\right) \\ &+\tau^2\left(\left(\frac{\rho_{\mathrm{f}}}{2} - C\delta_4\right)\|\partial_\tau\boldsymbol{\theta}_h^n\|^2_{0,\Omega} + \frac{\rho_s\epsilon}{2}\|\partial_\tau\dot{\boldsymbol{\xi}}_h^n\|^2_{0,\Sigma} + \frac{1}{2}\|\partial_\tau\boldsymbol{\xi}_h^n\|^2_{\mathrm{s}}\right) \\ &+\tau\left(2\mu - \delta_3 C_K - C\delta_4\right)\|\boldsymbol{\epsilon}(\boldsymbol{\theta}_h^n)\|^2_{0,\Omega} + \tau(1 - C\delta_4)|\varphi_h^n|^2_{s_h} \\ &\leq \\ &\left(\frac{\rho_{\mathrm{f}}}{2\delta_1}\|\partial_t\boldsymbol{\theta}_\pi\|^2_{L^2(t_{n-1},t_n;L^2(\Omega)^d)} + \frac{\rho_s\epsilon}{2\delta_2}\|\partial_t\dot{\boldsymbol{\xi}}_\pi\|^2_{L^2(t_{n-1},t_n;L^2(\Sigma)^d)} + C\delta_4\|\partial_t\boldsymbol{\theta}_\pi\|^2_{L^2(t_{n-1},t_n;L^2(\Omega)^d)}\right. \\ &\left.+\frac{C}{\delta_4}\|\dot{\boldsymbol{\xi}}_\pi\|^2_{L^\infty(0,T;H^{\frac{1}{2}}(\Sigma)^d)} + \frac{C}{\delta_4}\|\boldsymbol{\theta}_\pi\|^2_{L^\infty(0,T;H^1(\Omega)^d)}\right) \\ &+\tau\left(\frac{1}{\delta_3}\|\boldsymbol{\omega}_\pi\|^2_{L^\infty(0,T;\boldsymbol{\Lambda})} + C\delta_4\|\boldsymbol{\omega}_\pi\|^2_{L^\infty(0,T;\boldsymbol{\Lambda})} + \frac{C}{\delta_4}\|\dot{\boldsymbol{\xi}}_\pi\|^2_{L^\infty(0,T;H^{\frac{1}{2}}(\Sigma)^d)} + \frac{C}{\delta_4}\|\boldsymbol{\theta}_\pi\|^2_{L^\infty(0,T;H^1(\Omega)^d)}\right) \\ &+\tau^2\left(\frac{1}{2\delta_0}\|\partial_{tt}\boldsymbol{d}\|^2_{L^2(t_{n-1},t_n;\boldsymbol{W})} + \frac{\rho_{\mathrm{f}}}{2\delta_1}\|\partial_{tt}\boldsymbol{u}\|^2_{L^2(t_{n-1},t_n;L^2(\Omega)^d)}\right. \\ &\left.+\frac{\rho_s\epsilon}{2\delta_2}\|\partial_{tt}\dot{\boldsymbol{d}}\|^2_{L^2(t_{n-1},t_n;L^2(\Sigma)^d)} + C\delta_4\|\partial_{tt}\boldsymbol{u}\|^2_{L^2(t_{n-1},t_n;L^2(\Omega)^d)}\right) \\ &+\rho_{\mathrm{f}}\tau\delta_1\|\boldsymbol{\theta}_h^n\|^2_{0,\Omega} + \frac{\delta_0\tau}{2}\|\boldsymbol{\xi}_h^n\|^2_{\mathrm{s}} + \rho_s\epsilon\tau\delta_2\|\dot{\boldsymbol{\xi}}_h^n\|^2_{0,\Sigma}.\end{aligned}$$

Considering the following values of the parameters

$$\delta_0 = \delta_1 = \delta_2 = \frac{1}{2}, \qquad \delta_3 = \frac{\mu}{C}, \qquad \delta_4 = \min\left\{\frac{1}{C}, \frac{\rho_f}{2}, \mu\right\}, \tag{3.55}$$

and summing over $m = 1..n$ we obtain

$$\begin{aligned}&\underbrace{\left(\frac{\rho_{\mathrm{f}}}{2}\|\boldsymbol{\theta}_h^n\|^2_{0,\Omega} + \frac{\rho_s\epsilon}{2}\|\dot{\boldsymbol{\xi}}_h^n\|^2_{0,\Sigma} + \frac{1}{2}\|\boldsymbol{\xi}_h^n\|^2_{\mathrm{s}}\right)}_{\mathscr{E}_h^n} - \underbrace{\left(\frac{\rho_{\mathrm{f}}}{2}\|\boldsymbol{\theta}_h^0\|^2_{0,\Omega} + \frac{\rho_s\epsilon}{2}\|\dot{\boldsymbol{\xi}}_h^0\|^2_{0,\Sigma} + \frac{1}{2}\|\boldsymbol{\xi}_h^0\|^2_{\mathrm{s}}\right)}_{\mathscr{E}_h^0} \\ &\tau\sum_{m=1}^n\left[C\|\boldsymbol{\epsilon}(\boldsymbol{\theta}_h^m)\|^2_{0,\Omega} + C|\varphi_h^m|^2_{s_h}\right] + \tau^2\sum_{m=1}^n\left(C\|\partial_\tau\boldsymbol{\theta}_h^m\|^2_{0,\Omega} + \frac{\rho_s\epsilon}{2}\|\partial_\tau\dot{\boldsymbol{\xi}}_h^m\|^2_{0,\Sigma} + \frac{1}{2}\|\partial_\tau\boldsymbol{\xi}_h^m\|^2_{\mathrm{s}}\right) \\ &\leq \underbrace{\sum_{m=1}^n C_m^0}_{A_0} + \tau\underbrace{\sum_{m=1}^n C_m^1}_{A_1} + \tau^2\underbrace{\sum_{m=1}^n C_m^2}_{A_2} + \tau\sum_{m=1}^n\underbrace{\left(\frac{\rho_{\mathrm{f}}}{2}\|\boldsymbol{\theta}_h^m\|^2_{0,\Omega} + \frac{\rho_s\epsilon}{2}\|\dot{\boldsymbol{\xi}}_h^m\|^2_{0,\Sigma} + \frac{1}{2}\|\boldsymbol{\xi}_h^m\|^2_{\mathrm{s}}\right)}_{\mathscr{E}_h^m}.\end{aligned} \tag{3.56}$$

where the terms in the right hand side of the inequality are defined ad analyzed in the following using results of Proposition 3.3.6.

- The terms $A_0$ can be bounded as follows

$$
\begin{aligned}
A_0 \stackrel{def}{=} C \sum_{n=1}^{m} \Big( & \|\partial_t \boldsymbol{\theta}_\pi\|^2_{L^2(t_{n-1},t_n;L^2(\Omega)^d)} + \|\partial_t \dot{\boldsymbol{\xi}}_\pi\|^2_{L^2(t_{n-1},t_n;L^2(\Sigma)^d)} \\
& + \|\dot{\boldsymbol{\xi}}_\pi\|^2_{L^\infty(0,T;H^{\frac{1}{2}}(\Sigma)^d)} + \|\boldsymbol{\theta}_\pi\|^2_{L^\infty(0,T;H^1(\Omega)^d)} \Big) \\
\leq & C\Big( \|\partial_t \boldsymbol{\theta}_\pi\|^2_{L^2(0,T;L^2(\Omega)^d)} + \|\partial_t \dot{\boldsymbol{\xi}}_\pi\|^2_{L^2(0,T;L^2(\Sigma)^d)} + \|\partial_t \boldsymbol{\theta}_\pi\|^2_{L^2(0,T;L^2(\Omega)^d)} \\
& + \|\dot{\boldsymbol{\xi}}_\pi\|^2_{L^\infty(0,T;H^{\frac{1}{2}}(\Sigma)^d)} + \|\boldsymbol{\theta}_\pi\|^2_{L^\infty(0,T;H^1(\Omega)^d)} \Big) \\
& \approx \mathcal{O}(h_f^{2l+2}) + \mathcal{O}(h_s^{2m+2}) + \mathcal{O}(h_f^{2l+2}) + \mathcal{O}(h_f^{2l}) + \mathcal{O}(h_s^{2m+1}).
\end{aligned}
\tag{3.57}
$$

- The term $A_1$ is defined and bounded as follows

$$
\begin{aligned}
A_1 \stackrel{def}{=} C \sum_{n=1}^{m} & \Big( \|\boldsymbol{\omega}_\pi\|^2_{L^\infty(0,T;\boldsymbol{\Lambda})} + \|\dot{\boldsymbol{\xi}}_\pi\|^2_{L^\infty(0,T;H^{\frac{1}{2}}(\Sigma)^d)} + \|\boldsymbol{\theta}_\pi\|^2_{L^\infty(0,T;H^1(\Omega)^d)} \Big) \\
& \approx \mathcal{O}(h_s^{2z+1}) + \mathcal{O}(h_s^{2m+1}) + \mathcal{O}(h_f^{2l}).
\end{aligned}
\tag{3.58}
$$

- The term $A_2$ is defined and bounded in the following

$$
\begin{aligned}
A_2 \stackrel{def}{=} C \sum_{n=1}^{m} & \Big( \|\partial_{tt} \boldsymbol{d}\|^2_{L^2(t_{n-1},t_n;\mathbf{W})} + \|\partial_{tt} \boldsymbol{u}\|^2_{L^2(t_{n-1},t_n;L^2(\Omega)^d)} + \|\partial_{tt} \dot{\boldsymbol{d}}\|^2_{L^2(t_{n-1},t_n;L^2(\Sigma)^d)} \Big) \\
\leq & C\Big( \|\partial_{tt} \boldsymbol{d}\|^2_{L^2(0,T;\mathbf{W})} + \|\partial_{tt} \boldsymbol{u}\|^2_{L^2(0,T;L^2(\Omega)^d)} + \|\partial_{tt} \dot{\boldsymbol{d}}\|^2_{L^2(0,T;L^2(\Sigma)^d)} \Big).
\end{aligned}
\tag{3.59}
$$

In sum the error estimate (3.56) can be written as follows

$$
\mathcal{E}_h^n + \tau \mathcal{D}_h^n \leq \tau \sum_{m=1}^{n} \mathcal{E}_h^m + A_0 + \tau A_1 + \tau^2 A_2
\tag{3.60}
$$

where $\tau\mathcal{D}$ represents the dissipation term and is given by

$$
\tau \mathcal{D}_h^n = \tau \sum_{m=1}^{n} \Big[ C\|\boldsymbol{\epsilon}(\boldsymbol{\theta}_h^m)\|^2_{0,\Omega} + C|\varphi_h^m|^2_{s_h} \Big] + \tau^2 \sum_{m=1}^{n} \Big( C\|\partial_\tau \boldsymbol{\theta}_h^m\|^2_{0,\Omega} + \frac{\rho_s \epsilon}{2} \|\partial_\tau \dot{\boldsymbol{\xi}}_h^m\|^2_{0,\Sigma} + \frac{1}{2}\|\partial_\tau \boldsymbol{\xi}_h^m\|^2_{\mathrm{s}} \Big).
$$

Introducing the convergence rates with respect to $h_f$ and $h_s$ of the term $A_0, A_1, A_2$ given in (3.57) and (3.58), we have the desired result applying the discrete Gronwall Lemma 3.3.1, since by hypothesis we have $\tau < 1$

$$
\mathcal{E}_h^n + \tau \mathcal{D}_h^n \lesssim \exp\left( \sum_{m=1}^{n} \frac{\tau}{1-\tau} \right) \Bigg[ \underbrace{\mathcal{O}(h_f^{2l}) + \mathcal{O}(h_s^{2m+1})}_{A_0} + \tau \Bigg( \underbrace{\mathcal{O}(h_s^{2z+1}) + \mathcal{O}(h_s^{2m+1}) + \mathcal{O}(h_f^{2l})}_{A_1} \Bigg) + \tau^2 A_2 \Bigg].
\tag{3.61}
$$

□

## 3.4 Numerical experiments

In this section, we perform numerical tests in order to check numerically the rate of convergence of the monolithic scheme. All the numerical tests are performed using the classical benchmark problem of an ellipsoidal structure that evolves to a circular equilibrium position. The matrix form of the monolithic algorithm has already been presented in (2.47). In order to check the convergence rate of the algorithm, we consider a **reference solution**, obtained considering the following discretization parameters

$$h_f^{ref} = 1/256, \qquad h_s^{ref} = 1/256, \qquad \tau^{ref} = 0.00005s, \tag{3.62}$$

and we evaluate the errors with respect to this solution. In the following paragraphs we perform three types of convergence analysis, namely, we test the convergence with respect to time-step size when the fluid and solid mesh sizes are fixed, then we check the convergence with respect to the fluid and solid mesh sizes for fixed time-step size, lastly we test the convergence when time-step size and the spatial mesh sizes go to zero simultaneously.

**Convergence with respect to time step size.** The first test performed is suited in order to evaluate the convergence rate with respect to the time discretization. In this regard, we consider the test case of the elliptical structure that evolves to a circular equilibrium position. The discretization parameters are chosen in such a way that the effect of the error due to spatial discretization is negligible, that is we consider

$$h_f = 1/256, \qquad h_s = 1/256.$$

The time step $\tau$ has been varied in the range

$$\tau \in \{0.064, 0.032, 0.016, 0.008, 0.004\},$$

Using the convergence estimate (3.40) in the case $h_f = h_s$ we expect that the total error

$$\text{Total Error} \stackrel{def}{=} \mathscr{E}_h^n = C\left(\|\boldsymbol{u}_h^n - \boldsymbol{u}_{ref}^n\|_{0,\Omega}^2 + \|\boldsymbol{d}_h^n - \boldsymbol{d}_{ref}^n\|_{s,\Sigma}^2 + \|\dot{\boldsymbol{d}}_h^n - \dot{\boldsymbol{d}}_{ref}^n\|_{0,\Sigma}^2\right)^{\frac{1}{2}}$$

goes to zero as $\tau$.

In the following table, we report the error computed for the fluid velocity, the solid displacement, and the solid velocity; all the errors are evaluated at time $t = 0.064s$. The last column display the convergence rates of the total error which is in agreement with the theoretical results.

| $\tau$ | $\|\boldsymbol{u}_h^n - \boldsymbol{u}_{ref}^n\|_{0,\Omega}$ | $\|\boldsymbol{d}_h^n - \boldsymbol{d}_{ref}^n\|_{s,\Sigma}$ | $\|\dot{\boldsymbol{d}}_h^n - \dot{\boldsymbol{d}}_{ref}^n\|_{0,\Sigma}$ | Total Error | Total Rate |
|---|---|---|---|---|---|
| 0.064 | 0.0106936667 | 0.0554678227 | 0.0162457138 | 0.058778883 | |
| 0.032 | 0.0066462263 | 0.0312720022 | 0.0084203747 | 0.0330607495 | 8.30178513e-01 |
| 0.016 | 0.003750269 | 0.0164390978 | 0.0045079773 | 0.0174536619 | 9.21589675e-01 |
| 0.008 | 0.001959067 | 0.0082047341 | 0.0023974795 | 0.0087694648 | 9.92969053e-01 |
| 0.004 | 0.0009522398 | 0.0038842392 | 0.0011954948 | 0.0041741206 | 1.07101650e+00 |

Table 3.1: Convergence w.r.t time-step size of the Monolithic algorithm at time $t = 0.064s$.

Moreover, the total error is plotted in Figure 3.1, where it also possible to appreciate that the convergence rate is 1. This point is particularly interesting since the kinematic constraint is enforced in an explicit way, namely using the domain configuration at the previous time step. This choice, as it is shown in Table 3.1 and in Figure 3.1, does not affect the first order of convergence of the Euler scheme.

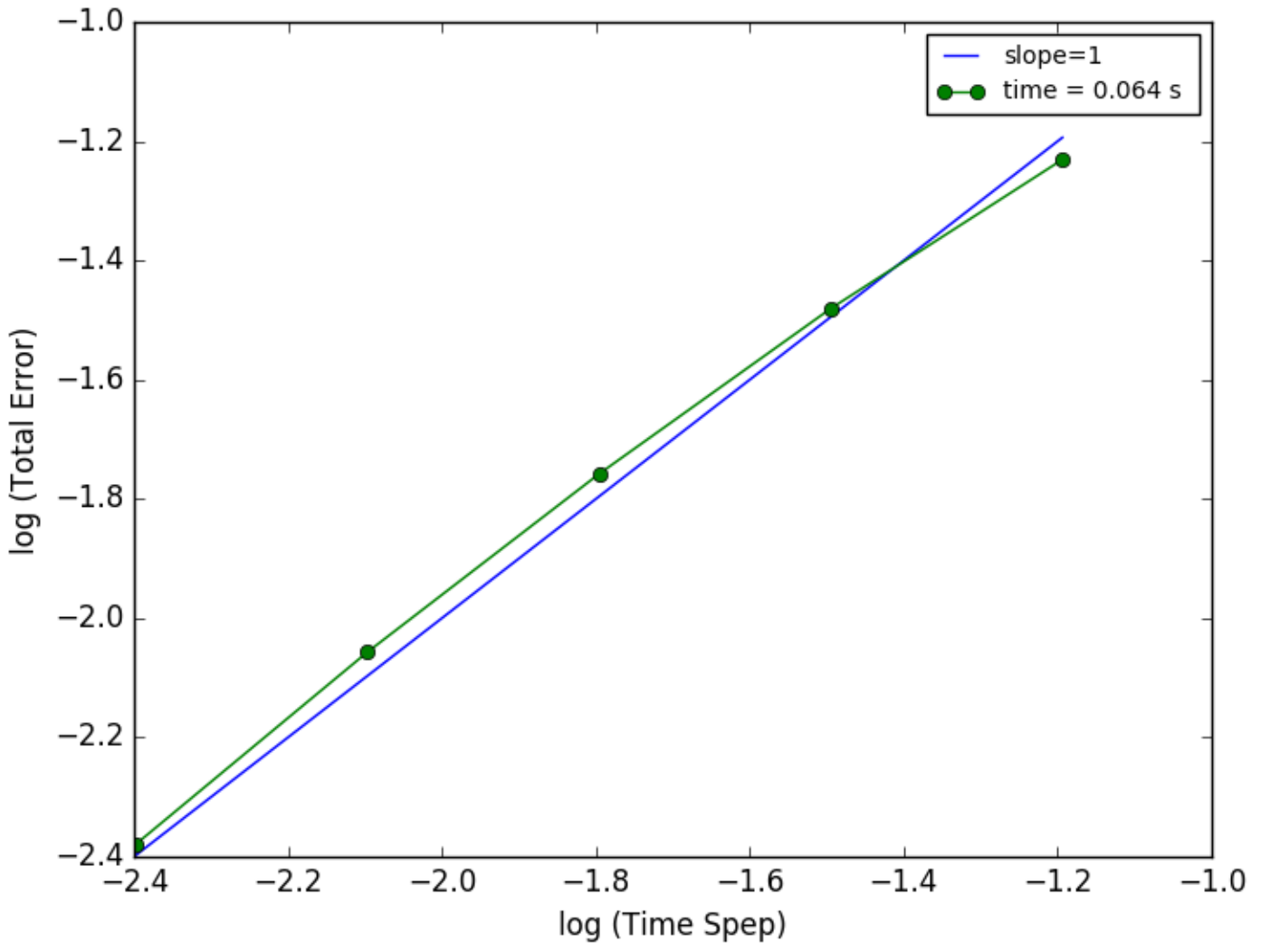

Figure 3.1: Total Error vs the time-step size of the Monolithic algorithm at time $t = 0.064s$.

**Convergence with respect to the spatial mesh sizes** In this paragraph, we analyze the error behavior with respect to the fluid and solid mesh sizes, $h_f$ and $h_s$ respectively for fixed value of the time step size. In the numerical experiments, we consider the following discretization parameters

$$h_f = h_s \in \{1/8, 1/16, 1/32, 1/64, 1/128\}, \qquad \tau = 0.00005s$$

We observe that the time step size is chosen in such a way to be negligible with respect to the spatial mesh sizes. Owing to (3.40), we expect that the total error goes to zero like $h_f^l$, with "$l$" close to $\frac{1}{2}$. In fact, in the following table, we report the computed errors that are in good agreement with the theoretical results. The simulations are performed considering a simple steady test where the initial position of the structure is a circle immersed in the fluid at rest, so that the asymptotic configuration remains unchanged.

| $h_f = h_s$ | $\\|\boldsymbol{u}_h^n - \boldsymbol{u}_{ref}^n\\|_{0,\Omega}$ | $\\|\boldsymbol{d}_h^n - \boldsymbol{d}_{ref}^n\\|_{s,\Sigma}$ | $\\|\dot{\boldsymbol{d}}_h^n - \dot{\boldsymbol{d}}_{ref}^n\\|_{0,\Sigma}$ | Total Error | Total Rate |
|---|---|---|---|---|---|
| 1/8 | 6.03323E-04 | 3.11070E-02 | 4.76608E-04 | 3.11165E-02 | |
| 1/16 | 5.32387E-04 | 1.61895E-02 | 4.52774E-04 | 1.62046E-02 | 9.41280E-01 |
| 1/32 | 4.58502E-04 | 8.55449E-03 | 3.57483E-04 | 8.57422E-03 | 9.18325E-01 |
| 1/64 | 3.43716E-04 | 4.76883E-03 | 2.41294E-04 | 4.78729E-03 | 8.40797E-01 |
| 1/128 | 1.92169E-04 | 2.86237E-03 | 1.35411E-04 | 2.87201E-03 | 7.37149E-01 |

Table 3.2: Convergence w.r.t $h_f = h_s$ of the Monolithic algorithm at time $t = 0.0001s$ for fixed time step size $\tau = 0.00005s$.

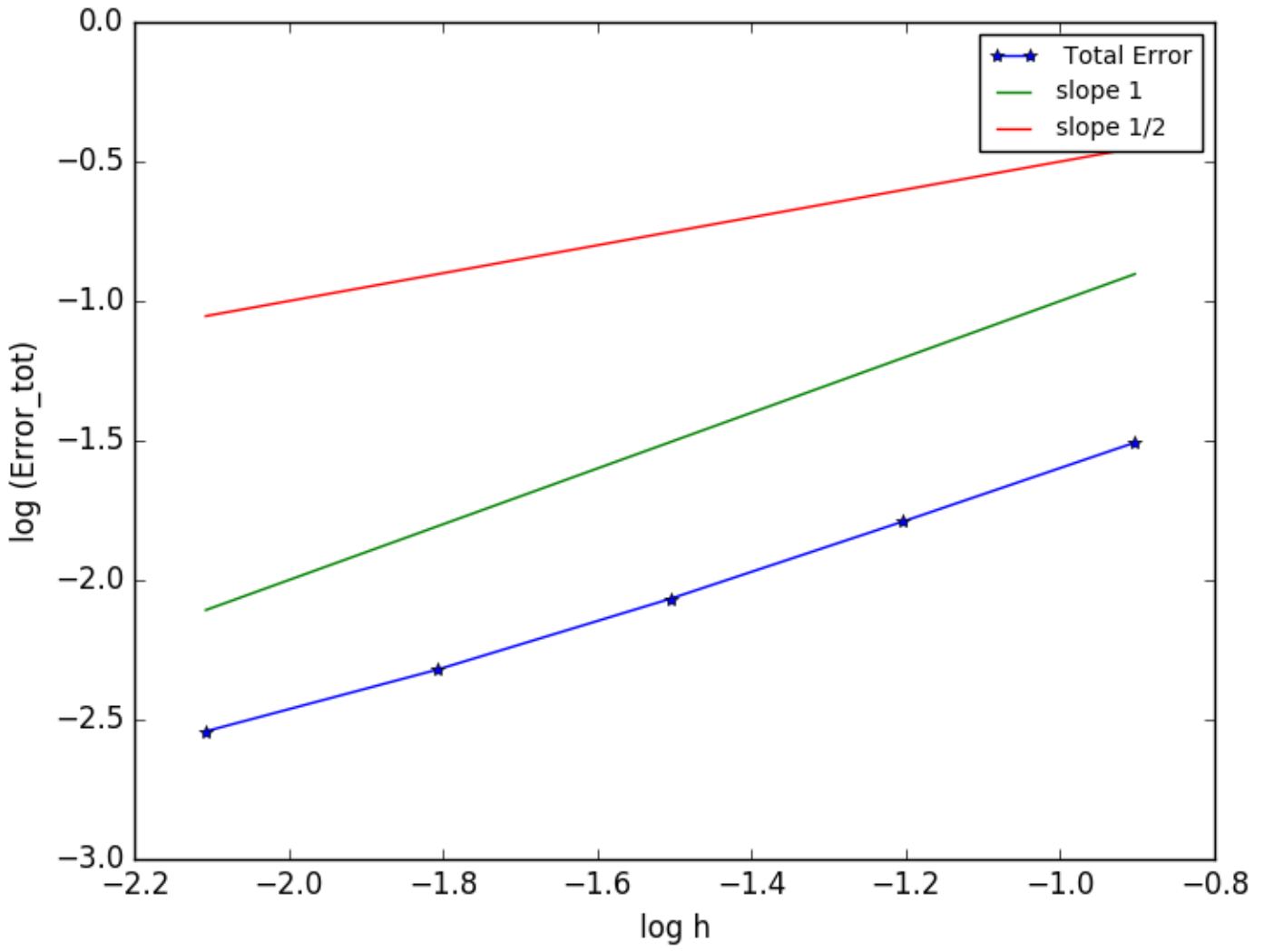

Figure 3.2: Total Error vs the mesh size of the Monolithic algorithm at time $t = 0.0001s$. for fixed time-step size $\tau = 0.00005s$

**Global Convergence rate** In this paragraph, we evaluate the global convergence rate of the monolithic algorithm considering the following discretization parameters

$$h_f = h_s = \frac{0.125}{2^j}, \qquad \tau = \frac{0.064}{2^j} \qquad \text{for} \qquad j = 0, 1, 2, 3, 4.$$

As in the previous cases, we consider a simple steady test where the initial position of the structure is a circle immersed in the fluid at rest, so that the asymptotic configuration remains unchanged. The reference solution is computed using the parameters given in (3.62). The following tables and figures report the rates of convergence which result in agreement with the supposed regularity of the solution. In particular in tables and figures 3.3 and 3.4 we report the convergence rates of the fluid velocity and pressure respectively. Since the pressure is discontinuous across the structure, the optimal convergence rate is $\frac{3}{2}$ for the fluid velocity in $L^2$ and $\frac{1}{2}$ for the pressure in $L^2$. Concerning the solid displacements and velocities, the convergence rates are reported in tables and figures 3.5 and 3.6 respectively. According with the theoretical results, the convergence rate is 1 for the solid displacement in the "energy norm" associated to the bilinear form $a^s(\cdot,\cdot)$, or equivalently in $H^1$ norm. The numerical results are in agreement with the theoretical expectation. The computed convergence rate for the solid velocity results to be 1 when evaluated in $L^2$ norm. In the last table we report the convergence rate of the total error evaluated with respect to fluid and solid mesh sizes. The numerical evaluation of the convergence rate in that case $h_f = h_s$ results in agreement with the theoretical result obtained in Theorem 3.3.7. In fact, for $h_f = h_s$ the estimate (3.40) gives a global rate close to 1 as reported in table and figure 3.7

| $h_f = h_s$ | $\|\boldsymbol{u}_h^n - \boldsymbol{u}_{ref}^n\|_{0,\Omega}$ | Rate |
|---|---|---|
| 1/8 | 1.11909E-02 | |
| 1/16 | 6.24951E-03 | 8.40505E-01 |
| 1/32 | 2.31694E-03 | 1.43152E+00 |
| 1/64 | 8.43140E-04 | 1.45838E+00 |
| 1/128 | 2.95381E-04 | 1.51320E+00 |

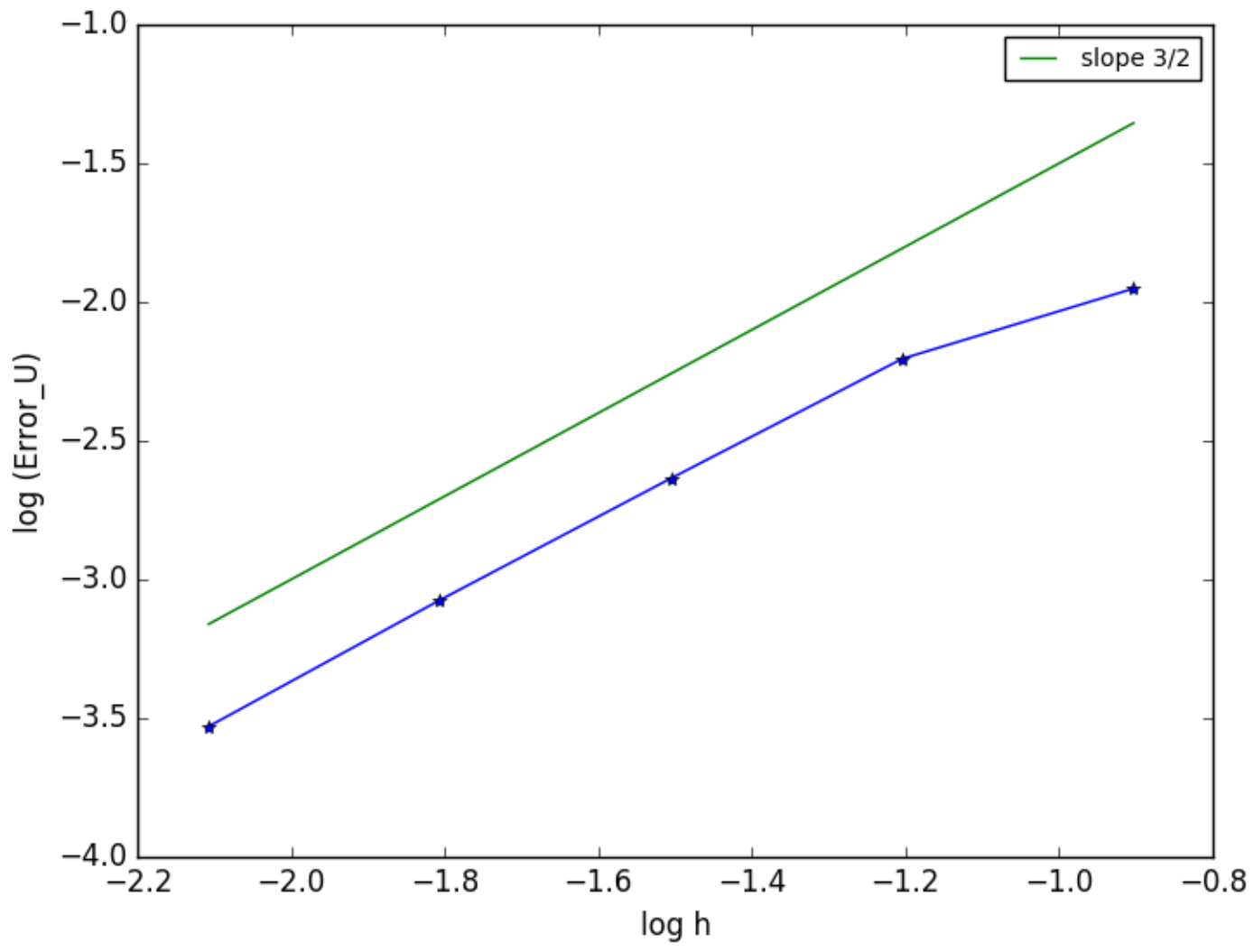

Table 3.3: Convergence rate of the fluid velocity with respect to the spatial mesh sizes $h_f$ and $h_s$ evaluated for decreasing values of $h_f, h_s, \tau$ at time $t = 0.064s$.

| $h_f = h_s$ | $\|\boldsymbol{p}_h^n - \boldsymbol{p}_{ref}^n\|_{0,\Omega}$ | Rate |
|---|---|---|
| 1/8 | 1.85807E+00 | |
| 1/16 | 1.32246E+00 | 4.90578E-01 |
| 1/32 | 9.70482E-01 | 4.46455E-01 |
| 1/64 | 7.36350E-01 | 3.98310E-01 |
| 1/128 | 5.88048E-01 | 3.24457E-01 |

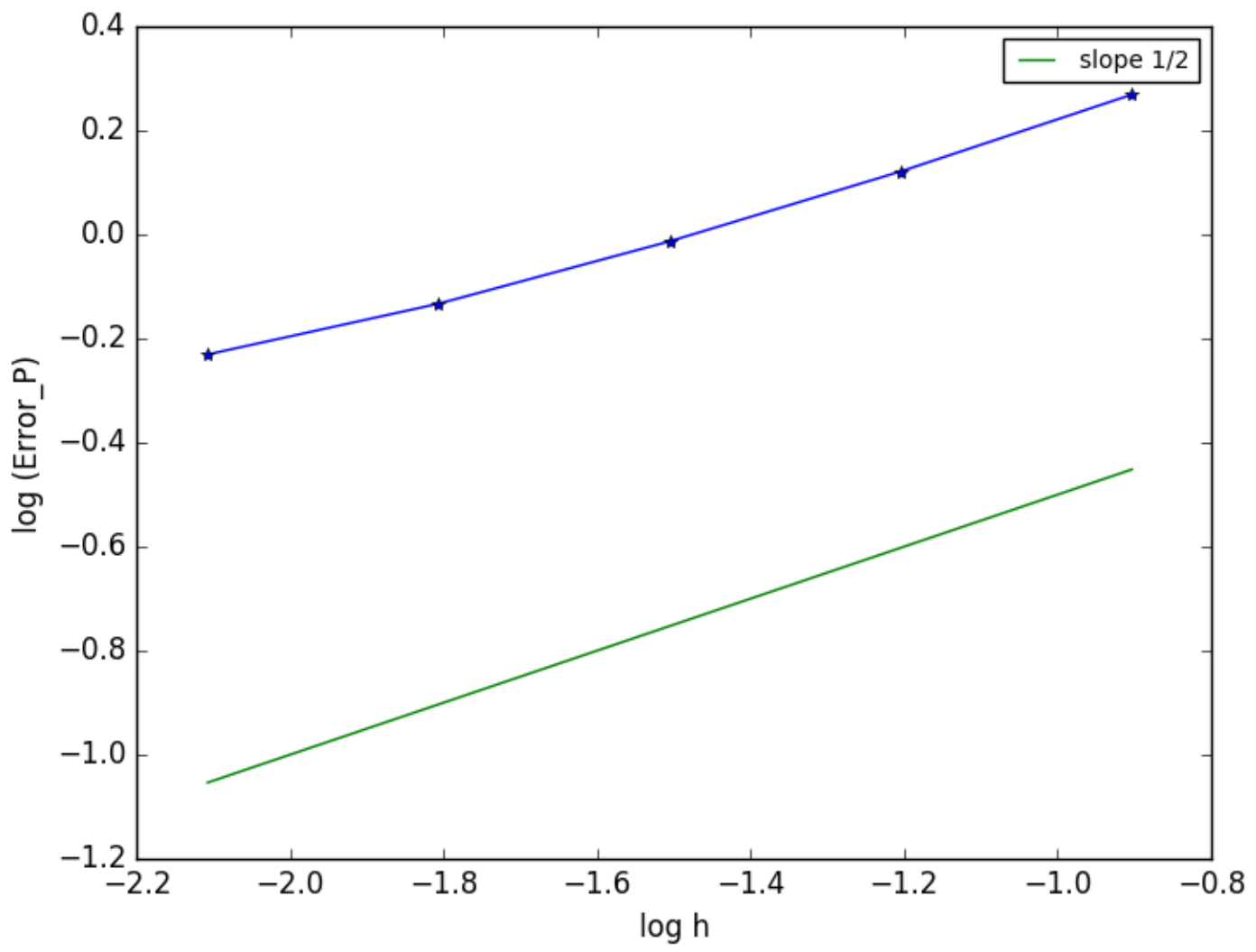

Table 3.4: Convergence rate of the fluid pressure with respect to the spatial mesh sizes $h_f$ and $h_s$ evaluated for decreasing values of $h_f, h_s, \tau$ at time $t = 0.064s$.

| $h_f = h_s$ | $\|\boldsymbol{d}_h^n - \boldsymbol{d}_{ref}^n\|_{s,\Sigma}$ | Rate |
|---|---|---|
| 1/8 | 3.04363E-02 | |
| 1/16 | 1.58093E-02 | 9.45024E-01 |
| 1/32 | 8.33133E-03 | 9.24151E-01 |
| 1/64 | 4.68751E-03 | 8.29726E-01 |
| 1/128 | 2.83502E-03 | 7.25464E-01 |

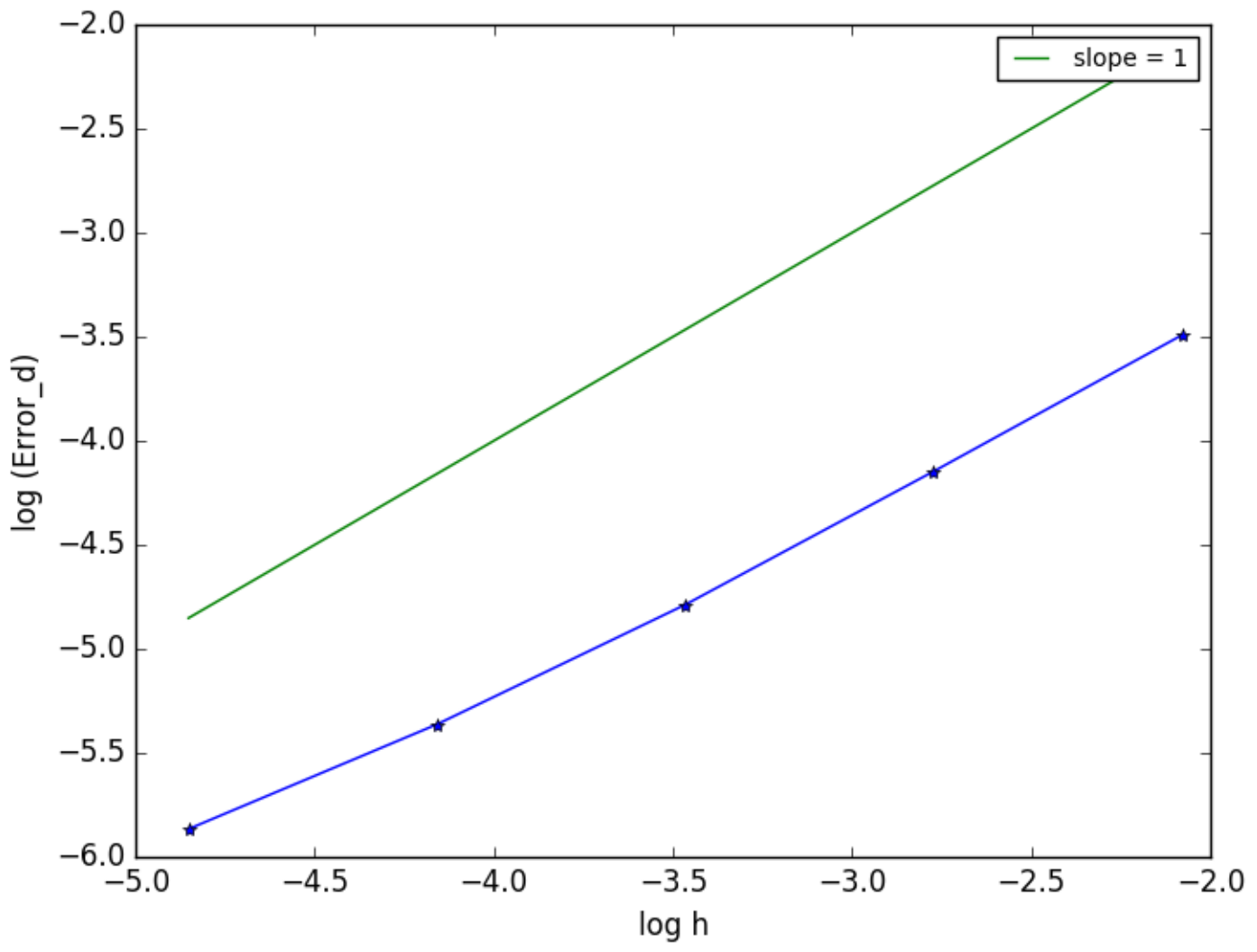

Table 3.5: Convergence rate of the solid displacements with respect to the spatial mesh sizes $h_f$ and $h_s$ evaluated for decreasing values of $h_f, h_s, \tau$ at time $t = 0.064s$.

| $h_f = h_s$ | $\|\dot{\boldsymbol{d}}_h^n - \dot{\boldsymbol{d}}_{ref}^n\|_{0,\Sigma}$ | Rate |
|---|---|---|
| 1/8 | 7.01531E-04 | |
| 1/16 | 6.10178E-04 | 2.01276E-01 |
| 1/32 | 3.72660E-04 | 7.11370E-01 |
| 1/64 | 1.72963E-04 | 1.10740E+00 |
| 1/128 | 7.86480E-05 | 1.13698E+00 |

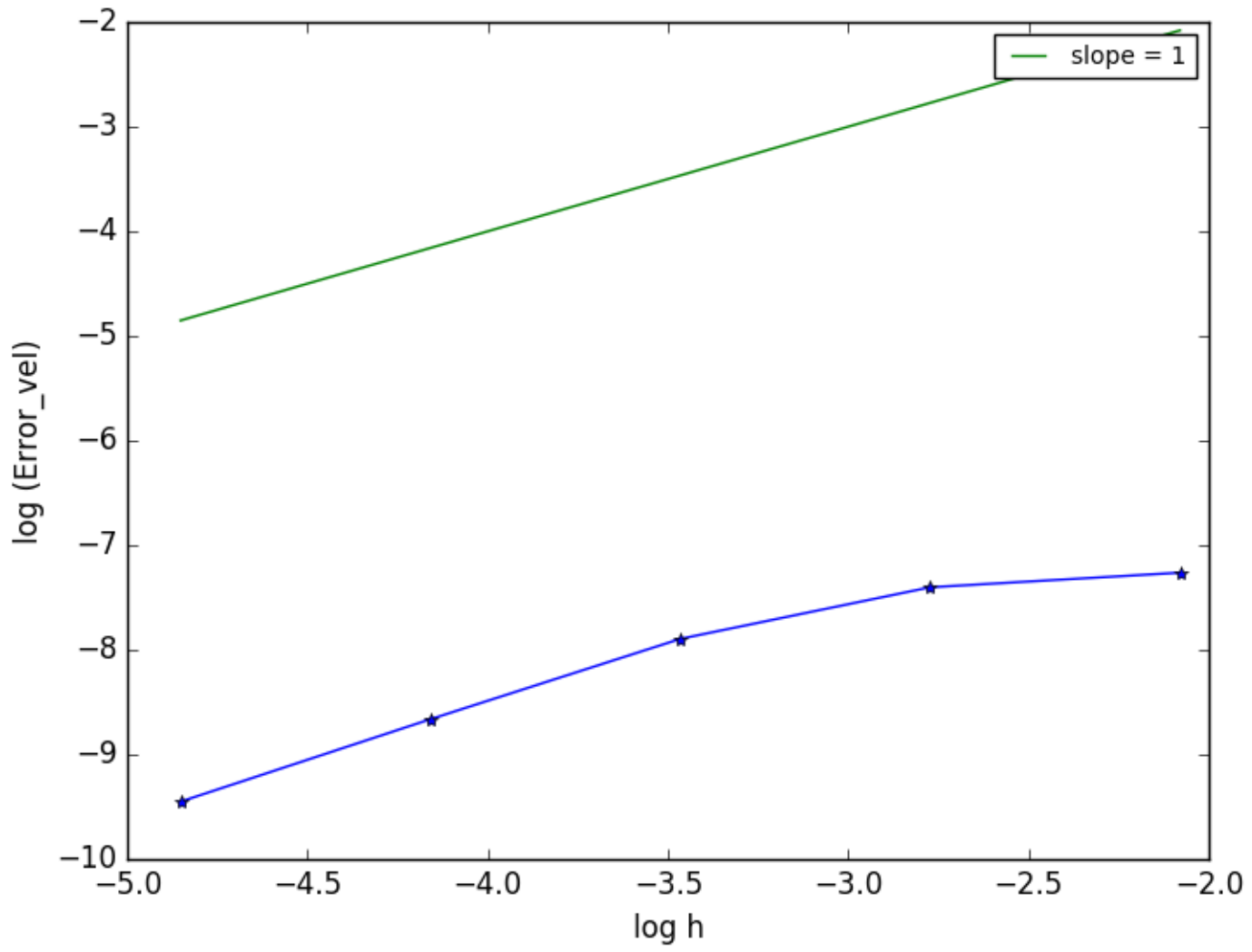

Table 3.6: Convergence rate of the solid velocity with respect to the spatial mesh sizes $h_f$ and $h_s$ evaluated for decreasing values of $h_f, h_s, \tau$ at time $t = 0.064s$.

| $h_f = h_s$ | Total Error | Total Rate |
|---|---|---|
| 1/8 | 3.24360E-02 | |
| 1/16 | 1.70106E-02 | 9.31162E-01 |
| 1/32 | 8.65553E-03 | 9.74742E-01 |
| 1/64 | 4.76587E-03 | 8.60882E-01 |
| 1/128 | 2.85145E-03 | 7.41045E-01 |

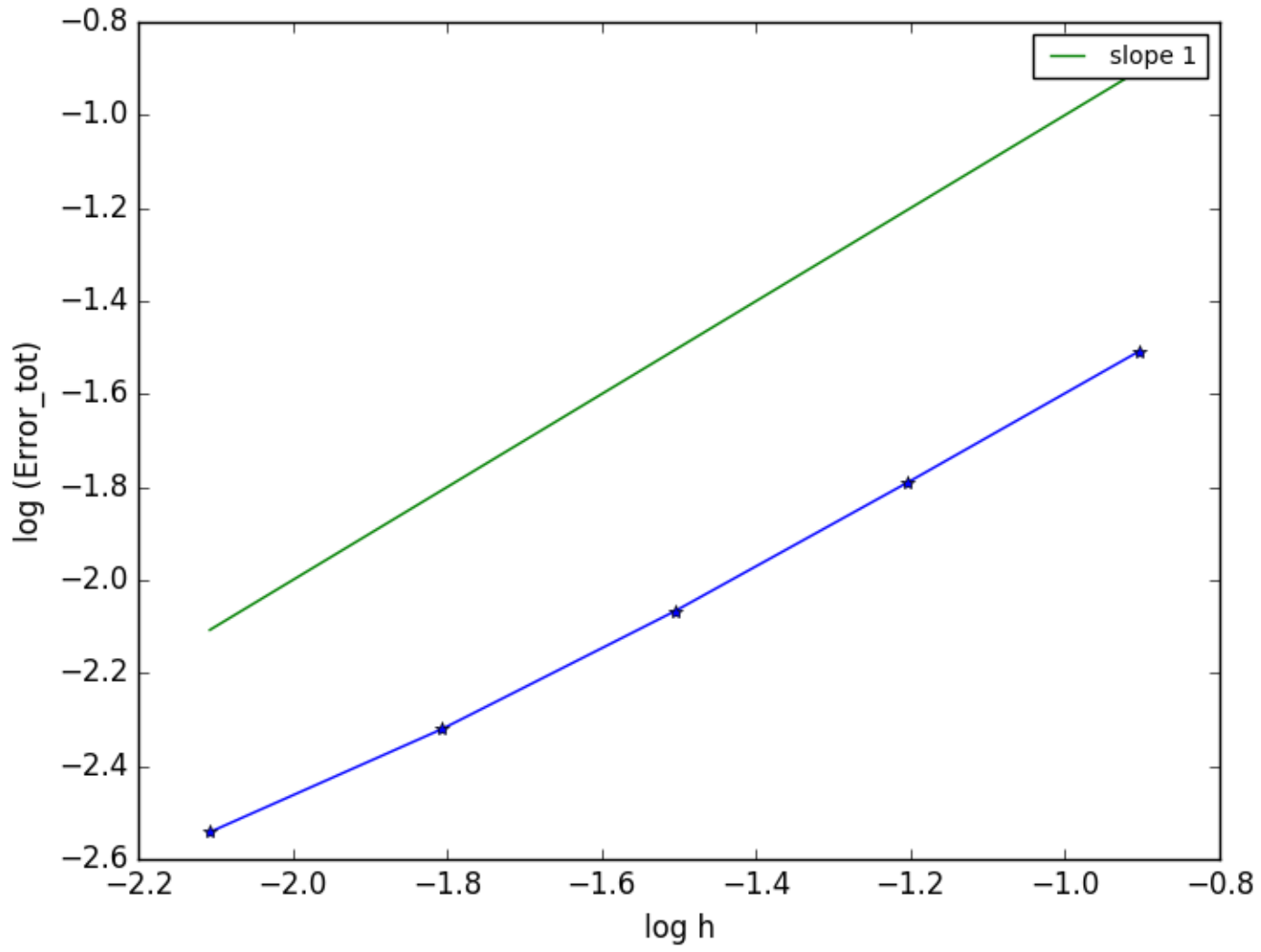

Table 3.7: Convergence rate of the Total Error with respect to the spatial mesh sizes $h_f$ and $h_s$ evaluated for decreasing values of $h_f, h_s, \tau$ at time $t = 0.064s$.

# Chapter 4

# Partitioned schemes for Problem 2

## 4.1 Introduction

In the previous chapters, we analyzed the monolithic algorithm given by problem 7 of Chapter 2. We showed the unconditional stability and the error estimate in the linearized case. Despite the monolithic algorithm express the very desirable features already mentioned, the numerical implementation requires the writing of a specific code and not allows for the reuse of existing and well-tested solver for the fluid and the solid problems. Moreover, the monolithic paradigm is very often quite expensive from the computational point of view and in some cases, the conditioning of the whole system is poorest than the conditioning of the solid and fluid problems. Given all these issues, the possibility of using existing and optimized solver for the subproblems is very attractive and researchers have proposed many possibilities to split monolithic algorithms in order to obtain partitioned schemes. These advantages are not cost free. In Chapter 1 we presented the principal difficulties encountered in designing of partitioned algorithms, here we recall that the partitioned approach requires careful formulation and implementation to avoid serious degradation in stability and accuracy, moreover gains in computational efficiency over a monolithic approach are not guaranteed. The aim of this chapter is to introduce and analyze two partitioned algorithms for the fluid-structure interaction problem that stems from a split of the solid subproblem.

In order to introduce the partitioned algorithms, we observe that, in the monolithic algorithm, the forces that the solid applies to the fluid can be divided into inertial forces, due to the acceleration of the solid mass, and elastic forces due to the solid deformation. The main idea behind the splitting technique proposed is to treat separately the two force terms. We propose, following [60], to couple the fluid problem implicitly with the solid inertial forcing term and explicitly with the elastic forcing term. In this way, the fluid problem, at each time step, is coupled with a rigid solid problem whose solution gives a guess for the solid velocity to be used in a second stage of the time-step, to advance the solid in the case of Algorithm 2, or to solve a residual solid problem in Algorithm 3. The solid velocity computed in the first stage is indicated as the fractional step solid velocity $\dot{\boldsymbol{d}}_h^{n-\frac{1}{2}}$. In the following Figure 4.1, we represent the conceptual scheme just described.

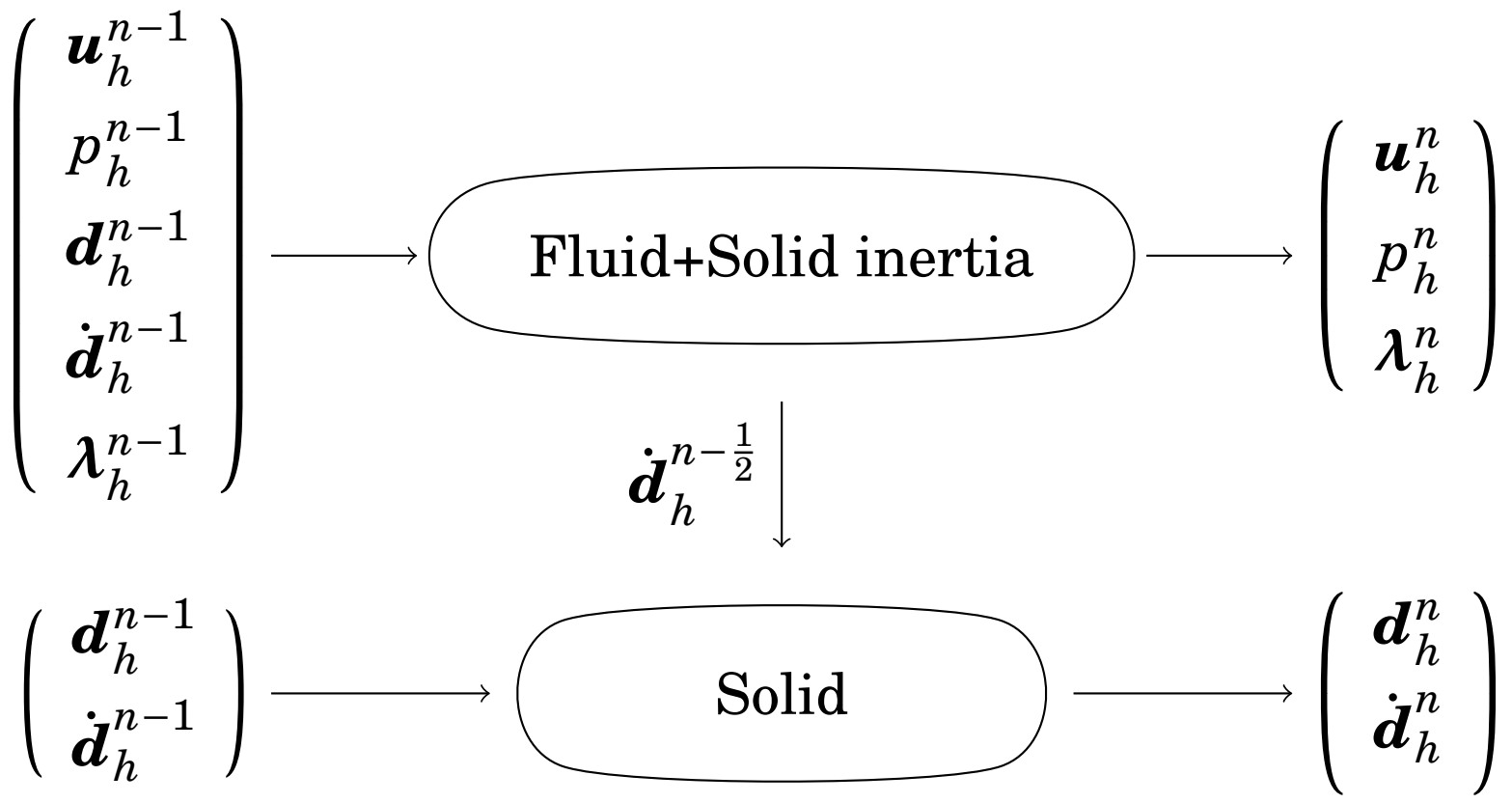

Figure 4.1: Partitioned Algorithm

This procedure has been successfully applied in fluid structure interaction problem in which the coupling conditions are imposed in weak form using the Nitsche's method [60]. In the following we propose two algorithms, both of them splits the problem in a fluid sub-problem and a solid sub-problem; in both the cases, in the fluid sub-problem we have the coupling of fluid and solid in the simplified manner described above. The difference between the two algorithms is in the solid sub-problem, in fact in Algorithm 2 we simply advance the solid position using the velocity computed in the fluid sub-problem, while in Algorithm 3 we compute the solid position and a correction of the solid velocity solving a residual problem.

For both the proposed algorithms, the explicit coupling of the fluid equations with the solid elastic forces, is realized introducing an extrapolation of the quantities of interest. We introduce the following notation to indicate the $r$th order explicit extrapolation at time $t_n$

$$x^{n\star} = \begin{cases} 0 & \text{if } r = 0, \\ x^{n-1} & \text{if } r = 1, . \\ x^{n-1} + \tau \dot{x}^{n-1} & \text{if } r = 2, \end{cases} \tag{4.1}$$

**Algorithm 2** The first algorithm proposed is reported in the following. As already specified, Algorithm 2 does not require a correction sub-step for the solid velocity computed in the fluid sub problem, then **we don't need** the introduction of the notation $\dot{\boldsymbol{d}}_h^{n-\frac{1}{2}}$. We present two version of this scheme differing for the order of extrapolation of the solid displacements $\boldsymbol{d}_h^{n\star}$. For $r = 1,2$ (see (4.1))

**Algorithm 2** Semi-implicit scheme based on [20] (not strongly coupled).

Given $\boldsymbol{\phi}_0, \boldsymbol{u}_h^0, \boldsymbol{d}_h^0, \dot{\boldsymbol{d}}_h^0$, for $n \geq 1$

1. Fluid with solid inertia step: find $(\boldsymbol{u}_h^n, p_h^n, \dot{\boldsymbol{d}}_h^n, \boldsymbol{\lambda}_h^n) \in \boldsymbol{V}_h \times Q_h \times \boldsymbol{W}_h \times \boldsymbol{\Lambda}_h$ such that

$$\rho_{\mathrm{f}}\big(\partial_\tau \boldsymbol{u}_h^n, \boldsymbol{v}_h\big) + a_h^{\mathrm{f}}\big((\boldsymbol{u}_h^n, p_h^n),(\boldsymbol{v}_h, q_h)\big) + c\big(\boldsymbol{\lambda}_h^n, \boldsymbol{v}_h \circ \boldsymbol{\phi}_h^{n-1} - \boldsymbol{w}_h\big)$$
$$- c\big(\boldsymbol{\mu}_h, \boldsymbol{u}_h^n \circ \boldsymbol{\phi}_h^{n-1} - \dot{\boldsymbol{d}}_h^n\big) + \frac{\rho_{\mathrm{s}}\epsilon}{\tau}\big(\dot{\boldsymbol{d}}_h^n - \dot{\boldsymbol{d}}_h^{n-1}, \boldsymbol{w}_h\big)_\Sigma = -a^{\mathrm{s}}\big(\boldsymbol{d}_h^{n,\star}, \boldsymbol{w}_h\big) \quad (4.2)$$

   for all $(\boldsymbol{v}_h, q_h, \boldsymbol{w}_h, \boldsymbol{\mu}_h) \in \boldsymbol{V}_h \times Q_h \times \boldsymbol{W}_h \times \boldsymbol{\Lambda}_h$.

2. Solid displacement update:
$$\boldsymbol{d}_h^n = \boldsymbol{d}_h^{n-1} + \tau \dot{\boldsymbol{d}}_h^n.$$

3. Interface update: $\boldsymbol{\phi}_h^n = \boldsymbol{\phi}_0 + \boldsymbol{d}_h^n$.

For extrapolation of order $r = 0$, that is for $\boldsymbol{d}_h^{n\star} = \boldsymbol{0}$, the elastic contribution of the solid is completely ignored, for this reason in Algorithm 2 we consider only the case with $r > 0$. Essentially Algorithm 2 is a simplified version of the monolithic algorithm given by Problem 7 and it doesn't require the solution of the solid problem. We will see in Theorem 4.3.1 that this strong simplification has a price, in fact Algorithm 2 is conditionally stable, while the monolithic algorithm is unconditionally stable.

**Algorithm 3** The second algorithm proposed is presented in the following scheme. This procedure, differently from Algorithm 2, is really a split scheme that gives a partitioned algorithm of the type of Figure 4.1. In fact, Algorithm 3 involves the solution of two staggered sub problems in each time step, namely the fluid problem that gives the fluid unknown $\boldsymbol{u}_h^n, p_h^n$, the Lagrange multiplier $\boldsymbol{\lambda}_h^n$ and a prediction of the solid velocity $\dot{\boldsymbol{d}}_h^{n-\frac{1}{2}}$. The solid residual problem, using the fractional step velocity, gives the solid displacement $\boldsymbol{d}_h^n$ and an updated solid velocity $\dot{\boldsymbol{d}}_h^n$ to be passed to the succeeding time step. The additional complexity due to the solid subproblem is able to guarantee the unconditional stability of the scheme, at least for the extrapolations of order 0 and 1. In any case Algorithm 3, with respect to the monolithic algorithm, presents the advantage to have a partitioned structure that allows for use of existing codes to solve the fluid and solid problems.

The splitting of the scheme is obtained working on the solid equation as follows. Let $\boldsymbol{d}_h^{n\star}$ be an extrapolation of the solid displacement, then we can introduce the following splitting of the solid problem: find $\dot{\boldsymbol{d}}_h^{n-\frac{1}{2}}$, $\dot{\boldsymbol{d}}_h^n$, $\boldsymbol{d}_h^n$ such that $\dot{\boldsymbol{d}}_h^n = \partial_\tau \boldsymbol{d}_h^n$ and for all $\boldsymbol{w}_h \in \boldsymbol{W}_h$

$$\frac{\rho_s \epsilon}{\tau}\left(\dot{\boldsymbol{d}}_h^{n-\frac{1}{2}} - \dot{\boldsymbol{d}}_h^{n-1}, \boldsymbol{w}_h\right)_{0,\Sigma} + a^s(\boldsymbol{d}_h^{n\star}, \boldsymbol{w}_h) = 0 \quad (4.3)$$

$$\frac{\rho_s \epsilon}{\tau}\left(\dot{\boldsymbol{d}}_h^{n} - \dot{\boldsymbol{d}}_h^{n-\frac{1}{2}}, \boldsymbol{w}_h\right)_{0,\Sigma} + a^s(\boldsymbol{d}_h^{n} - \boldsymbol{d}_h^{n\star}, \boldsymbol{w}_h) = 0. \quad (4.4)$$

We observe that the first equation does not depend on $\dot{\boldsymbol{d}}_h^n$ and $\boldsymbol{d}_h^n$ and can be solved giving $\dot{\boldsymbol{d}}_h^{n-\frac{1}{2}}$ that, used in the second equation, allows for the computation of $\dot{\boldsymbol{d}}_h^n$ and $\boldsymbol{d}_h^n$. The splitting introduced for the solid problem allows for the partitioning of the coupled problem that can be written now as

**Algorithm 3** Semi-implicit scheme (not strongly coupled).

Given $\boldsymbol{\phi}_0, \boldsymbol{u}_h^0, \boldsymbol{d}_h^0, \dot{\boldsymbol{d}}_h^0$, for $n \geq 1$

1. Fluid with solid inertia step: find $(\boldsymbol{u}_h^n, p_h^n, \dot{\boldsymbol{d}}_h^{n-\frac{1}{2}}, \boldsymbol{\lambda}_h^n) \in \boldsymbol{V}_h \times Q_h \times \boldsymbol{W}_h \times \boldsymbol{\Lambda}_h$ such that

$$
\begin{aligned}
\rho_{\mathrm{f}}\big(\partial_\tau \boldsymbol{u}_h^n, \boldsymbol{v}_h\big) + a_h^{\mathrm{f}}\big((\boldsymbol{u}_h^n, p_h^n), (\boldsymbol{v}_h, q_h)\big) + c\big(\boldsymbol{\lambda}_h^n, \boldsymbol{v}_h \circ \boldsymbol{\phi}_h^{n-1} - \boldsymbol{w}_h\big) \\
- c\big(\boldsymbol{\mu}_h, \boldsymbol{u}_h^n \circ \boldsymbol{\phi}_h^{n-1} - \dot{\boldsymbol{d}}_h^{n-\frac{1}{2}}\big) + \frac{\rho_{\mathrm{s}}\epsilon}{\tau}\big(\dot{\boldsymbol{d}}_h^{n-\frac{1}{2}} - \dot{\boldsymbol{d}}_h^{n-1}, \boldsymbol{w}_h\big)_\Sigma = -a_{\mathrm{s}}\big(\boldsymbol{d}_h^{n,\star}, \boldsymbol{w}_h\big)
\end{aligned}
\tag{4.5}
$$

for all $(\boldsymbol{v}_h, q_h, \boldsymbol{w}_h, \boldsymbol{\mu}_h) \in \boldsymbol{V}_h \times Q_h \times \boldsymbol{W}_h \times \boldsymbol{\Lambda}_h$.

2. Solid update: find $(\boldsymbol{d}_h^n, \dot{\boldsymbol{d}}_h^n,) \in \boldsymbol{W}_h \times \boldsymbol{W}_h$ such that

$$
\begin{cases}
\partial_\tau \boldsymbol{d}_h^n = \dot{\boldsymbol{d}}_h^n, \\
\dfrac{\rho_s \epsilon}{\tau}\left(\dot{\boldsymbol{d}}_h^n - \dot{\boldsymbol{d}}_h^{n-\frac{1}{2}}, \boldsymbol{z}_h\right) + a^s(\boldsymbol{d}_h^n - \boldsymbol{d}_h^{n\star}, \boldsymbol{z}_h) = 0,
\end{cases}
\tag{4.6}
$$

for all $\boldsymbol{z}_h \in \boldsymbol{W}_h$.

3. Interface update: $\boldsymbol{\phi}_h^n = \boldsymbol{\phi}_0 + \boldsymbol{d}_h^n$.

## 4.2 Well-posedness of Algorithms 2 and 3

In this section we study the well-posedness of the time-step problems associated to the proposed algorithms. In this context the same type of analysis carried out for the monolithic algorithm applies. In particular we consider the well-posedness of the time-step problems.

- **Well-posedness of Algorithm 2.** We introduce the matrix form of Algorithm 2. Given the initial data $\boldsymbol{\phi}_0, \boldsymbol{u}_h^0, \boldsymbol{d}_h^0, \dot{\boldsymbol{d}}_h^0$, for $n \geq 1$

  1. find $(\boldsymbol{u}_h^n, \dot{\boldsymbol{d}}_h^n, \boldsymbol{\lambda}_h^n, p_h^n) \in \boldsymbol{V}_h \times \boldsymbol{W}_h \times \boldsymbol{\Lambda}_h \times Q_h$ such that

$$
\left(\begin{array}{ccc|c}
A_{fh} & 0 & (L_{fh}^n)^\top & B_h^\top \\
0 & A_{sh} & -L_{sh}^\top & 0 \\
L_{fh}^n & -L_{sh} & 0 & 0 \\
\hline
B_h & 0 & 0 & S_h
\end{array}\right)
\left(\begin{array}{c}
\boldsymbol{u}_h^n \\
\dot{\boldsymbol{d}}_h^n \\
\boldsymbol{\lambda}_h^n \\
\hline
p_h^n
\end{array}\right)
=
\left(\begin{array}{c}
\boldsymbol{f}^{n-1} \\
\boldsymbol{g}^{n-1} \\
\mathbf{m}^{n-1} \\
\hline
0
\end{array}\right),
\tag{4.7}
$$

2. Solid displacement update:
$$\boldsymbol{d}_h^n = \boldsymbol{d}_h^{n-1} + \tau \dot{\boldsymbol{d}}_h^n.$$
3. Interface update: $\boldsymbol{\phi}_h^n = \boldsymbol{\phi}_0 + \boldsymbol{d}_h^n$.

Denoting by $\varphi, \psi, \chi$ and $\zeta$ the basis functions respectively in $\boldsymbol{V}_h, Q_h, \boldsymbol{W}_h, \boldsymbol{\Lambda}_h$, the submatrices in 4.7 have the following expressions

$$
\begin{aligned}
&A_{fh} = \frac{\rho_{\mathrm{f}}}{\tau} M_{\mathrm{f}} + K_{\mathrm{f}}, \quad \text{with} \quad (M_{\mathrm{f}})_{ij} = \left(\varphi_j, \varphi_i\right)_{0,\Omega}, \quad (K_{\mathrm{f}})_{ij} = \left(\boldsymbol{\epsilon}\varphi_j, \boldsymbol{\epsilon}\varphi_i\right)_{0,\Omega}, \\
&(B_h)_{ij} = -\left(\operatorname{div}\varphi_j, \psi_i\right)_{0,\Omega}, \\
&(S_h)_{ij} = \gamma \sum_{K\in\mathcal{T}_h} h_K^2 \left(\nabla\psi_j, \nabla\psi_i\right)_{0,K}, \\
&A_{sh} = \frac{\rho_{\mathrm{s}}\epsilon}{\tau} M_{\mathrm{s}}, \quad \text{with} \quad (M_{\mathrm{s}})_{ij} = \left(\chi_j, \chi_i\right)_{0,\Sigma}, \\
&(L_{fh}^n)_{ij} = c(\boldsymbol{\zeta}_j, \boldsymbol{\varphi}_i \circ \boldsymbol{\phi}^{n-1}) = (\boldsymbol{\zeta}_j, \boldsymbol{\varphi}_i \circ \boldsymbol{\phi}^{n-1})_{0,\Omega}, \\
&(L_{sh})_{ij} = c(\boldsymbol{\zeta}_j, \boldsymbol{\chi}_i) = (\boldsymbol{\zeta}_j, \boldsymbol{\chi}_i)_{0,\Sigma}, \\
&\boldsymbol{f}^{n-1} = \frac{\rho_{\mathrm{f}}}{\tau} M_{\mathrm{f}} \boldsymbol{u}_h^{n-1}, \\
&\boldsymbol{g}^{n-1} = \frac{\rho_{\mathrm{s}}\epsilon}{\tau} M_{\mathrm{s}} \dot{\boldsymbol{d}}_h^{n-1} - K_{\mathrm{s}} \boldsymbol{d}_h^{n\star}, \qquad \text{with} \qquad (K_{\mathrm{s}})_{ij} = a^s\left(\chi_j, \chi_i\right), \\
&\boldsymbol{m}^{n-1} = -\frac{1}{\tau} L_{\mathrm{s}} \boldsymbol{d}_h^{n-1}.
\end{aligned}
\tag{4.8}
$$

For any fixed $n \geq 1$ it is easy to recognize in (4.7) the same structure of the steady saddle point problem that was analyzed to show the well-posedness for the monolithic algorithm, the only difference being the solid matrix $A_{sh}$. Following the same steps of Proposition 2.4.3 it is possible to prove the well-posedness on the function space

$$\mathbb{V}_h = \boldsymbol{V}_h \times \boldsymbol{W}_h \times \boldsymbol{\Lambda}_h \times Q_h$$

with norm $\|\boldsymbol{U}_h\|_{\mathbb{V}_h} = \left(\|\boldsymbol{u}_h\|_{1,\Omega}^2 + \|\dot{\boldsymbol{d}}_h\|_{0,\Sigma}^2 + \|\boldsymbol{\lambda}_h\|_{\Lambda}^2 + \|p_h\|_{0,\Omega}^2\right)^{\frac{1}{2}}$. Also in this case the key result to prove the well-posedness is (2.40) that remains valid in this new context.

- **Well-posedness of Algorithm 3.** The matrix form of Algorithm 3 suitable for the well-posedness analysis is reported in the following. Given the initial data $\boldsymbol{\phi}_0, \boldsymbol{u}_h^0$, $\boldsymbol{d}_h^0$, $\dot{\boldsymbol{d}}_h^0$, for $n \geq 1$

  1. find $(\boldsymbol{u}_h^n, \dot{\boldsymbol{d}}_h^{n-\frac{1}{2}}, \boldsymbol{\lambda}_h^n, p_h^n) \in \boldsymbol{V}_h \times \boldsymbol{W}_h \times \boldsymbol{\Lambda}_h \times Q_h$ such that

$$
\left(\begin{array}{ccc|c}
A_{fh} & 0 & (L_{fh}^n)^\top & B_h^\top \\
0 & A_{sh} & -L_{sh}^\top & 0 \\
L_{fh}^n & -L_{sh} & 0 & 0 \\
\hline
B_h & 0 & 0 & S_h
\end{array}\right)
\left(\begin{array}{c}
\boldsymbol{u}_h^n \\
\dot{\boldsymbol{d}}_h^{n-\frac{1}{2}} \\
\boldsymbol{\lambda}_h^n \\
\hline
p_h^n
\end{array}\right)
=
\left(\begin{array}{c}
\boldsymbol{f}^{n-1} \\
\boldsymbol{g}^{n-1} \\
\mathbf{m}^{n-1} \\
\hline
0
\end{array}\right),
\tag{4.9}
$$

2. Solid displacement update:

$$\begin{aligned}
&\boldsymbol{d}_h^n = \boldsymbol{d}_h^{n-1} + \tau \dot{\boldsymbol{d}}_h^n, \\
&\frac{\rho_s \epsilon}{\tau} M_s \dot{\boldsymbol{d}}_h^n + K_s \boldsymbol{d}_h^n = \mathbf{r}^{n-1}
\end{aligned} \tag{4.10}$$

3. Interface update: $\boldsymbol{\phi}_h^n = \boldsymbol{\phi}_0 + \boldsymbol{d}_h^n$.

The sub-matrices in (4.9) are the same used in Algorithm 2, and are defined in (4.8). Instead, the stiffness matrix introduced in the second step is defined as

$$(K_{\mathrm{s}})_{ij} = a^s\left(\chi_j, \chi_i\right)$$

where $\chi_j$ are the basis functions for $\boldsymbol{W}_h$. The right hand side of the second substep has the expression

$$\boldsymbol{r}^{n-1} = \frac{\rho_s \epsilon}{\tau} M_s \dot{\boldsymbol{d}}_h^{n-1} - K_s \boldsymbol{d}_h^{n\star}$$

. Then the first step (4.9) is well-posed on

$$\mathbb{V}_h = \boldsymbol{V}_h \times \boldsymbol{W}_h \times \boldsymbol{\Lambda}_h \times \boldsymbol{Q}_h$$

with norm $\|\boldsymbol{U}_h\|_{\mathbb{V}_h} = \left(\|\boldsymbol{u}_h\|_{1,\Omega}^2 + \|\dot{\boldsymbol{d}}_h\|_{0,\Sigma}^2 + \|\boldsymbol{\lambda}_h\|_{\Lambda}^2 + \|p_h\|_{0,\Omega}^2\right)^{\frac{1}{2}}$. This result can be obtained in the same way of the similar result for Algorithm 2. The second step of the algorithm is a pure elastic problem whose well-posedness is well established owing to the coercivity of $a^s(\cdot,\cdot)$.

## 4.3 Stability analysis of Algorithms 2 and 3

In this section we analyze the stability property of the proposed Algorithms. In the following we shall apply the well-known discrete Gronwall Lemma 3.3.1 already used in chapter 3 and we will use the **discrete energy** of the coupled system at time $t_n$, defined as

$$E_h^n := \rho_{\mathrm{f}} \|\boldsymbol{u}_h^n\|_{0,\Omega}^2 + \rho_{\mathrm{s}} \epsilon \|\dot{\boldsymbol{d}}_h^n\|_{0,\Sigma}^2 + \|\boldsymbol{d}_h^n\|_{\mathrm{s}}^2, \quad \text{for} \quad n \geq 0, \tag{4.11}$$

and the **discrete dissipation** of the coupled system at time $t_n$ defined as

$$\begin{aligned}
\tilde{D}_h^n := &\tau \sum_{k=1}^n \left( \|\epsilon(\boldsymbol{u}_h^k)\|_{0,\Omega}^2 + |p_h^k|_{s_h}^2 \right) \\
&+ \tau^2 \sum_{k=1}^n \left( \rho_{\mathrm{f}} \|\partial_\tau \boldsymbol{u}_h^k\|_{0,\Omega}^2 + \rho_{\mathrm{s}} \epsilon \|\partial_\tau \dot{\boldsymbol{d}}_h^k\|_{0,\Sigma_0}^2 \right) \quad \text{for} \quad n \geq 1.
\end{aligned} \tag{4.12}$$

We are now ready to present the stability results.

**Algorithm 2** The stability analysis of Algorithm 2 is presented in the following theorem

**Theorem 4.3.1.** *Energy Stability for Algorithm 2*
*With reference to* (4.11) *and* (4.12)*, let* $(\boldsymbol{u}_h^0, \dot{\boldsymbol{d}}_h^0, \boldsymbol{d}_h^0) \in \boldsymbol{V}_h \times \boldsymbol{W}_h \times \boldsymbol{W}_h$ *be initial values for the fluid and solid velocities and for solid displacement; for* $n \geq 1$*, let* $(\boldsymbol{u}_h^n, p_h^n, \dot{\boldsymbol{d}}_h^n, \boldsymbol{d}_h^n, \boldsymbol{\lambda}_h^n) \in \boldsymbol{V}_h \times Q_h \times \boldsymbol{W}_h \times \boldsymbol{W}_h \times \boldsymbol{\Lambda}_h$ *be the solution of Algorithm 2 with* $r = 1$ *or* $r = 2$ *, let* $\tau$ *and* $h_s$ *be such that there exists* $\alpha > 0$ *such that*

$$\frac{\tau \beta_{\mathrm{s}} C_I^2}{h_s^2} \leq \alpha, \qquad \tau \alpha < \rho_{\mathrm{s}} \epsilon, \tag{4.13}$$

*then for all* $n \geq 1$

$$E_h^n + \tilde{D}_h^n \lesssim exp\left(\tau \sum_{m=1}^{n} \frac{\alpha}{\rho_{\mathrm{s}}\epsilon - \tau\alpha}\right) E_h^0.$$

*Proof.* Let us consider the fully discrete problem of Algorithm 2, and let us take as test functions in (4.2)

$$(\boldsymbol{v}_h, q_h) = \tau(\boldsymbol{u}_h^n, p_h^n), \qquad \boldsymbol{\mu}_h = \tau \boldsymbol{\lambda}_h^n, \qquad \boldsymbol{w}_h = \tau \dot{\boldsymbol{d}}_h^n,$$

then, we obtain

$$\begin{aligned} &\rho_{\mathrm{f}}\left(\boldsymbol{u}_h^n - \boldsymbol{u}_h^{n-1}, \boldsymbol{u}_h^n\right)_{0,\Omega} + 2\mu\tau\|\boldsymbol{\epsilon}(\boldsymbol{u}_h^n)\|_{0,\Omega}^2 + \tau|p_h^n|_{s_h}^2 \\ &+ \rho_{\mathrm{s}}\epsilon\left(\dot{\boldsymbol{d}}_h^n - \dot{\boldsymbol{d}}_h^{n-1}, \dot{\boldsymbol{d}}_h^n\right)_{0,\Sigma} + \tau a_{\mathrm{s}}\left(\boldsymbol{d}_h^{n\star}, \dot{\boldsymbol{d}}_h^n\right) = 0, \end{aligned}$$

hence, after some computation,

$$\begin{aligned} &\frac{\rho_{\mathrm{f}}}{2}\left(\|\boldsymbol{u}_h^n\|_{0,\Omega}^2 - \|\boldsymbol{u}_h^{n-1}\|_{0,\Omega}^2 + \|\boldsymbol{u}_h^n - \boldsymbol{u}_h^{n-1}\|_{0,\Omega}^2\right) + 2\mu\tau\|\boldsymbol{\epsilon}(\boldsymbol{u}_h^n)\|_{0,\Omega}^2 + \tau|p_h^n|_{s_h}^2 \\ &+ \frac{\rho_{\mathrm{s}}\epsilon}{2}\left(\|\dot{\boldsymbol{d}}_h^n\|_{0,\Sigma}^2 - \|\dot{\boldsymbol{d}}_h^{n-1}\|_{0,\Sigma}^2 + \|\dot{\boldsymbol{d}}_h^n - \dot{\boldsymbol{d}}_h^{n-1}\|_{0,\Sigma}^2\right) + \underbrace{\tau a_{\mathrm{s}}\left(\boldsymbol{d}_h^{n\star}, \dot{\boldsymbol{d}}_h^n\right)}_{T_0} = 0, \end{aligned} \tag{4.14}$$

where we have used the following relation valid for all $\boldsymbol{a}, \boldsymbol{b} \in L^2(G)^d$

$$2(\boldsymbol{a}, \boldsymbol{b})_{0,G} = \|\boldsymbol{a}\|_{0,G}^2 + \|\boldsymbol{b}\|_{0,G}^2 - \|\boldsymbol{a} - \boldsymbol{b}\|_{0,G}^2. \tag{4.15}$$

Equation (4.14) represents the discrete form of the energy equation (2.2) for the continuum problem 2; in the following we derive expressions for the term $T_0$ which depends on the extrapolation order.

- **Scheme with** $r = 1$, $\boldsymbol{d}_h^{n\star} = \boldsymbol{d}_h^{n-1}$

Using the inverse inequality (3.30) we obtain the following estimate from below for the term $T_0$,

$$\begin{aligned}
T_0 &= \tau a_{\mathrm{s}}(\boldsymbol{d}_h^{n-1}, \dot{\boldsymbol{d}}_h^n) = a_{\mathrm{s}}(\boldsymbol{d}_h^{n-1}, \boldsymbol{d}_h^n - \boldsymbol{d}_h^{n-1}) \\
&= \frac{1}{2}\|\boldsymbol{d}_h^n\|_{\mathrm{s}}^2 - \frac{1}{2}\|\boldsymbol{d}_h^{n-1}\|_{\mathrm{s}}^2 - \frac{1}{2}\|\boldsymbol{d}_h^n - \boldsymbol{d}_h^{n-1}\|_{\mathrm{s}}^2 \\
&= \frac{1}{2}\left(\|\boldsymbol{d}_h^n\|_{\mathrm{s}}^2 - \|\boldsymbol{d}_h^{n-1}\|_{\mathrm{s}}^2\right) - \frac{\tau^2}{2}\|\dot{\boldsymbol{d}}_h^n\|_{\mathrm{s}}^2 \geq \frac{1}{2}\left(\|\boldsymbol{d}_h^n\|_{\mathrm{s}}^2 - \|\boldsymbol{d}_h^{n-1}\|_{\mathrm{s}}^2\right) - \frac{\tau^2\beta_{\mathrm{s}}C_I^2}{2h_s^2}\|\dot{\boldsymbol{d}}_h^n\|_{0,\Sigma}^2 \\
&\geq \frac{1}{2}\left(\|\boldsymbol{d}_h^n\|_{\mathrm{s}}^2 - \|\boldsymbol{d}_h^{n-1}\|_{\mathrm{s}}^2\right) - \frac{\tau\alpha}{2}\|\dot{\boldsymbol{d}}_h^n\|_{0,\Sigma}^2,
\end{aligned}$$

where $\alpha$ is the constant given in (4.13) for fixed $\tau$ and $h_s$.

Using this bound for $T_0$ in the energy equation (4.14), and summing over $m = 1, \ldots, n$, we obtain

$$\begin{aligned}
&\frac{\rho_{\mathrm{f}}}{2}\|\boldsymbol{u}_h^n\|_{0,\Omega}^2 + \frac{1}{2}\|\boldsymbol{d}_h^n\|_{\mathrm{s}}^2 + \underbrace{\frac{\rho_{\mathrm{s}}\epsilon}{2}\|\dot{\boldsymbol{d}}_h^n\|_{0,\Sigma}^2}_{a_n} \\
&+ \tau\sum_{m=1}^{n}\underbrace{\left[2\mu\|\boldsymbol{\epsilon}(\boldsymbol{u}_h^m)\|_{0,\Omega}^2 + |p_h^m|_{s_h}^2 + \tau\left(\frac{\rho_{\mathrm{f}}}{2}\|\partial_\tau\boldsymbol{u}_h^m\|_{0,\Omega}^2 + \frac{\rho_{\mathrm{s}}\epsilon}{2}\|\partial_\tau\dot{\boldsymbol{d}}_h^m\|_{0,\Sigma}^2\right)\right]}_{b_m} \\
&\leq \tau\sum_{m=1}^{n}\underbrace{\frac{\alpha}{\rho_{\mathrm{s}}\epsilon}}_{\gamma_m}\underbrace{\frac{\rho_{\mathrm{s}}\epsilon}{2}\|\dot{\boldsymbol{d}}_h^m\|_{0,\Sigma}^2}_{a_m} + \underbrace{\left(\frac{\rho_{\mathrm{f}}}{2}\|\boldsymbol{u}_h^0\|_{0,\Omega}^2 + \frac{1}{2}\|\boldsymbol{d}_h^0\|_{\mathrm{s}}^2 + \left(\frac{\rho_{\mathrm{s}}\epsilon}{2} - \tau\alpha\right)\|\dot{\boldsymbol{d}}_h^0\|_{0,\Sigma}^2\right)}_{B=E_h^0}.
\end{aligned}$$

Since, by hypothesis, we have $\tau\alpha < \rho_{\mathrm{s}}\epsilon$, we can apply the discrete Gronwall Lemma 3.3.1 obtaining

$$E_h^n + \tilde{D}_h^n \lesssim exp\left(\tau\sum_{m=1}^{n}\frac{\alpha}{\rho_{\mathrm{s}}\epsilon - \tau\alpha}\right)E_h^0.$$

- **Scheme with** $r = 2$, $\boldsymbol{d}_h^{n\star} = \boldsymbol{d}_h^{n-1} + \tau\dot{\boldsymbol{d}}_h^{n-1}$

The term $T_0$ is controlled as follows,

$$T_0 = \tau a_{\mathrm{s}}\left((\boldsymbol{d}_h^{n-1} + \tau\dot{\boldsymbol{d}}_h^{n-1}), \dot{\boldsymbol{d}}_h^n\right) = \underbrace{a_{\mathrm{s}}\left(\boldsymbol{d}_h^{n-1}, \boldsymbol{d}_h^n - \boldsymbol{d}_h^{n-1}\right)}_{A} + \underbrace{\tau^2 a_{\mathrm{s}}\left(\dot{\boldsymbol{d}}_h^{n-1}, \dot{\boldsymbol{d}}_h^n\right)}_{B}.$$

The term A is bounded using the inequality obtained in the case $r = 1$, with $\alpha$ being the constant given in (4.13) for fixed $\tau$ and $h_s$.

$$\begin{aligned}
A = a_{\mathrm{s}}\left(\boldsymbol{d}_h^{n-1}, \boldsymbol{d}_h^n - \boldsymbol{d}_h^{n-1}\right) &= \frac{1}{2}\left(\|\boldsymbol{d}_h^n\|_{\mathrm{s}}^2 - \|\boldsymbol{d}_h^{n-1}\|_{\mathrm{s}}^2\right) - \frac{1}{2}\|\boldsymbol{d}_h^n - \boldsymbol{d}_h^{n-1}\|_{\mathrm{s}}^2 \\
&\geq \frac{1}{2}\left(\|\boldsymbol{d}_h^n\|_{\mathrm{s}}^2 - \|\boldsymbol{d}_h^{n-1}\|_{\mathrm{s}}^2\right) - \frac{\tau\alpha}{2}\|\dot{\boldsymbol{d}}_h^n\|_{0,\Sigma}^2.
\end{aligned}$$

The estimate of term B is given in the following

$$
\begin{aligned}
B = \tau^2 a_s\left(\dot{\boldsymbol{d}}_h^{n-1}, \dot{\boldsymbol{d}}_h^n\right) &= \frac{\tau^2}{2}\|\dot{\boldsymbol{d}}_h^n\|_s^2 + \frac{\tau^2}{2}\|\dot{\boldsymbol{d}}_h^{n-1}\|_s^2 - \frac{\tau^2}{2}\|\dot{\boldsymbol{d}}_h^n - \dot{\boldsymbol{d}}_h^{n-1}\|_s^2 \\
&\geq -\frac{\tau^2}{2}\|\dot{\boldsymbol{d}}_h^n - \dot{\boldsymbol{d}}_h^{n-1}\|_s^2 \\
&\geq -\frac{\tau^2\beta_s C_I^2}{2h_s^2}\|\dot{\boldsymbol{d}}_h^n - \dot{\boldsymbol{d}}_h^{n-1}\|_{0,\Sigma}^2 \geq -\frac{\tau\alpha}{2}\|\dot{\boldsymbol{d}}_h^n - \dot{\boldsymbol{d}}_h^{n-1}\|_{0,\Sigma}^2,
\end{aligned}
$$

where $\alpha$ is the constant given in (4.13) for fixed $\tau$ and $h_s$. Using the estimates of A and B we have

$$
T_0 \geq \frac{1}{2}\left(\|\boldsymbol{d}_h^n\|_s^2 - \|\boldsymbol{d}_h^{n-1}\|_s^2\right) - \frac{\tau\alpha}{2}\|\dot{\boldsymbol{d}}_h^n\|_{0,\Sigma}^2 - \frac{\tau\alpha}{2}\|\dot{\boldsymbol{d}}_h^n - \dot{\boldsymbol{d}}_h^{n-1}\|_{0,\Sigma}^2. \tag{4.16}
$$

Using this bound for $T_0$ in the energy equation (4.14), and summing over $m = 1,\dots,n$, we obtain

$$
\begin{aligned}
&\underbrace{\frac{\rho_f}{2}\|\boldsymbol{u}_h^n\|_{0,\Omega}^2 + \frac{1}{2}\|\boldsymbol{d}_h^n\|_s^2 + \frac{\rho_s\epsilon}{2}\|\dot{\boldsymbol{d}}_h^n\|_{0,\Sigma}^2}_{a_n} \\
&+\tau\sum_{m=1}^{n}\underbrace{\left[2\mu\|\boldsymbol{\epsilon}(\boldsymbol{u}_h^m)\|_{0,\Omega}^2 + |p_h^m|_{s_h}^2 + \tau\left(\frac{\rho_f}{2}\|\partial_\tau \boldsymbol{u}_h^m\|_{0,\Omega}^2 + \left(\frac{\rho_s\epsilon}{2} - \tau\alpha\right)\|\partial_\tau\dot{\boldsymbol{d}}_h^m\|_{0,\Sigma}^2\right)\right]}_{b_m} \\
&\leq \tau\sum_{m=1}^{n}\underbrace{\frac{\alpha}{\rho_s\epsilon}}_{\gamma_m}\underbrace{\frac{\rho_s\epsilon}{2}\|\dot{\boldsymbol{d}}_h^m\|_{0,\Sigma}^2}_{a_m} + \underbrace{\left(\frac{\rho_f}{2}\|\boldsymbol{u}_h^0\|_{0,\Omega}^2 + \frac{1}{2}\|\boldsymbol{d}_h^0\|_s^2 + \frac{\rho_s\epsilon}{2}\|\dot{\boldsymbol{d}}_h^0\|_{0,\Sigma}^2\right)}_{B=E_h^0}.
\end{aligned}
$$

Since, by hypothesis, we have $\tau\alpha < \rho_s\epsilon$, we can apply the discrete Gronwall Lemma 3.3.1 obtaining

$$
E_h^n + \tilde{D}_h^n \lesssim exp\left(\tau\sum_{m=1}^{n}\frac{\alpha}{\rho_s\epsilon - \tau\alpha}\right)E_h^0.
$$

□

Theorem 4.3.1, shows conditional stability of Algorithm 2 for all extrapolations analyzed; this feature in not surprising because of the strong simplification introduced with respect to the monolithic algorithm. In fact, in Algorithm 2 we never solve the elastic problem associated to the solid, in contrast to Problem 7.

**Algorithm 3** The stability analysis of Algorithm 3 is based on the following result that gives a relation between the fractional step solid velocity $\dot{\boldsymbol{d}}_h^{n-\frac{1}{2}}$ and the solid velocity $\dot{\boldsymbol{d}}_h^n$, the solid displacement $\boldsymbol{d}_h^n$ and the extrapolation $\boldsymbol{d}_h^{n\star}$. Moreover, Lemma 4.3.2 states that the splitting scheme can be interpreted as a kinematic perturbation of Problem 7.

**Lemma 4.3.2.** *For $n \geq 1$, let $(\boldsymbol{u}_h^n, p_h^n, \dot{\boldsymbol{d}}_h^{n-\frac{1}{2}}, \dot{\boldsymbol{d}}_h^n, \boldsymbol{d}_h^n, \boldsymbol{\lambda}_h^n) \in \boldsymbol{V}_h \times Q_h \times \boldsymbol{W}_h \times \boldsymbol{W}_h \times \boldsymbol{\Lambda}_h$ be the solution of Algorithm 3. Then, for $n \geq 1$, we have:*

- *Step 2 of Algorithm 3 gives the following expression of the fractional step velocity*

$$\dot{\boldsymbol{d}}_h^{n-\frac{1}{2}} = \dot{\boldsymbol{d}}_h^n + \frac{\tau}{\rho_s \epsilon} \boldsymbol{L}_h^s \left( \boldsymbol{d}_h^n - \boldsymbol{d}_h^{n\star} \right), \tag{4.17}$$

*being* $\boldsymbol{L}_h^s$ *the discrete elastic operator defined in* (3.25).

- *Algorithm 3 implies the following expression:*

$$\begin{aligned} &\forall (\boldsymbol{v}_h, q_h, \boldsymbol{w}_h, \boldsymbol{\mu}_h) \in \boldsymbol{V}_h \times Q_h \times \boldsymbol{W}_h \times \boldsymbol{\Lambda}_h \\ &\qquad \rho_f \left( \partial_\tau \boldsymbol{u}_h^n, \boldsymbol{v}_h \right) + a^f((\boldsymbol{u}_n^n, p_h^n); (\boldsymbol{v}_h, q_h)) + s_h(p_h^n, q_h) + c\left( \boldsymbol{\lambda}_h^n, \boldsymbol{v}_h \circ \boldsymbol{\phi}_h^{n-1} - \boldsymbol{w}_h \right) \\ &\qquad\qquad - c\left( \boldsymbol{\mu}_h, \boldsymbol{u}_h^n \circ \boldsymbol{\phi}_h^{n-1} - \dot{\boldsymbol{d}}_h^{n-\frac{1}{2}} \right) + \rho_s \epsilon \left( \partial_\tau \dot{\boldsymbol{d}}_h^n, \boldsymbol{w}_h \right) + a_s(\boldsymbol{d}_h^n, \boldsymbol{w}_h) = 0. \end{aligned} \tag{4.18}$$

$$\partial_\tau \boldsymbol{d}_h^n = \dot{\boldsymbol{d}}_h^n \tag{4.19}$$

*Proof.* The first claim is consequence of equation (4.6) of Algorithm 3 which, using the discrete elastic operator (3.25), gives

$$\dot{\boldsymbol{d}}_h^{n-\frac{1}{2}} = \dot{\boldsymbol{d}}_h^n + \frac{\tau}{\rho_s \epsilon} \boldsymbol{L}_h^s \left( \boldsymbol{d}_h^n - \boldsymbol{d}_h^{n\star} \right). \tag{4.20}$$

The second claim follows observing that the expression given in (4.18) is obtained summing equations (4.5) and (4.6) tested with $\boldsymbol{z}_h = \boldsymbol{w}_h \in \boldsymbol{W}_h$. □

**Remark 2.** *Lemma 4.3.2 leads to two remarkable consequences. The relation 4.18 indicates that Algorithm 3 treats the dynamic coupling condition* $(1.9)_1$ *as does Problem 7, in fact the term* $c\left( \boldsymbol{\lambda}_h^n, \boldsymbol{v}_h \circ \boldsymbol{\phi}_h^{n-1} - \boldsymbol{w}_h \right)$ *is present in both the algorithms in the same form. The difference between the schemes lies on the treatment of the kinematic coupling condition* $(1.9)_2$. *Indeed, Problem 7 imposes that the* $L^2$ *projection of the trace of the fluid velocity on the solid ( that is* $\boldsymbol{u}_h^n \circ \boldsymbol{\phi}^{n-1}$*) onto* $\boldsymbol{W}_h$ *is equal to the solid velocity, namely*

$$c\left( \boldsymbol{w}_h, \boldsymbol{u}_h^n \circ \boldsymbol{\phi}^{n-1} - \dot{\boldsymbol{d}}_h^n \right) = \left( \boldsymbol{w}_h, \boldsymbol{u}_h^n \circ \boldsymbol{\phi}^{n-1} - \dot{\boldsymbol{d}}_h^n \right)_{0,\Sigma} = 0 \qquad \forall \boldsymbol{w}_h \in \boldsymbol{W}_h \subset \boldsymbol{\Lambda}_h,$$

*whereas Algorithm 3 enforces on the interface the following condition*

$$c\left( \boldsymbol{w}_h, \boldsymbol{u}_h^n \circ \boldsymbol{\phi}^{n-1} - \dot{\boldsymbol{d}}_h^{n-\frac{1}{2}} \right) = \left( \boldsymbol{w}_h, \boldsymbol{u}_h^n \circ \boldsymbol{\phi}^{n-1} - \dot{\boldsymbol{d}}_h^{n-\frac{1}{2}} \right)_{0,\Sigma} = 0 \qquad \forall \boldsymbol{w}_h \in \boldsymbol{W}_h \subset \boldsymbol{\Lambda}_h,$$

*then it follows, considering expression* (4.17)

$$\left( \boldsymbol{u}_h^n \circ \boldsymbol{\phi}^{n-1}, \boldsymbol{w}_h \right)_{0,\Sigma} = \left( \dot{\boldsymbol{d}}_h^{n-\frac{1}{2}}, \boldsymbol{w}_h \right)_{0,\Sigma} = \left( \dot{\boldsymbol{d}}_h^n, \boldsymbol{w}_h \right)_{0,\Sigma} + \frac{\tau}{\rho_s \epsilon} \left( \boldsymbol{L}_h^s \left( \boldsymbol{d}_h^n - \boldsymbol{d}_h^{n\star} \right), \boldsymbol{w}_h \right)_{0,\Sigma}, \qquad \forall \boldsymbol{w}_h \in \boldsymbol{W}_h \subset \boldsymbol{\Lambda}_h.$$

*In other words, Algorithm 3 is nothing but a kinematic perturbation of Problem 7, with the perturbation being given by* $\frac{\tau}{\rho_s \epsilon} \boldsymbol{L}_h^s \left( \boldsymbol{d}_h^n - \boldsymbol{d}_h^{n\star} \right)$. *Note that this term depends on the extrapolation* $\boldsymbol{d}_h^{n\star}$ *and on the time step* $\tau$.

In the following stability analysis it is useful to define a new expression of the dissipation at time $t_n$ as

$$\begin{aligned} D_h^n :=& \tau \sum_{k=1}^{n} \left( \|\epsilon(\boldsymbol{u}_h^k)\|_{0,\Omega}^2 + |p_h^k|_{s_h}^2 \right) \\ &+ \tau^2 \sum_{k=1}^{n} \left( \rho_{\mathrm{f}} \|\partial_\tau \boldsymbol{u}_h^k\|_{0,\Omega}^2 + \rho_{\mathrm{s}} \epsilon \|\partial_\tau \dot{\boldsymbol{d}}_h^k\|_{0,\Sigma_0}^2 + \|\partial_\tau \boldsymbol{d}_h^k\|_{\mathrm{s}}^2 \right) \quad \text{for} \quad n \geq 1. \end{aligned} \tag{4.21}$$

**Theorem 4.3.3.** *Energy Stability for Algorithm 3*
*With reference to* (4.11) *and* (4.21)*, let* $(\boldsymbol{u}_h^0, \dot{\boldsymbol{d}}_h^0, \boldsymbol{d}_h^0) \in \boldsymbol{V}_h \times \boldsymbol{W}_h \times \boldsymbol{W}_h$ *be initial values for the fluid and solid velocities and for solid displacement; for* $n \geq 1$*, let* $(\boldsymbol{u}_h^n, p_h^n, \dot{\boldsymbol{d}}_h^n, \boldsymbol{d}_h^n, \boldsymbol{\lambda}_h^n) \in \boldsymbol{V}_h \times Q_h \times \boldsymbol{W}_h \times \boldsymbol{W}_h \times \boldsymbol{\Lambda}_h$ *the solution of Algorithm 3, then the following energy estimates holds:*

- ***Scheme with*** $r = 0$, $\boldsymbol{d}_h^{n\star} = 0$:

$$E_h^n + D_h^n + \frac{\tau^2}{4\rho_{\mathrm{s}}\epsilon} \sum_{m=1}^{n} \|\boldsymbol{L}_h^{\mathrm{s}} \boldsymbol{d}_h^m\|_{0,\Sigma}^2 \lesssim E_h^0. \tag{4.22}$$

- ***Scheme with*** $r = 1$, $\boldsymbol{d}_h^{n\star} = \boldsymbol{d}_h^{n-1}$:

$$\begin{aligned} E_h^n + D_h^n + \frac{\tau^2}{2\rho_{\mathrm{s}}\epsilon} \sum_{m=1}^{n} \|\boldsymbol{L}_h^{\mathrm{s}} \left( \boldsymbol{d}_h^m - \boldsymbol{d}_h^{m-1} \right)\|_{0,\Sigma}^2 & \\ + \frac{\tau^2}{2} \sum_{m=1}^{n} \|\dot{\boldsymbol{d}}_h^m - \dot{\boldsymbol{d}}_h^{m-1}\|_{\mathrm{s}}^2 + \frac{\tau^2}{2\rho_{\mathrm{s}}\epsilon} \|\boldsymbol{L}_h^{\mathrm{s}} \boldsymbol{d}_h^n\|_{0,\Sigma}^2 + \frac{\tau^2}{2} \|\dot{\boldsymbol{d}}_h^n\|_{\mathrm{s}}^2 & \\ \lesssim E_h^0 + \frac{\tau^2}{2} \|\dot{\boldsymbol{d}}_h^0\|_{\mathrm{s}}^2 + \frac{\tau^2}{2\rho_{\mathrm{s}}\epsilon} \|\boldsymbol{L}_h^{\mathrm{s}} \boldsymbol{d}_h^0\|_{0,\Sigma}^2. & \end{aligned} \tag{4.23}$$

- ***Scheme with*** $r = 2$, $\boldsymbol{d}_h^{n\star} = \boldsymbol{d}_h^{n-1} + \tau \dot{\boldsymbol{d}}_h^{n-1}$:
  *let* $\tau$ *and* $h_s$ *be such that there exist* $\alpha > 0$ *such that*

$$\frac{\tau^5}{h_{\mathrm{s}}^6} \frac{\beta_{\mathrm{s}}^3 C_I^6}{(\rho_{\mathrm{s}}\epsilon)^3} \leq \alpha^5, \qquad 2\tau\alpha^5 < 1, \tag{4.24}$$

  *then for* $n \geq 1$

$$\tau^2 \sum_{m=1}^{n} \|\dot{\boldsymbol{d}}_h^m - \dot{\boldsymbol{d}}_h^{m-1}\|_{\mathrm{s}}^2 + E_h^n + D_h^n \lesssim exp\left( \tau \sum_{m=1}^{n} \frac{2\alpha^5}{1 - 2\tau\alpha^5} \right) E_h^0 \tag{4.25}$$

*Proof.* Let us consider the fully discrete problem of Algorithm 3 as expressed in (4.18) with test functions

$$(\boldsymbol{v}_h, q_h) = \tau(\boldsymbol{u}_h^n, p_h^n), \qquad \boldsymbol{\mu}_h = \tau \boldsymbol{\lambda}_h^n, \qquad \boldsymbol{w}_h = \tau \dot{\boldsymbol{d}}_h^{n-\frac{1}{2}},$$

we obtain the following equation which represents the discrete form of the continuum energy equation (2.1)

$$\rho_{\mathrm{f}}\big(\boldsymbol{u}_h^n-\boldsymbol{u}_h^{n-1},\boldsymbol{u}_h^n\big)_{0,\Omega}^2+2\mu\tau\|\boldsymbol{\epsilon}(\boldsymbol{u}_h^n)\|_{0,\Omega}^2+\tau|p_h^n|_{s_h}^2 \\ +\rho_{\mathrm{s}}\epsilon\left(\dot{\boldsymbol{d}}_h^n-\dot{\boldsymbol{d}}_h^{n-1},\dot{\boldsymbol{d}}_h^{n-\frac{1}{2}}\right)_{0,\Sigma}+\tau a_{\mathrm{s}}\left(\boldsymbol{d}_h^n,\dot{\boldsymbol{d}}_h^{n-\frac{1}{2}}\right)=0. \quad (4.26)$$

Thanks to the definition of the discrete solid operator (3.25) and the expression of $\dot{\boldsymbol{d}}_h^{n-\frac{1}{2}}$ given in equation (4.17), we have

$$\frac{\rho_{\mathrm{f}}}{2}\left(\|\boldsymbol{u}_h^n\|_{0,\Omega}^2-\|\boldsymbol{u}_h^{n-1}\|_{0,\Omega}^2+\|\boldsymbol{u}_h^n-\boldsymbol{u}_h^{n-1}\|_{0,\Omega}^2\right)+2\mu\tau\|\boldsymbol{\epsilon}(\boldsymbol{u}_h^n)\|_{0,\Omega}^2+\tau|p_h^n|_{s_h}^2 \\ +\rho_{\mathrm{s}}\epsilon\left(\dot{\boldsymbol{d}}_h^n-\dot{\boldsymbol{d}}_h^{n-1},\dot{\boldsymbol{d}}_h^n\right)_{0,\Sigma}+\tau\left(\dot{\boldsymbol{d}}_h^n-\dot{\boldsymbol{d}}_h^{n-1},\boldsymbol{L}_h^{\mathrm{s}}\left(\boldsymbol{d}_h^n-\boldsymbol{d}_h^{n\star}\right)\right)_{0,\Sigma} \\ +\tau\left(\boldsymbol{L}_h^{\mathrm{s}}\boldsymbol{d}_h^n,\dot{\boldsymbol{d}}_h^n\right)_{0,\Sigma}+\frac{\tau^2}{\rho_{\mathrm{s}}\epsilon}\left(\boldsymbol{L}_h^{\mathrm{s}}\boldsymbol{d}_h^n,\boldsymbol{L}_h^{\mathrm{s}}\left(\boldsymbol{d}_h^n-\boldsymbol{d}_h^{n\star}\right)\right)_{0,\Sigma}=0. \quad (4.27)$$

Thanks to the the expressions

$$\dot{\boldsymbol{d}}_h^n=\frac{\boldsymbol{d}_h^n-\boldsymbol{d}_h^{n-1}}{\tau}, \qquad \partial_\tau\dot{\boldsymbol{d}}_h^n=\frac{\dot{\boldsymbol{d}}_h^n-\dot{\boldsymbol{d}}_h^{n-1}}{\tau},$$

after some computation, we have

$$\left(\frac{\rho_{\mathrm{f}}}{2}\|\boldsymbol{u}_h^n\|_{0,\Omega}^2+\frac{\rho_{\mathrm{s}}\epsilon}{2}\|\dot{\boldsymbol{d}}_h^n\|_{0,\Sigma}^2+\frac{1}{2}\|\boldsymbol{d}_h^n\|_{\mathrm{s}}^2\right)-\left(\frac{\rho_{\mathrm{f}}}{2}\|\boldsymbol{u}_h^{n-1}\|_{0,\Omega}^2+\frac{\rho_{\mathrm{s}}\epsilon}{2}\|\dot{\boldsymbol{d}}_h^{n-1}\|_{0,\Sigma}^2+\frac{1}{2}\|\boldsymbol{d}_h^{n-1}\|_{\mathrm{s}}^2\right) \\ +2\mu\tau\|\boldsymbol{\epsilon}(\boldsymbol{u}_h^n)\|_{0,\Omega}^2+\tau|p_h^n|_{s_h}^2+\frac{\rho_{\mathrm{f}}}{2}\tau^2\|\partial_\tau\boldsymbol{u}_h^n\|_{0,\Omega}^2+\frac{\rho_{\mathrm{s}}\epsilon}{2}\tau^2\|\partial_\tau\dot{\boldsymbol{d}}_h^n\|_{0,\Sigma}^2+\frac{1}{2}\tau^2\|\partial_\tau\boldsymbol{d}_h^n\|_{\mathrm{s}}^2 \\ +\underbrace{\tau^2\left(\partial_\tau\dot{\boldsymbol{d}}_h^n,\boldsymbol{L}_h^{\mathrm{s}}\left(\boldsymbol{d}_h^n-\boldsymbol{d}_h^{n\star}\right)\right)_{0,\Sigma}}_{T_1}+\underbrace{\frac{\tau^2}{\rho_{\mathrm{s}}\epsilon}\left(\boldsymbol{L}_h^{\mathrm{s}}\boldsymbol{d}_h^n,\boldsymbol{L}_h^{\mathrm{s}}\left(\boldsymbol{d}_h^n-\boldsymbol{d}_h^{n\star}\right)\right)_{0,\Sigma}}_{T_2}=0, \quad (4.28)$$

Hence, we have obtained the energy estimate that depends on the extrapolation $\boldsymbol{d}_h^{n\star}$, in the following we analyze separately the particular expressions for each extrapolations.

- **Scheme with** $r=0$, $\boldsymbol{d}_h^{n\star}=0$:
  for $n\geq 1$

  $$T_1+T_2=\left(\partial_\tau\dot{\boldsymbol{d}}_h^n,\tau^2\boldsymbol{L}_h^{\mathrm{s}}\boldsymbol{d}_h^n\right)_{0,\Sigma}+\frac{\tau^2}{\rho_{\mathrm{s}}\epsilon}\|\boldsymbol{L}_h^{\mathrm{s}}\boldsymbol{d}_h^n\|_{0,\Sigma}^2.$$

  We use the generalized Young's inequality to obtain an estimate from below of the term $T_1$

  $$T_1=\left(\partial_\tau\dot{\boldsymbol{d}}_h^n,\tau^2\boldsymbol{L}_h^{\mathrm{s}}\boldsymbol{d}_h^n\right)_{0,\Sigma}\geq-\frac{\rho_{\mathrm{s}}\epsilon}{2\delta}\|\dot{\boldsymbol{d}}_h^n-\dot{\boldsymbol{d}}_h^{n-1}\|_{0,\Sigma}^2-\frac{\tau^2\delta}{2\rho_{\mathrm{s}}\epsilon}\|\boldsymbol{L}_h^{\mathrm{s}}\boldsymbol{d}_h^n\|_{0,\Sigma}^2 \qquad \forall\delta>0.$$

  hence, with $\delta=3/2$,

  $$T_1+T_2\geq-\frac{\rho_{\mathrm{s}}\epsilon}{3}\|\dot{\boldsymbol{d}}_h^n-\dot{\boldsymbol{d}}_h^{n-1}\|_{0,\Sigma}^2+\frac{\tau^2}{4\rho_{\mathrm{s}}\epsilon}\|\boldsymbol{L}_h^{\mathrm{s}}\boldsymbol{d}_h^n\|_{0,\Sigma}^2.$$

Considering this expression in the energy equation (4.28), and summing over $m = 1,...,n$, we obtain

$$\left(\frac{\rho_{\mathrm{f}}}{2}\|\boldsymbol{u}_h^n\|_{0,\Omega}^2+\frac{1}{2}\|\boldsymbol{d}_h^n\|_{\mathrm{s}}^2+\frac{\rho_{\mathrm{s}}\epsilon}{2}\|\dot{\boldsymbol{d}}_h^n\|_{0,\Sigma}^2\right)-\left(\frac{\rho_{\mathrm{f}}}{2}\|\boldsymbol{u}_h^0\|_{0,\Omega}^2+\frac{1}{2}\|\boldsymbol{d}_h^0\|_{\mathrm{s}}^2+\frac{\rho_{\mathrm{s}}\epsilon}{2}\|\dot{\boldsymbol{d}}_h^0\|_{0,\Sigma}^2\right)$$
$$+\tau\sum_{m=1}^{n}\left(2\mu\|\boldsymbol{\epsilon}(\boldsymbol{u}_h^m)\|_{0,\Omega}^2+|p_h^m|_{s_h}^2\right)+\tau^2\sum_{m=1}^{n}\left(\frac{\rho_{\mathrm{f}}}{2}\|\partial_\tau\boldsymbol{u}_h^m\|_{0,\Omega}^2+\frac{\rho_{\mathrm{s}}\epsilon}{6}\|\partial_\tau\dot{\boldsymbol{d}}_h^m\|_{0,\Sigma}^2+\frac{1}{2}\|\partial_\tau\boldsymbol{d}_h^m\|_{\mathrm{s}}^2\right)$$
$$+\frac{\tau^2}{4\rho_{\mathrm{s}}\epsilon}\sum_{m=1}^{n}\|\boldsymbol{L}_h^{\mathrm{s}}\boldsymbol{d}_h^m\|_{0,\Sigma}^2\leq 0.$$

Recalling the quantities defined in (4.11) and (4.21) we obtain the desired result

$$E_h^n+D_h^n+\frac{\tau^2}{4\rho_{\mathrm{s}}\epsilon}\sum_{m=1}^{n}\|\boldsymbol{L}_h^{\mathrm{s}}\boldsymbol{d}_h^m\|_{0,\Sigma}^2\lesssim E_h^0.$$

- **Scheme with $r=1$, $\boldsymbol{d}_h^{n\star}=\boldsymbol{d}_h^{n-1}$:**
  for $n\geq 1$ we can rewrite the terms $T_1$ and $T_2$ as follow

$$T_1=\left(\partial_\tau\dot{\boldsymbol{d}}_h^n,\tau^2\boldsymbol{L}_h^{\mathrm{s}}\left(\boldsymbol{d}_h^n-\boldsymbol{d}_h^{n-1}\right)\right)_{0,\Sigma}=\tau^2 a_{\mathrm{s}}\left(\dot{\boldsymbol{d}}_h^n-\dot{\boldsymbol{d}}_h^{n-1},\dot{\boldsymbol{d}}_h^n\right)$$
$$=\frac{\tau^2}{2}\left(\|\dot{\boldsymbol{d}}_h^n\|_{\mathrm{s}}^2-\|\dot{\boldsymbol{d}}_h^{n-1}\|_{\mathrm{s}}^2+\|\dot{\boldsymbol{d}}_h^n-\dot{\boldsymbol{d}}_h^{n-1}\|_{\mathrm{s}}^2\right)_{0,\Sigma}$$

$$T_2=\frac{\tau^2}{2\rho_{\mathrm{s}}\epsilon}\left(\|\boldsymbol{L}_h^{\mathrm{s}}\boldsymbol{d}_h^n\|_{0,\Sigma}^2-\|\boldsymbol{L}_h^{\mathrm{s}}\boldsymbol{d}_h^{n-1}\|_{0,\Sigma}^2+\|\boldsymbol{L}_h^{\mathrm{s}}\left(\boldsymbol{d}_h^n-\boldsymbol{d}_h^{n-1}\right)\|_{0,\Sigma}^2\right)$$

Considering these expressions in the energy equation (4.28), and summing over $m = 1,...,n$, we obtain

$$\left(\frac{\rho_{\mathrm{f}}}{2}\|\boldsymbol{u}_h^n\|_{0,\Omega}^2+\frac{1}{2}\|\boldsymbol{d}_h^n\|_{\mathrm{s}}^2+\frac{\rho_{\mathrm{s}}\epsilon}{2}\|\dot{\boldsymbol{d}}_h^n\|_{0,\Sigma}^2\right)-\left(\frac{\rho_{\mathrm{f}}}{2}\|\boldsymbol{u}_h^0\|_{0,\Omega}^2+\frac{1}{2}\|\boldsymbol{d}_h^0\|_{\mathrm{s}}^2+\frac{\rho_{\mathrm{s}}\epsilon}{2}\|\dot{\boldsymbol{d}}_h^0\|_{0,\Sigma}^2\right)$$
$$+\tau\sum_{m=1}^{n}\left(2\mu\|\boldsymbol{\epsilon}(\boldsymbol{u}_h^m)\|_{0,\Omega}^2+|p_h^m|_{s_h}^2\right)+\tau^2\sum_{m=1}^{n}\left(\frac{\rho_{\mathrm{f}}}{2}\|\partial_\tau\boldsymbol{u}_h^m\|_{0,\Omega}^2+\frac{\rho_{\mathrm{s}}\epsilon}{2}\|\partial_\tau\dot{\boldsymbol{d}}_h^m\|_{0,\Sigma}^2+\frac{1}{2}\|\partial_\tau\boldsymbol{d}_h^m\|_{\mathrm{s}}^2\right)$$
$$\frac{\tau^2}{2}\|\dot{\boldsymbol{d}}_h^n\|_{\mathrm{s}}^2-\frac{\tau^2}{2}\|\dot{\boldsymbol{d}}_h^0\|_{\mathrm{s}}^2+\frac{\tau^2}{2}\sum_{m=1}^{n}\|\dot{\boldsymbol{d}}_h^m-\dot{\boldsymbol{d}}_h^{m-1}\|_{\mathrm{s}}^2$$
$$+\frac{\tau^2}{2\rho_{\mathrm{s}}\epsilon}\|\boldsymbol{L}_h^{\mathrm{s}}\boldsymbol{d}_h^n\|_{0,\Sigma}^2-\frac{\tau^2}{2\rho_{\mathrm{s}}\epsilon}\|\boldsymbol{L}_h^{\mathrm{s}}\boldsymbol{d}_h^0\|_{0,\Sigma}^2+\frac{\tau^2}{2\rho_{\mathrm{s}}\epsilon}\sum_{m=1}^{n}\|\boldsymbol{L}_h^{\mathrm{s}}\left(\boldsymbol{d}_h^m-\boldsymbol{d}_h^{m-1}\right)\|_{0,\Sigma}^2=0.$$

Recalling the quantities defined in (4.11) and (4.21) we obtain the desired result

$$E_h^n+D_h^n+\frac{\tau^2}{2\rho_{\mathrm{s}}\epsilon}\sum_{m=1}^{n}\|\boldsymbol{L}_h^{\mathrm{s}}\left(\boldsymbol{d}_h^m-\boldsymbol{d}_h^{m-1}\right)\|_{0,\Sigma}^2+\frac{\tau^2}{2}\sum_{m=1}^{n}\|\dot{\boldsymbol{d}}_h^m-\dot{\boldsymbol{d}}_h^{m-1}\|_{\mathrm{s}}^2$$
$$+\frac{\tau^2}{2\rho_{\mathrm{s}}\epsilon}\|\boldsymbol{L}_h^{\mathrm{s}}\boldsymbol{d}_h^n\|_{0,\Sigma}^2+\frac{\tau^2}{2}\|\dot{\boldsymbol{d}}_h^n\|_{\mathrm{s}}^2\lesssim E_h^0+\frac{\tau^2}{2}\|\dot{\boldsymbol{d}}_h^0\|_{\mathrm{s}}^2+\frac{\tau^2}{2\rho_{\mathrm{s}}\epsilon}\|\boldsymbol{L}_h^{\mathrm{s}}\boldsymbol{d}_h^0\|_{0,\Sigma}^2.$$

- **Scheme with** $r = 2$, $\boldsymbol{d}_h^{n\star} = \boldsymbol{d}_h^{n-1} + \tau\dot{\boldsymbol{d}}_h^{n-1}$:

$$T_1 = \left(\partial_\tau \dot{\boldsymbol{d}}_h^n, \tau^2 \boldsymbol{L}_h^{\mathrm{s}}\left(\boldsymbol{d}_h^n - \boldsymbol{d}_h^{n-1} - \tau\dot{\boldsymbol{d}}_h^{n-1}\right)\right)_{0,\Sigma}$$
$$= \tau^2 a_{\mathrm{s}}\left(\dot{\boldsymbol{d}}_h^n - \dot{\boldsymbol{d}}_h^{n-1}, \dot{\boldsymbol{d}}_h^n - \dot{\boldsymbol{d}}_h^{n-1}\right) = \tau^2 \|\dot{\boldsymbol{d}}_h^n - \dot{\boldsymbol{d}}_h^{n-1}\|_{\mathrm{s}}^2. \quad (4.29)$$

$$T_2 = \frac{\tau^2}{\rho_{\mathrm{s}}\epsilon}\left(\boldsymbol{L}_h^{\mathrm{s}}\boldsymbol{d}_h^n, \boldsymbol{L}_h^{\mathrm{s}}\left(\boldsymbol{d}_h^n - \boldsymbol{d}_h^{n-1} - \tau\dot{\boldsymbol{d}}_h^{n-1}\right)\right)_{0,\Sigma} = \frac{\tau^3}{\rho_{\mathrm{s}}\epsilon}\left(\boldsymbol{L}_h^{\mathrm{s}}\boldsymbol{d}_h^n, \boldsymbol{L}_h^{\mathrm{s}}\left(\dot{\boldsymbol{d}}_h^n - \dot{\boldsymbol{d}}_h^{n-1}\right)\right)_{0,\Sigma}$$
$$= a_{\mathrm{s}}\left(\boldsymbol{L}_h^{\mathrm{s}}\boldsymbol{d}_h^n, \dot{\boldsymbol{d}}_h^n - \dot{\boldsymbol{d}}_h^{n-1}\right) \geq -\frac{\tau^3}{\rho_{\mathrm{s}}\epsilon}\|\boldsymbol{L}_h^{\mathrm{s}}\boldsymbol{d}_h^n\|_{\mathrm{s}}\|\dot{\boldsymbol{d}}_h^n - \dot{\boldsymbol{d}}_h^{n-1}\|_{\mathrm{s}}$$

where we have used the Cauchy-Schwartz's inequality, now by inverse inequalities (3.28) and (3.30) and Young's inequality we obtain an estimate from below of the term $T_2$

$$\begin{aligned} T_2 &\geq -\frac{\tau^3}{\rho_{\mathrm{s}}\epsilon}\left(\frac{\beta_{\mathrm{s}}C_I^2}{h_s^2}\|\boldsymbol{d}_h^n\|_{\mathrm{s}}\right)\left(\frac{\beta_{\mathrm{s}}^{1/2}C_I}{h_s}\|\dot{\boldsymbol{d}}_h^n - \dot{\boldsymbol{d}}_h^{n-1}\|_{0,\Sigma}\right) \\ &\geq -\left(\frac{\tau^3}{(\rho_{\mathrm{s}}\epsilon)^{3/2}}\frac{\beta_{\mathrm{s}}^{3/2}C_I^3}{h_s^3}\|\boldsymbol{d}_h^n\|_{\mathrm{s}}\right)\left((\rho_{\mathrm{s}}\epsilon)^{1/2}\|\dot{\boldsymbol{d}}_h^n - \dot{\boldsymbol{d}}_h^{n-1}\|_{0,\Sigma}\right) \\ &\geq -\frac{\tau^6}{h_s^6}\frac{\beta_{\mathrm{s}}^3 C_I^6}{(\rho_{\mathrm{s}}\epsilon)^3}\|\boldsymbol{d}_h^n\|_{\mathrm{s}}^2 - \frac{\rho_{\mathrm{s}}\epsilon}{4}\|\dot{\boldsymbol{d}}_h^n - \dot{\boldsymbol{d}}_h^{n-1}\|_{0,\Sigma}^2 \geq -\tau\alpha^5\|\boldsymbol{d}_h^n\|_{\mathrm{s}}^2 - \frac{\rho_{\mathrm{s}}\epsilon}{4}\|\dot{\boldsymbol{d}}_h^n - \dot{\boldsymbol{d}}_h^{n-1}\|_{0,\Sigma}^2, \end{aligned} \quad (4.30)$$

where $\alpha$ is the constant given in the statement for the fixed $\tau$ and $h_s$. Considering expressions (4.29) and (4.30) in the energy equation (4.28) and summing over $m = 1,\ldots,n$, we obtain

$$\begin{aligned} &\tau^2\sum_{m=1}^n \|\dot{\boldsymbol{d}}_h^m - \dot{\boldsymbol{d}}_h^{m-1}\|_{\mathrm{s}}^2 + \frac{\rho_{\mathrm{f}}}{2}\|\boldsymbol{u}_h^n\|_{0,\Omega}^2 + \underbrace{\frac{1}{2}\|\boldsymbol{d}_h^n\|_{\mathrm{s}}^2}_{a_n} + \frac{\rho_{\mathrm{s}}\epsilon}{2}\|\dot{\boldsymbol{d}}_h^n\|_{0,\Sigma}^2 \\ &+ \tau\sum_{m=1}^n \underbrace{\left[2\mu\|\boldsymbol{\epsilon}(\boldsymbol{u}_h^m)\|_{0,\Omega}^2 + |p_h^m|_{s_h}^2 + \tau\left(\frac{\rho_{\mathrm{f}}}{2}\|\partial_\tau\boldsymbol{u}_h^m\|_{0,\Omega}^2 + \frac{\rho_{\mathrm{s}}\epsilon}{4}\|\partial_\tau\dot{\boldsymbol{d}}_h^m\|_{0,\Sigma}^2 + \frac{1}{2}\|\partial_\tau\boldsymbol{d}_h^m\|_{\mathrm{s}}^2\right)\right]}_{b_m} \\ &\leq \tau\sum_{m=1}^n \underbrace{2\alpha^5}_{\gamma_m}\underbrace{\frac{1}{2}\|\boldsymbol{d}_h^m\|_{\mathrm{s}}^2}_{a_m} + \underbrace{\left(\frac{\rho_{\mathrm{f}}}{2}\|\boldsymbol{u}_h^0\|_{0,\Omega}^2 + \frac{1}{2}\|\boldsymbol{d}_h^0\|_{\mathrm{s}}^2 + \frac{\rho_{\mathrm{s}}\epsilon}{2}\|\dot{\boldsymbol{d}}_h^0\|_{0,\Sigma}^2\right)}_{B=E_h^0}. \end{aligned}$$

Since we have, by hypothesis, $2\tau\alpha^5 < 1$, we can apply the discrete Gronwall Lemma 3.3.1 obtaining

$$\tau^2\sum_{m=1}^n \|\dot{\boldsymbol{d}}_h^m - \dot{\boldsymbol{d}}_h^{m-1}\|_{\mathrm{s}}^2 + E_h^n + D_h^n \lesssim exp\left(\tau\sum_{m=1}^n \frac{2\alpha^5}{1-2\tau\alpha^5}\right)E_h^0$$

□

Theorem 4.3.3 shows that Algorithm 3 with r=0,1 is unconditionally stable, whereas the variant $r = 2$ is conditionally stable, in fact for $r = 2$, in order to obtain stability we assume a CFL condition, that is

$$\tau \leq C h_s^{\frac{6}{5}},$$

where $C = C(\rho_s, \epsilon, C_I, \beta_s) > 0$ is a constant that depends on the physical parameters of the solid ($\rho_s$ and $\epsilon$), on the inverse inequality constant $C_I$ and on the continuity constant $\beta_s$ of $a_s(\cdot,\cdot)$.

## 4.4 Convergence analysis of Algorithm 3 - linearized case.

As for the monolithic scheme, in order to perform a convergence analysis of the proposed Algorithms, we simplify the setting of the fluid-solid interaction problem considering a linearized version of Problem 2 obtained assuming small solid displacements. This assumption allows for considering the actual solid configuration coincident with the reference configuration at each instant of time, that is $\boldsymbol{\Sigma}(t) = \boldsymbol{\Sigma}$ for all $T \geq t \geq 0$. We report here, for convenience of the reader, the linearized problem.

**Problem 10.** *For $t \geq 0$, find $\big(\boldsymbol{u}(t), p(t), \boldsymbol{d}(t), \dot{\boldsymbol{d}}(t), \boldsymbol{\lambda}(t)\big) \in \boldsymbol{V} \times Q \times \boldsymbol{W} \times \boldsymbol{W} \times \boldsymbol{\Lambda}$ such that*

- $\partial_t \boldsymbol{d} = \dot{\boldsymbol{d}}$
- *for all* $(\boldsymbol{v}, q, \boldsymbol{w}, \boldsymbol{\mu}) \in \boldsymbol{V} \times Q \times \boldsymbol{W} \times \boldsymbol{\Lambda}$

$$\rho_f\big(\partial_t \boldsymbol{u}, \boldsymbol{v}\big)_{0,\Omega} + a^f\big((\boldsymbol{u}, p), (\boldsymbol{v}, q)\big) + c\big(\boldsymbol{\lambda}, \boldsymbol{v}|_\Sigma - \boldsymbol{w}\big) \\ - c\big(\boldsymbol{\mu}, \boldsymbol{u}|_\Sigma - \dot{\boldsymbol{d}}\big) + \rho_s \epsilon \big(\partial_t \dot{\boldsymbol{d}}, \boldsymbol{w}\big)_{0,\Sigma} + a_s\big(\boldsymbol{d}, \boldsymbol{w}\big) = 0 \quad (4.31)$$

- $\boldsymbol{u}(0) = \boldsymbol{u}_0$ *in* $\Omega$, $\dot{\boldsymbol{d}}(0) = \dot{\boldsymbol{d}}_0$, $\boldsymbol{d}(0) = \boldsymbol{d}_0$ *in* $\Sigma$.

### 4.4.1 Convergence analysis

As in the case of the monolithic algorithm, we will study the behavior of the error when the discretization parameters $\tau$, $h_s$ and $h_f$ go to zero. We split the error $\boldsymbol{e} = \boldsymbol{a} - \boldsymbol{a}_h$ as follows

$$\boldsymbol{e} = \boldsymbol{a} - \boldsymbol{a}_h = \underbrace{\boldsymbol{a} - \boldsymbol{a}_\Pi}_{\boldsymbol{e}_\Pi} + \underbrace{\boldsymbol{a}_\Pi - \boldsymbol{a}_h}_{\boldsymbol{e}_h},$$

where $\boldsymbol{a}_\Pi$ is the Ritz projection of $\boldsymbol{a}$ on the discrete space. The term $\boldsymbol{e}_\Pi$ is called projection error and the part $\boldsymbol{e}_h$ is the discrete error. We use the same notation of Definition 2 of Chapter 3, moreover we introduce the following notation

$$\boldsymbol{\chi}_h^n \overset{def}{=} \dot{\boldsymbol{d}}_\Pi^n - \dot{\boldsymbol{d}}_h^{n-\frac{1}{2}}. \quad (4.32)$$

In the following proposition we give an useful expression of $\boldsymbol{\chi}_h^n$.

**Proposition 4.4.1.** *With reference to Definition 2, we have*

$$\chi_h^n = \dot{\xi}_h^n + \frac{\tau}{\rho_s\epsilon} L_h^s(\xi_h^n - \xi_h^{n\star}) - \frac{\tau}{\rho_s\epsilon} L_h^s(d^n - d^{n\star}), \tag{4.33}$$

*Proof.* The claim follows from the definition of $\chi_h^n$, equation (4.17) and the relation $L_h^s \circ \Pi_W = L_h^s$ (see (3.27)),

$$\begin{aligned}
\chi_h^n = \dot{d}_\Pi^n - \dot{d}_h^{n-\frac{1}{2}} = \dot{d}_\Pi^n - \dot{d}_h^n + \dot{d}_h^n - \dot{d}_h^{n-\frac{1}{2}} = \dot{\xi}_h^n - \frac{\tau}{\rho_s\epsilon} L_h^s(d_h^n - d_h^{n\star}) \\
= \dot{\xi}_h^n - \frac{\tau}{\rho_s\epsilon} L_h^s(d_h^n - \Pi_W d^n + \Pi_W d^{n\star} - d_h^{n\star}) - \frac{\tau}{\rho_s\epsilon} L_h^s(\Pi_W(d^n - d^{n\star})) \\
= \dot{\xi}_h^n + \frac{\tau}{\rho_s\epsilon} L_h^s(\xi_h^n - \xi_h^{n\star}) - \frac{\tau}{\rho_s\epsilon} L_h^s(d^n - d^{n\star}).
\end{aligned}$$

□

**Convergence of Algorithm 3** In this paragraph we give the result about the convergence of Algorithm 3. The proof of the Theorem is given for the extrapolation orders $r = 0, 1, 2$ and it uses in all the cases the discrete Gronwall Lemma 3.3.1.

**Theorem 4.4.2.** *Convergence of Algorithm 3*

*Let $(u, p, \dot{d}, d, \lambda) \in V \times Q \times W \times W \times \Lambda$ satisfy Problem 10 with the regularity assumptions (3.39) and $\{(u_h^n, p_h^n, \dot{d}_h^n, \dot{d}_h^{n-\frac{1}{2}}, d_h^n, \lambda_h^n)\}_{n>0} \subset V_h \times Q_h \times W_h \times W_h \times W_h \times \Lambda_h$ be given by Algorithm 3 with the stability hypothesis of Theorem 4.3.3, then we have the following error estimates. For $n > 0$ and $n\tau \le T$*

- **Scheme with** $r = 0$, *if* $\tau < 1$

$$\begin{aligned}
\mathscr{E}_h^n \stackrel{def}{=} \frac{\rho_f}{2}\|\theta_h^n\|_{0,\Omega}^2 + \frac{\rho_s\epsilon}{2}\|\dot{\xi}_h^n\|_{0,\Sigma}^2 + \frac{1}{2}\|\xi_h^n\|_s^2 \\
\lesssim \exp\left(\sum_{m=1}^{n} \frac{\tau}{1-\tau}\right)\Big[\underbrace{\mathcal{O}(h_f^{2l}) + \mathcal{O}(h_s^{2m+1})}_{A_0} \\
+ \tau\underbrace{\Big(\mathcal{O}(h_s^{2m}) + \mathcal{O}(h_s^{2z+1}) + \mathcal{O}(h_f^{2l})\Big)}_{A_1} + \tau^2 A_2 + \tau^3 A_3 + \tau^4 A_4\Big].
\end{aligned} \tag{4.34}$$

- **Scheme with** $r = 1$, *if* $\tau < 1$ *then*

$$\mathscr{E}_h^n \stackrel{def}{=} \frac{\rho_f}{2}\|\theta_h^n\|_{0,\Omega}^2 + \frac{\rho_s\epsilon}{2}\|\dot{\xi}_h^n\|_{0,\Sigma}^2 + \frac{1}{2}\|\xi_h^n\|_s^2 + \frac{\tau^2}{2\rho_s\epsilon}\|L_h^s\xi_h^n\|_{0,\Sigma}^2 + \frac{\tau^2}{2}\|\dot{\xi}_h^n\|_s^2$$

$$\lesssim \exp\left(\sum_{m=1}^{n} \frac{\tau}{1-\tau}\right) \Big[\underbrace{\mathcal{O}(h_f^{2l}) + \mathcal{O}(h_s^{2m+1})}_{A_0} + \tau\Big(\underbrace{\mathcal{O}(h_s^{2m}) + \mathcal{O}(h_s^{2z+1}) + \mathcal{O}(h_f^{2l})}_{A_1}\Big) + \tau^2 A_2 + \tau^3 A_3\Big]. \quad (4.35)$$

- **Scheme with** $r = 2$, *for* $m = 1, \ldots, n$, *let* $\gamma_m = \max\{1, 4\alpha^5\}$, *if* $\tau\gamma_m < 1$ *then*

$$\mathcal{E}_h^n \stackrel{def}{=} \frac{\rho_{\mathrm{f}}}{2}\|\boldsymbol{\theta}_h^n\|_{0,\Omega}^2 + \frac{\rho_s\epsilon}{2}\|\dot{\boldsymbol{\xi}}_h^n\|_{0,\Sigma}^2 + \frac{1}{2}\|\boldsymbol{\xi}_h^n\|_{\mathrm{s}}^2$$
$$\lesssim \exp\left(\sum_{m=1}^{n} \frac{\gamma_m\tau}{1-\gamma_m\tau}\right) \Big[\underbrace{\mathcal{O}(h_f^{2l}) + \mathcal{O}(h_s^{2m+1})}_{A_0} + \tau\Big(\underbrace{\mathcal{O}(h_s^{2m+1}) + \mathcal{O}(h_s^{2z+1}) + \mathcal{O}(h_f^{2l})}_{A_1}\Big) + \tau^2 A_2 + \tau^3 A_3\Big]. \quad (4.36)$$

*Proof.* Let us consider Problem 10 for $t = t^n$ and test functions taken in the discrete spaces

$$\rho_{\mathrm{f}}\big(\partial_t \boldsymbol{u}^n, \boldsymbol{v}_h\big)_{0,\Omega} + a^{\mathrm{f}}\big((\boldsymbol{u}^n, p^n), (\boldsymbol{v}_h, q_h)\big) + c\left(\boldsymbol{\lambda}^n, \boldsymbol{v}_h|_\Sigma - \boldsymbol{w}_h\right) - c\left(\boldsymbol{\mu}_h, \boldsymbol{u}^n|_\Sigma - \dot{\boldsymbol{d}}^n\right) + \rho_s\epsilon\left(\partial_t \dot{\boldsymbol{d}}^n, \boldsymbol{w}_h\right)_{0,\Sigma} + a_{\mathrm{s}}\left(\boldsymbol{d}^n, \boldsymbol{w}_h\right) = 0, \quad (4.37)$$

for all $(\boldsymbol{v}_h, q_h, \boldsymbol{w}_h, \boldsymbol{\mu}_h) \in \boldsymbol{V}_h \times Q_h \times \boldsymbol{W}_h \times \boldsymbol{\Lambda}_h$.

The following error equation is obtained subtracting the expression of Algorithm 3 given in (4.18) from (4.37) and adding and subtracting the terms $\rho_{\mathrm{f}}(\partial_\tau \boldsymbol{u}^n, \boldsymbol{v}_h)_{0,\Omega}$ and $\rho_s\epsilon\left(\partial_\tau \dot{\boldsymbol{d}}^n \boldsymbol{w}_h\right)_{0,\Sigma}$

$$\rho_{\mathrm{f}}\big(\partial_\tau(\boldsymbol{u}^n - \boldsymbol{u}_h^n), \boldsymbol{v}_h\big)_{0,\Omega} + a^{\mathrm{f}}\big((\boldsymbol{u}^n - \boldsymbol{u}_h^n, p^n - p_h^n), (\boldsymbol{v}_h, q_h)\big) + c\left(\boldsymbol{\lambda}^n - \boldsymbol{\lambda}_h^n, \boldsymbol{v}_h|_\Sigma - \boldsymbol{w}_h\right)$$
$$- c\left(\boldsymbol{\mu}_h, \boldsymbol{u}^n|_\Sigma - \boldsymbol{u}_h^n|_\Sigma\right) + c\left(\boldsymbol{\mu}_h, \dot{\boldsymbol{d}}^n - \dot{\boldsymbol{d}}_h^{n-\frac{1}{2}}\right) + \rho_s\epsilon\left(\partial_\tau(\dot{\boldsymbol{d}}^n - \dot{\boldsymbol{d}}_h^n), \boldsymbol{w}_h\right)_{0,\Sigma} + a_{\mathrm{s}}\left(\boldsymbol{d}^n - \boldsymbol{d}_h^n, \boldsymbol{w}_h\right)$$
$$= s_h(p_h^n, q_h) + \rho_{\mathrm{f}}\big((\partial_\tau - \partial_t)\boldsymbol{u}^n, \boldsymbol{v}_h\big)_{0,\Omega} + \rho_s\epsilon\left((\partial_\tau - \partial_t)\dot{\boldsymbol{d}}^n, \boldsymbol{w}_h\right)_{0,\Sigma} \quad (4.38)$$

for all $(\boldsymbol{v}_h, q_h, \boldsymbol{w}_h, \boldsymbol{\mu}_h) \in \boldsymbol{V}_h \times Q_h \times \boldsymbol{W}_h \times \boldsymbol{\Lambda}_h$.

Using the notations of Definition 2 of Chapter 3 and (4.32), we obtain

$$\rho_{\mathrm{f}}\big(\partial_\tau \boldsymbol{\theta}_h^n, \boldsymbol{v}_h\big)_{0,\Omega} + a^{\mathrm{f}}\big((\boldsymbol{\theta}_h^n, \varphi_h^n), (\boldsymbol{v}_h, q_h)\big) + \rho_s\epsilon\left(\partial_\tau \dot{\boldsymbol{\xi}}_h^n, \boldsymbol{w}_h\right)_{0,\Sigma} + a_{\mathrm{s}}\left(\boldsymbol{\xi}_h^n, \boldsymbol{w}_h\right) =$$
$$\rho_{\mathrm{f}}\big((\partial_\tau - \partial_t)\boldsymbol{u}^n, \boldsymbol{v}_h\big)_{0,\Omega} - \rho_{\mathrm{f}}\big(\partial_\tau \boldsymbol{\theta}_\pi^n, \boldsymbol{v}_h\big)_{0,\Omega} + \rho_s\epsilon\left((\partial_\tau - \partial_t)\dot{\boldsymbol{d}}^n, \boldsymbol{w}_h\right)_{0,\Sigma} - \rho_s\epsilon\left(\partial_\tau \dot{\boldsymbol{\xi}}_\pi^n, \boldsymbol{w}_h\right)_{0,\Omega}$$
$$- c\left(\boldsymbol{\omega}_h^n, \boldsymbol{v}_h|_\Sigma - \boldsymbol{w}_h\right) + c\big(\boldsymbol{\mu}_h, \boldsymbol{\theta}_h^n|_\Sigma - \boldsymbol{\chi}_h^n\big) + s_h(p_h^n, q_h) + c\left(\boldsymbol{\mu}_h, \boldsymbol{\theta}_\pi^n|_\Sigma - \dot{\boldsymbol{\xi}}_\pi^n\right)$$

$$-\left[a^{\mathrm{f}}\big((\boldsymbol{\theta}_\pi^n,\varphi_\pi^n),(\boldsymbol{v}_h,q_h)\big)+a_{\mathrm{s}}\left(\boldsymbol{\xi}_\pi^n,\boldsymbol{w}_h\right)+c\left(\boldsymbol{\omega}_\pi^n,\boldsymbol{v}_h|_\Sigma-\boldsymbol{w}_h\right)\right], \quad (4.39)$$

for all $(\boldsymbol{v}_h,q_h,\boldsymbol{w}_h,\boldsymbol{\mu}_h)\in \boldsymbol{V}_h\times \boldsymbol{Q}_h\times \boldsymbol{W}_h\times \boldsymbol{\Lambda}_h$.

From the definition of the projection operator given in Definition 1 of Chapter 3, we have

$$\begin{array}{lll}
\text{from (3.7)} & a_{\mathrm{s}}\left(\boldsymbol{\xi}_\pi^n,\boldsymbol{w}_h\right)=0 & \forall \boldsymbol{w}_h\in \boldsymbol{W}_h,\\
\text{from (3.8)} & a^{\mathrm{f}}\big((\boldsymbol{\theta}_\pi^n,\varphi_\pi^n),(\boldsymbol{v}_h,q_h)\big)=s_h(p_\Pi^n,q_h) & \forall(\boldsymbol{v}_h,q_h)\in \boldsymbol{V}_h\times \boldsymbol{Q}_h,\\
\text{from (3.9)} & c\left(\boldsymbol{\omega}_\pi^n,\boldsymbol{v}_h|_\Sigma-\boldsymbol{w}_h\right)=c\left(\boldsymbol{\omega}_\pi^n,\boldsymbol{v}_h|_\Sigma\right) & \forall(\boldsymbol{v}_h,\boldsymbol{w}_h)\in \boldsymbol{V}_h\times \boldsymbol{W}_h,
\end{array} \quad (4.40)$$

Hence, using (4.40) and testing with $(\boldsymbol{v}_h,q_h,\boldsymbol{\mu}_h,\boldsymbol{w}_h)=\tau(\boldsymbol{\theta}_h^n,\varphi_h^n,\boldsymbol{\omega}_h^n,\boldsymbol{\chi}_h^n)$ we have the following error expression

$$\begin{aligned}
&\rho_{\mathrm{f}}\tau\left(\partial_\tau\boldsymbol{\theta}_h^n,\boldsymbol{\theta}_h^n\right)_{0,\Omega}+\tau a^{\mathrm{f}}\big((\boldsymbol{\theta}_h^n,\varphi_h^n),(\boldsymbol{\theta}_h^n,\varphi_h^n)\big)+\rho_s\epsilon\tau\left(\partial_\tau\dot{\boldsymbol{\xi}}_h^n,\boldsymbol{\chi}_h^n\right)_{0,\Sigma}\\
&\quad+\tau a_{\mathrm{s}}\left(\boldsymbol{\xi}_h^n,\boldsymbol{\chi}_h^n\right)+\tau s_h(\varphi_h^n,\varphi_h^n)=\rho_{\mathrm{f}}\tau\left((\partial_\tau-\partial_t)\boldsymbol{u}^n,\boldsymbol{\theta}_h^n\right)_{0,\Omega}-\rho_{\mathrm{f}}\tau\left(\partial_\tau\boldsymbol{\theta}_\pi^n,\boldsymbol{\theta}_h^n\right)_{0,\Omega}\\
&\quad+\rho_s\epsilon\tau\left((\partial_\tau-\partial_t)\dot{\boldsymbol{d}}^n,\boldsymbol{\chi}_h^n\right)_{0,\Sigma}-\rho_s\epsilon\tau\left(\partial_\tau\dot{\boldsymbol{\xi}}_\pi^n,\boldsymbol{\chi}_h^n\right)_{0,\Sigma}-\tau c\left(\boldsymbol{\omega}_\pi^n,\boldsymbol{\theta}_h|_\Sigma\right)+\tau c\left(\boldsymbol{\omega}_h^n,\boldsymbol{\theta}_\pi^n|_\Sigma-\dot{\boldsymbol{\xi}}_\pi^n\right)
\end{aligned} \quad (4.41)$$

Using the results of Proposition 4.4.1 about $\boldsymbol{\chi}_h^h$ and Proposition 3.3.5 about $\dot{\boldsymbol{\xi}}_h^n$, we have the following expression where the terms $T_0-T_{10}$ are presented in (4.44).

$$\begin{aligned}
&\rho_{\mathrm{f}}\tau\left(\partial_\tau\boldsymbol{\theta}_h^n,\boldsymbol{\theta}_h^n\right)_{0,\Omega}+\tau a^{\mathrm{f}}\big((\boldsymbol{\theta}_h^n,\varphi_h^n),(\boldsymbol{\theta}_h^n,\varphi_h^n)\big)+\rho_s\epsilon\tau\left(\partial_\tau\dot{\boldsymbol{\xi}}_h^n,\dot{\boldsymbol{\xi}}_h^n\right)_{0,\Sigma}\\
&\qquad+\tau a_{\mathrm{s}}\left(\boldsymbol{\xi}_h^n,\partial_\tau\boldsymbol{\xi}_h^n\right)+\tau s_h(\varphi_h^n,\varphi_h^n)=\sum_{i=0}^{10}T_i.
\end{aligned} \quad (4.42)$$

computing the left hand side of equation (4.42), the relation can be rewritten as follows

$$\begin{aligned}
&\left(\frac{\rho_{\mathrm{f}}}{2}\|\boldsymbol{\theta}_h^n\|_{0,\Omega}^2+\frac{\rho_s\epsilon}{2}\|\dot{\boldsymbol{\xi}}_h^n\|_{0,\Sigma}^2+\frac{1}{2}\|\boldsymbol{\xi}_h^n\|_{\mathrm{s}}^2\right)-\left(\frac{\rho_{\mathrm{f}}}{2}\|\boldsymbol{\theta}_h^{n-1}\|_{0,\Omega}^2+\frac{\rho_s\epsilon}{2}\|\dot{\boldsymbol{\xi}}_h^{n-1}\|_{0,\Sigma}^2+\frac{1}{2}\|\boldsymbol{\xi}_h^{n-1}\|_{\mathrm{s}}^2\right)\\
&\quad+\tau^2\left(\frac{\rho_{\mathrm{f}}}{2}\|\partial_\tau\boldsymbol{\theta}_h^n\|_{0,\Omega}^2+\frac{\rho_s\epsilon}{2}\|\partial_\tau\dot{\boldsymbol{\xi}}_h^n\|_{0,\Sigma}^2+\frac{1}{2}\|\partial_\tau\boldsymbol{\xi}_h^n\|_{\mathrm{s}}^2\right)+\tau\left(2\mu\|\boldsymbol{\epsilon}(\boldsymbol{\theta}_h^n)\|_{0,\Omega}^2+|\varphi_h^n|_{s_h}^2\right)=\sum_{i=0}^{10}T_i.
\end{aligned} \quad (4.43)$$

The terms $T_0 - T_{10}$ are reported in the following

$$
\begin{aligned}
T_0 &= -\tau a_{\mathrm{s}} \left( \boldsymbol{\xi}_h^n, \dot{\mathbf{\Pi}}_{\boldsymbol{W}}^n - \mathbf{\Pi}_{\boldsymbol{W}} \partial_\tau \boldsymbol{d}^n \right) \\
T_1 &= \rho_{\mathrm{f}} \tau \left( (\partial_\tau - \partial_t) \boldsymbol{u}^n, \boldsymbol{\theta}_h^n \right)_{0,\Omega} - \rho_{\mathrm{f}} \tau \left( \partial_\tau \boldsymbol{\theta}_\pi^n, \boldsymbol{\theta}_h^n \right)_{0,\Omega} \\
T_2 &= \rho_s \epsilon \tau \left( (\partial_\tau - \partial_t) \dot{\boldsymbol{d}}^n, \dot{\boldsymbol{\xi}}_h^n \right)_{0,\Sigma} - \rho_s \epsilon \tau \left( \partial_\tau \dot{\boldsymbol{\xi}}_\pi^n, \dot{\boldsymbol{\xi}}_h^n \right)_{0,\Sigma} \\
T_3 &= -\tau c \left( \boldsymbol{\omega}_\pi^n, \boldsymbol{\theta}_h^n |_\Sigma \right) \\
T_4 &= c \left( \boldsymbol{\mu}_h, \boldsymbol{\theta}_\pi^n |_\Sigma - \dot{\boldsymbol{\xi}}_\pi^n \right) \\
T_5 &= \rho_s \epsilon \tau \left( (\partial_\tau - \partial_t) \dot{\boldsymbol{d}}^n, \frac{\tau}{\rho_{\mathrm{s}} \epsilon} \boldsymbol{L}_h^{\mathrm{s}} (\boldsymbol{\xi}_h^n - \boldsymbol{\xi}_h^{n\star}) \right)_{0,\Sigma} - \rho_s \epsilon \tau \left( \partial_\tau \dot{\boldsymbol{\xi}}_\pi^n, \frac{\tau}{\rho_{\mathrm{s}} \epsilon} \boldsymbol{L}_h^{\mathrm{s}} (\boldsymbol{\xi}_h^n - \boldsymbol{\xi}_h^{n\star}) \right)_{0,\Sigma} \\
T_6 &= -\rho_s \epsilon \tau \left( (\partial_\tau - \partial_t) \dot{\boldsymbol{d}}^n, \frac{\tau}{\rho_{\mathrm{s}} \epsilon} \boldsymbol{L}_h^{\mathrm{s}} (\boldsymbol{d}^n - \boldsymbol{d}^{n\star}) \right)_{0,\Sigma} + \rho_s \epsilon \tau \left( \partial_\tau \dot{\boldsymbol{\xi}}_\pi^n, \frac{\tau}{\rho_{\mathrm{s}} \epsilon} \boldsymbol{L}_h^{\mathrm{s}} (\boldsymbol{d}^n - \boldsymbol{d}^{n\star}) \right)_{0,\Sigma} \\
T_7 &= \tau a_{\mathrm{s}} \left( \boldsymbol{\xi}_h^n, \frac{\tau}{\rho_{\mathrm{s}} \epsilon} \boldsymbol{L}_h^{\mathrm{s}} (\boldsymbol{d}^n - \boldsymbol{d}^{n\star}) \right) \\
T_8 &= \rho_s \epsilon \tau \left( \partial_\tau \dot{\boldsymbol{\xi}}_h^n, \frac{\tau}{\rho_{\mathrm{s}} \epsilon} \boldsymbol{L}_h^{\mathrm{s}} (\boldsymbol{d}^n - \boldsymbol{d}^{n\star}) \right)_{0,\Sigma} \\
T_9 &= -\rho_s \epsilon \tau \left( \partial_\tau \dot{\boldsymbol{\xi}}_h^n, \frac{\tau}{\rho_{\mathrm{s}} \epsilon} \boldsymbol{L}_h^{\mathrm{s}} (\boldsymbol{\xi}_h^n - \boldsymbol{\xi}_h^{n\star}) \right)_{0,\Sigma} \\
T_{10} &= -\tau a_{\mathrm{s}} \left( \boldsymbol{\xi}_h^n, \frac{\tau}{\rho_{\mathrm{s}} \epsilon} \boldsymbol{L}_h^{\mathrm{s}} (\boldsymbol{\xi}_h^n - \boldsymbol{\xi}_h^{n\star}) \right)
\end{aligned}
\tag{4.44}
$$

The terms $T_0$, $T_1$, $T_2$, $T_3$ and $T_4$ do not depend on the extrapolation order and they can be estimated using the relations reported in Proposition 3.3.2 or the inf-sup condition for the bilinear form $c(\cdot,\cdot)$. We observe that we found the same terms in the analysis of the monolithic algorithm performed in Theorem 3.3.7. We report here the final estimates and refer to (3.47), (3.48), (3.49) and (3.50) respectively, for the details.

$$
T_0 \leq \frac{\delta_0 \tau}{2} \|\boldsymbol{\xi}_h^n\|_{\mathrm{s}}^2 + \frac{\tau^2}{2\delta_0} \|\partial_{tt} \boldsymbol{d}\|^2_{L^2(t_{n-1},t_n;\boldsymbol{W})}, \tag{4.45}
$$

$$
T_1 \leq \frac{\rho_{\mathrm{f}} \tau^2}{2\delta_1} \|\partial_{tt} \boldsymbol{u}\|^2_{L^2(t_{n-1},t_n;L^2(\Omega)^d)} + \frac{\rho_{\mathrm{f}}}{2\delta_1} \|\partial_t \boldsymbol{\theta}_\pi\|^2_{L^2(t_{n-1},t_n;L^2(\Omega)^d)} + \rho_{\mathrm{f}} \tau \delta_1 \|\boldsymbol{\theta}_h^n\|^2_{0,\Omega}, \tag{4.46}
$$

$$
T_2 \leq \frac{\rho_s \epsilon \tau^2}{2\delta_2} \|\partial_{tt} \dot{\boldsymbol{d}}\|^2_{L^2(t_{n-1},t_n;L^2(\Sigma)^d)} + \frac{\rho_s \epsilon}{2\delta_2} \|\partial_t \dot{\boldsymbol{\xi}}_\pi\|^2_{L^2(t_{n-1},t_n;L^2(\Sigma)^d)} + \rho_s \epsilon \tau \delta_2 \|\dot{\boldsymbol{\xi}}_h^n\|^2_{0,\Sigma}, \tag{4.47}
$$

$$
T_3 \leq \frac{\tau}{\delta_3} \|\boldsymbol{\omega}_\pi\|^2_{L^\infty(0,T;\boldsymbol{\Lambda})} + \tau \delta_3 C_K \|\boldsymbol{\epsilon}(\boldsymbol{\theta}_h^n)\|^2_{,\Omega}, \tag{4.48}
$$

$$
\begin{aligned}
T_4 \leq & C\delta_4 \|\partial_t \boldsymbol{\theta}_\pi\|^2_{L^2(t_{n-1},t_n;L^2(\Omega)^d)} + \frac{C}{\delta_4}\left[\|\dot{\boldsymbol{\xi}}_\pi\|^2_{L^\infty(0,T;H^{\frac{1}{2}}(\Sigma)^d)} + \|\boldsymbol{\theta}_\pi\|^2_{L^\infty(0,T;H^1(\Omega)^d)}\right] \\
& + C\tau\delta_4\|\boldsymbol{\omega}_\pi\|^2_{L^\infty(0,T;\boldsymbol{\Lambda})} + \frac{C\tau}{\delta_4}\left[\|\dot{\boldsymbol{\xi}}_\pi\|^2_{L^\infty(0,T;H^{\frac{1}{2}}(\Sigma)^d)} + \|\boldsymbol{\theta}_\pi\|^2_{L^\infty(0,T;H^1(\Omega)^d)}\right] \\
& + C\tau^2\delta_4\|\partial_{tt}\boldsymbol{u}\|^2_{L^2(t_{n-1},t_n;L^2(\Omega)^d)} + C\delta_4\|\tau\partial_\tau\boldsymbol{\theta}_h^n\|^2_{0,\Omega} + C\tau\delta_4\|\boldsymbol{\epsilon}(\boldsymbol{\theta}_h^n)\|^2_{0,\Omega} + \tau C\delta_4|\varphi_h^n|^2_{s_h},
\end{aligned}
\tag{4.49}
$$

with $\delta_0 > 0$, $\delta_1 > 0$, $\delta_2 > 0$, $\delta_3 > 0$ and $\delta_4 > 0$ to be fixed.

In order to estimate the terms $T_5 - T_{10}$ in (4.43) we introduce the extrapolations.

- **Scheme with** $r = 0$ implies $\boldsymbol{d}^{n\star} = \boldsymbol{0}$ e $\boldsymbol{\xi}_h^{n\star} = \boldsymbol{0}$;
  In the following we obtain the estimates for the terms from $T_5$ to $T_{10}$; such estimates are achieved using results of Proposition 3.3.2, in particular for terms $T_5$ and $T_6$ we have

$$
\begin{aligned}
T_5 = & \rho_s\epsilon\tau\left((\partial_\tau - \partial_t)\dot{\boldsymbol{d}}^n, \frac{\tau}{\rho_s\epsilon}\boldsymbol{L}_h^{\mathrm{s}}\boldsymbol{\xi}_h^n\right)_{0,\Sigma} - \rho_s\epsilon\tau\left(\partial_\tau\dot{\boldsymbol{\xi}}_\pi^n, \frac{\tau}{\rho_s\epsilon}\boldsymbol{L}_h^{\mathrm{s}}\boldsymbol{\xi}_h^n\right)_{0,\Sigma} \\
& \leq \frac{\tau^4}{2\delta_5}\|\partial_{tt}\dot{\boldsymbol{d}}\|^2_{L^2(t_{n-1},t_n;\boldsymbol{W})} + \frac{\tau^2}{2\delta_5}\|\partial_t\dot{\boldsymbol{\xi}}_\pi\|^2_{L^2(t_{n-1},t_n;\boldsymbol{W})} + \tau\delta_5\|\boldsymbol{\xi}_h^n\|^2_{\mathrm{s}},
\end{aligned}
\tag{4.50}
$$

$$
\begin{aligned}
T_6 = & -\rho_s\epsilon\tau\left((\partial_\tau - \partial_t)\dot{\boldsymbol{d}}^n, \frac{\tau}{\rho_s\epsilon}\boldsymbol{L}_h^{\mathrm{s}}\boldsymbol{d}^n\right)_{0,\Sigma} + \rho_s\epsilon\tau\left(\partial_\tau\dot{\boldsymbol{\xi}}_\pi^n, \frac{\tau}{\rho_s\epsilon}\boldsymbol{L}_h^{\mathrm{s}}\boldsymbol{d}^n\right)_{0,\Sigma} \\
& \leq \frac{\tau^3}{2\delta_6}\|\partial_{tt}\dot{\boldsymbol{d}}\|^2_{L^2(t_{n-1},t_n;L^2(\Sigma)^d)} + \frac{\tau}{2\delta_6}\|\partial_t\dot{\boldsymbol{\xi}}_\pi\|^2_{L^2(t_{n-1},t_n;L^2(\Sigma)^d)} + \tau^2\delta_6\|\boldsymbol{L}_h^{\mathrm{s}}\boldsymbol{d}^n\|^2_{0,\Sigma}.
\end{aligned}
\tag{4.51}
$$

The estimate for the term $T_7$ uses the definition of the solid discrete operator $\boldsymbol{L}_h^{\mathrm{s}}$ given in (3.25)

$$
\begin{aligned}
T_7 = \tau a_{\mathrm{s}}\left(\boldsymbol{\xi}_h^n, \frac{\tau}{\rho_s\epsilon}\boldsymbol{L}_h^{\mathrm{s}}\boldsymbol{d}^n\right) = & \frac{\tau^2}{\rho_s\epsilon}\left(\boldsymbol{L}_h^{\mathrm{s}}\boldsymbol{\xi}_h^n, \boldsymbol{L}_h^{\mathrm{s}}\boldsymbol{d}^n\right)_{0,\Sigma} \\
& \leq \frac{\tau^2\delta_7}{2\rho_s\epsilon}\|\boldsymbol{L}_h^{\mathrm{s}}\boldsymbol{d}^n\|^2_{0,\Sigma} + \frac{\tau^2}{2\rho_s\epsilon\delta_7}\|\boldsymbol{L}_h^{\mathrm{s}}\boldsymbol{\xi}_h^n\|^2_{0,\Sigma}.
\end{aligned}
\tag{4.52}
$$

The estimates for terms $T_8$ and $T_9$ are obtained by a straight application of Cauchy-Schwarz inequality

$$
T_8 = \rho_s\epsilon\tau\left(\partial_\tau\dot{\boldsymbol{\xi}}_h^n, \frac{\tau}{\rho_s\epsilon}\boldsymbol{L}_h^{\mathrm{s}}\boldsymbol{d}^n\right)_{0,\Sigma} \leq \frac{\tau^2}{2\delta_8}\|\boldsymbol{L}_h^{\mathrm{s}}\boldsymbol{d}^n\|^2_{0,\Sigma} + \frac{\tau^2\delta_8}{2}\|\partial_\tau\dot{\boldsymbol{\xi}}_h^n\|^2_{0,\Sigma},
\tag{4.53}
$$

$$T_9 = -\rho_s\epsilon\tau\left(\partial_\tau\dot{\boldsymbol{\xi}}_h^n, \frac{\tau}{\rho_s\epsilon}\boldsymbol{L}_h^s\boldsymbol{\xi}_h^n\right)_{0,\Sigma} \geq -\frac{\tau^2}{2\delta_9}\|\boldsymbol{L}_h^s\boldsymbol{\xi}_h^n\|_{0,\Sigma}^2 - \frac{\tau^2\delta_9}{2}\|\partial_\tau\dot{\boldsymbol{\xi}}_h^n\|_{0,\Sigma}^2. \quad (4.54)$$

The estimate of the term $T_{10}$, like the estimate of the term $T_7$ uses the definition of the discrete solid operator (3.25)

$$T_{10} = -\tau a_s\left(\boldsymbol{\xi}_h^n, \frac{\tau}{\rho_s\epsilon}\boldsymbol{L}_h^s\boldsymbol{\xi}_h^n\right) = -\frac{\tau^2}{\rho_s\epsilon}\left(\boldsymbol{L}_h^s\boldsymbol{\xi}_h^n, \boldsymbol{L}_h^s\boldsymbol{\xi}_h^n\right)_{0,\Sigma} = -\frac{\tau^2}{\rho_s\epsilon}\|\boldsymbol{L}_h^s\boldsymbol{\xi}_h^n\|_{0,\Sigma}^2. \quad (4.55)$$

Using the estimates of terms $T_0 - T_{10}$ in (4.43) and collecting together the similar terms, we obtain the following inequality that gives a bound for the discrete error at each instant of time $t_n$ in the case of extrapolation order $r = 0$.

$$\begin{aligned}
&\left(\frac{\rho_f}{2}\|\boldsymbol{\theta}_h^n\|_{0,\Omega}^2 + \frac{\rho_s\epsilon}{2}\|\dot{\boldsymbol{\xi}}_h^n\|_{0,\Sigma}^2 + \frac{1}{2}\|\boldsymbol{\xi}_h^n\|_s^2\right) - \left(\frac{\rho_f}{2}\|\boldsymbol{\theta}_h^{n-1}\|_{0,\Omega}^2 + \frac{\rho_s\epsilon}{2}\|\dot{\boldsymbol{\xi}}_h^{n-1}\|_{0,\Sigma}^2 + \frac{1}{2}\|\boldsymbol{\xi}_h^{n-1}\|_s^2\right)\\
&+\tau^2\left(\left(\frac{\rho_f}{2} - C\delta_4\right)\|\partial_\tau\boldsymbol{\theta}_h^n\|_{0,\Omega}^2 + \left(\frac{\rho_s\epsilon}{2} - \frac{\delta_9}{2} - \frac{\delta_8}{2}\right)\|\partial_\tau\dot{\boldsymbol{\xi}}_h^n\|_{0,\Sigma}^2 + \frac{1}{2}\|\partial_\tau\boldsymbol{\xi}_h^n\|_s^2\right)\\
&+\tau\left[\left(2\mu - \delta_3 C_K - C\delta_4\right)\|\boldsymbol{\epsilon}(\boldsymbol{\theta}_h^n)\|_{0,\Omega}^2 + (1 - C\delta_4)|\varphi_h^n|_{s_h}^2\right]\\
&+\left(\frac{\tau^2}{\rho_s\epsilon} - \frac{\tau^2}{2\delta_9} - \frac{\tau^2}{2\rho_s\epsilon\delta_7}\right)\|\boldsymbol{L}_h^s\boldsymbol{\xi}_h^n\|_{0,\Sigma}^2\\
\leq &\left[\frac{\rho_s\epsilon}{2\delta_2}\|\partial_t\dot{\boldsymbol{\xi}}_\pi\|_{L^2(t_{n-1},t_n;L^2(\Sigma)^d)}^2 + \frac{\rho_f}{2\delta_1}\|\partial_t\boldsymbol{\theta}_\pi\|_{L^2(t_{n-1},t_n;L^2(\Omega)^d)}^2 + C\delta_4\|\partial_t\boldsymbol{\theta}_\pi\|_{L^2(t_{n-1},t_n;L^2(\Omega)^d)}^2\right.\\
&\left.+\frac{C}{\delta_4}\|\dot{\boldsymbol{\xi}}_\pi\|_{L^\infty(0,T;H^{\frac{1}{2}}(\Sigma)^d)}^2 + \frac{C}{\delta_4}\|\boldsymbol{\theta}_\pi\|_{L^\infty(0,T;H^1(\Omega)^d)}^2\right]\\
&+\tau\left[\frac{1}{2\delta_6}\|\partial_t\dot{\boldsymbol{\xi}}_\pi\|_{L^2(t_{n-1},t_n;L^2(\Sigma)^d)}^2 + \frac{1}{2\delta_0}\|\dot{\boldsymbol{\xi}}_\pi\|_{L^\infty(0,T;\boldsymbol{W})}^2 + \frac{1}{\delta_3}\|\boldsymbol{\omega}_\pi\|_{L^\infty(0,T;\boldsymbol{\Lambda})}^2\right.\\
&\left.+C\delta_4\|\boldsymbol{\omega}_\pi\|_{L^\infty(0,T;\boldsymbol{\Lambda})}^2 + \frac{C}{\delta_4}\|\dot{\boldsymbol{\xi}}_\pi\|_{L^\infty(0,T;H^{\frac{1}{2}}(\Sigma)^d)}^2 + \frac{C}{\delta_4}\|\boldsymbol{\theta}_\pi\|_{L^\infty(0,T;H^1(\Omega)^d)}^2\right]\\
&+\tau^2\left[\left(\frac{1}{2\delta_8} + \frac{\delta_7}{2\rho_s\epsilon} + \delta_6\right)\|\boldsymbol{L}_h^s\boldsymbol{d}^n\|_{0,\Sigma}^2 + \frac{\rho_s\epsilon\tau^2}{2\delta_2}\|\partial_{tt}\dot{\boldsymbol{d}}\|_{L^2(t_{n-1},t_n;L^2(\Sigma)^d)}^2\right.\\
&+\frac{1}{2\delta_0}\|\partial_{tt}\boldsymbol{d}\|_{L^2(t_{n-1},t_n;\boldsymbol{W})}^2 + \frac{\rho_f}{2\delta_1}\|\partial_{tt}\boldsymbol{u}\|_{L^2(t_{n-1},t_n;L^2(\Omega)^d)}^2\\
&\left.+\frac{1}{2\delta_5}\|\partial_t\dot{\boldsymbol{\xi}}_\pi\|_{L^2(t_{n-1},t_n;\boldsymbol{W})}^2 + C\delta_4\|\partial_{tt}\boldsymbol{u}\|_{L^2(t_{n-1},t_n;L^2(\Omega)^d)}^2\right]\\
&+\frac{\tau^3}{2\delta_6}\|\partial_{tt}\dot{\boldsymbol{d}}\|_{L^2(t_{n-1},t_n;L^2(\Sigma)^d)}^2 + \frac{\tau^4}{2\delta_5}\|\partial_{tt}\dot{\boldsymbol{d}}\|_{L^2(t_{n-1},t_n;\boldsymbol{W})}^2\\
&+\left[\rho_f\tau\delta_1\|\boldsymbol{\theta}_h^n\|_{0,\Omega}^2 + \rho_s\epsilon\tau\delta_2\|\dot{\boldsymbol{\xi}}_h^n\|_{0,\Sigma}^2 + \left(\tau\delta_5 + \frac{\delta_0\tau}{2}\right)\|\boldsymbol{\xi}_h^n\|_s^2\right]. \quad (4.56)
\end{aligned}$$

In order to obtain the cumulative error after "$n$" time steps, we sum over $m = 1,\dots,n$, we set the parameters in such a way that all the quantities at the left hand side result positives

$$\delta_0 = \frac{1}{4}, \qquad \delta_1 = \delta_2 = \frac{1}{2}, \qquad \delta_3 = \frac{\mu}{C_K}, \qquad \delta_4 = \min\left\{\frac{1}{C}, \frac{\rho_{\rm f}}{2C}, \frac{\mu}{2C}\right\}, \qquad \delta_5 = \frac{1}{4},$$
$$\delta_6 = 1, \qquad \delta_7 = \frac{10}{8}, \qquad \delta_8 = \frac{\rho_s\epsilon}{10}, \qquad \delta_9 = \frac{13\rho_s\epsilon}{15}.$$

Using this values we obtain the following expression of the cumulative error at time step "$n$"

$$\begin{aligned}
&\underbrace{\left(\frac{\rho_{\rm f}}{2}\|\boldsymbol{\theta}_h^n\|_{0,\Omega}^2 + \frac{\rho_s\epsilon}{2}\|\dot{\boldsymbol{\xi}}_h^n\|_{0,\Sigma}^2 + \frac{1}{2}\|\boldsymbol{\xi}_h^n\|_{\rm s}^2\right)}_{\mathscr{E}_h^n} - \underbrace{\left(\frac{\rho_{\rm f}}{2}\|\boldsymbol{\theta}_h^0\|_{0,\Omega}^2 + \frac{\rho_s\epsilon}{2}\|\dot{\boldsymbol{\xi}}_h^0\|_{0,\Sigma}^2 + \frac{1}{2}\|\boldsymbol{\xi}_h^0\|_{\rm s}^2\right)}_{\mathscr{E}_h^0} \\
&\quad + \tau\sum_{m=1}^{n}\left[C\|\boldsymbol{\epsilon}(\boldsymbol{\theta}_h^m)\|_{0,\Omega}^2 + C|\varphi_h^m|_{s_h}^2\right] \\
&\quad + \tau^2\sum_{m=1}^{n}\left[\frac{\rho_{\rm f}}{2}\|\partial_\tau\boldsymbol{\theta}_h^m\|_{0,\Omega}^2 + C\|\partial_\tau\dot{\boldsymbol{\xi}}_h^m\|_{0,\Sigma}^2 + \frac{1}{2}\|\partial_\tau\boldsymbol{\xi}_h^m\|_{\rm s}^2 + C\|\boldsymbol{L}_h^{\rm s}\boldsymbol{\xi}_h^m\|_{0,\Sigma}^2\right] \\
&\le \underbrace{\sum_{m=1}^{n} C_m^0}_{A_0} + \tau\underbrace{\sum_{m=1}^{n} C_m^1}_{A_1} + \tau^2\underbrace{\sum_{m=1}^{n} C_m^2}_{A_2} + \tau^3\underbrace{\sum_{m=1}^{n} C_m^3}_{A_3} + \tau^4\underbrace{\sum_{m=1}^{n} C_m^4}_{A_4} \\
&\quad + \tau\sum_{m=1}^{n}\underbrace{\left[\frac{\rho_{\rm f}}{2}\|\boldsymbol{\theta}_h^m\|_{0,\Omega}^2 + \frac{\rho_s\epsilon}{2}\|\dot{\boldsymbol{\xi}}_h^m\|_{0,\Sigma}^2 + \frac{1}{2}\|\boldsymbol{\xi}_h^m\|_{\rm s}^2\right]}_{\mathscr{E}_h^m},
\end{aligned} \tag{4.57}$$

where the terms in the right hand side of the inequality are defined ad analyzed in the following using results of Proposition 3.3.6.

- The terms $A_0$ can be bounded as follows

$$\begin{aligned}
A_0 \overset{def}{=} C\sum_{n=1}^{m}&\Big(\|\partial_t\dot{\boldsymbol{\xi}}_\pi\|_{L^2(t_{n-1},t_n;L^2(\Sigma)^d)}^2 + \|\partial_t\boldsymbol{\theta}_\pi\|_{L^2(t_{n-1},t_n;L^2(\Omega)^d)}^2 \\
&+ \|\dot{\boldsymbol{\xi}}_\pi\|_{L^\infty(0,T;H^{\frac{1}{2}}(\Sigma)^d)}^2 + \|\boldsymbol{\theta}_\pi\|_{L^\infty(0,T;H^1(\Omega)^d)}^2\Big) \\
\le C&\Big(\|\partial_t\dot{\boldsymbol{\xi}}_\pi\|_{L^2(0,T;L^2(\Sigma)^d)}^2 + \|\partial_t\boldsymbol{\theta}_\pi\|_{L^2(0,T;L^2(\Omega)^d)}^2 \\
&+ \|\dot{\boldsymbol{\xi}}_\pi\|_{L^\infty(0,T;H^{\frac{1}{2}}(\Sigma)^d)}^2 + \|\boldsymbol{\theta}_\pi\|_{L^\infty(0,T;H^1(\Omega)^d)}^2\Big) \\
\approx\;& \mathcal{O}(h_s^{2m+2}) + \mathcal{O}(h_f^{2l+2}) + \mathcal{O}(h_s^{2m+1}) + \mathcal{O}(h_f^{2l}).
\end{aligned} \tag{4.58}$$

- The term $A_1$ is defined and bounded as follows

$$
\begin{aligned}
A_1 &\stackrel{def}{=} C \sum_{n=1}^{m} \Big( \|\partial_t \dot{\boldsymbol{\xi}}_\pi\|^2_{L^2(t_{n-1},t_n;L^2(\Sigma)^d)} + \|\dot{\boldsymbol{\xi}}_\pi\|^2_{L^\infty(0,T;\boldsymbol{W})} \\
&\qquad + \|\boldsymbol{\omega}_\pi\|^2_{L^\infty(0,T;\boldsymbol{\Lambda})} + \|\dot{\boldsymbol{\xi}}_\pi\|^2_{L^\infty(0,T;H^{\frac{1}{2}}(\Sigma)^d)} + \|\boldsymbol{\theta}_\pi\|^2_{L^\infty(0,T;H^1(\Omega)^d)} \Big) \\
&\leq C\Big( \|\partial_t \dot{\boldsymbol{\xi}}_\pi\|^2_{L^2(0,T;L^2(\Sigma)^d)} + \|\dot{\boldsymbol{\xi}}_\pi\|^2_{L^\infty(0,T;\boldsymbol{W})} \\
&\qquad + \|\boldsymbol{\omega}_\pi\|^2_{L^\infty(0,T;\boldsymbol{\Lambda})} + \|\dot{\boldsymbol{\xi}}_\pi\|^2_{L^\infty(0,T;H^{\frac{1}{2}}(\Sigma)^d)} + \|\boldsymbol{\theta}_\pi\|^2_{L^\infty(0,T;H^1(\Omega)^d)} \Big) \\
&\approx \mathcal{O}(h_s^{2m+2}) + \mathcal{O}(h_s^{2m}) + \mathcal{O}(h_s^{2z+1}) + \mathcal{O}(h_s^{2m+1}) + \mathcal{O}(h_f^{2l}).
\end{aligned}
\tag{4.59}
$$

– The term $A_2$ is defined and bounded in the following using also the properties of the solid operator $\|\boldsymbol{L}_h^{\mathrm{s}}\boldsymbol{d}^n\|_{0,\Sigma} \leq C\|\boldsymbol{L}^{\mathrm{s}}\boldsymbol{d}^n\|_{0,\Sigma}$

$$
\begin{aligned}
A_2 &\stackrel{def}{=} C \sum_{n=1}^{m} \Big( \|\boldsymbol{L}_h^{\mathrm{s}}\boldsymbol{d}^n\|^2_{0,\Sigma} + \|\partial_{tt}\dot{\boldsymbol{d}}\|^2_{L^2(t_{n-1},t_n;L^2(\Sigma)^d)} \\
&\qquad + \|\partial_{tt}\boldsymbol{d}\|^2_{L^2(t_{n-1},t_n;\boldsymbol{W})} + \|\partial_{tt}\boldsymbol{u}\|^2_{L^2(t_{n-1},t_n;L^2(\Omega)^d)} \\
&\qquad + \|\partial_t \dot{\boldsymbol{\xi}}_\pi\|^2_{L^2(t_{n-1},t_n;\boldsymbol{W})} + \|\partial_{tt}\boldsymbol{u}\|^2_{L^2(t_{n-1},t_n;L^2(\Omega)^d)} \Big) \\
&\leq C\Big( \|\boldsymbol{L}^{\mathrm{s}}\boldsymbol{d}\|^2_{L^\infty(0,T;L^2(\Sigma)^d)} + \|\partial_{tt}\dot{\boldsymbol{d}}\|^2_{L^2(0,T;L^2(\Sigma)^d)} \\
&\qquad + \|\partial_{tt}\boldsymbol{d}\|^2_{L^2(0,T;\boldsymbol{W})} + \|\partial_{tt}\boldsymbol{u}\|^2_{L^2(0,T;L^2(\Omega)^d)} \\
&\qquad + \|\partial_t \dot{\boldsymbol{\xi}}_\pi\|^2_{L^2(0,T;\boldsymbol{W})} + \|\partial_{tt}\boldsymbol{u}\|^2_{L^2(0,T;L^2(\Omega)^d)} \Big).
\end{aligned}
\tag{4.60}
$$

– The terms $A_3$ and $A_4$ are bounded as follows

$$
A_3 \stackrel{def}{=} C \sum_{n=1}^{m} \|\partial_{tt}\dot{\boldsymbol{d}}\|^2_{L^2(t_{n-1},t_n;L^2(\Sigma)^d)} \leq C\|\partial_{tt}\dot{\boldsymbol{d}}\|^2_{L^2(0,T;L^2(\Sigma)^d)} \tag{4.61}
$$

$$
A_4 \stackrel{def}{=} C \sum_{n=1}^{m} \|\partial_{tt}\dot{\boldsymbol{d}}\|^2_{L^2(t_{n-1},t_n;\boldsymbol{W})} \leq C\|\partial_{tt}\dot{\boldsymbol{d}}\|^2_{L^2(0,T;\boldsymbol{W})} \tag{4.62}
$$

In sum the error estimate (4.57) can be written as follows

$$
\mathcal{E}_h^n + \tau\mathcal{D}_h^n \leq \tau \sum_{m=1}^{n} \mathcal{E}_h^m + A_0 + \tau A_1 + \tau^2 A_2 + \tau^3 A_3 + \tau^4 A_4 \tag{4.63}
$$

where $\tau\mathcal{D}$ represents the dissipation term and is given by

$$
\begin{aligned}
\tau\mathcal{D}_h^n = \tau \sum_{m=1}^{n} &\Big[ C\|\boldsymbol{\epsilon}(\boldsymbol{\theta}_h^m)\|^2_{0,\Omega} + C|\varphi_h^m|^2_{s_h} \Big] \\
&+ \tau^2 \sum_{m=1}^{n} \left[ \frac{\rho_{\mathrm{f}}}{2}\|\partial_\tau \boldsymbol{\theta}_h^m\|^2_{0,\Omega} + C\|\partial_\tau \dot{\boldsymbol{\xi}}_h^m\|^2_{0,\Sigma} + \frac{1}{2}\|\partial_\tau \boldsymbol{\xi}_h^m\|^2_{\mathrm{s}} + C\|\boldsymbol{L}_h^{\mathrm{s}}\boldsymbol{\xi}_h^m\|^2_{0,\Sigma} \right].
\end{aligned}
$$

Introducing the convergence rates with respect to $h_f$ and $h_s$ of the term $A_0$,$A_1$,$A_2$ given in (4.58)-(4.60), we have the desired result applying the discrete Gronwall Lemma 3.3.1, since by hypothesis we have $\tau < 1$

$$\mathscr{E}_h^n + \tau \mathscr{D}_h^n \lesssim \exp\left(\sum_{m=1}^{n} \frac{\tau}{1-\tau}\right)\Big[\underbrace{\mathcal{O}(h_f^{2l}) + \mathcal{O}(h_s^{2m+1})}_{A_0} \tag{4.64}$$

$$+ \tau\Big(\underbrace{\mathcal{O}(h_s^{2m}) + \mathcal{O}(h_s^{2z+1}) + \mathcal{O}(h_f^{2l})}_{A_1}\Big) + \tau^2 A_2 + \tau^3 A_3 + \tau^4 A_4\Big]. \tag{4.65}$$

- **Scheme with** $r = 1$ implies $\boldsymbol{d}^{n\star} = \boldsymbol{d}^{n-1}$ e $\boldsymbol{\xi}_h^{n\star} = \boldsymbol{\xi}_h^{n-1}$.

In the following we obtain the estimates for the terms from $T_5$ to $T_{10}$; such estimates are obtained using, results of proposition 3.3.2, in particular for terms $T_5$ and $T_6$ we have

$$\begin{aligned} T_5 &= \rho_s \epsilon \tau \left((\partial_\tau - \partial_t)\dot{\boldsymbol{d}}^n, \frac{\tau}{\rho_s \epsilon}\boldsymbol{L}_h^{s}(\boldsymbol{\xi}_h^n - \boldsymbol{\xi}_h^{n-1})\right)_{0,\Sigma} - \rho_s \epsilon \tau \left(\partial_\tau \dot{\boldsymbol{\xi}}_\pi^n, \frac{\tau}{\rho_s \epsilon}\boldsymbol{L}_h^{s}(\boldsymbol{\xi}_h^n - \boldsymbol{\xi}_h^{n-1})\right)_{0,\Sigma} \\ &\le \frac{\tau^3 \delta_5}{2}\|\partial_{tt}\dot{\boldsymbol{d}}\|^2_{L^2(t_{n-1},t_n;L^2(\Sigma)^d)} + \frac{\tau \delta_5}{2}\|\partial_t \dot{\boldsymbol{\xi}}_\pi\|^2_{L^2(t_{n-1},t_n;L^2(\Sigma)^d)} \\ &\quad + \frac{\tau^2}{\delta_5}\|\boldsymbol{L}_h^{s}(\boldsymbol{\xi}_h^n - \boldsymbol{\xi}_h^{n-1})\|^2_{0,\Sigma}, \end{aligned} \tag{4.66}$$

$$\begin{aligned} T_6 &= -\rho_s \epsilon \tau \left((\partial_\tau - \partial_t)\dot{\boldsymbol{d}}^n, \frac{\tau}{\rho_s \epsilon}\boldsymbol{L}_h^{s}(\boldsymbol{d}^n - \boldsymbol{d}^{n-1})\right)_{0,\Sigma} + \rho_s \epsilon \tau \left(\partial_\tau \dot{\boldsymbol{\xi}}_\pi^n, \frac{\tau}{\rho_s \epsilon}\boldsymbol{L}_h^{s}(\boldsymbol{d}^n - \boldsymbol{d}^{n-1})\right)_{0,\Sigma} \\ &\le \frac{\tau^3}{2\delta_6}\|\partial_{tt}\dot{\boldsymbol{d}}\|^2_{L^2(t_{n-1},t_n;L^2(\Sigma)^d)} + \frac{\tau}{2\delta_6}\|\partial_t \dot{\boldsymbol{\xi}}_\pi\|^2_{L^2(t_{n-1},t_n;L^2(\Sigma)^d)} \\ &\quad + \tau^2 \delta_6 \|\boldsymbol{L}_h^{s}(\boldsymbol{d}^n - \boldsymbol{d}^{n-1})\|^2_{0,\Sigma}. \end{aligned} \tag{4.67}$$

The estimate for the term $T_7$ uses the definition of the solid discrete operator $\boldsymbol{L}_h^{s}$ given in (3.25)

$$T_7 = \tau a_{\rm s}\left(\boldsymbol{\xi}_h^n, \frac{\tau}{\rho_{\rm s}\epsilon}\boldsymbol{L}_h^{\rm s}(\boldsymbol{d}^n - \boldsymbol{d}^{n-1})\right) = \frac{\tau^2}{\rho_{\rm s}\epsilon}\left(\boldsymbol{L}_h^{\rm s}\boldsymbol{\xi}_h^n, \boldsymbol{L}_h^{\rm s}(\boldsymbol{d}^n - \boldsymbol{d}^{n-1})\right)_{0,\Sigma}$$
$$\leq \frac{\tau^2}{2\rho_{\rm s}\epsilon\delta_7}\|\boldsymbol{L}_h^{\rm s}(\boldsymbol{d}^n - \boldsymbol{d}^{n-1})\|_{0,\Sigma}^2 + \frac{\tau^2\delta_7}{2\rho_{\rm s}\epsilon}\|\boldsymbol{L}_h^{\rm s}\boldsymbol{\xi}_h^n\|_{0,\Sigma}^2. \quad (4.68)$$

The estimates for terms $T_8$ and $T_9$ are obtained by a straight application of Cauchy-Schwarz inequality

$$T_8 = \rho_s\epsilon\tau\left(\partial_\tau\dot{\boldsymbol{\xi}}_h^n, \frac{\tau}{\rho_{\rm s}\epsilon}\boldsymbol{L}_h^{\rm s}(\boldsymbol{d}^n - \boldsymbol{d}^{n-1})\right)_{0,\Sigma} \leq \frac{\tau^2}{2\delta_8}\|\boldsymbol{L}_h^{\rm s}(\boldsymbol{d}^n - \boldsymbol{d}^{n-1})\|_{0,\Sigma}^2 + \frac{\tau^2\delta_8}{2}\|\partial_\tau\dot{\boldsymbol{\xi}}_h^n\|_{0,\Sigma}^2, \quad (4.69)$$

$$T_9 = -\rho_s\epsilon\tau\left(\partial_\tau\dot{\boldsymbol{\xi}}_h^n, \frac{\tau}{\rho_{\rm s}\epsilon}\boldsymbol{L}_h^{\rm s}(\boldsymbol{\xi}_h^n - \boldsymbol{\xi}_h^{n-1})\right)_{0,\Sigma} = -\tau^2\left(\tau\partial_\tau\dot{\boldsymbol{\xi}}_h^n, \boldsymbol{L}_h^{\rm s}(\dot{\boldsymbol{\xi}}_h^n)\right)_{0,\Sigma}$$
$$= -\frac{\tau^2}{2}\left(\|\dot{\boldsymbol{\xi}}_h^n\|_s^2 - \|\dot{\boldsymbol{\xi}}_h^{n-1}\|_s^2 + \|\dot{\boldsymbol{\xi}}_h^n - \dot{\boldsymbol{\xi}}_h^{n-1}\|_s^2\right) \quad (4.70)$$

In the estimate of the term $T_{10}$ we use once again the definition of the solid discrete operator (3.25)

$$T_{10} = -\tau a_{\rm s}\left(\boldsymbol{\xi}_h^n, \frac{\tau}{\rho_{\rm s}\epsilon}\boldsymbol{L}_h^{\rm s}(\boldsymbol{\xi}_h^n - \boldsymbol{\xi}_h^{n-1})\right) = -\frac{\tau^2}{\rho_{\rm s}\epsilon}\left(\boldsymbol{L}_h^{\rm s}\boldsymbol{\xi}_h^n, \boldsymbol{L}_h^{\rm s}(\boldsymbol{\xi}_h^n - \boldsymbol{\xi}_h^{n-1})\right)_{0,\Sigma}$$
$$= -\frac{\tau^2}{2\rho_{\rm s}\epsilon}\left(\|\boldsymbol{L}_h^{\rm s}\boldsymbol{\xi}_h^n\|_{0,\Sigma}^2 - \|\boldsymbol{L}_h^{\rm s}\boldsymbol{\xi}_h^{n-1}\|_{0,\Sigma}^2 + \|\boldsymbol{L}_h^{\rm s}(\boldsymbol{\xi}_h^n - \boldsymbol{\xi}_h^{n-1})\|_{0,\Sigma}^2\right). \quad (4.71)$$

Using the estimates of terms from $T_0$ to $T_{10}$ in (4.43) and collecting together the similar terms, we obtain the following inequality that gives a bound for the discrete error at each instant of time $t_n$ in the case of extrapolation order $r = 1$.

$$
\begin{aligned}
&\left(\frac{\rho_{\mathrm{f}}}{2}\|\boldsymbol{\theta}_h^n\|_{0,\Omega}^2 + \frac{\rho_s\epsilon}{2}\|\dot{\boldsymbol{\xi}}_h^n\|_{0,\Sigma}^2 + \frac{1}{2}\|\boldsymbol{\xi}_h^n\|_{\mathrm{s}}^2 + \frac{\tau^2}{2\rho_{\mathrm{s}}\epsilon}\|\boldsymbol{L}_h^{\mathrm{s}}\boldsymbol{\xi}_h^n\|_{0,\Sigma}^2 + \frac{\tau^2}{2}\|\dot{\boldsymbol{\xi}}_h^n\|_s^2\right)\\
&-\left(\frac{\rho_{\mathrm{f}}}{2}\|\boldsymbol{\theta}_h^{n-1}\|_{0,\Omega}^2 + \frac{\rho_s\epsilon}{2}\|\dot{\boldsymbol{\xi}}_h^{n-1}\|_{0,\Sigma}^2 + \frac{1}{2}\|\boldsymbol{\xi}_h^{n-1}\|_{\mathrm{s}}^2 + \frac{\tau^2}{2\rho_{\mathrm{s}}\epsilon}\|\boldsymbol{L}_h^{\mathrm{s}}\boldsymbol{\xi}_h^{n-1}\|_{0,\Sigma}^2 + \frac{\tau^2}{2}\|\dot{\boldsymbol{\xi}}_h^{n-1}\|_s^2\right)\\
&+\tau^2\left(\left(\frac{\rho_{\mathrm{f}}}{2} - C\delta_4\right)\|\partial_\tau\boldsymbol{\theta}_h^n\|_{0,\Omega}^2 + \left(\frac{\rho_s\epsilon}{2} - \frac{\delta_8}{2}\right)\|\partial_\tau\dot{\boldsymbol{\xi}}_h^n\|_{0,\Sigma}^2 + \frac{1}{2}\|\partial_\tau\boldsymbol{\xi}_h^n\|_{\mathrm{s}}^2 + \frac{1}{2}\|\partial_\tau\dot{\boldsymbol{\xi}}_h^n\|_s^2\right)\\
&+\tau\left((2\mu - C\delta_4 - \delta_3 C_K)\|\boldsymbol{\epsilon}(\boldsymbol{\theta}_h^n)\|_{0,\Omega}^2 + \tau(1 - C\delta_4)|\varphi_h^n|_{s_h}^2\right)\\
&+\left(\frac{\tau^2}{2\rho_{\mathrm{s}}\epsilon} - \frac{\tau^2}{\delta_5}\right)\|\boldsymbol{L}_h^{\mathrm{s}}(\boldsymbol{\xi}_h^n - \boldsymbol{\xi}_h^{n-1})\|_{0,\Sigma}^2\\
&\le\\
&\Big[\frac{\rho_s\epsilon}{2\delta_2}\|\partial_t\dot{\boldsymbol{\xi}}_\pi\|_{L^2(t_{n-1},t_n;L^2(\Sigma)^d)}^2 + \frac{\rho_{\mathrm{f}}}{2\delta_1}\|\partial_t\boldsymbol{\theta}_\pi\|_{L^2(t_{n-1},t_n;L^2(\Omega)^d)}^2 + C\delta_4\|\partial_t\boldsymbol{\theta}_\pi\|_{L^2(t_{n-1},t_n;L^2(\Omega)^d)}^2\\
&+\frac{C}{\delta_4}\|\dot{\boldsymbol{\xi}}_\pi\|_{L^\infty(0,T;H^{\frac{1}{2}}(\Sigma)^d)}^2 + \frac{C}{\delta_4}\|\boldsymbol{\theta}_\pi\|_{L^\infty(0,T;H^1(\Omega)^d)}^2\Big]\\
&+\tau\Big[\frac{1}{2\delta_6}\|\partial_t\dot{\boldsymbol{\xi}}_\pi\|_{L^2(t_{n-1},t_n;L^2(\Sigma)^d)}^2 + \frac{\delta_5}{2}\|\partial_t\dot{\boldsymbol{\xi}}_\pi\|_{L^2(t_{n-1},t_n;L^2(\Sigma)^d)}^2 + \frac{1}{2\delta_0}\|\dot{\boldsymbol{\xi}}_\pi\|_{L^\infty(0,T;\boldsymbol{W})}^2\\
&+\frac{1}{\delta_3}\|\boldsymbol{\omega}_\pi\|_{L^\infty(0,T;\boldsymbol{\Lambda})}^2 + C\delta_4\|\boldsymbol{\omega}_\pi\|_{L^\infty(0,T;\boldsymbol{\Lambda})}^2 + \frac{C}{\delta_4}\|\dot{\boldsymbol{\xi}}_\pi\|_{L^\infty(0,T;H^{\frac{1}{2}}(\Sigma)^d)}^2 + \frac{C}{\delta_4}\|\boldsymbol{\theta}_\pi\|_{L^\infty(0,T;H^1(\Omega)^d)}^2\Big]\\
&+\tau^2\Big[\left(\frac{1}{2\delta_8} + \frac{1}{2\rho_{\mathrm{s}}\epsilon\delta_7} + \delta_6\right)\|\boldsymbol{L}_h^{\mathrm{s}}(\boldsymbol{d}^n - \boldsymbol{d}^{n-1})\|_{0,\Sigma}^2 + \frac{1}{2\delta_0}\|\partial_{tt}\boldsymbol{d}\|_{L^2(t_{n-1},t_n;\boldsymbol{W})}^2 + \frac{\rho_s\epsilon}{2\delta_2}\|\partial_{tt}\dot{\boldsymbol{d}}\|_{L^2(t_{n-1},t_n;L^2(\Sigma)^d)}^2\\
&+\frac{\rho_{\mathrm{f}}}{2\delta_1}\|\partial_{tt}\boldsymbol{u}\|_{L^2(t_{n-1},t_n;L^2(\Omega)^d)}^2 + C\delta_4\|\partial_{tt}\boldsymbol{u}\|_{L^2(t_{n-1},t_n;L^2(\Omega)^d)}^2\Big]\\
&+\frac{\tau^3\delta_5}{2}\|\partial_{tt}\dot{\boldsymbol{d}}\|_{L^2(t_{n-1},t_n;L^2(\Sigma)^d)}^2\\
&+\left[\tau\rho_{\mathrm{f}}\delta_1\|\boldsymbol{\theta}_h^n\|_{0,\Omega}^2 + \rho_s\epsilon\tau\delta_2\|\dot{\boldsymbol{\xi}}_h^n\|_{0,\Sigma}^2 + \frac{\tau\delta_0}{2}\|\boldsymbol{\xi}_h^n\|_{\mathrm{s}}^2 + \frac{\tau^2\delta_7}{2\rho_{\mathrm{s}}\epsilon}\|\boldsymbol{L}_h^{\mathrm{s}}\boldsymbol{\xi}_h^n\|_{0,\Sigma}^2\right].
\end{aligned}
\tag{4.72}
$$

In order to obtain the cumulative error after "$n$" time steps, we sum over $m = 1\dots n$, we set the parameters in such a way that all the quantities at the left hand side result positives

$$
\delta_0 = 1, \qquad \delta_1 = \delta_2 = \frac{1}{2}, \qquad \delta_3 = \frac{\mu}{C_K}, \qquad \delta_4 = \min\left\{\frac{1}{C}, \frac{\rho_{\mathrm{f}}}{2C}, \frac{\mu}{2C}\right\}, \qquad \delta_5 = 4\rho_s\epsilon,
$$
$$
\delta_6 = 1, \qquad \delta_7 = 1, \qquad \delta_8 = \frac{\rho_s\epsilon}{2}.
$$

Using this values we obtain the following expression of the cumulative error at time step "$n$"

$$
\underbrace{\left(\frac{\rho_{\mathrm{f}}}{2}\|\boldsymbol{\theta}_h^n\|_{0,\Omega}^2+\frac{\rho_s\epsilon}{2}\|\dot{\boldsymbol{\xi}}_h^n\|_{0,\Sigma}^2+\frac{1}{2}\|\boldsymbol{\xi}_h^n\|_{\mathrm{s}}^2+\frac{\tau^2}{2\rho_{\mathrm{s}}\epsilon}\|\boldsymbol{L}_h^{\mathrm{s}}\boldsymbol{\xi}_h^n\|_{0,\Sigma}^2+\frac{\tau^2}{2}\|\dot{\boldsymbol{\xi}}_h^n\|_s^2\right)}_{\mathscr{E}_h^n}
$$

$$
-\underbrace{\left(\frac{\rho_{\mathrm{f}}}{2}\|\boldsymbol{\theta}_h^{n-1}\|_{0,\Omega}^2+\frac{\rho_s\epsilon}{2}\|\dot{\boldsymbol{\xi}}_h^{n-1}\|_{0,\Sigma}^2+\frac{1}{2}\|\boldsymbol{\xi}_h^{n-1}\|_{\mathrm{s}}^2+\frac{\tau^2}{2\rho_{\mathrm{s}}\epsilon}\|\boldsymbol{L}_h^{\mathrm{s}}\boldsymbol{\xi}_h^{n-1}\|_{0,\Sigma}^2+\frac{\tau^2}{2}\|\dot{\boldsymbol{\xi}}_h^{n-1}\|_s^2\right)}_{\mathscr{E}_h^0}
$$

$$
+\tau\sum_{m=1}^{n}\left[C\|\boldsymbol{\epsilon}(\boldsymbol{\theta}_h^m)\|_{0,\Omega}^2+C|\varphi_h^m|_{s_h}^2\right]
$$

$$
+\tau^2\sum_{m=1}^{n}\left[C\|\partial_\tau\boldsymbol{\theta}_h^m\|_{0,\Omega}^2+\frac{\rho_s\epsilon}{4}\|\partial_\tau\dot{\boldsymbol{\xi}}_h^m\|_{0,\Sigma}^2+\frac{1}{2}\|\partial_\tau\boldsymbol{\xi}_h^m\|_{\mathrm{s}}^2+\frac{1}{2}\|\partial_\tau\dot{\boldsymbol{\xi}}_h^m\|_s^2+\frac{1}{4\rho_{\mathrm{s}}\epsilon}\|\boldsymbol{L}_h^{\mathrm{s}}(\partial_\tau\boldsymbol{\xi}_h^m)\|_{0,\Sigma}^2\right]
$$

$$
\leq
$$

$$
\underbrace{\sum_{m=1}^{n}C_m^0}_{A_0}+\tau\underbrace{\sum_{m=1}^{n}C_m^1}_{A_1}+\tau^2\underbrace{\sum_{m=1}^{n}C_m^2}_{A_2}+\tau^3\underbrace{\sum_{m=1}^{n}C_m^3}_{A_3}
$$

$$
+\tau\sum_{m=1}^{n}\underbrace{\left[\frac{\rho_{\mathrm{f}}}{2}\|\boldsymbol{\theta}_h^m\|_{0,\Omega}^2+\frac{\rho_s\epsilon}{2}\|\dot{\boldsymbol{\xi}}_h^m\|_{0,\Sigma}^2+\frac{1}{2}\|\boldsymbol{\xi}_h^m\|_{\mathrm{s}}^2+\frac{\tau^2}{2\rho_{\mathrm{s}}\epsilon}\|\boldsymbol{L}_h^{\mathrm{s}}\boldsymbol{\xi}_h^n\|_{0,\Sigma}^2\right]}_{\mathscr{E}_h^m}. \quad (4.73)
$$

where the terms in the right hand side of the inequality are defined ad analyzed in the following using results of Proposition 3.3.6.

- The terms $A_0$ can be bounded as follows

$$
\begin{aligned}
A_0 \stackrel{def}{=} C\sum_{n=1}^{m}&\Big(\|\partial_t\dot{\boldsymbol{\xi}}_\pi\|_{L^2(t_{n-1},t_n;L^2(\Sigma)^d)}^2+\|\partial_t\boldsymbol{\theta}_\pi\|_{L^2(t_{n-1},t_n;L^2(\Omega)^d)}^2\\
&+\|\dot{\boldsymbol{\xi}}_\pi\|_{L^\infty(0,T;H^{\frac{1}{2}}(\Sigma)^d)}^2+\|\boldsymbol{\theta}_\pi\|_{L^\infty(0,T;H^1(\Omega)^d)}^2\Big)\\
\leq& C\Big(\|\partial_t\dot{\boldsymbol{\xi}}_\pi\|_{L^2(0,T;L^2(\Sigma)^d)}^2+\|\partial_t\boldsymbol{\theta}_\pi\|_{L^2(0,T;L^2(\Omega)^d)}^2\\
&+\|\dot{\boldsymbol{\xi}}_\pi\|_{L^\infty(0,T;H^{\frac{1}{2}}(\Sigma)^d)}^2+\|\boldsymbol{\theta}_\pi\|_{L^\infty(0,T;H^1(\Omega)^d)}^2\Big)\\
\approx& \mathcal{O}(h_s^{2m+2})+\mathcal{O}(h_f^{2l+2})+\mathcal{O}(h_s^{2m+1})+\mathcal{O}(h_f^{2l}).
\end{aligned} \quad (4.74)
$$

- The term $A_1$ is defined and bounded as follows

$$
\begin{aligned}
A_1 \stackrel{def}{=} & C \sum_{n=1}^{m} \Big( \|\partial_t \dot{\boldsymbol{\xi}}_\pi\|^2_{L^2(t_{n-1},t_n;L^2(\Sigma)^d)} + \|\dot{\boldsymbol{\xi}}_\pi\|^2_{L^\infty(0,T;\boldsymbol{W})} \\
& + \|\boldsymbol{\omega}_\pi\|^2_{L^\infty(0,T;\boldsymbol{\Lambda})} + \|\dot{\boldsymbol{\xi}}_\pi\|^2_{L^\infty(0,T;H^{\frac{1}{2}}(\Sigma)^d)} + \|\boldsymbol{\theta}_\pi\|^2_{L^\infty(0,T;H^1(\Omega)^d)} \Big) \\
\le & C\Big( \|\partial_t \dot{\boldsymbol{\xi}}_\pi\|^2_{L^2(0,T;L^2(\Sigma)^d)} + \|\dot{\boldsymbol{\xi}}_\pi\|^2_{L^\infty(0,T;\boldsymbol{W})} \\
& + \|\boldsymbol{\omega}_\pi\|^2_{L^\infty(0,T;\boldsymbol{\Lambda})} + \|\dot{\boldsymbol{\xi}}_\pi\|^2_{L^\infty(0,T;H^{\frac{1}{2}}(\Sigma)^d)} + \|\boldsymbol{\theta}_\pi\|^2_{L^\infty(0,T;H^1(\Omega)^d)} \Big) \\
\approx & \mathcal{O}(h_s^{2m+2}) + \mathcal{O}(h_s^{2m}) + \mathcal{O}(h_s^{2z+1}) + \mathcal{O}(h_s^{2m+1}) + \mathcal{O}(h_f^{2l}).
\end{aligned}
\tag{4.75}
$$

– The term $A_2$ is defined and bounded in the following using also the properties of the solid operator $\|\boldsymbol{L}_h^{\mathrm{s}} \boldsymbol{d}^n\|_{0,\Sigma} \le C\|\boldsymbol{L}^{\mathrm{s}} \boldsymbol{d}^n\|_{0,\Sigma}$

$$
\begin{aligned}
A_2 \stackrel{def}{=} & C \sum_{n=1}^{m} \Big( \|\boldsymbol{L}_h^{\mathrm{s}}(\boldsymbol{d}^n - \boldsymbol{d}^{n-1})\|^2_{0,\Sigma} + \|\partial_{tt}\boldsymbol{d}\|^2_{L^2(t_{n-1},t_n;\boldsymbol{W})} \\
& + \|\partial_{tt}\dot{\boldsymbol{d}}\|^2_{L^2(t_{n-1},t_n;L^2(\Sigma)^d)} + \|\partial_{tt}\boldsymbol{u}\|^2_{L^2(t_{n-1},t_n;L^2(\Omega)^d)} \Big) \\
\le & C\Big( \|\partial_t \boldsymbol{L}^{\mathrm{s}}(\boldsymbol{d})\|^2_{L^\infty(0,T;L^2(\Sigma)^d)} + \|\partial_{tt}\boldsymbol{d}\|^2_{L^2(0,T;\boldsymbol{W})} \\
& + \|\partial_{tt}\dot{\boldsymbol{d}}\|^2_{L^2(0,T;L^2(\Sigma)^d)} + \|\partial_{tt}\boldsymbol{u}\|^2_{L^2(0,T;L^2(\Omega)^d)} \Big).
\end{aligned}
\tag{4.76}
$$

– The terms $A_3$ is bounded as follows

$$
A_3 \stackrel{def}{=} C \sum_{n=1}^{m} \|\partial_{tt}\dot{\boldsymbol{d}}\|^2_{L^2(t_{n-1},t_n;L^2(\Sigma)^d)} \le C\|\partial_{tt}\dot{\boldsymbol{d}}\|^2_{L^2(0,T;L^2(\Sigma)^d)}
\tag{4.77}
$$

In sum the error estimate (4.73) can be written as follows

$$
\mathcal{E}_h^n + \tau\mathcal{D}_h^n \le \tau \sum_{m=1}^{n} \mathcal{E}_h^m + A_0 + \tau A_1 + \tau^2 A_2 + \tau^3 A_3
\tag{4.78}
$$

where $\tau\mathcal{D}_h^n$ represents the dissipation term and is given by

$$
\begin{aligned}
\tau\mathcal{D}_h^n = & \tau \sum_{m=1}^{n} \Big[ C\|\boldsymbol{\epsilon}(\boldsymbol{\theta}_h^m)\|^2_{0,\Omega} + C|\varphi_h^m|^2_{s_h} \Big] \\
& + \tau^2 \sum_{m=1}^{n} \Big[ C\|\partial_\tau \boldsymbol{\theta}_h^m\|^2_{0,\Omega} + \frac{\rho_s\epsilon}{4}\|\partial_\tau \dot{\boldsymbol{\xi}}_h^m\|^2_{0,\Sigma} + \frac{1}{2}\|\partial_\tau \boldsymbol{\xi}_h^m\|^2_{\mathrm{s}} + \frac{1}{2}\|\partial_\tau \dot{\boldsymbol{\xi}}_h^m\|^2_{s} + \frac{1}{4\rho_{\mathrm{s}}\epsilon}\|\boldsymbol{L}_h^{\mathrm{s}}(\partial_\tau \boldsymbol{\xi}_h^m)\|^2_{0,\Sigma} \Big].
\end{aligned}
$$

Introducing the convergence rates with respect to $h_f$ and $h_s$ of the term $A_0, A_1, A_2, A_3$ given in (4.74)-(4.75), we have the desired result applying the discrete Gronwall Lemma 3.3.1, since by hypothesis we have $\tau < 1$

$$
\mathcal{E}_h^n + \tau\mathcal{D}_h^n \lesssim \exp\left( \sum_{m=1}^{n} \frac{\tau}{1-\tau} \right) \Big[ \underbrace{\mathcal{O}(h_f^{2l}) + \mathcal{O}(h_s^{2m+1})}_{A_0}
\tag{4.79}
$$

$$+\tau\Big(\underbrace{\mathcal{O}(h_s^{2m})+\mathcal{O}(h_s^{2z+1})+\mathcal{O}(h_f^{2l})}_{A_1}\Big)+\tau^2A_2+\tau^3A_3\Big]. \quad (4.80)$$

- **Scheme with** $r=2$ implies $\boldsymbol{d}^{n\star}=\boldsymbol{d}^{n-1}+\tau\dot{\boldsymbol{d}}^{n-1}$; in this case the error extrapolation has the following expression

$$\begin{aligned}\boldsymbol{\xi}_h^{n\star}=\boldsymbol{d}_\Pi^{n\star}-\boldsymbol{d}_h^{n\star}&=(\boldsymbol{d}_\Pi^{n-1}+\tau\dot{\boldsymbol{d}}_\Pi^{n-1})-\boldsymbol{d}_h^{n-1}-\tau\dot{\boldsymbol{d}}_h^{n-1}\\&=(\boldsymbol{d}_\Pi^{n-1}-\boldsymbol{d}_h^{n-1})+\tau(\dot{\boldsymbol{d}}_\Pi^{n-1}-\dot{\boldsymbol{d}}_h^{n-1})\\&=\boldsymbol{\xi}_h^{n-1}+\tau\dot{\boldsymbol{\xi}}_h^{n-1},\end{aligned} \quad (4.81)$$

then, using the relation between the error $\dot{\boldsymbol{\xi}}_h^n$ and the error $\boldsymbol{\xi}_h^n$ expressed by (3.31), we have

$$\begin{aligned}\boldsymbol{\xi}_h^n-\boldsymbol{\xi}_h^{n\star}=\boldsymbol{\xi}_h^n-\boldsymbol{\xi}_h^{n-1}-\tau\dot{\boldsymbol{\xi}}_h^{n-1}&=\tau\partial_\tau\boldsymbol{\xi}_h^n-\tau\dot{\boldsymbol{\xi}}_h^{n-1}=\tau(\dot{\boldsymbol{\xi}}_h^n-\dot{\boldsymbol{\xi}}_h^{n-1})+\tau\left(\partial_\tau\boldsymbol{d}_\Pi^n-\dot{\boldsymbol{d}}_\Pi^n\right)\\&=\tau(\dot{\boldsymbol{\xi}}_h^n-\dot{\boldsymbol{\xi}}_h^{n-1})+\tau\left(\partial_\tau(\boldsymbol{d}_\Pi^n-\boldsymbol{d}^n)+(\partial_\tau-\partial_t)\boldsymbol{d}^n+(\dot{\boldsymbol{d}}^n-\dot{\boldsymbol{d}}_\Pi^n)\right)\\&=\tau(\dot{\boldsymbol{\xi}}_h^n-\dot{\boldsymbol{\xi}}_h^{n-1})-\tau\partial_\tau\boldsymbol{\xi}_\pi^n+\tau\dot{\boldsymbol{\xi}}_\pi^n+\tau(\partial_\tau-\partial_t)\boldsymbol{d}^n\end{aligned} \quad (4.82)$$

In order to obtain the required estimate for the terms from $T_5$ to $T_{10}$ we use results of proposition 3.3.2, in particular for terms $T_5$ and $T_6$ we have

$$\begin{aligned}T_5&=\rho_s\epsilon\tau\left((\partial_\tau-\partial_t)\dot{\boldsymbol{d}}^n-\partial_\tau\dot{\boldsymbol{\xi}}_\pi^n,\frac{\tau}{\rho_s\epsilon}\boldsymbol{L}_h^s(\boldsymbol{\xi}_h^n-\boldsymbol{\xi}_h^{n\star})\right)_{0,\Sigma}\\&=\tau^2\left(\tau(\partial_\tau-\partial_t)\dot{\boldsymbol{d}}^n-\tau\partial_\tau\dot{\boldsymbol{\xi}}_\pi^n,\boldsymbol{L}_h^s\left((\dot{\boldsymbol{\xi}}_h^n-\dot{\boldsymbol{\xi}}_h^{n-1})-\partial_\tau\boldsymbol{\xi}_\pi^n+\dot{\boldsymbol{\xi}}_\pi^n+(\partial_\tau-\partial_t)\boldsymbol{d}^n\right)\right)_{0,\Sigma}\\&\le\tau^2\left(\|\tau(\partial_\tau-\partial_t)\dot{\boldsymbol{d}}^n\|_{0,\Sigma}+\|\tau\partial_\tau\dot{\boldsymbol{\xi}}_\pi^n\|_{0,\Sigma}\right)\left\|\boldsymbol{L}_h^s\left((\dot{\boldsymbol{\xi}}_h^n-\dot{\boldsymbol{\xi}}_h^{n-1})-\partial_\tau\boldsymbol{\xi}_\pi^n+\dot{\boldsymbol{\xi}}_\pi^n+(\partial_\tau-\partial_t)\boldsymbol{d}^n\right)\right\|_{0,\Sigma}\\&\le\frac{\tau^{\frac{7}{2}}\beta_s^{\frac{3}{2}}C_I^3}{h_s^3}\|\partial_{tt}\boldsymbol{d}\|_{L^2(t_{n-1},t_n;L^2(\Sigma)^d)}\|(\dot{\boldsymbol{\xi}}_h^n-\dot{\boldsymbol{\xi}}_h^{n-1})-\partial_\tau\boldsymbol{\xi}_\pi^n+\dot{\boldsymbol{\xi}}_\pi^n+(\partial_\tau-\partial_t)\boldsymbol{d}^n\|_{0,\Sigma}\\&+\frac{\tau^{\frac{5}{2}}\beta_s^{\frac{3}{2}}C_I^3}{h_s^3}\|\partial_t\dot{\boldsymbol{\xi}}_\pi\|_{L^2(t_{n-1},t_n;L^2(\Sigma)^d)}\|(\dot{\boldsymbol{\xi}}_h^n-\dot{\boldsymbol{\xi}}_h^{n-1})-\partial_\tau\boldsymbol{\xi}_\pi^n+\dot{\boldsymbol{\xi}}_\pi^n+(\partial_\tau-\partial_t)\boldsymbol{d}^n\|_{0,\Sigma}\\&\le\frac{\tau^2\alpha^5}{2\delta_5}\|\partial_{tt}\boldsymbol{d}\|^2_{L^2(t_{n-1},t_n;L^2(\Sigma)^d)}+\frac{\alpha^5}{2\delta_5}\|\partial_t\dot{\boldsymbol{\xi}}_\pi\|^2_{L^2(t_{n-1},t_n;L^2(\Sigma)^d)}\\&\quad+(\rho_s\epsilon)^3\delta_5\|\dot{\boldsymbol{\xi}}_h^n-\dot{\boldsymbol{\xi}}_h^{n-1}\|^2_{0,\Sigma}+(\rho_s\epsilon)^3\delta_5\|\partial_\tau\boldsymbol{\xi}_\pi^n\|^2_{0,\Sigma}\\&\quad+(\rho_s\epsilon)^3\delta_5\|\dot{\boldsymbol{\xi}}_\pi^n\|^2_{0,\Sigma}+(\rho_s\epsilon)^3\delta_5\|(\partial_\tau-\partial_t)\boldsymbol{d}^n\|^2_{0,\Sigma}\\&\le\frac{\tau^2\alpha^5}{2\delta_5}\|\partial_{tt}\boldsymbol{d}\|^2_{L^2(t_{n-1},t_n;L^2(\Sigma)^d)}+\frac{\alpha^5}{2\delta_5}\|\partial_t\dot{\boldsymbol{\xi}}_\pi\|^2_{L^2(t_{n-1},t_n;L^2(\Sigma)^d)}\end{aligned}$$

$$+(\rho_s\epsilon)^3\delta_5\|\partial_t\boldsymbol{\xi}_\pi\|^2_{L^\infty(0,T;L^2(\Sigma)^d)}+(\rho_s\epsilon)^3\delta_5\|\dot{\boldsymbol{\xi}}_\pi\|^2_{L^\infty(0,T;L^2(\Sigma)^d)}$$
$$+\tau^2(\rho_s\epsilon)^3\delta_5\|\partial_{tt}\boldsymbol{d}\|^2_{L^\infty(0,T;L^2(\Sigma)^d)}+(\rho_s\epsilon)^3\delta_5\|\dot{\boldsymbol{\xi}}^n_h-\dot{\boldsymbol{\xi}}^{n-1}_h\|^2_{0,\Sigma}. \quad (4.83)$$

$$\begin{aligned}T_6&=\rho_s\epsilon\tau\left(\partial_\tau\dot{\boldsymbol{\xi}}^n_\pi-(\partial_\tau-\partial_t)\dot{\boldsymbol{d}}^n,\frac{\tau}{\rho_s\epsilon}\boldsymbol{L}^s_h(\boldsymbol{\Pi_W}(\boldsymbol{d}^n-\boldsymbol{d}^{n-1}-\tau\dot{\boldsymbol{d}}^{n-1}))\right)_{0,\Sigma}\\
&=\tau^2\left(\tau\partial_\tau\dot{\boldsymbol{\xi}}^n_\pi-\tau(\partial_\tau-\partial_t)\dot{\boldsymbol{d}}^n,\boldsymbol{L}^s_h(\partial_\tau\boldsymbol{d}^n-\dot{\boldsymbol{d}}^{n-1}))\right)_{0,\Sigma}\\
&\le\left(\frac{\tau^2\beta_s^{\frac{3}{2}}C_I^3}{h_s^3}\|\tau(\partial_\tau-\partial_t)\dot{\boldsymbol{d}}^n\|_{0,\Sigma}+\frac{\tau^2\beta_s^{\frac{3}{2}}C_I^3}{h_s^3}\|\tau\partial_\tau\dot{\boldsymbol{\xi}}^n_\pi\|_{0,\Sigma}\right)\|\partial_\tau\boldsymbol{d}^n-\dot{\boldsymbol{d}}^{n-1}\|_{0,\Sigma}\\
&\le\left(\frac{\tau^{\frac{7}{2}}\beta_s^{\frac{3}{2}}C_I^3}{h_s^3}\|\partial_{tt}\dot{\boldsymbol{d}}\|_{L^2(t_{n-1},t_n;L^2(\Sigma)^d)}+\frac{\tau^{\frac{5}{2}}\beta_s^{\frac{3}{2}}C_I^3}{h_s^3}\|\partial_t\dot{\boldsymbol{\xi}}_\pi\|_{L^2(t_{n-1},t_n;L^2(\Sigma)^d)}\right)\tau\|\partial_{tt}\boldsymbol{d}\|_{L^\infty(0,T;L^2(\Sigma)^d)}\\
&\le\frac{\tau^2\alpha^5\delta_5}{2}\|\partial_{tt}\dot{\boldsymbol{d}}\|^2_{L^2(t_{n-1},t_n;L^2(\Sigma)^d)}+\frac{\alpha^5\delta_6}{2}\|\partial_t\dot{\boldsymbol{\xi}}_\pi\|^2_{L^2(t_{n-1},t_n;L^2(\Sigma)^d)}+\frac{\tau^2(\rho_s\epsilon)^3}{\delta_6}\|\partial_{tt}\boldsymbol{d}\|_{L^\infty(0,T;L^2(\Sigma)^d)}.\end{aligned} \quad (4.84)$$

The estimates for the terms $T_7$, $T_8$ and $T_9$ use the definition of the solid discrete operator $\boldsymbol{L}^s_h$ given in (3.25) and the CFL condition 4.24.

$$\begin{aligned}T_7&=\tau a_s\left(\boldsymbol{\xi}^n_h,\frac{\tau}{\rho_s\epsilon}\boldsymbol{L}^s_h(\boldsymbol{\Pi_W}(\boldsymbol{d}^n-\boldsymbol{d}^{n-1}-\tau\dot{\boldsymbol{d}}^{n-1}))\right)=\frac{\tau^3}{\rho_s\epsilon}a_s\left(\boldsymbol{\xi}^n_h,\boldsymbol{L}^s_h(\partial_\tau\boldsymbol{d}^n-\dot{\boldsymbol{d}}^{n-1})\right)\\
&\le\frac{\tau^3}{\rho_s\epsilon}\|\boldsymbol{\xi}^n_h\|_s\|\boldsymbol{L}^s_h(\partial_\tau\boldsymbol{d}^n-\dot{\boldsymbol{d}}^{n-1})\|_s\le\frac{\tau^3\beta_s^{\frac{3}{2}}C_I^3}{\rho_s\epsilon h_s^3}\|\boldsymbol{\xi}^n_h\|_s\|(\partial_\tau\boldsymbol{d}^n-\dot{\boldsymbol{d}}^{n-1})\|_{0,\Sigma}\\
&\le\frac{\tau\alpha^5\delta_7}{2}\|\boldsymbol{\xi}^n_h\|^2_s+\frac{\tau^2\rho_s\epsilon}{2\delta_7}\|\partial_{tt}\boldsymbol{d}\|^2_{L^\infty(0,T;L^2(\Sigma)^d)}.\end{aligned} \quad (4.85)$$

In order to estimate the term $T_8$ we use the CFL condition $\tau\alpha^5=\frac{\tau^6\beta_s^3C_I^6}{h_s^6(\rho_s\epsilon)^{\frac{3}{2}}}<1$, then

$$\begin{aligned}T_8&=\rho_s\epsilon\tau\left(\partial_\tau\dot{\boldsymbol{\xi}}^n_h,\frac{\tau}{\rho_s\epsilon}\boldsymbol{L}^s_h(\boldsymbol{\Pi_W}(\boldsymbol{d}^n-\boldsymbol{d}^{n-1}-\tau\dot{\boldsymbol{d}}^{n-1}))\right)_{0,\Sigma}=\tau^2\left(\tau\partial_\tau\dot{\boldsymbol{\xi}}^n_h,\boldsymbol{L}^s_h(\partial_\tau\boldsymbol{d}^n-\dot{\boldsymbol{d}}^{n-1})\right)_{0,\Sigma}\\
&\le\tau^2\|\boldsymbol{L}^s_h(\tau\partial_\tau\dot{\boldsymbol{\xi}}^n_h)\|_{0,\Sigma}\|\partial_\tau\boldsymbol{d}^n-\dot{\boldsymbol{d}}^{n-1})\|_{0,\Sigma}\le\tau\underbrace{\left(\frac{\tau\beta_s^{\frac{1}{2}}C_I}{h_s(\rho_s\epsilon)^{\frac{1}{2}}}\right)}_{\le1}\|\tau\partial_\tau\dot{\boldsymbol{\xi}}^n_h\|_s\tau(\rho_s\epsilon)^{\frac{1}{2}}\|\partial_{tt}\boldsymbol{d}\|_{L^\infty(0,T;L^2(\Sigma)^d)}\\
&\le\frac{\tau^2\delta_8}{2}\|\tau\partial_\tau\dot{\boldsymbol{\xi}}^n_h\|^2_s+\frac{\tau^2\rho_s\epsilon}{2\delta_8}\|\partial_{tt}\boldsymbol{d}\|^2_{L^\infty(0,T;L^2(\Sigma)^d)}.\end{aligned} \quad (4.86)$$

$$
\begin{aligned}
T_9 = \rho_s \epsilon \tau & \left( \partial_\tau \dot{\boldsymbol{\xi}}_h^n, \frac{\tau}{\rho_{\rm s}\epsilon} \boldsymbol{L}_h^{\rm s} \left( \tau(\dot{\boldsymbol{\xi}}_h^n - \dot{\boldsymbol{\xi}}_h^{n-1}) - \tau \partial_\tau \boldsymbol{\xi}_\pi^n + \tau \dot{\boldsymbol{\xi}}_\pi^n + \tau(\partial_\tau - \partial_t) \boldsymbol{d}^n \right) \right)_{0,\Sigma} \\
&= \tau^2 \left( \boldsymbol{L}_h^{\rm s}(\tau \partial_\tau \dot{\boldsymbol{\xi}}_h^n), \tau \partial_\tau \dot{\boldsymbol{\xi}}_h^n \right)_{0,\Sigma} - \tau^2 \left( \boldsymbol{L}_h^{\rm s}(\tau \partial_\tau \dot{\boldsymbol{\xi}}_h^n), \partial_\tau \boldsymbol{\xi}_\pi^n \right)_{0,\Sigma} \\
&+ \tau^2 \left( \boldsymbol{L}_h^{\rm s}(\tau \partial_\tau \dot{\boldsymbol{\xi}}_h^n), \dot{\boldsymbol{\xi}}_\pi^n \right)_{0,\Sigma} - \tau^2 \left( \boldsymbol{L}_h^{\rm s}(\tau \partial_\tau \dot{\boldsymbol{\xi}}_h^n), (\partial_\tau - \partial_t) \boldsymbol{d}^n \right)_{0,\Sigma} \\
&\geq \tau^2 \|\tau \partial_\tau \dot{\boldsymbol{\xi}}_h^n\|_{\rm s}^2 - \frac{3\tau^2 \delta_9}{2} \|\tau \partial_\tau \dot{\boldsymbol{\xi}}_h^n\|_{\rm s}^2 - \frac{\tau^2}{2\delta_9} \|\partial_\tau \boldsymbol{\xi}_\pi^n\|_{\rm s}^2 - \frac{\tau^2}{2\delta_9} \|\dot{\boldsymbol{\xi}}_\pi^n\|_{\rm s}^2 - \frac{\tau^2}{2\delta_9} \|(\partial_\tau - \partial_t) \boldsymbol{d}^n\|_{\rm s}^2 \\
&\geq \tau^2 \left( 1 - \frac{3\delta_9}{2} \right) \|\tau \partial_\tau \dot{\boldsymbol{\xi}}_h^n\|_{\rm s}^2 - \frac{\tau^2}{2\delta_9} \frac{\beta_{\rm s} C_I^2}{h_{\rm s}^2} \left( \|\partial_\tau \boldsymbol{\xi}_\pi^n\|_{0,\Sigma}^2 + \|\dot{\boldsymbol{\xi}}_\pi^n\|_{0,\Sigma}^2 + \|(\partial_\tau - \partial_t) \boldsymbol{d}^n\|_{0,\Sigma}^2 \right) \\
&\geq \tau^2 \left( 1 - \frac{3\delta_9}{2} \right) \|\tau \partial_\tau \dot{\boldsymbol{\xi}}_h^n\|_{\rm s}^2 \\
&- \frac{\alpha^{\frac{5}{3}}}{2\delta_9} \left( \|\partial_t \boldsymbol{\xi}_\pi\|_{L^\infty(0,T;L^2(\Sigma)^d)}^2 + \|\dot{\boldsymbol{\xi}}_\pi\|_{L^\infty(0,T;L^2(\Sigma)^d)}^2 + \tau^2 \|\partial_{tt} \boldsymbol{d}\|_{L^\infty(0,T;L^2(\Sigma)^d)}^2 \right)
\end{aligned} \tag{4.87}
$$

The estimate of the term $T_{10}$ is performed in two steps; both of them use the CFL condition (4.24) stated in the stability analysis of algorithm 3 for the case of extrapolation of order $r = 2$

$$
\begin{aligned}
T_{10} &= \frac{\tau^2}{\rho_{\rm s}\epsilon} a^{\rm s} \left( \boldsymbol{\xi}_h^n, \boldsymbol{L}_h^{\rm s} \left( \tau(\dot{\boldsymbol{\xi}}_h^n - \dot{\boldsymbol{\xi}}_h^{n-1}) - \tau \partial_\tau \boldsymbol{\xi}_\pi^n + \tau \dot{\boldsymbol{\xi}}_\pi^n + \tau(\partial_\tau - \partial_t) \boldsymbol{d}^n \right) \right) \\
&= \underbrace{\frac{\tau^3}{\rho_{\rm s}\epsilon} \left( \boldsymbol{L}_h^{\rm s} \boldsymbol{\xi}_h^n, \boldsymbol{L}_h^{\rm s} (\dot{\boldsymbol{\xi}}_h^n - \dot{\boldsymbol{\xi}}_h^{n-1}) \right)_{0,\Sigma}}_{T_{10}^A} + \underbrace{\frac{\tau^3}{\rho_{\rm s}\epsilon} \left( \boldsymbol{L}_h^{\rm s} \boldsymbol{\xi}_h^n, \boldsymbol{L}_h^{\rm s} ((\partial_\tau - \partial_t) \boldsymbol{d}^n) \right)_{0,\Sigma}}_{T_{10}^B} \\
&+ \underbrace{\frac{\tau^3}{\rho_{\rm s}\epsilon} a^s \left( \boldsymbol{L}_h^{\rm s} \boldsymbol{\xi}_h^n, \dot{\boldsymbol{\xi}}_\pi^n - \partial_\tau \boldsymbol{\xi}_\pi^n \right)_{0,\Sigma}}_{=0 \text{ by the definitions of } \boldsymbol{\xi}_\pi^n, \dot{\boldsymbol{\xi}}_\pi^n} = T_{10}^A + T_{10}^B.
\end{aligned} \tag{4.88}
$$

The term $T_{10}^A$ is bounded as follows using the inverse estimates given in Lemma 3.3.4

$$
\begin{aligned}
T_{10}^A &= \frac{\tau^3}{\rho_{\rm s}\epsilon} a_{\rm s} (\boldsymbol{L}_h^{\rm s} \boldsymbol{\xi}_h^n, \dot{\boldsymbol{\xi}}_h^n - \dot{\boldsymbol{\xi}}_h^{n-1}) \leq \frac{\tau^3}{\rho_{\rm s}\epsilon} \|\boldsymbol{L}_h^{\rm s} \boldsymbol{\xi}_h^n\|_{\rm s} \|\dot{\boldsymbol{\xi}}_h^n - \dot{\boldsymbol{\xi}}_h^{n-1}\|_{\rm s} \leq \frac{\tau^3}{\rho_{\rm s}\epsilon} \frac{\beta_{\rm s}^{\frac{1}{2}} C_I}{h} \|\boldsymbol{L}_h^{\rm s} \boldsymbol{\xi}_h^n\|_{\rm s} \|\dot{\boldsymbol{\xi}}_h^n - \dot{\boldsymbol{\xi}}_h^{n-1}\|_{0,\Sigma} \\
&\leq \frac{\tau^3}{(\rho_{\rm s}\epsilon)^{\frac{3}{2}}} \frac{\beta_{\rm s}^{\frac{3}{2}} C_I^3}{h^3} \|\boldsymbol{\xi}_h^n\|_{\rm s} (\rho_{\rm s}\epsilon)^{\frac{1}{2}} \|\dot{\boldsymbol{\xi}}_h^n - \dot{\boldsymbol{\xi}}_h^{n-1}\|_{0,\Sigma} \leq \frac{\tau^6}{(\rho_{\rm s}\epsilon)^3} \frac{\beta_{\rm s}^3 C_I^6}{h^6} \|\boldsymbol{\xi}_h^n\|_{\rm s}^2 + \frac{\rho_{\rm s}\epsilon}{4} \|\dot{\boldsymbol{\xi}}_h^n - \dot{\boldsymbol{\xi}}_h^{n-1}\|_{0,\Sigma}^2 \\
&\leq \tau \alpha^5 \|\boldsymbol{\xi}_h^n\|_{\rm s}^2 + \frac{\rho_{\rm s}\epsilon}{4} \|\dot{\boldsymbol{\xi}}_h^n - \dot{\boldsymbol{\xi}}_h^{n-1}\|_{0,\Sigma}^2.
\end{aligned} \tag{4.89}
$$

Also the estimate of term $T_{10}^B$ is obtained using inequalities given in Lemma 3.3.4

$$
T_{10}^B = \frac{\tau^3}{\rho_{\rm s}\epsilon} \left( \boldsymbol{L}_h^{\rm s} \boldsymbol{\xi}_h^n, \boldsymbol{L}_h^{\rm s} ((\partial_\tau - \partial_t) \boldsymbol{d}^n) \right)_{0,\Sigma} \leq \frac{\tau^3}{\rho_{\rm s}\epsilon} \|\boldsymbol{L}_h^{\rm s} \boldsymbol{\xi}_h^n\|_{\rm s} \|\boldsymbol{L}_h^{\rm s} ((\partial_\tau - \partial_t) \boldsymbol{d}^n)\|_{\rm s}
$$

$$\leq \frac{\tau^3}{\rho_s \epsilon}\left(\frac{\beta_s^{\frac{1}{2}} C_I}{h_s}\|\boldsymbol{\xi}_h^n\|_s\right)\left(\frac{\beta_s C_I^2}{h_s^2}\|(\partial_\tau - \partial_t)\boldsymbol{d}^n)\|_{0,\Sigma}\right) \leq$$

$$\leq \frac{\delta_{10}^B \tau \alpha^5}{2}\|\boldsymbol{\xi}_h^n\|_s^2 + \frac{\tau^2 \rho_s \epsilon}{2\delta_{10^B}}\|\partial_{tt}\boldsymbol{d}\|^2_{L^\infty(0,T;L^2(\Sigma)^d)}. \quad (4.90)$$

Then we have

$$T_{10} \;\leq\; \tau\alpha^5\left(1+\frac{\delta_{10}^B}{2}\right)\|\boldsymbol{\xi}_h^n\|_s^2 \;+\; \frac{\rho_s\epsilon}{4}\|\dot{\boldsymbol{\xi}}_h^n \;-\; \dot{\boldsymbol{\xi}}_h^{n-1}\|^2_{0,\Sigma} \;+\; \frac{\tau^2\rho_s\epsilon}{2\delta_{10^B}}\|\partial_{tt}\boldsymbol{d}\|^2_{L^\infty(0,T;L^2(\Sigma)^d)}. \quad (4.91)$$

Using the estimates of terms from $T_0$ to $T_{10}$ in (4.43) and collecting together the similar terms, we obtain the following inequality that gives a bound for the discrete error at any instant of time $t_n$ in the case of extrapolation order $r = 2$.

$$\begin{aligned}
&\left(\frac{\rho_f}{2}\|\boldsymbol{\theta}_h^n\|^2_{0,\Omega} + \frac{\rho_s\epsilon}{2}\|\dot{\boldsymbol{\xi}}_h^n\|^2_{0,\Sigma} + \frac{1}{2}\|\boldsymbol{\xi}_h^n\|_s^2\right) - \left(\frac{\rho_f}{2}\|\boldsymbol{\theta}_h^{n-1}\|^2_{0,\Omega} + \frac{\rho_s\epsilon}{2}\|\dot{\boldsymbol{\xi}}_h^{n-1}\|^2_{0,\Sigma} + \frac{1}{2}\|\boldsymbol{\xi}_h^{n-1}\|_s^2\right) \\
&\quad + \tau^2\left(\left(\frac{\rho_f}{2} - C\delta_4\right)\|\boldsymbol{\theta}_h^n - \boldsymbol{\theta}_h^{n-1}\|^2_{0,\Omega} + \left(\frac{\rho_s\epsilon}{4} - (\rho_s\epsilon)^3\delta_5\right)\|\dot{\boldsymbol{\xi}}_h^n - \dot{\boldsymbol{\xi}}_h^{n-1}\|^2_{0,\Sigma} + \frac{1}{2}\|\boldsymbol{\xi}_h^n - \boldsymbol{\xi}_h^{n-1}\|_s^2\right) \\
&\quad + \tau\left[\left(2\mu - \delta_3 C_K - C\delta_4\right)\|\boldsymbol{\epsilon}(\boldsymbol{\theta}_h^n)\|^2_{0,\Omega} + (1 - C\delta_4)|\varphi_h^n|^2_{s_h}\right] \\
&\quad + \tau^2\left(1 - \frac{3\delta_9}{2} - \frac{\delta_8}{2}\right)\|\dot{\boldsymbol{\xi}}_h^n - \dot{\boldsymbol{\xi}}_h^{n-1}\|_s^2 \\
&\leq \\
&\Big[\frac{\rho_f}{2\delta_1}\|\partial_t\boldsymbol{\theta}_\pi\|^2_{L^2(t_{n-1},t_n;L^2(\Omega)^d)} + \frac{\rho_s\epsilon}{2\delta_2}\|\partial_t\dot{\boldsymbol{\xi}}_\pi\|^2_{L^2(t_{n-1},t_n;L^2(\Sigma)^d)} \\
&\quad + C\delta_4\|\partial_t\boldsymbol{\theta}_\pi\|^2_{L^2(t_{n-1},t_n;L^2(\Omega)^d)} + \frac{C}{\delta_4}\|\dot{\boldsymbol{\xi}}_\pi\|_{L^\infty(0,T;H^{\frac{1}{2}}(\Sigma)^d)} + \frac{C}{\delta_4}\|\boldsymbol{\theta}_\pi\|_{L^\infty(0,T;H^1(\Omega)^d)} \\
&\quad + \frac{\alpha^{\frac{5}{3}}}{2\delta_9}\|\partial_t\boldsymbol{\xi}_\pi\|^2_{L^\infty(0,T;L^2(\Sigma)^d)} + \frac{\alpha^{\frac{5}{3}}}{2\delta_9}\|\dot{\boldsymbol{\xi}}_\pi\|^2_{L^\infty(0,T;L^2(\Sigma)^d)} \\
&\quad + \frac{\alpha^5\delta_6}{2}\|\partial_t\dot{\boldsymbol{\xi}}_\pi\|^2_{L^2(t_{n-1},t_n;L^2(\Sigma)^d)} + \frac{\alpha^5}{2\delta_5}\|\partial_t\dot{\boldsymbol{\xi}}_\pi\|^2_{L^2(t_{n-1},t_n;L^2(\Sigma)^d)} \\
&\quad + (\rho_s\epsilon)^3\delta_5\|\partial_t\boldsymbol{\xi}_\pi\|^2_{L^\infty(0,T;L^2(\Sigma)^d)} + (\rho_s\epsilon)^3\delta_5\|\dot{\boldsymbol{\xi}}_\pi\|^2_{L^\infty(0,T;L^2(\Sigma)^d)}\Big] \\
&+ \tau\Big[\frac{1}{\delta_3}\|\boldsymbol{\omega}_\pi\|^2_{L^\infty(0,T;\boldsymbol{\Lambda})} + C\delta_4\|\boldsymbol{\omega}_\pi\|^2_{L^\infty(0,T;\boldsymbol{\Lambda})} + \frac{C}{\delta_4}\|\dot{\boldsymbol{\xi}}_\pi\|^2_{L^\infty(0,T;H^{\frac{1}{2}}(\Sigma)^d)} + \frac{C}{\delta_4}\|\boldsymbol{\theta}_\pi\|^2_{L^\infty(0,T;H^1(\Omega)^d)}\Big]
\end{aligned}$$

$$
\begin{aligned}
&+\tau^2\Big[\frac{\alpha^{\frac{5}{3}}}{2\delta_9}\|\partial_{tt}\boldsymbol{d}\|^2_{L^\infty(0,T;L^2(\Sigma)^d)}+\frac{\rho_{\mathrm{s}}\epsilon}{2\delta_{10^B}}\|\partial_{tt}\boldsymbol{d}\|^2_{L^\infty(0,T;L^2(\Sigma)^d)}+\frac{\alpha^5}{2\delta_5}\|\partial_{tt}\dot{\boldsymbol{d}}\|^2_{L^2(t_{n-1},t_n;L^2(\Sigma)^d)}\\
&\qquad+\frac{\alpha^5}{2\delta_5}\|\partial_{tt}\dot{\boldsymbol{d}}\|^2_{L^2(t_{n-1},t_n;L^2(\Sigma)^d)}+(\rho_{\mathrm{s}}\epsilon)^3\delta_5\|\partial_{tt}\boldsymbol{d}\|^2_{L^\infty(0,T;L^2(\Sigma)^d)}\\
&\qquad+\frac{1}{2\delta_0}\|\partial_{tt}\boldsymbol{d}\|^2_{L^2(t_{n-1},t_n;\boldsymbol{W})}+\frac{\rho_{\mathrm{f}}}{2\delta_1}\|\partial_{tt}\boldsymbol{u}\|^2_{L^2(t_{n-1},t_n;L^2(\Omega)^d)}\\
&\qquad+\frac{\rho_s\epsilon}{2\delta_2}\|\partial_{tt}\dot{\boldsymbol{d}}\|^2_{L^2(t_{n-1},t_n;L^2(\Sigma)^d)}+C\delta_4\|\partial_{tt}\boldsymbol{u}\|^2_{L^2(t_{n-1},t_n;L^2(\Omega)^d)}\\
&+\frac{\rho_s\epsilon}{2\delta_8}\|\partial_{tt}\boldsymbol{d}\|^2_{L^\infty(0,T;L^2(\Sigma)^d)}+\frac{\rho_s\epsilon}{2\delta_7}\|\partial_{tt}\boldsymbol{d}\|^2_{L^\infty(0,T;L^2(\Sigma)^d)}+\frac{(\rho_s\epsilon)^3}{\delta_6}\|\partial_{tt}\boldsymbol{d}\|^2_{L^\infty(0,T;L^2(\Sigma)^d)}\Big]\\
&\qquad+\tau\Big[\rho_{\mathrm{f}}\tau\delta_1\|\boldsymbol{\theta}_h^n\|^2_{0,\Omega}+\rho_s\epsilon\tau\delta_2\|\dot{\boldsymbol{\xi}}_h^n\|^2_{0,\Sigma}+\alpha^5\left(1+\frac{\delta_0}{2}+\frac{\delta_7}{2}+\frac{\delta_{10}^B}{2}\right)\|\boldsymbol{\xi}_h^n\|^2_{\mathrm{s}}\Big].
\end{aligned}\tag{4.92}
$$

In order to obtain the cumulative error after "$n$" time steps, we sum over $m=1..n$, we set the parameters in such a way that all the quantities at the left hand side result positives

$$
\delta_0=\frac{2}{3},\qquad \delta_1=\delta_2=\frac{1}{2},\qquad \delta_3=\frac{\mu}{C_K},\qquad \delta_4=\min\left\{\frac{1}{C},\frac{\rho_{\mathrm{f}}}{2C},\frac{\mu}{2C}\right\},\qquad \delta_5=\frac{1}{8(\rho_s\epsilon)^2},
$$
$$
\delta_6=1,\qquad \delta_7=\frac{2}{3},\qquad \delta_8=\frac{1}{2},\qquad \delta_9=\frac{1}{4},\qquad \delta_{10}^B=\frac{2}{3}.
$$

Using this values we obtain the following expression of the cumulative error at time step "$n$"

$$
\begin{aligned}
&\underbrace{\left(\frac{\rho_{\mathrm{f}}}{2}\|\boldsymbol{\theta}_h^n\|^2_{0,\Omega}+\frac{\rho_s\epsilon}{2}\|\dot{\boldsymbol{\xi}}_h^n\|^2_{0,\Sigma}+\frac{1}{2}\|\boldsymbol{\xi}_h^n\|^2_{\mathrm{s}}\right)}_{\mathscr{E}_h^n}-\underbrace{\left(\frac{\rho_{\mathrm{f}}}{2}\|\boldsymbol{\theta}_h^0\|^2_{0,\Omega}+\frac{\rho_s\epsilon}{2}\|\dot{\boldsymbol{\xi}}_h^0\|^2_{0,\Sigma}+\frac{1}{2}\|\boldsymbol{\xi}_h^0\|^2_{\mathrm{s}}\right)}_{\mathscr{E}_h^0}\\
&\qquad+\tau\sum_{m=1}^{n}\Big[C\|\boldsymbol{\epsilon}(\boldsymbol{\theta}_h^m)\|^2_{0,\Omega}+C|\varphi_h^m|^2_{s_h}\Big]\\
&+\sum_{m=1}^{n}C\|\tau\partial_\tau\boldsymbol{\theta}_h^m\|^2_{0,\Omega}+C\|\tau\partial_\tau\dot{\boldsymbol{\xi}}_h^m\|^2_{0,\Sigma}+\frac{1}{2}\|\tau\partial_\tau\boldsymbol{\xi}_h^m\|^2_{\mathrm{s}}+C\|\tau\partial_\tau\dot{\boldsymbol{\xi}}_h^m\|^2_{\mathrm{s}}\\
&\le\\
&\underbrace{\sum_{m=1}^{n}C_m^0}_{A_0}+\tau\underbrace{\sum_{m=1}^{n}C_m^1}_{A_1}+\tau^2\underbrace{\sum_{m=1}^{n}C_m^2}_{A_2}+\tau\underbrace{\sum_{m=1}^{n}\Big[\frac{\rho_{\mathrm{f}}}{2}\|\boldsymbol{\theta}_h^m\|^2_{,\Omega}+\frac{\rho_s\epsilon}{2}\|\dot{\boldsymbol{\xi}}_h^m\|^2_{0,\Sigma}+2\alpha^5\|\boldsymbol{\xi}_h^m\|^2_{\mathrm{s}}\Big]}_{A_3}.
\end{aligned}\tag{4.93}
$$

where the terms in the right hand side of the inequality are defined ad analyzed in the following using results of Proposition 3.3.6.

– The terms $A_0$ can be bounded as follows

$$
\begin{aligned}
A_0 \overset{def}{=} C \sum_{n=1}^{m} \Big( & \|\partial_t \boldsymbol{\theta}_\pi\|^2_{L^2(t_{n-1},t_n;L^2(\Omega)^d)} + \|\partial_t \dot{\boldsymbol{\xi}}_\pi\|^2_{L^2(t_{n-1},t_n;L^2(\Sigma)^d)} \\
& + \|\dot{\boldsymbol{\xi}}_\pi\|_{L^\infty(0,T;H^{\frac{1}{2}}(\Sigma)^d)} + \|\boldsymbol{\theta}_\pi\|_{L^\infty(0,T;H^1(\Omega)^d)} \\
& + \|\partial_t \boldsymbol{\xi}_\pi\|^2_{L^\infty(0,T;L^2(\Sigma)^d)} + \|\dot{\boldsymbol{\xi}}_\pi\|^2_{L^\infty(0,T;L^2(\Sigma)^d)} \Big) \\
\le C \Big( & \|\partial_t \boldsymbol{\theta}_\pi\|^2_{L^2(0,T;L^2(\Omega)^d)} + \|\partial_t \dot{\boldsymbol{\xi}}_\pi\|^2_{L^2(0,T;L^2(\Sigma)^d)} \\
& + \|\dot{\boldsymbol{\xi}}_\pi\|_{L^\infty(0,T;H^{\frac{1}{2}}(\Sigma)^d)} + \|\boldsymbol{\theta}_\pi\|_{L^\infty(0,T;H^1(\Omega)^d)} \\
& + \|\partial_t \boldsymbol{\xi}_\pi\|^2_{L^\infty(0,T;L^2(\Sigma)^d)} + \|\dot{\boldsymbol{\xi}}_\pi\|^2_{L^\infty(0,T;L^2(\Sigma)^d)} \Big) \\
& \approx \mathcal{O}(h_f^{2l+2}) + \mathcal{O}(h_s^{2m+2}) + \mathcal{O}(h_s^{2m+1}) + \mathcal{O}(h_f^{2l}) + \mathcal{O}(h_s^{2m+2}) + \mathcal{O}(h_s^{2m+2}).
\end{aligned}
\tag{4.94}
$$

– The term $A_1$ is defined and bounded as follows

$$
\begin{aligned}
A_1 \overset{def}{=} C \sum_{n=1}^{m} \Big( & \|\boldsymbol{\omega}_\pi\|^2_{L^\infty(0,T;\Lambda)} + \|\dot{\boldsymbol{\xi}}_\pi\|^2_{L^\infty(0,T;H^{\frac{1}{2}}(\Sigma)^d)} + \|\boldsymbol{\theta}_\pi\|^2_{L^\infty(0,T;H^1(\Omega)^d)} \Big) \\
& \approx \mathcal{O}(h_s^{2z+1}) + \mathcal{O}(h_s^{2m+1}) + \mathcal{O}(h_f^{2l}).
\end{aligned}
\tag{4.95}
$$

– The term $A_2$ is defined and bounded in the following using also the properties of the solid operator $\|\boldsymbol{L}_h^{\mathrm{s}} \boldsymbol{d}^n\|_{0,\Sigma} \le C \|\boldsymbol{L}^{\mathrm{s}} \boldsymbol{d}^n\|_{0,\Sigma}$

$$
\begin{aligned}
A_2 \overset{def}{=} C \sum_{n=1}^{m} \Big( & \|\partial_{tt} \boldsymbol{d}\|^2_{L^\infty(0,T;L^2(\Sigma)^d)} + \|\partial_{tt} \dot{\boldsymbol{d}}\|^2_{L^2(t_{n-1},t_n;L^2(\Sigma)^d)} \\
& + \|\partial_{tt} \boldsymbol{d}\|^2_{L^2(t_{n-1},t_n;\boldsymbol{W})} + \|\partial_{tt} \boldsymbol{u}\|^2_{L^2(t_{n-1},t_n;L^2(\Omega)^d)} \Big) \\
\le C \Big( & \|\partial_{tt} \boldsymbol{d}\|^2_{L^\infty(0,T;L^2(\Sigma)^d)} + \|\partial_{tt} \dot{\boldsymbol{d}}\|^2_{L^2(0,T;L^2(\Sigma)^d)} \\
& + \|\partial_{tt} \boldsymbol{d}\|^2_{L^2(0,T;\boldsymbol{W})} + \|\partial_{tt} \boldsymbol{u}\|^2_{L^2(0,T;L^2(\Omega)^d)} \Big).
\end{aligned}
\tag{4.96}
$$

– The terms $A_3$ is bounded as follows

$$
\begin{aligned}
A_3 &= \sum_{m=1}^{n} \left[ \frac{\rho_{\mathrm{f}}}{2} \|\boldsymbol{\theta}_h^m\|^2_{0,\Omega} + \frac{\rho_s \epsilon}{2} \|\dot{\boldsymbol{\xi}}_h^m\|^2_{0,\Sigma} + 2\alpha^5 \|\boldsymbol{\xi}_h^m\|^2_{\mathrm{s}} \right] \\
&\le \sum_{m=1}^{n} \gamma_m \underbrace{\left( \frac{\rho_{\mathrm{f}}}{2} \|\boldsymbol{\theta}_h^m\|^2_{0,\Omega} + \frac{\rho_s \epsilon}{2} \|\dot{\boldsymbol{\xi}}_h^m\|^2_{0,\Sigma} + \frac{1}{2} \|\boldsymbol{\xi}_h^m\|^2_{\mathrm{s}} \right)}_{\mathscr{E}_h^m},
\end{aligned}
\tag{4.97}
$$

where $\gamma_m$ is given by

$$
\gamma_m = \max\{1, 4\alpha^5\}
$$

In sum the error estimate (4.93) can be written as follows

$$
\mathscr{E}_h^n + \tau \mathscr{D}_h^n \le \tau \sum_{m=1}^{n} \gamma_m \mathscr{E}_h^m + A_0 + \tau A_1 + \tau^2 A_2 + \tau^3 A_3 \tag{4.98}
$$

where $\tau\mathscr{D}_h^n$ represents the dissipation term and is given by

$$\tau\mathscr{D}_h^n = \tau \sum_{m=1}^{n} \Big[ C\|\boldsymbol{\epsilon}(\boldsymbol{\theta}_h^m)\|_{0,\Omega}^2 + C|\varphi_h^m|_{s_h}^2 \Big] + \sum_{m=1}^{n} C\|\tau\partial_\tau\boldsymbol{\theta}_h^m\|_{0,\Omega}^2 + C\|\tau\partial_\tau\dot{\boldsymbol{\xi}}_h^m\|_{0,\Sigma}^2 + \frac{1}{2}\|\tau\partial_\tau\boldsymbol{\xi}_h^m\|_{\mathrm{s}}^2 + C\|\tau\partial_\tau\dot{\boldsymbol{\xi}}_h^m\|_{\mathrm{s}}^2.$$

Introducing the convergence rates with respect to $h_f$ and $h_s$ of the term $A_0, A_1, A_2, A_3$ given in (4.94)-(4.95), we have the desired result applying the discrete Gronwall Lemma 3.3.1, since by hypothesis we have $\tau\gamma_m < 1$

$$\mathscr{E}_h^n + \tau\mathscr{D}_h^n \lesssim \exp\left(\sum_{m=1}^{n} \frac{\tau}{1-\tau}\right)\Big[\underbrace{\mathscr{O}(h_f^{2l}) + \mathscr{O}(h_s^{2m+1})}_{A_0} \tag{4.99}$$

$$+ \tau\Big(\underbrace{\mathscr{O}(h_s^{2m+1}) + \mathscr{O}(h_s^{2z+1}) + \mathscr{O}(h_f^{2l})}_{A_1}\Big) + \tau^2 A_2 + \tau^3 A_3\Big]. \tag{4.100}$$

□

## 4.5 Numerical experiments

In this section we perform numerical simulations in order to test the behavior of the partitioned Algorithm 3, whose theoretical analysis was given in the previous sections. All the numerical tests are performed using the classical benchmark problem of an ellipsoidal structure that evolves to a circular equilibrium position. At the end of this section, we present some snapshots showing the time evolution of the system in the cases of extrapolation zero, one, two (see Figures 4.11, 4.12 and 4.13). We point out that in the case of extrapolation of order zero, in order to have meaningful numerical results, it is necessary to use finer mesh and time-step sizes with respect to the other extrapolations. We begin testing the stability of the algorithm with the different extrapolations.

### 4.5.1 Stability Analysis

Our first objective is to test the stability of the partitioned scheme, in order to verify the results of Theorem 4.3.3 regarding the unconditional stability of Algorithm 3 for extrapolation order $r = 0$ and $r = 1$ and conditional stability for extrapolation of order $r = 2$. In the following figures we report the time evolution of the energy of the system for different values of the time-step size $\tau$ and for fixed fluid and solid mesh sizes $h_f = h_s = 1/64$.

- case $r = 0$

  In the case $\boldsymbol{d}_h^\star = \boldsymbol{0}$, Algorithm 3 was proved to be unconditionally stable. In fact, the numerical experiments are in agreement with this theoretical result, as it is possible to appreciate in the following Figure 4.2 showing the time evolution of the total energy of the system.

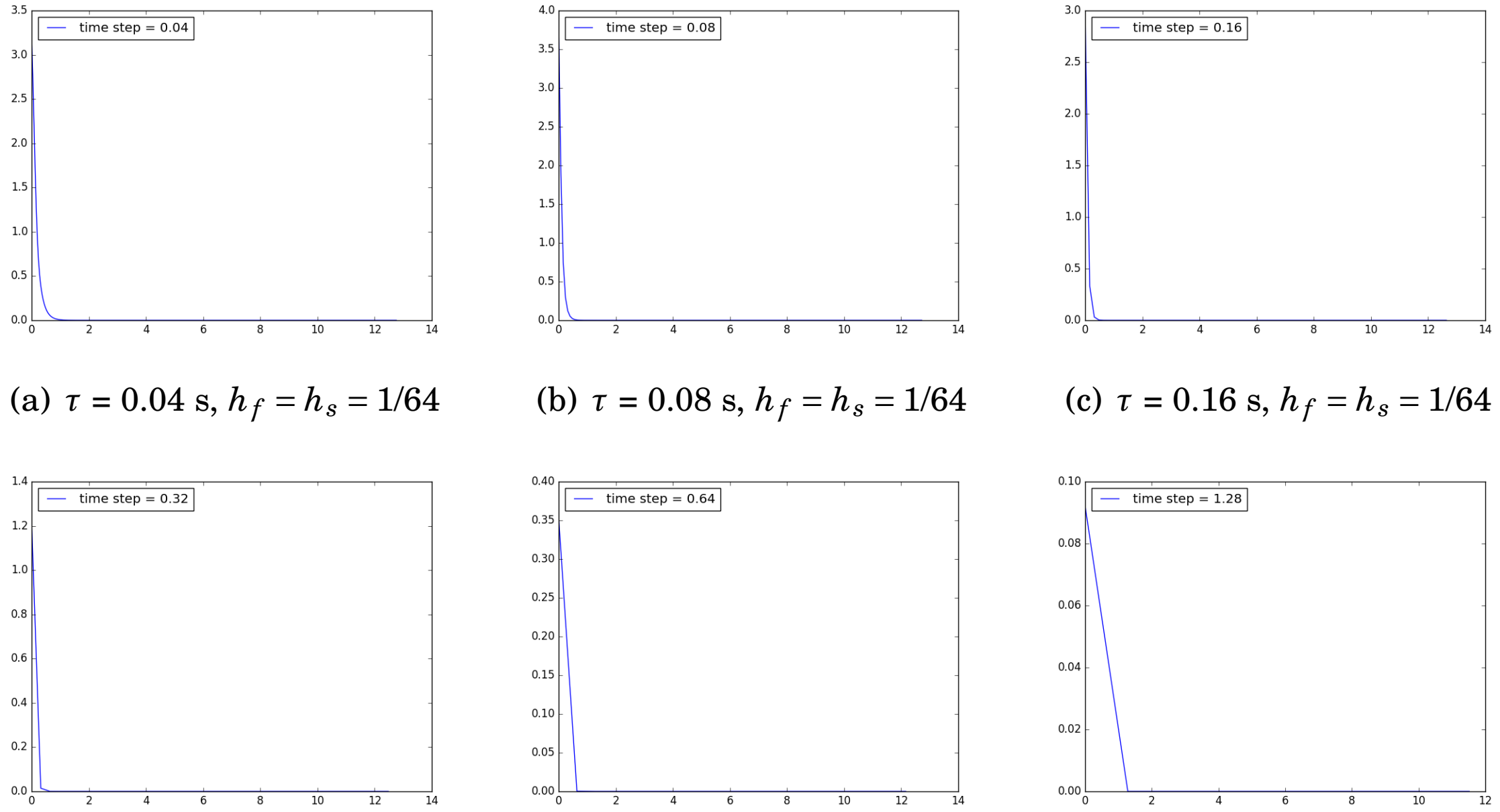

(a) $\tau$ = 0.04 s, $h_f = h_s = 1/64$ (b) $\tau$ = 0.08 s, $h_f = h_s = 1/64$ (c) $\tau$ = 0.16 s, $h_f = h_s = 1/64$

(d) $\tau$ = 0.32 s, $h_f = h_s = 1/64$ (e) $\tau$ = 0.64 s, $h_f = h_s = 1/64$ (f) $\tau$ = 1.28 s, $h_f = h_s = 1/64$

Figure 4.2: Stability of Algorithm 3, extrapolation $\boldsymbol{d}_h^\star = \boldsymbol{0}$.

- case $r = 1$

In the case $\boldsymbol{d}_h^\star = \boldsymbol{d}_h^{n-1}$ Algorithm 3 was proved to be unconditionally stable and the numerical experiments are in agreement with this theoretical result. In the following Figure 4.3 we report the results of the simulation performed on the test case. It is possible to appreciate that the energy of the system decreases during the simulation.

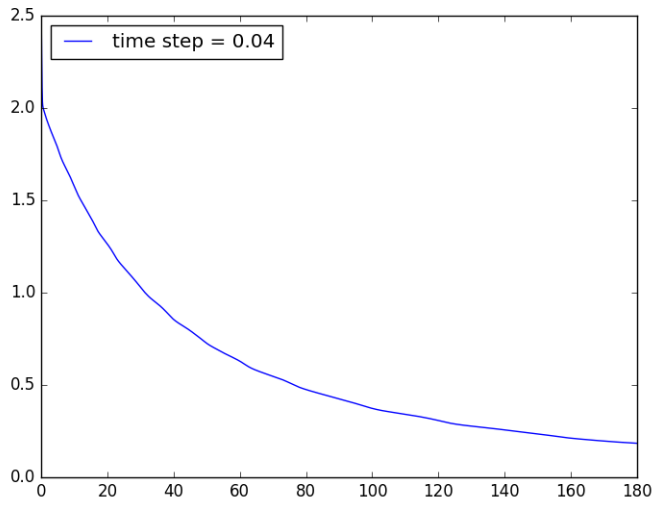

(a) $\tau$ = 0.04 s, $h_f = h_s = 1/64$

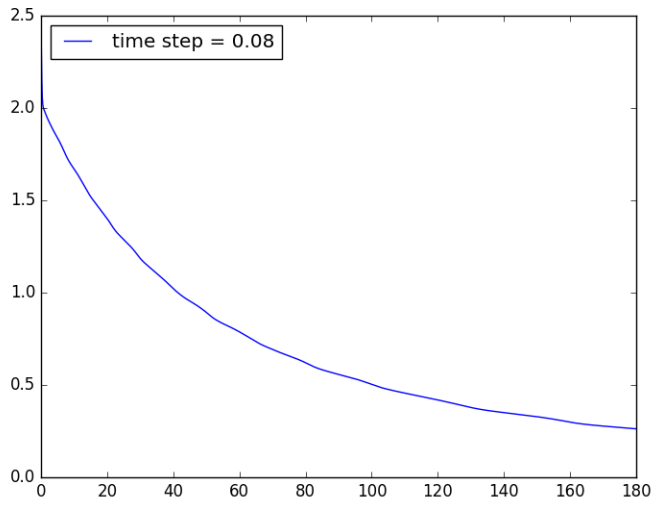

(b) $\tau$ = 0.08 s, $h_f = h_s = 1/64$

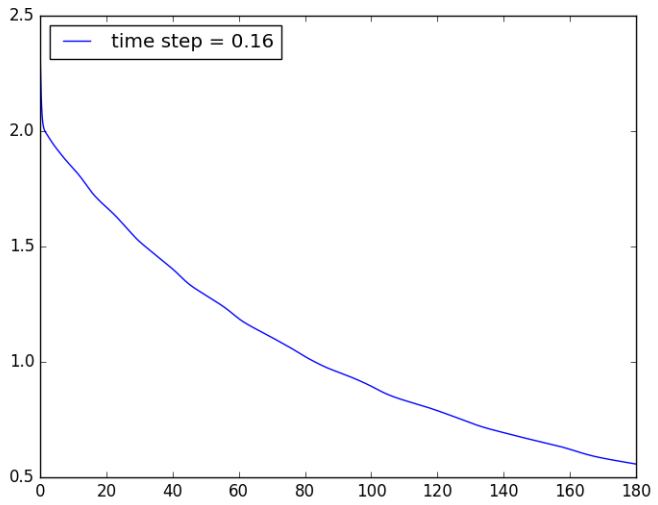

(c) $\tau$ = 0.16 s, $h_f = h_s = 1/64$

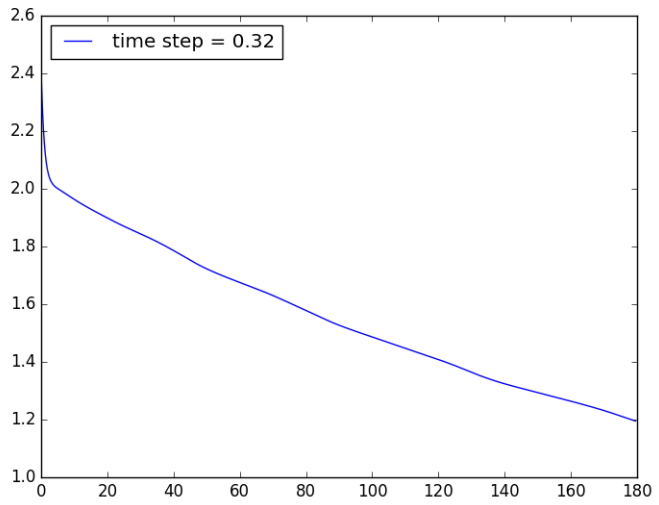

(d) $\tau$ = 0.32 s, $h_f = h_s = 1/64$

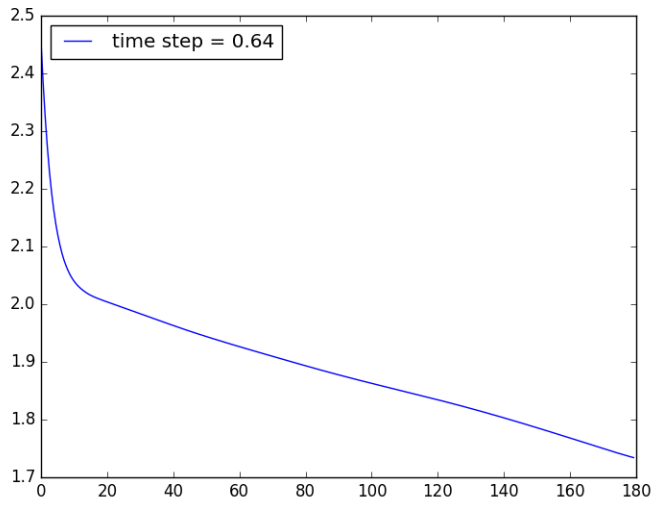

(e) $\tau$ = 0.64 s, $h_f = h_s = 1/64$

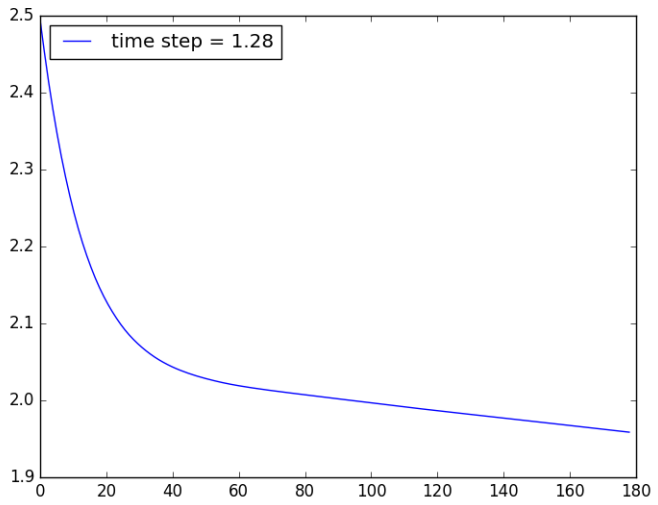

(f) $\tau$ = 1.28 s, $h_f = h_s = 1/64$

Figure 4.3: Stability of Algorithm 3, extrapolation $\boldsymbol{d}_h^\star = \boldsymbol{d}_h^{n-1}$.

- case $r = 2$ :

In the case $\boldsymbol{d}_h^\star = \boldsymbol{d}_h^{n-1} + \tau\dot{\boldsymbol{d}}_h^{n-1}$ Algorithm 3 was proved to be conditionally stable, in fact for time-step sizes $\tau = 1.28$ and $\tau = 0.64$ we observe fluctuations of the energy that are related to solutions that are in some sense not stable, meaning that in this case the system does not reach the rest position as in the other cases. Despite this behavior, that indicate non meaningful solution, the energy remain bounded as it is shown in Figure 4.4.

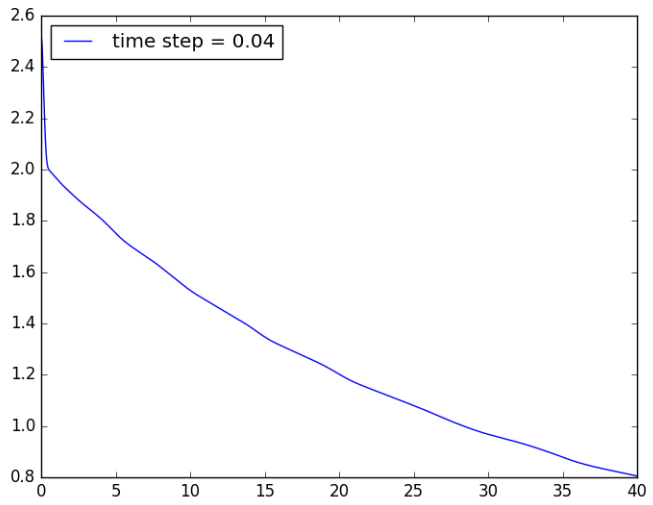

(a) $\tau$ = 0.04 s, $h_f = h_s = 1/64$

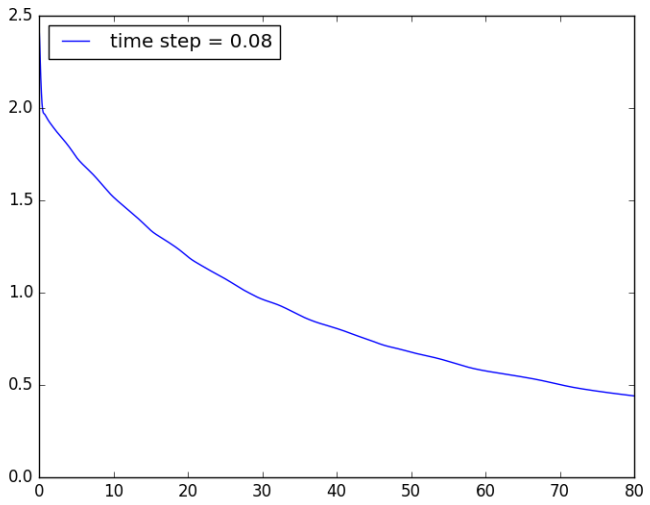

(b) $\tau$ = 0.08 s, $h_f = h_s = 1/64$

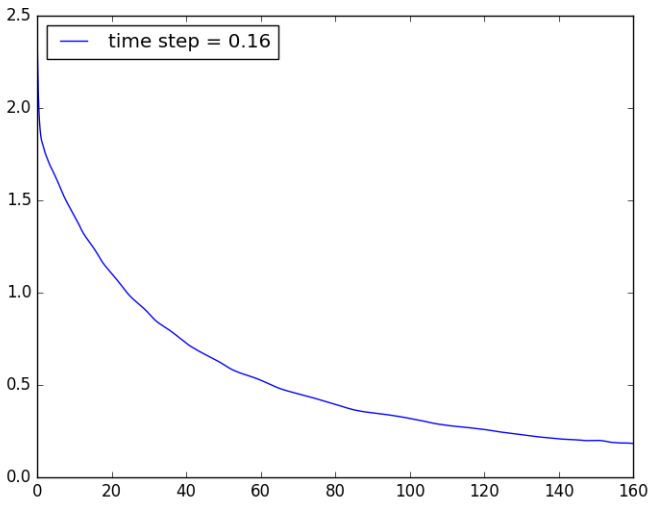

(c) $\tau$ = 0.16 s, $h_f = h_s = 1/64$

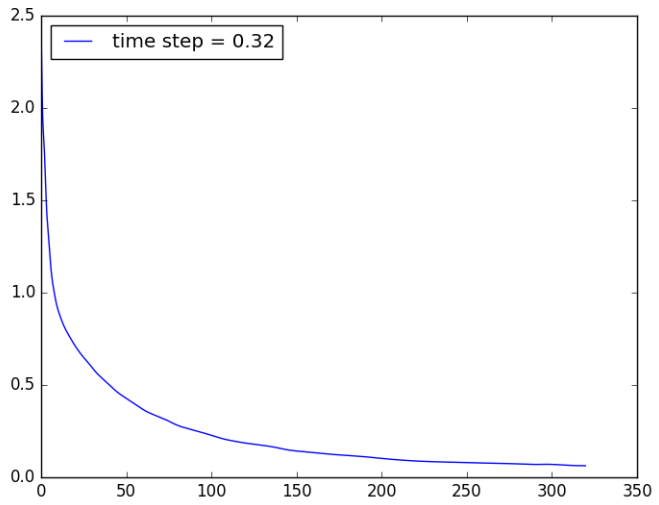

(d) $\tau$ = 0.32 s, $h_f = h_s = 1/64$

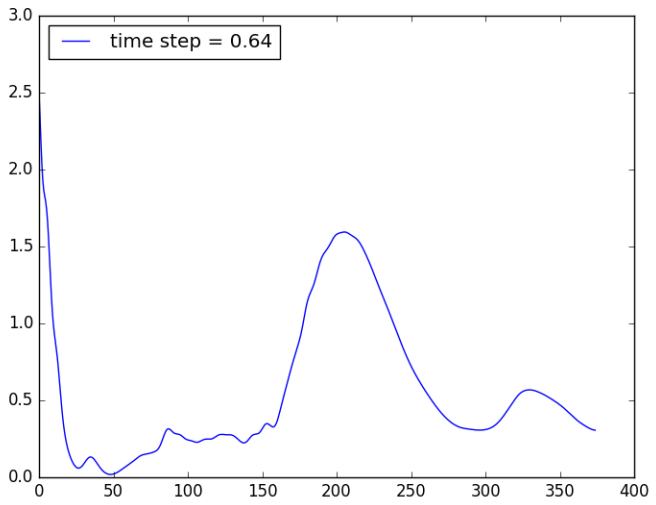

(e) $\tau$ = 0.64 s, $h_f = h_s = 1/64$

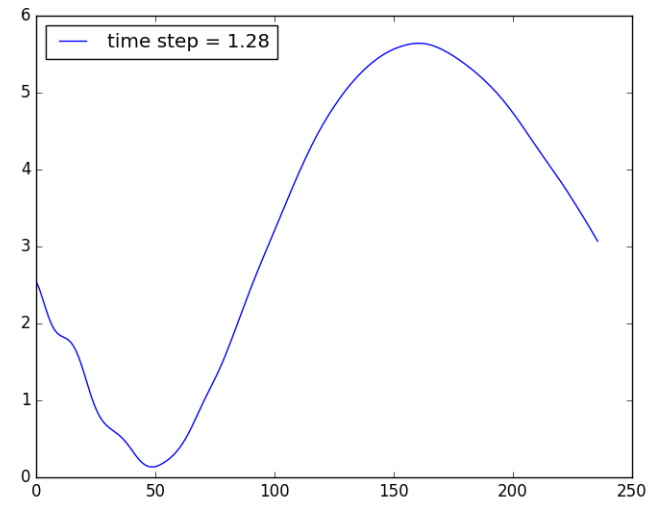

(f) $\tau$ = 1.28 s, $h_f = h_s = 1/64$

Figure 4.4: Stability of Algorithm 3, extrapolation $\boldsymbol{d}_h^\star = \boldsymbol{d}_h^{n-1} + \tau \dot{\boldsymbol{d}}_h^{n-1}$.

### 4.5.2 Convergence analysis

In this paragraph we consider the error analysis of the partitioned Algorithm 3. In order to verify the convergence estimates given in Theorem 4.4.2, we consider the errors evaluated with respect to reference solution obtained using the monolithic algorithm with discretization parameters 3.62.

**Convergence with respect to time** In order to investigate the convergence with respect to time-step size we consider decreasing values of $\tau$ and very fines spatial meshes for the solid and the fluid. The discretization parameters are chosen in such a way that the effect of the error due to spatial discretization is negligible.

$$h_{\rm s} = h_{\rm f} = 1/256, \qquad \tau \in \{0.004, 0.008, 0.016, 0.032, 0.064\}.$$

In the following tables we report the errors, with respect to the reference solution 3.62, computed at time $t = 0.064s$ for the various values of $\tau$ and fixed values of the spatial mesh sizes $h_s = h_f = 1/256$. We report moreover the figures in which we plot the total error as function of the time-step size.

Using the convergence estimates of Theorem 4.4.2 in the case $h_f = h_s$ we expect that the total error

$$\text{Total Error} \stackrel{def}{=} \mathscr{E}_h^n = C\left(\|\boldsymbol{u}_h^n - \boldsymbol{u}_{ref}^n\|_{0,\Omega}^2 + \|\boldsymbol{d}_h^n - \boldsymbol{d}_{ref}^n\|_{s,\Sigma}^2 + \|\dot{\boldsymbol{d}}_h^n - \dot{\boldsymbol{d}}_{ref}^n\|_{0,\Sigma}^2\right)^{\frac{1}{2}}$$

goes to zero as $\tau$.

- case $r = 0$

  For the extrapolation of order zero, meaning for $\boldsymbol{d}_h^{\star} = \mathbf{0}$, the convergence rate given by the numerical simulations corresponds to the one expected from the theoretical analysis

| $\tau$ | $\\|\boldsymbol{u}_h^n - \boldsymbol{u}_{ref}^n\\|_{0,\Omega}$ | $\\|\boldsymbol{d}_h^n - \boldsymbol{d}_{ref}^n\\|_{s,\Sigma}$ | $\\|\dot{\boldsymbol{d}}_h^n - \dot{\boldsymbol{d}}_{ref}^n\\|_{0,\Sigma}$ | Total Error | Total Rate |
|---|---|---|---|---|---|
| 0.064 | 0.0294989805 | 0.2649662054 | 0.5315079303 | 0.5946242173 | |
| 0.032 | 0.0165615557 | 0.1097510364 | 0.2998033689 | 0.3196899046 | 8.95E-01 |
| 0.016 | 0.0086397675 | 0.0473799536 | 0.1568161663 | 0.1640451632 | 9.63E-01 |
| 0.008 | 0.0043766559 | 0.0216450918 | 0.0798075811 | 0.082806492 | 9.86E-01 |
| 0.004 | 0.0022085527 | 0.0103633971 | 0.0402345623 | 0.0416064623 | 9.93E-01 |

Table 4.1: Convergence w.r.t time-step size of partitioned Algorithm 3 at time $t = 0.064s$, for $h_s = h_f = 1/256$.

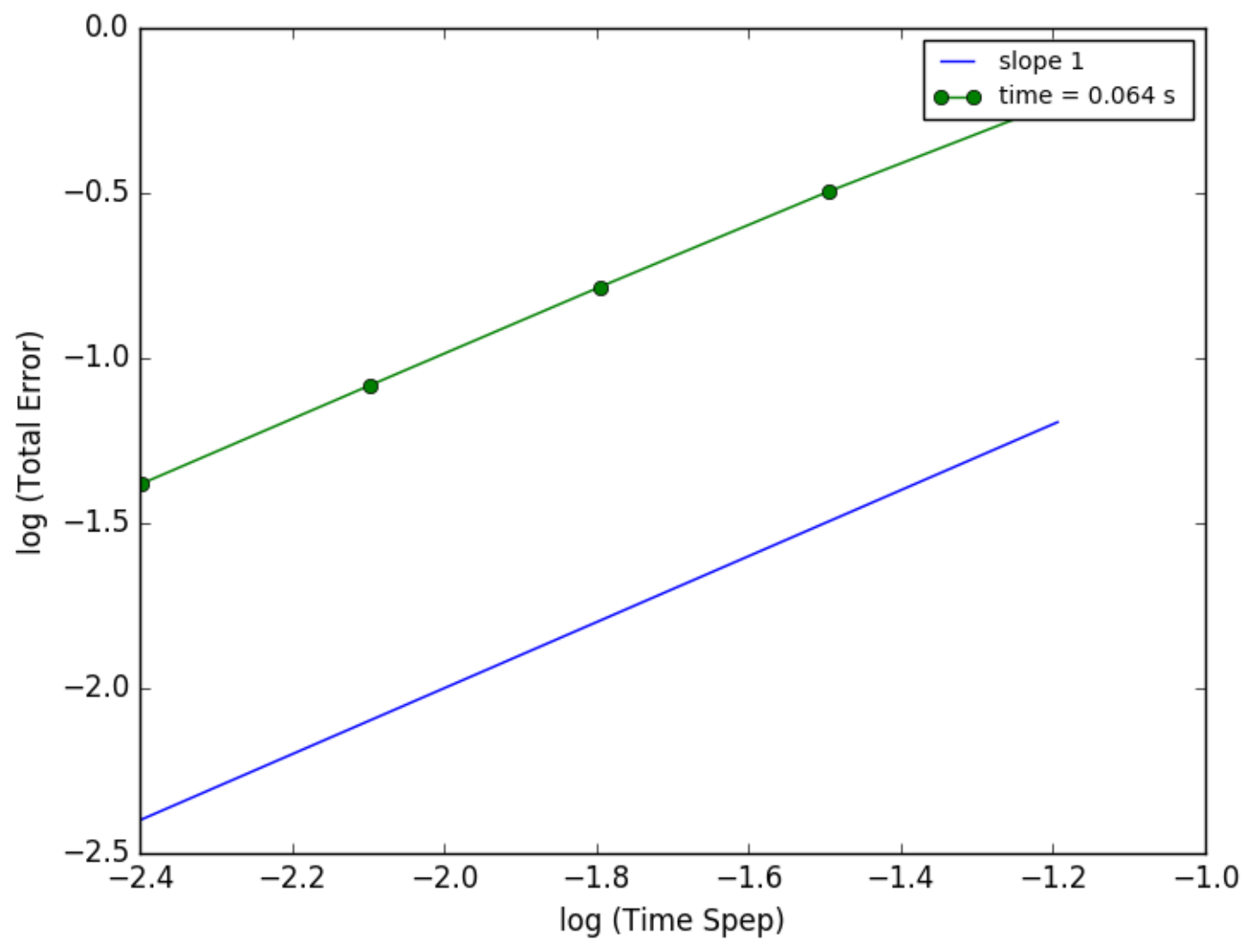

(a) $t = 0.064s$, $h_f = h_s = 1/256$

Figure 4.5: Convergence with respect to time-step size of Algorithm 3, extrapolation $\boldsymbol{d}_h^{\star} = \mathbf{0}$.

- Case $r = 1$

  For the extrapolation of order one, meaning $\boldsymbol{d}_h^{\star} = \boldsymbol{d}_h^{n-1}$, we expected a convergence rate, with respect to time-step size, of order 1; the numerical simulations seems do be in agreement with the theoretical analysis

| $\tau$ | $\lVert \boldsymbol{u}_h^n - \boldsymbol{u}_{ref}^n \rVert_{0,\Omega}$ | $\lVert \boldsymbol{d}_h^n - \boldsymbol{d}_{ref}^n \rVert_{s,\Sigma}$ | $\lVert \dot{\boldsymbol{d}}_h^n - \dot{\boldsymbol{d}}_{ref}^n \rVert_{0,\Sigma}$ | Total Error | Total Rate |
|---|---|---|---|---|---|
| 0.064 | 0.0106936667 | 0.0554678227 | 0.0162457138 | 0.058778883 | |
| 0.032 | 0.0066462263 | 0.0312720022 | 0.0084203747 | 0.0330607495 | 8.30E-01 |
| 0.016 | 0.003750269 | 0.0164390978 | 0.0045079773 | 0.0174536619 | 9.22E-01 |
| 0.008 | 0.001959067 | 0.0082047341 | 0.0023974795 | 0.0087694648 | 9.93E-01 |
| 0.004 | 0.0009522398 | 0.0038842392 | 0.0011954948 | 0.0041741206 | 1.07E+00 |

Table 4.2: Convergence w.r.t time-step size of partitioned Algorithm 3 at time $t = 0.064s$, for $h_s = h_f = 1/256$.

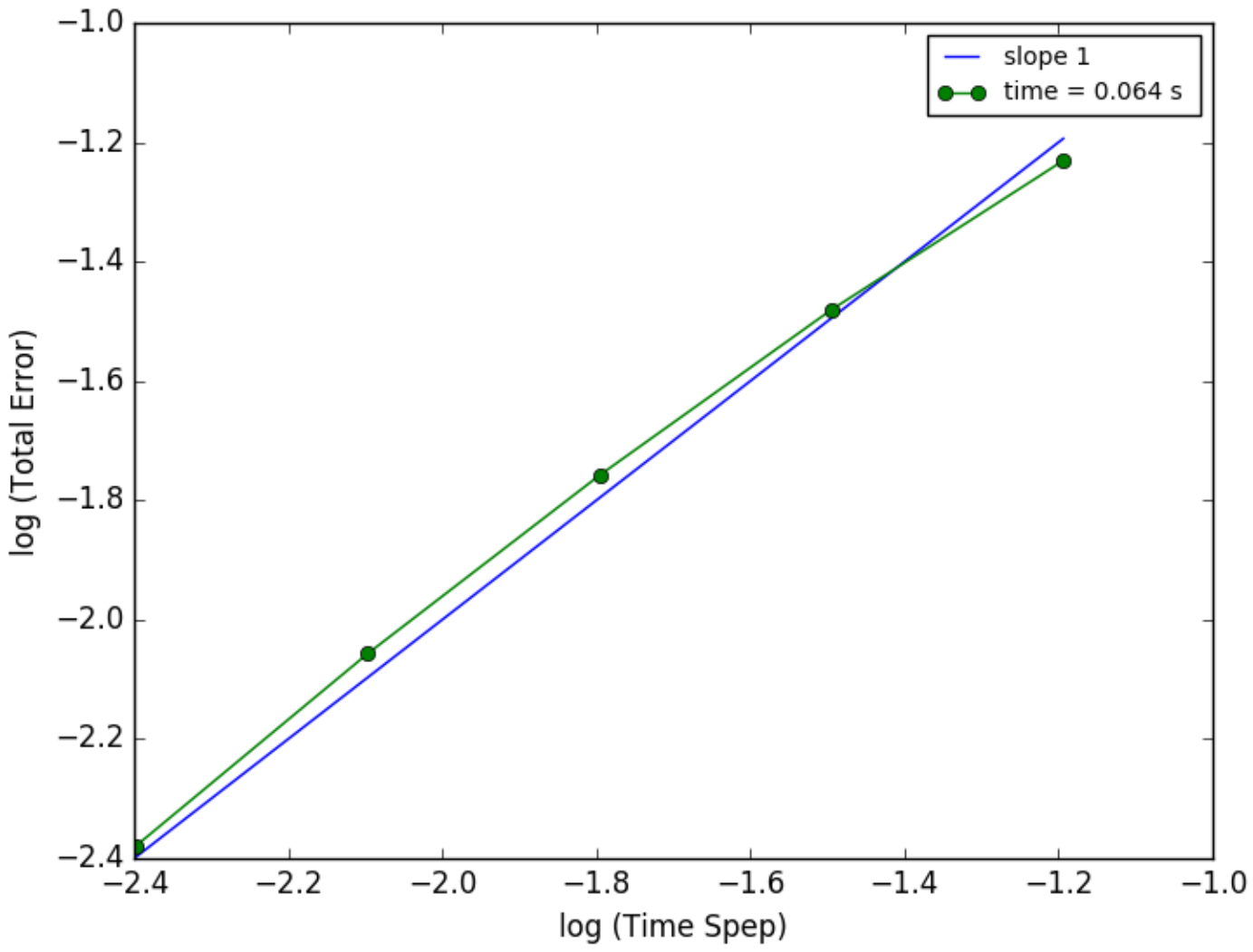

(a) $t = 0.064s$, $h_f = h_s = 1/256$

Figure 4.6: Convergence with respect to time-step size of Algorithm 3, extrapolation $\boldsymbol{d}_h^\star = \boldsymbol{d}_h^{n-1}$.

- Case $r = 2$

  For extrapolation $\boldsymbol{d}_h^\star = \boldsymbol{d}_h^{n-1} + \tau \dot{\boldsymbol{d}}_h^{n-1}$, the scheme was proven to be conditionally stable. In the case that the stability conditions are fulfilled, the theoretical analysis predicted a convergence rate, with respect to the time-step size, of order 1. The numerical experiments seems to confirm the theoretical analysis.

| $\tau$ | $\|\boldsymbol{u}_h^n - \boldsymbol{u}_{ref}^n\|_{0,\Omega}$ | $\|\boldsymbol{d}_h^n - \boldsymbol{d}_{ref}^n\|_{s,\Sigma}$ | $\|\dot{\boldsymbol{d}}_h^n - \dot{\boldsymbol{d}}_{ref}^n\|_{0,\Sigma}$ | Total Error | Total Rate |
|---|---|---|---|---|---|
| 0.064 | 0.0145145369 | 0.0617916977 | 0.0081849754 | 0.0639990586 | |
| 0.032 | 0.0081647849 | 0.0339125999 | 0.0053601233 | 0.0352910621 | 8.59E-01 |
| 0.016 | 0.0040807811 | 0.0171783216 | 0.0036681733 | 0.0180333858 | 9.69E-01 |
| 0.008 | 0.0020232439 | 0.0084183654 | 0.0021876433 | 0.0089301833 | 1.01E+00 |
| 0.004 | 0.0009661696 | 0.0039985616 | 0.0011535881 | 0.004272323 | 1.06E+00 |

Table 4.3: Convergence w.r.t time-step size of partitioned Algorithm 3 at time $t = 0.064s$, for $h_s = h_f = 1/256$.

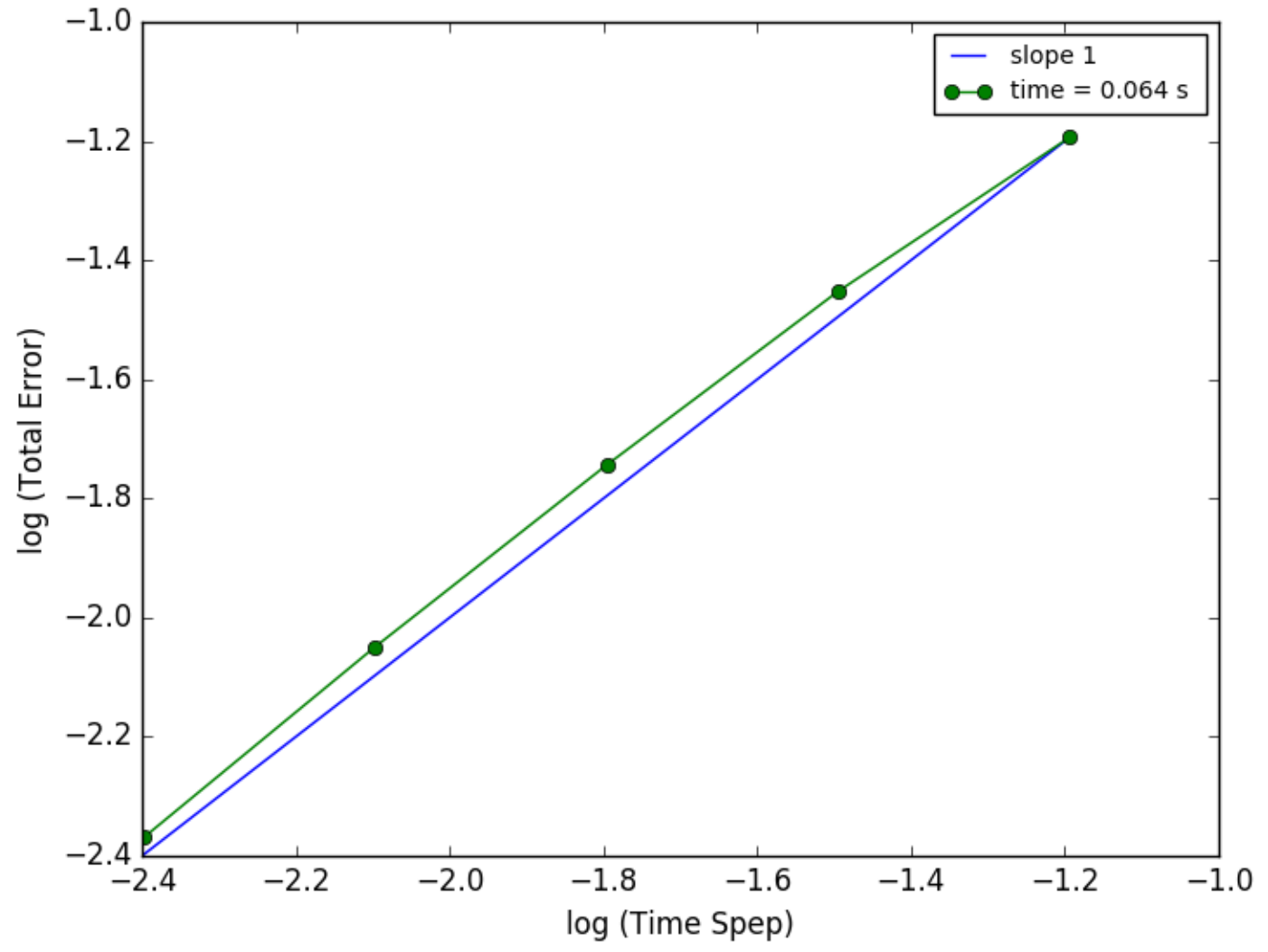

(a) $t = 0.064s$, $h_f = h_s = 1/256$

Figure 4.7: Convergence with respect to time-step size of Algorithm 3, extrapolation $\boldsymbol{d}_h^{\star} = \boldsymbol{d}_h^{n-1} + \tau\dot{\boldsymbol{d}}_h^{n-1}$.

In conclusion, the convergence analysis with respect to time-step size, of the partitioned Algorithm 3 seems to be in good agreement with the results of the numerical experiments. In the following paragraphs we will consider the convergence with respect to the spatial mesh sizes.

### 4.5.3 Convergence with respect to spatial mesh sizes

In this paragraph, we analyze the error behavior with respect to the fluid and solid mesh sizes, $h_f$ and $h_s$ respectively, for fixed value of the time mesh size. In order to perform the analysis, we consider the following discretization parameters

$$h_f = h_s \in \{1/8, 1/16, 1/32, 1/64, 1/128\}, \qquad \tau = 0.00005$$

We perform analysis considering the steady test case of a circular elastic string immersed in the fluid. The error is evaluated with respect to a reference solution obtained considering the same time step size $\tau = 0.00005$ and spatial mesh sizes $h_f = h_s = 1/256$. In the following we report the results of the numerical experiments evaluated for any extrapolation.

- Case $r = 0$

  In the following table we report the errors for the fluid velocity, the solid displacements and velocity, moreover we report the plot of the total error. From the table and figure it is possible to appreciate the total convergence rate that goes around "$h^l$", in good agreement with the theoretical analysis results.

| $h_f = h_s$ | $\|\boldsymbol{u}_h^n - \boldsymbol{u}_{ref}^n\|_{0,\Omega}$ | $\|\boldsymbol{d}_h^n - \boldsymbol{d}_{ref}^n\|_{s,\Sigma}$ | $\|\dot{\boldsymbol{d}}_h^n - \dot{\boldsymbol{d}}_{ref}^n\|_{0,\Sigma}$ | Total Error | Total Rate |
|---|---|---|---|---|---|
| 1/8 | 9.93998E-04 | 3.11651E-02 | 1.95490E-04 | 3.11816E-02 | |
| 1/16 | 8.67699E-04 | 1.62672E-02 | 1.94071E-04 | 1.62915E-02 | 9.36574E-01 |
| 1/32 | 7.61955E-04 | 8.64093E-03 | 1.97383E-04 | 8.67671E-03 | 9.08900E-01 |
| 1/64 | 6.07689E-04 | 4.80470E-03 | 1.88353E-04 | 4.84664E-03 | 8.40163E-01 |
| 1/128 | 3.63058E-04 | 2.88349E-03 | 1.65336E-04 | 2.91096E-03 | 7.35492E-01 |

Table 4.4: Convergence rate with respect to the spatial mesh size evaluated for $\boldsymbol{d}_h^\star = \boldsymbol{0}$, $\tau = 0.00005s$ at time $t = 0.0001s$.

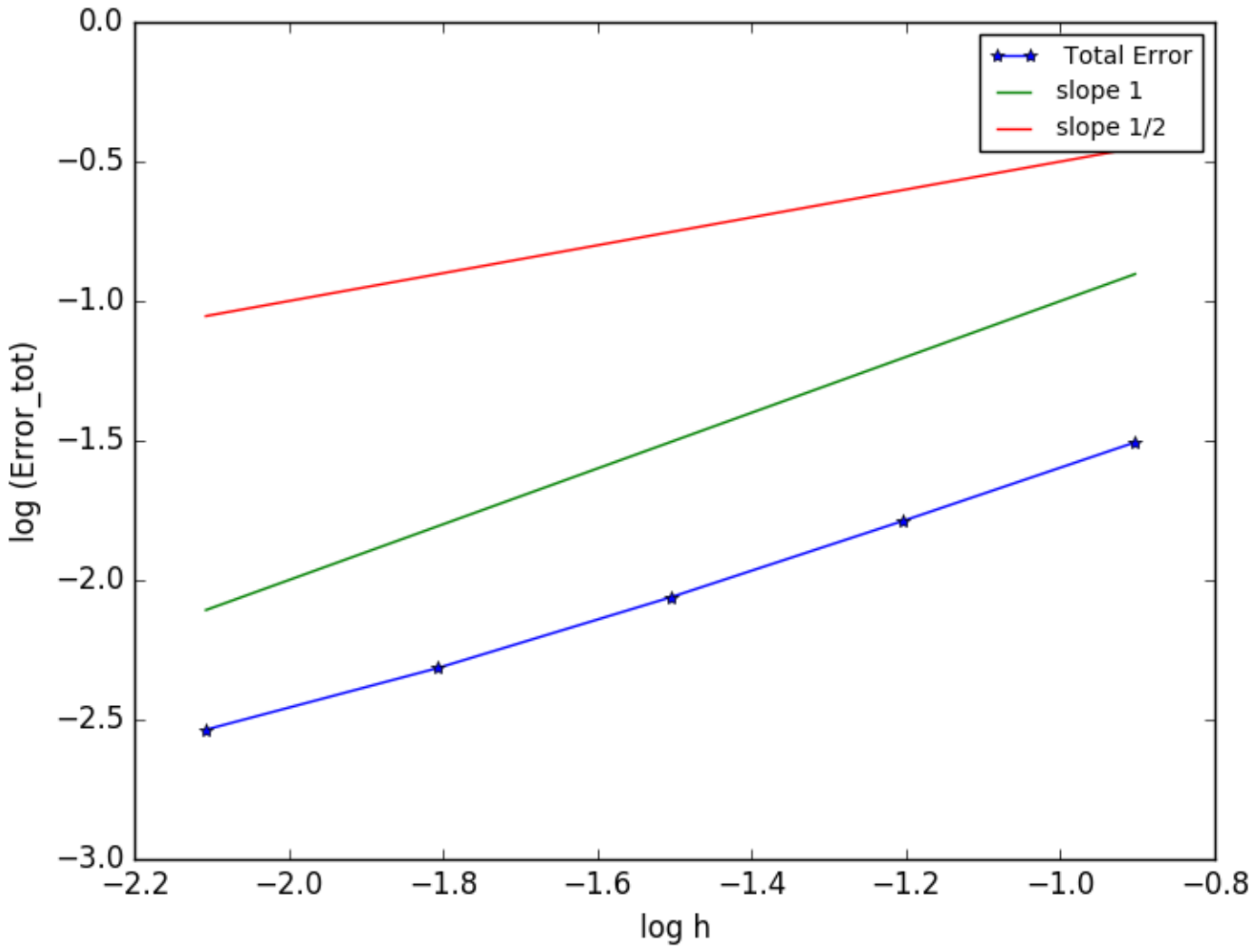

(a) $t = 0.0001ss$, $h_f = h_s$, $\tau = 0.00005s$

Figure 4.8: Convergence with respect to the spatial mesh size of Algorithm 3, extrapolation $\boldsymbol{d}_h^\star = \boldsymbol{0}$ evaluated at $t = 0.0001s$.

- Case r=1

  In the following table we report the errors for the fluid velocity, the solid displacements and velocity. In figure we plot the total error. It is possible to appreciate the total convergence rate, that, in agreement with the theoretical analysis results around "$h^l$".

| $h_f = h_s$ | $\|\boldsymbol{u}_h^n - \boldsymbol{u}_{ref}^n\|_{0,\Omega}$ | $\|\boldsymbol{d}_h^n - \boldsymbol{d}_{ref}^n\|_{s,\Sigma}$ | $\|\dot{\boldsymbol{d}}_h^n - \dot{\boldsymbol{d}}_{ref}^n\|_{0,\Sigma}$ | Total Error | Total Rate |
|---|---|---|---|---|---|
| 1/8 | 1.48533E-03 | 3.11651E-02 | 1.64340E-04 | 3.12009E-02 | |
| 1/16 | 1.33734E-03 | 1.62672E-02 | 1.64468E-04 | 1.63229E-02 | 9.34690E-01 |
| 1/32 | 1.14503E-03 | 8.64091E-03 | 1.68944E-04 | 8.71809E-03 | 9.04815E-01 |
| 1/64 | 8.59875E-04 | 4.80470E-03 | 1.50629E-04 | 4.88336E-03 | 8.36138E-01 |
| 1/128 | 4.72989E-04 | 2.88348E-03 | 1.11350E-04 | 2.92414E-03 | 7.39861E-01 |

Table 4.5: Convergence rate with respect to the spatial mesh size evaluated for $\boldsymbol{d}_h^\star = \boldsymbol{d}_h^{n-1}$, $\tau = 0.00005s$ at time $t = 0.0001s$.

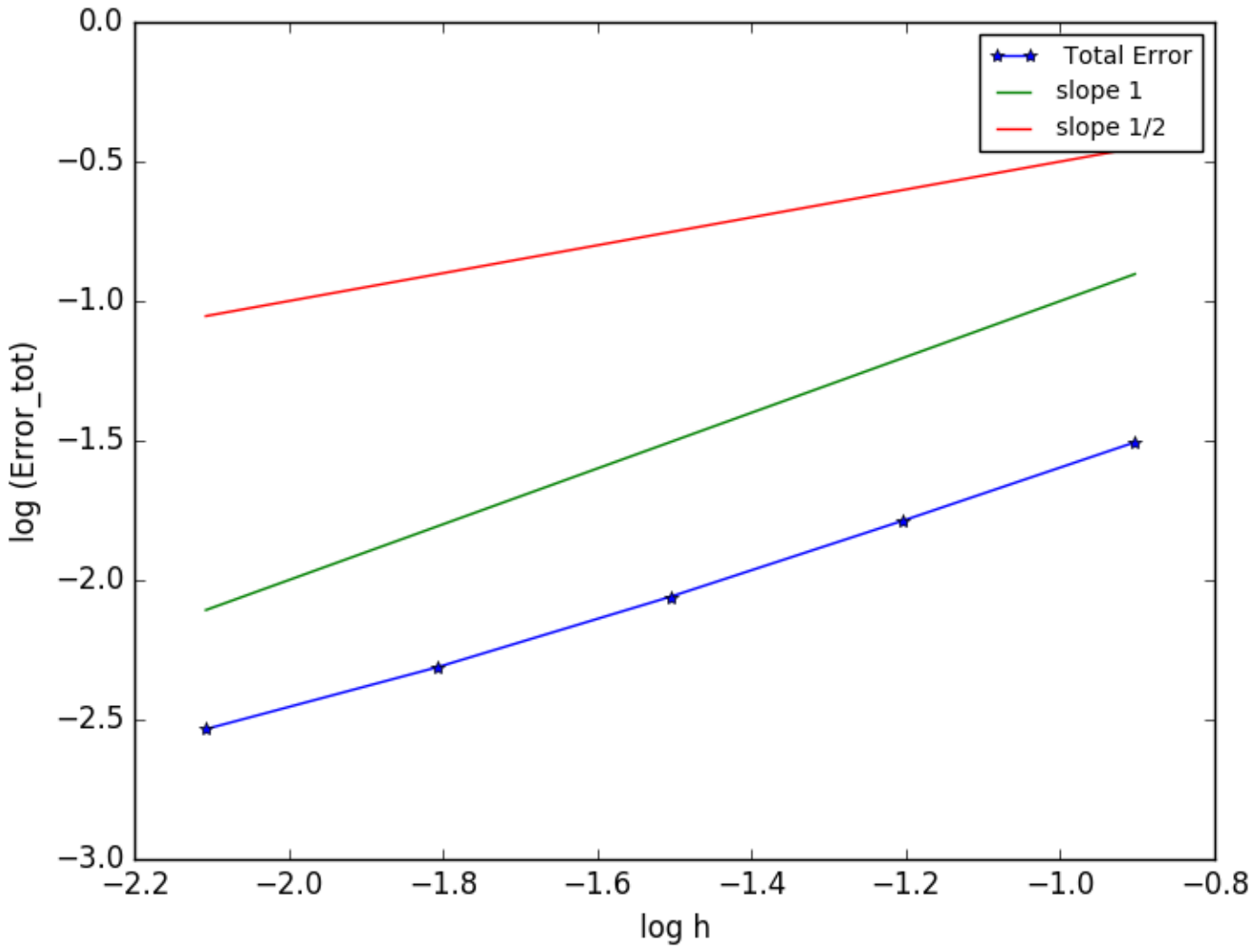

(a) $t = 0.0001ss$, $h_f = h_s$, $\tau = 0.00005s$

Figure 4.9: Convergence with respect to the spatial mesh size of Algorithm 3, extrapolation $\boldsymbol{d}_h^\star = \boldsymbol{d}_h^{n-1}$ evaluated at $t = 0.0001s$.

- Case r=2

  In the following table we report the errors for the fluid velocity, the solid displacements and velocity, moreover we plot the total error. In the tables and the figure it is possible to appreciate the total convergence rate. As expected from the theoretical analysis, the total error goes to zero about as "$h^l$".

| $h_f = h_s$ | $\|\boldsymbol{u}_h^n - \boldsymbol{u}_{ref}^n\|_{0,\Omega}$ | $\|\boldsymbol{d}_h^n - \boldsymbol{d}_{ref}^n\|_{s,\Sigma}$ | $\|\dot{\boldsymbol{d}}_h^n - \dot{\boldsymbol{d}}_{ref}^n\|_{0,\Sigma}$ | Total Error | Total Rate |
|---|---|---|---|---|---|
| 1/8 | 1.48533E-03 | 3.11651E-02 | 1.64340E-04 | 3.12009E-02 | |
| 1/16 | 1.33734E-03 | 1.62672E-02 | 1.64468E-04 | 1.63229E-02 | 9.34690E-01 |
| 1/32 | 1.14503E-03 | 8.64091E-03 | 1.68944E-04 | 8.71809E-03 | 9.04815E-01 |
| 1/64 | 8.59875E-04 | 4.80470E-03 | 1.50629E-04 | 4.88336E-03 | 8.36138E-01 |
| 1/128 | 4.72990E-04 | 2.88348E-03 | 1.11350E-04 | 2.92414E-03 | 7.39861E-01 |

Table 4.6: Convergence rate with respect to the spatial mesh size evaluated for $\boldsymbol{d}_h^\star = \boldsymbol{d}_h^{n-1} + \tau\dot{\boldsymbol{d}}_h^{n-1}$, $\tau = 0.00005s$ at time $t = 0.0001s$.

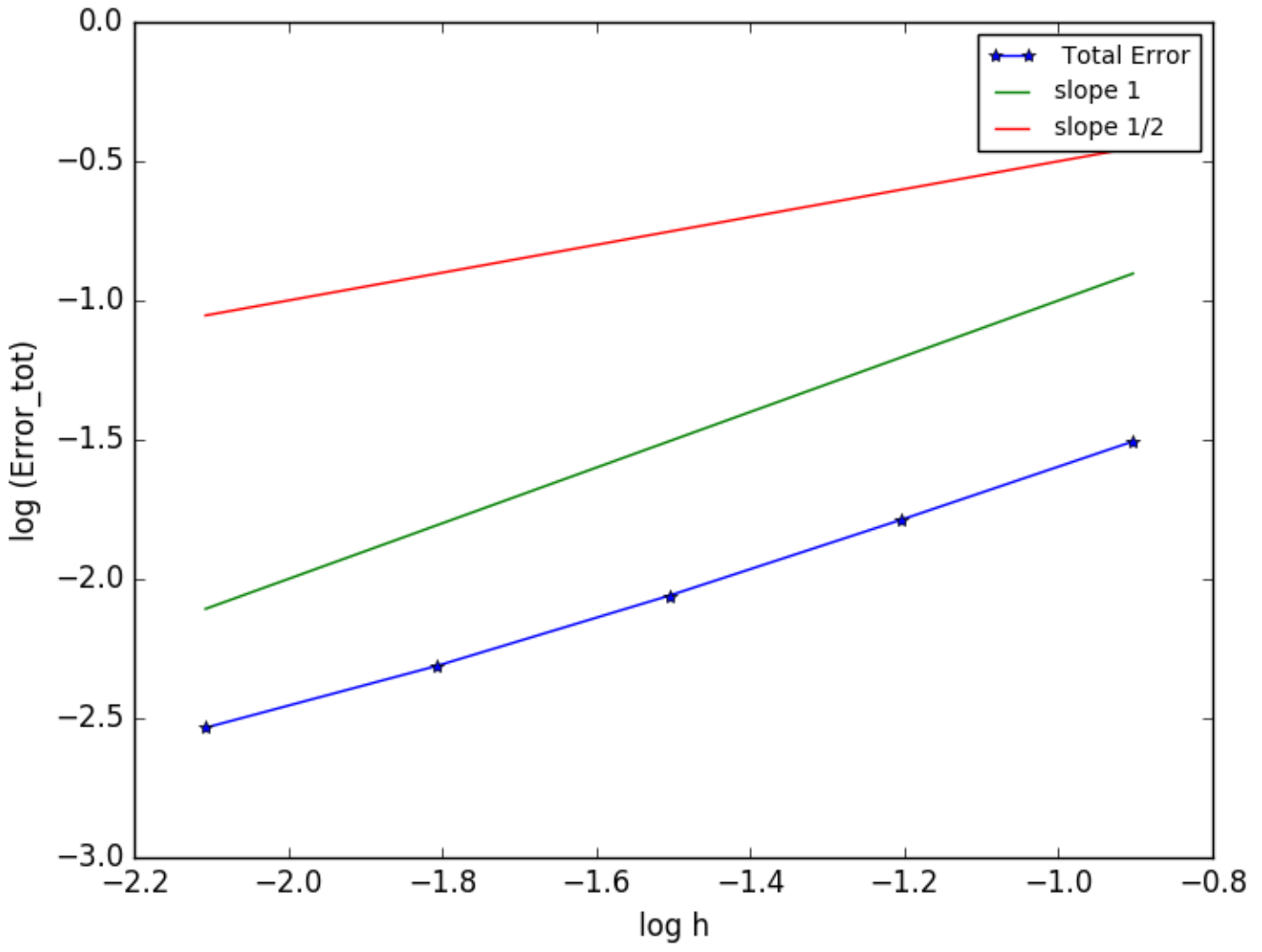

(a) $t = 0.0001ss$, $h_f = h_s$, $\tau = 0.00005s$

Figure 4.10: Convergence with respect to the spatial mesh size of Algorithm 3, extrapolation $\boldsymbol{d}_h^\star = \boldsymbol{d}_h^{n-1} + \tau\dot{\boldsymbol{d}}_h^{n-1}$ evaluated at $t = 0.0001s$.

### 4.5.4 Global convergence rate

In this paragraph, we analyze the error behavior when the discretization parameters $h_f, h_s, \tau$ go to zero simultaneously. In order to perform the analysis, we consider the following discretization parameters.

$$h_f = h_s = \frac{0.125}{2^j}, \qquad \tau = \frac{0.064}{2^j} \qquad \text{for} \qquad j = 0, 1, 2, 3, 4.$$

The error is evaluated with respect to a reference solution obtained using the monolithic algorithm with discretization parameters 3.62. For all the cases, the convergence rates for the fluid velocity is in agreement with the supposed regularity of the solution. In fact, since the pressure is discontinuous across the structure, the optimal convergence rate is $\frac{3}{2}$ for the fluid velocity in $L^2$. For the solid displacement the computed convergence rate close to 1 in the energy norm which is equivalent to the $H^1$-norm. This result is in agreement with the supposed regularity for the solid displacement. The convergence rate of the solid velocity results 1 evaluated in the $L^2$ norm. We report also the convergence rate of the total error that reaches in all the cases a value around 1, in agreement with the theoretical results.

- Case $r = 0$

  In this case we have $\boldsymbol{d}^{n\star} = \boldsymbol{0}$ and $\boldsymbol{d}_h^{n\star} = \boldsymbol{0}$. The accuracy of the split scheme with extrapolation $r = 0$ is evaluated in the following solving the steady test case of equilibrium of a circular elastic string immersed in the fluid.

| $h_f = h_s$ | $\\|\boldsymbol{u}_h^n - \boldsymbol{u}_{ref}^n\\|_{0,\Omega}$ | Rate |
|---|---|---|
| 1/8 | 1.41273E-02 | |
| 1/16 | 8.89585E-03 | 6.67283E-01 |
| 1/32 | 4.00575E-03 | 1.15106E+00 |
| 1/64 | 1.60436E-03 | 1.32007E+00 |
| 1/128 | 6.03243E-04 | 1.41119E+00 |

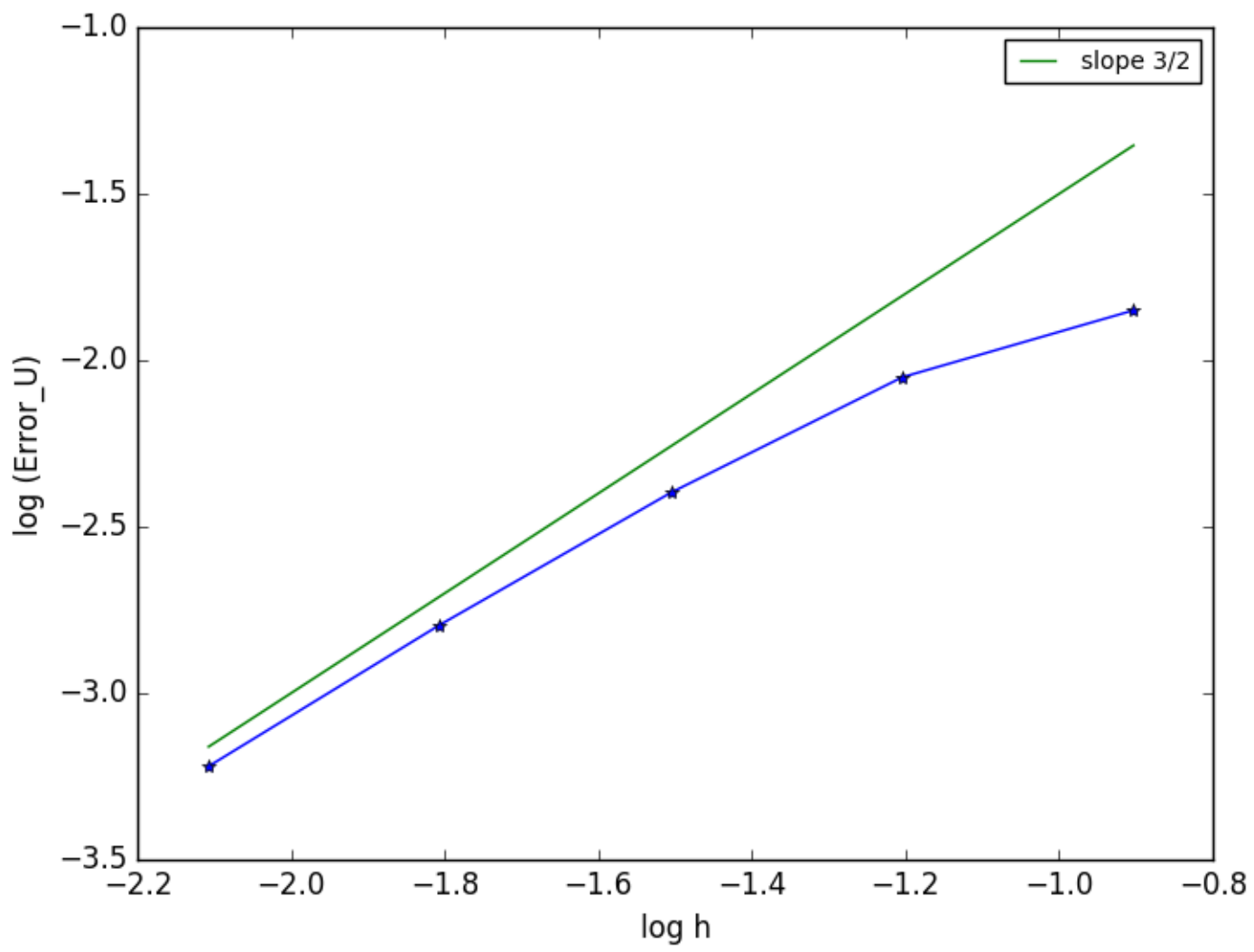

Table 4.7: Fluid velocity error in the case $r = 0$, evaluated at time $t = 0.0001s$.

| $h_f = h_s$ | $\\|\boldsymbol{d}_h^n - \boldsymbol{d}_{ref}^n\\|_s$ | Rate |
|---|---|---|
| 1/8 | 4.30822E-02 | |
| 1/16 | 2.24909E-02 | 9.37752E-01 |
| 1/32 | 1.13877E-02 | 9.81865E-01 |
| 1/64 | 5.95546E-03 | 9.35190E-01 |
| 1/128 | 3.34373E-03 | 8.32755E-01 |

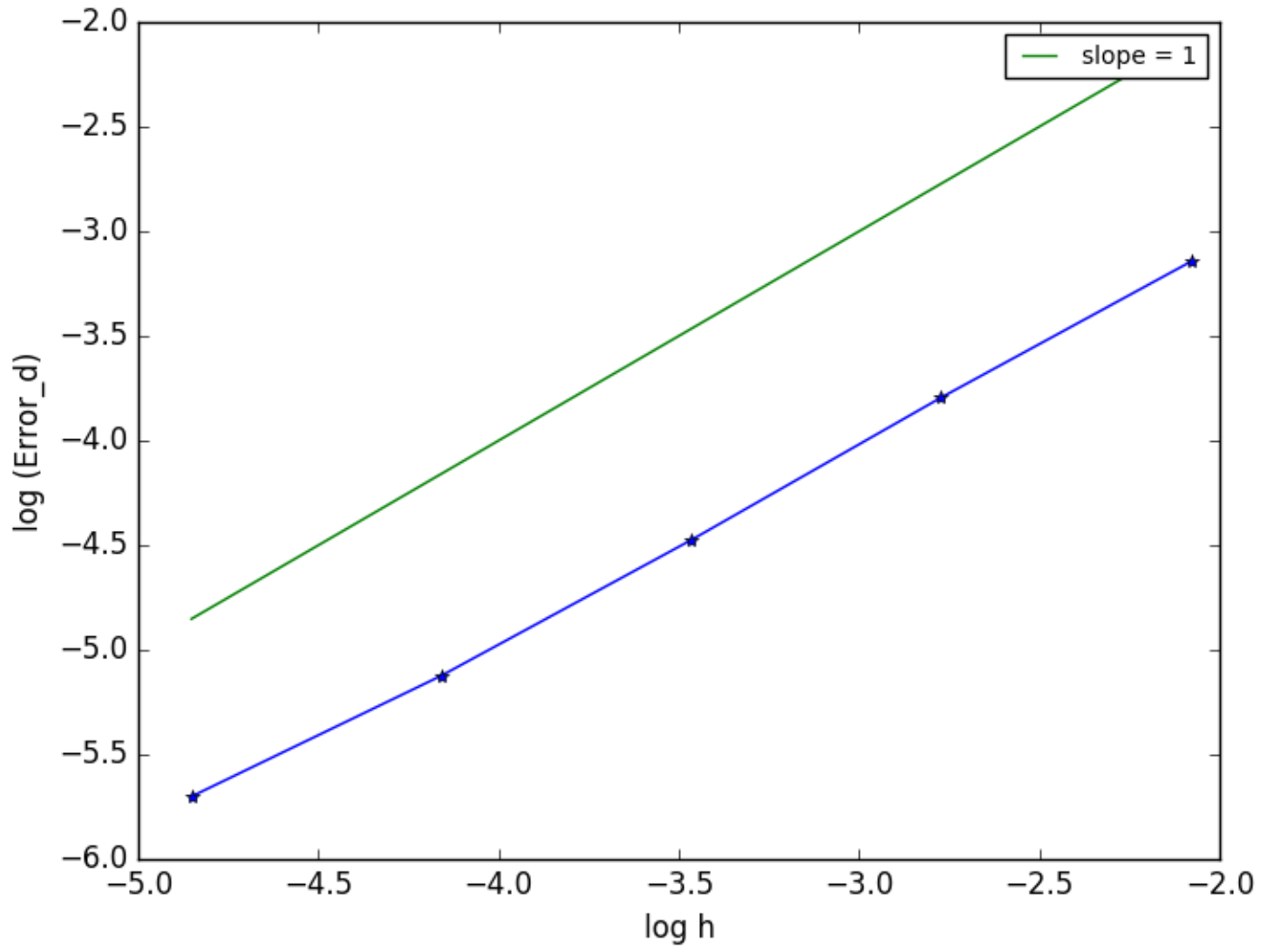

Table 4.8: Solid displacement error in the case $r = 0$, evaluated at time $t = 0.0001s$.

| $h_f = h_s$ | $\|\dot{\boldsymbol{d}}_h^n - \dot{\boldsymbol{d}}_{ref}^n\|_{0,\Sigma}$ | Rate |
|---|---|---|
| 1/8 | 5.78817E-02 | |
| 1/16 | 4.02564E-02 | 5.23887E-01 |
| 1/32 | 2.30805E-02 | 8.02543E-01 |
| 1/64 | 1.22144E-02 | 9.18090E-01 |
| 1/128 | 6.28019E-03 | 9.59706E-01 |

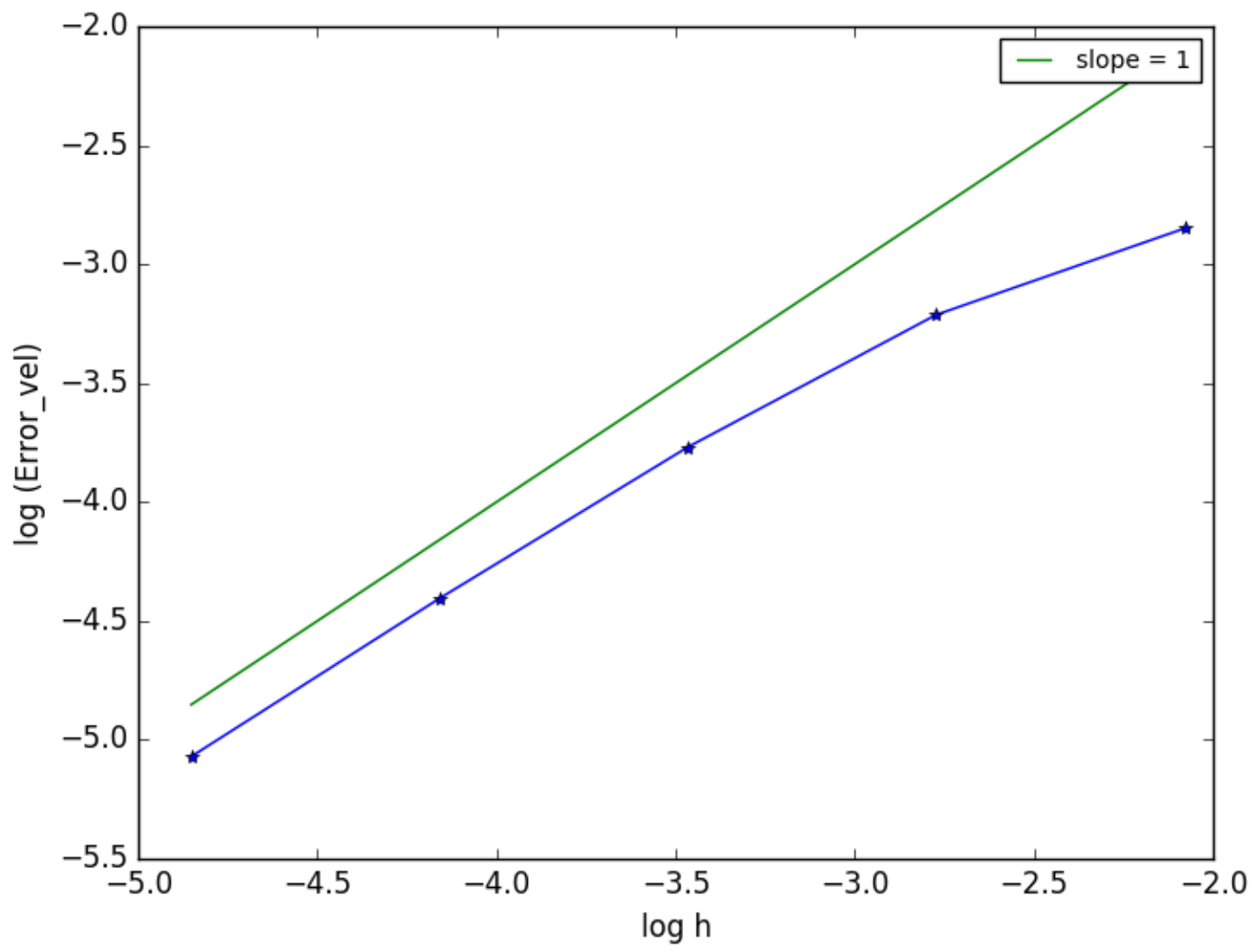

Table 4.9: Solid velocity error in the case $r = 0$, evaluated at time $t = 0.0001s$.

| $h_f = h_s$ | Total Error | Rate |
|---|---|---|
| 1/8 | 7.35251E-02 | |
| 1/16 | 4.69634E-02 | 6.46702E-01 |
| 1/32 | 2.60468E-02 | 8.50429E-01 |
| 1/64 | 1.36833E-02 | 9.28687E-01 |
| 1/128 | 7.14039E-03 | 9.38345E-01 |

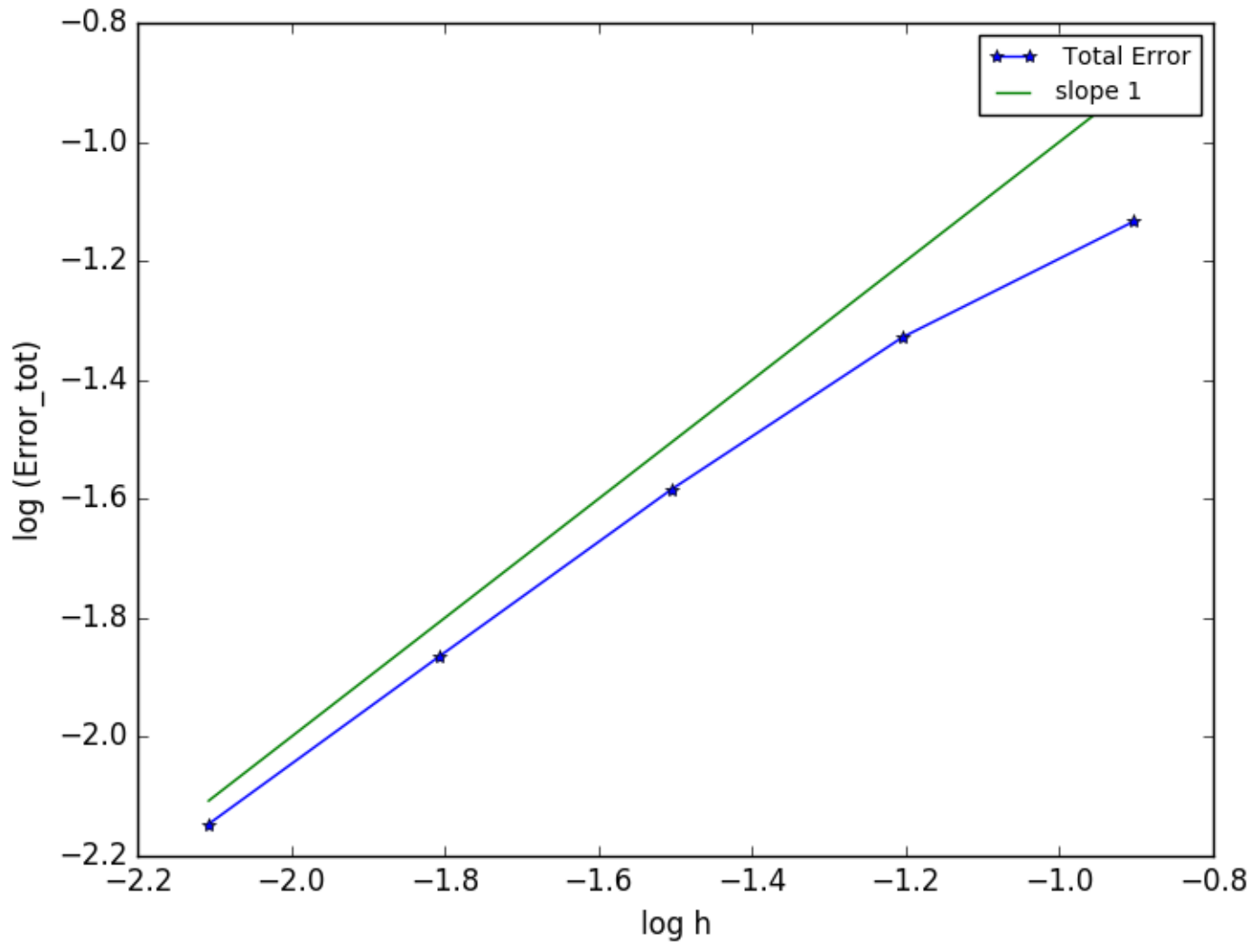

Table 4.10: Total error in the case $r = 0$, evaluated at time $t = 0.0001s$.

- Case $r = 1$

  In this case we have $\boldsymbol{d}^{n\star} = \boldsymbol{d}^{n-1}$ and $\boldsymbol{d}_h^{n\star} = \boldsymbol{d}_h^{n-1}$. As in the previous case, the fluid and velocity and pressure exhibit an optimal convergence rate according to the regularity assumptions. Also the accuracy od the approximation of the solid displacements results to be optimal in the energy norm or equivalently in $H^1$-norm. As in the previous case, the solid velocity exhibit a suboptimal behavior.

| $h_f = h_s$ | $\|\boldsymbol{u}_h^n - \boldsymbol{u}_{ref}^n\|_{0,\Omega}$ | Rate |
|---|---|---|
| 1/8 | 3.20977E-02 | |
| 1/16 | 1.50066E-02 | 1.09687E+00 |
| 1/32 | 5.63142E-03 | 1.41403E+00 |
| 1/64 | 2.20714E-03 | 1.35132E+00 |
| 1/128 | 7.69605E-04 | 1.51999E+00 |

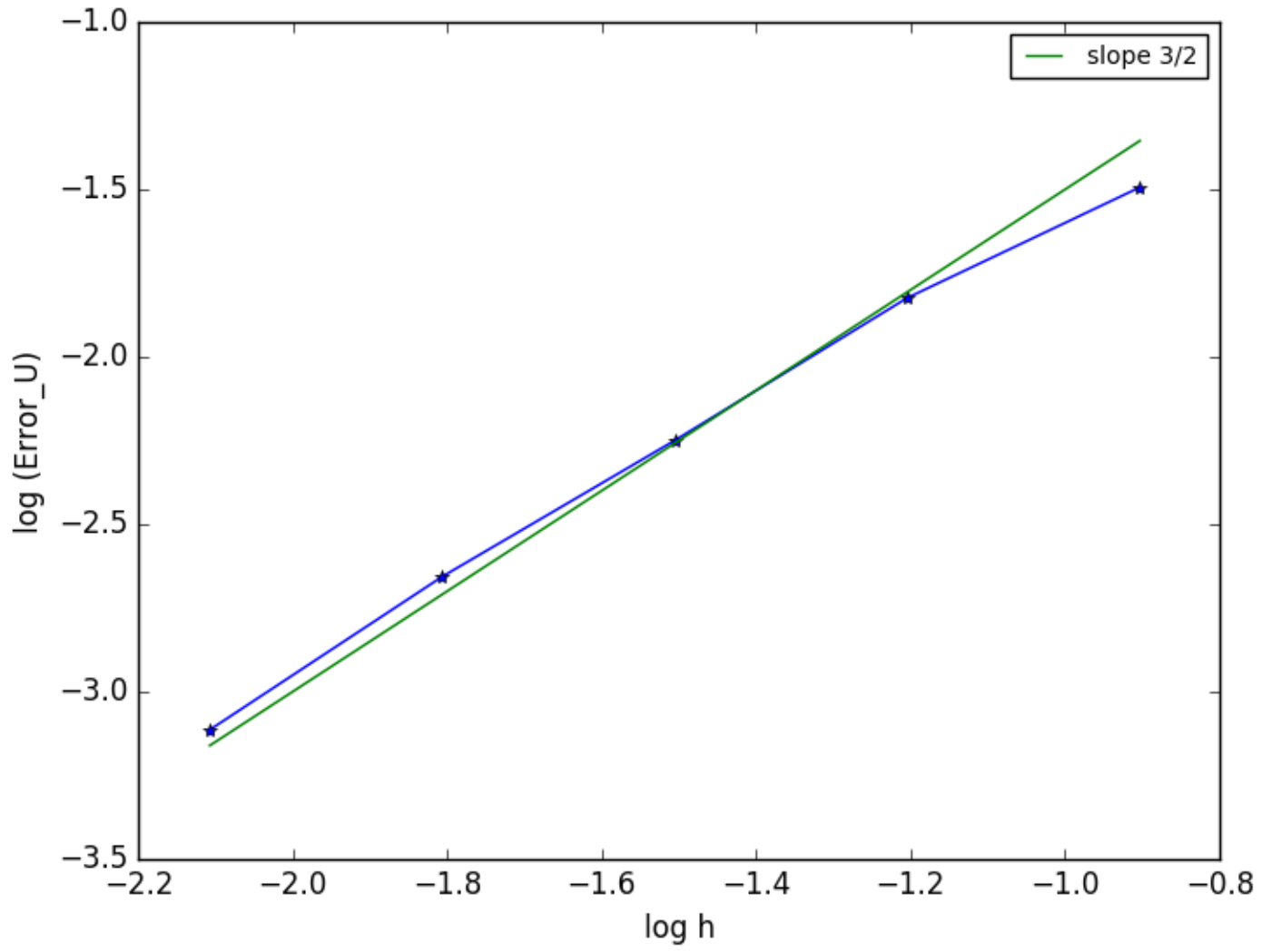

Table 4.11: Fluid velocity error in the case $r = 1$, evaluated at time $t = 0.0001s$

| $h_f = h_s$ | $\|\boldsymbol{d}_h^n - \boldsymbol{d}_{ref}^n\|_s$ | Rate |
|---|---|---|
| 1/8 | 3.10446E-02 | |
| 1/16 | 1.60457E-02 | 9.52157E-01 |
| 1/32 | 8.38743E-03 | 9.35883E-01 |
| 1/64 | 4.69113E-03 | 8.38293E-01 |
| 1/128 | 2.83775E-03 | 7.25185E-01 |

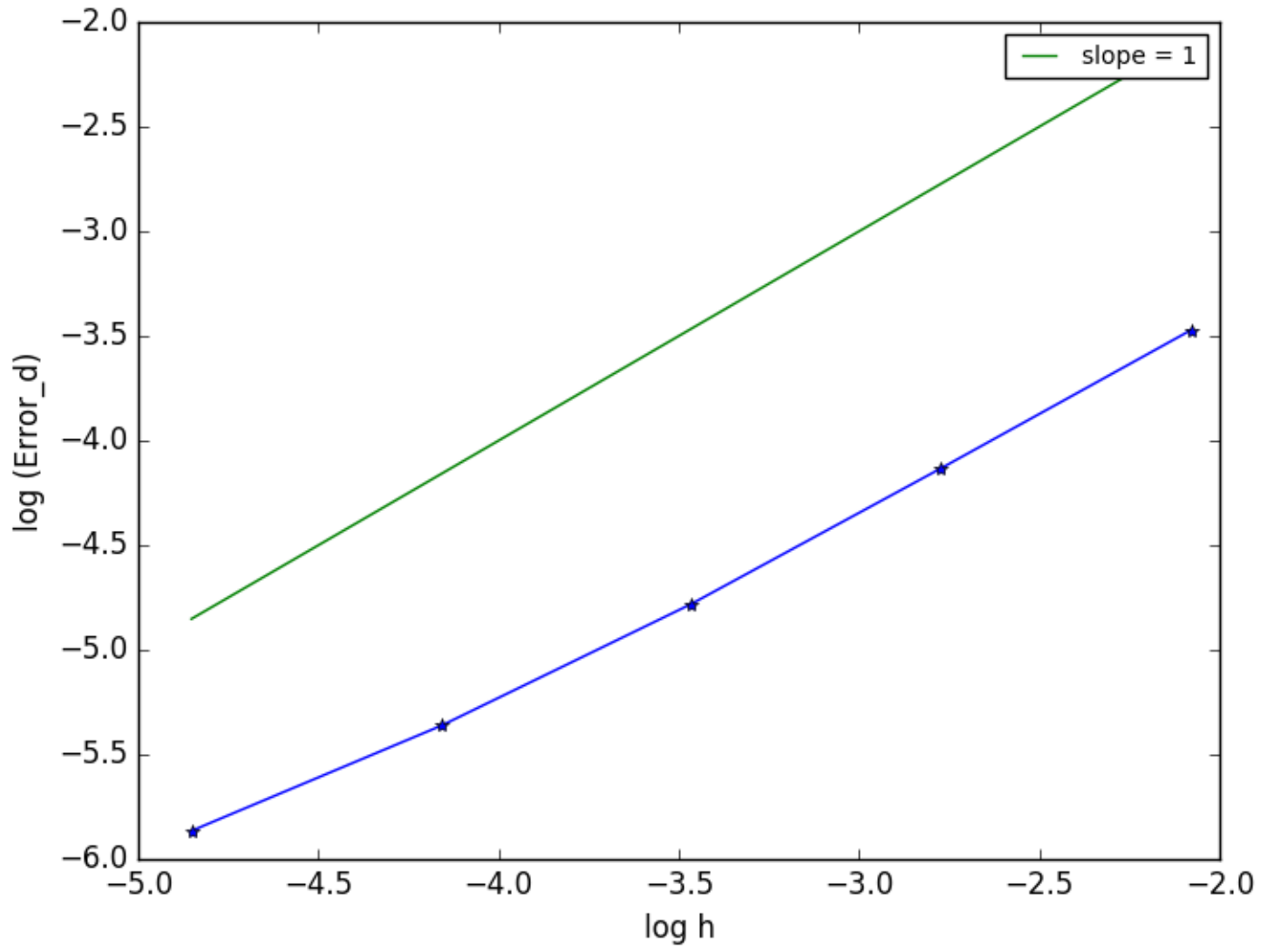

Table 4.12: Solid displacement error in the case $r = 1$, evaluated at time $t = 0.0001s$.

| $h_f = h_s$ | $\|\dot{\boldsymbol{d}}_h^n - \dot{\boldsymbol{d}}_{ref}^n\|_{0,\Sigma}$ | Rate |
|---|---|---|
| 1/8 | 9.84642E-04 | |
| 1/16 | 9.15567E-04 | 1.04934E-01 |
| 1/32 | 5.62710E-04 | 7.02273E-01 |
| 1/64 | 2.73571E-04 | 1.04048E+00 |
| 1/128 | 1.43604E-04 | 9.29822E-01 |

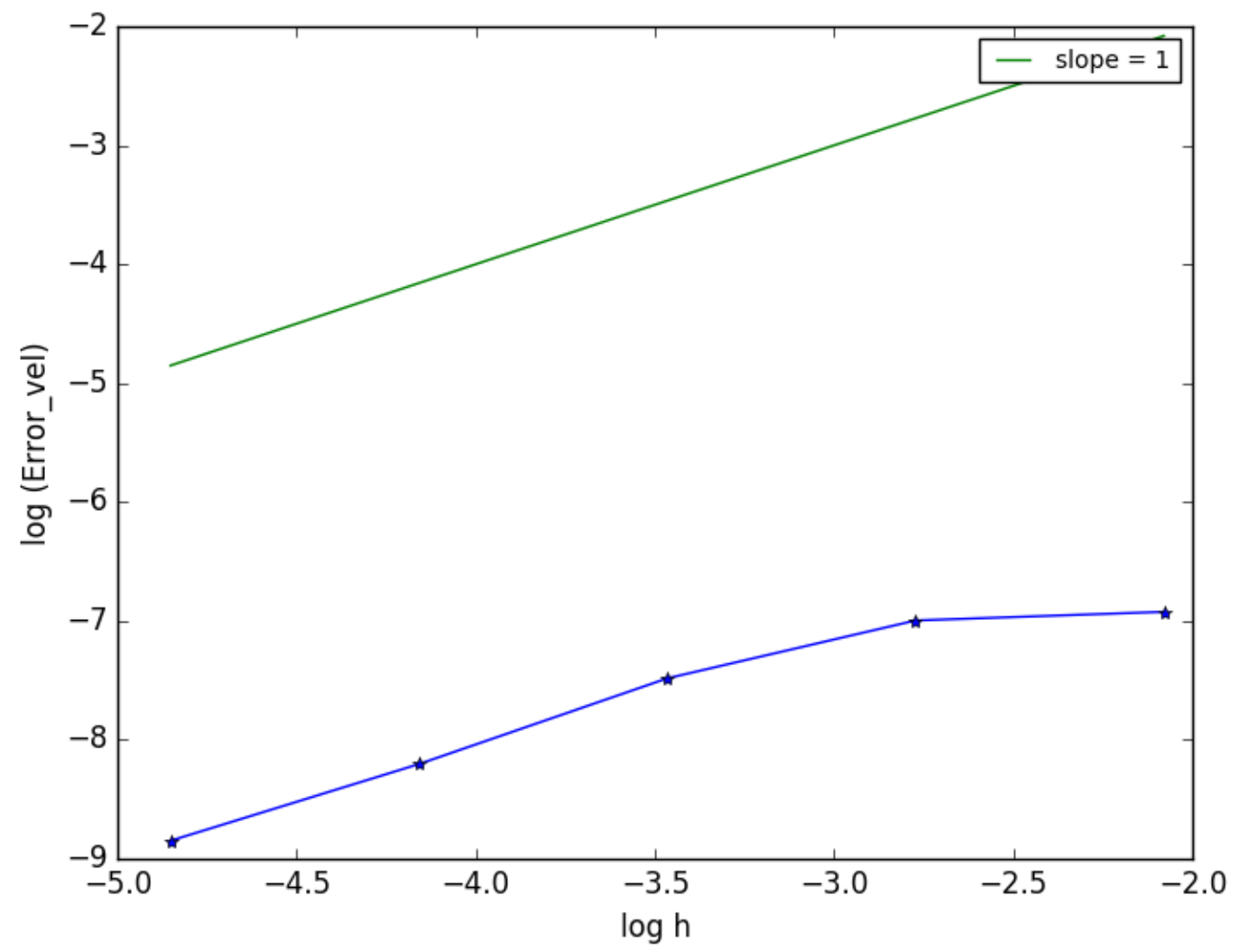

Table 4.13: Solid velocity error in the case $r = 1$, evaluated at time $t = 0.0001s$.

| $h_f = h_s$ | Total Error | Rate |
|---|---|---|
| 1/8 | 4.46654E-02 | |
| 1/16 | 2.19886E-02 | 1.02240E+00 |
| 1/32 | 1.01182E-02 | 1.11980E+00 |
| 1/64 | 5.19163E-03 | 9.62698E-01 |
| 1/128 | 2.94377E-03 | 8.18523E-01 |

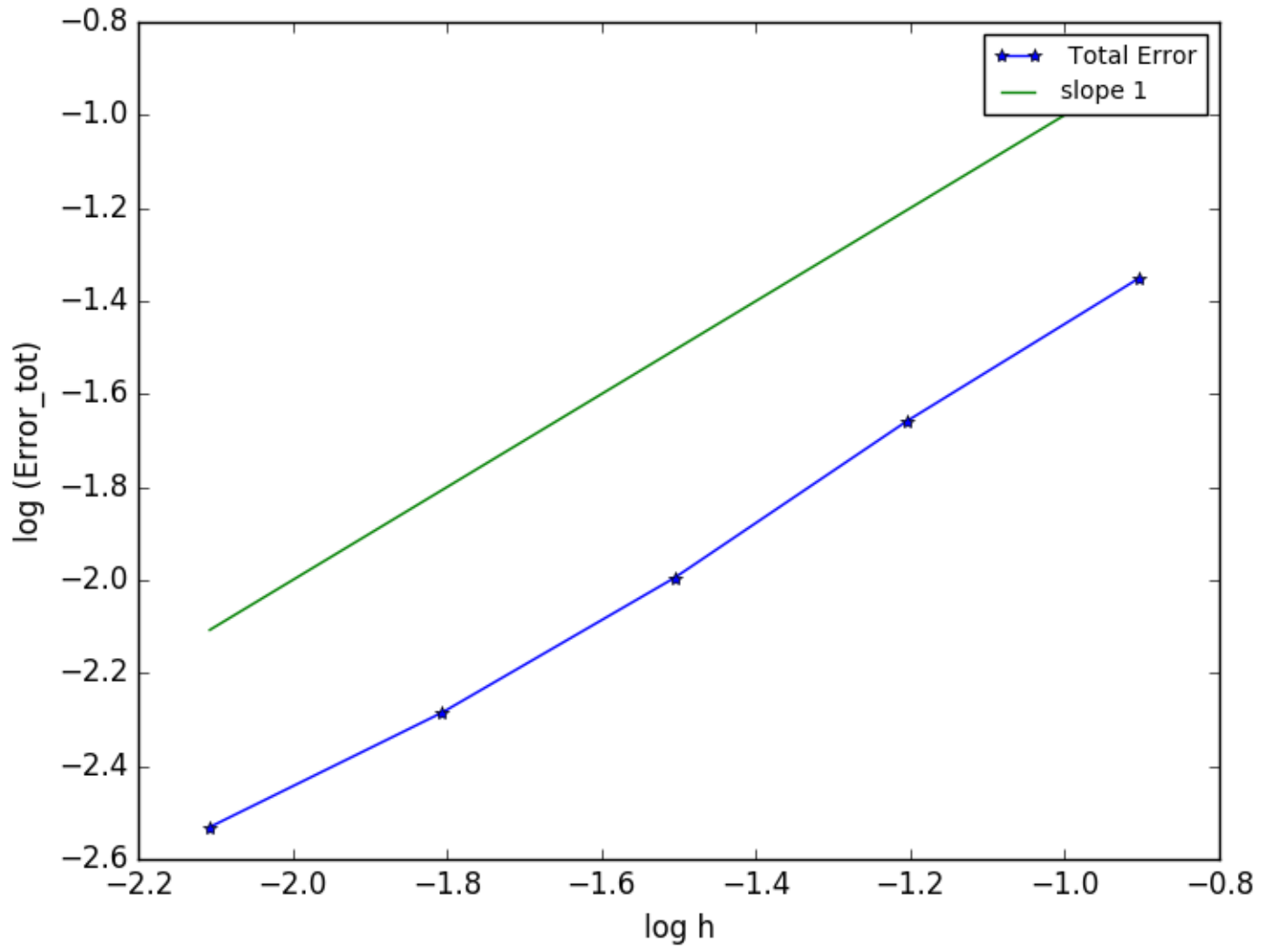

Table 4.14: Total error in the case $r = 1$, evaluated at time $t = 0.0001s$.

- Case $r = 2$

  In this case we have $\boldsymbol{d}^{n\star} = \boldsymbol{d}^{n-1} + \tau\dot{\boldsymbol{d}}^{n-1}$ and $\boldsymbol{d}_h^{n\star} = \boldsymbol{d}_h^{n-1} + \tau\dot{\boldsymbol{d}}_h^{n-1}$. The accuracy of the split scheme with extrapolation $r = 2$ is evaluated, as in the previous cases with reference to the steady test case of equilibrium of a circular elastic string immersed in the fluid solved using the monolithic algorithm with discretization parameters 3.62.

| $h_f = h_s$ | $\Vert \boldsymbol{u}_h^n - \boldsymbol{u}_{ref}^n \Vert_{0,\Omega}$ | Rate |
|---|---|---|
| 1/8 | 2.98843E-02 | |
| 1/16 | 1.52670E-02 | 9.68969E-01 |
| 1/32 | 6.22870E-03 | 1.29341E+00 |
| 1/64 | 2.28059E-03 | 1.44953E+00 |
| 1/128 | 8.00233E-04 | 1.51091E+00 |

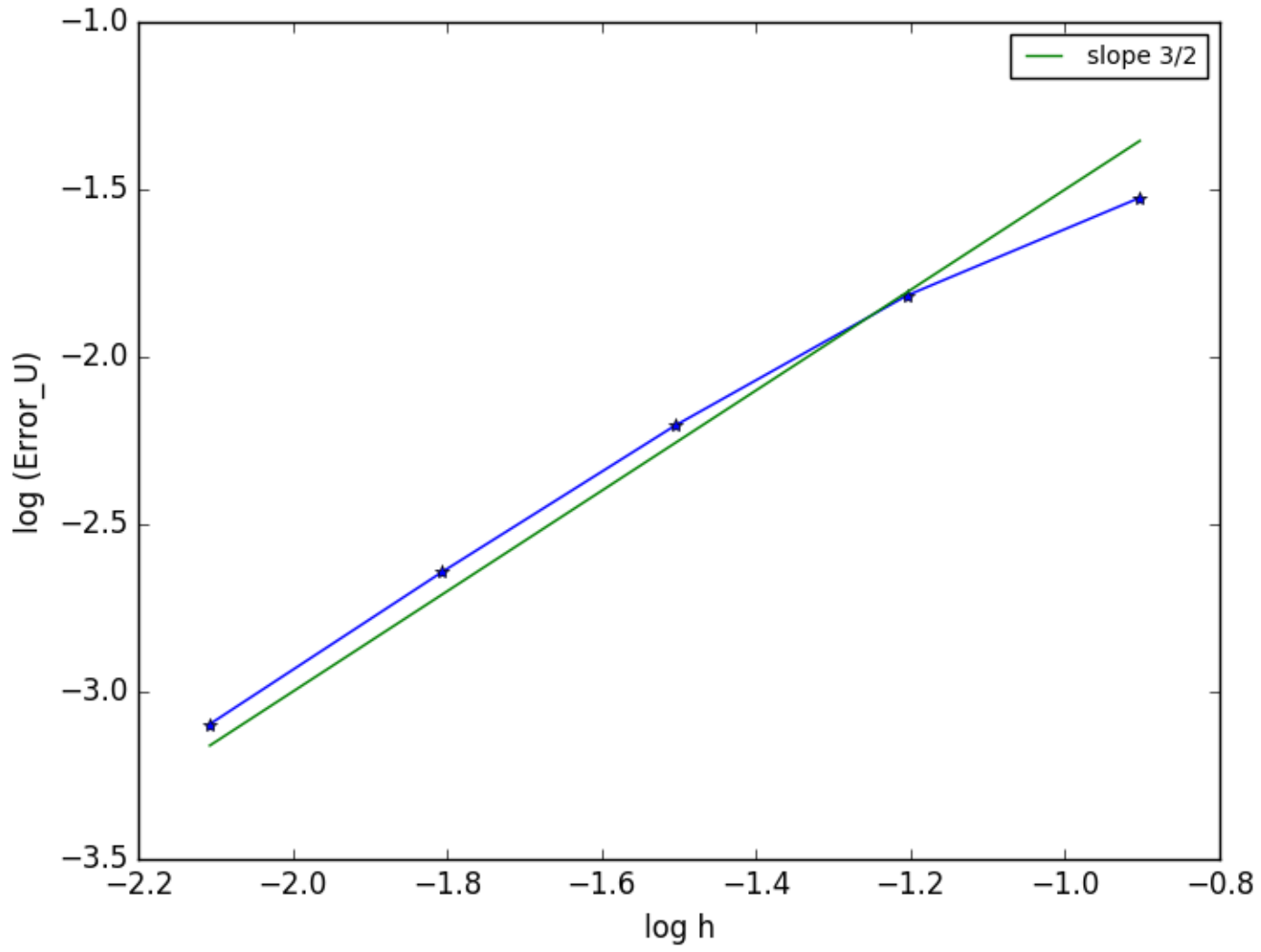

Table 4.15: Fluid velocity error in the case $r = 2$, evaluated at time $t = 0.0001s$.

| $h_f = h_s$ | $\Vert \boldsymbol{d}_h^n - \boldsymbol{d}_{ref}^n \Vert_s$ | Rate |
|---|---|---|
| 1/8 | 3.09276E-02 | |
| 1/16 | 1.58438E-02 | 9.64977E-01 |
| 1/32 | 8.19538E-03 | 9.51035E-01 |
| 1/64 | 4.70732E-03 | 7.99905E-01 |
| 1/128 | 2.83916E-03 | 7.29442E-01 |

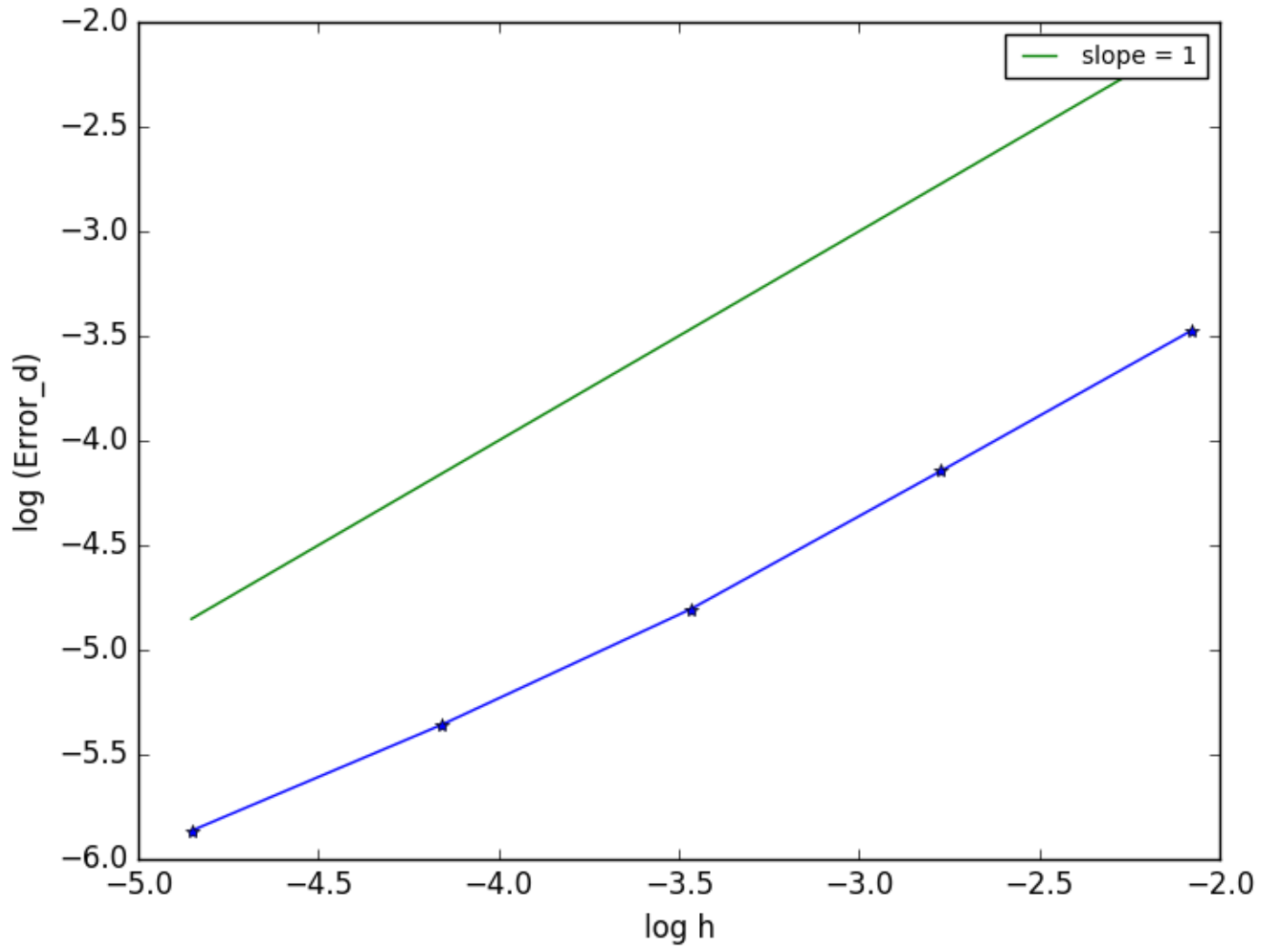

Table 4.16: Solid displacement error in the case $r = 2$, evaluated at time $t = 0.0001s$.

| $h_f = h_s$ | $\|\dot{\boldsymbol{d}}_h^n - \dot{\boldsymbol{d}}_{ref}^n\|_{0,\Sigma}$ | Rate |
|---|---|---|
| 1/8 | 1.47116E-03 | |
| 1/16 | 1.14310E-03 | 3.64010E-01 |
| 1/32 | 5.77172E-04 | 9.85875E-01 |
| 1/64 | 2.62340E-04 | 1.13756E+00 |
| 1/128 | 1.44076E-04 | 8.64605E-01 |

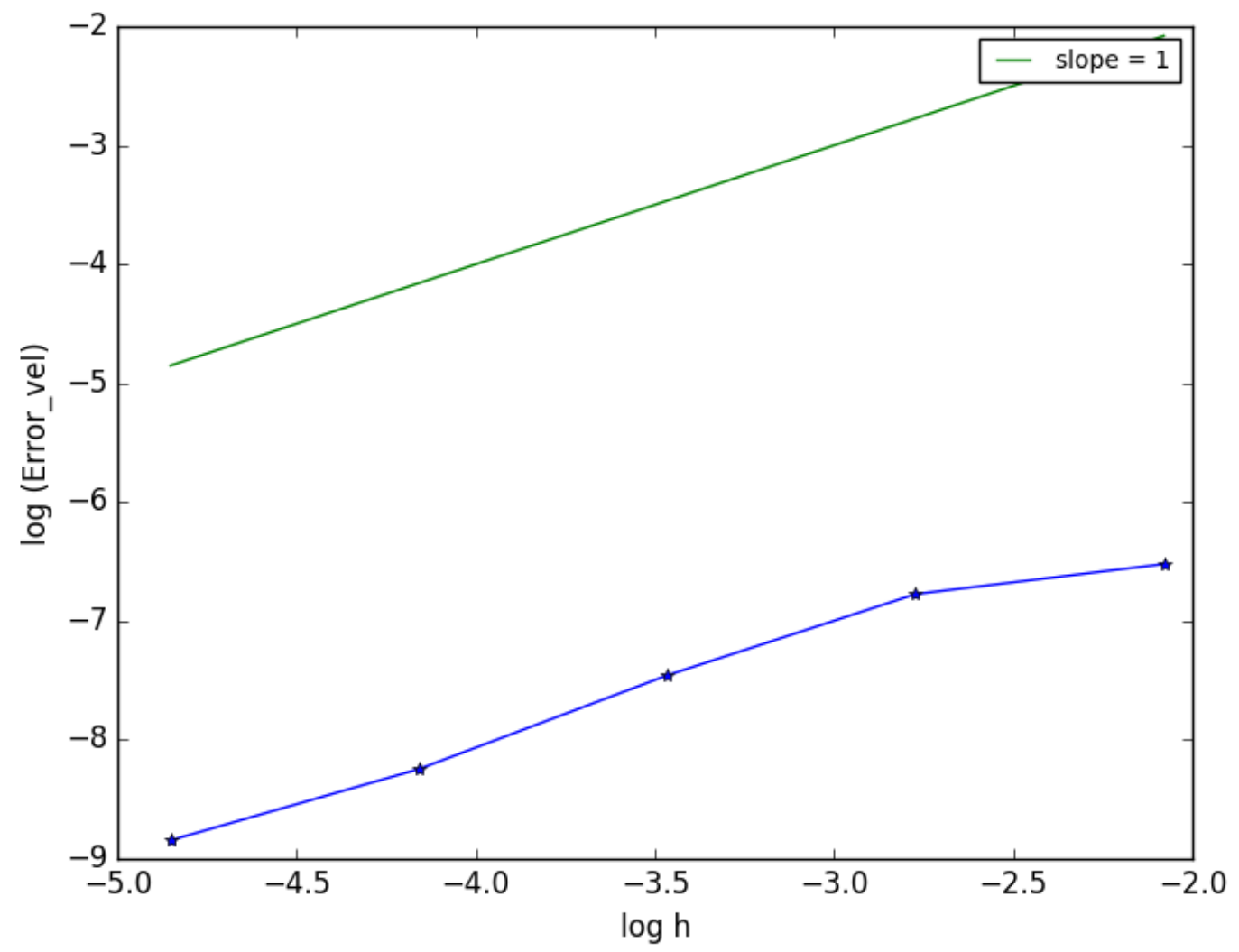

Table 4.17: Solid velocity error in the case $r = 2$, evaluated at time $t = 0.0001s$.

| $h_f = h_s$ | Total Error | Rate |
|---|---|---|
| 1/8 | 4.30320E-02 | |
| 1/16 | 2.20321E-02 | 9.65801E-01 |
| 1/32 | 1.03099E-02 | 1.09558E+00 |
| 1/64 | 5.23725E-03 | 9.77150E-01 |
| 1/128 | 2.95330E-03 | 8.26482E-01 |

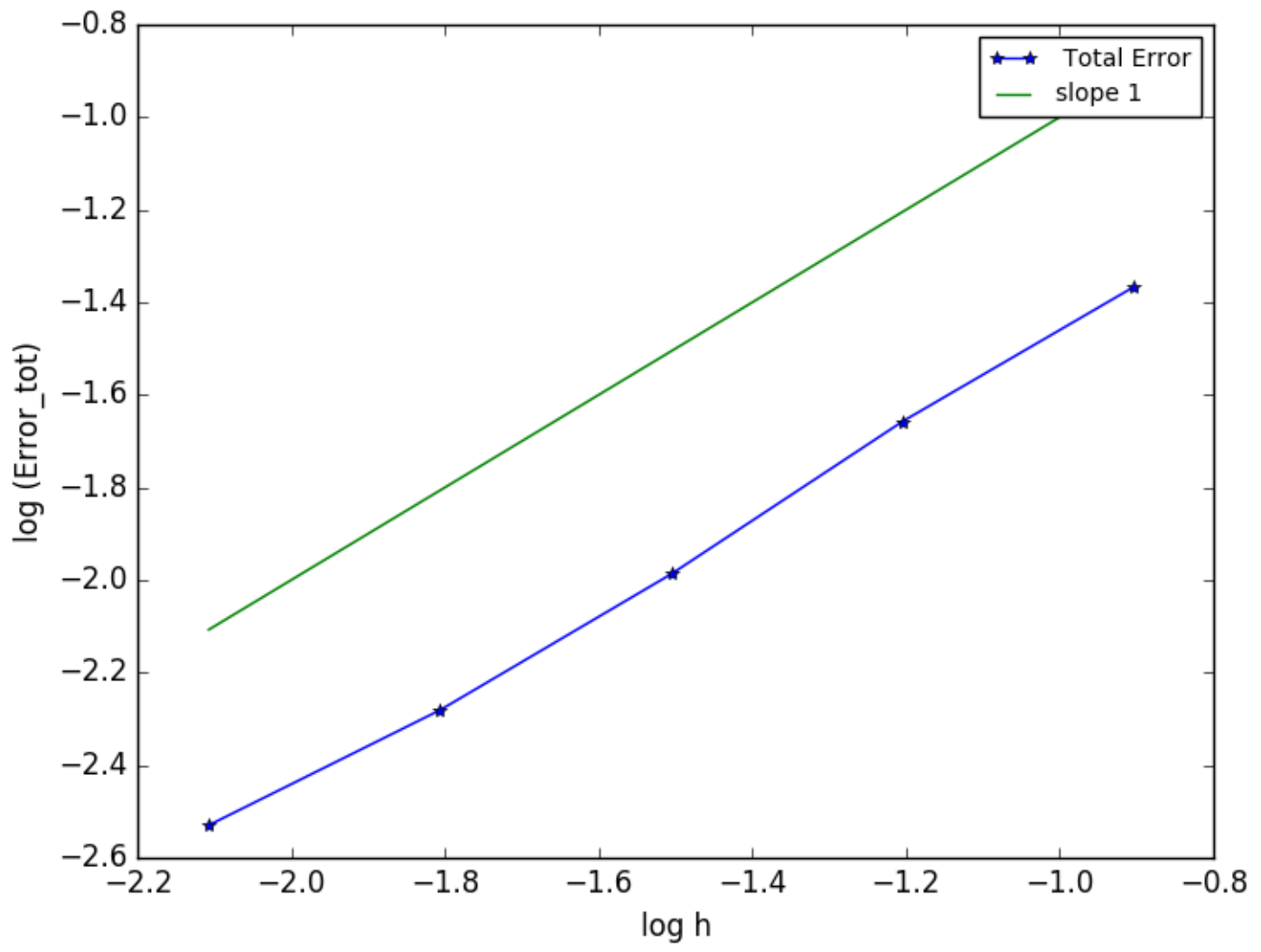

Table 4.18: Total error in the case $r = 2$, evaluated at time $t = 0.0001s$.

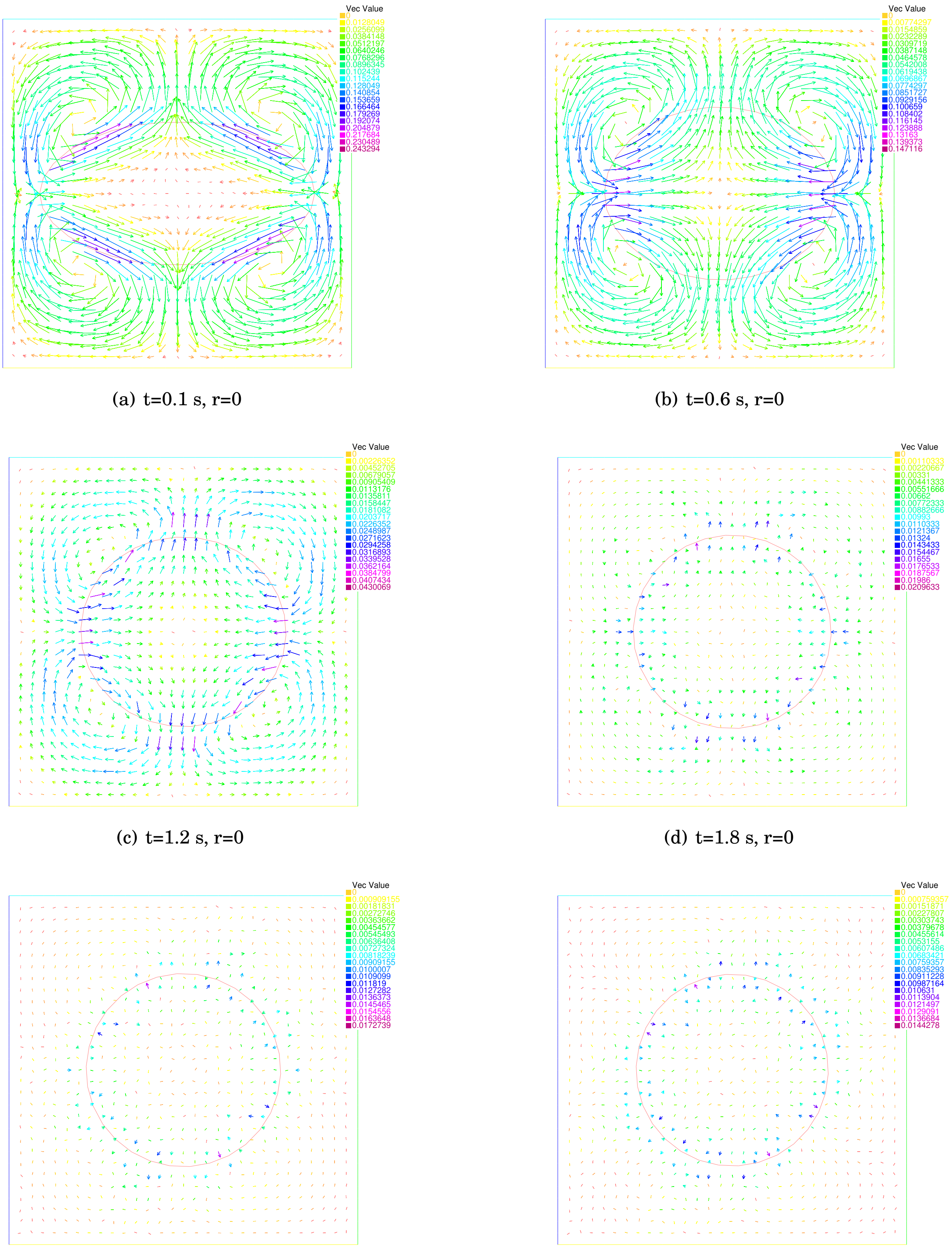

(a) t=0.1 s, r=0
(b) t=0.6 s, r=0
(c) t=1.2 s, r=0
(d) t=1.8 s, r=0
(e) t=2.4 s, r=0
(f) t=3.0 s, r=0

Figure 4.11: Algorithm 3, r=0

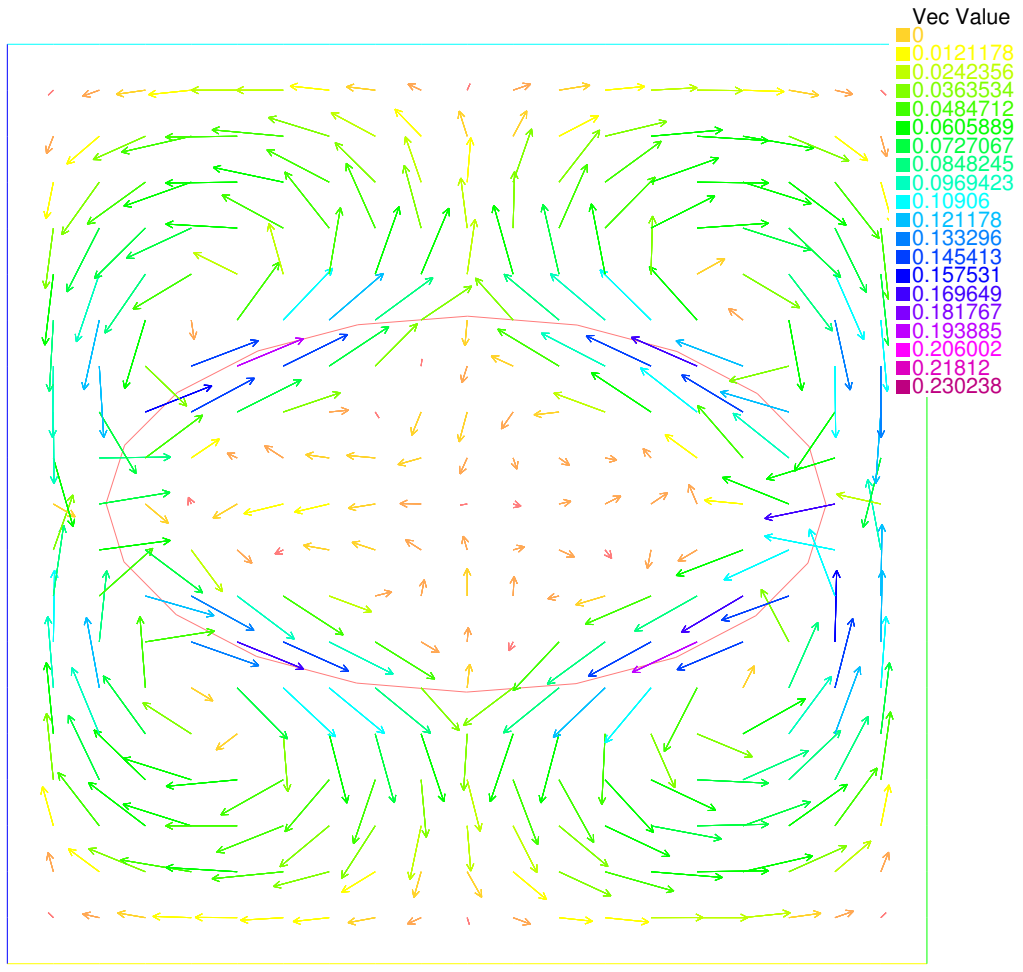

(a) t=0.1 s, r=1

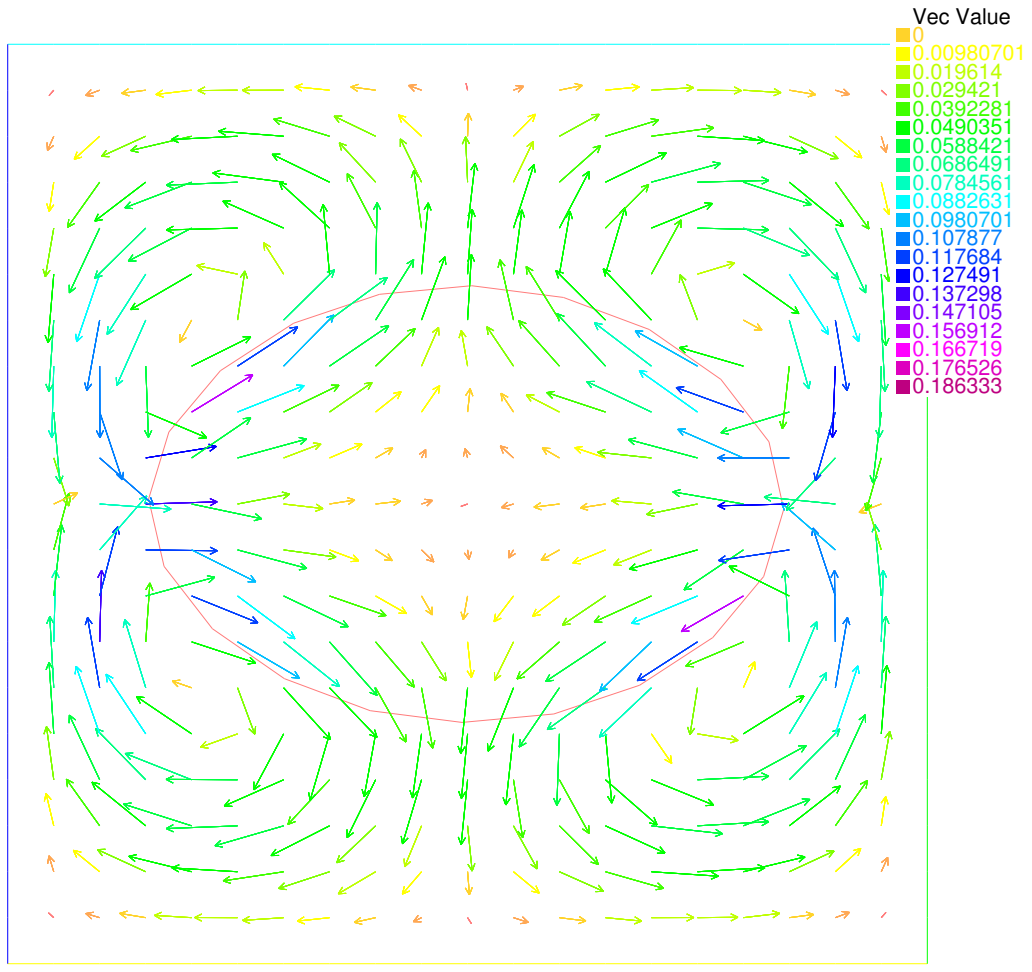

(b) t=0.6 s, r=1

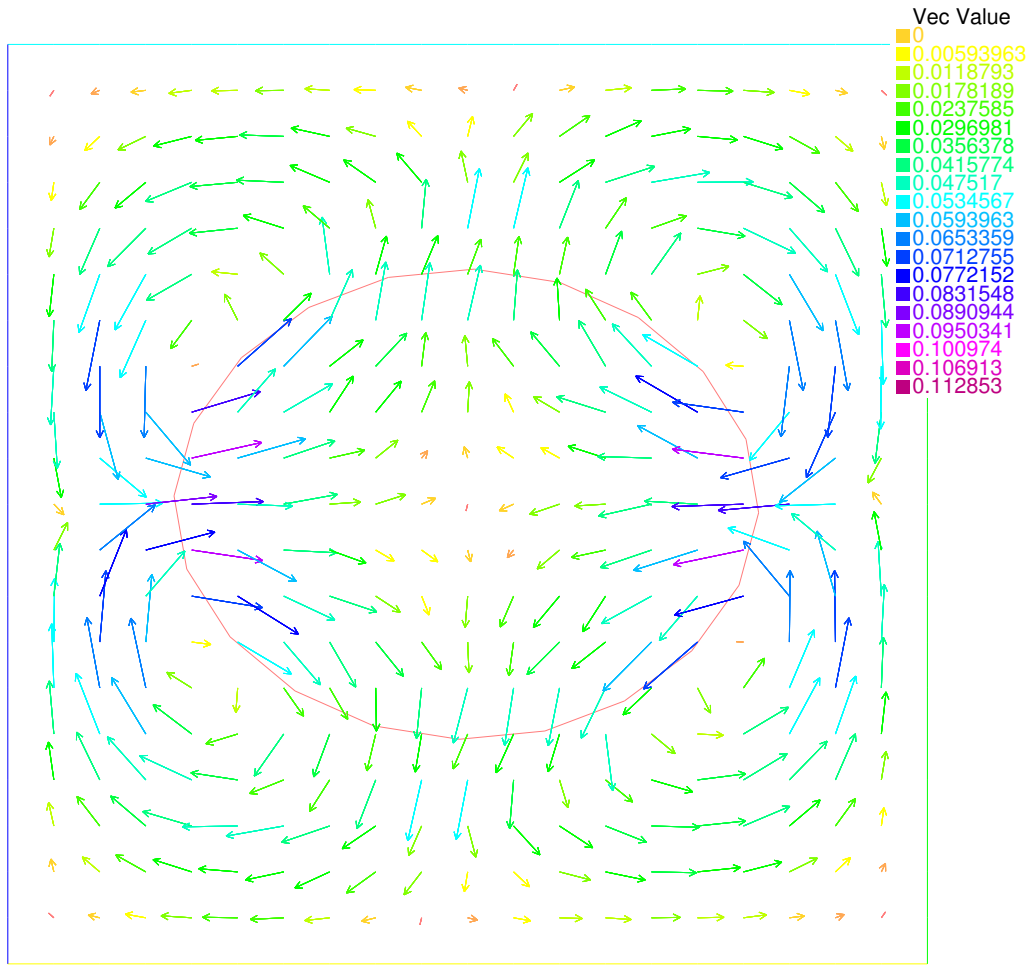

(c) t=1.2 s, r=1

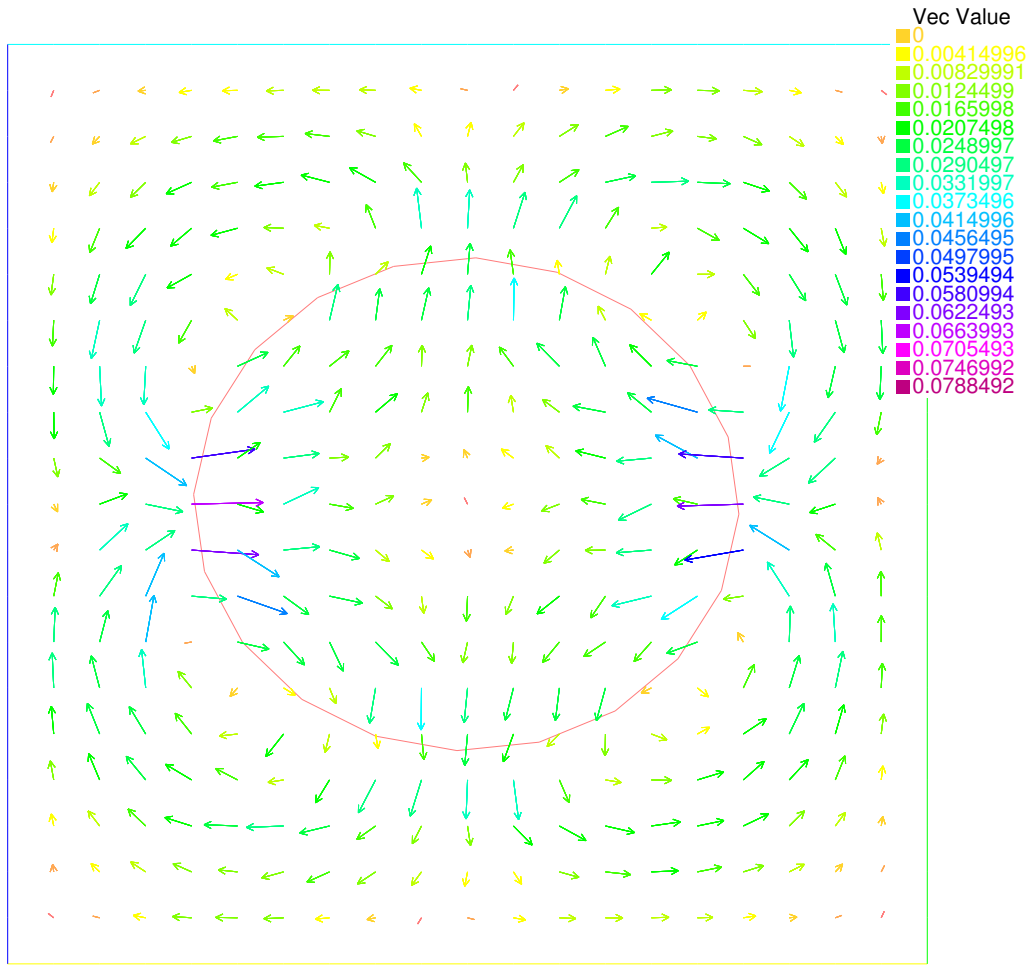

(d) t=1.8 s, r=1

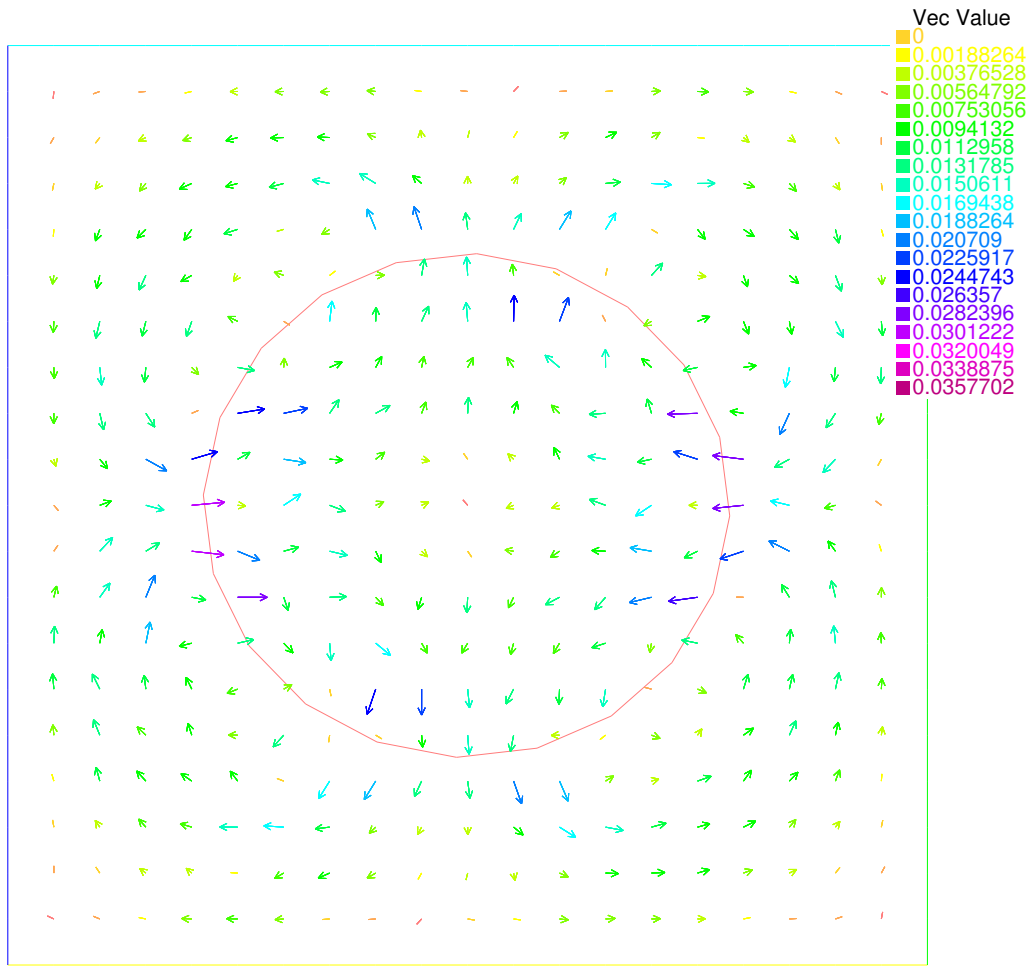

(e) t=2.4 s, r=1

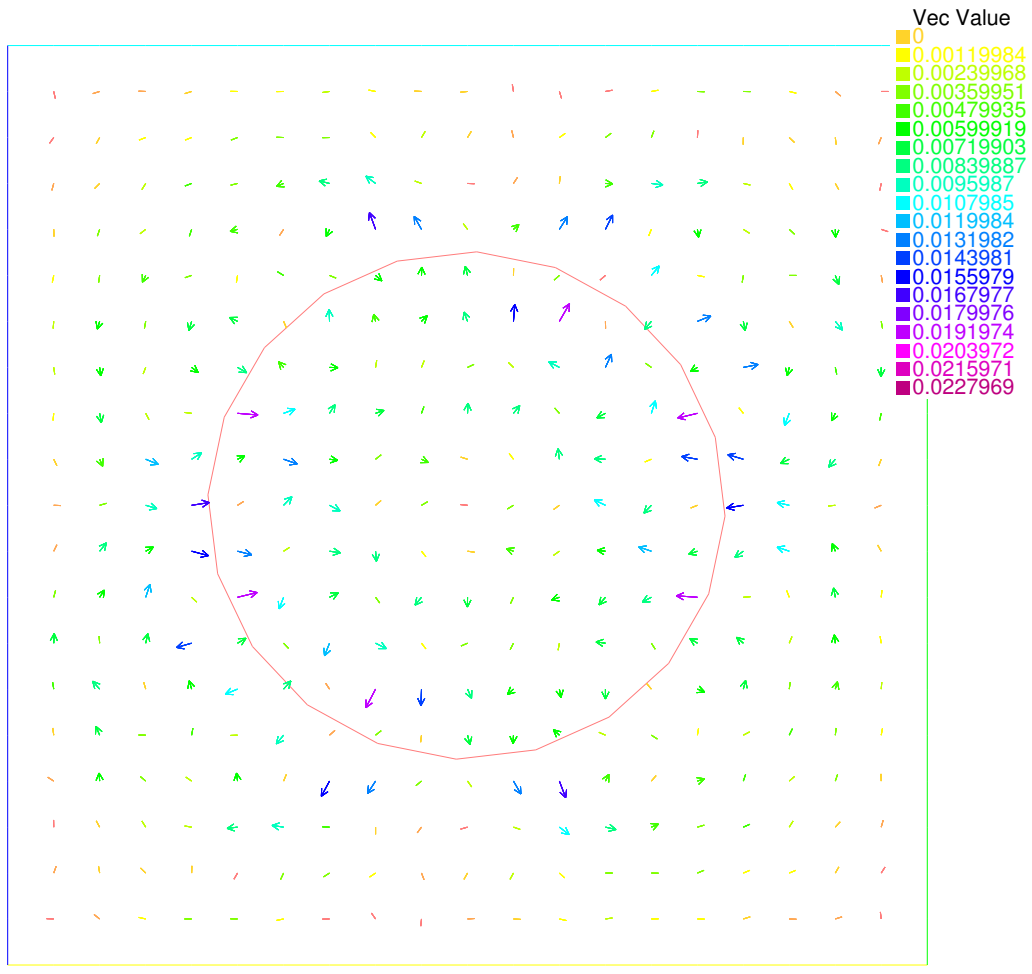

(f) t=3.0 s, r=1

Figure 4.12: Algorithm 3, r=1

ince

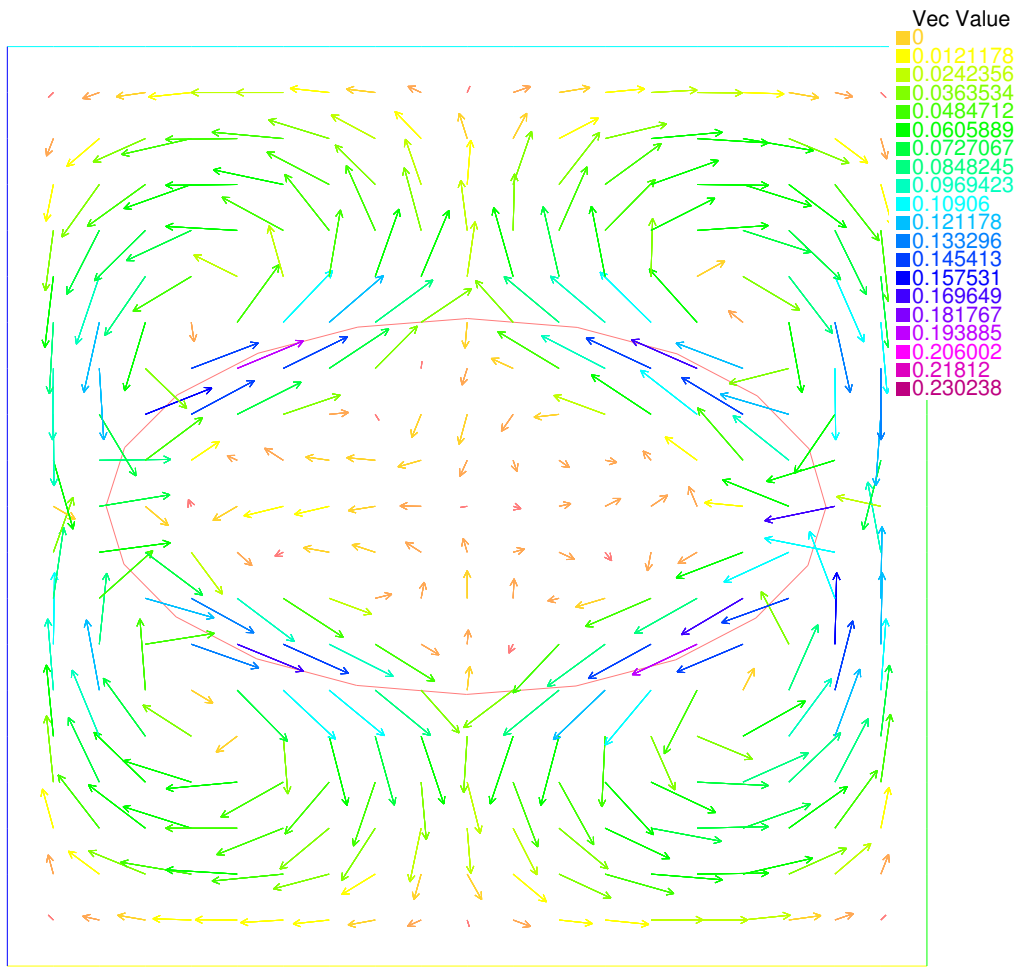

(a) t=0.1 s, r=2

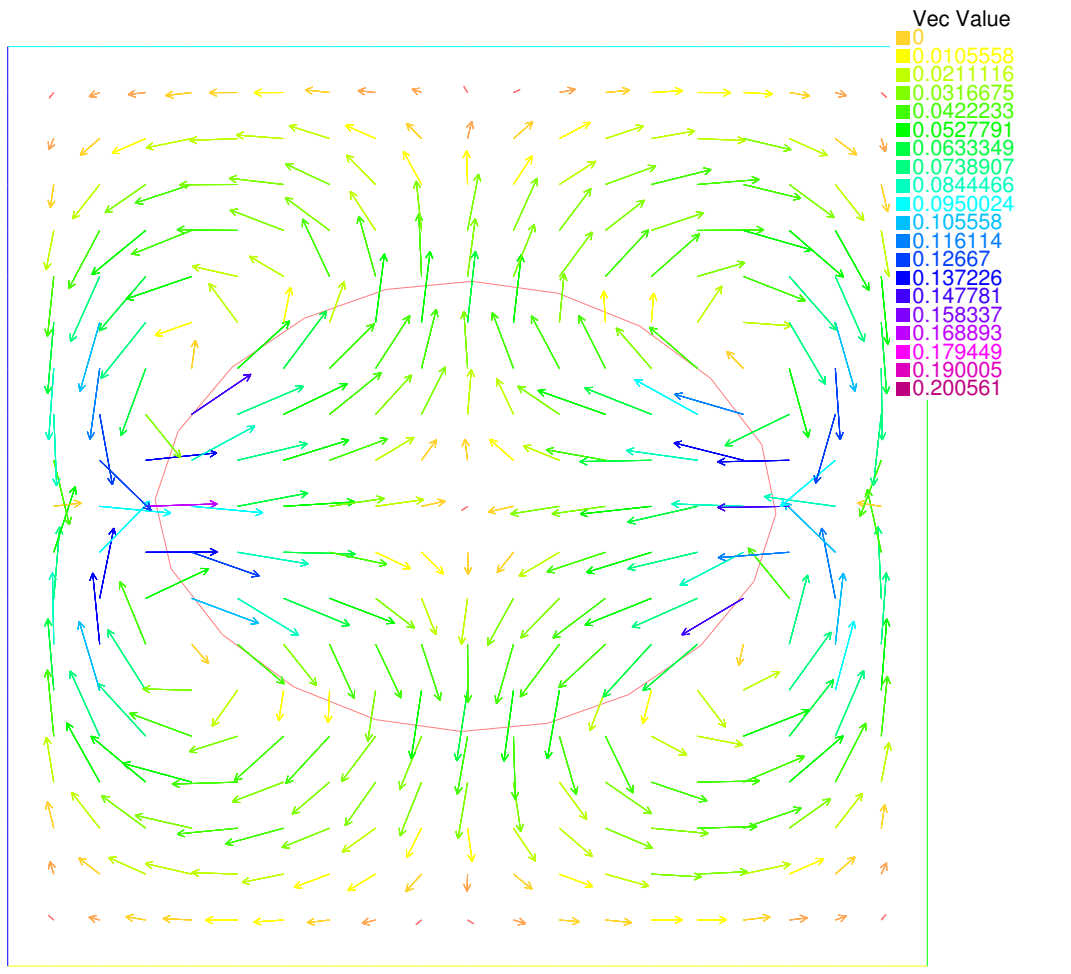

(b) t=0.6 s, r=2

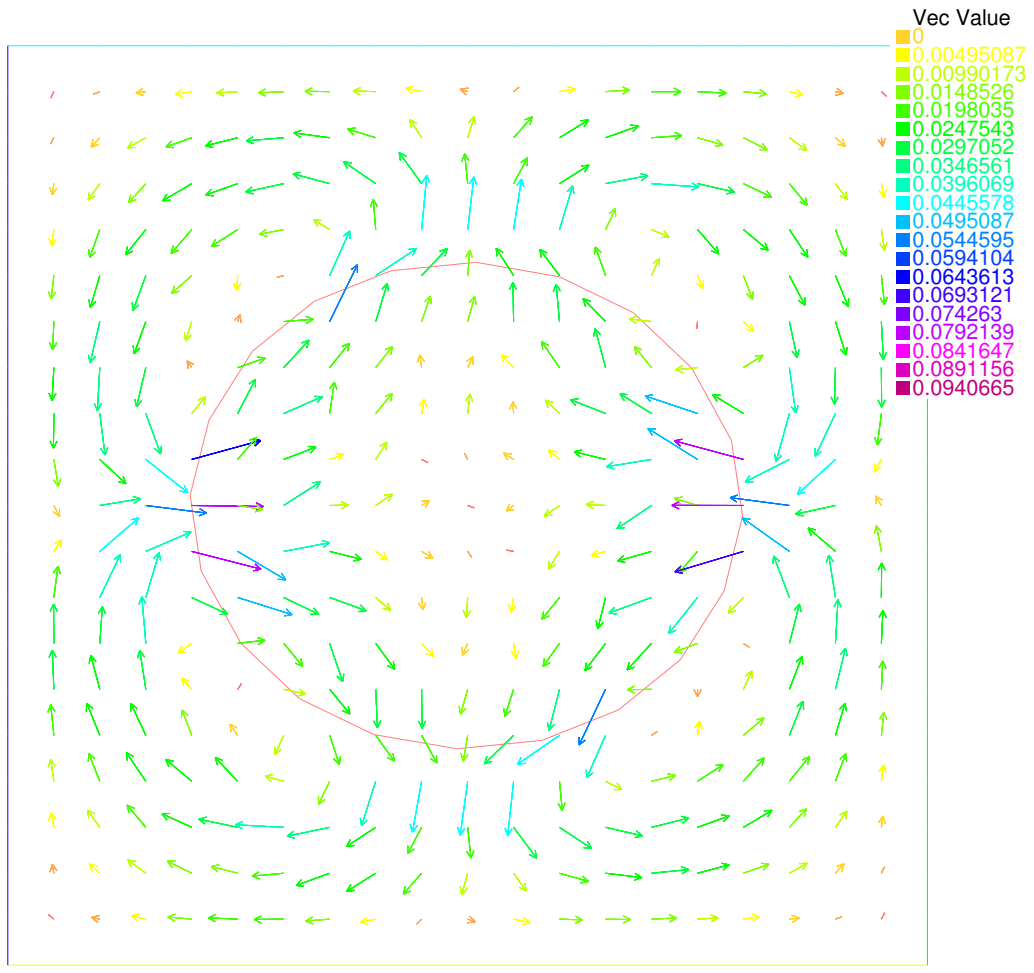

(c) t=1.2 s, r=2

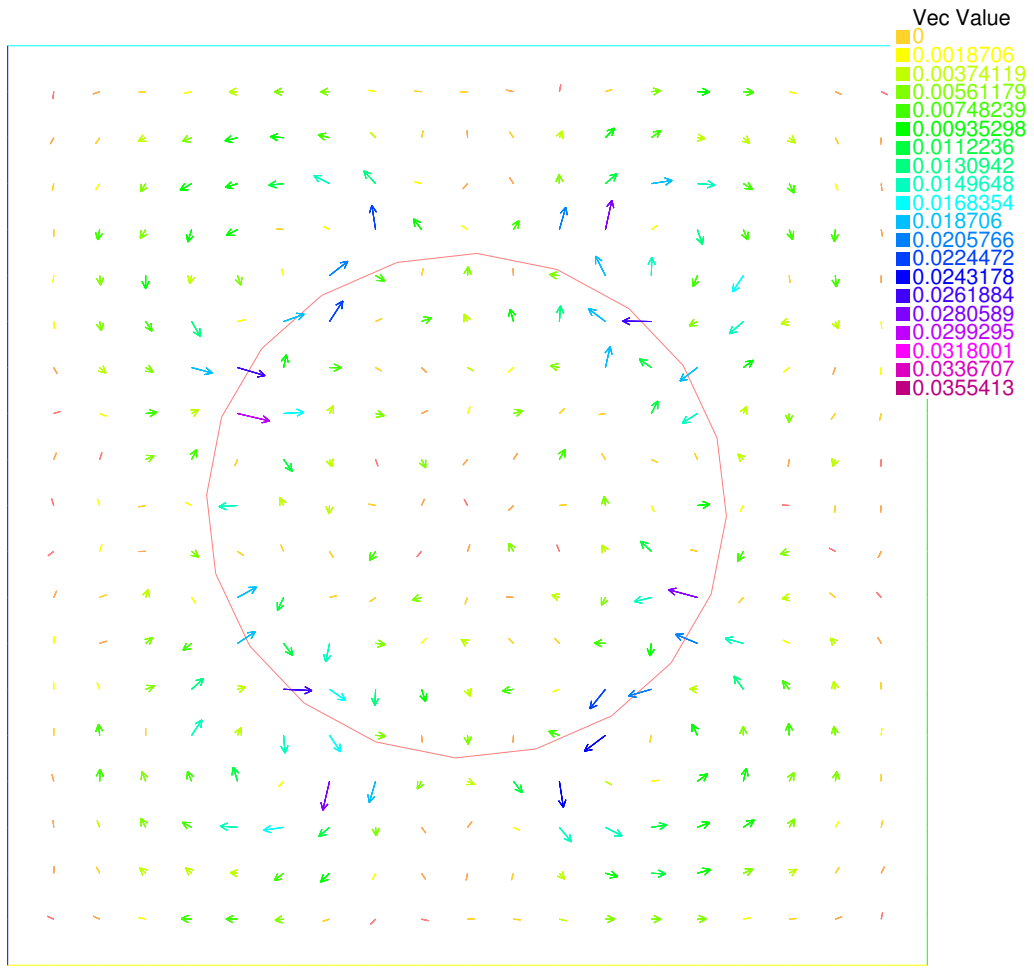

(d) t=1.8 s, r=2

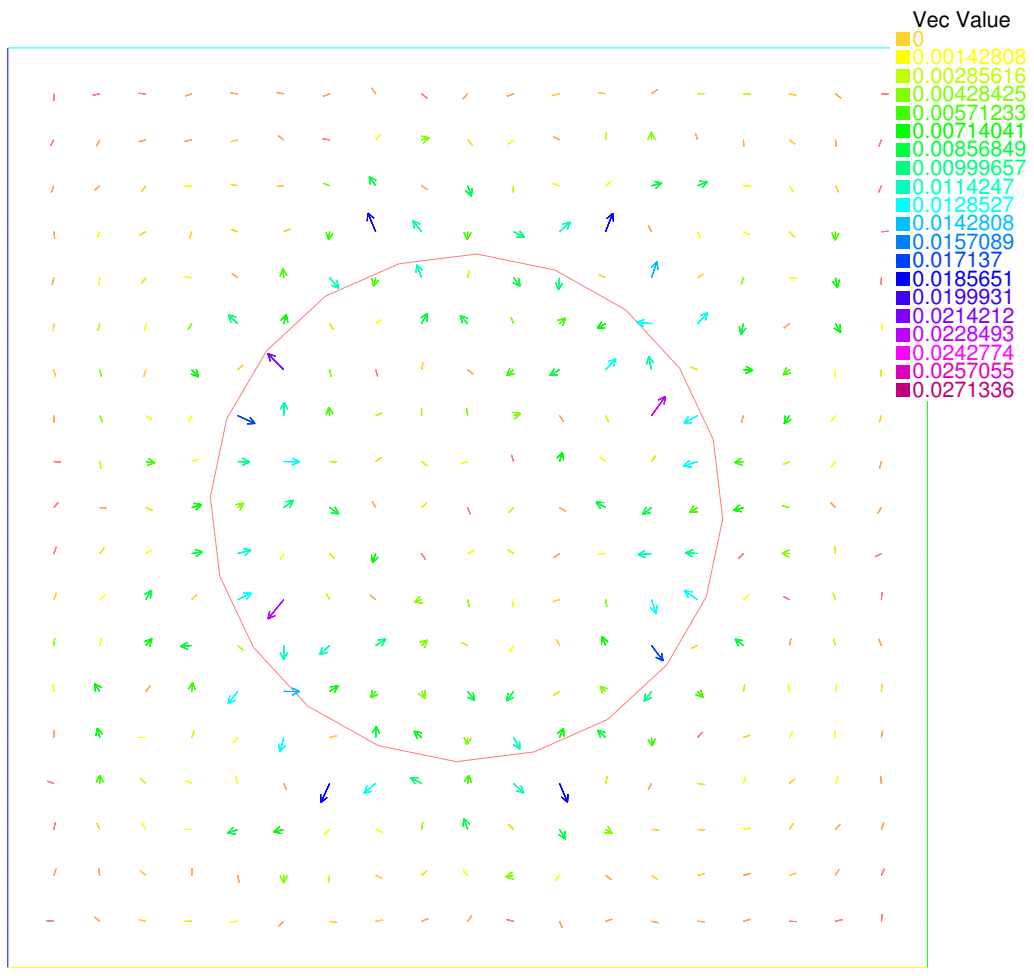

(e) t=2.4 s, r=2

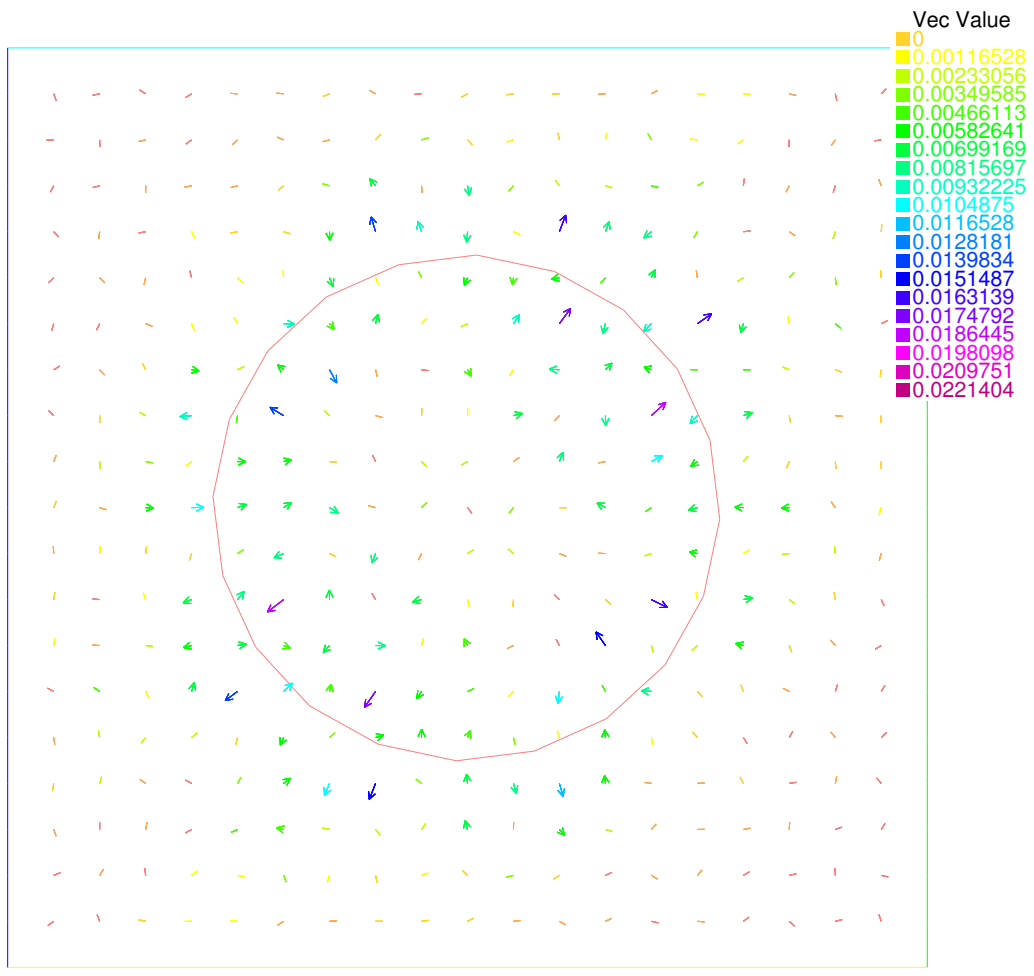

(f) t=3.0 s, r=2

Figure 4.13: Algorithm 3, r=2

# Conclusions

In conclusion, we have analyzed from the theoretical point of view the stability, the well-posedness and the convergence of a monolithic and two partitioned numerical schemes designed in order to solve the coupled problem of thin elastic solids immersed in incompressible fluids by using finite elements. Both the monolithic and partitioned paradigms use unfitted fluid and solid meshes.

The coupling of the fluid and solid is realized using the technique of Lagrange multipliers as introduced in [20, 23] where it was presented the monolithic scheme and it was shown its unconditional stablility. In [23] it was proven also the well-posedness of the time-semidiscrete problem and the well-posedness of the time-steps of the fully discrete system using Stokes stable finite elements. In this thesis, we extend these results to the case of stabilized $\mathbb{P}_1 - \mathbb{P}_1$ fluid finite elements. In Chapter 3 we analyzed the convergence of the monolithic scheme showing a first order convergence rate in time, as expected since we use the Euler scheme, and convergence rates in space that are in agreement with the supposed regularity of the solution of the continuous problem.

The design of the partitioned schemes has been inspired by the splitting paradigm of [60, 56, 1] that was shown to be able to avoid strong coupling without compromising the stability and the accuracy. In Chapter 4, we introduce the partitioned algorithm exploiting a solid splitting that uses a fractional-step time-marching for the computation of the solid velocity.

The analysis performed has shown that the partitioned schemes conceived give results comparables with the results of the monolithic scheme with the advantage to allow for the reuse of existing fluid and solid solvers.

The partitioned schemes introduced, in most of the cases, have the desirable features of unconditional stability, moreover, the theoretical analysis shows that the accuracy of the monolithic and the partitional schemes are comparable. The theoretical results are confirmed by the numerical experiments performed on the test case and are encouraging to further explore the possibilities offered by the splitting introduced.

It is well known that one of the critical points of the non-matching grid methods for incompressible fluid-structure interaction, is related to the imperfect imposition of the free divergence constraint for the fluid at the discrete level. As consequence, even in the case of finite elements that allow for local mass conservation, in the elements cut by the thin structure, there is the fluid that crosses the structure. The numerical experiments have put in evidence that the problem concerns to all the analyzed schemes and can be partially corrected refining the discretization parameters. A different approach could be represented by the possibility to use differents techniques, suited to obtain local conservation of mass, whose implementation requires local enrichment of the fluid finite element space. In the literature we can find solutions to this drawback, for example in [41] it has been proposed an enrichment of the pressure elements. There are also other techniques to cure the problem of nonconservation of mass, we cite here, for example, the method of XFem [1]. The introduction of such techniques could have also the merit to capture the jump of the pressure in the finite element solution. A possible future research topic can be the overcoming of the mentioned problem by using enhanced fluid finite element spaces.

An interesting point for further research could be the extension of the theoretical results to the more complex solid models described in Appendix B. Concerning the three-dimensional shell models, we note that they include higher-order through-the-thickness effects that have to be taken into account. Consequently, the formulation of the coupled problem in the case of immersed shell-type solids seems to be not straightforward.

For what concerns the analysis of the accuracy of the partitioned methods, we focused on the splitting scheme given in Algorithm 3 that shows the property of unconditional stability; in future, it would be interesting to evaluate the accuracy of Algorithm 2 performing also numerical experiments.

Another interesting point, not addressed in this thesis, is the convergence analysis in the nonlinear case of big solid displacements for both the monolithic and the Partitioned schemes. This is a delicate problem which has received little consideration in the literature.

# Appendix A

# Functional setting and finite element discretization

## A.1 Introduction

In this appendix we collect some definitions and results that will be useful in the following analysis of the fluid-structure coupled problem.

## A.2 Variational equations and finite dimensional approximation

In this section we give some definitions and theorems that we will use very often in the following chapters. We begin with the definition of some properties concerning bilinear forms

**Definition 3.** *Let* $\left(\boldsymbol{V},(\cdot,\cdot)_{\boldsymbol{V}}\right)$ *and* $\left(\boldsymbol{W},(\cdot,\cdot)_{\boldsymbol{W}}\right)$ *be Hilbert spaces on* $\mathbb{R}$*, let* $A:\boldsymbol{V}\times\boldsymbol{W}\to\mathbb{R}$ *be a* bilinear form*, then*

- $A(\cdot,\cdot)$ *is continuous if*

$$\|A\| \stackrel{def}{=} \sup_{\boldsymbol{0}\neq\boldsymbol{v}\in\boldsymbol{V},\boldsymbol{0}\neq\boldsymbol{w}\in\boldsymbol{W}} \frac{A(\boldsymbol{v},\boldsymbol{w})}{\|\boldsymbol{v}\|_{\boldsymbol{V}}\|\boldsymbol{w}\|_{\boldsymbol{W}}} < \infty,$$

- **if** $\left(\boldsymbol{V},(\cdot,\cdot)_{\boldsymbol{V}}\right) = \left(\boldsymbol{W},(\cdot,\cdot)_{\boldsymbol{W}}\right)$*,* $A(\cdot,\cdot)$ *is* **coercive** *if there exists* $C>0$ *independent of* $\boldsymbol{v}\in\boldsymbol{V}$ *such that*

$$A(\boldsymbol{v},\boldsymbol{v}) \geq C\|\boldsymbol{v}\|_{\boldsymbol{V}}^2 \qquad \forall\boldsymbol{v}\in\boldsymbol{V}.$$

In the present work we will analyze many variational problems whose general form is

**Problem 11.** *Let* $\left(\boldsymbol{V},(\cdot,\cdot)_{\boldsymbol{V}}\right)$ *and* $\left(\boldsymbol{W},(\cdot,\cdot)_{\boldsymbol{W}}\right)$ *be Hilbert spaces on* $\mathbb{R}$*, let* $A:\boldsymbol{V}\times\boldsymbol{W}\to\mathbb{R}$ *be a* continuous bilinear form *and* $\boldsymbol{F}\in\boldsymbol{W}'$*. Find* $\boldsymbol{v}\in\boldsymbol{V}$ *such that*

$$A(\boldsymbol{v},\boldsymbol{w}) = \langle\boldsymbol{F},\boldsymbol{w}\rangle \qquad \forall\boldsymbol{w}\in\boldsymbol{W}. \tag{A.1}$$

Necessary and sufficient conditions for existence of the solution of Problem 11 are given in the following **generalized Lax-Milgram theorem** [21, 48]

**Theorem A.2.1.** *Let* $(\boldsymbol{V},(\cdot,\cdot)_{\boldsymbol{V}})$ *and* $(\boldsymbol{W},(\cdot,\cdot)_{\boldsymbol{W}})$ *be Hilbert spaces on* $\mathbb{R}$*, let* $A : \boldsymbol{V} \times \boldsymbol{W} \to \mathbb{R}$ *be a* continuous bilinear form. *For all* $\boldsymbol{F} \in \boldsymbol{W}'$*, Problem 11 admits a unique solution if and only if*

- $\exists \alpha > 0$ *such that*

$$\inf_{0 \neq \boldsymbol{v} \in \boldsymbol{V}} \sup_{0 \neq \boldsymbol{w} \in \boldsymbol{W}} \frac{A(\boldsymbol{v},\boldsymbol{w})}{\|\boldsymbol{v}\|_{\boldsymbol{V}} \|\boldsymbol{w}\|_{\boldsymbol{W}}} \geq \alpha \tag{A.2}$$

- 

$$\forall \boldsymbol{v} \in \boldsymbol{V} \backslash \{\boldsymbol{0}\}, \qquad \sup_{\boldsymbol{w} \in \boldsymbol{W}} A(\boldsymbol{v},\boldsymbol{w}) > 0. \tag{A.3}$$

*Moreover, if* $\boldsymbol{v} \in \boldsymbol{V}$ *is the solution of Problem 11, we have the following estimate*

$$\|\boldsymbol{v}\|_{\boldsymbol{V}} \leq \frac{1}{\alpha} \|\boldsymbol{F}\|_{\boldsymbol{W}'}. \tag{A.4}$$

A sufficient condition for the existence of the solution is given by the **Lax-Milgram lemma** whose conditions are more simple to check with respect to the conditions of Theorem A.2.1, [21, 48]

**Lemma A.2.2.** *Let* $A : \boldsymbol{V} \times \boldsymbol{V} \to \mathbb{R}$ *be a continuous and coercive bilinear form on the Hilbert space* $(\boldsymbol{V},(\cdot,\cdot)_{\boldsymbol{V}})$*, then for all* $\boldsymbol{F} \in \boldsymbol{V}'$ *there exists a unique* $\boldsymbol{v} \in \boldsymbol{V}$ *such that*

$$A(\boldsymbol{v},\boldsymbol{w}) = \langle \boldsymbol{F},\boldsymbol{w} \rangle \qquad \forall \boldsymbol{w} \in \boldsymbol{V},$$

*moreover*

$$\|\boldsymbol{v}\|_{\boldsymbol{V}} \leq \frac{1}{\alpha} \|\boldsymbol{F}\|_{\boldsymbol{V}'}, \tag{A.5}$$

*where* $\alpha$ *is the coercivity constant of* $A(\cdot,\cdot)$*.*

**Finite dimensional approximation** In this paragraph we collect some basic facts about the approximation of the solution to weak problems using finite dimensional subspaces. We limit ourselves to the conforming case since this is the context in which we will work in the following chapters. Let $\{\boldsymbol{V}_h\}_{h>0}$ and $\{\boldsymbol{W}_h\}_{h>0}$ be families of finite dimensional spaces, such that, for all $h > 0$, $\boldsymbol{V}_h \subset \boldsymbol{V}$ and $\boldsymbol{W}_h \subset \boldsymbol{W}$ respectively. For all $h > 0$, we have a finite dimensional problem of the following type

**Problem 12.** *Let* $\boldsymbol{V}_h \subset \boldsymbol{V}$ *and* $\boldsymbol{W}_h \subset \boldsymbol{W}$*, let* $A : \boldsymbol{V} \times \boldsymbol{W} \to \mathbb{R}$ *be a* continuous bilinear form *and* $\boldsymbol{F} \in \boldsymbol{W}'$*, then we consider the following discrete variational problem: find* $\boldsymbol{v}_h \in \boldsymbol{V}_h$ *such that*

$$A(\boldsymbol{v}_h,\boldsymbol{w}_h) = \langle \boldsymbol{F},\boldsymbol{w}_h \rangle \qquad \forall \boldsymbol{w}_h \in \boldsymbol{W}_h. \tag{A.6}$$

In the general case, the well-posedness of Problem 12 is not assured by the well-posedness of Problem 11. In this regard, in order to have solution for all finite dimensional problems, we need to fulfill the requirements of the generalized Lax-Milgram Theorem A.2.1 for each discrete problem, that is for each $h$, there exists $\alpha_h > 0$ such that

$$\inf_{\mathbf{0}\neq \boldsymbol{w}_h\in \boldsymbol{W}_h} \sup_{\mathbf{0}\neq \boldsymbol{v}_h\in \boldsymbol{V}_h} \frac{A(\boldsymbol{v}_h,\boldsymbol{w}_h)}{\|\boldsymbol{v}_h\|_{\boldsymbol{V}}\|\boldsymbol{w}_h\|_{\boldsymbol{W}}} \geq \alpha_h, \qquad \sup_{\boldsymbol{w}_h\in \boldsymbol{W}_h} A(\boldsymbol{v}_h,\boldsymbol{w}_h) > 0 \qquad \forall \boldsymbol{v}_h \in \boldsymbol{V}_h\setminus\{\mathbf{0}\}. \tag{A.7}$$

In particular, for fixed $\boldsymbol{V}_h$ and $\boldsymbol{W}_h$ we have the following [21, 48]

**Proposition A.2.3.** *Let $\boldsymbol{F} \in \boldsymbol{W}'$, suppose that $A(\cdot,\cdot)$ satisfies the hypotheses of the generalized Lax-Milgram Theorem A.2.1 on $\boldsymbol{V}\times\boldsymbol{W}$ and $\boldsymbol{V}_h\times\boldsymbol{W}_h$ respectively, let $\boldsymbol{v}\in\boldsymbol{V}$ and $\boldsymbol{v}_h\in\boldsymbol{V}_h$ be the solutions to Problems 11 and 12 respectively, then*

- *The error $\boldsymbol{v}-\boldsymbol{v}_h$ fulfills the following* **orthogonality relation**

$$A(\boldsymbol{v}-\boldsymbol{v}_h,\boldsymbol{w}_h) = 0 \qquad \forall \boldsymbol{w}_h \in \boldsymbol{W}_h.$$

- *Assume* $\dim \boldsymbol{V}_h = \dim \boldsymbol{W}_h$, *then it holds the following* **Céa inequality**

$$\|\boldsymbol{v}-\boldsymbol{v}_h\|_{\boldsymbol{V}} \leq \left(1+\frac{\|A\|}{\alpha_h}\right) \inf_{\boldsymbol{z}_h\in \boldsymbol{V}_h} \|\boldsymbol{v}-\boldsymbol{z}_h\|_{\boldsymbol{V}},$$

*where $\alpha_h > 0$ is the inf-sup constant of $A(\cdot,\cdot)$ on $(\boldsymbol{V}_h,\boldsymbol{W}_h)$ from* (A.7).

**Remark 3.** *In the Céa Lemma given in Proposition A.2.3 the requirement* $\dim \boldsymbol{V}_h = \dim \boldsymbol{W}_h$ *assures that the two conditions of the generalized Lax-Milgram Theorem A.2.1 are in fact equivalent in the finite dimensional case.*

In order to have convergence of the finite element solutions $\{\boldsymbol{v}_h\}_{h>0}$ to the continuum solution $\boldsymbol{v}$ of Problem 11, we require that the inf-sup conditions are fulfilled uniformly in $h$, namely

$$\inf_{h>0} \alpha_h = \alpha > 0.$$

In that case we have the following **error estimate**

$$\|\boldsymbol{v}-\boldsymbol{v}_h\|_{\boldsymbol{V}} \leq \left(1+\frac{\|A\|}{\alpha}\right) \inf_{\boldsymbol{z}_h\in \boldsymbol{V}_h} \|\boldsymbol{v}-\boldsymbol{z}_h\|_{\boldsymbol{V}}.$$

## A.3 Function spaces

In this thesis we use many Hilbert function spaces; elements of such spaces are defined on some open set $U\subset\mathbb{R}^k$ with values in $\mathbb{R}^d$. We refer to [21, 28, 87] for the details.

**Definition 4.** *Let $\boldsymbol{f}$ be a measurable function defined on an open set $U\subset\mathbb{R}^k$ ($k\in\mathbb{N}$) with values in $\mathbb{R}^d$, then*

- $$\boldsymbol{f} \in L^2(U)^d \quad \overset{def}{\Longleftrightarrow} \quad \int_U |f_1|^2 + \cdots + |f_d|^2 d\boldsymbol{x} < \infty, \tag{A.8}$$

  *for all* $\boldsymbol{f}$ *and* $\boldsymbol{g}$ *in* $L^2(U)^d$ *the scalar product and the associated norm are*

  $$(\boldsymbol{f},\boldsymbol{g})_{0,U} \overset{def}{=} \int_U \boldsymbol{f}\cdot\boldsymbol{g}d\boldsymbol{x}; \qquad \|\boldsymbol{f}\|_{0,U} \overset{def}{=} \sqrt{(\boldsymbol{f},\boldsymbol{f})_{0,U}}. \tag{A.9}$$

- *For* $k \in \mathbb{N}$ *and* $\alpha \in \mathbb{N}^d$

  $$\boldsymbol{f} \in H^k(U)^d \quad \overset{def}{\Longleftrightarrow} \quad \partial_\alpha \boldsymbol{f} \in L^2(U)^d \quad \textit{for all} \quad 0 \le |\alpha| \le k, \tag{A.10}$$

  *where* $\partial_\alpha \boldsymbol{f}$ *stands for the distributional derivative of* $\boldsymbol{f}$. *Moreover we have, for all* $\boldsymbol{f}$ *and* $\boldsymbol{g}$ *in* $H^k(U)^d$

  $$(\boldsymbol{f},\boldsymbol{g})_{k,U} \overset{def}{=} \sum_{0\le|\alpha|\le k} \int_U \partial_\alpha\boldsymbol{f}\cdot\partial_\alpha\boldsymbol{g}d\boldsymbol{x}; \qquad \|\boldsymbol{f}\|_{k,U} \overset{def}{=} \sqrt{(\boldsymbol{f},\boldsymbol{f})_{k,U}}. \tag{A.11}$$

- *For* $0 < \mu < 1$

  $$\boldsymbol{f} \in H^\mu(U)^d \quad \overset{def}{\Longleftrightarrow} \quad \boldsymbol{f} \in L^2(U)^d \quad \textit{and} \quad \left(\int_U\int_U \frac{|\boldsymbol{f}(\boldsymbol{x})-\boldsymbol{f}(\boldsymbol{y})|^2}{|\boldsymbol{x}-\boldsymbol{y}|^{d+2\mu}} d\boldsymbol{x}d\boldsymbol{y}\right)^{\frac{1}{2}} < \infty, \tag{A.12}$$

  *moreover, for all* $\boldsymbol{f}$ *and* $\boldsymbol{g}$ *in* $H^\mu(U)^d$ *the scalar product and the norm are*

  $$(\boldsymbol{f},\boldsymbol{g})_{\mu,U} \overset{def}{=} \int_U\int_U \frac{|\boldsymbol{f}(\boldsymbol{x})-\boldsymbol{f}(\boldsymbol{y})||\boldsymbol{g}(\boldsymbol{x})-\boldsymbol{g}(\boldsymbol{y})|}{|\boldsymbol{x}-\boldsymbol{y}|^{d+2\mu}} d\boldsymbol{x}d\boldsymbol{y}; \qquad \|\boldsymbol{f}\|_{\mu,U} \overset{def}{=} \sqrt{(\boldsymbol{f},\boldsymbol{f})_{\mu,U}} \tag{A.13}$$

- *For* $s = k + \mu$ *with* $k \in \mathbb{N}$, $0 < \mu < 1$ *and* $\alpha \in \mathbb{N}^d$

  $$\boldsymbol{f} \in H^s(U)^d \quad \overset{def}{\Longleftrightarrow} \quad \boldsymbol{f} \in H^k(U)^d \quad \textit{and} \quad \partial_\alpha \boldsymbol{f} \in H^\mu(U)^d \quad \textit{for all} \quad |\alpha| = k. \tag{A.14}$$

  *moreover, for all* $\boldsymbol{f}$ *and* $\boldsymbol{g}$ *in* $H^s(U)^d$ *the scalar product and the norm are*

  $$(\boldsymbol{f},\boldsymbol{g})_{s,U} \overset{def}{=} (\boldsymbol{f},\boldsymbol{g})_{k,U} + \sum_{|\alpha|=k} (\partial_\alpha\boldsymbol{f},\partial_\alpha\boldsymbol{g})_{\mu,U}, \qquad \|\boldsymbol{f}\|^2_{s,U} = \|\boldsymbol{f}\|^2_{k,U} + \sum_{|\alpha|=k} \|\partial_\alpha\boldsymbol{f}\|^2_{\mu,U}. \tag{A.15}$$

It is well known that $L^2(U)$ and $H^s(U)^d$ for all $s > 0$ are separable Hilbert spaces. Another useful space is $H^1_0(U)^d$, defined as the closure of the space of $C^\infty_0(U)^d$ with respect to the $H^1$-norm

$$H^1_0(U)^d \overset{def}{=} \overline{C^\infty_0(U)^d}^{\|\cdot\|_{1,\Omega}}. \tag{A.16}$$

For the space $H^1_0(U)^d$ we have a useful result, namely the Poincaré Lemma

**Lemma A.3.1.** *Let $U \subset \mathbb{R}^k$ a bounded open set, then for any $\boldsymbol{u} \in H_0^1(U)^d$*

$$\|\boldsymbol{u}\|_{0,\Omega} \leq C\|\nabla \boldsymbol{u}\|_{0,\Omega} \tag{A.17}$$

*with C independent of $\boldsymbol{u}$.*

Poincaré Lemma can be generalized, first of all it holds for functions in $H_\Gamma^1(U)^d$, that is functions in $H^1(U)^d$ with "zero trace" on a part of the boundary $\Gamma$ with positive (d-1)-measure (the meaningful definition of trace will be given in Proposition A.3.2).

$$\|\boldsymbol{u}\|_{0,\Omega} \leq C\|\nabla \boldsymbol{u}\|_{0,\Omega}, \qquad \forall \boldsymbol{u} \in H_\Gamma^1(U)^d. \tag{A.18}$$

Another generalization of Lemma A.3.1, valid for all functions in $H^1(U)^d$, is the following estimate

$$\|\boldsymbol{u} - \overline{\boldsymbol{u}}\|_{0,\Omega} \leq C\|\nabla \boldsymbol{u}\|_{0,\Omega}, \tag{A.19}$$

where

$$\overline{\boldsymbol{u}} = \frac{1}{|U|}\int_U \boldsymbol{u}\, d\boldsymbol{x}. \tag{A.20}$$

Poincaré Lemma implies that the $H^1$-semi-norm ($\|\nabla \boldsymbol{u}\|_{0,\Omega}$) is, in fact, a norm on $H_\Gamma^1(U)^d$ equivalent to the $H^1$-norm, that is there exist constants $C_p$ and $C$ such that

$$C_p\|\boldsymbol{u}\|_{1,\Omega} \leq \|\nabla \boldsymbol{u}\|_{0,\Omega} \leq C\|\boldsymbol{u}\|_{1,\Omega} \qquad \forall \boldsymbol{u} \in H_\Gamma^1(U)^d. \tag{A.21}$$

**Remark 4.** *In all the definitions and theorems we will give from now on we will consider* **convex polygonal domain** $\Omega \subset \mathbb{R}^d$ *with $d = 2$ or $d = 3$ and the* **one-dimensional interval** $\Sigma = [0, L] \subset \mathbb{R}$. *This setting is enough for our purposes but many results presented are valid also for more general open sets. The precise meaning of polygonal domain is given in the following*

**Definition 5.** *Let $d \geq 1$ and $\Omega \subset \mathbb{R}^d$, then*

- $\Omega$ *is a* **domain** *if $\Omega$ is a bounded, connected and open set,*
- $\Omega$ *is a* **Lipschitz domain** *if $\Omega$ is a domain and its boundary is "locally the graph of a Lipschitz function" (see [40])*
- $\Omega$ *is a* **polygonal domain** *if $\Omega$ is a domain and $\overline{\Omega}$ is union of finitely many polyhedra.*

We observe that polygonal domains are in fact Lipschitz domains. Function spaces defined above present very useful properties when their domain of definition is a Lipschitz domain, and in particular a polygonal domain. We list some of these properties in the following proposition (see [87])

**Proposition A.3.2.** *Let $\Omega \subset \mathbb{R}^d$ be a Lipschitz domain (in particular a polygonal domain), then*

1. *$C^\infty(\overline{\Omega})^d$ is dense in $H^s(\Omega)^d$ for $s \geq 0$;*

2. *for $s_1 \geq s_2 \geq 0$, $H^{s_1}(\Omega)^d \hookrightarrow H^{s_2}(\Omega)^d$ with dense immersion;*

3. *the* **trace operator** *$\gamma : C^\infty(\overline{\Omega})^d \to C^\infty(\partial\Omega)^d$, defined by $\gamma \boldsymbol{u} = \boldsymbol{u}\big|_{\partial\Omega}$, has unique extension to a bounded linear operator*

$$\gamma : H^s(\Omega)^d \to H^{s-\frac{1}{2}}(\partial\Omega)^d \qquad \text{for} \qquad \frac{1}{2} < s < \frac{3}{2}, \tag{A.22}$$

   *and this extension has a continuous right inverse if $\frac{1}{2} < s \leq 1$;*

4. *if $0 \leq s \leq \frac{1}{2}$, then $H^s_0(\Omega) = H^s(\Omega)$;*

5. *if $\frac{1}{2} < s \leq 1$, then $H^s_0(\Omega)^d = \{\boldsymbol{u} \in H^s(\Omega)^d : \gamma(\boldsymbol{u}) = \boldsymbol{0}\}$;*

6. *if $k \in \mathbb{N}$ and $2k > d$ then for all $\boldsymbol{v} \in H^k(\Omega)^d$ there is a continuous bounded function $\tilde{\boldsymbol{v}} \in C_b(\Omega)^d$ such that*

$$\tilde{\boldsymbol{v}} \overset{a.e.}{=} \boldsymbol{v} \qquad \text{and} \qquad \sup_{\Omega} \|\tilde{\boldsymbol{v}}\|_E \leq C\|\boldsymbol{v}\|_{k,\Omega}. \tag{A.23}$$

Point 3 in Proposition A.3.2 implies the existence of a constant $C$ such that for all $\boldsymbol{v} \in H^s(\Omega)^d$, $\|\gamma(\boldsymbol{v})\|_{s-\frac{1}{2},\partial\Omega} \leq C\|\boldsymbol{v}\|_{s,\Omega}$. Regarding the trace of functions in Sobolev spaces, for Lipschitz domain, we will use also the following theorem (see [28, Theorem 1.6.6])

**Theorem A.3.3.** *Suppose that $\Omega$ has a Lipschitz boundary, then there exists a constant $C > 0$, such that*

$$\|\gamma v\|_{0,\partial\Omega} \leq C\|v\|_{0,\Omega}^{\frac{1}{2}}\|v\|_{1,\Omega}^{\frac{1}{2}} \qquad \forall v \in H^1(\Omega) \tag{A.24}$$

In particular a consequence of Theorem A.3.3 is that there is a constant, $C > 0$, such that $\|v\|_{0,\partial\Omega} \leq C\|v\|_{1,\Omega}$ for all $v \in H^1(\Omega)$.

Another theorem used quite often in the following is the **Korn's Theorem** [40] in which we introduce the symmetric gradient defined as

$$\boldsymbol{\epsilon}(\boldsymbol{v}) \overset{def}{=} \frac{1}{2}\left(\nabla\boldsymbol{v} + \nabla\boldsymbol{v}^\top\right)$$

**Theorem A.3.4.** *Let $\Omega \subset \mathbb{R}^d$ be a Lipschitz domain, then*

- **First Korn's inequality** *Let $\Gamma \subset \partial\Omega$ with positive measure and $H^1_\Gamma(\Omega)^d = \{\boldsymbol{v} \in H^1(\Omega)^d : \boldsymbol{v}\big|_\Gamma = \boldsymbol{0}\}$, then there exist $C > 0$ such that*

$$\forall \boldsymbol{v} \in H^1_\Gamma(\Omega)^d, \qquad \|\boldsymbol{v}\|_{1,\Omega} \leq C\|\boldsymbol{\epsilon}(\boldsymbol{v})\|_{0,\Omega}.$$

- **Second Korn's inequality** *There exist $C > 0$ such that*

$$\forall \boldsymbol{v} \in H^1(\Omega)^d, \qquad \|\boldsymbol{v}\|_{1,\Omega} \leq C\left(\|\boldsymbol{\epsilon}(\boldsymbol{v})\|^2_{0,\Omega\cdot} + \|\boldsymbol{v}\|^2_{0,\Omega}\right)^{\frac{1}{2}}.$$

In the thesis we use the fact that *fractional Sobolev spaces*, are, in fact, interpolation spaces between Sobolev spaces of integer index. A very important property of interpolation spaces is expressed in the following proposition and concerns linear operators defined on Sobolev spaces [28, Proposition 14.1.5]

**Proposition A.3.5.** *Let $T$ be a linear operator, such that $T : H^k(\Omega) \to H^m(\Omega)$ and $T : H^{k+l}(\Omega) \to H^{m+t}(\Omega)$, suppose that $T$ is continuous in both the cases, then for any $0 \le \theta \le 1$, $T$ is a continuous linear operator from $H^{k+\theta l}(\Omega)$ to $H^{m+\theta t}(\Omega)$, moreover*

$$\|T\|_{\mathscr{L}(H^{k+\theta l}(\Omega),H^{m+\theta t}(\Omega))} \le \|T\|^{1-\theta}_{\mathscr{L}(H^k(\Omega),H^m(\Omega))}\|T\|^{\theta}_{\mathscr{L}(H^{k+l}(\Omega),H^{m+t}(\Omega))}$$

The last proposition is very useful in the applications since it allows for the extension of properties proved for integer index Sobolev spaces to real index Sobolev spaces.

## A.4 Finite element convergence estimates in Sobolev spaces

### A.4.1 Triangulation

If $\overline{\Omega} \subset \mathbb{R}^d$ (d=2,3) is a polygonal domain we can consider a finite decomposition into polyhedra

$$\overline{\Omega} = \cup_{K\in\mathscr{T}_h} K \tag{A.25}$$

such that

1. each polyhedron is such that $\overset{\circ}{K} \neq \emptyset$,
2. $\overset{\circ}{K}_1 \cap \overset{\circ}{K}_2 = \emptyset$ for $K_1$ and $K_2$ in $\mathscr{T}_h$ distinct,
3. If $F = K_1 \cap K_2 \neq \emptyset$ and $K_1 \neq K_2$, then F is a common face, side or vertex of $K_1$ and $K_2$,
4. diam(K) $\le$ h for each $K \in \mathscr{T}_h$.

$\mathscr{T}_h$ called a *triangulation* of $\overline{\Omega}$. We assume moreover, that each polyhedron is obtained as $K = \boldsymbol{T}_K(\hat{K})$ with $\hat{K}$ a reference polyhedron and $\boldsymbol{T}_K$ affine map, namely $\boldsymbol{T}_K(\hat{\boldsymbol{x}}) = \boldsymbol{B}_K\hat{\boldsymbol{x}} + \boldsymbol{b}_K$ with denoting a non singular matrix $\boldsymbol{B}_K$ and a vector $\boldsymbol{b}_K$. In order to evaluate the convergence of the finite element method, we are interested to family of triangulations $\{\mathscr{T}_h\}_{h>0}$. we give the following definition

**Definition 6.** *Let $\{\mathscr{T}_h\}_{h>0}$ be a family of triangulations of $\overline{\Omega}$, then*

- *given a triangulation $\mathscr{T}_h \in \{\mathscr{T}_h\}_{h>0}$, we define*

$$\textit{for all} \quad K \in \mathscr{T}_h \quad h_K \overset{def}{=} \mathrm{diam(K)}, \qquad \rho_{\mathrm{K}} \overset{\mathrm{def}}{=} \max_{\mathrm{B}\subset\mathrm{K}} \mathrm{diam(B)} \quad \textit{with} \quad \boldsymbol{B} \quad \textit{ball}$$

$$h \overset{def}{=} \max_{K\in\mathscr{T}_h} h_K$$

- $\{\mathscr{T}_h\}_{h>0}$ *is a* **regular family** *if there exists* $\sigma > 0$ *such that*

$$\max_{K\in\mathscr{T}_h} \frac{h_K}{\rho_K} \le \sigma \qquad \forall h > 0.$$

- $\{\mathscr{T}_h\}_{h>0}$ *is a* **quasi-uniform family** *if*
  1. $\{\mathscr{T}_h\}_{h>0}$ *is regular,*
  2. *there exists a constant* $\sigma_1 > 0$ *such that*

$$\min_{K\in\mathscr{T}_h} h_K \ge \sigma_1 h, \qquad \forall h > 0$$

.

### A.4.2 Piecewise-Polynomial Subspaces

A basic aspect of finite element method is that we try to approximate the space of the solution of a given problem using finite dimensional spaces of piecewise-polynomial functions. In order to construct piecewise polynomial functions we will consider the space $\mathbb{P}_k$ of polynomials of degree less that or equal to some natural number $k$ in the spatial variables $x_1, \cdots, x_d$. Then we define the following finite dimensional space. Given a triangulation $\mathscr{T}_h$ of $\overline{\Omega}$

$$\boldsymbol{X}_h^k = \{v_h \in C^0(\overline{\Omega}) : v_h|_K \in \mathbb{P}_k, \forall K \in \mathscr{T}_h\}. \tag{A.26}$$

A very well known property (see [103, Proposition 3.2.1]) of $\boldsymbol{X}_h^k$ is that for all $k \ge 1$

$$\boldsymbol{X}_h^k \subset H^1(\Omega). \tag{A.27}$$

It is now necessary to construct a basis for the space $\boldsymbol{X}_h^k$. A good starting point is to give a basis for $\mathbb{P}_k$ on each element $K$ of $\mathscr{T}_h$. Following the standard definition of *finite element*, in order to assign a basis for a finite dimensional function space $\mathscr{P}$ (with dim$\mathscr{P}$ = n < ∞), we assign a basis for the dual space, namely a family of linearly independent functionals $\mathscr{N} = \{\Phi_1, \cdots, \Phi_n\}$ in $\mathscr{P}'$, then the basis for $\mathscr{P}$ is the dual basis to $\mathscr{N}$, namely the set $\{\phi_1, \cdots, \phi_n\} \subset \mathscr{P}$ such that $\Phi_i(\phi_j) = \delta_{ij}$; then each element $\boldsymbol{v}_h \in \mathscr{P}$ is expressed as

$$\boldsymbol{v}_h = \sum_{i=1}^n v_h^i \phi_i, \qquad \text{with} \qquad v_h^i = \Phi_i(\boldsymbol{v}_h).$$

In our work we use *Lagrange finite elements*, namely elements for which the functionals in $\mathscr{N}$ are the evaluation at given points $\{\boldsymbol{a}_j\}_{j=1}^n$ named *nodes*. Then, in order to assign a basis for $\mathbb{P}_k$ on each triangle K, we have to assign the evaluation points on each triangle, namely

$$\{\boldsymbol{a}_j\}_{j=1}^n \subset K, \qquad \text{with} \qquad \dim \mathbb{P}_k = n,$$

such that the associated evaluation functionals $\{\Phi_1, \cdots, \Phi_n\}$ are independent. Since the functions in $\boldsymbol{X}_{\boldsymbol{h}}^{\boldsymbol{k}}$ are globally continuous, we are forced to choose some nodes on the faces of the elements, in such a way that functions obtained gluing the polynomials computed on adjacent elements are continuous on the union of these elements. Denoting by $N_h$ the total number of nodes in $\overline{\Omega}$, a basis for $\boldsymbol{X}_{\boldsymbol{h}}^{\boldsymbol{k}}$ is now easily constructed considering the set $\{\boldsymbol{a}_j\}_{j=1}^{N_h} \subset \overline{\Omega}$ and the **shape functions** $\varphi_i \in \boldsymbol{X}_h^k$ such that

$$\varphi_i(\boldsymbol{a}_j) = \delta_{ij}, \qquad i, j = 1, \cdots, N_h. \tag{A.28}$$

### A.4.3 Interpolation operator

In the error analysis of finite element approximation it is important to control the error given by the choice of the space $\boldsymbol{X}_h^k$, namely given a function in some Hilbert space $v \in \boldsymbol{H}$ with $\boldsymbol{X}_h^k \subset \boldsymbol{H}$, we want to control the quantity

$$\inf_{v_h \in \boldsymbol{X}_h^k} \|v - v_h\|_{\boldsymbol{H}}.$$

We reach the goal if we are able to find a continuous linear operator $\Pi_h : \boldsymbol{H} \to \boldsymbol{X}_h^k$ ; in fact, in this case we have

$$\inf_{v_h \in \boldsymbol{X}_h^k} \|v - v_h\|_{\boldsymbol{H}} \leq \|v - \Pi_h(v)\|_{\boldsymbol{H}} \leq \|I - \Pi_h\|_{\mathscr{L}(\boldsymbol{H},\boldsymbol{H})} \|v\|_{\boldsymbol{H}}.$$

**Interpolation of continuous functions** If the functions in $\boldsymbol{H}$ are continuous (in the case of Sobolev spaces $H^k(\Omega)$ this is the case when $2k > d$ see Proposition A.3.2) we can define

$$\mathscr{I}_h(v) \stackrel{def}{=} \sum_{i=1}^{N_h} v(\boldsymbol{a}_i)\phi_i \qquad \forall v \in C^0(\overline{\Omega}). \tag{A.29}$$

then we have the following result (see [103, Theorem 3.4.2])

**Theorem A.4.1.** *Let $\mathscr{T}_h$ be a regular family of triangulations, assume that $m = 0, 1$ and $n = \min\{k, r-1\} \geq 1$ are all integers, then there exist a constant $C > 0$, independent of $h$, such that*

$$|v - \mathscr{I}_h(v)|_{m,\Omega} \leq C h^{n+1-m} |v|_{n+1,\Omega} \qquad \forall v \in H^r(\Omega) \tag{A.30}$$

**Projection operators** For functions that are not continuous, the interpolation operator introduced in (A.29) is not well defined. In these cases, since we work in Hilbert spaces, in order to obtain optimal error estimates, we can use the orthogonal projection operator. Let $\boldsymbol{S}$ be a closed subspace of the Hilbert space $\boldsymbol{H}$, the projection of $v \in \boldsymbol{H}$ over $\boldsymbol{S}$ is defined as

$$P_S(v) \in S \quad : \quad (P_S(v), w)_H = (v, w)_H, \qquad \forall w \in \boldsymbol{S} \tag{A.31}$$

It is characterized by the property

$$\|v - P_S(v)\|_H = \min_{w \in H} \|v - w\|_H. \tag{A.32}$$

In particular, we define the following projection operators

$$P_h^0 : L^2(\Omega) \to \boldsymbol{X}_h^k, \tag{A.33}$$

$$P_h^1 : H^1(\Omega) \to \boldsymbol{X}_h^k. \tag{A.34}$$

We have the following error estimates [103], for $v \in H^{l+1}(\Omega)$ with $l$ **integer**. Let $k$ be the polynomial degree, then

$$\|v - P_h^0(v)\|_{0,\Omega} \leq C h^{l+1} |v|_{l+1,\Omega}, \qquad 1 \leq l \leq k, \tag{A.35}$$

$$\|v - P_h^1(v)\|_{1,\Omega} \leq C h^{l} |v|_{l+1,\Omega}, \qquad 1 \leq l \leq k. \tag{A.36}$$

**Clément interpolation** Another Interpolation operator that will be used many times in the following is the Clément interpolation operator that allows for the approximation of functions in $L^1(\Omega)$ using $H^1$-conforming Lagrange finite elements. The definition of this operator is not useful for our purposes, then we skip this point and refer to [48]. The important point about the Clément operator are the stability and interpolation properties stated in the following proposition [48]

**Proposition A.4.2.** *Let $\mathscr{C}_h : L^1(\Omega) \to \boldsymbol{X}_h^k$ be the Clément interpolation operator, then*

- *Stability: for $0 \le s \le 1$, there exists a constant $C > 0$, independent of h, such that*

$$\|\mathscr{C}_h \boldsymbol{v}\|_{s,\Omega} \le C\|\boldsymbol{v}\|_{s,\Omega} \qquad \forall \boldsymbol{v} \in H^s(\Omega). \tag{A.37}$$

- *Approximation: for $K \in \mathscr{T}_h$, denote by $\Delta_K$ the union of elements in $\mathscr{T}_h$ sharing at least one vertex with K. Let $0 \le s \le l \le k+1$, then there exists a constant $C > 0$, independent of K and h, such that*

$$\|\boldsymbol{v} - \mathscr{C}_h \boldsymbol{v}\|_{s,K} \le C h_K^{l-s} \|\boldsymbol{v}\|_{l,\Delta_K} \qquad \forall \boldsymbol{v} \in H^l(\Delta_K). \tag{A.38}$$

**Scott-Zhang interpolation** Clément interpolation is not a projection and it does not preserve boundary conditions. Scott-Zhang interpolation operator gives a solution to this two issues. In the case of our interest Scott-Zhang interpolation is defined for function in $H^l(\Omega)$ with $l > \frac{1}{2}$. The definition of this operator is not useful for our purposes, then we skip this point and refer to [48]. Instead we give the stability and interpolation properties of Scott-Zhang interpolation stated in the following proposition [48]

**Proposition A.4.3.** *Let $\mathscr{S}\mathscr{Z}_h : H^l(\Omega) \to \boldsymbol{X}_h^k$ be the Scott-Zhang interpolation operator, then*

- *Projection property and boundary values:*

$$\forall \boldsymbol{v}_h \in \boldsymbol{X}_h^k, \qquad \mathscr{S}\mathscr{Z}_h \boldsymbol{v}_h = \boldsymbol{v}_h; \qquad \forall \boldsymbol{v} \in H^l(\Omega), \qquad \boldsymbol{v}\big|_{\partial\Omega} = \boldsymbol{0} \Rightarrow (\mathscr{S}\mathscr{Z}_h \boldsymbol{v})\big|_{\partial\Omega} = \boldsymbol{0}. \tag{A.39}$$

- *Stability: let $0 \le s \le \min\{1, l\}$, then there exists a constant $C > 0$, independent of h, such that*

$$\|\mathscr{S}\mathscr{Z}_h \boldsymbol{v}\|_{s,\Omega} \le C\|\boldsymbol{v}\|_{l,\Omega} \qquad \forall \boldsymbol{v} \in H^l(\Omega). \tag{A.40}$$

- *Approximation: for $K \in \mathscr{T}_h$, denote by $\Delta_K$ the union of elements in $\mathscr{T}_h$ sharing at least one vertex with K. Let $0 \le s \le l \le k+1$, then there exists a constant $C > 0$, independent of K and h, such that*

$$\|\boldsymbol{v} - \mathscr{S}\mathscr{Z}_h \boldsymbol{v}\|_{s,K} \le C h_K^{l-s} \|\boldsymbol{v}\|_{l,\Delta_K} \qquad \forall \boldsymbol{v} \in H^l(\Delta_K). \tag{A.41}$$

### A.4.4 Inverse inequalities

It is well known that in finite dimensional spaces all norms are equivalents; we have a useful result given in the following lemma (see [48, 103])

**Lemma A.4.4.** *Let* $\{\mathcal{T}_h\}_{1\geq h>0}$ *be a family of triangulations of* $\Omega$, *let* $X_h^k \subset H^m(\Omega) \subset H^n(\Omega)$ *with integers* $0 \leq n \leq m$, *then*

- *if* $\{\mathcal{T}_h\}_{h>0}$ *is regular, there exist a constants* $C > 0$ *independent on* $h$ *such that for all* $v_h \in X_h^k$ *and for all* $K \in \mathcal{T}_h$

$$\|v_h\|_{m,K} \leq C \frac{1}{h_K^{m-n}} \|v_h\|_{n,K}, \tag{A.42}$$

- *if* $\{\mathcal{T}_h\}_{h>0}$ *is quasi-uniform, then there exist a constants* $C > 0$ *independent on* $h$ *such that for all* $v_h \in X_h^k$

$$\|v_h\|_{m,\Omega} \leq C \frac{1}{h^{m-n}} \|v_h\|_{n,\Omega} \tag{A.43}$$

### A.4.5 Estimates and inverse inequality in fractional Sobolev spaces

In many results in the thesis use estimates in fractional Sobolev spaces $H^s(\Omega)$ with $s \geq 0$. In order to obtain the desired estimates we use Proposition A.3.5

**Proposition A.4.5.** *Assume that the family of triangulations of* $\Omega$ $\{\mathcal{T}_h\}_{h>0}$ *is regular, let* $X_h^k \subset H^l(\Omega) \subset H^m(\Omega)$ *with* $0 \leq m < l$ *integers, let* $s \in [0,\infty)$ *such that* $m \leq s \leq l$ , *let* $I_h : H^m(\Omega) \to \boldsymbol{X}_h^k$ *a linear interpolation operator with the following properties*

- $\|I_h v\|_{m,\Omega} \leq C\|v\|_{m,\Omega} \qquad \forall v \in H^m(\Omega),$
- $\|v - I_h v\|_{m,\Omega} \leq Ch^{l-m}\|v\|_{l,\Omega} \qquad \forall v \in H^l(\Omega),$

*then*

- *for all* $v \in H^s(\Omega)$

$$\|v - I_h(v)\|_{m,\Omega} \leq Ch^{s-m}\|v\|_{s,\Omega}, \tag{A.44}$$

- *for all* $K \in \mathcal{T}_h$, *for all* $v_h \in \boldsymbol{X}_h^k$

$$\|v_h\|_{l,K} \leq C \frac{1}{h_K^{l-s}} \|v_h\|_{s,K}, \tag{A.45}$$

*Proof.* The first point follows if we apply Proposition A.3.5 to the operator $\boldsymbol{T} : H^m(\Omega) \to H^m(\Omega)$ defined, for all $v \in H^m(\Omega)$ as

$$\boldsymbol{T}(v) \stackrel{def}{=} v - I_h(v), \tag{A.46}$$

then,

$$\|T\|_{\mathscr{L}(H^m(\Omega),H^m(\Omega))} \leq C, \qquad \text{and} \qquad \|T\|_{\mathscr{L}(H^l(\Omega),H^m(\Omega))} \leq Ch^{l-m}$$

applying Proposition A.3.5 we have

$$\|T\|_{\mathscr{L}(H^s(\Omega),H^m(\Omega))} \leq Ch^{s-m}.$$

The definition of the operator norm gives, for all $v \in H^s(\Omega)$

$$\|v - I_h(v)\|_{m,\Omega} \leq Ch^{s-m}\|v\|_{s,\Omega}.$$

The second point comes applying Proposition A.3.5 to the injection operator $\boldsymbol{j}$ such that $\boldsymbol{j} : \boldsymbol{X}_h^k \to H^l(\Omega)$ and $\boldsymbol{j} : \boldsymbol{X}_h^k \to H^m(\Omega)$. For all $v_h \in \boldsymbol{X}_h^k$ by local inverse inequality (A.42) there holds:

$$\|\boldsymbol{j}(v_h)\|_{m,\Omega} = \|v_h\|_{m,\Omega}, \qquad \text{and} \qquad \|\boldsymbol{j}(v_h)\|_{l,\Omega} = \|v_h\|_{l,\Omega} \leq C\frac{1}{h_K^{l-m}}\|v_h\|_{m,\Omega},$$

hence, thanks to Proposition A.3.5 with $\theta = \frac{l-s}{l-m}$,

$$\|v_h\|_{l,\Omega} = \|\boldsymbol{j}(v_h)\|_{l,\Omega} \leq C\frac{1}{h_K^{l-s}}\|v_h\|_{s,\Omega}$$

□

# Appendix B

# Fluids and solids

## B.1 Introduction

In this appendix we rewiev some preliminary notion on solid and fluid equations. The main goal is to establish the notations and to review important facts about the physical and mathematical aspects of the problem.

## B.2 Kinematics

In this section we rewiev the main ideas related to the kinematic of a continuum body useful for the description of solids and fluid motion, we refer to [95, 69].

### B.2.1 Notations

We consider a body, initially in a location denoted by $\Omega$, undergoing a time-dependent motion $\boldsymbol{\phi} : \overline{\Omega} \times [0,T] \to \mathbb{R}^d$ with $(d = 2,3)$, that describes its trajectory through the ambient space $\mathbb{R}^d$. The set $\Omega$ is called the *reference configuration* and can be thought of as consisting of points $\boldsymbol{s}$ that serve as labels for the material points. For this reason the coordinates $\boldsymbol{s}$ are often called reference or material or Lagrangian coordinates. If we fix the time argument of $\boldsymbol{\phi}$, we obtain a configuration mapping $\boldsymbol{\phi}_t = \boldsymbol{\phi}(\cdot, t)$ which gives us the location of the body, at time $t$, namely $\Omega(t) \overset{def}{=} \boldsymbol{\phi}_t(\Omega)$; we call $\Omega(t)$ *spatial configuration* and the coordinates of points in $\Omega(t)$, denoted by $\boldsymbol{x}$, are called *Eulerian coordinates*. We assume that the boundary $\partial\Omega$ of $\Omega$ can be decomposed into subsets $\Gamma_\sigma$ and $\Gamma_u$, obeying the following restrictions

$$\Gamma_u \cup \Gamma_\sigma = \partial\Omega, \qquad \Gamma_u \cap \Gamma_\sigma = \emptyset. \tag{B.1}$$

Traction boundary conditions (Neumann type B.C.) will be imposed on $\Gamma_\sigma$, and displacement boundary conditions (Dirichlet type B.C.) will be imposed on $\Gamma_u$.

### B.2.2 Lagrangian and Eulerian Descriptions

Equations of continuum mechanics can be written using reference coordinates $\boldsymbol{s}$ or spatial coordinates $\boldsymbol{x}$; the choice is highly dependent on the physical system to be studied. Gen-erally, for a fluid flowing through a fixed physical region the identification of individual particle trajectories is not of primary interest; very often quantities as pressure, velocity, temperature etc. at particular locations in the flow field are more desirable. In such problems it is most appropriate to express equations in the current configuration, namely in the Eulerian coordinates $\boldsymbol{x}$. A system described in this way uses the Eulerian description. In most solid mechanics applications, by contrast, the relative position of material particles is of central interest. For example, the elastic response of materials is related to the change of distance between material points. Then the predominant approach to solid mechanics systems is to write all equations in terms of material coordinates or to use the Lagrangian frame of reference.

### B.2.3 Motion

As already specified, the body moves through the space over a period of time $[0,T]$ and occupies the current configuration $\Omega(t)$ at time $t$; a material point that occupies the position $\boldsymbol{s}$ in the reference configuration $\Omega$, will be at the position $\boldsymbol{x} = \boldsymbol{\phi}(\boldsymbol{s},t)$ at time $t$. We require that $\boldsymbol{\phi}$ be **differentiable with continuity, invertible, and orientation preserving**. The **displacement field** from the reference configuration is defined as

$$\boldsymbol{d}(\boldsymbol{s},t) \stackrel{def}{=} \boldsymbol{\phi}(\boldsymbol{s},t) - \boldsymbol{s}. \tag{B.2}$$

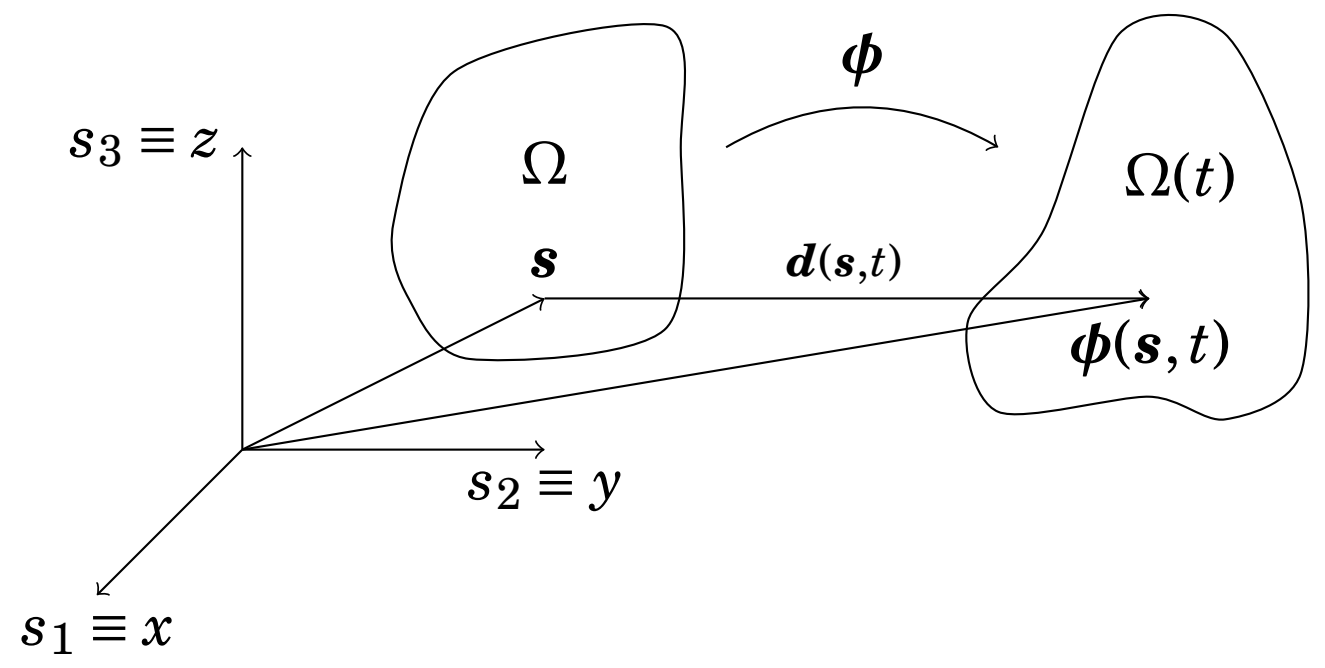

Figure B.1: Body motion

The gradient of $\boldsymbol{\phi}$ is called **deformation gradient**

$$\boldsymbol{F}(\boldsymbol{s},t) \stackrel{def}{=} \nabla_{\boldsymbol{s}}\boldsymbol{\phi}(\boldsymbol{s},t). \tag{B.3}$$

Clearly,

$$\boldsymbol{F}(\boldsymbol{s},t) = \boldsymbol{I} + \nabla_{\boldsymbol{s}}\boldsymbol{d}(\boldsymbol{s},t), \tag{B.4}$$

where $\boldsymbol{I}$ is the identity tensor and $\nabla_{\boldsymbol{s}}\boldsymbol{d}$ is the displacement gradient. The **deformation gradient** is useful for the local description of the motion, in fact, the relative motion of two points can be described using the Taylor expansion of $\boldsymbol{\phi}$; given $\boldsymbol{s}_1$ and $\boldsymbol{s}_2$ in $\Omega$, we can write

$$\boldsymbol{x}_2 - \boldsymbol{x}_1 = \boldsymbol{\phi}(\boldsymbol{s}_2,t) - \boldsymbol{\phi}(\boldsymbol{s}_1,t) = \boldsymbol{F}(\boldsymbol{s}_1,t)(\boldsymbol{s}_2 - \boldsymbol{s}_1) + o(\|\boldsymbol{s}_2 - \boldsymbol{s}_1\|). \tag{B.5}$$

It is possible to relate the deformation gradient $\boldsymbol{F}$ to the change of length of an infinitesimal material segment, indeed, let $\boldsymbol{\delta x} = \boldsymbol{x}_2 - \boldsymbol{x}_1$ be the current configuration of $\boldsymbol{\delta s} = \boldsymbol{s}_2 - \boldsymbol{s}_1$, then, neglecting infinitesimal of higher order,

$$\|\boldsymbol{\delta x}\|^2 = \boldsymbol{\delta x}^T \boldsymbol{\delta x} = \boldsymbol{\delta s}^T \underbrace{\boldsymbol{F}^T \boldsymbol{F}}_{\boldsymbol{C}} \boldsymbol{\delta s}. \tag{B.6}$$

The tensor $\boldsymbol{C} = \boldsymbol{F^T F}$ is called the **right Cauchy-Green deformation tensor** and it gives a possible measure of the deformation state in the body. Related to the Cauchy-Green deformation tensor is the *Green-Lagrange strain tensor* defined as

$$\boldsymbol{E} \stackrel{def}{=} \frac{1}{2}(\boldsymbol{F}^T \boldsymbol{F} - \boldsymbol{I}) = \frac{1}{2}(\nabla \boldsymbol{d} + \nabla \boldsymbol{d}^T + \nabla \boldsymbol{d}^T \nabla \boldsymbol{d}). \tag{B.7}$$

It represents the change in the length of an infinitesimal fiber $\delta \boldsymbol{s}$, namely

$$\frac{\|\boldsymbol{\delta x}\|^2 - \|\boldsymbol{\delta s}\|^2}{2} = \boldsymbol{\delta s}^T \boldsymbol{E} \boldsymbol{\delta s}.$$

## B.2.4 Rates of Motion

If $\boldsymbol{\phi}(\boldsymbol{s}, t)$ is the position of $\boldsymbol{s}$ at time $t$, i.e.,

$$\boldsymbol{x} = \boldsymbol{\phi}(\boldsymbol{s}, t),$$

then

$$\dot{\boldsymbol{x}} \stackrel{def}{=} \partial_t \boldsymbol{\phi}(\boldsymbol{s}, t),$$

is the velocity and

$$\ddot{\boldsymbol{x}} \stackrel{def}{=} \partial_{tt} \boldsymbol{\phi}(\boldsymbol{s}, t),$$

is the acceleration of the material point $\boldsymbol{s}$ at time $t$. Since $\boldsymbol{\phi}$ is bijective, we can also describe the velocity as a function of the position $\boldsymbol{x} \in \mathbb{R}^d$ and time $t$:

$$\boldsymbol{v} = \boldsymbol{v}(\boldsymbol{x}, t) \stackrel{def}{=} \dot{\boldsymbol{x}}(\boldsymbol{\phi}^{-1}(\boldsymbol{x}, t), t). \tag{B.8}$$

This is the **spatial description of the velocity**.

## B.2.5 Material derivative

There are differences in the way rates of change appear in the Lagrangian and Eulerian formulations. Let $\psi(\boldsymbol{s}, t)$ be a function in the Lagrangian description; its total time derivative $\frac{d}{dt}\psi$ coincides with its partial time derivative, because the material coordinate $\boldsymbol{s}$ is not time dependent

$$\frac{d}{dt}\psi(\boldsymbol{s}, t) = \partial_t \psi(\boldsymbol{s}, t).$$

On the contrary, if $\psi(\boldsymbol{x},t)$ is a function in the Eulerian description, the total time derivative of the function evaluated on a material particle $\boldsymbol{s}$ moving with trajectory $\boldsymbol{x} = \boldsymbol{\phi}(\boldsymbol{s},t)$, is expressed using the so called **material derivative**

$$\begin{aligned}\frac{d}{dt}\psi(\boldsymbol{x},t)\big|_{\boldsymbol{x}=\boldsymbol{\phi}(\boldsymbol{s},t)} &= \partial_t\psi(\boldsymbol{x},t)\big|_{\boldsymbol{x}=\boldsymbol{\phi}(\boldsymbol{s},t)} + \nabla_{\boldsymbol{x}}\psi(\boldsymbol{x},t)\big|_{\boldsymbol{x}=\boldsymbol{\phi}(\boldsymbol{s},t)}\cdot\partial_t\boldsymbol{\phi}(\boldsymbol{s},t)\big|_{\boldsymbol{s}=\boldsymbol{\phi}^{-1}(\boldsymbol{x},t)} \\ &= \big(\partial_t\psi(\boldsymbol{x},t) + \boldsymbol{v}(\boldsymbol{x},t)\cdot\nabla_{\boldsymbol{x}}\psi(\boldsymbol{x},t)\big)\Big|_{\boldsymbol{x}=\boldsymbol{\phi}(\boldsymbol{s},t)}.\end{aligned}$$

In the literature it is common to indicate the **material derivative operator** of an Eulerian function using the notation

$$\frac{D}{Dt} = \partial_t + \boldsymbol{v}\cdot\nabla_{\boldsymbol{x}}. \tag{B.9}$$

Thus, in the Eulerian formulation, the acceleration is the "material time derivative" of the Eulerian velocity $\boldsymbol{v}$,

$$\boldsymbol{a} = \frac{D}{Dt}\boldsymbol{v} = \partial_t\boldsymbol{v} + (\boldsymbol{v}\cdot\nabla_{\boldsymbol{x}})\boldsymbol{v}. \tag{B.10}$$

### B.2.6 Rates of Deformation

Given the Eulerian velocity field $\boldsymbol{v}(x,t)$, we introduce the spatial (Eulerian) field called **velocity gradient**

$$\boldsymbol{L}(\boldsymbol{x},t) \stackrel{def}{=} \nabla_{\boldsymbol{x}}\boldsymbol{v}(\boldsymbol{x},t). \tag{B.11}$$

Using the notation $\hat{\boldsymbol{L}} \stackrel{def}{=} \nabla_{\boldsymbol{x}}\boldsymbol{v}(\boldsymbol{x},t)\big|_{\boldsymbol{x}=\boldsymbol{\phi}(\boldsymbol{s},t)}$, the velocity gradient is related to the time derivative of the deformation gradient as follows

$$\dot{\boldsymbol{F}} = \partial_t\nabla_{\boldsymbol{s}}\boldsymbol{\phi}(\boldsymbol{s},t) = \nabla_{\boldsymbol{s}}\partial_t\boldsymbol{\phi}(\boldsymbol{s},t) = \nabla_{\boldsymbol{s}}\boldsymbol{v}(\boldsymbol{x},t)\big|_{\boldsymbol{x}=\boldsymbol{\phi}(\boldsymbol{s},t)} = \nabla_{\boldsymbol{x}}\boldsymbol{v}(\boldsymbol{x},t)\big|_{\boldsymbol{x}=\boldsymbol{\phi}(\boldsymbol{s},t)}\nabla_{\boldsymbol{s}}\boldsymbol{\phi}(\boldsymbol{s},t) = \hat{\boldsymbol{L}}\boldsymbol{F}. \tag{B.12}$$

It is standard practice to write $\boldsymbol{L}$ in terms of its symmetric and skew-symmetric parts

$$\boldsymbol{L} = \boldsymbol{\epsilon}(\boldsymbol{v}) + \boldsymbol{W}. \tag{B.13}$$

where

$$\boldsymbol{\epsilon}(\boldsymbol{v}) = \frac{1}{2}(\nabla_{\boldsymbol{x}}\boldsymbol{v} + \nabla_{\boldsymbol{x}}\boldsymbol{v}^T) \qquad \boldsymbol{W} = \frac{1}{2}(\nabla_{\boldsymbol{x}}\boldsymbol{v} - \nabla_{\boldsymbol{x}}\boldsymbol{v}^T), \tag{B.14}$$

are the **deformation rate tensor** and the **spin tensor** respectively. The spin tensor is related to the rigid rotation of the continuum, in fact it can be proved that if $\boldsymbol{v}$ is the velocity field,

$$\boldsymbol{W}\boldsymbol{v} = \frac{\mathbf{1}}{\mathbf{2}}\boldsymbol{\omega}\wedge\boldsymbol{v} \tag{B.15}$$

where $\boldsymbol{\omega}$ is the vorticity

$$\boldsymbol{\omega} = \mathbf{curl}\,\boldsymbol{v}. \tag{B.16}$$

On the other hand the deformation rate tensor $\boldsymbol{\epsilon}(\boldsymbol{v})$ is responsible of the changes of shape of material volumes.

### B.2.7 Piola Trasformation

Suppose that a subdomain $\omega$ of the reference configuration $\Omega$ of a body, with boundary $\partial\omega$ and unit exterior vector $\hat{\boldsymbol{n}}$ normal to the surface-area element $d\hat{\boldsymbol{\sigma}}$, is mapped by the motion into a subdomain $\omega(t)$ of the current configuration with boundary $\partial\omega(t)$ with unit exterior vector $\boldsymbol{n}$ normal to the "deformed" surface area $d\boldsymbol{\sigma}$. Let $T = T(\boldsymbol{x})$ denote a tensor field defined on $\omega(t)$ and $T(\boldsymbol{x})\boldsymbol{n}(\boldsymbol{x})$ the flux of $T$ across $\partial\omega(t)$. Corresponding to $T$, a tensor field $\hat{T} = \hat{T}(\boldsymbol{s})$ is defined on $\omega$ and we can consider its flux across $\partial\omega$, $\hat{T}(\boldsymbol{s})\hat{\boldsymbol{n}}(\boldsymbol{s})$. The **Piola Transformation** is a relationship between $\hat{T}(\boldsymbol{s})$ and $T(\boldsymbol{x})\big|_{\boldsymbol{x}=\boldsymbol{\phi}_t(\boldsymbol{s})}$, that preserves the total flux through the surfaces $\partial\omega$ and $\partial\omega(t)$, that is

$$\int_{\partial\omega} \hat{T}\hat{\boldsymbol{n}}d\hat{\boldsymbol{\sigma}} = \int_{\partial\omega(t)} T\boldsymbol{n}d\boldsymbol{\sigma}. \tag{B.17}$$

We have the following proposition.

**Proposition B.2.1.** *The property expressed by* (B.17) *is true if*

$$\hat{T}(\boldsymbol{s}) = \det(\boldsymbol{F}(\boldsymbol{s}))\,T(\boldsymbol{\phi}(\boldsymbol{s}))\boldsymbol{F}^{-T}(\boldsymbol{s}). \tag{B.18}$$

*Proof.* See [95, Proposition 1.1]. □

The Piola transform has also other properties other than (B.17), in particular it provides fundamental relationships between differential surface areas and their orientations in the reference and current configurations; we will list them without proof (see [95]).

1. Since $\omega$ can be arbitrary we have the following symbolic relation

$$\hat{T}\hat{\boldsymbol{n}}d\hat{\boldsymbol{\sigma}} = T\boldsymbol{n}d\boldsymbol{\sigma}. \tag{B.19}$$

2. If $T(\boldsymbol{x},t) = I$ then we have the following relation between normals

$$\det(\boldsymbol{F})\boldsymbol{F}^{-T}\hat{\boldsymbol{n}}d\hat{\boldsymbol{\sigma}} = \boldsymbol{n}d\boldsymbol{\sigma}. \tag{B.20}$$

3. From the previous point, since the euclidean norm $\|\boldsymbol{n}\|_{\mathbb{R}^d} = 1$, it comes the Nanson's formula

$$d\boldsymbol{\sigma} = \det(\boldsymbol{F})\,\|\boldsymbol{F}^{-T}\hat{\boldsymbol{n}}\|_{\mathbb{R}^d}d\hat{\boldsymbol{\sigma}} \tag{B.21}$$

## B.3 Fundamental balance laws

We will establish the conservation laws in the Eulerian framework, then let $\Omega \subset \mathbb{R}^d$ be a region in two or three dimensional space filled with a continuum mean. We call $\boldsymbol{v}$ the (spatial) velocity field. For each time $t$, we assume that the continuum has a well-defined positive mass density $\rho(\boldsymbol{x},t)$. In order to derive the equations describing the continuum motion, we assume that the velocity field and the density field are smooth enough to perform all the required calculus operations. Thus, if $\omega$ is any subregion of $\Omega$, the mass contained in $\omega$ at time $t$ is given by

$$m(\omega,t) = \int_{\omega} \rho(\boldsymbol{x},t)d\boldsymbol{x}. \tag{B.22}$$

The derivation of the conservation laws is based on three basic principles:

1. mass is neither created nor destroyed ;
2. the rate of change of momentum of a portion of the continuum equals the force applied to it (Newton's second law);
3. energy is neither created nor destroyed.

### B.3.1 Conservation of Mass

Let $\omega$ be a **fixed** subregion of $\Omega$; the rate of change of mass in $\omega$ is given by

$$\frac{d}{dt}m(\omega,t)=\frac{d}{dt}\int_{\omega}\rho(\boldsymbol{x},t)d\boldsymbol{x}=\int_{\omega}\partial_t\rho(\boldsymbol{x},t)d\boldsymbol{x}. \tag{B.23}$$

where the symbol $\partial_t$ denote the partial derivative with respect to time (see [95, 69]). Let $\partial\omega$ denote the boundary of $\omega$ , assumed to be smooth; let $\boldsymbol{n}$ denote the unit outward normal defined at points of $\partial\omega$ ; and let $d\sigma$ denote the area element on $\partial\omega$ . The volume flow rate across $\partial\omega$ per unit area is $\boldsymbol{v}\cdot\boldsymbol{n}$ and the mass flow rate per unit area is $\rho\boldsymbol{v}\cdot\boldsymbol{n}$. The principle of conservation of mass can be more precisely stated as follows:

**the rate of increase of mass in $\omega$ equals the rate at which mass is crossing $\partial\omega$ in the inward direction**;

$$\frac{d}{dt}\int_{\omega}\rho d\boldsymbol{x}=-\int_{\partial\omega}\rho\boldsymbol{v}\cdot\boldsymbol{n}d\boldsymbol{\sigma}. \tag{B.24}$$

This is the **integral form of the law of conservation of mass**. Applying the divergence theorem to the right hand side and using the definition of the rate of change of mass in $\omega$ given in (B.23), this statement is equivalent to

$$\int_{\omega}\partial_t\rho+\mathrm{div}(\rho\boldsymbol{v})d\boldsymbol{x}=0. \tag{B.25}$$

If $\rho$ and $\boldsymbol{v}$ are smooth enough to apply the mean value theorem, by the arbitrariness of $\omega$, we can derive the **differential form of the law of conservation of mass**, also known as the **continuity equation**

$$\partial_t\rho+\mathrm{div}(\rho\boldsymbol{v})=0 \tag{B.26}$$

### B.3.2 Transport theorem

Let $\Omega$ denote the region in which the fluid is moving in time. In many situations we are led to compute the total time derivative of integrals of a field $f(\boldsymbol{x},t)$ defined over a volume depending on time $\omega(t)\subset\Omega$

$$\frac{d}{dt}\int_{\omega(t)}f(\boldsymbol{x},t)d\boldsymbol{x}.$$

In this case we have to consider that the domain of integration is changing in time. Let $\omega\subset\Omega$ be the material volume corresponding to $\omega(t)$, namely $\omega(t)=\boldsymbol{\phi}_t(\omega)$. Let $J(\boldsymbol{s},t)=\det(\boldsymbol{F}(\boldsymbol{s},t))$ the Jacobian determinant of the map $\boldsymbol{\phi}_t$, then it holds

**Lemma B.3.1.**

$$\partial_t J(\boldsymbol{s},t)\Big|_{\boldsymbol{s}=\boldsymbol{\phi}_t^{-1}(\boldsymbol{x})} = J(\boldsymbol{s},t)\Big|_{\boldsymbol{s}=\boldsymbol{\phi}_t^{-1}(\boldsymbol{x})}\left[\operatorname{div}\boldsymbol{v}(\boldsymbol{\phi}(\boldsymbol{s},t),t)\right]. \tag{B.27}$$

The proof of this lemma can be found in [69, 95]. The important consequence of this lemma is the following theorem.

**Theorem B.3.2 (Transport Theorem).** *For any sufficiently regular function $f(\boldsymbol{x},t)$ defined in the Eulerian framework, denoting by $\rho$ the density of the continuum body, we have*

$$\frac{d}{dt}\int_{\omega(t)} \rho f(\boldsymbol{x},t)d\boldsymbol{x} = \int_{\omega(t)} \rho\frac{D}{Dt}f(\boldsymbol{x},t)d\boldsymbol{x} \tag{B.28}$$

*Proof.* Let assume that $\omega(t) = \boldsymbol{\phi}_t(\omega)$, then, by a change of variables we can write

$$\frac{d}{dt}\int_{\omega(t)} \rho(\boldsymbol{x},t)f(\boldsymbol{x},t)d\boldsymbol{x} = \frac{d}{dt}\int_{\omega} \rho(\boldsymbol{\phi}_t(\boldsymbol{s}),t)f(\boldsymbol{\phi}_t(\boldsymbol{s}),t)J(\boldsymbol{s},t)d\boldsymbol{s}. \tag{B.29}$$

Because the domain $\omega$ is fixed and we assume regularity for the integrand, we may differentiate under the integral sign. Note that

$$\frac{d}{dt}\rho(\boldsymbol{\phi}_t(\boldsymbol{s}),t)f(\boldsymbol{\phi}_t(\boldsymbol{s}),t)J(\boldsymbol{s},t) = \frac{D}{Dt}\left(\rho f\right)J(\boldsymbol{s},t) + \left(\rho f\right)\partial_t J(\boldsymbol{s},t). \tag{B.30}$$

Applying the previous lemma B.3.1 and the conservation of mass (B.26), we obtain

$$\begin{aligned}\int_\omega \frac{D}{Dt}\left(\rho f\right)J(\boldsymbol{s},t) + \left(\rho f\right)\partial_t J(\boldsymbol{s},t)d\boldsymbol{s} &= \int_\omega \left[\frac{D}{Dt}\left(\rho f\right) + \left(\rho f\right)\operatorname{div}\boldsymbol{v}\right]J(\boldsymbol{s},t)d\boldsymbol{s} \\ &= \int_{\omega(t)} f\underbrace{\left(\frac{D}{Dt}\rho + \rho\operatorname{div}\boldsymbol{v}\right)}_{\partial_t\rho+\operatorname{div}(\rho\boldsymbol{v})=0} + \rho\frac{D}{Dt}f d\boldsymbol{x} = \int_{\omega(t)} \rho\frac{D}{Dt}f d\boldsymbol{x}.\end{aligned} \tag{B.31}$$

□

### B.3.3 Forces

For any continuum, forces acting on a material region $\omega \subset \Omega$ are of two types (see [95]). First, there are external, or body, forces such as gravity or a magnetic field, which exert a force per unit volume on the continuum. We suppose the existence of a mass density of body forces $\boldsymbol{b}(\boldsymbol{x},t)$ defined on $\omega(t)$. Hence the resultant of body forces acting on $\omega(t)$ is given by

$$B(\omega(t)) = \int_{\omega(t)} \rho\boldsymbol{b}d\boldsymbol{x}. \tag{B.32}$$

Second, there are forces of stress acting on $\omega(t)$ across its boundary $\partial\omega(t)$. These forces are given by the rest of the continuum or are applied forces on $\partial\Omega$ (in the case $\partial\omega \cap \partial\Omega(t) \neq \emptyset$). We suppose the existence of a surface density of stress named $\boldsymbol{t}(\boldsymbol{x},t,\boldsymbol{n}_t)$ that depend with continuity from the point $\boldsymbol{x} \in \partial\omega(t)$ and from the unit normal $\boldsymbol{n}_t$ to $\partial\omega(t)$. According to this assumption we can express the resultant of the forces acting on $\omega(t)$ across the boundary $\partial\omega(t)$ as follows

$$F(\omega(t)) = \int_{\partial\omega(t)} \boldsymbol{t} d\boldsymbol{\sigma}. \tag{B.33}$$

It is of fundamental importance the fact that the vector field $\boldsymbol{t}(\boldsymbol{x}, t, \boldsymbol{n})$ depends linearly from $\boldsymbol{n}$, in fact it holds the following Theorem due to Cauchy ([95]).

**Theorem B.3.3.** *There exist a tensor field $T(\boldsymbol{x}, t)$ defined on $\Omega \times \mathbb{R}^+$ such that for any unit normal vector $\boldsymbol{n}$*

$$\boldsymbol{t}(\boldsymbol{x}, t, \boldsymbol{n}) = T(\boldsymbol{x}, t)\boldsymbol{n}. \tag{B.34}$$

$T(\boldsymbol{x}, t)$ as above is called **Cauchy stress tensor**.

### B.3.4 Balance of Momentum

Using the expression of forces (B.32) and (B.33), the **conservation of the linear momentum** can be expressed as follows [95, 69]

$$\frac{d}{dt}\int_{\omega(t)} \rho \boldsymbol{v} d\boldsymbol{x} = \int_{\omega(t)} \rho \boldsymbol{b} d\boldsymbol{x} + \int_{\partial\omega(t)} \boldsymbol{t} d\boldsymbol{\sigma}. \tag{B.35}$$

Equation (B.35) expresses the fact that the rate of change of the linear momentum is equal to the resultant of the external forces acting on material contained in $\omega(t)$. Applying the Cauchy theorem B.3.3 and the Transport Theorem B.3.2 we obtain

$$\int_{\omega(t)} \rho \frac{D}{Dt} \boldsymbol{v} d\boldsymbol{x} = \int_{\omega(t)} \rho \boldsymbol{b} d\boldsymbol{x} + \int_{\partial\omega(t)} T\boldsymbol{n} d\boldsymbol{\sigma} = \int_{\omega(t)} \rho \boldsymbol{b} + \operatorname{div} T d\boldsymbol{x}. \tag{B.36}$$

From the arbitrariness of the choice of $\omega(t)$ we are led to the **differential form of conservation of the linear momentum**.

$$\rho \frac{D}{Dt} \boldsymbol{v} = \rho \boldsymbol{b} + \operatorname{div} T. \tag{B.37}$$

We will assume also the **balance of the angular momentum** that is expressed in the following equation

$$\frac{d}{dt}\int_{\omega(t)} \rho \boldsymbol{v} \wedge \boldsymbol{x} d\boldsymbol{x} = \int_{\omega(t)} \rho \boldsymbol{b} \wedge \boldsymbol{x} d\boldsymbol{x} + \int_{\partial\omega(t)} \boldsymbol{t} \wedge \boldsymbol{x} d\boldsymbol{\sigma}. \tag{B.38}$$

It can be proved that (B.38) is equivalent to the symmetry of the Cauchy stress tensor $T(\boldsymbol{x}, t)$ (see [95, 69]).

### B.3.5 Balance of energy

We will recall here only the balance of the kinetic energy since, in the chapters that follow we are not interested in the thermodynamic aspect of continuum theory. We define the **kinetic energy** of a continuum material region $\omega(t)$ as follows [95, 69]

$$K(\omega(t)) = \frac{1}{2}\int_{\omega(t)} \rho \boldsymbol{v} \cdot \boldsymbol{v} d\boldsymbol{x}, \tag{B.39}$$

where $\boldsymbol{u}$ is the velocity of the continuum. Then it holds the following **theorem of conservation of the kinetic energy**.

**Theorem B.3.4.** *Suppose that the balance laws of mass, linear momentum and angular momentum are satisfied, then*

$$\frac{d}{dt}K(\omega(t)) = \int_{\omega(t)} \rho \boldsymbol{b}\cdot\boldsymbol{v}d\boldsymbol{x} + \int_{\partial\omega(t)} \boldsymbol{t}\cdot\boldsymbol{v}d\boldsymbol{\sigma} - \int_{\omega(t)} T : \boldsymbol{\epsilon}(\boldsymbol{v})d\boldsymbol{x}. \tag{B.40}$$

*where $T$ is the Cauchy stress tensor and $\boldsymbol{\epsilon}(\boldsymbol{v})$ is the tensor of velocity of deformation defined in* (B.14).

## B.4 Viscous fluids and Navier-Stokes equations.

In the following chapters we are interested in the analysis of interaction of solid and viscous incompressible fluids. The incompressibility, in general, is a property of the flow, indeed we say that a **flow is incompressible in** $\Omega$ if all material volumes $\omega \subset \Omega$, preserve their volume during the motion, namely

$$\frac{d}{dt}\int_{\omega(t)} d\boldsymbol{x} = 0. \tag{B.41}$$

Using Lemma B.3.1 we obtain

$$0 = \frac{d}{dt}\int_{\omega(t)} d\boldsymbol{x} = \int_{\omega} \partial_t J(\boldsymbol{s},t)d\boldsymbol{s} = \int_{\omega} \operatorname{div}\boldsymbol{v}Jd\boldsymbol{s} = \int_{\omega(t)} \operatorname{div}\boldsymbol{v}d\boldsymbol{x}. \tag{B.42}$$

Since the last equation is true for all $\omega(t)$, we can state that the **flow is incompressible in the region $\Omega$ occupied by the fluid if and only if** $\operatorname{div}\boldsymbol{u} = 0$ **in** $\Omega$. We can introduce now the *Newtonian fluids* for which the Cauchy stress tensor is expressed as follows

$$T = (-p + \lambda \operatorname{div}\boldsymbol{v})\boldsymbol{I} + 2\mu\boldsymbol{\epsilon}(\boldsymbol{v}),$$

hence, for *incompressible Newtonian fluids* we have

$$T = -p\boldsymbol{I} + 2\mu\boldsymbol{\epsilon}(\boldsymbol{v}). \tag{B.43}$$

Using the constitutive relation (B.43) in equation of conservation of the linear momentum (B.37), we obtain

$$\rho\frac{D}{Dt}\boldsymbol{v} = \rho\boldsymbol{b} - \nabla p + 2\mu\operatorname{div}(\boldsymbol{\epsilon}(\boldsymbol{v})). \tag{B.44}$$

We observe that the stress tensor (B.43) is symmetric, hence it satisfies the conservation of angular momentum; It follows that the equations that describe the motion of an incompressible Newtonian fluid are the **Navier-Stokes equations**

$$\begin{aligned} \rho\frac{D}{Dt}\boldsymbol{v} = \rho\boldsymbol{b} - \nabla p + 2\mu\operatorname{div}(\boldsymbol{\epsilon}(\boldsymbol{v})) \qquad &\text{in} \quad \Omega, \\ \operatorname{div}\boldsymbol{v} = 0 \qquad &\text{in} \quad \Omega, \end{aligned} \tag{B.45}$$

completed with boundary and initial conditions.

## B.5 Solid Equations and Thin structural models

Solid problems, as said, are better treated in reference configuration; the reason is that for solid constitutive laws, generally, it is important to evaluate the changes in the reciprocal positions of the material particles. For this reason we will write the continuum equations in the Lagrangian framework.

### B.5.1 The Equations of Motion in the Lagrangian framework: The Piola-Kirchhoff Stress Tensors

The equation of conservation of the linear momentum (B.37) is written in the current configuration, hence we apply a change of framework in order to write the equation (B.37) in the reference configuration. We introduce the following quantities

$$\hat{\boldsymbol{b}}(\boldsymbol{s},t) = \boldsymbol{b}(\boldsymbol{\phi}(\boldsymbol{s},t),t)\det(\boldsymbol{F}), \tag{B.46}$$

$$\hat{\rho}(\boldsymbol{s},t) = \rho(\boldsymbol{\phi}(\boldsymbol{s},t),t)\det(\boldsymbol{F}), \tag{B.47}$$

$$\boldsymbol{P}(\boldsymbol{s},t) = \det(\boldsymbol{F})\,T(\boldsymbol{\phi}(\boldsymbol{s},t),t)\boldsymbol{F}^{-T}. \tag{B.48}$$

The tensor $\boldsymbol{P}(\boldsymbol{s},t)$ is called the **first Piola-Kirchhoff Stress Tensor**. We observe that, using the repeated indexes convention,

$$\mathrm{Div}\left(\boldsymbol{F}^{-\mathrm{T}}\det(\boldsymbol{F})\right) \stackrel{\text{def}}{=} \frac{\partial}{\partial \boldsymbol{s}_{\mathrm{j}}}\left(\boldsymbol{F}^{-\mathrm{T}}\det(\boldsymbol{F})\right)_{\mathrm{kj}} = 0$$

this fact can be proved by direct computation; we refer for proof to [95, Proposition 1.1]. Then we have

$$\mathrm{Div}\boldsymbol{P} \stackrel{\text{def}}{=} \frac{\partial \boldsymbol{P}_{\mathrm{ij}}}{\partial \boldsymbol{s}_{\mathrm{j}}} = \frac{\partial \mathrm{T}_{\mathrm{ik}}}{\partial \boldsymbol{x}_{\mathrm{k}}}\frac{\partial \boldsymbol{x}_{\mathrm{k}}}{\partial \boldsymbol{s}_{\mathrm{j}}}\left(\boldsymbol{F}^{-\mathrm{T}}\det\boldsymbol{F}\right)_{\mathrm{kj}} + \mathrm{T}_{\mathrm{ik}}\underbrace{\frac{\partial}{\partial \boldsymbol{s}_{j}}\left(\boldsymbol{F}^{-T}\det\boldsymbol{F}\right)_{kj}}_{=0} = \det\boldsymbol{F}\,\mathrm{div}\mathrm{T} \tag{B.49}$$

Hence, denoting by $\boldsymbol{v}$ the Eulerian velocity of the solid and by $\boldsymbol{d}$ the solid displacement in Lagrangian coordinate, from the local form of conservation of linear and angular momentum, referring to (B.49), we obtain

$$\mathrm{div}\,T + \boldsymbol{b} = \rho\frac{D}{Dt}\boldsymbol{v} \quad \Rightarrow \quad \frac{1}{\det\boldsymbol{F}}\left(\mathrm{Div}\boldsymbol{P} + \hat{\boldsymbol{b}}\right) = \frac{\hat{\rho}}{\det\boldsymbol{F}}\partial_{tt}\boldsymbol{d}, \tag{B.50}$$

$$T = T^{T} \quad \Rightarrow \quad \frac{1}{\det\boldsymbol{F}}\boldsymbol{P}\boldsymbol{F}^{T} = \frac{1}{\det\boldsymbol{F}}\boldsymbol{F}\boldsymbol{P}^{T}. \tag{B.51}$$

Note that $\boldsymbol{P}$ is not symmetric; however, $\boldsymbol{P}\boldsymbol{F}^{T} = \boldsymbol{F}\boldsymbol{P}^{T}$ since $T$ is symmetric. Thus, the equations of motion (linear and angular momentum) in the reference configuration are

$$\begin{aligned} \mathrm{Div}\boldsymbol{P} + \hat{\boldsymbol{b}} &= \hat{\rho}\partial_{\mathrm{tt}}\boldsymbol{d}, \\ \boldsymbol{P}\boldsymbol{F}^{T} &= \boldsymbol{F}\boldsymbol{P}^{T}. \end{aligned} \tag{B.52}$$

Introducing the **second Piola-Kirchhoff Stress Tensor** $\boldsymbol{S} = \boldsymbol{F}^{-1}\boldsymbol{P}$, we obtain a symmetric stress tensor and the equilibrium equations become

$$\mathrm{Div}\boldsymbol{F}\boldsymbol{S} + \hat{\boldsymbol{b}} = \hat{\rho}\partial_{\mathrm{tt}}\boldsymbol{d}, \tag{B.53}$$

$$\boldsymbol{S} = \boldsymbol{S}^{T} \tag{B.54}$$

### B.5.2 Hyperelastic materials

The constitutive law for a hyperelastic material is defined assigning a stress-strain relationship that derives from a strain energy density function. Often the constitutive law is better expressed in terms of the *second Piola-Kirchhoff stress tensor*, defined before as the symmetric tensor $\boldsymbol{S} = \boldsymbol{F}^{-1}\boldsymbol{P}$, and the *Green-Lagrange strain tensor* given in (B.7) and reported here

$$\boldsymbol{E} \stackrel{def}{=} \frac{1}{2}(\boldsymbol{F}^T\boldsymbol{F} - \boldsymbol{I}) = \frac{1}{2}(\nabla \boldsymbol{d} + \nabla \boldsymbol{d}^T + \nabla \boldsymbol{d} \nabla \boldsymbol{d}^T). \tag{B.55}$$

Assigning the elastic potential $W(\boldsymbol{E}) \geq 0$ we can obtain each component of the Second Piola-Kirchhoff stress tensor deriving the elastic potential with respect to the corresponding component of the strain tensor $\boldsymbol{E}$ [95, 69],

$$S_{ij} = \frac{\partial W}{\partial E_{ij}}. \tag{B.56}$$

Using the relation between $\boldsymbol{S}$ and $\boldsymbol{P}$, we derive also

$$P_{ij} = F_{ik} S_{kj} = F_{ik} \frac{\partial W}{\partial E_{kj}}. \tag{B.57}$$

The simplest hyper-elastic constitutive law for an isotropic, homogeneous material is the *Saint Venant-Kirchhoff model* which is just an extension of the linear elastic material model to the nonlinear regime. In practice the *Saint Venant-Kirchhoff elastic potential* is

$$W(\boldsymbol{E}) = \frac{L_1}{2}(tr(\boldsymbol{E}))^2 + L_2 tr(\boldsymbol{E}^2), \tag{B.58}$$

where $L_1 > 0$ and $L_2 > 0$ are the *Lamé constants* of the material. Usually the elastic properties of a material are given in terms of its *Young modulus* $E$ and *Poisson ratio* $\nu$ that are simple to be measured during experimental tests on materials. Young modulus $E$ represent the ratio of the stress (force per unit area) along an axis to the strain (ratio of deformation over initial length) along the same axis; the Poisson ratio $\nu$ is the negative ratio of transverse to axial strain along an axis and for isotropic materials it does not depend on the axis chosen for his measurement. The relation between Lamé constants and Young modulus and Poisson ratio is

$$L_1 = \frac{E\nu}{(1+\nu)(1-2\nu)}, \qquad L_2 = \frac{E}{2(1+\nu)}. \tag{B.59}$$

Then the stress tensor $\boldsymbol{S}$ is

$$\boldsymbol{S}(\boldsymbol{E}) = L_1 tr(\boldsymbol{E})\boldsymbol{I} + 2L_2\boldsymbol{E}. \tag{B.60}$$

### B.5.3 Thin-walled structures

The main results of this work concern the behavior of coupled systems composed by fluid and thin walled structure; in order to describe the dynamic of the structural part we will use structural models for solids of co-dimension one (shell, plate or membrane models). The dynamic of such solids is described tracking the position of their *mid-surface* $\Sigma(t)$, which is parametrized using a **smooth, injective and orientation preserving** *configuration map*

$$\boldsymbol{\phi} : \Sigma \times [0,T] \to \mathbb{R}^d, \qquad \boldsymbol{\phi}(\boldsymbol{s},t) = \boldsymbol{x}, \tag{B.61}$$

where $\Sigma \subset \mathbb{R}^{d-1}$ represents the reference configuration of the mid-surface of the thin solid that, in general, can be different from the initial configuration. Related to the initial configuration of the mid-surface $\boldsymbol{\phi}_0 = \boldsymbol{\phi}(s_1, s_2, 0)$ is the *local covariant basis*, whose vectors are given by $(\boldsymbol{a}_1 = \partial_1 \boldsymbol{\phi}_0, \boldsymbol{a}_2 = \partial_2 \boldsymbol{\phi}_0)$, where $\partial_1 \overset{def}{=} \frac{\partial}{\partial s_1}$ and $\partial_2 \overset{def}{=} \frac{\partial}{\partial s_2}$ are the derivative w.r.t. to $s_1$ and $s_2$ respectively. The covariant base vectors $(\boldsymbol{a}_1(s_1,s_2), \boldsymbol{a}_2(s_1,s_2))$ lay in the tangential plane which touches the surface at the point with the position vector $\boldsymbol{\phi}_0(s_1,s_2)$. We introduce, moreover, a normal vector to the couple $(\boldsymbol{a}_1, \boldsymbol{a}_2)$, which is denoted $\boldsymbol{a}_3$ and is given by

$$\boldsymbol{a}_3 \overset{def}{=} \frac{\boldsymbol{a}_1 \wedge \boldsymbol{a}_2}{\|\boldsymbol{a}_1 \wedge \boldsymbol{a}_2\|_{\mathbb{R}^d}}. \tag{B.62}$$

In the following we will use also the *contravariant base vectors* denoted as $(\boldsymbol{a}^1, \boldsymbol{a}^2)$ lying in the tangent plane and defined by the following equations

$$\boldsymbol{a}^\alpha \cdot \boldsymbol{a}_\beta = \delta^\alpha_\beta, \qquad \alpha, \beta = 1,2. \tag{B.63}$$

In the derivation of the shell equations the following tensors associated to the covariant and contravariant base vectors are also very important

- *first fundamental form* $\underline{\underline{\boldsymbol{a}}}$ whose covariant and contravariant components are respectively

$$\text{for} \quad \alpha,\beta = 1,2 \qquad a_{\alpha\beta} \overset{def}{=} \boldsymbol{a}_\alpha \cdot \boldsymbol{a}_\beta, \qquad a^{\alpha\beta} \overset{def}{=} \boldsymbol{a}^\alpha \cdot \boldsymbol{a}^\beta, \tag{B.64}$$

- the *second fundamental form* $\underline{\underline{\boldsymbol{b}}}$ whose covariant and contravariant components are, respectively

$$\text{for} \quad \alpha,\beta = 1,2 \qquad b_{\alpha\beta} \overset{def}{=} \boldsymbol{a}_3 \cdot \partial_\beta \boldsymbol{a}_\alpha, \qquad b^\alpha_\beta \overset{def}{=} \boldsymbol{a}_3 \cdot \partial_\beta \boldsymbol{a}^\alpha, \tag{B.65}$$

  since $\boldsymbol{a}_\alpha \cdot \boldsymbol{a}_3 = \boldsymbol{a}^\alpha \cdot \boldsymbol{a}_3 = 0$ for $\alpha = 1,2$, we have the following relations

$$\text{for} \quad \alpha,\beta = 1,2 \qquad b_{\alpha\beta} = \boldsymbol{a}_3 \cdot \partial_\beta \boldsymbol{a}_\alpha = -\boldsymbol{a}_\alpha \cdot \partial_\beta \boldsymbol{a}_3, \qquad b^\alpha_\beta = \boldsymbol{a}_3 \cdot \partial_\beta \boldsymbol{a}^\alpha = -\boldsymbol{a}^\alpha \cdot \partial_\beta \boldsymbol{a}_3. \tag{B.66}$$

**Remark 5.** *We observe that the two fundamental tensors are symmetric, in fact this is evident for* $\underline{\underline{\boldsymbol{a}}}$*. The symmetry of* $\underline{\underline{\boldsymbol{b}}}$ *derives from the following observations*

$$\partial_\beta \boldsymbol{a}_\alpha = \partial_{\beta\alpha} \boldsymbol{\phi}_0 = \partial_{\alpha\beta} \boldsymbol{\phi}_0 = \partial_\alpha \boldsymbol{a}_\beta$$

$$\partial_\beta \boldsymbol{a}^\alpha = \partial_\beta (\boldsymbol{a}^\alpha \cdot \boldsymbol{a}_\gamma) \boldsymbol{a}_\gamma = \partial_\beta \boldsymbol{a}_\alpha = \partial_\alpha \boldsymbol{a}_\beta = \partial_\alpha (\boldsymbol{a}_\beta \cdot \boldsymbol{a}^\gamma) \boldsymbol{a}^\gamma = \partial_\alpha \boldsymbol{a}^\beta.$$

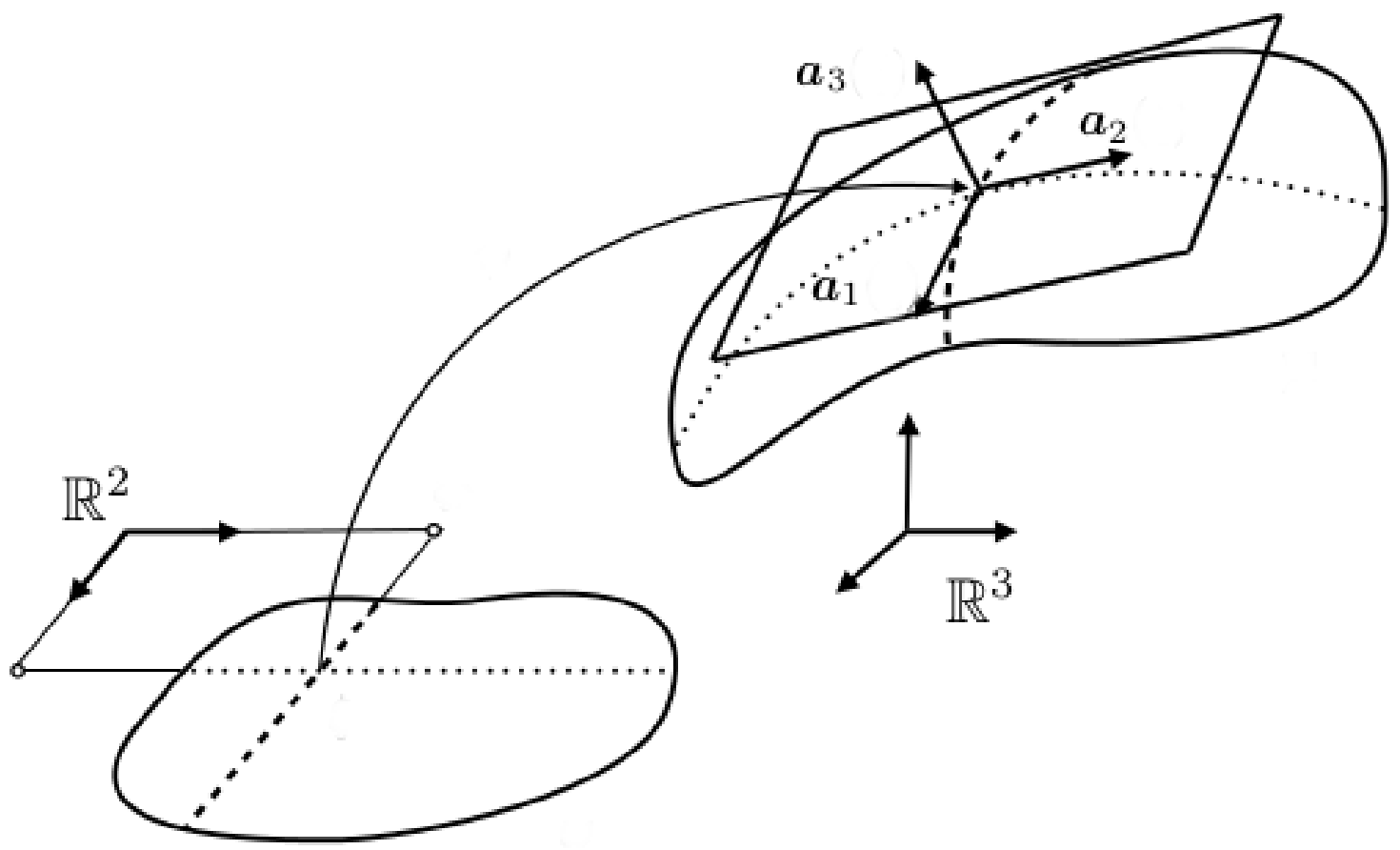

Figure B.2: Mid-surface geometry

Finally, it is important to consider also reference configuration of the thin solid itself $\mathscr{B} \stackrel{def}{=} \Sigma \times [-\frac{\epsilon}{2}, \frac{\epsilon}{2}] \subset \mathbb{R}^d$ that represent the region in $\mathbb{R}^d$ occupied by the solid when its mid-surface coincides with $\Sigma$. In relation to $\mathscr{B}$ we introduce the position vector of each point in the solid,

$$\boldsymbol{R}(s_1, s_2, s_3, t) = \boldsymbol{\phi}(s_1, s_2, t) + s_3 \boldsymbol{\Lambda}(s_1, s_2, t) \boldsymbol{a}_3(s_1, s_2), \tag{B.67}$$

where $\boldsymbol{\Lambda}$ is a rotation tensor. The initial position of the solid is then given by the map

$$\boldsymbol{R}_0 = \boldsymbol{R}(s_1, s_2, s_3, 0) = \boldsymbol{\phi}_0(s_1, s_2) + s_3 \boldsymbol{a}_3(s_1, s_2). \tag{B.68}$$

It is also useful to introduce the covariant and contravariant base vectors associated to $\boldsymbol{R}_0$

- covariant base vectors

$$\boldsymbol{g}_i(s_1, s_2, s_3) \stackrel{def}{=} \frac{\partial \boldsymbol{R}_0}{\partial s_i}(s_1, s_2, s_3), \qquad \text{for} \quad i = 1,2,3; \tag{B.69}$$

- contravariant base vectors for $i, j = 1,2,3$

$$\boldsymbol{g}^i \in \mathbb{R}^d \quad \text{such that} \quad \boldsymbol{g}^i \cdot \boldsymbol{g}_j = \delta^i_j. \tag{B.70}$$

We observe that for $s_3 = 0$ we have $\boldsymbol{g}_i = \boldsymbol{a}_i$ for $i = 1,2,3$; if $s_3 \neq 0$ we have the following relations between $\boldsymbol{g}_i$ and $\boldsymbol{a}_i$

$$\boldsymbol{g}_\alpha = \boldsymbol{a}_\alpha + s_3 \partial_\alpha \boldsymbol{a}_3 = (\delta^\lambda_\alpha - s_3 b^\lambda_\alpha) \boldsymbol{a}_\lambda \quad \text{for} \quad \alpha, \lambda = 1,2; \qquad \boldsymbol{g}_3 = \boldsymbol{a}_3. \tag{B.71}$$

The components of the metric tensor associated to the base vectors $(\boldsymbol{g}_i)_{i=1}^3$ and $(\boldsymbol{g}^i)_{i=1}^3$ are

$$g_{ij} = \boldsymbol{g}_i \cdot \boldsymbol{g}_j, \qquad g^{ij} = \boldsymbol{g}^i \cdot \boldsymbol{g}^j. \tag{B.72}$$

The displacement field in the solid is the difference between the actual position and the initial position and is given by

$$\boldsymbol{D}(s_1,s_2,s_3,t) = \boldsymbol{d}(s_1,s_2,t) + s_3\boldsymbol{\Theta}(s_1,s_2,t)\boldsymbol{a}_3(s_1,s_2), \tag{B.73}$$

where we have used the following notations

$$\boldsymbol{d}(s_1,s_2,t) \stackrel{def}{=} \boldsymbol{\phi}(s_1,s_2,t) - \boldsymbol{\phi}(s_1,s_2,0), \qquad \boldsymbol{\Theta}(s_1,s_2,t) \stackrel{def}{=} \boldsymbol{\Lambda}(s_1,s_2,t) - \boldsymbol{\Lambda}(s_1,s_2,0). \tag{B.74}$$

to denote the displacement field of the mid-surface and the relative rotation tensor, respectively. The hypothesis of small thickness allows for the following approximation

$$s_3\boldsymbol{\Theta}\boldsymbol{a}_3 \approx s_3(\theta_1\boldsymbol{a}^1 + \theta_2\boldsymbol{a}^2), \tag{B.75}$$

in fact, the rotation of an infinitely small straight line orthogonal to the mid-surface is uniquely defined by a rotation vector lying in the tangent plane. We can represent this rotation vector by the contravariant surface tensor $\boldsymbol{\theta} = \theta_1\boldsymbol{a}^1 + \theta_2\boldsymbol{a}^2$. The last assumption is called the *Reissner-Mindlin kinematic assumption* and allows for writing the displacement field as follows

$$\boldsymbol{D}(s_1,s_2,s_3,t) = \boldsymbol{d}(s_1,s_2,t) + s_3(\theta_1\boldsymbol{a}^1 + \theta_2\boldsymbol{a}^2), \tag{B.76}$$

In the following we assume that the the components of the Green-Lagrange strain tensor are small, that is we can consider

$$\boldsymbol{E} \approx \frac{1}{2}(\nabla\boldsymbol{D} + \nabla\boldsymbol{D}^T).$$

Expressing the displacement field in the base $(\boldsymbol{a}^1, \boldsymbol{a}^2, \boldsymbol{a_3})$, we have, using the Einstein convention of sum on repeated index, the following expression of (B.76)

$$\boldsymbol{D} = (d_\alpha + s_3\theta_\alpha)\boldsymbol{a}^\alpha + d_3\boldsymbol{a_3}, \qquad \alpha = 1,2. \tag{B.77}$$

In order to compute the component of the strain tensor in curvilinear coordinate, we perform the derivatives of $D$. We observe that, for all vectors $\boldsymbol{w} \in \mathbb{R}^d$ it holds $\boldsymbol{w}\cdot\boldsymbol{a}^\beta = \boldsymbol{w}\cdot\boldsymbol{a}_\beta$ for $\beta = 1,2$, hence we have the following identity

$$\boldsymbol{w} = (\boldsymbol{w}\cdot\boldsymbol{a}_\beta)\boldsymbol{a}^\beta + (\boldsymbol{w}\cdot\boldsymbol{a}_3)\boldsymbol{a}_3. \tag{B.78}$$

Using (B.78) with $\boldsymbol{w} = \partial_\alpha\boldsymbol{a}^\beta$, we derive the component of the strain tensor

$$\begin{aligned}\partial_\alpha(d_\beta\boldsymbol{a}^\beta) = \partial_\alpha d_\beta\boldsymbol{a}^\beta + d_\beta\partial_\alpha\boldsymbol{a}^\beta &= \partial_\alpha d_\beta\boldsymbol{a}^\beta + d_\beta\left[\left(\partial_\alpha\boldsymbol{a}^\beta\cdot\boldsymbol{a}_\lambda\right)\boldsymbol{a}^\lambda + \left(\partial_\alpha\boldsymbol{a}^\beta\cdot\boldsymbol{a}_3\right)\boldsymbol{a}_3\right]\\ &= \underbrace{\left(\partial_\alpha d_\beta + d_\lambda\partial_\alpha\boldsymbol{a}^\lambda\cdot\boldsymbol{a}_\beta\right)}_{d_{\beta|\alpha}}\boldsymbol{a}^\beta + d_\beta b^\beta_\alpha\boldsymbol{a}_3.\end{aligned} \tag{B.79}$$

The quantity $d_{\beta|\alpha}$ is the covariant surface derivative of the component $d_\beta$ of $\boldsymbol{d}$; then the derivative of the mid-surface displacement, for $\alpha = 1,2$, is

$$\partial_\alpha\boldsymbol{d} = d_{\beta|\alpha}\boldsymbol{a}^\beta + d_\beta b^\beta_\alpha\boldsymbol{a}_3 + \partial_\alpha d_3\boldsymbol{a}_3 + d_3\partial_\alpha\boldsymbol{a}_3 = \left(d_{\beta|\alpha} - b^\beta_\alpha d_3\right)\boldsymbol{a}^\beta + \left(d_\beta b^\beta_\alpha + \partial_\alpha d_3\right)\boldsymbol{a}_3. \tag{B.80}$$

In a similar way we compute also the derivatives of the term $s_3\theta_\beta \boldsymbol{a}^\beta$

$$\partial_\alpha(s_3\theta_\beta \boldsymbol{a}^\beta) = s_3\partial_\alpha\theta_\beta \boldsymbol{a}^\beta + s_3\theta_\beta\partial_\alpha \boldsymbol{a}^\beta = s_3\theta_{\beta|\alpha}\boldsymbol{a}^\beta + s_3\theta_\beta b^\beta_\alpha \boldsymbol{a}_3 \tag{B.81}$$

therefore

$$\begin{aligned}\partial_\alpha \boldsymbol{D} = \partial_\alpha \boldsymbol{d} + \partial_\alpha(s_3\theta_\beta \boldsymbol{a}^\beta) &= \left(d_{\beta|\alpha} - b^\beta_\alpha d_3\right)\boldsymbol{a}^\beta + \left(d_\beta b^\beta_\alpha + \partial_\alpha d_3\right)\boldsymbol{a}_3 + s_3\theta_{\beta|\alpha}\boldsymbol{a}^\beta + s_3\theta_\beta b^\beta_\alpha \boldsymbol{a}_3\\ &= \left(d_{\beta|\alpha} + s_3\theta_{\beta|\alpha} - b^\beta_\alpha d_3\right)\boldsymbol{a}^\beta + \left(d_\beta b^\beta_\alpha + s_3\theta_\beta b^\beta_\alpha + \partial_\alpha d_3\right)\boldsymbol{a}_3. \end{aligned}\tag{B.82}$$

Moreover

$$\partial_3 \boldsymbol{D} = \theta_\alpha \boldsymbol{a}^\alpha. \tag{B.83}$$

Before computing the components of the strain tensor we express the previous derivative in the directions $\boldsymbol{g}_1$, $\boldsymbol{g}_2$ and $\boldsymbol{g}_3$; thanks to the relations (B.71) we have

$$\begin{aligned}\partial_\alpha \boldsymbol{D}\cdot \boldsymbol{g}_\gamma &= \left(d_{\beta|\alpha} + s_3\theta_{\beta|\alpha} - b^\beta_\alpha d_3\right)(\delta^\lambda_\gamma - s_3 b^\lambda_\gamma)\delta^\beta_\lambda = \left(d_{\beta|\alpha} + s_3\theta_{\beta|\alpha} - b^\beta_\alpha d_3\right)(\delta^\beta_\gamma - s_3 b^\beta_\gamma)\\ &= \left(d_{\gamma|\alpha} - b^\gamma_\alpha d_3\right) + s_3\left(\theta_{\gamma|\alpha} - d_{\beta|\alpha}b^\beta_\gamma - b^\beta_\alpha b^\beta_\gamma d_3\right) - s_3^2\theta_{\beta|\alpha}b^\beta_\gamma, \qquad \alpha,\gamma = 1,2, \end{aligned}\tag{B.84}$$

$$\partial_\alpha \boldsymbol{D}\cdot \boldsymbol{g}_3 = \left(d_\beta b^\beta_\alpha + \partial_\alpha d_3\right) + s_3\theta_\beta b^\beta_\alpha \qquad \alpha = 1,2, \tag{B.85}$$

$$\partial_3 \boldsymbol{D}\cdot \boldsymbol{g}_\gamma = \theta_\alpha(\delta^\lambda_\gamma - s_3 b^\lambda_\gamma)\delta^\alpha_\lambda = \theta_\alpha(\delta^\alpha_\gamma - s_3 b^\alpha_\gamma) = \theta_\gamma - s_3\theta_\alpha b^\alpha_\gamma \qquad \gamma = 1,2, \tag{B.86}$$

$$\partial_3 \boldsymbol{D}\cdot \boldsymbol{g}_3 = 0. \tag{B.87}$$

Since the linear part of the strain tensor is $E_{ij} = \frac{1}{2}\left(\partial_i \boldsymbol{D}\cdot \boldsymbol{g}_j + \partial_j \boldsymbol{D}\cdot \boldsymbol{g}_i\right)$ we obtain

$$\begin{aligned} &E_{\alpha\beta} = \gamma_{\alpha\beta}(\boldsymbol{d}) + s_3\chi_{\alpha\beta}(\boldsymbol{d},\boldsymbol{\theta}) - s_3^2\kappa_{\alpha\beta}(\boldsymbol{\theta}), && \alpha,\beta = 1,2,\\ &E_{\alpha 3} = \zeta_\alpha(\boldsymbol{d},\boldsymbol{\theta}), && \alpha = 1,2,\\ &E_{33} = 0, \end{aligned}\tag{B.88}$$

where

$$\begin{aligned} \gamma_{\alpha\beta}(\boldsymbol{d}) &= \frac{1}{2}\left(d_{\alpha|\beta} + d_{\beta|\alpha}\right) - b^\beta_\alpha d_3\\ \chi_{\alpha\beta}(\boldsymbol{d},\boldsymbol{\theta}) &= \frac{1}{2}\left(\theta_{\alpha|\beta} + \theta_{\beta|\alpha} - b^\lambda_\beta d_{\lambda|\alpha} - b^\lambda_\alpha d_{\lambda|\alpha}\right) - b^\lambda_\alpha b^\lambda_\beta d_3\\ \kappa_{\alpha\beta}(\boldsymbol{\theta}) &= \frac{1}{2}\left(b^\lambda_\beta\theta_{\lambda|\alpha} + b^\lambda_\alpha\theta_{\lambda|\beta}\right)\\ \zeta_\alpha(\boldsymbol{d},\boldsymbol{\theta}) &= \frac{1}{2}\left(\theta_\alpha + \partial_\alpha d_3 + b^\lambda_\alpha d_\lambda\right). \end{aligned}\tag{B.89}$$

The quantities $\gamma_{\alpha\beta}(\boldsymbol{d})$, $\chi_{\alpha\beta}(\boldsymbol{d},\boldsymbol{\theta})$ and $\zeta_\alpha(\boldsymbol{d},\boldsymbol{\theta})$ represent the covariant components of the *membrane*, *bending* and *shear* strain tensors of the shell, respectively. Thanks to that expression of the Green-Lagrange strain tensor and the constitutive relation, we obtain the first Piola-Kirchhoff stress tensor. As usual in shell modeling, we make the assumption of plane stresses, that is stresses vanish along the normal direction. The stress-strain relationship, if the material is considered linear elastic and isotropic, is given by Hooke's law

$$S^{ij} = H^{ijkl} E_{kl}, \qquad i,j,k,l = 1,2,3. \tag{B.90}$$

In the previous formula $\boldsymbol{S}$ is the second Piola-Kirchhoff stress tensor and $\boldsymbol{H}$ is the elasticity tensor expressed in curvilinear coordinates

$$H^{ijkl} = L_1 g^{ij} g^{kl} + L_2 (g^{ik} g^{jl} + g^{il} g^{jk}). \tag{B.91}$$

The Lamé constants $L_1$ and $L_2$ are related to Young and Poisson moduli as in (B.59). In the following statements we will use the notation $\partial_{tt}$ in order to indicate the second partial derivative with respect to time and we denote by $\mathscr{BC}$ the set

$$\mathscr{BC} = \{\boldsymbol{v} : \Sigma \to \mathbb{R}^d : \boldsymbol{v} \text{ is an admissible function that satisfy the boundary conditions}\}$$

Using Hooke's law and equations (B.88) and (B.89) we obtain the following weak form of equilibrium equation (B.52)

**Problem 13.** *Let $\boldsymbol{f}(t) \in L^2(\Sigma)$, Let be $\mathscr{C} = \{(\boldsymbol{v}, \boldsymbol{\eta}) \in H^1(\Sigma)^d \times H^1(\Sigma)^d\} \cap \mathscr{BC}$ be the set of admissible solution of the problem (function that respect the boundary conditions $\mathscr{BC}$ and are sufficiently regular) with norm $\|(\boldsymbol{v}, \boldsymbol{\eta})\|^2_{\mathscr{C}} \stackrel{def}{=} \|\boldsymbol{v}\|^2_{1,\Sigma} + \|\boldsymbol{\eta}\|^2_{1,\Sigma}$, then find $(\boldsymbol{d}, \boldsymbol{\theta}) \in \mathscr{C}$ such that, for all $(\boldsymbol{w}, \boldsymbol{\eta}) \in \mathscr{C}$*

$$\rho_s \epsilon (\partial_{tt} \boldsymbol{d}, \boldsymbol{w})_{0,\Sigma \times [-\frac{\epsilon}{2}, \frac{\epsilon}{2}]} + A((\boldsymbol{d}, \boldsymbol{\theta}), (\boldsymbol{v}, \boldsymbol{\eta})) = (\boldsymbol{f}, \boldsymbol{w})_{0,\Sigma \times [-\frac{\epsilon}{2}, \frac{\epsilon}{2}]}, \tag{B.92}$$

*where*

$$A((\boldsymbol{d}, \boldsymbol{\theta}), (\boldsymbol{v}, \boldsymbol{\eta})) \stackrel{def}{=} \left( C^{\alpha\beta\lambda\delta} E_{\lambda\delta}(\boldsymbol{d}, \boldsymbol{\theta}), E_{\alpha\beta}(\boldsymbol{w}, \boldsymbol{\eta}) + D^{\lambda\delta} E_{\lambda 3}(\boldsymbol{d}, \boldsymbol{\theta}), E_{\delta 3}(\boldsymbol{w}, \boldsymbol{\eta}) \right)_{0,\Sigma \times [-\frac{\epsilon}{2}, \frac{\epsilon}{2}]} \tag{B.93}$$

$$C^{\alpha\beta\lambda\delta} = \frac{E}{2(1+\nu)} \left( g^{\alpha\lambda} g^{\beta\delta} + g^{\alpha\delta} g^{\beta\lambda} + \frac{2\nu}{1-\nu} g^{\alpha\beta} g^{\lambda\delta} \right), \qquad D^{\lambda\delta} = \frac{2E}{1+\nu} g^{\lambda\delta}. \tag{B.94}$$

The bilinear form (B.93) of Problem 13 is continuous and coercive on the space $(\mathscr{C}, \|\cdot\|_{\mathscr{C}})$ if the boundary condition are such that rigid motions of the solid is not allowed, then the stationary version of Problem 13 admits a unique solution in $\mathscr{C}$. This result is stated in the following Proposition and its proof can be found in [39, Proposition 4.3.5]

**Proposition B.5.1.** *Suppose that the space $(\mathscr{C}, \|\cdot\|_{\mathscr{C}})$ does not contains rigid body motions, then the bilinear form* (B.93)

$$A((\boldsymbol{d}, \boldsymbol{\theta}), (\boldsymbol{v}, \boldsymbol{\eta})) : \mathscr{C} \times \mathscr{C} \to \mathbb{R}$$

*of Problem 13 is continuous and coercive on the space $(\mathscr{C}, \|\cdot\|_{\mathscr{C}})$*

**Membrane shell** A simplified model of shell theory, namely the membrane shell model, is useful especially when the stress momentum and shearing forces are neglected compared with the stress resultants. This particular stress state is generally realized in very thin shell or, if the thickness is not negligible, in regions away from edges of the shell. The membrane shell theory is a consequence of the general shell theory outlined in the previous paragraph, and in particular it is obtained considering a small thickness $\epsilon$, then the part of the strain tensor $s_3^2\underline{\underline{\kappa}}$ can be neglected (since $-\frac{\epsilon}{2} < s_3 < \frac{\epsilon}{2}$). In this case the quantities defined in (B.89) become

$$\begin{aligned} E_{\alpha\beta} &= \gamma_{\alpha\beta}(\boldsymbol{d}) + s_3\chi_{\alpha\beta}(\boldsymbol{d},\boldsymbol{\theta}) \qquad && \alpha,\beta = 1,2, \\ E_{\alpha 3} &= \zeta_\alpha(\boldsymbol{d},\boldsymbol{\theta}) && \alpha = 1,2, \\ E_{33} &= 0. \end{aligned} \tag{B.95}$$

Moreover, since the thickness is very small, we can assume that $\boldsymbol{a}_i \approx \boldsymbol{g}_i$ for $i = 1,2,3$, then

$$C_m^{\alpha\beta\lambda\delta} = \frac{E}{2(1+\nu)}\left(a^{\alpha\lambda}a^{\beta\delta} + a^{\alpha\delta}a^{\beta\lambda} + \frac{2\nu}{1-\nu}a^{\alpha\beta}a^{\lambda\delta}\right), \qquad D_m^{\lambda\delta} = \frac{2E}{1+\nu}a^{\lambda\delta}. \tag{B.96}$$

In order to eliminate the dependence of the strain component on the rotation $\theta$, we assume the Kirchhoff-Love kinematic hypothesis: *material lines orthogonal to the midsurface in the undeformed configuration remain straight, unstretched, and always orthogonal to the midsurface during the deformations.* It can be proved (see [39, Proposition 4.2.1]) that if a displacement field satisfies the Kirchhoff-Love hypothesis then

$$\zeta_\alpha(\boldsymbol{d},\boldsymbol{\theta}) = 0 \quad \Leftrightarrow \quad \theta_\alpha = -\partial_\alpha d_3 - b_\alpha^\lambda d_\lambda \qquad \text{for} \qquad \alpha = 1,2.$$

Under this restriction, the bending tensor $\underline{\underline{\chi}}$ does not depend on $\theta$ but it is of second order with respect to the third component of displacement $d_3$. The weak formulation of the membrane shell theory reads

**Problem 14.** *Find* $(d_1(t), d_2(t), d_3(t)) \in \boldsymbol{W} = (H^1(\Sigma)^2 \times H^2(\Sigma)) \cap \mathscr{BC}$

$$\rho_s\epsilon(\partial_{tt}\boldsymbol{d},\boldsymbol{w})_{0,\Sigma} + \epsilon\left(C_m^{\alpha\beta\lambda\delta}\gamma_{\lambda\delta}(\boldsymbol{d}),\gamma_{\alpha\beta}(\boldsymbol{w})\right)_{0,\Sigma} + \frac{\epsilon^3}{12}\left(C_m^{\alpha\beta\lambda\delta}\chi_{\lambda\delta}(\boldsymbol{d}),\chi_{\alpha\beta}(\boldsymbol{w})\right)_{0,\Sigma} = \epsilon(\boldsymbol{f},\boldsymbol{w})_{0,\Sigma}, \quad \forall \boldsymbol{w} \in \boldsymbol{W}. \tag{B.97}$$

Introducing the parameter $\lambda \stackrel{def}{=} \frac{\epsilon}{L}$ where $\epsilon$ is the thickness and $L$ is an overall characteristic dimension of the shell structure, we can see that the general form of that variational Problem 15 is

$$\lambda A_m(\boldsymbol{d}^\lambda,\boldsymbol{w}) + \lambda^3 A_b(\boldsymbol{d}^\lambda,\boldsymbol{w}) = \lambda(\boldsymbol{F},\boldsymbol{w})_{0,\Sigma}, \tag{B.98}$$

where, denoting $\partial_{tt}$ the second partial derivative with respect to time,

$$A_m = \left(LC_m^{\alpha\beta\lambda\delta}\gamma_{\lambda\delta}(\boldsymbol{d}),\gamma_{\alpha\beta}(\boldsymbol{w})\right)_{0,\Sigma}, \quad A_b = \frac{1}{12}\left(LC_m^{\alpha\beta\lambda\delta}\chi_{\lambda\delta}(\boldsymbol{d}),\chi_{\alpha\beta}(\boldsymbol{w})\right)_{0,\Sigma}, \quad \boldsymbol{F} = \boldsymbol{f} - \rho_s\partial_{tt}\boldsymbol{d}.$$

This means that the behavior of the solid comes from the superposition of a pure membrane model ($A_m$) and a pure bending model ($A_b$). In many situations it is interesting to consider very thin structures, namely structure for which the parameter $\lambda$ goes to zero. Structures of this type can exhibit a membrane-dominated behavior or a bending-dominated behavior depending on the constraint on the solid and on the applied loads.

**The simplified problem** In any case, the analysis of the coupling of this type of structures with fluids presents many challenging difficulties. For this reason in this work, we will often use a simplified structural model, namely the generalized string model, which is obtained by considering longitudinal sections of a 2D cylindrical surface of radius R and length L, then the mid-surface in this case is flat (see [62, Chapter 4]); in weak form we have

$$\boldsymbol{d} \in \boldsymbol{W}: \qquad \rho_s \epsilon \int_\Sigma \partial_{tt} \boldsymbol{d} \cdot \boldsymbol{w} + a_s(\boldsymbol{d}, \boldsymbol{w}) = \int_\Sigma \boldsymbol{f} \cdot \boldsymbol{w} \quad \forall \boldsymbol{w} \in \boldsymbol{W}, \tag{B.99}$$

where $\boldsymbol{W}$ is the subspace $H^1(\Sigma)^2 \cap \mathscr{BC}$ of $H^1(\Sigma)^2$ consisting of displacements and velocities orthogonal to the mid-surface that satisfy boundary conditions. In 2D, with the longitudinal and perpendicular direction lying on the x-axis and y-axis respectively, the model is retrieved considering

$$\boldsymbol{d} = \begin{pmatrix} 0 \\ \eta \end{pmatrix}, \qquad a_s(\boldsymbol{d}, \boldsymbol{w}) = \int_\Sigma \left( \lambda_0 \partial_s \eta \partial_s w + \lambda_1 \eta w \right) ds, \quad \forall \boldsymbol{w} = \begin{pmatrix} 0 \\ w \end{pmatrix} \in \boldsymbol{W} \subset H^1(\Sigma)^2 \tag{B.100}$$

with

$$\lambda_0 = \frac{E\epsilon}{2(1+\nu)}, \qquad \lambda_1 = \frac{E\epsilon}{R^2(1-\nu^2)}. \tag{B.101}$$

We observe explicitly that the bilinear form $a_s(\cdot,\cdot)$ defined in (B.100) and (B.101) is symmetric, continuous and coercive on the space $W = \{\boldsymbol{w} = (0,w)^T \in H^1(\Sigma)^2\}$ endowed with the $H^1$-norm: there exist positive constants $C_1$ and $C_2$, such that for all $\boldsymbol{d}, \boldsymbol{w} \in \boldsymbol{W}$ it holds

- **Continuity**

$$\begin{aligned} a_s(\boldsymbol{d}, \boldsymbol{w}) = \int_\Sigma \left( \lambda_0 \partial_s \eta \partial_s w + \lambda_1 \eta w \right) ds &\le C \left( \|\partial_s \eta\|_{0,\Sigma} \|\partial_s w\|_{0,\Sigma} + \|\eta\|_{0,\Sigma} \|w\|_{0,\Sigma} \right) \\ &\le C_1 \left( \|\boldsymbol{d}\|_{1,\Sigma} + \|\boldsymbol{w}\|_{1,\Sigma} \right), \end{aligned} \tag{B.102}$$

- **Coercivity**

$$a_s(\boldsymbol{w}, \boldsymbol{w}) = \int_\Sigma \left( \lambda_0 (\partial_s w)^2 + \lambda_1 w^2 \right) ds \ge C_2 \left( \|\partial_s w\|_{0,\Sigma}^2 + \|w\|_{0,\Sigma}^2 \right) = C_2 \|\boldsymbol{w}\|_{1,\Sigma}^2. \tag{B.103}$$

We derive the solid model used in the following starting by considering the *membrane shell* model presented in Problem 15 of Appendix B in which it is neglected the bending term. In developing the solid model we proceed in two steps; in the first step we develop a solid model in the case that the reference configuration of the solid is stress-free, in the second step we enrich the model supposing that the reference position is not stress-free; in the latter case we have to add a prestress term to the equations.

### B.5.4 Solid model for stress-free reference configuration

In the case $\Sigma = \Sigma(0)$, let $\boldsymbol{d}$ be the displacement from the reference position and let $\boldsymbol{W}$ the space of admissible displacements, then the problem that we are considering can be written in weak form as follows.

**Problem 15.** *Given a forcing term $\boldsymbol{f}$, find $\boldsymbol{d}(\cdot,t) \in \boldsymbol{W}$ such that*

$$\rho_s \epsilon \int_\Sigma \partial_{tt}\boldsymbol{d} \cdot \boldsymbol{w} ds + \epsilon \int_\Sigma C_m^{\alpha\beta\lambda\delta} \gamma_{\lambda\delta}(\boldsymbol{d}) \gamma_{\alpha\beta}(\boldsymbol{w}) ds = \epsilon \int_\Sigma \boldsymbol{f} \cdot \boldsymbol{w} ds \qquad \forall \boldsymbol{w} \in \boldsymbol{W}. \tag{B.104}$$

In the previous equation $\rho_s$ represents the solid density and $\epsilon$ is the solid thickness. The meaning of the other terms is presented in Appendix B; we recall only that $C_m^{\alpha\beta\lambda\delta}$ is the elastic tensor and $\gamma_{\alpha\beta}(\boldsymbol{d})$ is the membrane strain tensor whose expressions are:

$$C_m^{\alpha\beta\lambda\delta} = \frac{E}{2(1+\nu)} \left( a^{\alpha\lambda} a^{\beta\delta} + a^{\alpha\delta} a^{\beta\lambda} + \frac{2\nu}{1-\nu} a^{\alpha\beta} a^{\lambda\delta} \right), \tag{B.105}$$

$$\gamma_{\alpha\beta}(\boldsymbol{d}) = \frac{1}{2} \left( d_{\alpha|\beta} + d_{\beta|\alpha} \right) - b_\alpha^\beta d_3. \tag{B.106}$$

The proper choice of the functional space $\boldsymbol{W}$ depends on the boundary conditions imposed on the displacement $\boldsymbol{d}$. We point out that for structural problems it is a common practice to adopt a Lagrangian framework, indeed we write the equations of the membrane model in the reference configuration $\Sigma$.

We point out that the last assumption of our model imposes that the displacements from the reference configuration $\Sigma$ are normal, meaning that, with reference to the local basis $(\boldsymbol{a}^1, \boldsymbol{a}^2, \boldsymbol{a}_3)$, we consider displacements of the type $\boldsymbol{d} = (0, 0, d_3)$. Using this expression in the equation (B.106) we derive $\gamma_{\alpha\beta}(\boldsymbol{d}) = -b_\alpha^\beta d_3$. Then the model can be further simplified, in fact we have

$$C_m^{\alpha\beta\lambda\delta} \gamma_{\lambda\delta}(\boldsymbol{d}) \gamma_{\alpha\beta}(\boldsymbol{w}) = C_m^{\alpha\beta\lambda\delta} b_\lambda^\delta b_\alpha^\beta d_3 w_3 = K d_3 w_3,$$

with K a constant depending from the elastic properties of the material and from the geometry of the structure. Using the last expression in (B.104), we derive easily the following structural model

$$\begin{cases} d_1 = 0, \qquad d_2 = 0 & \text{in } \Sigma \times [0,T], \\ \rho_s \epsilon \partial_{tt}\boldsymbol{d} + K\boldsymbol{d} = \boldsymbol{f} & \text{in } \Sigma \times [0,T], \\ \boldsymbol{d}(0) = \boldsymbol{d}_0 & \text{in } \Sigma, \\ \partial_t \boldsymbol{d}(0) = \dot{\boldsymbol{d}}_0 & \text{in } \Sigma. \end{cases} \tag{B.107}$$

### B.5.5 Solid model for general reference configurations

The very simple model developed in the previous subsection can be enriched adding a prestress term that is useful if we are working with reference configurations that are not stress-free. A typical case in which this term is important is the one of the fluid-structure interaction in hemodynamics, in fact, the arteries and blood vessels are in a prestressed state. Starting from the nonlinear elasticity equations for a general shell-type stress-free domain $\hat{\mathscr{S}} = \hat{\Sigma} \times [-\hat{\epsilon}/2, \hat{\epsilon}/2]$ and linearizing around a deformed configuration $\mathscr{S} = \Sigma \times [-\epsilon/2, \epsilon/2]$, we obtain a model that describes the small deformations around $\mathscr{S}$. We point out that $\mathscr{S}$ is not, in general, an equilibrium position and it is obtained from the reference one by the displacement field $\hat{\boldsymbol{\eta}}$, that is, $\mathscr{S} = \hat{\mathscr{S}} + \hat{\boldsymbol{\eta}}$.

The Lagrangian weak formulation of the static elastic problem is written in the reference configuration $\hat{\mathscr{S}}$ using the Second Piola-Kirchhoff stress tensor $\boldsymbol{S}$ that is related to the Green-Lagrange strain tensor $\hat{\boldsymbol{E}}$ as follows

$$\boldsymbol{S} = \boldsymbol{H} : \hat{\boldsymbol{E}},$$

where $\boldsymbol{H}$ is the fourth order tensor of the elastic constants. Denoting by $\hat{\boldsymbol{W}}$ the functional space of admissible displacements defined on the reference configuration $\hat{\mathscr{S}}$, and $\hat{\boldsymbol{W}}_0$ the space of virtual displacements from $\hat{\mathscr{S}}$.

**Problem 16.** *find* $\hat{\boldsymbol{d}} \in \hat{\boldsymbol{W}}$ *such that, for all* $\hat{\boldsymbol{w}} \in \hat{\boldsymbol{W}}_0$,

$$\hat{W}(\hat{\boldsymbol{d}}, \hat{\boldsymbol{w}}) \stackrel{def}{=} \rho_s \int_{\hat{\mathscr{S}}} \partial_{tt}\hat{\boldsymbol{d}} \cdot \hat{\boldsymbol{w}} d\boldsymbol{x} + \int_{\hat{\mathscr{S}}} \hat{\boldsymbol{E}}(\hat{\boldsymbol{d}}) : \boldsymbol{H} : \hat{\boldsymbol{E}}(\hat{\boldsymbol{w}}) d\boldsymbol{x} - \int_{\hat{\mathscr{S}}} \hat{\boldsymbol{b}} \cdot \hat{\boldsymbol{w}} d\boldsymbol{x} - \int_{\partial\hat{\mathscr{S}}} \hat{\boldsymbol{t}} \cdot \hat{\boldsymbol{w}} d\boldsymbol{x} = 0. \quad \text{(B.108)}$$

Since we are considering the problem of small displacements around the position $\mathscr{S}$, we consider the linear part of the functional $\hat{W}(\cdot,\cdot)$ near the position $\mathscr{S}$. For "small" displacements from $\mathscr{S}$, $\boldsymbol{d} : \mathscr{S} \to \mathbb{R}^d$, denoted as $\hat{\boldsymbol{d}}$ in the configuration $\hat{\mathscr{S}}$, we can write

$$\hat{W}(\hat{\boldsymbol{\eta}} + \hat{\boldsymbol{d}}, \hat{\boldsymbol{w}}) \cong \hat{W}(\hat{\boldsymbol{\eta}}, \hat{\boldsymbol{w}}) + D\hat{W}(\hat{\boldsymbol{\eta}}, \hat{\boldsymbol{w}})[\hat{\boldsymbol{d}}]. \quad \text{(B.109)}$$

The displacement field $\hat{\boldsymbol{d}}$ is defined as $\hat{\boldsymbol{d}} \stackrel{def}{=} \boldsymbol{d} \circ \hat{\boldsymbol{\phi}}$ where $\hat{\boldsymbol{\phi}} : \hat{\mathscr{S}} \to \mathscr{S}$ represents the deformation map.

$$\hat{W}(\hat{\boldsymbol{\eta}} + \hat{\boldsymbol{d}}, \hat{\boldsymbol{w}}) \cong \hat{W}(\hat{\boldsymbol{\eta}}, \hat{\boldsymbol{w}}) + D\hat{W}(\hat{\boldsymbol{\eta}}, \hat{\boldsymbol{w}})[\hat{\boldsymbol{d}}] \quad \text{(B.110)}$$

where $D\hat{W}(\hat{\boldsymbol{\eta}}, \hat{\boldsymbol{w}})[\hat{\boldsymbol{d}}]$ is the derivative of $\hat{W}(\cdot,\cdot)$ evaluated in $(\hat{\boldsymbol{\eta}}, \hat{\boldsymbol{w}})$ in the direction of $\hat{\boldsymbol{d}}$ and its expression is given in the following.

$$\begin{aligned} D\hat{W}(\hat{\boldsymbol{\eta}}, \hat{\boldsymbol{w}})[\hat{\boldsymbol{d}}] &= \lim_{\epsilon \to 0} \frac{\hat{W}(\hat{\boldsymbol{\eta}} + \epsilon\hat{\boldsymbol{d}}, \hat{\boldsymbol{w}}) - \hat{W}(\hat{\boldsymbol{\eta}}, \hat{\boldsymbol{w}})}{\epsilon} \\ &= \rho_s \int_{\hat{\mathscr{S}}} \partial_{tt}\hat{\boldsymbol{d}} \cdot \hat{\boldsymbol{w}} d\boldsymbol{x} + \int_{\hat{\mathscr{S}}} \hat{\boldsymbol{\epsilon}}(\hat{\boldsymbol{d}}) : \boldsymbol{H} : \hat{\boldsymbol{E}}(\hat{\boldsymbol{w}}) d\boldsymbol{x} + \int_{\hat{\mathscr{S}}} (\nabla\hat{\boldsymbol{\eta}})^T \nabla\hat{\boldsymbol{d}} : \boldsymbol{H} : \hat{\boldsymbol{E}}(\hat{\boldsymbol{w}}) d\boldsymbol{x}, \quad \text{(B.111)} \end{aligned}$$

where $\hat{\boldsymbol{\epsilon}}(\hat{\boldsymbol{d}}) = \frac{1}{2}(\nabla\hat{\boldsymbol{d}} + (\nabla\hat{\boldsymbol{d}})^T)$ is the linear part of the Green-Lagrange strain tensor. Neglecting the non linear part of $\hat{\boldsymbol{E}}(\hat{\boldsymbol{w}})$ and assuming that, for some $\boldsymbol{f} \in \hat{\boldsymbol{W}}_0$, we can write $W(\hat{\boldsymbol{\eta}}, \hat{\boldsymbol{w}}) + \int_{\hat{\mathscr{S}}} \boldsymbol{f} \cdot \hat{\boldsymbol{w}} = 0, \quad \forall \hat{\boldsymbol{w}} \in \hat{\boldsymbol{W}}_0$, the linearized weak formulation is: find $\hat{\boldsymbol{d}} \in \hat{\boldsymbol{W}}$ such that

$$\begin{aligned} \widetilde{W}(\hat{\boldsymbol{d}}, \hat{\boldsymbol{w}}) \stackrel{def}{=} \rho_s \int_{\hat{\mathscr{S}}} \partial_{tt}\hat{\boldsymbol{d}} \cdot \hat{\boldsymbol{w}} d\boldsymbol{x} + \int_{\hat{\mathscr{S}}} \hat{\boldsymbol{\epsilon}}(\boldsymbol{d}) : \boldsymbol{H} : \hat{\boldsymbol{\epsilon}}(\hat{\boldsymbol{w}}) d\boldsymbol{x} \\ + \int_{\hat{\mathscr{S}}} (\nabla\hat{\boldsymbol{\eta}})^T \nabla\boldsymbol{d} : \boldsymbol{H} : \hat{\boldsymbol{\epsilon}}(\hat{\boldsymbol{w}}) d\boldsymbol{x} - \int_{\hat{\mathscr{S}}} \boldsymbol{f} \cdot \hat{\boldsymbol{w}} = 0 \qquad \forall \hat{\boldsymbol{w}} \in \hat{\boldsymbol{W}}_0. \quad \text{(B.112)} \end{aligned}$$

The tensor $T \stackrel{def}{=} \left(\nabla\hat{\boldsymbol{\eta}}\right)^T : \boldsymbol{H}$ represents the prestress tensor. In particular, since we are interested in deriving a membrane model, we consider a displacement field from the configuration $\mathscr{S}$ that, in the local basis $(\boldsymbol{a}_1, \boldsymbol{a}_2, \boldsymbol{a}_3)$ has the form

$$\boldsymbol{d} = d_3\boldsymbol{a}_3. \tag{B.113}$$

Introducing the following quantities

$$\Theta_i \stackrel{def}{=} (\nabla\hat{\boldsymbol{\eta}})^{-1}_{i3}, \qquad \hat{d}_3 \stackrel{def}{=} d_3(\hat{\boldsymbol{\eta}}), \qquad i = 1,2,3,$$

this displacement model, in the configuration $\hat{\mathscr{S}}$ is represented as

$$\hat{d}_i = (\nabla\hat{\boldsymbol{\eta}})^{-1}_{i3} d_3(\hat{\boldsymbol{\eta}}) = \Theta_i\hat{d}_3, \qquad i = 1,2,3. \tag{B.114}$$

Denoting by $\mathscr{BC}$ the set of functions $\boldsymbol{v} : \hat{\Sigma} \to \mathbb{R}^d$ that satisfy the essential boundary conditions, and by $\mathscr{BC}_0$ the set of functions that satisfy the homogeneous essential boundary conditions, then space of admissible and virtual displacements can be specified as

$$\hat{\boldsymbol{W}} = \left\{\hat{\boldsymbol{v}} \in H^1(\hat{\Sigma})^{d-1} : \hat{\boldsymbol{v}} = \boldsymbol{\Theta} v, v \in H^1(\hat{\Sigma})\right\} \cap \mathscr{BC},$$
$$\hat{\boldsymbol{W}}_0 = \left\{\hat{\boldsymbol{v}} \in H^1(\hat{\Sigma})^{d-1} : \hat{\boldsymbol{v}} = \boldsymbol{\Theta} v, v \in H^1(\hat{\Sigma})\right\} \cap \mathscr{BC}_0.$$

We observe that for displacement fields of type (B.114), we have

$$\hat{\epsilon}_{ij}(\hat{\boldsymbol{d}}) = \hat{\epsilon}_{ij}(\boldsymbol{\Theta})\hat{d}_3 + \frac{1}{2}\Theta_i\partial_j\hat{d}_3 + \frac{1}{2}\Theta_j\partial_i\hat{d}_3$$

Using this expression in B.112 and assuming that all te quantities are constant across the thickness, which is assumed very tiny, we can write, $\forall\hat{\boldsymbol{w}} \in \hat{\boldsymbol{W}}_0$

$$\begin{aligned}
\widetilde{W}(\hat{\boldsymbol{d}}, \hat{\boldsymbol{w}}) =& \rho_s\hat{\epsilon}\int_{\hat{\Sigma}} \boldsymbol{\Theta}\cdot\boldsymbol{\Theta}\partial_{tt}\hat{d}_3\hat{w}_3 d\boldsymbol{x} + \hat{\epsilon}\int_{\hat{\Sigma}} \left[\hat{\epsilon}_{ij}(\boldsymbol{\Theta})H^{ijkl}\hat{\epsilon}_{kl}(\boldsymbol{\Theta})\right]\hat{d}_3\hat{w}_3 d\boldsymbol{x} \\
&+ \frac{1}{2}\hat{\epsilon}\int_{\hat{\Sigma}} \left[\Theta_j H^{ijkl}\hat{\epsilon}_{kl}(\boldsymbol{\Theta})\right]\partial_i\hat{d}_3\hat{w}_3 d\boldsymbol{x} + \frac{1}{2}\hat{\epsilon}\int_{\hat{\Sigma}} \left[\Theta_i H^{ijkl}\hat{\epsilon}_{kl}(\boldsymbol{\Theta})\right]\partial_j\hat{d}_3\hat{w}_3 d\boldsymbol{x} \\
&+ \frac{1}{2}\hat{\epsilon}\int_{\hat{\Sigma}} \left[\hat{\epsilon}_{ij}(\boldsymbol{\Theta})H^{ijkl}\Theta_l\right]\partial_k\hat{w}_3\hat{d}_3 d\boldsymbol{x} + \frac{1}{2}\hat{\epsilon}\int_{\hat{\Sigma}} \left[\hat{\epsilon}_{ij}(\boldsymbol{\Theta})H^{ijkl}\Theta_k\right]\partial_l\hat{w}_3\hat{d}_3 d\boldsymbol{x} \\
&+ \hat{\epsilon}\int_{\hat{\Sigma}} \left[\partial_j\hat{\eta}_i\partial_i\Theta_k H^{jklm}\hat{\epsilon}_{lm}(\boldsymbol{\Theta})\right]\hat{d}_3\hat{w}_3 d\boldsymbol{x} \\
&+ \frac{1}{2}\hat{\epsilon}\int_{\hat{\Sigma}} \left[\partial_j\hat{\eta}_i\partial_i\Theta_k H^{jklm}\Theta_m\right]\partial_l\hat{w}_3\hat{d}_3 d\boldsymbol{x} + \frac{1}{2}\hat{\epsilon}\int_{\hat{\Sigma}} \left[\partial_j\hat{\eta}_i\partial_i\Theta_k H^{jklm}\Theta_l\right]\partial_m\hat{w}_3\hat{d}_3 d\boldsymbol{x} \\
&+ \hat{\epsilon}\int_{\hat{\Sigma}} \left[\partial_j\hat{\eta}_i\Theta_k H^{jklm}\hat{\epsilon}_{lm}(\boldsymbol{\Theta})\right]\partial_i\hat{d}_3\hat{w}_3 d\boldsymbol{x} \\
&+ \frac{1}{2}\hat{\epsilon}\int_{\hat{\Sigma}} \left[\partial_j\hat{\eta}_i\Theta_k H^{jklm}\Theta_m\right]\partial_l\hat{w}_3\partial_i\hat{d}_3 d\boldsymbol{x} + \frac{1}{2}\hat{\epsilon}\int_{\hat{\Sigma}} \left[\partial_j\hat{\eta}_i\Theta_k H^{jklm}\Theta_l\right]\partial_m\hat{w}_3\partial_i\hat{d}_3 d\boldsymbol{x} \\
&- \int_{\hat{\Sigma}} \boldsymbol{f}\cdot\hat{\boldsymbol{w}} = 0.
\end{aligned} \tag{B.115}$$

Rearranging the terms and integrating by part the integrals that contain the products $\hat{d}_3\partial_j\hat{w}_3$ for $j = 1,2,3$, equation (B.115) can be rewritten in the following form

$$\widetilde{W}(\hat{\boldsymbol{d}},\hat{\boldsymbol{w}}) = \rho_s\hat{\epsilon}\int_{\hat{\Sigma}} K_0\partial_{tt}\hat{d}_3\hat{w}_3 d\boldsymbol{x} + \hat{\epsilon}\int_{\hat{\Sigma}} K_1^{ij}\partial_i\hat{d}_3\partial_j\hat{w}_3 d\boldsymbol{x} + \hat{\epsilon}\int_{\hat{\Sigma}} K_2^j\partial_j\hat{d}_3\hat{w}_3 d\boldsymbol{x}$$
$$+\hat{\epsilon}\int_{\hat{\Sigma}} K_3\hat{d}_3\hat{w}_3 d\boldsymbol{x} + \hat{\epsilon}\int_{\partial\hat{\Sigma}} K_4^j n_j\hat{d}_3\hat{w}_3 d\boldsymbol{x} - \int_{\hat{\Sigma}} \boldsymbol{f}\cdot\hat{\boldsymbol{w}} = 0 \qquad \forall\hat{\boldsymbol{w}}\in\hat{\boldsymbol{W}}_0. \quad \text{(B.116)}$$

The last equation imply the following strong form of the solid problem

$$\begin{cases} d_1 = \Theta_1\hat{d}_3, \qquad d_2 = \Theta_2\hat{d}_3 & \text{in} \quad \Sigma\times[0,T], \\ K_0\partial_{tt}\hat{d}_3 + K_1^{ij}\partial_{ij}\hat{d}_3 + K_2^j\partial_j\hat{d}_3 + K_3\hat{d}_3 = \boldsymbol{f}\cdot\boldsymbol{\Theta} & \text{in} \quad \Sigma\times[0,T], \\ K_4^j n_j\hat{d}_3\hat{w}_3 = 0 & \text{on} \quad \partial\Sigma, \\ \boldsymbol{d}(0) = \boldsymbol{d}_0 & \text{in} \quad \Sigma, \\ \partial_t\boldsymbol{d}(0) = \dot{\boldsymbol{d}}_0 & \text{in} \quad \Sigma, \end{cases} \quad \text{(B.117)}$$

with $\hat{\boldsymbol{d}} = \boldsymbol{\Theta}\hat{d}_3$ satisfying the essential boundary conditions that can prescribe the value of $\boldsymbol{d}\big|_{\partial\Sigma}$ and $\nabla\boldsymbol{d}\cdot\boldsymbol{v}$, where $\boldsymbol{v}$ is the normal unit outward vector to $\partial\Sigma$ as shown in Figure 1.1. In the equations we have also the constant quantities $K_0, K_1^{ij}, K_2^j$ and $K_3$ that depends on the material and the geometry of the configuration $\mathscr{S}$.

# Appendix C

# F.E. analysis of the Stokes Problem

In this appendix we recall the analysis of the Stokes Problem. The results of this analysis are useful for the definition of the projection operators used in the study of convergence of the coupled problem.

## C.1 Finite element solution of Stokes equations

Before the introduction of the finite element discretization of the Stokes problem, we recall the theoretical analysis of the continuous problem.

### C.1.1 Stokes Problem

Stokes problem is associated to the steady motion of an incompressible fluid with no convective term. The weak formulation is reported in the following Problem 17

**Problem 17.** *Let $\Omega$ be a Lipschitz domain, given $\boldsymbol{b} \in L^2(\Omega)^d$, find $\boldsymbol{u} \in H_0^1(\Omega)^d$ and $p \in L_0^2(\Omega)$ such that*

$$2\mu(\boldsymbol{\epsilon}(\boldsymbol{u}), \boldsymbol{\epsilon}(\boldsymbol{v}))_{0,\Omega} - (p, \operatorname{div} \boldsymbol{v})_{0,\Omega} = \rho(\boldsymbol{b}, \boldsymbol{v})_{0,\Omega} \qquad \forall \boldsymbol{v} \in H_0^1(\Omega)^d, \tag{C.1}$$

$$(q, \operatorname{div} \boldsymbol{u})_{0,\Omega} = 0 \qquad \forall q \in L_0^2(\Omega). \tag{C.2}$$

We recall that $L_0^2(\Omega) = \{q \in L^2(\Omega) : \int_\Omega q\, d\boldsymbol{x} = 0\}$. A key result in the solution of the weak formulation of Stokes Problem 17 is the following lemma in which we use the function space $\boldsymbol{V}_{\text{div}} = \{\boldsymbol{v} \in H_0^1(\Omega)^d : \operatorname{div} \boldsymbol{v} = 0\}$

**Lemma C.1.1.** *Let $\Omega \subset \mathbb{R}^d$ be a Lipschitz domain and let $L \in H^{-1}(\Omega)^d$, then $L$ vanishes identically on $\boldsymbol{V}_{\text{div}}$ if and only if there exist a unique $q_{\boldsymbol{v}} \in L_0^2(\Omega)$ such that*

$$\langle L, \boldsymbol{v} \rangle = \int_\Omega q \operatorname{div} \boldsymbol{v} d\boldsymbol{x} \quad \forall \boldsymbol{v} \in H_0^1(\Omega)^d. \tag{C.3}$$

The proof of this result can be found in [66, Lemma 2.1]. We report instead the proof of the following result that states the existence and uniqueness of the solution to Problem 17

**Theorem C.1.2.** *Problem 17 admits a unique solution*

*Proof.* The proof of the theorem is performed in two steps; the first step consists in the proof of existence and uniqueness of solution to the following auxiliary problem

$$\boldsymbol{u}^{\star} \in \boldsymbol{V}_{\mathrm{div}} : \quad 2\mu\left(\boldsymbol{\epsilon}(\boldsymbol{u}^{\star}), \boldsymbol{\epsilon}(\boldsymbol{v})\right)_{0,\Omega} = \rho(\boldsymbol{b}, \boldsymbol{v})_{0,\Omega} \qquad \forall \boldsymbol{v} \in \boldsymbol{V}_{\mathrm{div}}. \tag{C.4}$$

Problem (C.4) admits a unique solution by Lax-Milgram Lemma A.2.2, in fact $a(\boldsymbol{u},\boldsymbol{v}) \overset{def}{=} 2\mu(\boldsymbol{\epsilon}(\boldsymbol{u}),\boldsymbol{\epsilon}(\boldsymbol{v}))_{0,\Omega}$ is coercive and continuous on $\boldsymbol{V}_{\mathrm{div}}$ as consequences of the Poincaré theorem and Korn's inequality: for all $\boldsymbol{u},\boldsymbol{v} \in \boldsymbol{V}_{\mathrm{div}}$

$$a(\boldsymbol{u},\boldsymbol{v}) = 2\mu(\boldsymbol{\epsilon}(\boldsymbol{u}),\boldsymbol{\epsilon}(\boldsymbol{v}))_{0,\Omega} \le C_1|\boldsymbol{u}|_{1,\Omega}|\boldsymbol{v}|_{1,\Omega} \le C_1\|\boldsymbol{u}\|_{1,\Omega}\|\boldsymbol{v}\|_{1,\Omega} \tag{C.5}$$

and for all $\boldsymbol{v} \in \boldsymbol{V}_{\mathrm{div}}$

$$a(\boldsymbol{v},\boldsymbol{v}) = 2\mu(\boldsymbol{\epsilon}(\boldsymbol{u}),\boldsymbol{\epsilon}(\boldsymbol{v}))_{0,\Omega} \ge C_2|\boldsymbol{v}|^2_{1,\Omega} \ge C_2\|\boldsymbol{v}\|^2_{1,\Omega} \ge C_2\left(\|\boldsymbol{u}\|^2_{0,\Omega} + \|\operatorname{div}\boldsymbol{v}\|^2_{0,\Omega}\right) = C_2\|\boldsymbol{v}\|^2_{\mathrm{div}}. \tag{C.6}$$

The second step of the proof is the existence and uniqueness of the pressure field $p^{\star}$ such that the couple $(\boldsymbol{u}^{\star}, p^{\star})$ solves Problem 17. In order to find $p^{\star}$ we define the following continuous linear functional

$$\langle L, \boldsymbol{v}\rangle \overset{def}{=} 2\mu\left(\boldsymbol{\epsilon}(\boldsymbol{u}^{\star}),\boldsymbol{\epsilon}(\boldsymbol{v})\right)_{0,\Omega} - \rho(\boldsymbol{b},\boldsymbol{v})_{0,\Omega} \qquad \forall \boldsymbol{v} \in H^1_0(\Omega)^d, \tag{C.7}$$

then

$$\langle L, \boldsymbol{v}\rangle = 0 \quad \forall \boldsymbol{v} \in \boldsymbol{V}_{\mathrm{div}}. \tag{C.8}$$

By Lemma C.1.1 this last property is equivalent to the existence and uniqueness of $p^{\star} \in L^2_0(\Omega)$ such that $2\mu\left(\boldsymbol{\epsilon}(\boldsymbol{u}^{\star}),\boldsymbol{\epsilon}(\boldsymbol{v})\right)_{0,\Omega} - \rho(\boldsymbol{b},\boldsymbol{v})_{0,\Omega} = \left(p^{\star}, \operatorname{div}\boldsymbol{v}\right)_{0,\Omega} \quad \forall \boldsymbol{v} \in H^1_0(\Omega)^d$. □

The regularity question for the solution of Stokes Problem 17 will be exploited in the analysis of convergence for the coupled problem, for this reason we report here the following result [48, Lemma 4.17]

**Lemma C.1.3.** *Let $(\boldsymbol{u}, p)$ solution of Stokes Problem 17 and suppose that*

- $\boldsymbol{d = 2}$ *$\Omega$ is a convex polygonal domain or*

- $\boldsymbol{d = 2,3}$ *$\Omega$ is a domain of class $C^{1,1}$, namely its boundary is "locally the graph of a function $C^{1,1}$",*

*then there exist $C > 0$ such that*

$$(\boldsymbol{u}, p) \in \left(H^2(\Omega)^d \cap H^1_0(\Omega)^d\right) \times \left(H^1(\Omega) \cap L^2_0(\Omega)\right) \quad \textit{and} \quad \|\boldsymbol{u}\|_{2,\Omega} + \|p\|_{1,\Omega} \le C\|\boldsymbol{b}\|_{0,\Omega}. \tag{C.9}$$

### C.1.2 Finite elements for Stokes Problem

The finite element approximation of solution to Stokes problem is obtained considering a family of regular meshes $\{\mathscr{T}_h\}_{h>0}$ of the domain $\overline{\Omega}$. Associated to the family of meshes we consider two families of finite element spaces $\{\boldsymbol{V}_h\}_h \subset H_0^1(\Omega)^d$ and $\{Q_h\}_h \subset L_0^2(\Omega)$ approximating the continuous spaces $H_0^1(\Omega)^d$ and $L_0^2(\Omega)$, respectively, as the family parameter $h > 0$ goes to zero. For each mesh of the family we have to solve the finite dimensional problem

**Problem 18.** *Find $\boldsymbol{u}_h \in \boldsymbol{V}_h$ and $p_h \in Q_h$ such that*

$$2\mu(\boldsymbol{\epsilon}(\boldsymbol{u}_h),\boldsymbol{\epsilon}(\boldsymbol{v}_h))_{0,\Omega} - (p_h, \operatorname{div} \boldsymbol{v}_h)_{0,\Omega} = \rho(\boldsymbol{b},\boldsymbol{v}_h)_{0,\Omega} \qquad \forall \boldsymbol{v}_h \in \boldsymbol{V}_h, \tag{C.10}$$

$$(q_h, \operatorname{div} \boldsymbol{u}_h)_{0,\Omega} = 0 \qquad \forall q_h \in Q_h. \tag{C.11}$$

Brezzi's Theorem applies to this problem; in order to express the result in a more concise way, we introduce the following notation

$$a(\boldsymbol{u}_h,\boldsymbol{v}_h) \stackrel{def}{=} 2\mu(\boldsymbol{\epsilon}(\boldsymbol{u}_h),\boldsymbol{\epsilon}(\boldsymbol{v}_h))_{0,\Omega}, \qquad \forall(\boldsymbol{u}_h,\boldsymbol{v}_h) \in \boldsymbol{V}_h \times \boldsymbol{V_h} \tag{C.12}$$

$$b(\boldsymbol{u}_h,p_h) \stackrel{def}{=} -(p_h, \operatorname{div} \boldsymbol{u}_h)_{0,\Omega} \qquad \forall(\boldsymbol{u}_h,p_h) \in \boldsymbol{V}_h \times Q_h \tag{C.13}$$

then the result is the following

**Theorem C.1.4.** ***(Brezzi)*** *If the bilinear forms $a(\cdot,\cdot)$ and $b(\cdot,\cdot)$ are continuous on $\boldsymbol{V}_h \times \boldsymbol{V}_h$ and $\boldsymbol{V}_h \times Q_h$ respectively and*

$$\exists \alpha_h > 0: \quad a(\boldsymbol{z}_h,\boldsymbol{z}_h) \geq \alpha_h \|\boldsymbol{z}_h\|_{1,\Omega} \quad \forall \boldsymbol{z}_h \in \boldsymbol{Z}_h = \left\{\boldsymbol{v}_h \in \boldsymbol{V}_h : b(\boldsymbol{v}_h,q_h) = 0 \quad \forall q_h \in Q_h\right\}, \tag{C.14}$$

$$\exists \beta_h > 0: \qquad \inf_{0 \neq q_h \in Q_h} \sup_{0 \neq \boldsymbol{v}_h \in \boldsymbol{V}_h} \frac{b(\boldsymbol{v}_h,q_h)}{\|\boldsymbol{v}_h\|_{1,\Omega}\|q_h\|_{0,\Omega}} \geq \beta_h. \tag{C.15}$$

*Problem* (18) *has a unique solution. If* (C.14) *and* (C.15) *are verified for all spaces in the family and* $\inf_h \alpha_h \geq \alpha_0 > 0$ *and* $\inf_h \beta_h \geq \beta_0 > 0$, *supposing that $\boldsymbol{u}$ and $p$ are the solution of the continuous Problem 17, we have the following estimate*

$$\|\boldsymbol{u} - \boldsymbol{u}_h\|_{1,\Omega} + \|p - p_h\|_{0,\Omega} \lesssim \inf_{(\boldsymbol{v}_h,q_h) \in \boldsymbol{V}_h \times Q_h} \left(\|\boldsymbol{u} - \boldsymbol{v}_h\|_{1,\Omega} + \|p - q_h\|_{0,\Omega}\right). \tag{C.16}$$

Theorem C.1.4, gives sufficient conditions for the existence, uniqueness and stability of the solution to Stokes problem. In our particular case, the bilinear form $a(\cdot,\cdot)$ is coercive on the entire space $H_0^1(\Omega)^d$, then, since $\boldsymbol{V}_h \subset H_0^1(\Omega)^d$, condition (C.14) is automatically satisfied uniformly in $h$; condition (C.15) poses a limitation on the choice of the finite element spaces $\boldsymbol{V}_h$ and $Q_h$ that can not be arbitrary. Couples of finite element spaces $\boldsymbol{V}_h$ and $Q_h$ that satisfy (C.15) are said *inf-sup stables*.

**Stabilization** Condition (C.15) is quite restrictive and in practice, very convenient couples of spaces are not allowed. To overcome this inconvenience, given a particular couple of finite element spaces, it is possible to modify the discrete variational problem introducing a perturbation term for the bilinear form $b(\cdot,\cdot)$ that vanishes when the family parameter $h$ goes to zero and such that the choice of finite element spaces becames stable for the perturbed bilinear form.

**$P_1 - P_1$ Elements** In this work, we use extensively piecewise linear and globally continuous approximations for fluid velocity and pressure. Hence, considering the following spaces of functions associated to the mesh $\mathcal{T}_h$

$$X_h = \{v_h \in C^0(\overline{\Omega}) : v_h\big|_K \in \mathbb{P}_1^d, \forall K \in \mathcal{T}_h\},$$

we are lead to the following finite element spaces

$$\boldsymbol{V}_h \stackrel{def}{=} X_h^d \cap H_0^1(\Omega)^d, \qquad Q_h \stackrel{def}{=} X_h \cap L_0^2(\Omega). \tag{C.17}$$

It easy to see that the couple $(\boldsymbol{V}_h, \boldsymbol{Q}_h)$ is not inf-sup stable, in fact, if we consider $\Omega$ to be the hexagon showed in figure, and the function $\overline{q}_h$ illustrated in Figure C.1.

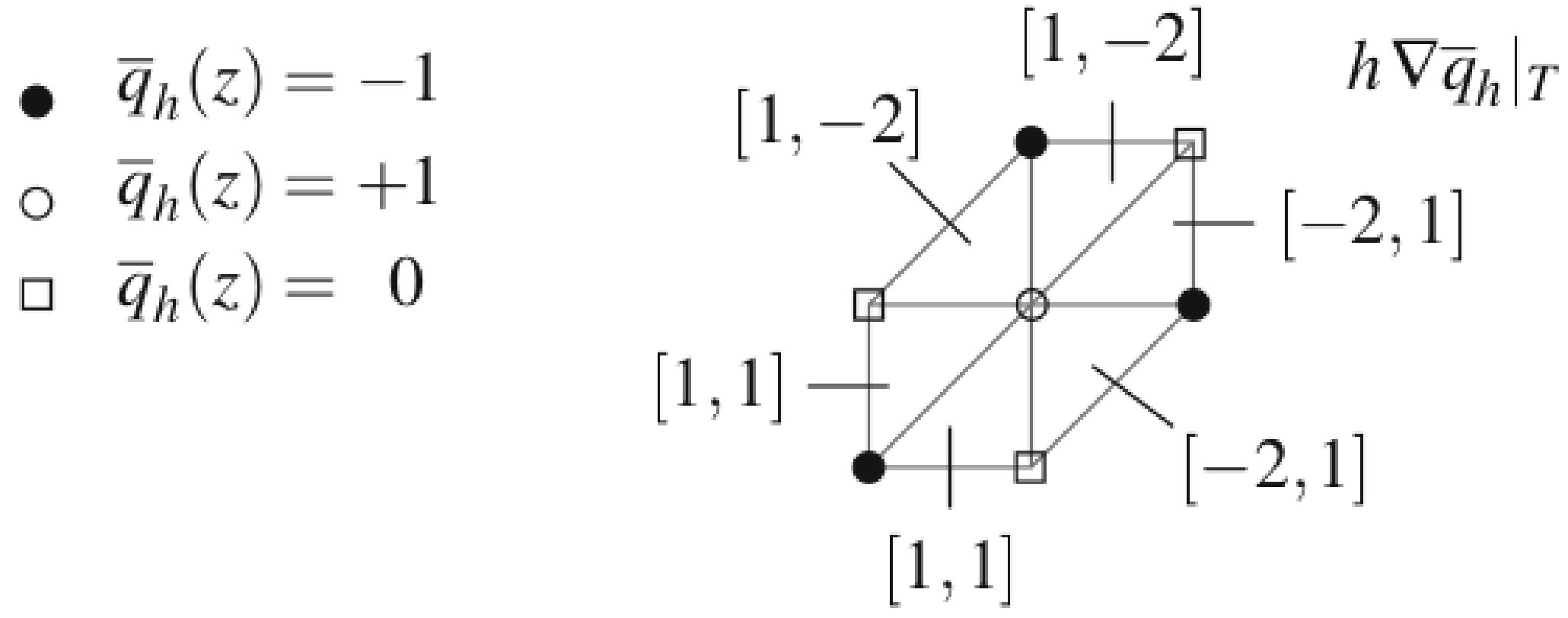

Figure C.1: $P_1 - P_1$ Spurious pressure mode

Let $\{\boldsymbol{e}_j\}_{j=1}^2$ be the canonical basis of $\mathbb{R}^2$, then for any $\boldsymbol{v}_h \in \boldsymbol{V}_h$, we can write $\boldsymbol{v}_h = \alpha_1\phi_h\boldsymbol{e}_1 + \alpha_2\phi_h\boldsymbol{e}_1$, computing the $L^2$ scalar product

$$h\left(\operatorname{div}\boldsymbol{v}_h, \overline{q}_h\right)_{0,\Omega} = -h\left(\boldsymbol{v}_h, \nabla\overline{q}_h\right)_{0,\Omega} = -h\left(\alpha_1\phi_h \sum_{Triangles}(\nabla\overline{q}_h)_1 + \alpha_2\phi_h \sum_{Triangles}(\nabla\overline{q}_h)_2\right) = 0. \tag{C.18}$$

Namely, there exists a non zero pressure field that is orthogonal to the divergence of all velocity fields, then we can not satisfy the discrete inf-sup condition. Even if the $\boldsymbol{P_1} - \boldsymbol{P_1}$ finite elements are not stables it is possible to introduce a stabilization term; that possibility is illustrated by the following well known generalized inf-sup relation that holds for every couple of spaces of continuous piecewise polynomial functions

**Lemma C.1.5.** *Let $Q_h \subset C^0(\overline{\Omega})$, then there exists $\beta' > 0$ and $c_p > 0$ such that for all $q_h \in Q_h$ it holds*

$$\sup_{0 \neq \boldsymbol{v}_h \in \boldsymbol{V}_h} \frac{b(\boldsymbol{v}_h, q_h)}{\|\boldsymbol{v}_h\|_{1,\Omega}} \geq \beta'\|q_h\|_{0,\Omega} - c_p\left(\sum_{K \in \mathcal{T}_h} h_K^2 \|\nabla q_h\|_{0,K}^2\right)^{\frac{1}{2}}. \tag{C.19}$$

*Proof.* The proof of this result is not reported here, it can be found, for example, in [21, Lemma 8.5.1]. □

**Pressure stabilization of $P_1$-$P_1$ elements** Lemma C.1.5 motivates the introduction of a stabilizing term to compensate for the failure of the inf-sup condition. Let $s_h : Q_h \times Q_h \to \mathbb{R}$ be an appropriately chosen bilinear form, we consider the *stabilized discrete Stokes problem*

**Problem 19.** *Find $(\boldsymbol{u}_h, p_h) \in \boldsymbol{V}_h \times Q_h$ such that*

$$a(\boldsymbol{u}_h, \boldsymbol{v}_h) + b(p_h, \boldsymbol{v}_h) = (\boldsymbol{b}, \boldsymbol{v}_h)_{0,\Omega} \qquad \forall \boldsymbol{v} \in \boldsymbol{V}_h, \tag{C.20}$$

$$b(q_h, \boldsymbol{u}_h) - s_h(p_h, q_h) = 0 \qquad \forall q_h \in Q_h. \tag{C.21}$$

As anticipated the bilinear form $s_h(\cdot,\cdot)$ has to be appropriately chosen, in the sense that it has to control the negative term in the perturbed inf-sup condition (C.19). In our work we choose a **pressure stabilization term** introduced by Brezzi and Pitkaranta [29], namely

$$s_h(p_h, q_h) \stackrel{def}{=} \gamma \sum_{K \in \mathcal{T}_h} h_K^2 \int_K \nabla p_h \cdot \nabla q_h d\boldsymbol{x}, \tag{C.22}$$

where $\gamma > 0$ is a constant to be chosen appropriately. We introduce also the notation

$$|q_h|_{s_h}^2 \stackrel{def}{=} s_h(q_h, q_h) \qquad \forall q_h \in Q_h. \tag{C.23}$$

We observe that if the *family of triangulations is quasi-uniform* (see Definition 6), we have, for some $C_1 > 0$ and $C_2 > 0$ independent on h

$$C_1 h^2 \|\nabla q_h\|_{0,\Omega}^2 \geq |q_h|_{s_h}^2 \geq C_2 h^2 \|\nabla q_h\|_{0,\Omega}^2 \qquad \forall q_h \in Q_h. \tag{C.24}$$

With the introduction of the stabilization term (C.22), Stokes Problem 19 is well posed.

**Proposition C.1.6.** *Let $\{\mathcal{T}_h\}_{1 \geq h > 0}$ a regular triangulation family, then Problem 19 with stabilization term* (C.22) *admits a unique solution.*

*Proof.* The existence and uniqueness of solution to the Stabilized Stokes Problem 19 is a consequence of the generalized Lax-Milgram Theorem A.2.1. The cited result will be applied to the bilinear form defined in the following

$$a_h^{\mathrm{f}}((\boldsymbol{u}_h, p_h);(\boldsymbol{v}_h, q_h)) \stackrel{def}{=} 2\mu(\boldsymbol{\epsilon}(\boldsymbol{u}_h), \boldsymbol{\epsilon}(\boldsymbol{v}_h))_{0,\Omega} - (\operatorname{div} \boldsymbol{v}_h, p_h)_{0,\Omega} + (\operatorname{div} \boldsymbol{u}_h, q_h)_{0,\Omega} + s_h(p_h, q_h)$$
$$\forall \boldsymbol{u}_h, \boldsymbol{v}_h \in \boldsymbol{V}_h, \forall p_h, q_h \in Q_h. \tag{C.25}$$

Given $(\boldsymbol{u}_h, p_h) \in \boldsymbol{V}_h \times Q_h$, according to the generalized inf-sup inequality (C.19), there exist $\tilde{\boldsymbol{w}}_h \in \boldsymbol{V}_h$ such that

$$\frac{(div \tilde{\boldsymbol{w}}_h, p_h)_{0,\Omega}}{\|\tilde{\boldsymbol{w}}_h\|_{1,\Omega}} \geq \beta' \|p_h\|_{0,\Omega} - c_p \left( \sum_{K \in \mathcal{T}_h} h_K^2 \|\nabla p_h\|_{0,K}^2 \right)^{\frac{1}{2}}, \qquad \|\tilde{\boldsymbol{w}}_h\|_{1,\Omega} = \|p_h\|_{0,\Omega}. \tag{C.26}$$

Then, let $\delta > 0$, taking $(\boldsymbol{v}_h, q_h) = (\boldsymbol{u}_h - \delta \tilde{\boldsymbol{w}}_h, p_h) \in \boldsymbol{V}_h \times Q_h$ we have

$$a_h^{\mathrm{f}}((\boldsymbol{u}_h, p_h);(\boldsymbol{v}_h, q_h)) = 2\mu(\boldsymbol{\epsilon}(\boldsymbol{u}_h), \boldsymbol{\epsilon}(\boldsymbol{u}_h))_{0,\Omega} - 2\mu\delta(\boldsymbol{\epsilon}(\boldsymbol{u}_h), \boldsymbol{\epsilon}(\tilde{\boldsymbol{w}}_h))_{0,\Omega} + \delta(\operatorname{div} \tilde{\boldsymbol{w}}_h, p_h)_{0,\Omega} + |p_h|_{s_h}^2$$

$$\geq \underbrace{2\mu\|\boldsymbol{\epsilon}(\boldsymbol{u}_h)\|_{0,\Omega}^2}_{A} \underbrace{-2\mu\delta\|\boldsymbol{\epsilon}(\boldsymbol{u}_h)\|_{0,\Omega}\|\boldsymbol{\epsilon}(\tilde{\boldsymbol{w}}_h)\|_{0,\Omega}}_{B} + \underbrace{\delta(\operatorname{div}\tilde{\boldsymbol{w}}_h, p_h)_{0,\Omega} + |p_h|_{s_h}^2}_{C}, \quad (C.27)$$

using a combination of Korn's and Poincaré's inequalities , since $\boldsymbol{u}_h \in H_0^1(\Omega)^d$ we can write

$$A = 2\mu\|\boldsymbol{\epsilon}(\boldsymbol{u}_h)\|_{0,\Omega}^2 \geq 2\mu C\|\nabla\boldsymbol{u}_h\|_{0,\Omega}^2 \geq 2\mu C_K\|\boldsymbol{u}_h\|_{1,\Omega}^2, \quad (C.28)$$

The term $B$ is bounded using Young inequality and considering that $\|\tilde{\boldsymbol{w}}_h\|_{1,\Omega} = \|p_h\|_{0,\Omega}$

$$\begin{aligned} B &= -2\mu\delta\|\boldsymbol{\epsilon}(\boldsymbol{u}_h)\|_{0,\Omega}\|\boldsymbol{\epsilon}(\tilde{\boldsymbol{w}}_h)\|_{0,\Omega} \\ &\geq -\frac{\mu\delta}{\eta}\|\boldsymbol{\epsilon}(\boldsymbol{u}_h)\|_{0,\Omega}^2 - \mu\delta\eta\|\boldsymbol{\epsilon}(\tilde{\boldsymbol{w}}_h)\|_{0,\Omega}^2 \\ &\geq -C_K^2\frac{\mu\delta}{\eta}\|\boldsymbol{u}_h\|_{1,\Omega}^2 - C_K^2\mu\delta\eta\|\tilde{\boldsymbol{w}}_h\|_{1,\Omega}^2 \\ &= -C_K^2\frac{\mu\delta}{\eta}\|\boldsymbol{u}_h\|_{1,\Omega}^2 - C_K^2\mu\delta\eta\|p_h\|_{0,\Omega}^2, \end{aligned} \quad (C.29)$$

using relations in (C.26) and considering that $\|\tilde{\boldsymbol{w}}_h\|_{1,\Omega} = \|p_h\|_{0,\Omega}$, we obtain

$$\begin{aligned} C &= \delta(\operatorname{div}\tilde{\boldsymbol{w}}_h, p_h)_{0,\Omega} + |p_h|_{s_h}^2 \\ &\geq \delta\beta'\|p_h\|_{0,\Omega}\|\tilde{\boldsymbol{w}}_h\|_{1,\Omega} - \delta c_p\left(\sum_{K\in\mathscr{T}_h} h_K^2\|\nabla p_h\|_{0,K}^2\right)^{\frac{1}{2}}\|\tilde{\boldsymbol{w}}_h\|_{1,\Omega} + |p_h|_{s_h}^2 \\ &= \delta\beta'\|p_h\|_{0,\Omega}^2 - \delta c_p\left(\sum_{K\in\mathscr{T}_h} h_K^2\|\nabla p_h\|_{0,K}^2\right)^{\frac{1}{2}}\|p_h\|_{0,\Omega} + \gamma\sum_{K\in\mathscr{T}_h} h_K^2\|\nabla p_h\|_{0,K}^2 \\ &\geq \frac{\delta\beta'}{2}\|p_h\|_{0,\Omega}^2 + \left(\gamma - \frac{\delta c_p^2}{2\beta'}\right)\sum_{K\in\mathscr{T}_h} h_K^2\|\nabla p_h\|_{0,K}^2, \end{aligned} \quad (C.30)$$

In summary, using expressions (C.28), (C.29) and (C.30) in (C.27) we obtain

$$\begin{aligned} a_h^{\mathrm{f}}((\boldsymbol{u}_h, p_h);(\boldsymbol{v}_h, q_h)) &\geq \underbrace{\mu C_K^2\left(2 - \frac{\delta}{\eta}\right)}_{A(\delta,\eta)}\|\boldsymbol{u}_h\|_{1,\Omega}^2 \\ &\quad + \delta\underbrace{\left(\frac{\beta'}{2} - \eta C_K^2\mu\right)}_{B(\eta)}\|p_h\|_{0,\Omega}^2 + \underbrace{\left(\gamma - \frac{\delta c_p^2}{2\beta'}\right)}_{C(\delta)}\sum_{K\in\mathscr{T}_h} h_K^2\|\nabla p_h\|_{0,K}^2. \end{aligned} \quad (C.31)$$

If we choose $\delta > 0, \gamma > 0, \eta > 0$, such that $C(\delta) > 0$, $A(\delta,\eta) \geq \alpha_0 > 0$ and $B(\eta) \geq \alpha_0 > 0$, for some $\alpha_0 > 0$, we have

$$a_h^{\mathrm{f}}((\boldsymbol{u}_h, p_h);(\boldsymbol{v}_h, q_h)) \geq \alpha_0\left(\|\boldsymbol{u}_h\|_{1,\Omega}^2 + \|p_h\|_{0,\Omega}^2\right). \quad (C.32)$$

This inequality is fulfilled if

$$\eta = \frac{\beta'}{4\mu C_K^2}, \qquad \delta = \frac{1}{2}\min\left\{\frac{2\beta'\gamma}{c_p^2}, \frac{\beta'}{4\mu C_K^2}\right\}$$

Moreover we observe that for any $(\boldsymbol{u}_h, q_h) \in \boldsymbol{V}_h \times Q_h$ the couple $(\boldsymbol{v}_h, p_h) = (\boldsymbol{u}_h - \delta\tilde{\boldsymbol{w}}_h, p_h)$ satisfies

$$\begin{aligned}\|\boldsymbol{v}_h\|_{1,\Omega}^2 + \|q_h\|_{0,\Omega}^2 &\le 2\|\boldsymbol{u}_h\|_{1,\Omega}^2 + 2\delta^2\|\tilde{\boldsymbol{w}}_h\|_{1,\Omega}^2 + \|p_h\|_{0,\Omega}^2 \\ &= 2\|\boldsymbol{u}_h\|_{1,\Omega}^2 + \left(2\delta^2+1\right)\|p_h\|_{0,\Omega}^2 \le C(\|\boldsymbol{u}_h\|_{1,\Omega}^2 + \|p_h\|_{0,\Omega}^2), \end{aligned} \quad \text{(C.33)}$$

then, using (C.33) in (C.32) we arrive at the following result

$$\forall(\boldsymbol{u}_h, p_h) \in \boldsymbol{V}_h \times Q_h, \quad \exists(\boldsymbol{v}_h, q_h) \in \boldsymbol{V}_h \times Q_h, \qquad \frac{a_h^{\mathrm{f}}((\boldsymbol{u}_h, p_h);(\boldsymbol{v}_h, q_h))}{\sqrt{\|\boldsymbol{v}_h\|_{1,\Omega}^2 + \|q_h\|_{0,\Omega}^2}} \ge \alpha_1\sqrt{\|\boldsymbol{u}_h\|_{1,\Omega}^2 + \|p_h\|_{0,\Omega}^2}. \quad \text{(C.34)}$$

for some $\alpha_1 > 0$. □

**Theorem C.1.7.** *Assume that the conditions of Proposition C.1.6 are satisfied. Let* $(\boldsymbol{u}, p) \in H_0^1(\Omega)^d \times L_0^2(\Omega)$ *be the solution of Problem 17, and* $(\boldsymbol{u}_h, p_h) \in \boldsymbol{V}_h \times Q_h$ *be the solution of Problem 19, then we have*

$$\|\boldsymbol{u} - \boldsymbol{u}_h\|_{1,\Omega} + \|p - p_h\|_{0,\Omega} \le C\left[\inf_{\boldsymbol{v}_h \in \boldsymbol{V}_h}\|\boldsymbol{u} - \boldsymbol{v}_h\|_{1,\Omega} + \inf_{q_h \in Q_h}\left(\|p - q_h\|_{0,\Omega} + |q_h|_{s_h}\right)\right] \quad \text{(C.35)}$$

*Proof.* Let $(\boldsymbol{u}_h, p_h) \in \boldsymbol{V}_h \times Q_h$ be the solution of Problem 19, in order to prove the error inequality (C.35), we consider any $\left(\tilde{\boldsymbol{u}}_h, \tilde{q}_h\right) \in \boldsymbol{V}_h \times Q_h$, then, using the inf-sup condition (C.34), we obtain the existence of $(\boldsymbol{v}_h, q_h) \in \boldsymbol{V}_h \times Q_h$ such that $\|\boldsymbol{v}_h\|_{1,\Omega}^2 + \|q_h\|_{0,\Omega}^2 = 1$ and

$$a^{\mathrm{f}}((\tilde{\boldsymbol{u}}_h - \boldsymbol{u}_h, \tilde{p}_h - p_h);(\boldsymbol{v}_h, q_h)) + s_h(\tilde{p}_h - p_h, q_h) \ \ge\ \alpha_1\sqrt{\|\tilde{\boldsymbol{u}}_h - \boldsymbol{u}_h\|_{1,\Omega}^2 + \|\tilde{p}_h - p_h\|_{0,\Sigma}^2}. \quad \text{(C.36)}$$

Moreover we have, for all $(\tilde{\boldsymbol{u}}_h, \tilde{p}_h) \in \boldsymbol{V}_h \times Q_h$

$$\begin{aligned} &a^{\mathrm{f}}((\tilde{\boldsymbol{u}}_h - \boldsymbol{u}, \tilde{p}_h - p);(\boldsymbol{v}_h, q_h)) + s_h(\tilde{p}_h, q_h) \\ &\quad = 2\mu(\epsilon(\tilde{\boldsymbol{u}}_h - \boldsymbol{u}), \epsilon(\boldsymbol{v}_h))_{0,\Omega} - (\operatorname{div}\boldsymbol{v}_h, \tilde{p}_h - p)_{0,\Omega} + (\operatorname{div}(\tilde{\boldsymbol{u}}_h - \boldsymbol{u}), q_h)_{0,\Omega} + s_h(\tilde{p}_h, q_h) \\ &\quad\quad \lesssim \|\tilde{\boldsymbol{u}}_h - \boldsymbol{u}\|_{1,\Omega}\|\boldsymbol{v}_h\|_{1,\Omega} + \|\tilde{p}_h - p\|_{0,\Omega}\|\operatorname{div}\boldsymbol{v}_h\|_{0,\Omega} + \|\operatorname{div}(\tilde{\boldsymbol{u}}_h - \boldsymbol{u})\|_{0,\Omega}\|q_h\|_{0,\Omega} \\ &\qquad\qquad + \gamma\sum_{K \in \mathscr{T}_h} h_K^2\int_K \nabla\tilde{p}_h\nabla q_h d\boldsymbol{x}, \end{aligned}$$

since $\|\boldsymbol{v}_h\|_{1,\Omega}^2 + \|q_h\|_{0,\Omega}^2 = 1$, and

$$\begin{aligned}\gamma\sum_{K \in \mathscr{T}_h} h_K^2\int_K \nabla\tilde{p}_h\nabla q_h d\boldsymbol{x} &\le C\gamma\left(\sum_{K \in \mathscr{T}_h} h_K^2\|\nabla\tilde{p}_h\|_{0,K}^2\right)^{\frac{1}{2}}\left(\sum_{K \in \mathscr{T}_h} h_K^2\|\nabla q_h\|_{0,\Omega}^2\right)^{\frac{1}{2}} \\ &\le C|\tilde{p}_h|_{s_h}\|q_h\|_{0,\Omega} \le C|\tilde{p}_h|_{s_h}, \end{aligned}$$

we have

$$a^{\mathrm{f}}((\tilde{\boldsymbol{u}}_h - \boldsymbol{u}, \tilde{p}_h - p);(\boldsymbol{v}_h, q_h)) + s_h(\tilde{p}_h, q_h) \lesssim \|\tilde{\boldsymbol{u}}_h - \boldsymbol{u}\|_{1,\Omega} + \|\tilde{p}_h - p\|_{0,\Omega} + C|\tilde{p}_h|_{s_h}. \quad \text{(C.37)}$$

Considering (C.36) and (C.37) we obtain, for all $(\tilde{\boldsymbol{u}}_h, \tilde{p}_h) \in \boldsymbol{V}_h \times Q_h$

$$\begin{aligned}
&\|\boldsymbol{u}-\boldsymbol{u}_h\|_{1,\Omega} + \|p-p_h\|_{0,\Omega} \\
&\qquad \leq \|\boldsymbol{u}-\tilde{\boldsymbol{u}}_h\|_{1,\Omega} + \|\tilde{\boldsymbol{u}}_h-\boldsymbol{u}_h\|_{1,\Omega} + \|p-\tilde{p}_h\|_{0,\Omega} + \|\tilde{p}_h-p_h\|_{0,\Sigma} \\
&\qquad\qquad \lesssim \|\boldsymbol{u}-\tilde{\boldsymbol{u}}_h\|_{1,\Omega} + \|p-\tilde{p}_h\|_{0,\Omega} + C|\tilde{p}_h|_{s_h}. \qquad \text{(C.38)}
\end{aligned}$$

Hence

$$\|\boldsymbol{u}-\boldsymbol{u}_h\|_{1,\Omega} + \|p-p_h\|_{0,\Omega} \lesssim \inf_{\tilde{\boldsymbol{u}}_h \in \boldsymbol{V}_h} \|\boldsymbol{u}-\tilde{\boldsymbol{u}}_h\|_{1,\Omega} + \inf_{\tilde{p}_h \in Q_h} \left(\|p-\tilde{p}_h\|_{0,\Omega} + C|\tilde{p}_h|_{s_h}\right) \qquad \text{(C.39)}$$

□

We will use also the $L^2$ estimate for the error in the fluid velocities.

**Lemma C.1.8.** *Let $\Omega$ as in the hypothesis of Lemma C.1.3. If $(\boldsymbol{u}_h, p_h) \in \boldsymbol{V}_h \times Q_h$ is the solution of Problem 18, and $(\boldsymbol{u}, p) \in \boldsymbol{V} \times Q$ is the solution Problem 17, then*

$$\|\boldsymbol{u}-\boldsymbol{u}_h\|_{0,\Omega} \leq Ch\left(\|\boldsymbol{u}-\boldsymbol{u}_h\|_{1,\Omega} + \|p-p_h\|_{0,\Omega} + |p_h|_{s_h}\right). \qquad \text{(C.40)}$$

*Proof.* We consider the dual Stokes Problem: find $(\boldsymbol{u}^\star, p^\star) \in \boldsymbol{V} \times Q$ such that

$$\left(\epsilon(\boldsymbol{v}), \epsilon(\boldsymbol{u}^\star)\right) - \left(\operatorname{div}\boldsymbol{v}, p^\star\right) = (\boldsymbol{u}-\boldsymbol{u}_h, \boldsymbol{v}) \qquad \forall \boldsymbol{v} \in \boldsymbol{V} \qquad \text{(C.41)}$$

$$\left(\operatorname{div}\boldsymbol{u}^\star, q\right) = 0 \qquad \forall q \in Q. \qquad \text{(C.42)}$$

Since $\Omega$ satisfies the hypotheses of Lemma C.1.3, we have the following regularity estimate for the solution of the dual Stokes Problem

$$\|\boldsymbol{u}^\star\|_{2,\Omega} + \|p^\star\|_{1,\Omega} \lesssim \|\boldsymbol{u}-\boldsymbol{u}_h\|_{0,\Omega}. \qquad \text{(C.43)}$$

Let $\boldsymbol{u}_h^\star \in \boldsymbol{V}_h$ be the Lagrange interpolation of $\boldsymbol{u}^\star$ and let $p_h^\star \in Q_h$ be the Clément interpolation of $p^\star$, then equation (C.41) yields

$$\begin{aligned}
\|\boldsymbol{u}-\boldsymbol{u}_h\|_{0,\Omega}^2 &= \left(\epsilon(\boldsymbol{u}-\boldsymbol{u}_h), \epsilon(\boldsymbol{u}^\star)\right)_{0,\Omega} - \left(\operatorname{div}(\boldsymbol{u}-\boldsymbol{u}_h), , p^\star\right)_{0,\Omega} \\
&= \underbrace{\left(\epsilon(\boldsymbol{u}-\boldsymbol{u}_h), \epsilon(\boldsymbol{u}^\star-\boldsymbol{u}_h^\star)\right)_{0,\Omega} - \left(\operatorname{div}(\boldsymbol{u}-\boldsymbol{u}_h,), p^\star-p_h^\star\right)_{0,\Omega}}_{A} \\
&\qquad + \underbrace{\left(\epsilon(\boldsymbol{u}-\boldsymbol{u}_h), \epsilon(\boldsymbol{u}_h^\star)\right)_{0,\Omega}}_{B} \underbrace{- \left(\operatorname{div}(\boldsymbol{u}-\boldsymbol{u}_h), p_h^\star\right)_{0,\Omega}}_{C} \qquad \text{(C.44)}
\end{aligned}$$

We estimate in the following the various terms of relation (C.44) using the regularity (C.43) and the estimates for Lagrange and Clément interpolations (see (A.30) and (A.38) respectively), for $m = 0, 1$

$$\|\boldsymbol{u}^\star-\boldsymbol{u}_h^\star\|_{m,\Omega} \leq Ch^{2-m}\|\boldsymbol{u}^\star\|_{2,\Omega} \leq Ch^{2-m}\|\boldsymbol{u}-\boldsymbol{u}_h\|_{0,\Omega} \qquad \text{(C.45)}$$

$$\|p^\star-p_h^\star\|_{0,\Omega} \leq Ch\|p^\star\|_{1,\Omega} \leq Ch\|\boldsymbol{u}-\boldsymbol{u}_h\|_{0,\Omega}. \qquad \text{(C.46)}$$

Then we have, using (C.45) and (C.46)

$$\begin{aligned} A &\le \left| \left(\epsilon(\boldsymbol{u}-\boldsymbol{u}_h), \epsilon(\boldsymbol{u}^\star - \boldsymbol{u}_h^\star)\right)_{0,\Omega} - \left(\operatorname{div}(\boldsymbol{u}-\boldsymbol{u}_h), p^\star - p_h^\star\right)_{0,\Omega} \right| \\ &\le C\|\boldsymbol{u}-\boldsymbol{u}_h\|_{1,\Omega}\left(\|\boldsymbol{u}^\star - \boldsymbol{u}_h^\star\|_{1,\Omega} + \|p^\star - p_h^\star\|_{0,\Omega}\right) \\ &\le Ch\|\boldsymbol{u}-\boldsymbol{u}_h\|_{1,\Omega}\left(\|\boldsymbol{u}^\star\|_{2,\Omega} + \|p^\star\|_{1,\Omega}\right) \\ &\le Ch\|\boldsymbol{u}-\boldsymbol{u}_h\|_{1,\Omega}\|\boldsymbol{u}-\boldsymbol{u}_h\|_{0,\Omega}, \end{aligned} \tag{C.47}$$

the term $B$ is controlled as follows using (C.45)

$$\begin{aligned} B &= \left(\epsilon(\boldsymbol{u}-\boldsymbol{u}_h), \epsilon(\boldsymbol{u}_h^\star)\right)_{0,\Omega} = -\left(\operatorname{div}\boldsymbol{u}_h^\star, p - p_h\right)_{0,\Omega} = \left(\operatorname{div}\left(\boldsymbol{u}^\star - \boldsymbol{u}_h^\star\right), p - p_h\right)_{0,\Omega} \\ &\le C\|\boldsymbol{u}^\star - \boldsymbol{u}_h^\star\|_{1,\Omega}\|p - p_h\|_{0,\Omega} \le Ch\|\boldsymbol{u}^\star\|_{2,\Omega}\|p - p_h\|_{0,\Omega} \\ &\le Ch\|\boldsymbol{u}-\boldsymbol{u}_h\|_{0,\Omega}\|p - p_h\|_{0,\Omega}, \end{aligned} \tag{C.48}$$

moreover, using the inverse inequality (A.42) on each element $K \in \mathscr{T}_h$ and (C.46)

$$\begin{aligned} C &= -\left(\operatorname{div}(\boldsymbol{u}-\boldsymbol{u}_h), p_h^\star\right)_{0,\Omega} = -s_h(p_h, p_h^\star) \le \gamma \sum_{K\in\mathscr{T}_h} h_K^2 \int_K \nabla p_h \nabla p_h^\star d\boldsymbol{x} \\ &\le \gamma C \left(\sum_{K\in\mathscr{T}_h} h_K^2 \|\nabla p_h\|_{0,K}^2\right)^{\frac{1}{2}} \left(\sum_{K\in\mathscr{T}_h} h_K^2 \|\nabla p_h^\star\|_{0,\Omega}^2\right)^{\frac{1}{2}} \le Ch|p_h|_{s_h}\|\nabla p_h^\star\|_{0,\Omega} \\ &\le Ch|p_h|_{s_h}\left(\|\nabla p_h^\star - \nabla p^\star\|_{0,\Omega} + \|\nabla p^\star\|_{0,\Omega}\right) \\ &\le Ch|p_h|_{s_h}\|\boldsymbol{u}-\boldsymbol{u}_h\|_{0,\Omega} \end{aligned} \tag{C.49}$$

Hence, using the estimates (C.47), (C.48) and (C.49) in (C.44)

$$\|\boldsymbol{u}-\boldsymbol{u}_h\|_{0,\Omega} \le Ch\left(\|\boldsymbol{u}-\boldsymbol{u}_h\|_{1,\Omega} + \|p - p_h\|_{0,\Omega}\right) + Ch|p_h|_{s_h}. \tag{C.50}$$

□

# Appendix D

# Matrices for the numerical simulations

## D.1 Introduction

In this appendix, we report the computation needed for the implementation of the numerical simulations. The computations are performed on the test case of an elliptical structure that evolves to the rest circular position.

## D.2 Matrix expression of the Algorithms and computation

For the numerical implementation we refer to the matrix expression of the algorithms that heve been presented in (2.47) for what regards the monolithic algorithm, in (4.7) for Algorithm 2 and in (4.9) and (4.10) for Algorithm 3 with the definitions of submatrices given there.

In all the algorithms presented, we have to compute some sub-matrices in order to obtain the system of linear equations. In particular, in all the cases we have to derive expressions for the following matrices

- *fluid matrices*: $M_{\mathrm{f}}$,$K_{\mathrm{f}}$,$B$,$S$;
- *solid matrices*: $M_{\mathrm{s}}$,$K_{\mathrm{s}}$;
- *matrices coupling the fluid and the structure*: $L_{\mathrm{f}}(\boldsymbol{\phi}_h^{n-1})$,$L_{\mathrm{s}}$.

In the following paragraphs we illustrate the construction of these matrices associated to the finite element model.

### D.2.1 Fluid matrices

The matrices related to the fluid equations are the fluid mass matrix $M_{\mathrm{f}}$, the fluid "stiffness matrix" $K_{\mathrm{f}}$, the matrix $B$ associated to the divergence operator and the matrix $S$ associated to the stabilizing term. The computation of these matrices is standard; given the transformation $T_K : \hat{K} \to K$ from the reference element $\hat{K}$ to the generic element $K \in \mathcal{T}_h$ and denoting by $\boldsymbol{J}_K = |\det \nabla T_K|$, we have

$$(M_\mathrm{f})_{ij} = \left(\varphi_j, \varphi_i\right)_{0,\Omega} = \int_\Omega \varphi_j \cdot \varphi_i d\boldsymbol{x} = \sum_{K \in \mathcal{T}_h} \int_K \varphi_j \cdot \varphi_i d\boldsymbol{x}$$
$$= \sum_{K \in \mathcal{T}_h} \int_{\hat{K}} \left[\varphi_j \circ T_K(\hat{\boldsymbol{x}})\right] \cdot \left[\varphi_i \circ T_K(\hat{\boldsymbol{x}})\right] \boldsymbol{J}_K d\hat{\boldsymbol{x}}. \quad \text{(D.1)}$$

In the same way we compute the fluid "stiffness matrix"

$$(K_\mathrm{f})_{ij} = \left(\boldsymbol{\epsilon}\varphi_j, \boldsymbol{\epsilon}\varphi_i\right)_{0,\Omega} = \int_\Omega \boldsymbol{\epsilon}\varphi_j : \boldsymbol{\epsilon}\varphi_i d\boldsymbol{x} = \sum_{K \in \mathcal{T}_h} \int_K \boldsymbol{\epsilon}\varphi_j : \boldsymbol{\epsilon}\varphi_i d\boldsymbol{x}$$
$$= \sum_{K \in \mathcal{T}_h} \int_{\hat{K}} \left[\boldsymbol{\epsilon}\varphi_j \circ T_K(\hat{\boldsymbol{x}})\right] : \left[\boldsymbol{\epsilon}\varphi_i \circ T_K(\hat{\boldsymbol{x}})\right] \boldsymbol{J}_K d\hat{\boldsymbol{x}}. \quad \text{(D.2)}$$

Now we compute the matrix associated to the divergence operator

$$(B)_{ij} = -\left(\operatorname{div}\varphi_j, \psi_i\right)_{0,\Omega} = -\int_\Omega \operatorname{div}\varphi_j \psi_i d\boldsymbol{x} = -\int_{\partial\Omega} \varphi_j \psi_i \boldsymbol{n} d\boldsymbol{\sigma} + \int_\Omega \varphi_j \nabla\psi_i d\boldsymbol{x}$$
$$= \sum_{K \in \mathcal{T}_h} \int_K \varphi_j \nabla\psi_i d\boldsymbol{x} = \sum_{K \in \mathcal{T}_h} \int_{\hat{K}} \left[\varphi_j \circ T_K(\hat{\boldsymbol{x}})\right] \left[\nabla\psi_i \circ T_K(\hat{\boldsymbol{x}})\right] \boldsymbol{J}_K d\hat{\boldsymbol{x}} \quad \text{(D.3)}$$

It remains to compute the matrix associated to the pressure stabilization

$$(S)_{ij} = \sum_{K \in \mathcal{T}_h} h_K^2 \left(\nabla\psi_j, \nabla\psi_i\right)_{0,K} = \sum_{K \in \mathcal{T}_h} h_K^2 \int_K \nabla\psi_j \nabla\psi_i d\boldsymbol{x}$$
$$= \sum_{K \in \mathcal{T}_h} h_K^2 \int_{\hat{K}} \left[\nabla\psi_j \circ T_K(\hat{\boldsymbol{x}})\right] \left[\nabla\psi_i \circ T_K(\hat{\boldsymbol{x}})\right] \boldsymbol{J}_K d\hat{\boldsymbol{x}} \quad \text{(D.4)}$$

### D.2.2 Solid matrices

For the actual computation of the matrices related to the solid problem we need to introduce a parametrization of the reference configuration $\Sigma$. The following computations are performed for the case of a 1-dimensional solid immersed in a 2-dimensional fluid domain, In practice we need only a parametrization of each element (see Figure D.1). In the case of linear elements of extrema $A$ and $B$, this parametrization is

$$\boldsymbol{\alpha}(\varsigma) = \vec{\boldsymbol{A}} + (\vec{\boldsymbol{B}} - \vec{\boldsymbol{A}})\frac{\varsigma}{l} \qquad \text{with} \qquad \varsigma \in [0, l], \quad \text{(D.5)}$$

where $\vec{\boldsymbol{A}}$ and $\vec{\boldsymbol{B}}$, are the position vectors of the points $A$ and $B$ respectively and $l = \|\vec{\boldsymbol{B}} - \vec{\boldsymbol{A}}\|_{\mathbb{R}^d}$ is the length of the element. Using the parametrization introduced in (D.5) for each solid element $K \in \mathcal{S}_h$, solid mass matrix has the following expression

$$
\begin{aligned}
(M_{\mathrm{s}})_{ij} = \left(\chi_j, \chi_i\right)_{0,\Sigma} = \int_\Sigma \chi_j \cdot \chi_i d\boldsymbol{s} = \sum_{K \in \mathscr{S}_h} \int_K \chi_j \cdot \chi_i d\boldsymbol{s} \\
= \sum_{K \in \mathscr{S}_h} \int_0^{l_K} \left[\chi_j \circ \boldsymbol{\alpha}_K(\varsigma)\right] \cdot \left[\chi_i \circ \boldsymbol{\alpha}_K(\varsigma)\right] \left\| \frac{d\boldsymbol{\alpha}_K}{d\varsigma} \right\|_{\mathbb{R}^d} d\varsigma, \quad \text{(D.6)}
\end{aligned}
$$

with the following values of the integrands:

$$
\left\| \frac{d\boldsymbol{\alpha}_K}{d\hat{\varsigma}} \right\|_{\mathbb{R}^d} = 1, \qquad \chi_j \circ \boldsymbol{\alpha}_K(\varsigma) = \begin{cases} \left(1 - \frac{\varsigma}{l_K}, 1 - \frac{\varsigma}{l_K}\right) & \text{shape functions for node "0"} \\ \left(\frac{\varsigma}{l_K}, \frac{\varsigma}{l_K}\right) & \text{shape functions for node "1"}, \end{cases} \quad \text{(D.7)}
$$

where $\chi_j \circ \boldsymbol{\alpha}_K(\varsigma)$ are the shape functions on the element K.

In a similar way we can compute the components of the solid stiffness matrix

$$
\begin{aligned}
(K_{\mathrm{s}})_{ij} = \left(\boldsymbol{\epsilon}(\chi_j), \boldsymbol{\epsilon}(\chi_i)\right)_{0,\Sigma} = \int_\Sigma \boldsymbol{\epsilon}(\chi_j) \cdot \boldsymbol{\epsilon}(\chi_i) d\boldsymbol{s} = \sum_{K \in \mathscr{S}_h} \int_K \boldsymbol{\epsilon}(\chi_j) \cdot \boldsymbol{\epsilon}(\chi_i) d\boldsymbol{s} \\
= \sum_{K \in \mathscr{S}_h} \int_0^{l_K} \left[\boldsymbol{\epsilon}(\chi_j) \circ \boldsymbol{\alpha}_K(\varsigma)\right] \cdot \left[\boldsymbol{\epsilon}(\chi_i) \circ \boldsymbol{\alpha}_K(\varsigma)\right] \left\| \frac{d\boldsymbol{\alpha}_K}{d\varsigma} \right\|_{\mathbb{R}^d} d\varsigma, \quad \text{(D.8)}
\end{aligned}
$$

where the integrands are given in the following

$$
\left\| \frac{d\boldsymbol{\alpha}_K}{d\varsigma} \right\|_{\mathbb{R}^d} = 1, \qquad \boldsymbol{\epsilon}(\chi_j) \circ \boldsymbol{\alpha}_K(\varsigma) = \begin{cases} \left(-\frac{1}{l_K}, -\frac{1}{l_K}\right) & \text{for node "0"} \\ \left(\frac{1}{l_K}, \frac{1}{l_K}\right) & \text{for node "1"}. \end{cases} \quad \text{(D.9)}
$$

### D.2.3 Lagrange multiplier matrix

The computation of the matrix related to the Lagrange multiplier is performed in two stages, in fact we have to consider the matrix $L_{\mathrm{s}}$ and the matrix $L_{\mathrm{f}}$, the first one does not depend on the time, instead the second one depends on time. As in the previous section, for a 1-dimensional solid immersed in a 2-dimensional fluid domain, we use the parametrization of the generic linear element of extrema $A$ and $B$

$$
\boldsymbol{\alpha}(\varsigma) = \vec{\boldsymbol{A}} + (\vec{\boldsymbol{B}} - \vec{\boldsymbol{A}})\frac{\varsigma}{l} \qquad \text{with} \qquad \varsigma \in [0, l]. \quad \text{(D.10)}
$$

Hence the effective computation of the Lagrange multiplier matrix for the solid problem can be performed as follow

$$
(L_{\mathrm{s}})_{ij} = (\boldsymbol{\zeta}_j, \boldsymbol{\chi}_i)_{0,\Sigma} = \int_\Sigma \boldsymbol{\zeta}_j \cdot \boldsymbol{\chi}_i d\boldsymbol{s} = \sum_{K \in \mathscr{S}_h} \int_0^{l_K} \left[\boldsymbol{\zeta}_j \circ \boldsymbol{\alpha}_K(\varsigma)\right] \cdot \left[\boldsymbol{\chi}_i \circ \boldsymbol{\alpha}_K(\varsigma)\right] \left\| \frac{d\boldsymbol{\alpha}_K}{d\varsigma} \right\|_{\mathbb{R}^d} d\varsigma, \quad \text{(D.11)}
$$

hence the matrix $L_s$ is identical to the solid mass matrix $M_s$.

The computation of the Lagrange matrix for the fluid problem is more involved since it requires the evaluation of the fluid velocity on the current configuration of the structure; for the mentioned reason we introduce also the parametrization of each element of extrema $A_t$ and $B_t$ in the current configuration

$$\boldsymbol{\alpha}_{Kt}(\varsigma) = \overrightarrow{\boldsymbol{A}}_t + (\overrightarrow{B}_t - \overrightarrow{A}_t)\frac{\varsigma}{l} \qquad \text{with} \qquad \varsigma \in [0,l], \tag{D.12}$$

If $\boldsymbol{\phi}_t : \boldsymbol{\Sigma} \to \boldsymbol{\Sigma}_t$, represent the transformation from the reference configuration to the current configuration, we have that $\boldsymbol{\alpha}_{Kt}(\varsigma) = \boldsymbol{\phi}_t \circ \boldsymbol{\alpha}_K(\varsigma)$. In the following is computed the matrix $L_{\mathrm{f}}$ with reference to the situation presented in Figure D.1, the generalization to more complex situations is simple. Suppose that each solid element $[A_t, B_t]$ in the current configuration has intersection with two triangles, then

$$\begin{aligned}
(L_{\mathrm{f}})_{ij} = (\boldsymbol{\zeta}_j, \boldsymbol{\varphi}_i \circ \boldsymbol{\phi}^{n-1})_{0,\Sigma} &= \int_\Sigma \boldsymbol{\zeta}_j \cdot [\boldsymbol{\varphi}_i \circ \boldsymbol{\phi}^{n-1}]\, d\boldsymbol{s} \\
&= \sum_{K\in\mathscr{S}_h} \int_0^{l_K} [\boldsymbol{\zeta}_j \circ \boldsymbol{\alpha}_K(\varsigma)] \cdot [\boldsymbol{\varphi}_i \circ \boldsymbol{\phi}^{n-1} \circ \boldsymbol{\alpha}_K(\varsigma)] \left\| \frac{d\boldsymbol{\alpha}_K}{d\varsigma} \right\|_{\mathbb{R}^d} d\varsigma \\
= \sum_{K\in\mathscr{S}_h} \Bigg( \int_0^{l_K^C} [\boldsymbol{\zeta}_j \circ \boldsymbol{\alpha}_K(\varsigma)] \cdot [\boldsymbol{\varphi}_i \circ \boldsymbol{\alpha}_{Kt}(\varsigma)] &\left\| \frac{d\boldsymbol{\alpha}_K}{d\varsigma} \right\|_{\mathbb{R}^d} d\varsigma + \int_{l_K^C}^{l_K} [\boldsymbol{\zeta}_j \circ \boldsymbol{\alpha}_K(\varsigma)] \cdot [\boldsymbol{\varphi}_i \circ \boldsymbol{\alpha}_{Kt}(\varsigma)] \left\| \frac{d\boldsymbol{\alpha}_K}{d\varsigma} \right\|_{\mathbb{R}^d} d\varsigma \Bigg) \\
&= \sum_{K\in\mathscr{S}_h} \int_{-1}^{1} \left[\boldsymbol{\zeta}_j \circ \boldsymbol{\alpha}_K(\varsigma) \circ \boldsymbol{\psi}_1^K(\hat{\varsigma})\right] \cdot \left[\boldsymbol{\varphi}_i \circ \boldsymbol{\alpha}_{Kt}(\varsigma) \circ \boldsymbol{\psi}_1^K(\hat{\varsigma})\right] \left\| \frac{d\boldsymbol{\alpha}_K}{d\varsigma} \circ \boldsymbol{\psi}_1^K(\hat{\varsigma}) \right\|_{\mathbb{R}^d} \left| \frac{d\boldsymbol{\psi}_1^K}{d\hat{\varsigma}} \right| d\hat{\varsigma} \\
&\quad + \sum_{K\in\mathscr{S}_h} \int_{-1}^{1} \left[\boldsymbol{\zeta}_j \circ \boldsymbol{\alpha}_K \circ \boldsymbol{\psi}_2^K(\hat{\varsigma})\right] \cdot \left[\boldsymbol{\varphi}_i \circ \boldsymbol{\alpha}_{Kt} \circ \boldsymbol{\psi}_2^K(\hat{\varsigma})\right] \left\| \frac{d\boldsymbol{\alpha}_K}{d\varsigma} \circ \boldsymbol{\psi}_2^K(\hat{\varsigma}) \right\|_{\mathbb{R}^d} \left| \frac{d\boldsymbol{\psi}_2^K}{d\hat{\varsigma}} \right| d\hat{\varsigma},
\end{aligned} \tag{D.13}$$

where $\boldsymbol{\psi}_1^K : [-1,1] \to [0, l_K^C]$ and $\boldsymbol{\psi}_2^K : [-1,1] \to [l_K^C, l_K]$ are changes of variable from $[-1,1]$ onto the reference intervals defined as follows. neglecting the dependence from the element $K$, which correspond to consider the case that for all the element $l_K = l$, and using the notation $l_0 = 0, l_1 = l_C, l_2 = l$

$$\boldsymbol{\psi}_k(\hat{\varsigma}) = l_{k-1} + \frac{l_k - l_{k-1}}{2}(1+\hat{\varsigma}) \qquad k = 1,2, \tag{D.14}$$

then

$$\left| \frac{d\boldsymbol{\psi}_k}{d\hat{\varsigma}} \right| = \frac{l_k - l_{k-1}}{2} \tag{D.15}$$

and

$$\left\| \frac{d\boldsymbol{\alpha}_K}{d\varsigma} \circ \boldsymbol{\psi}_k(\hat{\varsigma}) \right\|_{\mathbb{R}^d} = \left\| \frac{1}{l}(\overrightarrow{B} - \overrightarrow{A}) \circ \boldsymbol{\psi}_k(\hat{\varsigma}) \right\|_{\mathbb{R}^d} = \frac{1}{l} \| \overrightarrow{B} - \overrightarrow{A} \|_{\mathbb{R}^d} = 1. \tag{D.16}$$

Using the previous expressions of $\boldsymbol{\alpha}_K$ and $\boldsymbol{\psi}_k$ we have

$$\zeta_j \circ \boldsymbol{\alpha}_K \circ \boldsymbol{\psi}_k(\hat{\varsigma}) = \begin{cases} \left(1 - \frac{1}{l}\left[l_{k-1} + \frac{l_k - l_{k-1}}{2}(1+\hat{\varsigma})\right], 1 - \frac{1}{l}\left[l_{k-1} + \frac{l_k - l_{k-1}}{2}(1+\hat{\varsigma})\right]\right) \\ \\ \left(\frac{1}{l}\left[l_{k-1} + \frac{l_k - L_{k-1}}{2}(1+\hat{\varsigma})\right], \frac{1}{l}\left[l_{k-1} + \frac{l_k - l_{k-1}}{2}(1+\hat{\varsigma})\right]\right) \end{cases} \tag{D.17}$$

In order to compute the integrals we use a Gauss quadrature formula, we can write

$$\int_{-1}^{1} f(\hat{\varsigma})d\hat{\varsigma} \approx \sum_{m=1}^{N} c_m f(\hat{\varsigma}_m).$$

Hence, in our case, for $k = 1,2$

$$(I_k)_{ij} \stackrel{def}{=} \int_{-1}^{1} \left[\boldsymbol{\zeta}_j \circ \boldsymbol{\alpha}_K \circ \boldsymbol{\psi}_k(\hat{\varsigma})\right] \cdot \left[\boldsymbol{\varphi}_i \circ \boldsymbol{\alpha}_{Kt} \circ \boldsymbol{\psi}_k(\hat{\varsigma})\right] \left\| \frac{d\boldsymbol{\alpha}_K}{d\varsigma} \circ \boldsymbol{\psi}_k(\hat{\varsigma}) \right\|_{\mathbb{R}^d} \left| \frac{d\boldsymbol{\psi}_k}{d\hat{\varsigma}} \right| d\hat{\varsigma} \tag{D.18}$$

$$\approx \sum_{m=1}^{N} c_m \left[\boldsymbol{\zeta}_j \circ \boldsymbol{\alpha}_K \circ \boldsymbol{\psi}_k(\hat{\varsigma}_m)\right] \left[\boldsymbol{\varphi}_i \circ \boldsymbol{\alpha}_{Kt} \circ \boldsymbol{\psi}_k(\hat{\varsigma}_m)\right] \frac{l_k - l_{k-1}}{2} \tag{D.19}$$

then, owing to (D.13), we have

$$(L_{\mathrm{f}})_{ij} = (I_1)_{ij} + (I_2)_{ij}$$

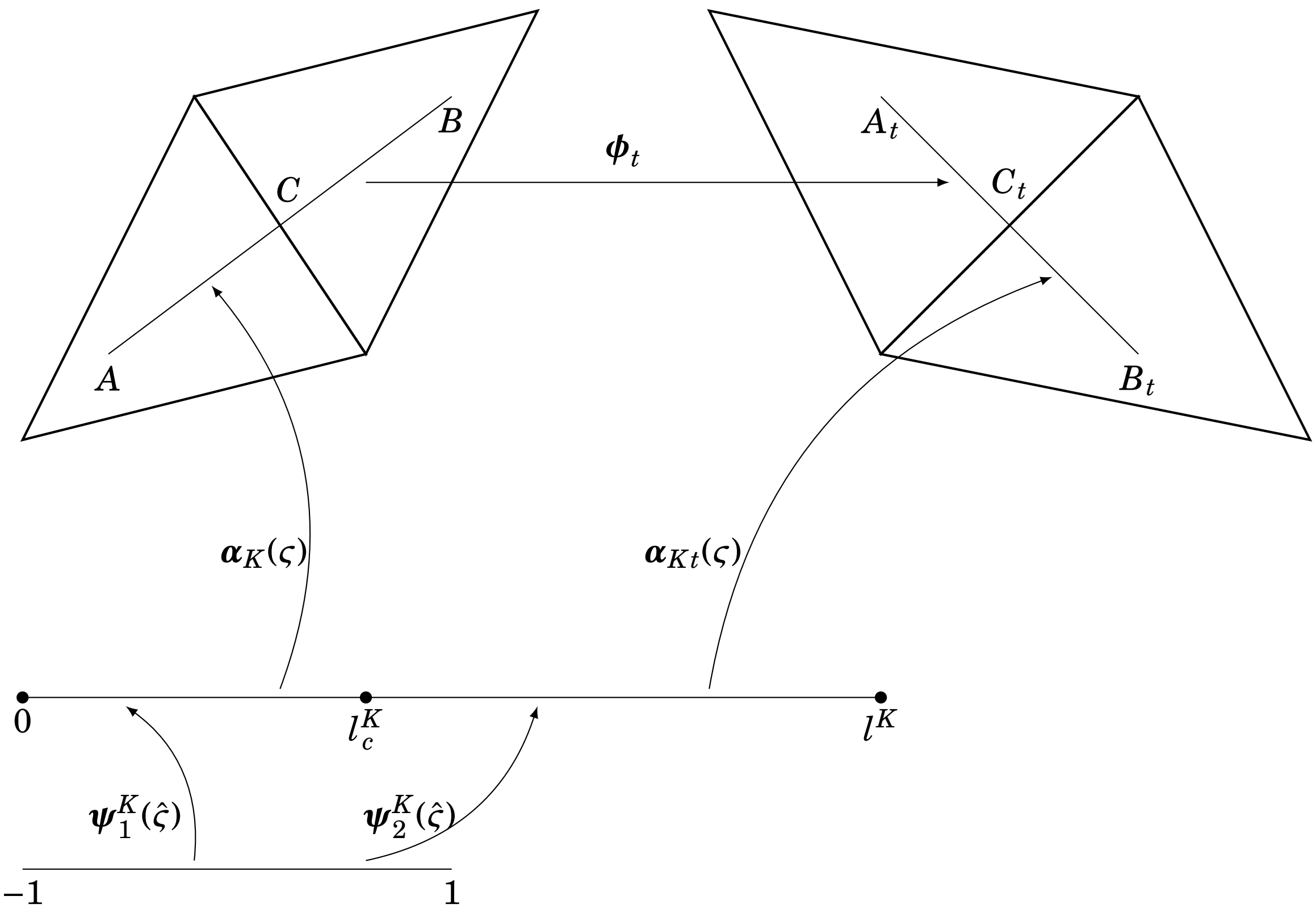

Figure D.1: Element parametrization